\documentclass[]{ucalgmthesis}

\usepackage{amssymb}
\usepackage{amsfonts}
\usepackage{amsthm}
\usepackage{amsmath}
\usepackage{mathtools}
\usepackage{bbm}
\usepackage{bussproofs}
\usepackage{mathrsfs}
\usepackage{tikz, tikz-cd}
\usepackage[all,2cell]{xy}
\usepackage{epigraph}
\usepackage{ stmaryrd } 
\usepackage{todonotes}
\usepackage[style = numeric]{biblatex}
\usepackage{imakeidx}
\tikzset{%
	symbol/.style={%
		draw=none,
		every to/.append style={%
			edge node={node [sloped, allow upside down, auto=false]{$#1$}}}
	}
}

\usetikzlibrary{matrix,arrows}

\newtheorem{Theorem}{Theorem}

\newtheorem{Corollary}[Theorem]{Corollary}
\newtheorem{conjecture}[Theorem]{Conjecture}

\newtheorem{proposition}[Theorem]{Proposition}
\newtheorem{lemma}[Theorem]{Lemma}
\newtheorem{corollary}[Theorem]{Corollary}

\theoremstyle{definition}
\newtheorem{example}[Theorem]{Example}
\newtheorem{remark}[Theorem]{Remark}
\newtheorem{definition}[Theorem]{Definition}
\newtheorem{desiderata}[Theorem]{Desiderata}

\DeclareMathOperator{\Ascr}{\mathscr{A}}

\DeclareMathOperator{\Cscr}{\mathscr{C}}
\DeclareMathOperator{\Dscr}{\mathscr{D}}
\DeclareMathOperator{\Escr}{\mathscr{E}}
\DeclareMathOperator{\Fscr}{\mathscr{F}}
\DeclareMathOperator{\Gscr}{\mathscr{G}}
\DeclareMathOperator{\Hscr}{\mathscr{H}}
\DeclareMathOperator{\Iscr}{\mathscr{I}}

\DeclareMathOperator{\Kscr}{\mathscr{K}}

\DeclareMathOperator{\Pscr}{\mathscr{P}}

\DeclareMathOperator{\Tscr}{\mathscr{T}}
\DeclareMathOperator{\Uscr}{\mathscr{U}}
\DeclareMathOperator{\Vscr}{\mathscr{V}}

\DeclareMathOperator{\Abb}{\mathbb{A}}
\DeclareMathOperator{\Cbb}{\mathbb{C}}
\DeclareMathOperator{\Dbb}{\mathbb{D}}

\DeclareMathOperator{\Gbb}{\mathbb{G}}
\DeclareMathOperator{\Hbb}{\mathbb{H}}

\DeclareMathOperator{\Sbb}{\mathbb{S}}
\DeclareMathOperator{\Tbb}{\mathbb{T}}

\DeclareMathOperator{\Xbb}{\mathbb{X}}

\DeclareMathOperator{\Acal}{\mathcal{A}}
\DeclareMathOperator{\Bcal}{\mathcal{B}}
\DeclareMathOperator{\Ccal}{\mathcal{C}}
\DeclareMathOperator{\Dcal}{\mathcal{D}}
\DeclareMathOperator{\Ecal}{\mathcal{E}}
\DeclareMathOperator{\Fcal}{\mathcal{F}}

\DeclareMathOperator{\Ical}{\mathcal{I}}
\DeclareMathOperator{\Jcal}{\mathcal{J}}

\DeclareMathOperator{\Lcal}{\mathcal{L}}
\DeclareMathOperator{\Mcal}{\mathcal{M}}

\DeclareMathOperator{\Ocal}{\mathcal{O}}
\DeclareMathOperator{\Pcal}{\mathcal{P}}

\DeclareMathOperator{\Scal}{\mathcal{S}}
\DeclareMathOperator{\Tcal}{\mathcal{T}}
\DeclareMathOperator{\Ucal}{\mathcal{U}}
\DeclareMathOperator{\Vcal}{\mathcal{V}}
\DeclareMathOperator{\Wcal}{\mathcal{W}}

\DeclareMathOperator{\afrak}{\mathfrak{a}}

\DeclareMathOperator{\Cfrak}{\mathfrak{C}}

\DeclareMathOperator{\Dfrak}{\mathfrak{D}}

\DeclareMathOperator{\Efrak}{\mathfrak{E}}

\DeclareMathOperator{\mfrak}{\mathfrak{m}}

\DeclareMathOperator{\Xfrak}{\mathfrak{X}}

\DeclareMathOperator{\yon}{\mathbf{y}}

\DeclareMathOperator{\QlVectfd}{\overline{\mathbb{Q}}_{\ell}-\mathbf{Vect}_{f.d.}}
\DeclareMathOperator{\QCoh}{\mathbf{QCoh}}
\DeclareMathOperator{\Shv}{\mathbf{Shv}}

\DeclareMathOperator{\id}{id}

\DeclareMathOperator{\Ker}{Ker}

\DeclareMathOperator{\Set}{\mathbf{Set}}

\DeclareMathOperator{\Topos}{\mathbf{Topos}}
\DeclareMathOperator{\R}{\mathbb{R}}
\DeclareMathOperator{\C}{\mathbb{C}}
\DeclareMathOperator{\Z}{\mathbb{Z}}
\DeclareMathOperator{\N}{\mathbb{N}}

\DeclareMathOperator{\RMod}{\mathbf{R-Mod}}

\DeclareMathOperator{\Ab}{\mathbf{Ab}}

\DeclareMathOperator{\Q}{\mathbb{Q}}

\DeclareMathOperator{\Ext}{Ext}

\DeclareMathOperator*{\colim}{colim}

\DeclareMathOperator{\CalO}{\mathcal{O}}

\DeclareMathOperator{\Spec}{Spec}
\DeclareMathOperator{\Sch}{\mathbf{Sch}}
\DeclareMathOperator{\Dom}{Dom}
\DeclareMathOperator{\FinSet}{\mathbf{FinSet}}

\DeclareMathOperator{\Forget}{Forget}
\DeclareMathOperator{\Cat}{\mathbf{Cat}}

\DeclareMathOperator{\Codom}{Codom}
\DeclareMathOperator{\codom}{Codom}
\DeclareMathOperator{\Iso}{\mathcal{I}so}
\DeclareMathOperator{\A}{\mathbb{A}}

\DeclareMathOperator{\op}{op}

\DeclareMathOperator{\GL}{GL}

\DeclareMathOperator{\SL}{SL}

\DeclareMathOperator{\Per}{\mathbf{Per}}

\DeclareMathOperator{\Sf}{\mathbf{Sf}}
\DeclareMathOperator{\GVar}{\mathnormal{G}-\mathbf{Var}}

\DeclareMathOperator{\OXMod}{\mathcal{O}_{\mathnormal{X}}-\mathbf{Mod}}

\DeclareMathOperator{\Ch}{\mathbf{Ch}}

\DeclareMathOperator{\per}{{}^{\mathnormal{p}}}

\DeclareMathOperator{\VSet}{\mathscr{V}-\mathbf{Set}}
\DeclareMathOperator{\VAb}{\mathscr{V}-\mathbf{Ab}}

\DeclareMathOperator{\fTopos}{\mathfrak{Topos}}
\DeclareMathOperator{\Point}{\mathbf{Point}}
\DeclareMathOperator{\tot}{tot}
\DeclareMathOperator{\Glue}{\mathsf{Glue}}
\DeclareMathOperator{\CAT}{\mathbf{CAT}}
\DeclareMathOperator{\AlgGrp}{\mathbf{AlgGrp}}
\DeclareMathOperator{\dom}{dom}
\DeclareMathOperator{\fCat}{\mathfrak{Katze}}
\DeclareMathOperator{\fCAT}{\mathfrak{KATZE}}
\DeclareMathOperator{\SfResl}{\mathbf{SfResl}}
\DeclareMathOperator{\quo}{\mathsf{quot}}

\DeclareMathOperator{\Loc}{\mathbf{Loc}}
\DeclareMathOperator{\DMod}{\mathcal{D}-\mathbf{Mod}}
\DeclareMathOperator{\DGMod}{\mathcal{D}_{\mathnormal{G}}-\mathbf{Mod}}
\DeclareMathOperator{\eq}{eq}
\DeclareMathOperator{\Eq}{Eq}
\DeclareMathOperator{\Coeq}{Coeq}
\DeclareMathOperator{\coeq}{coeq}
\DeclareMathOperator{\true}{\mathsf{true}}
\DeclareMathOperator{\twoCat}{2\mathbf{Cat}}

\DeclareMathOperator{\GEqCat}{\mathnormal{G}-\mathbf{EquivCat}}

\DeclareMathOperator{\fGEqCat}{\mathnormal{G}-\mathfrak{EquivCat}}

\DeclareMathOperator{\PreEq}{\mathbf{PreEq}}
\DeclareMathOperator{\fPreEq}{\mathfrak{PreEq}}
\DeclareMathOperator{\TShv}{\mathbf{TorShv}}
\DeclareMathOperator{\HVar}{\mathnormal{H}-\mathbf{Var}}
\DeclareMathOperator{\Var}{\mathbf{Var}}
\DeclareMathOperator{\Ad}{Ad}
\DeclareMathOperator{\inv}{inv}
\DeclareMathOperator{\tr}{tr}
\DeclareMathOperator{\Ran}{Ran}
\DeclareMathOperator{\cosk}{cosk}
\DeclareMathOperator{\tottop}{top}
\DeclareMathOperator{\cart}{Cart}	
\DeclareMathOperator{\Open}{\mathbf{Open}}
\DeclareMathOperator*{\twolim}{2lim}
\DeclareMathOperator{\eqstacksimplicial}{\mathnormal{D}_{\text{eq}}^{\mathnormal{b}}(\underline{[X/G]}_{\bullet};\overline{\Q}_{\ell})}

\let\emptyset\varnothing
\let\epsilon\varepsilon

\newcommand{\os}[2]{\overset{#1}{#2}}
\newcommand{\us}[2]{\underset{#1}{#2}}
\newcommand{\ds}[3]{\us{#2}{\os{#1}{#3}}}

\newcommand{\twoColim}[1]{\twolim_{\substack{\longrightarrow \\ {#1}}}}
\newcommand{\Dbeqsimp}[1]{\mathnormal{D}^{\mathnormal{b}}_{\text{eq}}(\underline{G \backslash {#1}}_{\bullet};\overline{\Q}_{\ell})}
\newcommand{\Dbeqstack}[1]{\mathnormal{D}^{\mathnormal{b}}_{\text{eq}}(\underline{[\mathnormal{G} \backslash \mathnormal{X}]}_{\bullet};\overline{\Q}_{\ell})}
\newcommand{\Pro}[1]{{\left[\mathbb{N}^{\op},{#1}\right]}}
\newcommand{\DbQl}[1]{\mathnormal{D}_{\mathnormal{c}}^{\mathnormal{b}}({#1};\overline{\Q}_{\ell})}
\newcommand{\DbeqQl}[1]{\mathnormal{D}_{\mathnormal{G}}^{\mathnormal{b}}({#1};\overline{\Q}_{\ell})}

\DeclareMathOperator{\fX}{\quot{\overline{\mathnormal{f}}}{\mathnormal{X}}}
\DeclareMathOperator{\fY}{\quot{\overline{\mathnormal{f}}}{\mathnormal{Y}}}
\DeclareMathOperator{\gX}{\quot{\overline{\mathnormal{g}}}{\mathnormal{X}}}
\DeclareMathOperator{\gY}{\quot{\overline{\mathnormal{g}}}{\mathnormal{Y}}}
\DeclareMathOperator{\hGamma}{\quot{\overline{\mathnormal{h}}}{\Gamma}}
\DeclareMathOperator{\hGammap}{\quot{\overline{\mathnormal{h}}}{\Gamma^{\prime}}}
\DeclareMathOperator{\hGammapp}{\quot{\overline{\mathnormal{h}}}{\Gamma^{\prime\prime}}}
\DeclareMathOperator{\XGamma}{\quot{\mathnormal{X}}{\Gamma}}
\DeclareMathOperator{\YGamma}{\quot{\mathnormal{Y}}{\Gamma}}
\DeclareMathOperator{\XGammap}{\quot{\mathnormal{X}}{\Gamma^{\prime}}}
\DeclareMathOperator{\YGammap}{\quot{\mathnormal{Y}}{\Gamma^{\prime}}}
\DeclareMathOperator{\XGammapp}{\quot{\mathnormal{X}}{\Gamma^{\prime\prime}}}
\DeclareMathOperator{\YGammapp}{\quot{\mathnormal{Y}}{\Gamma^{\prime\prime}}}
\DeclareMathOperator{\AGamma}{\quot{\mathnormal{A}}{\Gamma}}
\DeclareMathOperator{\AGammap}{\quot{\mathnormal{A}}{\Gamma^{\prime}}}

\DeclareMathOperator{\BGamma}{\quot{\mathnormal{B}}{\Gamma}}
\DeclareMathOperator{\BGammap}{\quot{\mathnormal{B}}{\Gamma^{\prime}}}

\DeclareMathOperator{\OtimesYGamma}{\ds{\Gamma}{Y}{\otimes}}
\DeclareMathOperator{\OtimesXGamma}{\ds{\Gamma}{X}{\otimes}}

\DeclareMathOperator{\const}{\mathsf{const}}
\DeclareMathOperator{\GSch}{\mathnormal{G}-\mathbf{Sch}}
\DeclareMathOperator{\opi}{\overline{\pi}}
\DeclareMathOperator{\of}{\overline{\mathnormal{f}}}
\DeclareMathOperator{\ogamma}{\overline{\gamma}}
\DeclareMathOperator{\oalphaGamma}{\quot{\overline{\alpha}_{\mathnormal{X}}}{\Gamma}}
\DeclareMathOperator{\oalphaGammap}{\quot{\overline{\alpha}_{\mathnormal{X}}}{\Gamma^{\prime}}}
\DeclareMathOperator{\fGX}{\quot{\of}{\mathnormal{G \times X}}}
\DeclareMathOperator{\muGamma}{\mu_{\Gamma}}
\DeclareMathOperator{\muGammap}{\mu_{\Gamma^{\prime}}}
\DeclareMathOperator{\opiGamma}{\quot{\opi_2}{\Gamma}}
\DeclareMathOperator{\opiGammap}{\quot{\opi_2}{\Gamma^{\prime}}}
\DeclareMathOperator{\pGamma}{\overline{\mathnormal{p}}_{\Gamma}}
\DeclareMathOperator{\pGammap}{\overline{\mathnormal{p}}_{\Gamma^{\prime}}}
\DeclareMathOperator{\aGamma}{\overline{\mathnormal{a}}_{\Gamma}}
\DeclareMathOperator{\aGammap}{\overline{\mathnormal{a}}_{\Gamma^{\prime}}}
\DeclareMathOperator{\rhoGamma}{\quot{\rho}{\Gamma}}

\let\emptyset\varnothing
\let\epsilon\varepsilon

\newcommand{\pushoutcorner}[1][dr]{\save*!/#1+1.2pc/#1:(1,-1)@^{|-}\restore} 
\newcommand{\pullbackcorner}[1][dr]{\save*!/#1-1.2pc/#1:(-1,1)@^{|-}\restore}
\newcommand{\quot}[2]{\,_{#2}{#1}}

\newcommand{\ABLDbeqQl}[1]{\quot{\mathnormal{D}_{\mathnormal{G}}^{\mathnormal{b}}({#1};\overline{\Q}_{\ell})}{ABL}}
\newcommand{\odGamma}[2]{\quot{\overline{\mathnormal{d}}_{#1}^{#2}}{{\Gamma}}}

\UseTwocells

\numberwithin{Theorem}{section}

\usepackage[colorlinks, allcolors = blue]{hyperref} 

\usepackage[numbered]{bookmark}

\author{Geoffrey Mark Vooys}

\title{Equivariant Functors and Sheaves}

\degree{Doctor of Philosophy}

\prog{Graduate Program in Mathematics and Statistics}

\monthname{AUGUST}

\thesisyear{2021}

\hypersetup{pdfinfo={Title={\thetitle},Author={\theauthor}}}

\makeindex[title = Index of Terminology, intoc] 
\makeindex[name = notation, title = Index of Notation, intoc] 

\begin{document}

\frontmatter



\makethesistitle



\chapter{Abstract}

In this thesis we study two main topics which culminate in a proof that four distinct definitions of the equivariant derived category of a smooth algebraic group $G$ acting on a variety $X$ are in fact equivalent. In the first part of this thesis we introduce and study equivariant categories on a quasi-projective variety $X$. These are a generalization of the equivariant derived category of Lusztig and are indexed by certain pseudofunctors that take values in the $2$-category of categories. This $2$-categorical generalization allow us to prove rigorously and carefully when such categories are additive, monoidal, triangulated, and admit $t$-structures, among other properties. We also define equivariant functors and natural transformations before using these to prove how to lift adjoints to the equivariant setting. We also give a careful foundation of how to manipulate $t$-structures on these equivariant categories for future use and with an eye towards future applications.

In the final part of this thesis we prove a four-way equivalence between the different formulations of the equivariant derived category of $\ell$-adic sheaves on a quasi-projective variety $X$. We show that the equivariant derived category of Lusztig is equivalent to the equivariant derived category of Bernstein-Lunts and the simplicial equivariant derived category. As a final step we show that these equivariant derived categories are equivalent to the derived $\ell$-adic category on the algebraic stack $[X/G]$ of Behrend. We also provide an isomorphism of the simplicial equivariant derived category on the variety $X$ with the simplicial equivariant derived category on the simplicial presentation of $[X/G]$, as well as prove explicit equivalences between the categories of equivariant $\ell$-adic sheaves, local systems, and perverse sheaves with the classical incarnations of such categories of equivariant sheaves.
%

\chapter{Preface}

\begin{quote}
	
	Chapters \ref{Chapter Intro EqCats}, \ref{Section: Equivariant Categories}, \ref{Chapter 4}, \ref{Chapter 5}, \ref{Chapter Intro Comparison}, \ref{Section: Lusztig EDC is equiva to ABL EDC}, \ref{Section: Lusztig EDC is the Simplicial EDC}, and \ref{Section: Simplicial EDC is the stcky EDC} of this thesis are original, unpublished, independent work by the author,
	Geoffrey M.\@ Vooys.
\end{quote}

\chapter{Acknowledgments}
I would like to first and foremost thank Clifton Cunningham for his support and helpful guidance as a supervisor, mentor, and friend throughout both my undergrad and PhD careers. I am especially grateful for the ability to find my own approach to research mathematics and the help Clifton has given me in both writing the abstract nonsense that I tend towards in a way that other people can understand, as well as helping me learn the ever important lesson of knowing my audience. I'm also grateful for Clifton for making math both funny and engaging, and of the most important lessons I've learned from him is that a sense of humor is invaluable (together with the fact that puns are the highest form of humor --- no contest).

I would like to thank Kristine Bauer for all the support and advice she's given me over the years. Her advice is always helpful and insightful, and the resources that she has guided me  towards have been invaluable in my research and for learning different perspectives on things.

I would like to thank Berndt Brenken for all his help, interesting discussions, wisdom, and friendship over the years. I learned a lot about teaching and grading both pure mathematics (while TAing his abstract algebra class) and applied mathematics (while grading his math for economics majors course); the algebra class was one of the best and most fun classes I've ever been a part of.

I would like to thank Larry Bates for all his guidance and the many interesting and illuminating conversations about life, geometry, music, and the universe over the years. I learned more about complex geometry, writing mathematics, and communicating mathematical ideas from conversations and courses taught by Larry than anywhere else. 

I would like to thank the $p$-attic and the Voganish group (especially James Steele, Alex Cameron, Andr{\'e} de Waal, Nicole Kitt, Qing Zhang, Sarah Dijols, Jerrod Smith, and Andrew Fiori) for attending the various talks I've given over the years, listening to me rant about (higher) categories and their use in arithmetic geometry, entertaining the questions I've asked (even when they were technical), and generally helping make my graduate student experience fantastic.

I would like to thank the MATRIX Institute for their hospitality during the conference ``Geometric and Categorical Representation Theory.'' Many of the results of the paper \cite{VooysGreenberg} were discovered during the intellectually stimulating environment of this two week conference; while this paper is not a part of this thesis, it was done during my PhD studies and I am grateful all the same.

I would like to thank the Banff International Research Station for their hospitality during the 2019 focused research group ``The Voganish Project.'' The work that began Chapter \ref{Chapter: EQCats} was born at this event, and many preliminary results involving triangulated equivariant categories were discovered after sessions that took place during this week.

I would like to thank Robin Cockett and the Peripatetic group at the University of Calgary for accepting me as a token categorical algebraic geometer. I learned a lot from the talks I attended (especially about how useful essentially algebraic theories, type theory, and fibrations are) and I'm very grateful for the opportunities to learn I've had over the years by being a part of the peripatetic group.

%

Finally, since this thesis was typeset using \LaTeX{} using the package \\\verb|ucalgarymthesis|, itself based on the \verb|memoir| class, I would like to thank
thank Don Knuth, Leslie Lamport, Peter Wilson, and Richard Zach.

\dedication{To Heather for always being there for me, my parents for being so wonderful and supportive of me throughout my career, my brother, and finally my pugs Petunia (she's a beautiful flower) and Cauchy (he leaves residues on poles).}

\tableofcontents






\mainmatter


\epigraph{La mer s'avance insensiblement et sans bruit, rien ne semble se casser rien ne bouge l'eau est si loin on l'entend {\`a} piene... Pourtant elle finit par entourer la substance r{\'e}tive celle-ci peu {\`a} peu devient une presqu'{\^\i}le, puis un {\^\i}lot, quit finit par {\^e}tre submerg{\'e} {\`a} son tour, comme si'il s'{\'e}tait finalement dissous dan l'oc{\'e}an s'{\'e}tendant {\`a} perte de vue..}{Alexandre Grothendieck, \emph{R{\'e}coltes et Semailles}}


\chapter{Introduction}


According to Alexandre Grothendieck in \cite{GorthendieckRecolte}, there are two approaches to performing mathematics and trying to learn the unknown. The original quote, as translated to English in \cite{mclarty2003}, paints the metaphor as someone attempting to get into a sealed walnut. The first approach takes a hammer and chisel and starts to smash away; they can reorient their hammer or chisel as much as needed, but they proceed by largely crashing through the walnut to get to the tantalizing centre of the walnut. This is the quickest and most direct way to get into the walnut, but it is provides the fewest tools for getting into other walnuts. If presented with a harder or differently shaped nut, there is no guarantee that the same methodology will allow us to crack the new nut.

The second approach to getting into the walnut of unknown knowledge is the approach that Grothendieck himself practiced. The trick here is to find and produce some sort of liquid in which to submerge the walnut. After letting it soak for long enough, adding more of the liquid if necessary, eventually the shell becomes supple and pops open. The benefit of this method is that the liquid may be used to soak other walnuts. While it may take considerable time and significant amounts of the liquid may need to be added to do a proper soaking, eventually the new walnut will part. Like a beachfront amidst the ocean, eventually the sea rises and covers the sand; even if the march of the sea seemed slow, perhaps even imperceptible, the sea eventually swallows and submerges the piece of land.

In my opinion, this thesis is much more of a giant bottle of liquid, a wave cresting on the beach, as opposed to a hammer with which to smash. If you are reading this thesis you are likely an examiner (because who else reads theses) and there is a good chance that you are holding a physical copy of the thing. In this case I give my apologies, both for writing such a long document of daunting size and weight and for taking the rising sea approach. This thesis is a marathon and not a sprint, but please let me be your guide through this too long swim.

This thesis is motivated by a study of equivariant derived categories (of $\ell$-adic or {\'e}tale sheaves) on a variety $X$. These are important categories in representation theory and algebraic geometry, but are much more difficult to understand and come across than someone might hope. The study of equivaraint sheaves becan in earnest in \cite{GIT}, when Mumford and Fogarty established the correct definition for a sheaf to be equivariant with respect to a group scheme $G$ acting on a scheme $X$. Their definition said roughly that a sheaf $\Fscr$ on a scheme $X$ is equivariant with respect to a group scheme $G$ acting on $X$ if there was an isomorphism, for group action $\alpha_X:G \times X \to X$ and projection $\pi_2:G \times X \to X$,
\[
\theta:\alpha_X^{\ast}\Fscr \to \pi_2^{\ast}\Fscr
\]
of sheaves on $G \times X$. These isomorphisms had to satisfy a certain cocycle condition that was equivalent to being compatible with the group action of $G$ on $X$ (cf.\@ Definition \ref{Defn: GIT Cocycle}), but this definition at least implies that for any point $(g,x) \in G \times X$ we have an isomorphism fitting into a commuting diagram
\[
\xymatrix{
\alpha_X^{\ast}\Fscr_{(g,x)} \ar[d]_{\cong} \ar[r]^-{\theta_{(g,x)}} & \pi_2^{\ast}\Fscr_{(g,x)} \ar[d]^{\cong} \\
\Fscr_{gx} \ar@{-->}[r]_-{\exists \cong} & \Fscr_{x}
}
\]
so the sheaf $\Fscr$ is suitably stable along the $G$-orbits of points of $X$. 

While this definition of equivariant sheaves is straightforward to understand and work with, it is poorly behaved with respect to the categorical techniques we use when studying cohomology of sheaves and the general six functor formalism that comes with derived algebraic geometry. In particular, in Chapter \ref{Chapter Intro Comparison} we give an explicit example due to Torsten Ekedahl (cf.\@ \cite{Torsten}) of a situation in which the collection of compelexes $A \in D^b_c(X)_0$ with an isomorphism $\theta:\alpha_X^{\ast} A \to \pi_2^{\ast}A$ as above need not even be a triangulated category. Let us spell out explicitly what we want out of an equivariant ($\ell$-adic) derived category:

\begin{desiderata}[{\cite[Desideratum 6.2.1]{PramodBook}}]\label{Desiderata}
For every smooth algebraic group $G$ over a field $K$ and for every quasi-projective left $G$-variety $X$, there should be a category $\DbeqQl{X}$ with the following properties:
\begin{itemize}
	\item The category $\DbeqQl{X}$ needs to be triangulated;
	\item There exist $t$-structures $t_{\text{stand}}$ and $t_{\text{per}}$ on $\DbeqQl{X}$ which induce equivalences of categories
	\[
	{}^{\text{stand}}\DbQl{\XGamma}^{\heartsuit} \simeq \Shv_G(X;\overline{\Q}_{\ell})
	\]
	and
	\[
	{}^{p}\DbeqQl{X}^{\heartsuit} \simeq \Per_G(X;\overline{\Q}_{\ell}).
	\]
	\item There are ``forgetful'' functors $\Forget:\DbeqQl{X} \to \DbQl{X}$,\\ ${}^{\text{stand}}\Forget:\Shv_G(X;\overline{\Q}_{\ell}) \to \Shv(X;\overline{\Q}_{\ell})$, and\\ ${}^{p}\Forget:\Per_G(X;\overline{\Q}_{\ell}) \to \Per(X;\overline{\Q}_{\ell})$ and the forgetful functor on the equivariant derived category restricts to either of these in the sense that the diagrams
	\[
	\begin{tikzcd}
	\DbeqQl{X} \ar[r, ""{name = U}]{}{{}^{\text{stand}}H^0} \ar[d, swap]{}{\Forget} & \Shv_G(X;\overline{\Q}_{\ell}) \ar[d]{}{{}^{\text{stand}}\Forget} \\
	\DbQl{X} \ar[r, swap, ""{name = L}]{}{H^0} & \Shv(X;\overline{\Q}_{\ell}) \ar[from = U, to = L, Rightarrow, shorten <= 4pt, shorten >= 4pt]{}{\cong}
	\end{tikzcd}	
	\]
	and
	\[
	\begin{tikzcd}
	\DbeqQl{X} \ar[r, ""{name = U}]{}{{}^{p}H^0} \ar[d, swap]{}{\Forget} & \Per_G(X;\overline{\Q}_{\ell}) \ar[d]{}{{}^{p}\Forget} \\
	\DbQl{X} \ar[r, swap, ""{name = L}]{}{{}^{p}H^0} & \Per(X;\overline{\Q}_{\ell}) \ar[from = U, to = L, Rightarrow, shorten <= 4pt, shorten >= 4pt]{}{\cong}
	\end{tikzcd}
	\]
	commute up to natural isomorphism.
	\item Given a $G$-equivariant morphism $h:X \to Y$ between quasi-projective left $G$-varieties, there are adjoints
	\[
	\begin{tikzcd}
	\DbeqQl{X} \ar[r, swap, bend right = 30, ""{name = L}]{}{Rh_{\ast}} & \DbeqQl{Y} \ar[l, bend right =30, swap, ""{name = U}]{}{h^{\ast}} \ar[from = U, to = L, symbol = \dashv]
	\end{tikzcd}
	\]
	\[
	\begin{tikzcd}
		\DbeqQl{X} \ar[r, bend left = 30, ""{name = L}]{}{Rh_{!}} & \DbeqQl{Y} \ar[l, bend left =30, ""{name = U}]{}{h^{!}} \ar[from = L, to = U, symbol = \dashv]
	\end{tikzcd}
	\]
	and, for any $A \in \DbeqQl{X}_0$,
	\[
	\begin{tikzcd}
	\DbeqQl{X} \ar[r, bend left =30, ""{name = U}]{}{A \os{L}{\otimes} (-)} & \DbeqQl{X} \ar[l, bend left = 30, ""{name =L}]{}{R[A,-]} \ar[from = U, to = L, symbol = \dashv]
	\end{tikzcd}
	\]
	which all commute with the forgetful functors up to isomorphism.
	\item For any $G$-equivariant open immersion $j:U \to X$ with closed complement $i:V \to X$ and for any object $A \in \DbeqQl{X}_0$ there is a distinguished triangle of the form
	\[
	\xymatrix{
		Rj_!\left(j^{!}A\right) \ar[r]^-{\epsilon_A} & A \ar[r]^-{\eta_A} & Ri_{\ast}\left(i^{\ast}A\right) \ar[r] & Rj_!\left(j^!A\right)[1]
	}
	\]
	where $\epsilon:Rj_{!} \circ j^{!} \Rightarrow \id_{\DbeqQl{X}}$ is the counit of adjunction and $\eta:\id_{\DbeqQl{X}} \Rightarrow Ri_{\ast} \circ i^{\ast}$ is the unit of adjunction.
	\item For any $G$-equivariant morphism $h:X\to Y$ and any objects $A \in \DbeqQl{X}_0, B \in \DbeqQl{Y}$ there are natural isomorphisms
	\begin{align*}
		\big(Rh_!A\big) \os{L}{\otimes}_Y B &\cong Rh_!\left(A \os{L}{\otimes}_X h^{\ast}B\right), \\
		R[Rh_!A,B]_Y &\cong Rh_{\ast}\left(R[A,h^{!}B]_X\right), \\
		h^{!}R[A,B]_Y &\cong R[h^{\ast}A,h^{!}B]_X.
	\end{align*}
	
\end{itemize}
\end{desiderata}

In light of Desiderata \ref{Desiderata}, it may seem like trying to realize $\DbQl{X}$ as a subcategory of $D^b_c(X;\overline{\Q}_{\ell})$ is a little na{\"i}ve and using a construction like $D^b(\Per_G(X;\overline{\Q}_{\ell}))$ might work. However, the construction of $\DbeqQl{X}$ resists even this approach. In \cite[Example 6.2.2]{PramodBook} it is shown that for the algebraic group $G = \SL(2, \C)$ and the varieties $X = \mathbb{P}_{\C}^1$ and $Y = \Spec \C$ the categories $D^b(\Per_{\SL(2,\C)}(-))$ fail to satisfy Desiderata \ref{Desiderata}. As such, a more delicate approach towards the equivariant derived category was required, which lead to the competing formats that were studied in \cite{DeligneHodge3} (the simplicial approach); \cite{PramodBook} and \cite{BernLun} (the Bernstein-Lunts approach); \cite{LusztigCuspidal2} (Lusztig's approach for use with graded Hecke algebras); and \cite{Behrend} (the stacky approach). We begin our examination of these categories in Part \ref{Chapter: EQCats} by setting the base camp of our journey at a generalization of Lusztig's equivariant derived category as it appears in \cite{LusztigCuspidal2}.

Our first step in the journey is Part \ref{Chapter: EQCats}, which is arguably the best avatar of the rising sea of Grothendieck within this thesis. In this chapter we define and systematically study equivariant categories. In this chapter we work with smooth algebraic groups $G$ and left $G$-varieties $X$. These equivariant categories, $F_G(X)$, are indexed by pseudofunctors
\[
F:\SfResl_G(X)^{\op} \to \fCat
\]
where $\SfResl_G(X)$ is a category of smooth free $G$-resolutions of $X$ and $\fCat$ is the $2$-category of categories. We develop and study the theory of these categories in order to have tools and methods to carefully and rigorously study the equivariant derived category of Lusztig as defined in \cite{LusztigCuspidal2}. In the beginning of this chapter we provide some basic definitions and sanity results that say, among other things, that the equivariant category over the trivial algebraic group $1$,  $F_{1}(X)$, is equivalent to the basic category $F(X)$. Additionally we construct a forgetful functor $F_G(X) \to F(X)$ and an equivariantification functor $F(X) \to F_G(X)$. After this we get into studying the category theory of $F_G(X)$ and give conditions in which $F_G(X)$ is additive, admits (co)limits of specified shape, and is (symmetric) monoidal. Once we have seen the monoidal structure on $F_G(X)$ we begin a study of equivariant functors between equivariant categories. We then prove that equivariant functors of the form $F_G(X) \to E_G(X)$ arise from pseudonatural transformations
\[
\begin{tikzcd}
\SfResl_G(X)^{\op} \ar[r, bend left = 30, ""{name = U}]{}{F} \ar[r, bend right = 30, swap, ""{name = L}]{}{E} & \fCat \ar[from = U, to = L, Rightarrow, shorten <= 4pt, shorten >= 4pt]{}{\alpha}
\end{tikzcd}
\] 
and use these pseudonatural transformations to give natural conditions when equivariant functors $\underline{\alpha}:F_G(X) \to E_G(X)$ are (co)complete, additive, and monoidal. After introducing equivariant functors, we show how natural transformations between equivariant functors arise from modifications $\eta:\alpha \to \beta$ between pseudonatural transformations. This is a crucial step for proving that we can lift adjoints along resolutions to adjoints between equivariant categories. 

Once we have introduced modifications and adjoints, we move to discuss and define equivariant functors of the form $F_G(X) \to F_G(Y)$. We give definitions for cases in which these functors exist and then prove when they are additive, preserve (co)limits, and are monoidal. This section of the chapter closes with a quick study of when changing fibres $F_G(X) \to E_G(X)$ and then changing spaces $E_G(X) \to E_G(Y)$ is naturally isomorphic to changing spaces $F_G(X) \to F_G(Y)$ and then changing fibres $F_G(Y) \to E_G(Y)$.

After studying the interchange of equivariant functors, we introduce equivariant functors of the form $F_G(X) \to F_H(X)$ when there is a morphism of smooth algebraic groups $H \to G$. We work to prove that they exist and then go on to prove how they interact with the other equivariant functors. As particular sanity checks, we show that changing groups $F_G(X) \to F_H(X)$ and $F_H(X) \to F_{L}(X)$ along two changes of groups induced by $H \to G$ and $L \to H$ is naturally isomorphic to the change of groups induced by the composite $L \to H \to G$ and that the forgetful functor $F_G(X) \to F(X)$ is naturally isomorphic to the change of groups $F_G(X) \to F_1(X)$ induced by the inclusion of the identity of $G$.

Finally, after studying the equivariant functors which change groups, we study when equivariant categories admit a triangulation and some of the $t$-structures that arise on the categories $F_G(X)$. In this section we are able to show that under mild commutativity assumptions that if $F_G(X)$ admits a $t$-structure and if each fibre category of $F$ over the resolution category admits a heart of a $t$-structure, then there is an equality of categories $(F_G(X))^{\heartsuit} \simeq (F^{\heartsuit})_G(X)$ that I like to think of as a change of heart equality. 

After getting to meet $t$-structures of equivariant categories, the next step in our journey is Part \ref{Chapter: EDC Comp}. In this part of the thesis we show that the equivariant category $\DbeqQl{X}$ of Part \ref{Chapter: EQCats} satisfies Desiderata \ref{Desiderata}. After this we prove a four-way equivalence of categories
\[
\DbeqQl{X} \simeq \quot{\DbeqQl{X}}{ABL} \simeq \Dbeqsimp{X} \simeq \DbQl{[G\backslash X]}
\]
for quasi-projective varieties $X$ and smooth affine algebraic groups $G$ acting on $X$; note that $\DbeqQl{X}$ is the equivariant category as it appears in Part \ref{Chapter: EQCats}, $\quot{\DbeqQl{X}}{ABL}$ is the equivariant derived category of Achar-Bernstein-Lunts (as developed by Bernstein-Lunts in \cite{BernLun} and extended to schemes by Achar in \cite{PramodBook}), $\Dbeqsimp{X}$ is the simplicial equivariant derived category of Deligne (cf.\@ \cite{DeligneHodge3}), and $\DbQl{[G \backslash X]}$ is the $\ell$-adic category as defined in \cite{Behrend}. This chapter has a focus on explicit computation and focuses to a large degree in illustrating how computations and manipulations with the simplicial equivariant derived category $\Dbeqsimp{X}$ work. However, we begin by proving that $\DbeqQl{X} \simeq \quot{\DbeqQl{X}}{ABL}$. We also develop an isomorphism of simplicial equivariant derived categories
\[
\Dbeqsimp{X} \cong \Dbeqstack{X},
\]
where $\underline{[G \backslash X]}_{\bullet}$ is a simplicial approximation of $[G\backslash X]$, by proving an isomorphism of $2$-coskeletal simplicial schemes
\[
\underline{G \backslash X}_{\bullet} \cong \underline{[G \backslash X]}_{\bullet}.
\]
These techniques allow us to prove the equivalence of categories 
\[
\Dbeqsimp{X} \simeq \DbQl{[G \backslash X]}
\] 
by changing to the stacky simplicial category $\Dbeqstack{X}$ and then proving that $\Dbeqstack{X} \simeq \DbQl{[G \backslash X]}$.

Finally, a quick word of the foundations for this thesis. We will frequently talk about left $G$-varieties and the category of left $G$-varieties for an algebraic group $G$. These are defined explicitly later on, but it is helpful to keep in mind that the category $\GVar$ of left $G$-varieties is the category of varieties that carry a left $G$-action and whose morphisms are all $G$-equivariant morphisms of varieties.

Now that my task as map guide is complete, let me set the mood for going forward with a quick joke. What do you call someone who reads a category theory paper? A coauthor! Please allow me to take us through this thesis then as my coauthor as we begin our journey with the small mountain of equivariant categories.

\part{Equivariant Functors and Sheaves}\label{Chapter: EQCats}
\chapter{Introduction}\label{Chapter Intro EqCats}
The equivariant derived category is an important object in arithmetic geometry and representation theory, as it provides the natural location in which the theories of equivariant ($\ell$-adic) sheaves, equivariant ($\ell$-adic) local systems, and equivariant ($\ell$-adic) perverse sheaves (cf.\@ \cite{LusztigIntersectionCohomComp}, \cite{LusztigCharacter1}, \cite{AcharEtAl}, \cite{CFMMX}, \cite{Juteauetal}, \cite{TensorProdsPervSheaves}, \cite{BernLun}, \cite{LusztigCuspidal2}). However, while the notion of equivariance is straightforward to define and follow the definition as given by Mumford in \cite{GIT} (namely that an object $\Fscr$ in a category $\Cscr(X)$ of sheaves of some sort on $X$ is equivariant if, for an action morphism $\alpha_X:G \times X \to X$ and projection $\pi_2:G \times X \to X$ for an algebraic group $G$ and a left $G$-variety $X$, there is an isomorphism $\theta$ from the pullback of $\Fscr$ along $\alpha_X$ and the pullback of $\Fscr$ along the projection map in $\Cscr(G \times X)$), but the notion of equivariance is much more subtle and difficult to extend to the full derived category level. In fact, following an example of Ekedahl \cite{Torsten}, the category of complexes in $D_c^b(\ast)$ with an equivariance in the style of Mumford can even fail to be triangulated.

The notion of the equivariant derived category was eventually figured out by Bernstein-Lunts in \cite{BernLun} and also, in an alternative fashion, by Lusztig in \cite{LusztigCuspidal2}. While the approach of Bernstein-Lunts has been used much more in practice (cf.\@ \cite{Geordie}, \cite{PramodBook}, \cite{CharacterSheavesUnipotentGroupsPositiveChar}), the equivariant derived category introduced in \cite{LusztigCuspidal2}, $D_G^b(X)$, has been used in various places. It was developed to use a type of equivariant cohomology to construct the graded Hecke algebra; this equivariant derived category and its connections to graded Hecke algebras has been used in, for example, \cite{PerverseSheavesLoop}, \cite{RepGradedHecke}, and crucially in \cite{CFMMX} to give an equivariant notion of Brylinski's $\mathsf{Ev}$ functor and equivariant nearby and vanishing cycles. However, the equivariant bounded derived category of Lusztig has had few of its formal properties proved and treated explicitly. At this point people must take the properties of this category largely on faith. Notably, it is unclear that this equivariant derived category is triangulated, even though this is suggested by the fact that it is a natural home for the equivariant cohomology giving rise to the graded Hecke algebra.

The main goal of this thesis started with a desire to rigourously and systematically establish and examine the desired properties of the equivariant derived category of \cite{LusztigCuspidal2}. However, in going about this it became clear that the equivariant derived category of \cite{LusztigCuspidal2} was the tip of an iceberg of pseudofunctorial equivariant descent. While the equivariant derived category is a natural and important object, it is conceivable that it is reasonable and useful to consider an equivariant category on a scheme or variety $X$ which is not given by equivariant descent through derived categories or hearts of various $t$-strucutres. For example, it may be of interest to study the category of (sheaves) of $G$-equivariant dg-algebras over in the ringed site $(\Shv(\Open(\lvert X\rvert)),\Ocal_X)$; here we will run into a very similar problem as in the equivariant derived categorical case.

A second goal for this work began when trying to understand what it means to have equivariant functors and equivariant adjunctions between any of the above well-known equivariant categories. For instance, even when trying to produce functors $\Per_G(X) \to \Shv_G(X)$, $\Loc_G(X) \to \Shv_G(X)$, or $\Shv_G(X) \to \Shv_G(Y)$ it quickly becomes difficult to ascertain if these functors arise as equivariant versions of ``global'' functors and if these global functors restrict appropriately. To see an explicit example of this, let $G$ be a smooth connected reductive group over $\Spec K$ (for $K$ an arbitrary field) and let $X$ be an arbitrary $G$-variety. There is then a non-equivariant isomorphism of varieties
\[
\psi:X \xrightarrow{\cong} G\backslash(G \times X)
\]
which gives rise to an equivalence of categories $\Shv(X) \simeq \Shv(G \backslash (G \times X))$. However, when trying to find equivariant functors
\[
\xymatrix{
\Shv_G(X) \ar@<1ex>[rr]^-{(\psi^{-1})^{\ast}}&  & \Shv_G(G \backslash (G \times X)) \ar@<1ex>[ll]^-{\psi^{\ast}}
}
\]
in the na{\"i}ve sense, the equivariance situation becomes significantly more subtle than it may first appear. These subtle issues led to a general principle that (combined with the notions in \cite{PramodBook}, \cite{BernLun}, and \cite{LusztigCuspidal2} that the equivariant derived category is controlled by descent through acyclic $G$-equivariant resolutions of $X$) it is simultaneously the descent theory and the exact categories through which we descend that control exactly how the adjoints and functors that appear between equivariant categories appear and behave.

When making the observation above that we can study equivariant functors, adjoints, and behaviours by studying acyclic descent through the specific categories at hand, we are naturally lead to a pseudofunctorial and $2$-categorical language so that we can carefully and precisely explain exactly how equivariant categories, functors, and adjoints arise. These lead to a notion of equivariant category $F_G(X)$ as indexed by a pseudofunctor $F$ defined on a resolution category over $X$ landing in the $2$-category of categories, $\fCat$ (cf.\@ Definition \ref{Defn: Equivariant Cats}). This paper is devoted to the study of these categories and what properties arise on $F_G(X)$ based on the nature of the pseudofunctor itself; perhaps unsurprisingly, the structure of categories like $D^b_G(X), \DbeqQl{X}, \Per_G(X), \Per_G(X;\overline{\Q}_{\ell}), \Shv_G(X;\overline{\Q}_{\ell})$, and others arise from the properties of the various pseudofunctors that induce each given category. In fact, Example \ref{Example: Section 2: Big List of Examples} gives a large (certainly non-exhaustive!) list of examples of various pseudofunctors that give rise to equivariant categories of interest to algebraic geometers, category theorists, and representation theorists alike.

\section{Summary of Key Results and Structure of Part I}
To see the flexiblity and strength of the formalism introduced in this paper, we will describe some of the key results of this paper here (while pointing to the corresponding result as it appears in context in the work simultaneously) to illustrate exactly how this structure is provided and induced. For this, however, it is worth pointing out some language/jargon that will help us navigate the results themselves. The type of pseudofunctor $F$ defined on our resolution category $R_G(X)$ of $X$ that gives rise to an equivariant category is called a \textit{pre-equivariant pseudofunctor on $X$} (cf.\@ Definition \ref{Defn: Prequivariant pseudofunctors}).

In Section \ref{Section: Preliminary and Background for EQCats} we review the basic construction of Lusztig from \cite{LusztigCuspidal2} by introducing the category $\Sf(G)$ of smooth free $G$-varieties of pure dimension. The only big result in this section is establishing that for a Grothendieck universe $\Vscr$ for which the categories $\DbQl{X}$ are all $\Vscr$-small for varieties $X$ over $\Spec K$, the equivariant category $F_G(X)$ is $\Vscr$-small as well. While those who do not worry about foundations will find this unnecessary detail, this is a necessary step for showing that $F_G(X)$ is additive in certain cases (cf.\@ Corollary \ref{Cor: Section 2: Additive pseudofunctor gives additive equivariat cat}), as we need to know that the hom-objects are enriched in $\Ab$ for whatever universe we have at hand.

Section \ref{Section: Definitions and Stuff} defines and constructs the equivariant category (cf.\@ Definition \ref{Defn: Equivariant Cats}) and prove some basic sanity results. In particular, we show that for the trivial group $1 = \Spec K$ there is always an equivalence of categories $F_{\Spec K}(X) \simeq F(X)$ (cf.\@ Proposition \ref{Prop: Section 2: Equivalence of K equivariant cat with F X}). For any smooth algebraic group $G$ and left $G$-variety $X$ we also provide functors $\Forget:F_G(X) \to F(X)$ and $\Eq:F(X) \to F_G(X)$ which make $\Forget \circ \Eq \cong \id_{F(X)}$.

In Section \ref{Section 2: Cat theory of ECat} we begin a systematic study of the category theory of the equivariant category $F_G(X)$. The first main theorems of this section are Theorems \ref{Thm: Section 2: Equivariant Cat has lims} and \ref{Theorem: Section 2: Equivariant Cat has colims}, which give straightforward conditions in which we can guarantee that $F_G(X)$ has limits of shape $I$ whenever each fibre category and fibre functor of the pre-equivariant pseudofunctor $F$ have and preserve all (co)limits of shape $I$. This allows us to deduce as corollaries when $F_G(X)$ is complete and cocomplete, and together with a foundational lemma when $F_G(X)$ is additive. Afterwards we study how to extend regularity conditions to $F_G(X)$ and deduce in Proposition \ref{Prop: Section 2: Abelian Descends} when the equivariant category $F_G(X)$ is Abelian. The last main result of this section is Theorem \ref{Theorem: Section 2: Monoidal preequivariant pseudofunctor gives monoidal equivariant cat}, which states that so long as all the fibre categories of the pseudofunctor are monoidal and all the fibre functors are strong monoidal, the equivariant category $F_G(X)$ is monoidal as well with monoidal functor induced from the monoidal functors in each fibre category. From this we deduce that whenever each of the fibre categories of $F$ are symmetric monoidal and the fibre functors are strong monoidal, then so is $F_G(X)$ (cf.\@ Proposition \ref{Prop: Section 2: Equivariant cat is symmetric monoidal}).

Section \ref{Section: Section 3: Change of Fibre} begins our study of equivariant functors (functors between equivariant categories induced by equivariant descent data) by studying what we call Change of Fibre functors: Functors of the form $F_G(X) \to E_G(X)$ for pre-equivariant pseudofunctors $E$ and $F$ on $X$. The first main result of this section is Theorem \ref{Thm: Section 3: Psuedonatural trans lift to equivariant functors}, which says that equivariant functors of the form $F_G(X) \to E_G(X)$ come exactly from pseudonatural transformations $\alpha:F \Rightarrow E$ of pseudofunctors. This then gives rise to Proposition \ref{Prop: Section 3: Change of fibre equivariant functor preserves limit of a specific shape}, which says that if we have pseudofunctors $E$ and $F$ which satisfy the hypotheses of Theorem \ref{Thm: Section 2: Equivariant Cat has lims} so that both $F_G(X)$ and $E_G(X)$ have limits of shape $I$ and if each translation functor of the pseudonatural transformation $\alpha:F \Rightarrow E$ preserves the limits of shape $I$ in each fibre category, then the functor $\alpha:F_G(X) \to E_G(X)$ induced by $\alpha$ preserves limits of shape $I$ as well. After showing that pseudonatural transformations give rise to equivariant functors, we show that modifications $\eta$ between pseudonatural transformations $\alpha,\beta:F \Rightarrow E$ as in the diagram
\[
\begin{tikzcd}
F \ar[rr, bend left = 30, ""{name = U}]{}{\alpha} \ar[rr, swap, bend right = 30, ""{name = L}]{}{\beta} & & E \ar[from = U, to = L, shorten <= 4pt, shorten >= 4pt]{}{\eta}
\end{tikzcd}
\]
give rise to equivariant natural transformations (cf.\@ Lemma \ref{Lemma: Modifications give equivariant natural transformations}):
\[
\begin{tikzcd}
F_G(X) \ar[rr, bend left = 30, ""{name = U}]{}{\alpha} \ar[rr, swap, bend right = 30, ""{name = L}]{}{\beta} & & E_G(X) \ar[from = U, to = L, Rightarrow, shorten <= 4pt, shorten >= 4pt]{}{\eta}
\end{tikzcd}
\]
Using this and the fact that 2-functors preserve adjoints, we are able to prove Theorem \ref{Thm: Section 3: Gamma-wise adjoints lift to equivariant adjoints}, which is arguably the most crucial theorem in this paper. It says that when we have pre-equivariant pseudofunctors $F$ and $E$ on $X$ with pseudonatural transformations $\alpha:F \to E$ and $\beta:E \to F$ for which we have adjunctions
\[
\begin{tikzcd}
F(R) \ar[r, bend left = 30, ""{name = U}]{}{\alpha_R} & E(R) \ar[l, bend left = 30, ""{name = L}]{}{\beta_R} \ar[from = U, to = L, symbol = \dashv]
\end{tikzcd}
\]
for every resolution $R$ in the resolution category $R_G(X)$, we get an adjuction between equivariant categories:
\[
\begin{tikzcd}
F_G(X) \ar[r, bend left = 30, ""{name = U}]{}{\alpha} & E_G(X) \ar[l, bend left = 30, ""{name = L}]{}{\beta} \ar[from = U, to = L, symbol = \dashv]
\end{tikzcd}
\]
After this, to close off the section we study the notion of equivariant localization of categories and prove Theorem \ref{Theorem: Section 3: Equivariant localization is a 2-colimit} which shows conditions when we can give an equivariant localization of categories and get a universal property between other equivariant categories.

In Section \ref{Section: Section 3: Change of Space} we continue our study of equivariant functors by studying when we have equivariant functors of the form $F_G(X) \to F_G(Y)$ or $F_G(Y) \to F_G(X)$ that we call Change of Space functors. To study these in detail we define what are called descent pushforwards and descent pullbacks (cf.\@ Definitions \ref{Defn: Section 3.2: Descent Pushforwards} and \ref{Defn: Section 3.2: Descent Pullbacks}, respectively) in order to state when we have pseudonatural transformation-like conditions between a pre-equivariant pseudofunctor $F$ which is defined over both $X$ and $Y$ at the same time. These definitions give conditions in which we can derive Theorems \ref{Thm: Section 3: Existence of pushforwards} and \ref{Thm: Section 3.2: Existence of descent pullbacks} which give conditions in which for any $G$-equivariant morphism $h:X \to Y$ we can induce functors $F_G(X) \to F_G(Y)$ or $F_G(Y) \to F_G(X)$; Theorem \ref{Thm: Section 3: Existence of pullback functors}, which comes earlier in the section, shows that the pseudofunctor $F$ itself always induces a functor $F_G(X) \to F_G(Y)$. After this we prove Theorem \ref{Theorem: Section 3.2: Descent adjunctions lift to equivariant ajdunction} which gives conditions for which we have adjoints between pushforward and pullback functors induced by Definitions \ref{Defn: Section 3.2: Descent Pushforwards} and \ref{Defn: Section 3.2: Descent Pullbacks} (which give tools to show formally why it is that there is an equivariant pullback/pushforward adjunction on $\DbeqQl{X}$ and $\DbeqQl{Y}$ when there is an equivariant morphism $h:X \to Y$).

Section \ref{Section: Section 3: Change of Groups}, the final section in which we study equivariant functors, is the section that studies Change of Groups functors: Functors of the form $F_H(X) \to F_G(X)$ when there is a morphism of smooth algebraic groups $\varphi:G \to H$. The first main theorem in this section is Theorem \ref{Thm: Section 3: Change of groups functor} which proves that for any morphism of smooth algebraic groups $\varphi:G \to H$, there is always an equivariant Change of Groups functor $F^{\sharp}:F_H(X) \to F_G(X)$ by following an idea of Lusztig in \cite{LusztigCuspidal2}. The remainder of this section proves some technical results before giving two sanity check propositions. The first says that for any two morphisms of smooth algebriac groups $\varphi:G_0 \to G_1$ and $\psi:G_1 \to G_2$, there are natural isomorphisms of functors
\[
(\psi \circ \varphi)^{\sharp} \cong \varphi^{\sharp} \circ \psi^{\sharp}.
\]
The second sanity check says that for the unit morphism $1_G:\Spec K \to G$ of any smooth algebraic group, the Change of Groups $1_G^{\sharp}:F_{G}(X) \to F_{\Spec K}(X)$ post-composed with the equivalence $F_{\Spec K}(X) \simeq F(X)$ of Section \ref{Section: Definitions and Stuff} is naturally isomorphic to the forgetful functor $\Forget:F_G(X) \to F(X)$.

In the final section of the paper, Section \ref{Section: Triangulated Cats}, we prove some basic results on triangulations and $t$-structures of equivariant categories. The first of these results is Theorem \ref{Theorem: Section Triangle: Equivariant triangulation}, which shows that when each fibre category and fibre functor of a pre-equivariant pseudofunctor $F$ is triangulated, so is the category $F_G(X)$. Unfortunately because of the $2$-categorical nature of triangulated categories, the proof of Theorem \ref{Theorem: Section Triangle: Equivariant triangulation} requires some use of the Axiom of Choice to make sure that certain isomorphisms satisfy some strict coherence conditions that need not be automatic from the weak limits and colimits that induce them. Afterwards, we move on to discuss when equivariant categories admit $t$-structures induced by $t$-structures on the fibre categories of the pseudofunctor; cf.\@ Theorem \ref{Theorem: Section Triangle: t-structure on our friend the ECat}. After this we give a ``change of heart'' result in Theorem \ref{Theorem: Heart of equivariant t-structure is equivariant heart of t-structure} by showing that if $F_G(X)$ is a triangulated category with $t$-structure induced via equivariant descent, and if $F^{\heartsuit}$ is the pre-equivariant pseudofunctor induced by just descending through the hearts, then there is a strict equality of categories
\[
(F^{\heartsuit})_G(X) = F_G(X)^{\heartsuit}.
\]
Finally we conclude by showing there are $t$-structures on $F_G(X)$ which are induced by $t$-structures on $F(X)$; cf.\@ Theorem \ref{Thm: Section Triangle: t-structures on Equivariant cat}.

\section{Conventions and Notations}
Let us begin by setting some conventions with regards to the formalism we use in this paper.
\begin{itemize}
	\item If $\Cscr$ is a category, we write $\Cscr_0$ for the (class of) objects in $\Cscr$ and $\Cscr_1$ for the (class of) morphisms in $\Cscr$. Similarly, for $n$-categories $\Cscr$, $\Cscr_k$ is the (class of) $k$-cells in $\Cscr$.
	\item Generically script fonts $(\Cscr, \Dscr, \Escr, \cdots)$ are used to denote $1$-categories and German fonts $(\Cfrak, \Dfrak, \Efrak, \cdots)$ are used to denote $2$-categories. In particular, a named category like $\Cat$ is a $1$-category while $\fCat$ is a $2$-category.
	\item The field $K$ is allowed to be arbitrary.
	\item In this thesis a variety over a field $K$ is a quasi-projective reduced separated scheme of finite type over $\Spec K$.\index{Variety}
	\item Throughout this part we assume that for any variety $X$, the prime $\ell$ is coprime to the characteristic of $X$.
	\item All algebraic groups $G$ over $\Spec K$ are assumed to be smooth. While we will make sure to point this out as much as possible, it is worth making this convention now.
	\item Whenever we have a fixed Grothendieck universe $\Uscr$, $\fCAT$\index[notation]{KATZE@$\fCAT$} and $\CAT$\index[notation]{CAT@$\CAT$} denote the (meta)$2$-category and (meta)$1$-category of large categories, while $\fCat$\index[notation]{Katze@$\fCat$} and $\Cat$\index[notation]{Cat@$\Cat$} denote the $2$-category and $1$-category of $\Uscr$-small categories.
	\item Throughout the paper, if not otherwise stated $G$ is a smooth algebraic group over $\Spec K$ with multiplication $\mu:G \times G \to G$ and $X$ is a left $G$-variety over $\Spec K$ with action morphism $\alpha_X:G \times X \to X$. Explicitly, this means that the diagrams
	\[
	\xymatrix{
	G \times G \times X \ar[rr]^-{\mu \times \id_X} \ar[d]_{\id_G\times \alpha_X} & & G \times X \ar[d]^{\alpha_X}  \\
	G \times X \ar[rr]_-{\alpha_X} & & X
	}
	\]
	and
	\[
	\xymatrix{
	X \ar@{=}[d] \ar[rr]^-{\cong} & & \top \times X \ar[d]^{1_G \times \id_X} \\
	X & & G \times X \ar[ll]^-{\alpha_X}
	}
	\]
	commute where $\top$ is the Terminal object of the category $\Var_{/\Spec K}$ of $K$-varieties. Similarly, a morphism $f:X \to Y$ of $G$-varieties is $G$-equivariant if the diagram
	\[
	\xymatrix{
	G \times X \ar[r]^-{\alpha_X} \ar[d]_{\id_G \times f} & X \ar[d]^{f} \\
	G \times Y \ar[r]_-{\alpha_Y} & Y
	}
	\]
	commutes.
	\item We write $\Sch$\index[notation]{Sch@$\Sch, \Sch_{/S}$} (or $\Sch_{/S}$ for the relative case) to denote the category of schemes (or $S$-schemes), $\RMod$\index[notation]{RM@$\RMod$} to denote the category of left modules for a unital ring $R$, $\Var_{/\Spec K}$ (or sometimes $\Var_{/K}$)\index[notation]{Var@$\Var_{/\Spec K}, \Var_{/K}$} to denote the category of varieties over $\Spec K$, and $\Set$\index[notation]{Set@$\Set$} to denote the category of sets.
	\item If $G$ is a smooth algebraic group over $K$, we write $\GVar$\index[notation]{GVar@$\GVar$} and $\GSch$\index[notation]{GSch@$\GSch$} to denote the categories of left $G$-varieties and left $G$-schemes, respectively. These are the categories whose objects are left $G$-varieties (respectively left $G$-schemes) with morphisms given by $G$-equivariant morphisms of varieties (respectively schemes).
	\item We write $\Ab$\index[notation]{Ab@$\Ab$} for the category of Abelian groups. If $\Cscr$ is a category with finite products, we write $\Ab(\Cscr)$\index[notation]{AbC@$\Ab(\Cscr)$} for the category of internal Abelian groups in $\Cscr$.
	\item If $X$ is a locally ringed space, we will write $X = (\lvert X \rvert, \CalO_X)$ where $\lvert X \rvert$ is the underlying space of $X$ and $\CalO_X$ is the structure sheaf of $X$.
	\item Given a left $G$-scheme $X$, if a quotient of $X$ by the action of $G$ (a coequalizer of the action map $\alpha_X:G \times X \to X$ and the projection map $\pi_2:G \times X \to X$) we will write $G \backslash X$ for this quotient in reference to the fact that we want to think of the underlying space of $G \backslash X$ as the space of left $G$-orbits $\lbrace Gx \; | \; x \in \lvert X \rvert \rbrace$ whenever possible.
\end{itemize}

\chapter{A First Look: Equivariant Categories}\label{Section: Equivariant Categories}

\section{Preliminary Notions and Foundational Background}\label{Section: Preliminary and Background for EQCats}
In both \cite{BernLun} and \cite{LusztigCuspidal2}, the equivariant derived category is defined by descent data taken over some category of schemes which resolve the $G$-action on $X$ in order to get around the fact that $X$ may or may not have singularities or other various poorly behaved aspects. For this we follow the scheme-theoretic approach of Lusztig and resolve our $G$-variety $X$ by smooth free principal $G$-varieties. This motivates the definition of the category $\Sf(G)$ below.

\begin{definition}\label{Defn: Section 1: Sf(G)}\index{S@$\Sf(G)$ The Category of Smooth Free $G$-Varieties}\index[notation]{SfG@$\Sf(G)$}
The category $\Sf(G)$ of smooth free $G$-varieties is defined as follows:
\begin{itemize}
	\item Objects: Smooth left $G$-varieties $\Gamma$ of pure dimension with principal $G$-fibrations $\Gamma \to G\backslash \Gamma$, where $G \backslash \Gamma$ is a geometric quotient (cf.\@ \cite[Defnition 0.6]{GIT}) of $\Gamma$ by $G$. That is, the action $\alpha_X$ is {\'e}tale locally trivializable with fibre $G$. More explicitly, there is an {\'e}tale cover $\lbrace U_i \; |\; i \in I \rbrace$ of $G \backslash \Gamma$ for which in each pullback
	\[
	\xymatrix{
	\quo_{\Gamma}^{-1}(U_i) \ar[r] \ar[d] \pullbackcorner &  \Gamma \ar[d]^{\quo} \\
	U_i \ar[r] & G \backslash \Gamma
	}
	\]
	we have that $\quo_{\Gamma}^{-1}(U_i) \cong G \times U_i$ for all $i \in I$.
	\item Morphisms: Smooth morphisms of schemes $f:\Gamma \to \Gamma^{\prime}$ whose fibres have constant dimension, i.e., for each $\gamma^{\prime},\gamma^{\prime\prime} \in \lvert \Gamma^{\prime}\rvert$ we have $\dim(f^{-1}(\gamma^{\prime})) = \dim(f^{-1}(\gamma^{\prime\prime}))$.
	\item Composition and Identities: As in $\Sch$.
\end{itemize}
\end{definition}

We now record some important first properties of smooth morphisms and quotients of smooth morphisms that we will need later.
\begin{proposition}\label{Prop: Smooth quotient from principal G var}
	If $\Gamma$ is a $G$-variety in $\Sf(G)_0$ then the quotient morphism $\quo_{\Gamma}:\Gamma \to G\backslash \Gamma$ is smooth and the diagram
	\[
	\xymatrix{
		G \times \Gamma \ar[r]^{\pi_2} \ar[d]_{\alpha_{\Gamma}} \pullbackcorner & \Gamma \ar[d]^{\quo_{\Gamma}} \\
		\Gamma \ar[r]_-{\quo_{\Gamma}} & G \backslash \Gamma
	}
	\]
	is a pullback square.
\end{proposition}
\begin{proof}
	This is \cite[Proposition 0.9]{GIT}.
\end{proof}
\begin{lemma}[{\cite[\href{https://stacks.math.columbia.edu/tag/02K5}{Tag 02K5}]{stacks-project}}]\label{Lemma: Section 2.1: Smooth morphism lemma}
	Consider the commutative diagram
	\[
	\xymatrix{
		X \ar[rr]^-{f} \ar[dr]_{h} & & Y \ar[dl]^{g} \\
		& Z
	}
	\]
	of schemes for which $f$ is surjective and smooth and $h$ is smooth. Then $g$ is smooth as well.
\end{lemma}
\begin{proposition}\label{Prop: Section 2.1: Smooth quotient map}
	Let $Y, Z$ be $G$-schemes with smooth quotient morphisms $\quo_Y:Y \to G \backslash Y$ and $\quo_Z:Z \to G \backslash Z$. Then if there is a smooth $G$-equivariant morphism $f:Y \to Z$ there exists a unique morphism $\of:G \backslash Y \to G \backslash Z$ making the diagram
	\[
	\xymatrix{
		Y \ar[r]^-{f} \ar[d]_{\quo_Y} & Z \ar[d]^{\quo_Z} \\
		G\backslash Y \ar[r]_-{\of} & G\backslash Z
	}
	\]
	commute. Moreover, $\of$ is smooth as well.
\end{proposition}
\begin{proof}
	The existence of the morphism $\of$ follows from the universal property of the quotient (as a coequalizer of $\alpha_Y, \pi_2:G \times Y \to Y$) applied to the morphism $(\quo_Z \circ f):Y \to G \backslash Z$. That $\of$ is smooth follows from Lemma \ref{Lemma: Section 2.1: Smooth morphism lemma} together with the observations that $\quo_{Y}$ is smooth and surjective while the composite $\quo_Z \circ f$ is smooth.
\end{proof}

The next proposition is a crucial technical result that gives as a corollary that for any $\Gamma \in \Sf(G)_0$ and any left $G$-variety $X$, the product $\Gamma \times X$ is free and admits a smooth quotient. This allows us to resolve the action of $G$ on $X$ by looking at $G$-resolutions of the form $\pi_2:\Gamma \times X \to X$. This is the fundamental technical tool that allows us to describe equivariant information by taking descent data indexed by a pseudofunctor defined on these resolutions.
\begin{proposition}\label{Prop: Section 2.1: The quotient prop}
Let $\Gamma$ be a smooth free $G$-variety and let $Y$ be a $G$-variety. Then the variety $\Gamma \times Y$ is free with smooth quotient morphism 
\[
\quo_{\Gamma \times Y}:\Gamma \times Y \to G \backslash (\Gamma \times Y).
\]
Furthermore, if there is a $G$-equivariant morphism of smooth free $G$-varieties $f:\Gamma \to \Gamma^{\prime}$ then there is a smooth morphism making the diagram
\[
\xymatrix{
	\Gamma \times Y \ar[r]^-{f \times \id_Y} \ar[d]_{q_\Gamma} & \Gamma^{\prime} \times Y \ar[d]^{q_{\Gamma^{\prime}}} \\
	G\backslash(\Gamma \times Y) \ar@{-->}[r]_-{\exists!\overline{f}} & G \backslash(\Gamma^{\prime} \times Y)
}
\]
commute in $\Var_{/\Spec K}$.
\end{proposition}
\begin{proof}
By \cite[Proposition 4]{SerreAlgebraicFibreSpaces} it follows that if we can show for all finite subsets $F \subseteq \lvert Y \rvert$ there is an open affine $U$ of $Y$ for which $F \subseteq \lvert U \rvert$, then it follows that $\Gamma \times Y$ is free. However, by \cite[Proposition 3.3.36]{LiuAGAC} because $Y$ is quasi-projective over $\Spec K$, for any finite set $F \subseteq \lvert Y \rvert$ we can find an open affine subscheme $U$ of $Y$ for which $F \subseteq \lvert U \rvert$ so it follows that $\Gamma \times Y$ is a free $G$-variety. From this it follows that the morphism
\[
\quo_{\Gamma \times Y}:\Gamma \times Y \to G \backslash (\Gamma \times Y)
\]
is smooth, which establishes the first claim. For the second claim we simply apply Proposition \ref{Prop: Section 2.1: Smooth quotient map} to the $G$-equivariant morphism $f \times \id_Y$, which is smooth because $\id_Y$ is {\'e}tale, $f$ is smooth by assumption, and because the product of smooth morphisms is smooth.
\end{proof}
Let us now see how the above propositions allow us to recover what Lusztig did in \cite{LusztigCuspidal2} (which also uses the case where $K = \C$). Our strategy from here is to use Proposition \ref{Prop: Section 2.1: The quotient prop} as a jumping off point for our category of smooth resolutions of $X$. Because for each $G$-variety $X$ and each $\Gamma \in \Sf(G)_0$ the quotient $G \backslash (\Gamma \times X)$ exists and these quotients are functorial in $\Gamma$, we can essentially quotient each such variety by $G$ and then use the category $\Sf(G)$ to resolve $X$. In particular, as we fix $X$ and vary $\Gamma$, by quotienting each product $\Gamma \times X$ by the $G$-action and then taking the derived category of $\C$-sheaves on each quotient variety, we get a family of categories $D^b(G \backslash (\Gamma \times X))$ which are pseudofunctorial in $\Sf(G)$ by the functoriality of the quotients. Thus by taking descent data through this pseudofunctorial assembly of categories $D^b(G \backslash (\Gamma \times X))$ we arrive at Lusztig's definition of the equivariant derived category. To make this work in a more clear and systematic fashion, however, we will introduce the category of smooth free $G$-resolutions of the $G$-variety $X$.

\begin{definition}\label{Defn: Section 1: SfReslX}\index{SfResl@$\SfResl_G(X)$ The Category of Smooth Free $G$-Resolutions of $X$}\index[notation]{SfResl@$\SfResl_G(X)$}
For a $G$-variety $X$ the category $\SfResl_G(X)$ of smooth $G$-resolutions of $X$ is defined as follows:
\begin{itemize}
	\item Objects: Products $\Gamma \times X$ for $\Gamma \in \Sf(G)_0$;
	\item Morphisms: Morphisms of the form $f \times \id_X:\Gamma \times X \to \Gamma^{\prime} \times X$ for $f \in \Sf(G)_1$;
	\item Composition and Identities: As in $\Sch$.
\end{itemize}
\end{definition}
\begin{definition}\index[notation]{XQuot@$\quot{X}{\Gamma}$}
Let $\Gamma \in \Sf(G)_0$ and let $X \in \GVar_0$. Then we write the quotient variety of $\Gamma \times X$ as
\[
\XGamma := G\backslash(\Gamma \times X).
\]
Similarly, if $f \in \Sf(G)_1$ with $f:\Gamma \to \Gamma^{\prime}$, we write $\of:\XGamma \to \XGammap$\index[notation]{f@$\overline{f}$} for the induced map on quotient varieties as in Proposition \ref{Prop: Section 2.1: Smooth quotient map}.
\end{definition}
\begin{remark}
	The various lemmas and constructions above show that there is a functor\index[notation]{Quotient functor @$\quo$}\index{Quotient Functor}
	\[
	\quo:\SfResl_G(X) \to \Var_{/\Spec K}
	\]
	given on objects by $\quo(\Gamma \times X) = \XGamma$ and on morphisms by $\quo(f \times \id_X) = \of$. This is a functor exactly because $\quo$ is defined via universal properties; consequently we have that $\quo(g \circ f) = \overline{g \circ f} = \overline{g} \circ \overline{f} = \quo(g) \circ \quo(f)$. We will use this functor frequently, so it is best to have it defined explicitly here.
\end{remark}

Because we will use the language of pseduofunctors shortly, we will now take a brief digression to define them and our conventions around them. While we will not generally need the full weight of the definition (generically we will care only about the case where a $1$-category $\Cscr$ is regarded as a $2$-category and maps to $\fCat$), we present the full definition in the interest of being thorough.
\begin{definition}\index{Pseudofunctor}
	Let $\Cfrak$ and $\Dfrak$ be $2$-categories. A \textit{pseudofunctor} $F:\Cfrak \to \Dfrak$ is given by:
	\begin{itemize}
		\item An assignment $F_0:\Cfrak_0 \to \Dfrak_0$;
		\item An assignment $F_1:\Cfrak_1 \to \Dfrak_1$ which preserves domains and codomains;
		\item An assignment $F_2:\Cfrak_2 \to \Dfrak_2$ which preserves all sources and targets;
		\item For all composable pairs of $1$-morphisms $X \xrightarrow{f} Y \xrightarrow{g} Z$ in $\Cfrak$, there exists a natural isomorphism of functors $\phi_{f,g}:F(g) \circ F(f) \to F(g \circ f)$ such that the following associativity cocycle condition is satisfied: Given any composable triple $X \xrightarrow{\rho} Y \xrightarrow{\varphi} Z \xrightarrow{\psi} W$ in $\Cfrak$, the pasting diagram
		\[
		\begin{tikzcd}
		& FY \ar[Rightarrow, d, shorten <= 4pt, shorten >=4pt]{}{\phi_{\rho,\varphi}} \ar[dr]{}{F\varphi} & \\
		FX \ar[ur]{}{F\rho} \ar[dr,swap]{}{F(\psi \circ \varphi \circ \rho)} \ar[rr, ""{name = M}]{}[description]{F(\varphi \circ \rho)} & {} & FZ \ar[dl]{}{F\psi} \\
		& FW & \ar[from = M, to = 3-2, Rightarrow, shorten <= 4pt, shorten >= 4pt]{}{\phi_{\varphi \circ \rho, \psi}}
		\end{tikzcd}
		\]
		is equal to the pasting diagram:
		\[
		\begin{tikzcd}
		& FY \ar[dd, ""{name = M}]{}[description]{F(\psi \circ \varphi)} \ar[dr]{}{F\varphi} & \\
		FX \ar[ur]{}{F\rho} \ar[dr, swap]{}{F(\psi \circ \varphi \circ \rho)} & & FZ \ar[from = 2-3, to = M, Rightarrow, shorten <= 5pt, shorten >= 10pt]{}{\phi_{\varphi,\psi}} \ar[dl]{}{F\psi} \\
		& FW & \ar[from = M, to = 2-1, Rightarrow, shorten <= 15pt, shorten >= 2pt]{}{\phi_{\rho, \psi \circ \varphi}}
		\end{tikzcd}
		\]
		\item For any $2$-cell
		\[
		\begin{tikzcd}
		X \ar[r, bend left = 30, ""{name = UL}]{}{f} \ar[r,bend right = 30, swap, ""{name = BL}]{}{g} & Y \ar[r, bend left = 30, ""{name = UR}]{}{k} \ar[r, bend right = 30, swap, ""{name = BR}]{}{k} & Z \ar[from = UL, to = BL, Rightarrow, shorten <= 2pt, shorten >= 2pt]{}{\varphi} \ar[from = UR, to = BR, Rightarrow, shorten <= 2pt, shorten >= 2pt]{}{\psi}
		\end{tikzcd}
		\]
		in $\Cfrak$, we have the equality
		\[
		\phi_{g,k} \circ F(\psi \ast \rho) = F(\psi \ast \rho) \circ \phi_{f,h},
		\]
		where $\ast$ is used to denote the horizontal composition of $2$-morphisms.
	\end{itemize}
\end{definition}
\begin{remark}\index{Pseudofunctor!Rigidified}
	The pseudofunctors we use in this paper are what I refer to as {rigidified pseudofunctors}, i.e.,  the natural isomorphism $\phi:\id_{F(-)} \to F(\id_{-})$ is taken to be the identity natural transformation. This may be done freely for any pseudofunctor without losing information or fundamentally changing the pseudofunctor (in the sense that there is always a pseudonatural isomorphism between a pseudofunctor and its rigidified form; cf.\@ \cite[Expos{\'e} VI.9, pg.\@ 180, 181]{sga1}), so we make this technical simplification.
\end{remark}
\begin{remark}
If $\Cscr$ is a $1$-category, it is always possible to regard $\Cscr$ as a $2$-category. We do this by defining the $2$-cells $\Cscr_2$ to be trivial, i.e., the only $2$-morphisms are the identity $2$-morphism on a $1$-morphism. In this way we can also regard any functor $F:\Cscr \to \Dscr$ as a strict pseudofunctor and also obtain pseudofunctors
\[
F:\Cscr \to \Cfrak
\]
where $\Cfrak$ is a $2$-category. In this thesis we will be primarily interested in the case where $\Cscr = \SfResl_G(X)^{\op}$ and $\Cfrak = \fCat$.
\end{remark}
\begin{remark}
	In this paper we will mostly be interested in pseudofunctors defined on opposite categories. When we do this we will abuse notation and still index the natural transformations by composable arrows in $\Cscr$, save that for a composable set of arrows $X \xrightarrow{f} Y \xrightarrow{g} Z$ in $\Cscr$ we get the natural composition morphism
	\[
	\phi_{f,g}:F(f) \circ F(g) \Rightarrow F(g \circ f)
	\]
	despite the fact that $g \circ f$ is not defined in $\Cscr^{\op}$.
\end{remark}
\begin{remark}
	We need a few remarks here about the logical foundations we use in this paper. For instance, we have not proved that $\Sf(G)_0$ is a set, nor that the bounded derived category $D^{b}(\quot{X}{\Gamma})$ is a small category. In particular, we will see that the indexingwe use to define the collections $A, T_A,$ for an object $(A,T_A) \in D_G^{b}(X)_0$ and for a morphism $P \in D^{b}_G(X)_1$ are not necessarily given in a default language of set theory; perhaps more problematic is that we will use a pseudofunctor that does not take values in $\fCat$ but instead in the (meta)-$2$-category $\fCAT$ of large categories. In order to avoid foundational issues involving the existence of $\CAT$ and $\fCAT$, we instead fix a Grothendieck universe $\Vscr$\index[notations]{V@$\Vscr$} indexed by a strongly inaccessible cardinal $\kappa$ such that $\Vscr$ admits an extension of ZFC set theory for which the category $\Sf(G)$ and the categories $D^{b}(Z)$ and $D^{b}(Z;\overline{\Q}_{\ell})$ of bounded derived ({\'e}tale) sheaves and bounded $\ell$-adic sheaves, repsectively, are $\kappa$-small for all varieties $Z$. This should not be seen as any major technical obstruction, however, as we simply use it to have a well-behaved theory of simplicial $\Vscr$-sets and to avoid the issue of whether or not certain hom-objects are ZFC sets.
\end{remark}

Before proceeding to the full level of generality we use in this paper, let us motivate things by recalling Lusztig's equivariant derived category and show how it relates to the pseudofunctorial formalism described in the introduction. Consider the (strict) pseudofunctor $\Dbb^{b}:\mathbf{Var}_{/\Spec \C}^{\op} \to \fCat$ defined by sending a variety $Y$ over $\C$ to the bounded derived category $D^{b}(Y)$ of ({\'e}tale) sheaves of $\C$-vector spaces on $Y$. The pseudofunctor $\Dbb^{b}$ sends a morphism $f \in \Var(Y,Z)$ to $f^{\ast}:D^{b}(Z) \to D^{b}(Y)$ and is strict, i.e., $\Dbb^{b}(g \circ f) = (g \circ f)^{\ast} = f^{\ast} \circ g^{\ast}$ and $\Dbb^{b}(\id_Y) = \id_{D^{b}(Y)}$. Now define the functor $\quo:\SfResl_G(X) \to \Var_{/\Spec \C}$ via $\quo(\Gamma \times X) := \quot{X}{\Gamma} := G\backslash(\Gamma \times X)$ on objects and $\quo(f \times \id_X) := \overline{f}$ on morphisms. Now regard each $1$-category $\SfResl_G(X)^{\op}$ and $\Var_{/\Spec \C}^{\op}$ as a $2$-category with trivial $2$-morphisms and note that $\quo^{\op}$ is a strict $2$-functor. We then have a pseudofunctor $D^{b}:\SfResl_G(X)^{\op} \to \fCat$ defined via the composite:
\[
\xymatrix{
\SfResl_G(X)^{\op} \ar[rr]^-{D^b} \ar[dr]_{\quo^{\op}} & & \fCat \\
 & \mathbf{Var}_{/\Spec \C} \ar[ur]_{\Dbb^b}
}
\]
This pseudofunctor is the one we use to define Lusztig's equivariant derived category. 
%
\begin{definition}[\cite{LusztigCuspidal2}]\label{Defn: Lusztig EDC}\index[notation]{DbG@$D_G^b(X)$} 
Let $G$ be a smooth algebraic group and let $X$ be a $G$-variety. The bounded equivariant derived category of $X$, $D_G^{b}(X)$, is defined by:
\begin{itemize}
	\item Objects: Pairs $(A,T_A)$ where $A$ is a class of objects
	\[
	A := \lbrace \AGamma \rbrace_{\Gamma \in \Sf(G)_0}
	\]
	where each $\AGamma \in \DbQl{\XGamma}_0$, and where $T_A$ is a class of transition isomorphisms
	\[
	T_A := \left\lbrace \tau_f^A:\of^{\ast}(\AGammap) \xrightarrow{\cong} \AGamma \right\rbrace_{f \in \Sf(G)_1, f:\Gamma \to \Gamma^{\prime}}
	\]
	that satisfy the following cocycle condition: Given a composable pair of arrows $\Gamma \xrightarrow{f} \Gamma^{\prime} \xrightarrow{g} \Gamma^{\prime\prime}$  in the category $\Sf(G)$, the diagram
	\[
	\xymatrix{
	\overline{f}^{\ast}(\overline{g}^{\ast}\quot{A}{\Gamma^{\prime\prime}}) \ar@{=}[d] \ar[rr]^-{\overline{f}^{\ast}\tau_{g}^{A}} & & \overline{f}^{\ast}\quot{A}{\Gamma^{\prime}} \ar[d]^{\tau_{f}^{A}} \\
	(\overline{g \circ f})^{\ast}\quot{A}{\Gamma^{\prime\prime}} \ar[rr]_-{\tau_{g \circ f}^{A}} & & \quot{A}{\Gamma}
	}
	\]
	commutes.
	\item Morphisms: A morphism $P:(A,T_A) \to (B,T_B)$ is an indexed collection of morphisms
	\[
	P := \lbrace \rhoGamma:\AGamma \to \BGamma \rbrace_{\Gamma \in \Sf(G)_0}
	\]
	such that for all $f \in \Sf(G)_1$ with $f:\Gamma \to \Gamma^{\prime}$, the diagram
	\[
	\xymatrix{
	\overline{f}^{\ast}\quot{A}{\Gamma^{\prime}} \ar[r]^{\overline{f}^{\ast}(\quot{\rho}{\Gamma^{\prime}})}  \ar[d]_{\tau_f^A} & \overline{f}^{\ast}\quot{B}{\Gamma^{\prime}} \ar[d]^{\tau_f^B} \\
	\quot{A}{\Gamma} \ar[r]_-{\quot{\rho}{\Gamma}} & \quot{B}{\Gamma}
	}
	\]
	commutes.
	\item Composition: Pointwise over $\Gamma$, i.e., if $\Phi:(A,T_A) \to (B,T_B)$ and if $\Psi:(B,T_B) \to (C,T_C)$ we define
	\[
	\Psi \circ \Phi := \lbrace \quot{\psi}{\Gamma} \circ \quot{\varphi}{\Gamma} \; | \; \Gamma \in \Sf(G)_0, \quot{\varphi}{\Gamma} \in \Phi, \quot{\psi}{\Gamma} \in \Psi \rbrace.
	\]
	\item Identities: The identity morphism $\id_{(A,T_A)}:(A,T_A) \to (A,T_A)$ is the collection
	\[
	\id_{(A,T_A)} := \lbrace \id_{\quot{A}{\Gamma}} \; | \; \Gamma \in \Sf(G)_0 \rbrace.
	\]
\end{itemize}
\end{definition}
\begin{remark}
It is perhaps worth remarking that the category $D_G^b(X)$ that appears above was just defined directly in \cite{LusztigCuspidal2} with no mention of the pseudofunctor $D^b$ anywhere. However, because our point of departure and generalization is to focus on the descent theory controlled by the pseudofunctor $D^b$, we have emphasized it here.
\end{remark}
\section{Equivariant Categories: Definitions, Examples, and Basic Functors}\label{Section: Definitions and Stuff}
Before proceeding to discuss the theory of equivariant categories, we record a crucial scheme-theoretic lemma on algebraic groups. In particular, we need to know that for any smooth algebraic group $G$, $G \in \Sf(G)_0$. While this is obvious enough, it is crucial to the constructions we will do later on (especially involving the change of group functors; cf.\@ Theorem \ref{Thm: Section 3: Change of groups functor}) so we make sure to provide it here as a first step.
\begin{lemma}\label{Lemma: Section 2: Smooth alg group G is in Sf(G)}
Let $G$ be a smooth algebraic group over a field $K$. Then $G \in \Sf(G)_0$.
\end{lemma}
\begin{proof}
The fact that the action $G \times G \to G$ given by the multiplication of the group is free in sense that it is {\'e}tale locally trivializable is given by taking the {\'e}tale cover $\lbrace G \xrightarrow{\id_G} G \rbrace$ and from the pullback diagram:
\[
\xymatrix{
G \times G \ar@{=}[r] \pullbackcorner \ar[d] & G \times G \ar[d]^{\mu} \\
G \ar@{=}[r] & G
}
\]
Moreover, by \cite[Lemma VIB.1.5]{SGA3}, we find that $G$ is a pure variety. Since $G$ is a subgroup scheme of itself, the quotient $G\backslash G$ exists as well by \cite[Theorem V.3.2]{SGA3} (cf.\@ also \cite{milneiAG}). Finally, since $G$ is smooth by assumption, we are done.
\end{proof}

Let us begin our discussion of equivariant categories by first defining what our equivariant categories actually are.

In order to proceed, fix a field $K$ and a smooth algebraic group $G$ over $\Spec K$, together with a (left) $G$-variety $X$. Furthermore fix a pseudofunctor $F:\SfResl_G(X)^{\op} \to \Cat$ factoring as
\[
\xymatrix{
\SfResl_G(X)^{\op} \ar[rr]^{F} \ar[dr]_{\quo^{\op}} & & \fCat \\
 & \Var_{/\Spec K}^{\op} \ar[ur]_{\tilde{F}}
}
\]
for a pseudofunctor that we also call $F:\Var_{/\Spec K}^{\op} \to \fCat$. In particular, we will use this factorization and make an abuse of notation to be able to write
\[
F(f \times \id_X) = \tilde{F}(\quo(f \times \id_X)) = \tilde{F}(\overline{f})
\]
for morphisms and
\[
F(\Gamma \times X) = \tilde{F}(\quo(\Gamma \times X)) = \tilde{F}(\quot{X}{\Gamma})
\]
for objects.
\begin{definition}\label{Defn: Prequivariant pseudofunctors}\index{Pre-equivariant Pseudofunctor}
Let $X$ be a left $G$-variety. We call a pseudofunctor $F:\SfResl_G(X)^{\op} \to \fCat$ factoring as
	\[
	\xymatrix{
		\SfResl_G(X)^{\op} \ar[r]^-{\quo^{\op}} \ar[dr]_{F} & \Var_{/\Spec K}^{\op} \ar[d]^{F} \\
		& \fCat
	}
	\]
	a {pre-equivariant pseudofunctor on (or over) $X$.}
\end{definition}
We are now in a good place to define the equivariant category of a pre-equivariant pseudofunctor $F$ on $X$.
\begin{definition}\label{Defn: Equivariant Cats}\index{Equivariant Category!Definition}
Let $X$ be a left $G$-variety and let $F$ be a pre-equivariant pseudofunctor on $X$. The equivariant category of $F$ on $X$ is the category $F_G(X)$\index[notations]{FGX@$F_G(X)$} defined by:
\begin{itemize}
	\item Objects: Pairs $(A,T_A)$ where $A$ is a collection of objects
	\[
	A = \lbrace \AGamma \rbrace_{\Gamma \in \Sf(G)_0}
	\]
	where each $\AGamma \in F(\XGamma)_0$, and where $T_A$ is a collection of transition isomorphisms
	\begin{align*}
	&\left\lbrace \tau_f^A:F(\of)(\AGammap) \xrightarrow{\cong} \AGamma \right\rbrace_{f \in \Sf(G)_1, f:\Gamma \to \Gamma^{\prime}}
	\end{align*}
	that satisfy the following descent/cocycle condition: Given a composable pair $\Gamma \xrightarrow{f} \Gamma^{\prime} \xrightarrow{g} \Gamma^{\prime\prime}$ in $\Sf(G)$, the diagram
	\[
	\xymatrix{
	F(\overline{f})(F(\overline{g})\quot{A}{\Gamma^{\prime\prime}}) \ar[r]^-{F(\overline{f})\tau_{g}^{A}} \ar[d]_{\phi_{f,g}^{\quot{A}{\Gamma^{\prime\prime}}}} & F(\overline{f})\quot{A}{\Gamma^{\prime}} \ar[d]^{\tau_f^{A}} \\
	F(\overline{g} \circ \overline{f})\quot{A}{\Gamma^{\prime\prime}} \ar[r]_-{\tau_{g \circ f}^{A}} & \quot{A}{\Gamma}
	}
	\]
	commutes in $F(\quot{X}{\Gamma})$.
	\item Morphisms: A morphism $P:(A,T_A) \to (B,T_B)$ is a collection
	\[
	P := \lbrace \rhoGamma:\AGamma \to \BGamma\rbrace_{\Gamma \in \Sf(G)_0}
	\]
	where each $\rhoGamma \in F(\XGamma)_1$ such that for all $f \in \Sf(G)_1$ with $f:\Gamma \to \Gamma^{\prime}$, the diagram
	\[
	\xymatrix{
	F(\overline{f})\quot{A}{\Gamma^{\prime}} \ar[d]_{\tau_f^A} \ar[rr]^-{F(\overline{f})(\quot{\rho}{\Gamma^{\prime}})} & & F(\overline{f})\quot{B}{\Gamma^{\prime}} \ar[d]^{\tau_f^B} \\
	\quot{A}{\Gamma} \ar[rr]_-{\quot{\rho}{\Gamma}} & & \quot{B}{\Gamma}
	}
	\]
	commutes.
	\item Composition: Given pointwise, i.e., if $\Phi:(A,T_A) \to (B,T_B)$ and $\Psi:(B,T_B) \to (C,T_C)$ are morphisms with $\Phi = \lbrace \quot{\phi}{\Gamma} \; | \; \Gamma \in \Sf(G)_0 \rbrace$ and if $\Psi = \lbrace \quot{\psi}{\Gamma} \; | \; \Gamma \in \Sf(G)_0 \rbrace$, then
	\[
	\Psi \circ \Phi := \lbrace \quot{\psi}{\Gamma} \circ \quot{\phi}{\Gamma} \; | \; \Gamma \in \Sf(G)_0 \rbrace.
	\]
	\item Identities: The identity on $(A,T_A)$ is
	\[
	\id_{(A,T_A)} := \lbrace \id_{\quot{A}{\Gamma}} \; | \; \Gamma \in \Sf(G)_0 \rbrace.
	\]
\end{itemize}
\end{definition}
\begin{remark}\label{Remark: Grothendieck Construction Trick}
The equivariant category $F_G(X)$ has an incarnation as an object defined instead of in terms of a pseudofunctor $F$, but in terms of a fibration instead. Indeed, if we have a pseudofunctor $F:\SfResl_G(X)^{\op} \to \fCat$ which gives rise to an equivariant category, let $p:\Fscr \to \SfResl_G(X)$ be the corresponding cloven fibration as induced by the Grothendieck construction. Now let $p_{\cart}:\Fscr_{\cart} \to \SfResl_G(X)$ be the corresponding Cartesian fibration and let $\Fscr_G(X)$ denote the category of sections of $p_{\cart}$, i.e., the category of functors
\[
\underline{A}:\SfResl_G(X) \to \Fscr_{\cart}
\]
for which the diagram
\[
\xymatrix{
\SfResl_G(X) \ar[r]^-{\underline{A}} \ar@{=}[dr] & \Fscr_{\cart} \ar[d]^{p_{\cart}} \\
 & \SfResl_G(X)
}
\]
commutes in $\Cat$. Then it can be shown that $\Fscr_G(X) \simeq F_G(X)$. However, in our study of $F_G(X)$, all the properties we use will be deduced from various properties of $F$, so we have chosen to use the pseudofunctors $F$ to construct and index equivariant categories.
\end{remark}
\begin{remark}
	If it is unlikely to cause confusion, we will refer to an object pair $(A,T_A)$ of $F_G(X)$ simply by $A$ and leave the transition isomorphisms implicit. This is primariy so that we can write things such as $P:A \to B$ or $P \in F_G(X)(A,B)$ for morphisms, as opposed to the more cumbersome $P:(A,T_A) \to (B,T_B)$ and $P \in F_G(X)((A,T_A),(B,T_B))$.
\end{remark}
\begin{remark}
There is a functor
\[
\tilde{\i}:F_G(X) \to \prod_{\Gamma \in \Sf(G)_0} F(\XGamma)
\]
given on objects by
\[
\tilde{\i}(A,T_A) := (\AGamma)_{\Gamma \in \Sf(G)_0}
\]
and on morphisms by
\[
\tilde{\i}(P) = (\rhoGamma)_{\Gamma \in \Sf(G)_0}.
\]
This functor is not generally full, however. To see an example of this let $F = D^b_c(-;\text{{\'e}t})$ be the pre-equivariant pseudofunctor which sends a variety $\Gamma \times X$ to the bounded derived category $D^b_c(\XGamma;\text{{\'e}t})$ of {\'e}tale sheaves on $\XGamma$ and which sends morphisms $f \times \id_X$ to functors $D_c^b(\of;\text{{\'e}t}) = \of^{\ast}$. Then the functor
\[
\tilde{\i}: D_G^b(X;\text{{\'e}t}) \to \prod_{\Gamma \in \Sf(G)_0} D_c^b(\XGamma;\text{{\'e}t})
\]
is not full. To see this consider the object $\mathbbm{1}^{\Z}$ in $D_G^b(X;\text{{\'e}t})$ (the equivariant constant sheaf with value $\Z$); it has object collection
\[
\mathbbm{1}^{\Z} = \left\lbrace (\mathbbm{1}_{\XGamma})^{\Z} \; | \; \Gamma \in \Sf(G)_0 \right\rbrace
\]
where each complex is $(\mathbbm{1}_{\XGamma})^{\Z}$ is the constant sheaf on $\XGamma$ with value $\Z$ embedded in degree $0$. The transition isomorphisms $T_{\mathbbm{1}}$ are constructed by defining each $\tau_f^A$ to be the unique isomorphism induced by the colimit and limit preservation isomorphisms
\[
\of^{\ast}\big((\mathbbm{1}_{\XGammap})^{\Z}\big) = \of^{\ast}\left(\coprod_{x \in \Z}\top_{\Gamma^{\prime}}\right) \cong \coprod_{x \in \Z} \of^{\ast}(\top_{\Gamma^{\prime}}) \cong \coprod_{x \in \Z} \top_{\Gamma} = (\mathbbm{1}_{\XGamma})^{\Z}
\]
where $\top_{\Gamma}, \top_{\Gamma^{\prime}}$ are the terminal objects in the {\'e}tale sheaf toposes $\Shv(\XGamma;\text{{\'e}t})$ and $\Shv(\XGammap;\text{{\'e}t})$, respectively. Note that these isomorphisms exist because $\of^{\ast}$ is a left adjoint and becasue since $\of$ is smooth, it also preserves all finite limits. We then compute that
\[
\tilde{\i}(\mathbbm{1}^{\Z}) = \left((\mathbbm{1}_{\XGamma})^{\Z}\right)_{\Gamma \in \Sf(G)_0}.
\]
Now let $m \in \Z$ with $m \ne 0$ and let $\lambda_m^{\Gamma}:(\mathbbm{1}_{\XGamma})^{\Z} \to (\mathbbm{1}_{\XGamma})^{\Z}$ be the endomorphism induced by left multiplication by $m$. Consider the morphism
\[
\alpha:\tilde{\i}(\mathbbm{1}^{\Z}) \to \tilde{\i}(\mathbbm{1}^{\Z})
\]
defined by
\[
\alpha_{\Gamma} = \begin{cases}
\lambda_m^{\Gamma} & \text{if}\, \Gamma \ne G;\\
0 = \lambda_0^{G} & \text{if}\, \Gamma = G;
\end{cases}
\]
where $0$ is the zero morphism on the constant sheaf. Note that $\alpha$ is a morphism from $\tilde{\i}(\mathbbm{1}^{\Z}) \to \tilde{\i}(\mathbbm{1}^{\Z})$ in $\prod_{\Gamma \in \Sf(G)_0} F(\XGamma)$. 

We claim that $\alpha$ cannot be the image of any $P \in D_c^b(X;\text{{\'e}t})(\mathbbm{1}^{\Z},\mathbbm{1}^{\Z})$; in fact, because of how the isomorphisms $\tau_f$ are constructed, a straightforward argument shows that if $P \in D_c^b(X;\text{{\'e}t})(\mathbbm{1}^{\Z},\mathbbm{1}^{\Z})$ is induced by multiplication by an integer $k$, then for all $\Gamma \in \Sf(G)_0$ we have $\rhoGamma = \lambda_k^\Gamma$. However, because $m \ne 0$ we cannot find a $P \in D_c^b(X;\text{{\'e}t})(\mathbbm{1}^{\Z},\mathbbm{1}^{\Z})$ for which $\rhoGamma = \lambda_m^{\Gamma}$ for $\Gamma \ne G$ and $\quot{\rho}{G} = \lambda_0^G$ when $\Gamma = G$. Thus we conclude that for all $P \in D^b_G(X;\text{{\'e}t})(\mathbbm{1}^{\Z}, \mathbbm{1}^{\Z})$, $\alpha \ne \tilde{\i}(P)$ so $\tilde{\i}$ is not full.
\end{remark}

\begin{example}\label{Example: Section 2: Big List of Examples}
	Here we list some examples of interest to us throughout this paper and for general constructions.
	\begin{itemize}
				\item Let $\Cscr$ be any category. Then define the constant pseudofunctor $F$ on $\SfResl_G(X)^{\op}$ by $F(\Gamma \times X) := \Cscr$ and $F(f \times \id_X) := \id_{\Cscr}$. Then $F_G(X) \simeq \Cscr$, which allows us to realize any category as an equivariant category of fibres over $X$ with respect to the trivial action of $G$ on $X$.
		\item  Let $F$ be the pseudofunctor $F:\SfResl_G(X)^{\op} \to \fCat$ defined by $F(\Gamma \times X) := \Shv(\quot{X}{\Gamma},\text{{\'e}t})$ which sends a scheme $\Gamma \times X$ to the topos of {\'e}tale sheaves of sets over the quotient $\quot{X}{\Gamma}$ and which sends morphisms $f \times \id_X:\Gamma \times X \to \Gamma^{\prime} \times X$ to the pullback
		\[
		\overline{f}^{\ast}:\Shv(\quot{X}{\Gamma^{\prime}},\text{{\'e}t}) \to\Shv(\quot{X}{\Gamma},\text{{\'e}t}).
		\]
		Then the category $F_{G}(X) := \Shv_{G}(X,\text{{\'e}t})$\index[notation]{Sheaves of Sets@$\Shv_G(X,\text{{\'e}t})$} defines the category of $G$-equivariant {\'e}tale sheaves.
		\item Let $F$ be the psuedofunctor defined by
		\[
		F(\Gamma \times X) := D^{b}_c(\quot{X}{\Gamma})
		\]
		and
		\[
		F(f \times \id_X:\Gamma \times X \to \Gamma^{\prime} \times X) := \of^{\ast}:D^{b}_c(\quot{X}{\Gamma^{\prime}}) \to D^{b}_c(\quot{X}{\Gamma}).
		\]
		Then the category $F_{G}(X) := D^{b}_{G,c}(X)$\index[notation]{DGbc@$D_{G,c}^b(X)$} is the equivariant bounded derived category of constructible sheaves on $X$.
		\item Let $F$ be the pseudofunctor defined by
		\[
		F(\Gamma \times X) := \Loc(\quot{X}{\Gamma},\overline{\Q}_{\ell}),
		\]
		where $\Loc(\quot{X}{\Gamma},\overline{\Q}_{\ell})$ is the category of $\ell$-adic ({\'e}tale) local systems of $\quot{X}{\Gamma}$ and
		\[
		F(f \times \id_X:\Gamma \times X \to \Gamma^{\prime} \times X):= \of^{\ast}: \Loc(\quot{X}{\Gamma^{\prime}},\overline{\Q}_{\ell}) \to \Loc(\quot{X}{\Gamma},\overline{\Q}_{\ell}).
		\]
		Then the category $F_{G}(X) := \Loc_{G}(X;\overline{\Q}_{\ell})$\index[notation]{LocLadic@$\Loc_G(X,\overline{\Q}_{\ell})$} is the category of $G$-equivariant $\ell$-adic ({\'e}tale) local systems on $X$.
		\item Let $F$ be the psuedofunctor defined by
		\[
		F(\Gamma \times X) := \QCoh(\quot{X}{\Gamma})
		\]
		and
		\[
		F(f \times \id_X:\Gamma \times X \to \Gamma^{\prime} \times X) := \of^{\ast}:\QCoh(\quot{X}{\Gamma^{\prime}}) \to \QCoh(\quot{X}{\Gamma}).
		\]
		Then the category $F_{G}(X) := \QCoh_{G}(X)$\index[notation]{Quasi-coherent sheaves@$\QCoh_G(X)$} is the category of equivariant quasi-coherent sheaves on $X$.
		\item Let $F$ be the psuedofunctor defined on objects by
		\[
		F(\Gamma \times X) := \Ab(\Shv)(\XGamma,\text{{\'e}t})
		\]
		where $\Ab(\Shv)(\XGamma,{\text{{\'e}t}})$ is the category of Abelian {\'e}tale sheaves on $\XGamma$, and on morphisms by
		\begin{align*}
		&F(f \times \id_X:\Gamma \times X \to \Gamma^{\prime} \times X) \\
		&:= \of^{\ast}:\Ab(\Shv)(\XGammap,\text{{\'e}t}) \to \Ab(\Shv)(\XGamma,\text{{\'e}t}).
		\end{align*}
		Then the category $F_{G}(X) := \Ab(\Shv)_G(X)$\index[notation]{AbSheavesG@$\Ab(\Shv)_G(X)$} is the category of equivariant Abelian sheaves on $X$.
		\item Let $F$ be the pseudofunctor defined by
		\[
		F(\Gamma \times X) := D(\XGamma),
		\]
		the derived category of $\XGamma$ sheaves and let $F$ be defined on morphisms by
		\[
		F(f \times \id_X) = F(\of).
		\]
		Then $F_G(X) = D_G(X)$\index[notation]{DGX@{$D_G(X)$}} is the full equivariant derived category of sheaves on $X$.
		\item Let $F$ be the pseudofunctor defined by
		\[
		F(\Gamma \times X) := \Per(\quot{X}{\Gamma})
		\]
		on objects, where $\Per(\quot{X}{\Gamma})$ is the category of perverse sheaves on $\quot{X}{\Gamma}$, and
		\[
		F(f \times \id_X:\Gamma \times X \to \Gamma^{\prime} \times X) := {}^{p}\of^{\ast}:\Per(\quot{X}{\Gamma^{\prime}}) \to \Per(\quot{X}{\Gamma})
		\]
		on morphisms; note that ${}^{p}\of^{\ast} = \of^{\ast}[d]$, where $d$ is the relative dimension of the morphism $\of$. Then $F_G(X) = \Per_G(X)$\index[notation]{PerGX@$\Per_G(X)$} defines the category of $G$-equivariant perverse sheaves on $X$. Following this same construction mutatis mutandis, save with the fibre categories taking values as 
		\[
		F(\Gamma \times X) = \Per(\quot{X}{\Gamma},\overline{\Q}_{\ell})
		\] 
		for the categories of $\ell$-adic perverse sheaves, we also get a $G$-equivariant category of $\ell$-adic perverse sheaves
		\[
		F_G(X) = \Per_G(X;\overline{\Q}_{\ell}).
		\]\index[notation]{PerGXladic@$\Per_G(X;\overline{\Q}_{\ell})$}
		\item Let $F$ be the pseudofunctor defined by \[
		F(\Gamma \times X) := \DMod(\quot{X}{\Gamma}) 
		\]
		on objects and 
		\[
		F(f \times \id):= \of^{\ast}:\DMod(\quot{X}{\Gamma^{\prime}}) \to \DMod(\XGamma)
		\] 
		on morphisms. Then the category $F_G(X) = \DGMod(X)$\index[notation]{DModulesG@$\Dcal$-$\mathbf{Mod}_{G}(X)$} is the category of equivariant $D$-modules on $X$.
	\end{itemize}
\end{example}
\begin{remark}
When defining and considering pre-equivariant pseudofunctors, while we have given them as pseudofunctors which factor through $\Var_{/\Spec K}^{\op}$, this is not strictly speaking a necessary condition. Because every quotient functor $\quo^{\op}:\SfResl_G(X)^{\op} \to \Sch_{/\Spec K}^{\op}$ factors through $\Var_{/\Spec K}^{\op}$, it suffices to define pre-equivariant pseudofunctors as pseudofunctors $F:\SfResl_{G}(X)^{\op} \to \fCat$ which factor as:
\[
\xymatrix{
\SfResl_G(X)^{\op} \ar[r]^-{\quo^{\op}} \ar[dr]_{F} & \Cscr^{\op} \ar[d]^{F} \\
 & \fCat
}
\]
where $\Cscr$ is any category $\Cscr \hookrightarrow \Sch_{/\Spec K}$ which contains the essential image of $\quo^{\op}$. This will suffice to give the same equivariant category $F_G(X)$ as the ones defined via the pre-equivariant pseudofunctors of Definition \ref{Defn: Prequivariant pseudofunctors}. When considering generalizations of the theory described in this paper to equivariant categories on schemes $X$ which are not varieites, it may be of interest to take categories like the categories of separated schemes over $\Spec K$, quasi-separated quasi-compact schemes over $\Spec K$, schemes over $\Spec K$, and the category of pro-varieties over $\Spec K$, or other such categories depending on the structure and properties of the schemes $X$ and $\XGamma$.
\end{remark}

We now give a short series of results that show that in the case $G = \Spec K$ is the trivial group over $K$, $F_{\Spec K}(X) \simeq F(X)$.
\begin{lemma}
If $G = \Spec K$, then $\Sf(G)$ is the category of smooth $K$-varieties of finite pure dimension with smooth morphisms with constant fibre dimension.
\end{lemma}
\begin{proof}
This is simply an unwinding of the definition of $\Sf(G)$ specialized to $G = \Spec K$. That any smooth $K$-variety $\Gamma$ of pure dimension is $\Spec K$-free follows from the fact that if $\lbrace U_i \; | \; i \in I \rbrace$ is any {\'e}tale cover of $\Gamma$, and if the action map, which is a canonical isomorphism $\alpha: \Spec K \times\Gamma \xrightarrow{\cong} \Gamma$, then the pullback diagram
\[
\xymatrix{
\alpha^{-1}(U_i) \ar[r] \ar[d]_{\cong} \pullbackcorner & \Spec K \times \Gamma \ar[d]^{\cong} \\
U_i \ar[r] & \Gamma
}
\]
gives an isomorphism $\alpha^{-1}(U_i) \cong U_i$ and hence a further isomorphism $\alpha^{-1}(U_i) \cong U_i \cong \Spec K \times U_i$.
\end{proof}
\begin{lemma}
The category $\Sf(\Spec K)$ has a terminal object.
\end{lemma}
\begin{proof}
The terminal object is $\Spec K$; note that if we verify that the variety $\Spec K \in \Sf(\Spec K)_0$, the lemma is immediate from the observations that $\Spec K$ is terminal in $\Var_{/\Spec K}$ and from $\Sf(\Spec K)$ being a subcategory of $\Var_{/\Spec K}$ with structure map $\nu:\Gamma \to \Spec K$ smooth of constant fibre dimension. Now the variety $\Spec K$ is smooth over itself (as the identity map $\id_{\Spec K}$ is always {\'e}tale) and $\Spec K$ is a $0$-dimensional irreducible variety over $\Spec K$. Thus $\Spec K \in \Sf(\Spec K)_0$.
\end{proof}
\begin{proposition}\label{Prop: Section 2: Equivalence of K equivariant cat with F X}
	Let $X$ be a $K$-variety and let $F:\SfResl_{\Spec K}(X)^{\op} \to \fCat$ be a pre-equivariant pseudofunctor (cf.\@ Definition \ref{Defn: Prequivariant pseudofunctors}). There is then an equivalence of categories $F(X) \simeq F_{\Spec K}(X)$.
\end{proposition}
\begin{proof}
First note that since the action of $\Spec K$ on $X$ is trivial, we have that there is an isomorphism
\[
\quot{X}{\Spec K} = \Spec K\backslash(\Spec K \times X) \cong \Spec K \backslash X \cong X.
\]
Let $\psi:\quot{X}{\Spec K} \xrightarrow{\cong} X$ be the induced isomorphism above. Recall also that since $F$ is pre-equivariant, it sits in a strictly commuting diagram
\[
\xymatrix{
\SfResl_{\Spec K}(X)^{\op} \ar[dr]_{\quo^{\op}} \ar[rr]^{F_X} & & \fCat \\
 & \Var_{/\Spec K}^{\op} \ar[ur]_{F}
}
\]
of $2$-categories where we abuse notation and write $F_X = F$. However, we have that the fibre category $F(X)$ is defined as well.

We now claim that $F(X) \simeq F(\quot{X}{\Spec K})$. This is immediate from the isomorphism $\psi$. In particular, because $\psi$ is an isomorphism and $\psi^{-1} \circ \psi = \id_{X}$, the pseudofunctoriality of $F$ gives us $2$-cells
\[
\begin{tikzcd}
 & F(\quot{X}{\Spec K}) \ar[dr]{}{F(\psi^{-1})} \ar[d, Rightarrow, shorten <= 4pt, shorten >= 4pt]{}{\cong}  &  \\
F(X) \ar[rr, equals]{}{} \ar[ur]{}{F(\psi)} & {} & F(X) 
\end{tikzcd}
\]
\[
\begin{tikzcd}
 & F(X) \ar[dr]{}{F(\psi)} \ar[d, Rightarrow, shorten <= 4pt, shorten >= 4pt]{}{\cong} \\
 F(\quot{X}{\Spec K}) \ar[rr, equals] \ar[ur]{}{F(\psi^{-1})} & {} & F(\quot{X}{\Spec K})
\end{tikzcd}
\]
where the natural isomorphisms are the compositors $\phi_{\psi,\psi^{-1}}$ and $\phi_{\psi^{-1},\psi}$, respectively. In particular, this shows that $F(\psi)$ is an equivalence.

We now define functors $L:F(X) \to F_{\Spec K}(X)$ and $R:F_{\Spec K}(X) \to F(X)$. For this first let $A \in F(X)_0$ and let $\rho \in F(X)_1$. For each $\Gamma \in \Sf(\Spec K)$, let $\nu_{\Gamma}:\Gamma \to \Spec K$ be the structure map induced by the fact that $\Spec K$ is terminal. We then define $L:F(X) \to F_{\Spec K}(X)$ by sending $A$ to the object $(LA, T_{LA})$, where 
\[
LA := \lbrace F(\overline{\nu}_{\Gamma})\big(F(\psi)A\big) \; | \; \Gamma \in \Sf(\Spec K)_0 \rbrace.
\]
To see how to define the transition isomorphisms $T_{LA}$, first note that for all $f \in \Sf(G)_1$ we have commuting diagrams:
\[
\xymatrix{
\Gamma \ar[rr]^{f} \ar[dr]_{\nu_{\Gamma}}  & & \Gamma^{\prime} \ar[dl]^{\nu_{\Gamma^{\prime}}} \\
 & \Spec K
}
\]
Applying $F$ gives rise to the pasting diagram
\[
\begin{tikzcd}
 & F(\XGammap) \ar[dr]{}{F(\overline{f})} \ar[d, Rightarrow, shorten <= 4pt, shorten <= 4pt]{}{\phi_{f,\nu_{\Gamma^{\prime}}}} \\
F(\quot{X}{\Spec K}) \ar[rr, swap]{}{F(\overline{\nu}_{\Gamma})} \ar[ur]{}{F(\overline{\nu}_{\Gamma^{\prime}})} & {} & F(\quot{X}{\Gamma})
\end{tikzcd}
\]
and hence induces the isomorphism 
\[
\phi_{f,\nu_{\Gamma^{\prime}}}^{F(\psi)A}:F(\overline{f})\big( F(\overline{\nu}_{\Gamma^{\prime}})\big(F(\psi)A\big)\big) \to F(\overline{\nu}_{\Gamma})\big(F(\psi)A\big).
\]
We thus define $\tau_f^{LA} := \phi_{f,\nu_{\Gamma^{\prime}}}^{F(\psi)A}$ and take $T_{LA} := \lbrace \tau_f^{LA} \; | \; f \in \Sf(\Spec K)_1 \rbrace$. That these satisfy the cocycle condition is routinely verified, as they are built up of the compositors of a pseudofunctor and pseudofunctors preserve composition up to compositors, which is exactly the cocycle condition. Finally, we define the functor $L$ on morphisms by setting
\[
L(\rho) := \lbrace F(\overline{\nu}_{\Gamma})\big(F(\psi)\rho\big) \; | \; \Gamma \in \Sf(\Spec K)_0 \rbrace.
\]
That this is an $F_{\Spec K}(X)$ morphism is also routinely verified and follows from the naturality of the compositors.

To define the functor $R:F_{\Spec K}(X) \to  F(X)$, let $(A,T_A)$ be an object with $A = \lbrace \AGamma \; | \; \Gamma \in \Sf(\Spec K)_0\rbrace$ and let $P$ be a morphism with $P = \lbrace \quot{\rho}{\Gamma} \; | \; \Gamma \in \Sf(\Spec K)_0 \rbrace$. We define
\[
RA := F(\psi^{-1})\quot{A}{\Spec K}
\]
and
\[
RP := F(\psi^{-1})\quot{\rho}{\Spec K}.
\]
That these assignments are functorial is trivially verified.

We now show that $L$ and $R$ are inverse equivalences of categories. For this we first verify that for all objects $A$ of $F(X)$,
\begin{align*}
(R \circ L)A &= R(LA) = R\left(\lbrace F({\overline{\nu}}_{\Gamma})\big(F(\psi)A\big) \; | \; \Gamma \in \Sf(\Spec K) \rbrace\right) \\
&= F(\psi^{-1})\big(F(\overline{\nu}_{\Spec K})\big(F(\psi)A\big)\big) = F(\psi^{-1})\big(F(\overline{\id_{\Spec K}})\big(F(\psi)A\big)\big) \\
&= F(\psi^{-1})\big(F(\psi)A\big) \cong A
\end{align*}
and similarly for morphisms we have that
\begin{align*}
(R \circ L)\rho &= R(L\rho) =R\big(\lbrace F\big(\overline{\nu}_{\Gamma}\big(F(\psi)\rho\big)\big) \; | \; \Gamma \in \Sf(\Spec K)_0 \rbrace\big) \\
&= F(\psi^{-1})\big(F(\psi)\rho\big) \cong \rho,
\end{align*}
where the naturality of the isomorphism follows from the pseudofunctoriality of $F$. Thus we get that $R \circ L \cong \id_{F(X)}$.

We now verify the other direction of the equivalence. Namely that $L \circ R$ is naturally isomorphic to the identity functor on $F_{\Spec K}(X)$. For this we consider first that for an object $(A,T_A)$ we have
\begin{align*}
(L \circ R)(A) &= L(R(A)) = L(F(\psi^{-1}\quot{A}{\Spec K})) \\
&= \lbrace F(\overline{\nu}_{\Gamma})\big(F(\psi)\big(F(\psi^{-1})\quot{A}{\Spec K}\big)\big) \; | \; \Gamma \in \Sf(\Spec K)_0 \rbrace
\end{align*}
while $T_{L(RA)}$ is induced as above. We first claim that for all $\Gamma \in \Sf(\Spec K)_0$, we have isomorphisms
\[
\quot{L(RA)}{\Gamma} \xrightarrow{\cong} \AGamma.
\]
To see this note that first there is a natural isomorphism $\theta$ which, when evaluated at $\quot{A}{\Spec K}$ gives an isomorphism
\[
\theta_{A}:\big(F(\psi) \circ F(\psi^{-1})\big)\quot{A}{\Spec K} \xrightarrow{\cong} \quot{A}{\Spec K}
\]
induced by the isomorphism $F(\psi) \circ F(\psi^{-1}) \cong \id_{F(\quot{X}{\Spec K})}$. From here note that the set of transition isomorphisms $T_A$ furnishes isomorphisms
\[
\tau_{\nu_{\Gamma}}^{A}:F(\overline{\nu}_{\Gamma})\quot{A}{\Spec K} \xrightarrow{\cong} \AGamma
\]
so by composing yields an isomorphism 
\[
\big(F(\overline{\nu}_{\Gamma}) \circ F(\psi) \circ F(\psi^{-1})\big)\quot{A}{\Spec K} \xrightarrow{F(\overline{\nu}_{\Gamma})\theta_A} F(\overline{\nu}_{\Gamma})\quot{A}{\Spec K} \xrightarrow{\tau_{\nu_{\Gamma}}^{A}} \AGamma.
\]
However, since $\quot{R(LA)}{\Gamma} = \big(F(\overline{\nu}_{\Gamma}) \circ F(\psi) \circ F(\psi^{-1})\big)\quot{A}{\Spec K}$, by setting 
\[
\quot{\zeta}{\Gamma} = \tau_{\nu_{\Gamma}} \circ F(\overline{\nu}_{\Gamma})
\]
we obtain an isomorphism
\[
\quot{\zeta}{\Gamma}:\quot{L(RA)}{\Gamma} \xrightarrow{\cong} \AGamma.
\]
We claim that the $\Vscr$-set $Z := \lbrace \quot{\zeta}{\Gamma} \; | \; \Gamma \in \Sf(\Spec K)_0 \rbrace$ is an $F_{\Spec K}(X)$-morphism. To prove this let $f \in \Sf(\Spec K)_1$ and write $\Dom f = \Gamma$, $\Codom f = \Gamma^{\prime}$. We will show that the diagram
\[
\xymatrix{
F(\overline{f})\quot{L(RA)}{\Gamma^{\prime}} \ar[rr]^-{F(\overline{f})\quot{\zeta}{\Gamma^{\prime}}} \ar[d]_{\tau_{f}^{L(RA)}} & & F(\overline{f})\quot{A}{\Gamma^{\prime}} \ar[d]^{\tau_f^A} \\
\quot{L(RA)}{\Gamma} \ar[rr]_-{\quot{\zeta}{\Gamma}} & & \AGamma
}
\]
commutes. For this note that
\begin{align*}
\tau_f^A \circ F(\overline{f})\quot{\zeta}{\Gamma^{\prime}} &= \tau_f^A \circ F(\overline{f})\big(\tau_{\nu_{\Gamma^{\prime}}}^{A} \circ F(\overline{\nu}_{\Gamma^{\prime}})\theta_A\big) \\
&= \tau_f^A \circ F(\overline{f})\tau_{\nu_{\Gamma^{\prime}}}^{A} \circ \big(F(\overline{f}) \circ F(\overline{\nu}_{\Gamma^{\prime}})\big)\theta_A  \\
&= \tau_{\nu_{\Gamma^{\prime}} \circ f}^{A} \circ \phi_{f,\nu_{\Gamma^{\prime}}} \circ \big(F(\overline{f}) \circ F(\overline{\nu}_{\Gamma^{\prime}})\big)\theta_A  \\
&= \tau_{\nu_{\Gamma}}^{A} \circ F(\overline{\nu}_{\Gamma})\theta_A \circ \phi_{f,\nu_{\Gamma}} = \tau_{\nu_{\Gamma}}^{A} \circ F(\overline{\nu}_{\Gamma})\theta_A \circ \tau_{f}^{L(RA)} \\
&= \quot{\zeta}{\Gamma} \circ \tau_{f}^{L(RA)}.
\end{align*}
We conclude that $Z$ is an isomorphism. Thus $L \circ R \cong \id_{F_{\Spec K}(X)}$ and so $F(X) \simeq F_{\Spec K}(X)$.
\end{proof}

As an example, we construct the functor $\Forget:F_G(X) \to F(X)$ which forgets the $G$-equivariance and returns the $F(X)$-component of the $G$-equivariant information.

\begin{example}\label{Example: Section 2: The forgetful functor}
Let $F$ be any pre-equivariant pseudofunctor on $X$. There is a functor
\[
\Forget:F_G(X) \to F(X)
\]
\index[notation]{ForgetfulFunctor@$\Forget$}given as follows: Note that there is a (non-equivariant) isomorphism of varieties
\[
\psi:X \xrightarrow{\cong} \quot{X}{G}
\]
which essentially trivializes the action of $G$ on $X$. We define our forgetful functor $\Forget:F_G(X) \to F(X)$ by
\[
\lbrace \AGamma \; | \; \Gamma \in \Sf(G)_0 \rbrace \mapsto F(\psi^{-1})(\quot{A}{G})
\]
on objects and
\[
\lbrace \quot{\rho}{\Gamma} \; | \; \Gamma \in \Sf(G) \rbrace \mapsto F(\psi^{-1})(\quot{\rho}{G})
\]
on morphisms. A braver author than I would call this functor $\text{Foret}:F_G(X) \to F(X)$ because it forgets the $G$.
\end{example}
\begin{proposition}
For any pre-equivariant pseudofunctor $F$, there is a functor $\Eq:F(X) \to F_G(X)$ and a natural isomorphism $\nu:\Forget \circ \Eq \xRightarrow{\cong} \id_{F(X)}$.
\end{proposition}
\begin{proof}
Begin by noting with $\psi:X \xrightarrow{\cong} \quot{X}{G}$ the isomorphism above, if $\Gamma \in \Sf(G)_0$ is arbitrary then there is a unique morphism $\overline{\Gamma}:\XGamma \to \quot{X}{G}$ induced by the universal property of the quotient in the following diagram:
\[
\xymatrix{
\Gamma \times X \ar[d]_{\quo_{\Gamma}} \ar[r]^-{\pi_2^{\Gamma}} & X \ar[d]^{\psi} \\
\XGamma \ar@{-->}[r]_-{\exists!\overline{\Gamma}} & \quot{X}{G}
}
\]
These maps $\overline{\Gamma}$ have the following property for any $f \in \Sf(G)_1$: The diagram
\[
\begin{tikzcd}
 & & & \quot{X}{G}\\
\\
 & X \ar[uurr]{}{\psi} & \XGamma \ar[rr, swap]{}{\of} \ar[uur, swap]{}{\overline{\Gamma}} & & \XGammap \ar[uul]{}{\overline{\Gamma}^{\prime}}\\
\\
\Gamma \times X  \ar[rr, swap]{}{f \times \id_X} \ar[uurr, swap, near start]{}{\quo_{\Gamma}} & & \Gamma^{\prime} \times X \ar[uurr, swap]{}{\quo_{\Gamma^{\prime}}}
\ar[from = 5-1, to = 3-2, crossing over]{}{\pi_2^{\Gamma}} \ar[from = 5-3, to = 3-2, swap, crossing over, near start]{}{\pi_2^{\Gamma^{\prime}}}
\end{tikzcd}
\]
commutes in $\Var_{\Spec K}$. From this we define, for any $A \in F(X)_0$ the object
\[
\Eq A := \lbrace F(\overline{\Gamma})F(\psi)A \; | \; \Gamma \in \Sf(G)_0 \rbrace
\]
with transition isomorphisms induced by the definition diagram
\[
\begin{tikzcd}
F(\of)F(\overline{\Gamma}^{\prime})F(\psi)A \ar[dd, equals] \ar[r]{}{\phi_{\of,\overline{\Gamma}}^{F(\psi)A}} & F(\overline{\Gamma}^{\prime} \circ \of)F(\psi)A \ar[d, equals] \\
 & F(\overline{\Gamma})F(\psi) \ar[d, equals]\\
F(\of)\quot{\Eq A}{\Gamma^{\prime}} \ar[r, swap]{}{\tau_f^{\Eq A}} & \quot{\Eq A}{\Gamma}
\end{tikzcd}
\]
for any $f \in \Sf(G)_1$. That these transition isomorphisms satisfy the cocycle condition is immediate from the pseudofuncoriality of $F$, so we do indeed have an object of$F_G(X)$. The assignment on morphisms is given similarly: For any $\alpha \in F(X)_1,$
\[
\Eq \alpha := \lbrace F(\overline{\Gamma})F(\psi)\alpha \; | \; \Gamma \in \Sf(G)_0 \rbrace
\]
and that these do indeed define morphisms follows from the fact that the compositors $\phi$ are natural isomorphisms. Thus we get our desired functor $\Eq$.

For the claim on the natural isomorphism we simply calculate that
\[
(\Forget \circ \Eq)A = F(\psi^{-1})F(\psi)A
\]
so setting $\nu := \phi_{\psi^{-1},\psi}$ proves the proposition.
\end{proof}

\section{Some Categorical Properties of Equivariant Categories}\label{Section 2: Cat theory of ECat}
To begin our study of equivariant categories, we first introduce some terminology (so that we can talk succinctly about pseudofunctors that universally take values in additive categories, triangulated categories, $\Vscr$-locally small categories, and other such things) before proceeding to develop some basic theory of equivariant categories.

\begin{definition}\label{Defn: Additive, triangulated, and locally small pseudofunctors}
Fix a pre-equivariant pseudofunctor $F:\SfResl_G(X)^{\op} \to \fCat$ (cf.\@ Definition \ref{Defn: Prequivariant pseudofunctors}). If for all $\Gamma \times X \in \SfResl_G(X)_0$ and $f \in \Sf(G)_1$ the category $F(\quot{X}{\Gamma})$ is additive,\index{Pre-equivariant Pseudofunctor! Additive} triangulated\index{Pre-equivariant Pseudofunctor!Triangulated} or $\Vscr$-locally small\index{Pre-equivariant Pseudofunctor! V@{$\Vscr$-Locally Small}}  and each functor $F(\overline{f})$ preserves the relevant structure, we call $F$ an additive, triangulated, or $\Vscr$-locally small (pre-equivariant) pseudofunctor.
\end{definition}
\begin{example}
The pseudofunctor $F:\SfResl_G(X)^{\op} \to \fCat$ given by
\[
F(\XGamma) := \DbQl{\XGamma}
\]
on objects and
\[
F(\of) := \of^{\ast}:\DbQl{\XGammap} \to \DbQl{\XGamma}
\]
on morphisms is an additive, triangulated, and $\Vscr$-locally small pseudofunctor.
\end{example}

We now show a quick smallness result. Explicitly, for our applications to additive and triangulated categories in mind, we need to know when we can define an addition on the hom-objects $F_G(X)(A,B)$ for all $A$ and $B$ without leaving the universe $\Vscr$. This leads us immediately to use $\Vscr$-locally small pre-equivariant pseudofunctors, as they give us sufficient structure to ensure that $F_G(X)(A,B)$ is always a $\Vscr$-set.

\begin{lemma}\label{Lemma: Section 2: Pseudofunctor values in Vlocally small means Equiv cat Vlocally small}
Let $F:\SfResl_G(X)^{\op} \to \fCat$ be a $\Vscr$-locally small pseudofunctor. Then $F_{G}(X)$ is $\Vscr$-locally small.
\end{lemma}
\begin{proof}
Fix two arbitrary objects $(A,T_A), (B,T_B) \in F_G(X)_0$; note that it suffices to show that $F_G(X)(A,B)$ is a $\Vscr$-set, i.e., that $\lvert F_G(X)(A,B)\rvert < \kappa$, where $\kappa$ is the strongly inaccessible cardinal defining the Grothendieck universe $\Vscr$.

Recall that the universe $\Vscr$ is defined with the property that $\Sf(G)_0$ and $\Sf(G)_1$ are both $\Vscr$-small, i.e., that $\lvert \Sf(G)_0 \rvert, \lvert \Sf(G)_1 \rvert < \kappa$. For what follows, consider the topos $\VSet$ of $\Vscr$-sets and let $\Omega \in \VSet_0$ be the corresponding subobject classifier. Because $\kappa$ is strongly inaccessible, the power objects $\Pcal_{\Vscr}(X) := [X,\Omega]$ in $\VSet_0$ satisfy $\lvert \Pcal_{\Vscr}(X)\rvert < \kappa$.

Consider now that since $P$ is defined as the set
\[
P = \lbrace \quot{\rho}{\Gamma} \; |\; \Gamma \in \Sf(G)_0 \rbrace = \bigcup_{\Gamma \in \Sf(G)_0} \lbrace \quot{\rho}{\Gamma} \rbrace
\]
we have that
\[
P \subseteq \bigcup_{\Gamma \in \Sf(G)_0} F(\quot{X}{\Gamma})(\quot{A}{\Gamma},\quot{B}{\Gamma}).
\]
Because each of the categories $F(\quot{X}{\Gamma})$ is locally $\Vscr$-small, each hom-collection $F(\quot{X}{\Gamma})(\quot{A}{\Gamma},\quot{B}{\Gamma})$ is a $\Vscr$-set. Thus, since $\Sf(G)_0$ is a $\Vscr$-set as well, we have that the union
\[
X := \bigcup_{\Gamma \in \Sf(G)_0} F(\quot{X}{\Gamma})(\quot{A}{\Gamma},\quot{B}{\Gamma})
\]
is a $\Vscr$-set as well.

Now note that for each $P \in F_G(X)(A,B)$, we have that $\lbrace P \rbrace \subseteq \Pcal_{\Vscr}(X)$. Thus it follows that
\[
\Pcal_{\Vscr}(X) \supseteq \bigcup_{P \in F_G(X)(A,B)}\lbrace P \rbrace = F_G(X)(A,B)
\]
and hence we have that $F_G(X)(A,B)$ is a $\Vscr$-set. This proves the lemma.
\end{proof}

Our main goal in the development of equivariant categories at the moment is to show that if each category $F(\quot{X}{\Gamma})$ admits either a zero object, finite products, or an $\Ab$-enrichment then so does $F_G(X)$ (with the $\Ab$-enrichment extended to a $\VAb$-enrichement as necessary). These together will allow us to conclude that the category $F_G(X)$ is additive (in the $\VAb$-enriched sense) whenever each category $F(\quot{X}{\Gamma})$ is and each of the functors $F(\overline{f})$ preserve the relevant structure. The first result we present below is admittedly a special case of Theorems \ref{Thm: Section 2: Equivariant Cat has lims} and \ref{Theorem: Section 2: Equivariant Cat has colims}, but we present it in isolation to give an explicit sense of how these proofs work before diving into maximal generality.

\begin{lemma}\label{Lemma: Section 2: Equivariant Cat has zero objects}
Let $F$ be a pre-equivariant pseudofunctor over $X$ and assume that for all $\Gamma \in \Sf(G)_0$ and for all $f \in \Sf(G)_1$ each category $F(\quot{X}{\Gamma})$ has a zero object $\quot{0}{\Gamma}$ and that each functor $F(\overline{f})$ preserves each corresponding zero object. Then the object
\[
0 = \lbrace \quot{0}{\Gamma} \; | \; \Gamma \in \Sf(G)_0\rbrace
\]
determines a zero object in $F_G(X)$.
\end{lemma}
\begin{proof}
Begin by defining $(0,T_0)$ by setting the collection $0$ as above and defining
\[
T_0 := \lbrace \tau_f^0:F(\overline{f})(\quot{0}{\Gamma^{\prime}}) \xrightarrow{\cong}  \quot{0}{\Gamma} \; | \; f \in \Sf(G)_1 \rbrace
\]
as the $\Vscr$-set of unique ismorphisms witnessing the fact that the $F(\overline{f})$ preserve the zero object. 

We will first verify that $0$ is terminal in $F_G(X)$; that it is initial follows by duality and will be omitted. First let $(A,T_A)$ be an aribtrary object in $F_G(X)_0$. Define the map $!_{A}:(A,T_A) \to (0,T_0)$ by taking the $\Gamma$-wise collection of the terminal maps $!$ on each $\quot{A}{\Gamma}$; explicitly
\[
!_{A} := \lbrace !_{\quot{A}{\Gamma}}:\quot{A}{\Gamma} \to \quot{0}{\Gamma} \; | \; \Gamma \in \Sf(G)_0 \rbrace.
\]
That $!_{A}$ is a morphism in $F_G(X)$ is clear, as the diagrams
\[
\xymatrix{
F(\overline{f})(\quot{A}{\Gamma^{\prime}}) \ar[r]^-{\tau_f^A} \ar[d]_{!_{F(\overline{f})(\quot{A}{\Gamma^{\prime}})}} & \quot{A}{\Gamma} \ar[d]^{!_{\quot{A}{\Gamma}}} \\
F(\overline{f})(\quot{0}{\Gamma^{\prime}}) \ar[r]_{\tau_f^0} & \quot{0}{\Gamma}
}
\]
all commute for all $f \in \Sf(G)_1$ and for all $\Gamma \in \Sf(G)_0$. That this map $!_{A}$ is unique follows from the fact that for any other morphism $P:A \to 0$, we would have by each $\quot{0}{\Gamma}$ being terminal that $\quot{\rho}{\Gamma} = !_{\quot{A}{\Gamma}}$ and hence that $!_{A} = P$. This proves the lemma.
\end{proof}
Following the proof of Lemma \ref{Lemma: Section 2: Equivariant Cat has zero objects} we get the following.
\begin{lemma}
If each category $F(\quot{X}{\Gamma})$ has a terminal (or initial) object and each functor $F(\overline{f})$ preserves it, then $F_G(X)$ has a terminal (or initial) object.
\end{lemma}

For what follows below we will prove when it is possible to deduce that the equivariant category has limits over diagrams of specific shape based on when limits of the same shape exist in each fibre category $F(\quot{X}{\Gamma})$. This will allow us to deduces, say in the cases of equivariant quasi-coherent sheaves, perverse sheaves, local systems, $\Ocal_X$-modules, and other examples where the fibre functors $F(\overline{f})$ are exact/left exact, that the category $F_G(X)$ has all finite limits (or limits of specific type; cf.\@ Corollaries \ref{Cor: Section 2: Limits of shape I in equiv cat} and \ref{Cor: Equivariant cat is finitely complete or cocomplete}).
\begin{Theorem}\label{Thm: Section 2: Equivariant Cat has lims}
Let $F$ be a pre-equivariant pseudofunctor on $X$ and let $I$ be an index category for which there is a diagram $d:I \to F_G(X)$. For each $\Gamma \in \Sf(G)_0$ let $d_{\Gamma}:I \to F(\XGamma)$ denote the diagram functor
\[
\xymatrix{
I \ar[drr]_{d_{\Gamma}} \ar[r]^-{d} & F_G(X) \ar[r]^-{\tilde{\i}} & \prod\limits_{\Gamma \in \Sf(G)_0} F(\XGamma) \ar[d]^{\pi_{\Gamma}} \\
 & & F(\XGamma)
}
\]
and assume that for each $\Gamma \in \Sf(G)_0$ the limit $\lim d_{\Gamma}$ exists in $F(\XGamma)$. Furthermore, assume that for every morphism $f \in \Sf(G)_1$ with $f:\Gamma \to \Gamma^{\prime}$ there is an isomorphism $\theta_f$ witnessing the limit preservation
\[
F(\overline{f})\left(\lim_{\substack{\longleftarrow \\ i \in I}}d_{\Gamma}(i)\right) \cong \lim_{\substack{\longleftarrow \\ i \in I}} F(\overline{f})(d_{\Gamma}(i)).
\]
Then the diagram $d$ admits a limit in $F_G(X)$.
\end{Theorem}
\begin{proof}
Let $\Gamma \xrightarrow{f} \Gamma^{\prime} \xrightarrow{g} \Gamma^{\prime\prime}$ be a pair of composable arrows in $\Sf(G)$. Recall the natural isomorphism $\phi_{f,g}:F(\overline{f}) \circ F(\overline{g}) \Rightarrow F(\overline{g}\circ \overline{f})$; we will abuse notation and write
\[
\phi_{f,g}^{A_i}:F(\overline{f})\big(F(\overline{g})d_{\Gamma^{\prime\prime}}(i)\big) \xrightarrow{\cong} F(\overline{g} \circ \overline{f})d_{\Gamma^{\prime\prime}}(i)
\]
in place of the more cumbersome $\phi_{f,g}^{d_{\Gamma^{\prime\prime}}(i)}$. Set the object $\Vscr$-set as
\[
A := \left\lbrace \lim_{\substack{\longleftarrow \\ i \in I}} d_{\Gamma}(i) \; | \; \Gamma \in \Sf(G)_0 \right\rbrace.
\]
To see how to define the transition isomorphisms $\tau_f^A$ for all $f \in \Sf(G)_1$, recall that for all $f \in \Sf(G)_1$ there are natural isomorphisms 
\[
\theta_f: F(\overline{f}) \circ \lim(d_{\Codom f}) \xrightarrow{\cong} \lim(F(\overline{f}) \circ d_{\Codom f}).
\] 
Then for all $i \in I_0$ we have commuting diagrams
\[
\xymatrix{
\lim(F(\overline{f})\circ d_{\Gamma^{\prime}}) \ar[d]_{\pi_{i}^{\prime}} \ar@{-->}[rr]^-{\exists!\lim(\tau_f^{A_i})} & & \lim(d_{\Gamma}) \ar[d]^{\pi_i} \\
F(\overline{f})(d_{\Gamma^{\prime}}(i)) \ar[rr]_-{\tau_f^{A_i}} & & d_{\Gamma(i)}
}
\]
and, because each $\tau_f^{A_i}$ is an isomoprhism, so is $\lim(\tau_f^{A_i})$. We thus define the transition isomorphism $\tau_f^{A}:F(\overline{f})\quot{A}{\Gamma^{\prime}} \xrightarrow{\cong} \quot{A}{\Gamma}$ via the composition:
\[
\xymatrix{
F(\overline{f})(\lim(d_{\Gamma^{\prime}})) \ar[r]^-{\theta_f} \ar[dr]_{\tau_f^A} & \lim(F(\overline{f}) \circ d_{\Gamma^{\prime}}) \ar[d]^{\lim(\tau_f^{A_i})} \\
 & \lim(d_{\Gamma})
}
\]
Note that because $\theta_f$ and $\lim(\tau_f^{A_i})$ are isomorphisms so is the composite $\tau_f^A$. 

We now need to verify the cocycle condition
\[
\tau_f^A \circ F(\overline{f})\tau_g^A = \tau_{g \circ f}^{A} \circ \phi_{f,g}^{\lim A_i}
\]
for the composable pair of morphisms $f$ and $g$. With this in mind we will decompose the isomorphism $\theta_{g \circ f}$ in order to get an algebraic identity that will allow us to prove the cocycle condition. Consider that
\[
\tau_{g \circ f}^{A} = \lim(\tau_{g \circ f}^{A_i}) \circ \theta_{g \circ f}
\]
and by the uniqueness of the ismorphism $\theta_{g \circ f}$ (as it is induced from the universal property of a limit) we find that the diagram (where we recall that $d_{\Gamma^{\prime\prime}}(i) = \quot{A_i}{\Gamma^{\prime\prime}}$ and $\lim \quot{A_i}{\Gamma^{\prime\prime}} = \quot{A}{\Gamma^{\prime\prime}}$)
\[
\xymatrix{
F(\overline{g} \circ \overline{f})(\lim(\quot{A_i}{\Gamma^{\prime\prime}})) \ar[rr]^-{\theta_{g \circ f}}_-{\cong} \ar[d]_{\left(\phi_{f,g}^{\lim A_i}\right)^{-1}}^{\cong} & & \lim(F(\overline{g} \circ \overline{f})(\quot{A_i}{\Gamma^{\prime\prime}})) \\
F(\overline{f})\big(F(\overline{g})(\lim(\quot{A_i}{\Gamma^{\prime\prime}}))\big) \ar[d]_{F(\overline{f})\theta_g}^{\cong} \\
F(\overline{f})\big(\lim F(\overline{g})(\quot{A_i}{\Gamma^{\prime\prime}})\big)\ar[d]_{F(\overline{f})\big(\lim \tau_g^{A_i}\big)}^{\cong} & & \lim F(\overline{f})\big(F(\overline{g})(\quot{A_i}{\Gamma^{\prime\prime}}) \big) \ar[uu]^{\cong}_{\lim \phi_{f,g}^{A_i}} \\
F(\overline{f})\big(\lim (\quot{A_i}{\Gamma^{\prime}})\big) \ar[rr]_-{\theta_f}^-{\cong} & & \lim\big( F(\overline{f})(\quot{A_i}{\Gamma^{\prime}})\big) \ar[u]^{\cong}_{\big(\lim(F(\overline{f})\tau_g^{A_i})\big)^{-1}}
}
\]
commutes in $F(\quot{X}{\Gamma})$.
We now consider the morphism $\psi:\lim\big(F(\overline{g} \circ \overline{f})\quot{A_i}{\Gamma^{\prime\prime}}\big) \to \lim\quot{A_i}{\Gamma}$ defined by the composition:
\[
\xymatrix{
\lim\big(F(\overline{g} \circ \overline{f})(\quot{A_i}{\Gamma^{\prime\prime}})\big) \ar[dd]_{\psi} \ar[rr]^-{\theta_{g \circ f}^{-1}} & & F(\overline{g} \circ \overline{f})\big(\lim (\quot{A_i}{\Gamma^{\prime\prime}}) \big) \ar[d]^{\big(\phi_{f,g}^{\lim A_i}\big)^{-1}} \\ 
 & & F(\overline{f})(F(\overline{g})\lim (\quot{A_i}{\Gamma^{\prime\prime}}) \ar[d]^{F(\overline{f})\tau_{g}^{A}})\\
\lim (\quot{A_i}{\Gamma}) & & F(\overline{f})\big(\lim (\quot{A_i}{\Gamma^{\prime}})\big) \ar[ll]^-{\tau_f^A}
}
\]
Substituting our calculation of $\theta^{-1}_{g \circ f}$ above into the definition of $\psi$ yields
\begin{align*}
\psi&=\tau_f^A \circ F(\overline{f})\tau_g^A \circ \big(\phi_{f,g}^{\lim A_i}\big)^{-1} \circ \theta_{g \circ f}^{-1} \\
&= \tau_f^A \circ F(\overline{f})\tau_g^A \circ F(\overline{f})\theta_g^{-1} \circ F(\overline{f})(\lim \tau_g^{A_i})^{-1} \circ \theta_f^{-1} \\
&\circ \lim(F(\overline{f})\tau_g^{A_i}) \circ \big(\lim \phi_{f,g}^{A_i}\big)^{-1} \\
&= \tau_f^A \circ F(\overline{f})\tau_g^A \circ F(\overline{f})\big(\theta_g^{-1} \circ (\lim \tau_g^{A_i})^{-1}\big) \circ \theta_f^{-1}\\
& \circ \lim(F(\overline{f})\tau_g^{A_i}) \circ \big(\lim \phi_{f,g}^{A_i}\big)^{-1} \\
&= \tau_f^A \circ F(\overline{f})\tau_g^A \circ F(\overline{f})\big(\tau_g^{A}\big)^{-1} \circ \theta_{f}^{-1} \circ \lim(F(\overline{f}) \tau_g^{A_i}) \circ \big(\lim \phi_{f,g}^{A_i}\big)^{-1} \\
&= \tau_f^A \circ \theta_f^{-1} \circ \lim(F(\overline{f})\tau_g^{A_i}) \circ \big(\lim \phi_{f,g}^{A_i}\big)^{-1} \\
&= \lim(\tau_f^{A_i}) \circ \theta_f \circ \theta_f^{-1} \circ \lim(F(\overline{f})\tau_g^{A_i}) \circ \big(\lim \phi_{f,g}^{A_i}\big)^{-1} \\
&= \lim(\tau_f^{A_i}) \circ \lim(F(\overline{f})\tau_g^{A_i}) \circ \big(\lim \phi_{f,g}^{A_i}\big)^{-1} \\
&= \lim\left(\tau_f^{A_i} \circ F(\overline{f})\tau_g^{A_i} \circ (\phi_{f,g}^{A_i})^{-1} \right) \\
&=\lim\left(\tau_{g \circ f}^{A_i} \circ \phi_{f,g}^{A_i} \circ (\phi_{f,g}^{A_i})^{-1} \right) \\
&= \lim\tau_{g \circ f}^{A_i},
\end{align*}
where we canceled the composite $\big(\phi_{f,g}^{\lim A_i}\big)^{-1} \circ \phi_{f,g}^{\lim A_i}$ without mention in the second equality and where we used the universal property of the limit to give the equality
\[
\lim(\tau_f^{A_i}) \circ \lim(F(\overline{f})\tau_g^{A_i}) \circ \big(\lim \phi_{f,g}^{A_i}\big)^{-1} = \lim\big(\tau_f^{A_i} \circ F(\overline{f})\tau_g^{A_i} \circ \phi_{f,g}^{A_i}\big).
\]
This gives that $\tau_f^A \circ F(\overline{f})\tau_g^A \circ \big(\phi_{f,g}^{\lim A_i}\big)^{-1} \circ \theta_{g \circ f}^{-1} = \lim\tau_{g \circ f}^{A_i}$, which in turn implies that
\[
\tau_f^A \circ F(\overline{f})\tau_g^A = \lim\left(\tau_{g \circ f}^{A_i}\right) \circ \theta_{g \circ f} \circ \phi_{f,g}^{\lim A_i} = \tau_{g \circ f}^{A} \circ \phi_{f,g}^{A}.
\]
Therefore it follows that $T_{A}$ satisfies the cocycle condition and hence the pair $(A,T_A)$ is an object in $F_G(X)$.

We now show that $(A,T_A)$ is the object vertex of a cone over the $(A_i, T_{A_i})$, i.e., that there are morphisms $P_i:A \to A_i$ for all $i \in I_0$ that make the diagrams
\[
\xymatrix{
 & A \ar[dr]^{P_j} \ar[dl]_{P_i} & \\
A_i \ar[rr]_{d(\alpha)} & & A_j
}
\]
commute whenever there is a morphism $\alpha \in I(i,j)$. For any $i \in I_0$ we define the $\Vscr$-set $P:A \to A_i$ by
\[
P := \lbrace \quot{\rho_i}{\Gamma}:\quot{A}{\Gamma} \to \quot{A_i}{\Gamma} \; | \; \Gamma \in \Sf(G)_0 \rbrace
\]
where each $\quot{\rho}{\Gamma} \in F(\quot{X}{\Gamma})(\quot{A}{\Gamma},\quot{A_i}{\Gamma})$ is the limit map. That this makes $A$ into a cone over the $A_i$ is clear if $P$ is a morphism in $F_G(X)$, as for each $\Gamma \in \Sf(G)_0$ the corresponding cone identity holds. Thus all we need to do is verify that $P$ is a morphism in $F_G(X)$, i.e., we need to show that for all $f \in \Sf(G)_1$, say with $\Dom f = \Gamma$ and $\Codom f = \Gamma^{\prime}$, the diagram
\[
\xymatrix{
F(\overline{f})(\quot{A}{\Gamma^{\prime}}) \ar[rr]^-{F(\overline{f})\quot{\rho_i}{\Gamma^{\prime}}} \ar[d]_{\tau_f^A}& & F(\overline{f})(\quot{A_i}{\Gamma^{\prime}}) \ar[d]^{\tau_f^{A_i}} \\
\quot{A}{\Gamma} \ar[rr]_-{\quot{\rho_i}{\Gamma}} & & \quot{A_i}{\Gamma}
}
\]
commutes. For this first recall that $\tau_f^A = \lim(\tau_f^{A_i}) \circ \theta_f$ and that the morphism $\lim(\tau_f^{A_i})$ is defined from the commuting diagram
\[
\xymatrix{
\lim\big(F(\overline{f})(\quot{A_i}{\Gamma^{\prime}})\big) \ar[rr]^-{\lim\tau_f^{A_i}} \ar[d]_{\pi_i} & & \lim(\quot{A_i}{\Gamma}) \ar[d]^{\quot{\rho_i}{\Gamma}} \\
F(\overline{f})(\quot{A_i}{\Gamma^{\prime}}) \ar[rr]_-{\tau_f^{A_i}} & & \quot{A_i}{\Gamma}
}
\]
in $F(\quot{X}{\Gamma})$, where $\pi_i:\lim\big(F(\overline{f})\quot{A_i}{\Gamma^{\prime}}\big) \to F(\overline{f})\quot{A_i}{\Gamma^{\prime}}$ is the limit map. Furthermore, from the isomorphism $\theta_f$ we also derive the commuting triangle
\[
\xymatrix{
F(\overline{f})\big(\lim (\quot{A_i}{\Gamma^{\prime}})\big) \ar[drr]_{F(\overline{f})\quot{\rho_i}{\Gamma^{\prime}}} \ar[rr]^-{\theta_f} & & \lim\big(F(\overline{f})(\quot{A_i}{\Gamma^{\prime}})\big) \ar[d]^{\pi_i} \\
	& & F(\overline{f})(\quot{A_i}{\Gamma^{\prime}})
}
\]
in $F(\quot{X}{\Gamma})$. These together give the commuting quadrilateral
\[
\xymatrix{
F(\overline{f})\big(\lim (\quot{A_i}{\Gamma^{\prime}})\big) \ar[drr]_{F(\overline{f})\quot{\rho_i}{\Gamma^{\prime}}} \ar[rr]^-{\theta_f} & &  \lim\big(F(\overline{f})(\quot{A_i}{\Gamma^{\prime}})\big) \ar[rr]^-{\lim\tau_f^{A_i}} \ar[d]_{\pi_i} & & \lim(\quot{A_i}{\Gamma}) \ar[d]^{\quot{\rho_i}{\Gamma}} \\
& & F(\overline{f})(\quot{A_i}{\Gamma^{\prime}}) \ar[rr]_-{\tau_f^{A_i}} & & \quot{A_i}{\Gamma}
}
\]
which gives the equality
\[
\tau_f^{A_i} \circ F(\overline{f})\quo{\rho_i}{\Gamma^{\prime}} = \tau_f^{A_i} \circ \pi_i \circ \theta_f = \quot{\rho_i}{\Gamma} \circ \lim\big(\tau_f^{A_i}\big) \circ \theta_f = \quot{\rho_i}{\Gamma} \circ \tau_f^A.
\]
Thus the diagram
\[
\xymatrix{
	F(\overline{f})(\quot{A}{\Gamma^{\prime}}) \ar[rr]^-{F(\overline{f})\quot{\rho_i}{\Gamma^{\prime}}} \ar[d]_{\tau_f^A}& & F(\overline{f})(\quot{A_i}{\Gamma^{\prime}}) \ar[d]^{\tau_f^{A_i}} \\
	\quot{A}{\Gamma} \ar[rr]_-{\quot{\rho_i}{\Gamma}} & & \quot{A_i}{\Gamma}
}
\]
commutes and so $P_i$ is a morphism in $F_G(X)$. This shows that $A$ is a cone over the $A_i$.

Finally we verify that $A$ is a universal cone over $F_G(X)$. For this let $(C,T_C)$ be a cone over the $A_i$ with morphisms $\Psi_i:C \to A_i$ for which 
\[
\Psi_i = \lbrace \quot{\psi_i}{\Gamma}:\quot{C}{\Gamma} \to \quot{A_i}{\Gamma} \; | \; \Gamma \in \Sf(G)_0\rbrace.
\] 
Because each $\quot{A}{\Gamma}$ is the limit in $F(\quot{X}{\Gamma})$ of the $\quot{A_i}{\Gamma}$, for all $\Gamma \in \Sf(G)_0$ we produce unique morphisms $\quot{\zeta}{\Gamma}:\quot{C}{\Gamma} \to \quot{A}{\Gamma}$ making the diagram
\[
\xymatrix{
 & \quot{C}{\Gamma} \ar@{-->}[d]^{\quot{\zeta}{\Gamma}} \ar@/^/[ddr]^{\quot{\psi_j}{\Gamma}} \ar@/_/[ddl]_{\quot{\psi_i}{\Gamma}} & \\
 & \quot{A}{\Gamma} \ar[dr]_{\quot{\rho_j}{\Gamma}} \ar[dl]^{\quot{\rho_i}{\Gamma}} \\
\quot{A_i}{\Gamma} \ar[rr]_-{d(\alpha)} & & \quot{A_j}{\Gamma}
}
\]
commute for all $i,j \in I_0$ and for all $\alpha \in I(i,j)$. We now will show that the $\Vscr$-set
\[
Z := \lbrace \quot{\zeta}{\Gamma} \; | \; \Gamma \in \Sf(G)_0 \rbrace
\]
is a morphism in $F_G(X)(C,A)$. Note that this will complete the proof of the proposition, as any other cone map $\Phi:C \to A$ will be component-wise equal to $Z$ by the uniqueness of the $\quot{\zeta}{\Gamma}$.

To see that $Z$ is a morphism we must verify that for any $f \in \Sf(G)_1$ the diagram
\[
\xymatrix{
F(\overline{f})(\quot{C}{\Gamma^{\prime}}) \ar[rr]^-{F(\overline{f})\quot{\zeta}{\Gamma^{\prime}}} \ar[d]_{\tau_f^C} & & F(\overline{f})(\quot{A}{\Gamma^{\prime}}) \ar[d]^{\tau_f^A} \\
\quot{C}{\Gamma} \ar[rr]_-{\quot{\zeta}{\Gamma}} & & \quot{A}{\Gamma}
}
\]
commutes where $\Dom f = \Gamma$ and $\Codom f = \Gamma^{\prime}$. To this end, let $i,j \in I_0$ and let $\alpha \in I(i,j)$. Then, because $C$ is a cone to the $A_i$ and because the assignment $d:I \to F_G(X)$ is a functor (so in particular $d(\alpha)$ is a morphism in $F_G(X)$), we find that the diagram
\[
\xymatrix{
 & F(\overline{f})(\quot{C}{\Gamma^{\prime}}) \ar@/^/[dr]^{F(\overline{f})\quot{\psi_j}{\Gamma^{\prime}}} \ar@/_/[dl]_{F(\overline{f})\quot{\psi_i}{\Gamma^{\prime}}} & \\
F(\overline{f})(\quot{A_i}{\Gamma^{\prime}}) \ar[rr]_-{F(\overline{f})\big(d_{\Gamma^{\prime}}(\alpha)\big)} \ar[d]_{\tau_f^{A_i}} & & F(\overline{f})(\quot{A_j}{\Gamma^{\prime}}) \ar[d]^{\tau_f^{A_j}} \\
\quot{A_i}{\Gamma} \ar[rr]_-{d_{\Gamma}(\alpha)} & & \quot{A_j}{\Gamma}
}
\]
commutes in $F(\quot{X}{\Gamma})$. Thus there exists a unique morphism $\xi_f$ making the diagram
\[
\xymatrix{
	& F(\overline{f})(\quot{C}{\Gamma^{\prime}}) \ar@{-->}[d]^{\xi_f} \ar@/^/[ddr]^{\tau_f^{A_j} \circ F(\overline{f})\quot{\psi_j}{\Gamma^{\prime}}} \ar@/_/[ddl]_{\tau_f^{A_i} \circ F(\overline{f})\quot{\psi_i}{\Gamma}^{\prime}} & \\
	& \quot{A}{\Gamma} \ar[dr]_{\quot{\rho_j}{\Gamma}} \ar[dl]^{\quot{\rho_i}{\Gamma}} \\
	\quot{A_i}{\Gamma} \ar[rr]_-{d(\alpha)} & & \quot{A_j}{\Gamma}
}
\]
commute. Using the above diagram together with the identity 
\[
\tau_f^{A_i} \circ F(\overline{f})\quot{\rho_i}{\Gamma^{\prime}} = \quot{\rho_i}{\Gamma} \circ \tau_f^A
\] 
gives that on one hand
\begin{align*}
\quot{\rho_i}{\Gamma} \circ \xi_f &= \tau_f^{A_i} \circ F(\overline{f})\quot{\psi_i}{\Gamma^{\prime}} = \tau_f^{A_i} \circ F(\overline{f})\big(\quot{\rho_i}{\Gamma^{\prime}} \circ \quot{\zeta}{\Gamma^{\prime}}\big) \\
&= \tau_f^{A_i} \circ F(\overline{f})\quot{\rho_i}{\Gamma^{\prime}} \circ F(\overline{f})\quot{\zeta}{\Gamma^{\prime}}  \\
&=\quot{\rho_i}{\Gamma} \circ \tau_f^A \circ F(\overline{f})\quot{\zeta}{\Gamma^{\prime}};
\end{align*}
following this same argument mutatis mutandis we derive also that 
\[
\quot{\rho_j}{\Gamma} \circ \xi_f = \quot{\rho_j}{\Gamma} \circ \tau_f^A \circ F(\overline{f})\quot{\zeta}{\Gamma^{\prime}}.
\] 
On the other hand, using that $\Psi_i$ is an $F_G(X)$-morphism gives that
\[
\quot{\rho_i}{\Gamma} \circ \quot{\zeta}{\Gamma} \circ \tau_f^C = \quot{\psi_i}{\Gamma} \circ \tau_f^C = \tau_f^{A_i} \circ F(\overline{f})\quot{\psi_i}{\Gamma^{\prime}} = \quot{\rho_i}{\Gamma} \circ \xi_f.
\]
As before, we also derive that $\quot{\rho_j}{\Gamma} \circ \quot{\zeta}{\Gamma} \circ \tau_f^C = \quot{\rho_j}{\Gamma} \circ \xi_f$. Putting all these identities together gives us that both $\quot{\zeta}{\Gamma} \circ \tau_f^C$ and $\tau_f^A \circ F(\overline{f})\quot{\zeta}{\Gamma^{\prime}}$ produce fillers for the universal cone diagram in the same way that $\xi_f$ does; consequently, by the universal property of $\quot{A}{\Gamma}$ as the limit over the $\quot{A_i}{\Gamma}$, we derive that
\[
\tau_f^{A} \circ F(\overline{f})\quot{\zeta}{\Gamma^{\prime}} = \xi_f = \quot{\zeta}{\Gamma} \circ \tau_f^C.
\]
This gives the commutativity of the diagram
\[
\xymatrix{
F(\overline{f})(\quot{C}{\Gamma^{\prime}}) \ar[rr]^-{F(\overline{f})\quot{\zeta}{\Gamma^{\prime}}} \ar[d]_{\tau_f^C} & & F(\overline{f})(\quot{A}{\Gamma^{\prime}}) \ar[d]^{\tau_f^A} \\
\quot{C}{\Gamma} \ar[rr]_-{\quot{\zeta}{\Gamma}} & & \quot{A}{\Gamma}
}
\]
and hence shows that $(A,T_A)$ is a universal cone over the $(A_i, T_{A_i})$. Because this shows that $d$ has a limit in $F_G(X)$, this completes the proof of the theorem.
\end{proof}
\begin{Theorem}\label{Theorem: Section 2: Equivariant Cat has colims}
Let $F$ be a pre-equivariant pseudofunctor on $X$ and let $I$ be an index category for which there is a diagram $d:I \to F_G(X)$. For each $\Gamma \in \Sf(G)_0$ let $d_{\Gamma}:I \to F(\XGamma)$ denote the diagram functor
\[
\xymatrix{
	I \ar[drr]_{d_{\Gamma}} \ar[r]^-{d} & F_G(X) \ar[r]^-{\tilde{\i}} & \prod\limits_{\Gamma \in \Sf(G)_0} F(\XGamma) \ar[d]^{\pi_{\Gamma}} \\
	& & F(\XGamma)
}
\]
and assume that for each $\Gamma \in \Sf(G)_0$ the colimit $\colim d_{\Gamma}$ exists in $F(\XGamma)$. Furthermore, assume that for every morphism $f \in \Sf(G)_1$ with $f:\Gamma \to \Gamma^{\prime}$ there is an isomorphism $\theta_f$ witnessing the limit preservation
\[
F(\overline{f})\left(\lim_{\substack{\longrightarrow \\ i \in I}}d_{\Gamma}(i)\right) \cong \lim_{\substack{\longrightarrow \\ i \in I}} F(\overline{f})(d_{\Gamma}(i)).
\]
Then the diagram $d$ admits a colimit in $F_G(X)$.
\end{Theorem}
\begin{proof}
Dualize the proof of Theorem \ref{Thm: Section 2: Equivariant Cat has lims}.
\end{proof}

\begin{corollary}\label{Cor: Section 2: Limits of shape I in equiv cat}
If $I$ is an index category and each category $F(\quot{X}{\Gamma})$ has all limits of shape $I$ and if each functor $F(\overline{f})$ preserves these limits, then $F_G(X)$ has all limits of shape $I$ as well. Dually, if each category $F(\quot{X}{\Gamma})$ has all colimits of shape $I$ and if each functor $F(\overline{f})$ preserves these colimits, then $F_G(X)$ has all colimits of shape $I$ as well.
\end{corollary}
\begin{proof}
Use Theorem \ref{Thm: Section 2: Equivariant Cat has lims} for every diagram $d:I \to \Sf(G)$, as such a diagram gives rise to a compatable family of diagrams $d_{\Gamma}:I \to \Sf(G)$ by way of projection to the $\Gamma$-th component.
\end{proof}
\begin{corollary}\label{Cor: Equivariant cat is finitely complete or cocomplete}
If each category $F(\quot{X}{\Gamma})$ is (finitely) complete and each functor $F(\overline{f})$ is (finitely) continuous then $F_G(X)$ is (finitely) complete. Dually, $F_G(X)$ is (finitely) cocomplete if each category $F(\quot{X}{\Gamma})$ is (finitely) cocomplete and if each functor $F(\overline{f})$ is finitely cocontinuous. 
\end{corollary}
\begin{proof}
Use Corollary \ref{Cor: Section 2: Limits of shape I in equiv cat} with the categories $I = \emptyset$ to give terminal objects, $I = \coprod_{s \in S} \mathbbm{1}$ (where $S$ is an indexing $\Vscr$-set) to give products indexed by $S$, and $I = \Lambda_2^2$ to give pullbacks. Invoking now the classical categorical result that a category $\Cscr$ admits all limits if and only if it admits products, terminal objects, and pullbacks we find that $F_G(X)$ is complete. This proves the corollary.
%
%
\end{proof}
\begin{remark}
The finiteness assumptions above are used exclusively when determining which indexing $\Vscr$-sets $S$ we are allowed to use. If one only knows of the existence of finite limits or colimits, we restrict our attention to those $\Vscr$-sets of finite cardinality.
\end{remark}
\begin{lemma}\label{Lemma: Section 2: Equiv cat is VAb enriched}
Let $F$ be a pre-equivariant pseudofunctor such that each category $F_G(X)$ is $\VAb$ enriched. Then $F_G(X)$ is $\VAb$ enriched as well.
\end{lemma}
\begin{proof}
Note that since each $F(\quot{X}{\Gamma})$ is $\VAb$ enriched, each such category is $\Vscr$-small; appealing to Lemma \ref{Lemma: Section 2: Pseudofunctor values in Vlocally small means Equiv cat Vlocally small} gives that $F_G(X)$ is locally $\Vscr$-small as well. Now fix two morphisms $\Phi, \Psi \in F_G(X)(A,B)$ for some objects $A, B \in F_G(X)_0$. We then define their addition via the equation
\[
\Phi + \Psi := \lbrace \quot{\varphi}{\Gamma} + \quot{\psi}{\Gamma} \; | \; \Gamma \in \Sf(G)_0 \rbrace,
\]
the identity element as
\[
0 := \lbrace \quot{0}{\Gamma} \; | \; \Gamma \in \Sf(G)_0 \rbrace
\]
and their negations via
\[
-\Phi := \lbrace -\quot{\varphi}{\Gamma} \; | \; \Gamma \in \Sf(G)_0 \rbrace.
\]
The group laws required of a $\Vscr$-Abelian group are then routinely verified using the fact that each category $F(\quot{X}{\Gamma})$ has the required algebraic laws at hand. Furthrmore, the bilinearity of composition is verified in the same way. This proves the lemma.
\end{proof}
\begin{corollary}\label{Cor: Section 2: Additive pseudofunctor gives additive equivariat cat}
If the pre-equivariant pseudofunctor $F:\SfResl_G(X)^{\op} \to \fCat$ is additive and locally $\Vscr$-small (cf.\@ Definition \ref{Defn: Additive, triangulated, and locally small pseudofunctors}) then the category $F_G(X)$ is a locally $\Vscr$-small additive category. 
\end{corollary}
\begin{proof}
Corollary \ref{Cor: Section 2: Limits of shape I in equiv cat} shows that $F_G(X)$ admits finite products and a zero object. Local $\Vscr$-smallness is Lemma \ref{Lemma: Section 2: Pseudofunctor values in Vlocally small means Equiv cat Vlocally small}, while the $\VAb$ enrichment is Lemma \ref{Lemma: Section 2: Equiv cat is VAb enriched}. Finally, appealing to a classical result in homological algebra (cf.\@ \cite[Proposition II.9.1]{HiltonStammHA}, for instance) gives that $F_G(X)$ is additive.
\end{proof}
\begin{corollary}\label{Cor: Section 2: Additive complete and cocomplete equivariant cats}
	If each of the categories $F(\quot{X}{\Gamma})$ are complete and cocomplete as well as additive, and if the functors $F(\overline{f})$ preserve limits and colimits, then $F_G(X)$ is additive, complete, and cocomplete as well.
\end{corollary}
\begin{proof}
Corollary \ref{Cor: Section 2: Additive pseudofunctor gives additive equivariat cat} gives the additive structure on $F_G(X)$ while Corollary \ref{Cor: Equivariant cat is finitely complete or cocomplete} gives the completion and cocompletion.
\end{proof}

We now present some lemmas that allow us to deduce some cases when the equivariant categories $F_G(X)$ are Abelian. Of particular use here is that this gives proofs that the categories of equivariant $\mathcal{D}$-modules, local systems, and quasi-coherent ({\'e}tale) sheaves on $X$ are Abelian (cf.\@ Corollary \ref{Cor: Section 2: Additive equivriant derived, local system, and qcoh cats} below). 

While illuminating the structure of Abelian equivariant categories, we also present a proposition about descending regular categories through these equivariant constructions. While this will not be crucial for us, it will show the use of equivariant categories in developing a notion of equivariant categorical logic.

We proceed by recalling that a category $\Cscr$ is Abelian if and only if it is additive, (finitely) complete, (finitely) cocomplete, and has the property that every monomorphism is a kernel and every epimorphism is a cokernel. We have already seen in Corollary \ref{Cor: Section 2: Additive complete and cocomplete equivariant cats} that under mild assumptions the category $F_G(X)$ is additive, complete, and cocomplete. Thus all we need to show are the regularity conditions, i.e., the fact that every monomorphism is a kernel and that every epimorphism is a cokernel.
\begin{lemma}\label{Lemma: Section 2: KErnels descend}
Assume that $F$ is a pre-equivariant pseudofunctor, that each functor $F(\overline{f})$ is continuous and cocontinuous, and that each category $F(\quot{X}{\Gamma})$ admits equalizers and cokernel pairs of morphisms. Then if every monic in each fibre category $F(\quot{X}{\Gamma})$ is a kernel, every monic in $F_G(X)$ is a kernel.
\end{lemma}
\begin{proof}
Begin by recalling that a morphism $\Psi$ is a kernel of a morphism $\Phi:A \to B$ if and only if there is an isomorphism $\Psi \cong \eq(I_0, I_1)$, where $I_0$ and $I_1$ are the cokernel pair of $\Phi$ that appear in the pushout below:
\[
\xymatrix{
\Dom \Phi \ar[rr]^-{\Phi} \ar[d]_{\Phi} & & \Codom \Phi \ar[d]^{I_1} \\
\Codom \Phi \ar[rr]_-{I_0} & & \Codom \Phi \coprod_{\Dom \Phi} \Codom \Phi
}
\]

To proceed we let $P = \lbrace \quot{\rho}{\Gamma} \; | \; \Gamma \in \Sf(G)_0\rbrace:A \to B$ be a monomorphism in $F_G(X)$ and we let the pair $(I_0,I_1)$ be the cokernel pair of $P$ arising from the pushout
\[
\xymatrix{
A \ar[r]^-{P} \ar[d]_{P} & B \ar[d]^{I_1} \\
B \ar[r]_-{I_0} & B \coprod_A B
}
\]
in $F_G(X)$. From Theorem \ref{Theorem: Section 2: Equivariant Cat has colims} we know that
\[
B \coprod_A B = \left\lbrace \quot{B}{\Gamma} \coprod_{\quot{A}{\Gamma}} \quot{B}{\Gamma} \; : \; \Gamma \in \Sf(G)_0 \right\rbrace
\]
and hence each $I_k$ has the form
\[
I_k = \left\lbrace \quot{i_k}{\Gamma}:\quot{B}{\Gamma} \to \quot{B}{\Gamma} \coprod_{\quot{A}{\Gamma}} \quot{B}{\Gamma} \; : \; \Gamma \in \Sf(G)_0 \right\rbrace,
\]
where $k \in \lbrace 0, 1 \rbrace$ and where each pair $(\quot{i_0}{\Gamma},\quot{i_1}{\Gamma})$ is the cokernel pair of $\quot{\rho}{\Gamma}$ in $F(\quot{X}{\Gamma})$. Using Theorem \ref{Thm: Section 2: Equivariant Cat has lims} we also find that the morphism $P$, as a monic in $F_G(X)$, has the property that each $\quot{\rho}{\Gamma}$ is monic in $F(\quot{X}{\Gamma})$. Because each monomorphism in each fibre category $F(\quot{X}{\Gamma})$ is a kernel, there is a unique isomorphism $\quot{\theta}{\Gamma}:\quot{X}{\Gamma} \xrightarrow{\cong} \Eq(\quot{i_0}{\Gamma},\quot{i_1}{\Gamma})$ making the diagram
\[
\xymatrix{
\quot{A}{\Gamma} \ar@{-->}[d]_{\exists!\quot{\theta}{\Gamma}} \ar[rr]^-{\quot{\rho}{\Gamma}} & & \quot{B}{\Gamma} \ar@{=}[d] \\
\Eq(\quot{i_0}{\Gamma},\quot{i_1}{\Gamma}) \ar[rr]_-{\eq(\quot{i_0}{\Gamma}, \quot{i_1}{\Gamma})} & & \quot{B}{\Gamma}
}
\]
commute in $F(\quot{X}{\Gamma})$; note that $\Eq(\quot{i_0}{\Gamma},\quot{i_1}{\Gamma})$ is the equalizer object of the $\quot{i_k}{\Gamma}$ and $\eq(\quot{i_0}{\Gamma},\quot{i_1}{\Gamma})$ is the canonical equalizer morphism. Using Theorem \ref{Thm: Section 2: Equivariant Cat has lims} again, we have that the collection
\[
\Eq(I_0,I_1) := \lbrace \Eq(\quot{i_0}{\Gamma},\quot{i_1}{\Gamma}) \; | \; \Gamma \in \Sf(G)_0 \rbrace
\]
determines the object assignment of the equalizer $\Eq(I_0,I_1)$ in $F_G(X)$ and the $\Vscr$-set
\[
\eq(I_0,I_1) := \lbrace \eq(\quot{i_0}{\Gamma},\quot{i_1}{\Gamma}) \;|\; \Gamma \in \Sf(G)_0\rbrace
\]
determines the equalizer morphism $\eq(I_0,I_1):\Eq(I_0,I_1) \to B$ in $F_G(X)$. Thus we will be done if we can show that the $\Vscr$-set
\[
\Theta := \lbrace \quot{\theta}{\Gamma} \; | \; \Gamma \in \Sf(G)_0 \rbrace
\]
is a morphism $\Theta:A \to \Eq(I_0, I_1)$ in $F_G(X)$, as then $\Theta$ is the unique morphism making the square
\[
\xymatrix{
A \ar[r]^{P} \ar@{-->}[d]_{\exists!\Theta} & B \ar@{=}[d] \\
\Eq(I_0,I_1) \ar[r]_-{\eq(I_0,I_1)} & B
}
\]
commute.

For this fix a morphism $f \in \Sf(G)_1$ and let $\Dom f = \Gamma$ and $\Codom f = \Gamma^{\prime}$. We now must show that $\quot{\theta}{\Gamma} \circ \tau_f^A \overset{?}{=} \tau_f^{\Eq(I_0,I_1)} \circ F(\overline{f})\eq(\quot{i_0}{\Gamma^{\prime}},\quot{i_1}{\Gamma^{\prime}})$. For this we calculate that
\begin{align*}
\eq(\quot{i_0}{\Gamma},\quot{i_1}{\Gamma}) \circ \quot{\theta}{\Gamma} \circ \tau_f^A &= \quot{\rho}{\Gamma} \circ \tau_f^A = \tau_f^B \circ F(\overline{f})\quot{\rho}{\Gamma^{\prime}} \\
&= \tau_f^B \circ F(\overline{f})\big(\eq(\quot{i_0}{\Gamma^{\prime}},\quot{i_1}{\Gamma^{\prime}}) \circ \quot{\theta}{\Gamma^{\prime}}\big) \\
&= \tau_f^B \circ F(\overline{f})\eq(\quot{i_0}{\Gamma^{\prime}},\quot{i_1}{\Gamma^{\prime}}) \circ F(\overline{f})\quot{\theta}{\Gamma{^\prime}} \\
&= \eq(\quot{i_0}{\Gamma},\quot{i_1}{\Gamma}) \circ \tau_f^{\Eq(I_0,I_1)} \circ F(\overline{f})\quot{\theta}{\Gamma^{\prime}}.
\end{align*}
Using that $\eq(\quot{i_0}{\Gamma},\quot{i_1}{\Gamma})$ is monic then allows us to deduce that
\[
\quot{\theta}{\Gamma} \circ \tau_f^A = \tau_f^{\Eq(I_0,I_1)} \circ F(\overline{f})\quot{\theta}{\Gamma},
\]
which proves that $\Theta \in F_G(X)(A,\Eq(I_0,I_1))$. This proves the lemma.
\end{proof}
By dualizing the above lemma, we also have the result below for cokernels and epimorphisms.
\begin{lemma}\label{Lemma: Section 2: Cokernels descend}
Assume that $F$ is a pre-equivariant pseudofunctor, that each functor $F(\overline{f})$ is continuous and cocontinuous, and that each category $F(\quot{X}{\Gamma})$ admits coequalizers and kernels pairs of morphisms. Then if every epimorphism in each fibre category $F(\quot{X}{\Gamma})$ is a cokernel, every epimorphism in $F_G(X)$ is a cokernel.
\end{lemma}
This leads us to the fact that Abelian categories descend through the equivariance, so if each fibre category $F(\quot{X}{\Gamma})$ is Abelian and if each fibre functor $F(\overline{f})$ is exact then $F_G(X)$ is Abelian.
\begin{proposition}\label{Prop: Section 2: Abelian Descends}
Let $F$ be a pre-equivariant pseudofunctor such that each category $F(\quot{X}{\Gamma})$ is Abelian and for which each functor $F(\overline{f})$ is exact. Then the category $F_G(X)$ is Abelian as well.
\end{proposition}
\begin{proof}
By Corollary \ref{Cor: Section 2: Additive complete and cocomplete equivariant cats} it follows that the category $F_G(X)$ is complete, cocomplete, and additive. To see that it has the necessary regularity property, simply use Lemmas \ref{Lemma: Section 2: KErnels descend} and \ref{Lemma: Section 2: Cokernels descend} to derive that every monic and epic in $F_G(X)$ is a kernel or cokernel, respectively.
\end{proof}

As a corollary, we have that the categories of equivariant $D$-modules, local systems, and quasicoherent sheaves are all Abelian.
\begin{corollary}\label{Cor: Section 2: Additive equivriant derived, local system, and qcoh cats}
The categories $D_G^b(X)$, $\DGMod(X)$, $\Loc_G(X)$, $\Per_G(X)$, and $\QCoh_G(X)$ are all additive. Moreover, the categories $\DGMod(X)$, $\Loc_G(X)$, $\Per_G(X)$, and $\QCoh_G(X)$ are Abelian.
\end{corollary}
\begin{proof}
That all four categories are additive follows from Corollary \ref{Cor: Section 2: Additive pseudofunctor gives additive equivariat cat}. Because each morphism $f \in \Sf(G)_1$ is smooth (and hence flat), so is the morphism $f \times \id_X: \Dom f \times X \to \Codom f \times X$; thus we get from a standard quotient argument that the morphism $\overline{f}:\quot{X}{\Gamma} \to \quot{X}{\Gamma^{\prime}}$ is smooth as well. This then implies that each map $\overline{f}$ is flat as well, and hence gives that each functor $\overline{f}^{\ast}$ is exact on the level of {\'e}tale sheaves. Thus we can apply Corollary \ref{Cor: Section 2: Additive complete and cocomplete equivariant cats} to give the finite completion and finite cocompletion of $\DGMod(X), \Loc_G(X)$, $\Per_G(X)$, and $\QCoh_G(X)$. Finally, Proposition \ref{Prop: Section 2: Abelian Descends} gives the Abelian structure for the categories $\DGMod(X), \Loc_G(X), \Per_G(X)$ and $\QCoh_G(X)$, as each of the categories $\DMod(\quot{X}{\Gamma}), \Loc(\quot{X}{\Gamma}), \Per_G(X)$, and $\QCoh(\quot{X}{\Gamma})$ are Abelian for all $\Gamma \in \Sf(G)_0$.
\end{proof}
The same proof may be used to show that the category $\Shv_G(X,R)$ of ({\'e}tale) sheaves of (left) modules over a fixed ring $R$ (or modules over a fixed sheaf of rings $\CalO_X$) is Abelian; similarly, the category $\Shv_G(X;\overline{\Q}_{\ell})$ of $\ell$-adic sheaves is Abelian. We record these results in the corollary below.
\begin{corollary}
The categories $\OXMod_G$\index[notation]{OXMod@$\OXMod_G$} of equivariant $\Ocal_X$-modules, $\Shv_G(X; R)$\index[notation]{SheavesofRmodules@$\Shv_G(X,R)$} of equivariant sheaves of $R$-modules, $\Shv_G(X;\overline{\Q}_{\ell})$ of equivariant $\ell$-adic sheaves, and $\Per_G(X;\overline{\Q}_{\ell})$ of equivariant $\ell$-adic perverse sheaves are Abelian.
\end{corollary}

We now begin providing cases when the equivariant categories $F_G(X)$ are monoidal. This will allow us to not only give a foundation for an equivariant monoidal theory, but also give us tools to study equivariant derived monoidal theory and equivariant dg-categories over the geometric objects about which we care. The main use that I can forsee is in the homological algebraic setting, as it gives us tools with which to construct the six-functor-formalism of the equivariant derived category. I expect us to be able to use this in noncommutative algebraic geometry (following the foundations and methodology of \cite{Keller}, \cite{Orlov_2018}, \cite{toën2014derived}, \cite{ToenDAGDQ}, among others, which define and work with noncommutative schemes in dg-categorical language) to give an equivariant notion of noncommutative (derived) algebraic geometry. For the moment, however, we will simply prove a class of cases that allow us to deduce that the equivariant category $F_G(X)$ is monoidal.
\begin{definition}\label{Defn: Section 2: Monoidal Preequivariant Pseudofunctor}
	Let $F:\SfResl_G(X)^{\op} \to \fCat$ be a pre-equivariant pseudofunctor. We say that $F$ is {monoidal}\index{Pre-equivariant Pseudofunctor! Monoidal} if the following hold:
	\begin{itemize}
		\item For all $\Gamma \in \Sf(G)_0$, each category $F(\quot{X}{\Gamma})$ is a monoidal category \\$\big(F(\quot{X}{\Gamma}), \os{\Gamma}{\otimes}, \quot{I}{\Gamma}, \quot{\alpha}{\Gamma}, \quot{\lambda}{\Gamma}, \quot{\rho}{\Gamma}\big)$;
		\item For all $f:\Gamma \to \Gamma^{\prime}$ in $\Sf(G)_1$, $f$ is monoidal in the sense that there are natural isomorphisms
		\[
		\theta_f^{A,B}:F(\overline{f})\big(A \os{\Gamma^{\prime}}{\otimes} B\big) \xrightarrow{\cong} F(\overline{f})A \os{\Gamma}{\otimes} F(\overline{f})B
		\]
		which are pseudofunctorial in the sense that for all composable pairs $\Gamma \xrightarrow{f} \Gamma^{\prime} \xrightarrow{g} \Gamma^{\prime\prime}$ in $\Sf(G)$, the pasting diagram
		\[
		\begin{tikzcd}
		F(\quot{X}{\Gamma^{\prime\prime}}) \times F(\quot{X}{\Gamma^{\prime\prime}}) \ar[d]{}{F(\overline{g}) \times F(\overline{g})} \ar[dd, bend right = 80, swap, ""{name = LL}]{}{F(\overline{g} \circ \overline{f}) \times F(\overline{g}\circ\overline{f})} \ar[rrr, ""{name = UU}]{}{\os{\Gamma^{\prime\prime}}{\otimes}} & & & F(\quot{X}{\Gamma^{\prime\prime}}) \ar[dd, bend left = 80, ""{name = RR}]{}{F(\overline{g} \circ \overline{f})} \ar[d,swap]{}{F(\overline{g})}\\
		F(\quot{X}{\Gamma^{\prime}}) \times F(\quot{X}{\Gamma^{\prime}}) \ar[d]{}{F(\overline{f}) \times F(\overline{f})} \ar[rrr, ""{name = M}]{}[description]{\os{\Gamma^{\prime}}{\otimes}} & & & F(\quot{X}{\Gamma^{\prime}}) \ar[d]{}{F(\overline{f})} \\
		F(\quot{X}{\Gamma}) \times F(\quot{X}{\Gamma}) \ar[rrr, swap, ""{name = BB}]{}{\os{\Gamma}{\otimes}} & & & F(\quot{X}{\Gamma})
		\ar[from = 2-4, to = RR, Rightarrow, shorten <=2pt, shorten >= 2pt]{}{\phi_{f,g}}
		\ar[from = LL, to = 2-1, swap, Rightarrow, shorten <= 2pt, shorten >= 2pt]{}{\phi_{f,g}^{-1} \times \phi_{f,g}^{-1}}
		\ar[from = UU, to = M, Rightarrow, shorten <= 4pt, shorten >= 4pt]{}{\theta_g}
		\ar[from = M, to = BB, Rightarrow, shorten <= 8pt, shorten >= 4pt]{}{\theta_f}
		\end{tikzcd}
		\]
		is equal to the $2$-cell:
		\[
		\begin{tikzcd}
		F(\quot{X}{\Gamma^{\prime\prime}}) \times F(\quot{X}{\Gamma^{\prime\prime}}) \ar[r, ""{name = U}]{}{\os{\Gamma^{\prime\prime}}{\otimes}} \ar[d, swap]{}{F(\overline{g}\circ\overline{f}) \times F(\overline{g} \circ \overline{f})} & F(\quot{X}{\Gamma^{\prime\prime}}) \ar[d]{}{F(\overline{g} \circ \overline{f})} \\
		F(\quot{X}{\Gamma}) \times F(\quot{X}{\Gamma}) \ar[r, swap, ""{name = L}]{}{\os{\Gamma}{\otimes}} & F(\quot{X}{\Gamma}) \ar[from = U, to = L, Rightarrow, shorten >= 4pt, shorten <= 4pt]{}{\theta_{g \circ f}}
		\end{tikzcd}
		\]
		\item The functors $F(\overline{f})$ preserve the units pseudofunctorially in the sense that there for all $f \in \Sf(G)_1$, there is an isomorphism
		\[
		\sigma_f:F(\overline{f})\quot{I}{\Gamma^{\prime}} \xrightarrow{\cong} \quot{I}{\Gamma}
		\]
		such that for all composable pairs $\Gamma \xrightarrow{f} \Gamma^{\prime} \xrightarrow{g} \Gamma^{\prime\prime}$ in $\Sf(G)$, we have that
		\[
		\sigma_f \circ F(\overline{f})\sigma_g = \sigma_{g \circ f} \circ \phi_{f,g}.
		\]
	\end{itemize}
\end{definition}
\begin{remark}
	The functors we have called monoidal in this paper are frequently called {strong monoidal}\index{Monoidal functor}\index{Strong monoidal functor|see {Monoidal Functor}} in the category-theoretic literature (cf.\@ \cite{AguiMahajan.}, \cite{KellyStreet}, \cite{Leinster}, \cite{Yetter}, for instance). We have chosen to drop the adjective ``strong'' for a few reasons: First, we will need the isomorphisms for each use of a monoidal functor in this paper; second, we do not want to confuse strong monoidal functors with monoidal functors that have tensorial strengths; and third, having monoidal functors be implicitly lax is not something that our descent formalism will allow. We need to be working with tensorial transition isomorphisms, not simply transition maps; as such, we will need the full strength of a monoidal functor which preserves the tensor and units up to natural isomorphism.
\end{remark}
\begin{Theorem}\label{Theorem: Section 2: Monoidal preequivariant pseudofunctor gives monoidal equivariant cat}
	Let $F:\SfResl_G(X)^{\op} \to \fCat$ be a monoidal pre-equivariant pseudofunctor. Then $F_G(X)$ is a monoidal category.
\end{Theorem}
\begin{proof}
	First note that if we can define the functor $\otimes,$ the object $I$, and the natural isomorphisms $\lambda, \rho,$ and $\alpha$ which $\Gamma$-wise agree with the functors, objects, and transformations in the fibre category
	\[
	\left(F(\quot{X}{\Gamma}), \os{\Gamma}{\otimes}, \quot{I}{\Gamma}, \quot{\lambda}{\Gamma}, \quot{\rho}{\Gamma}, \quot{\alpha}{\Gamma}\right)
	\]
	we will be done, as all triangle and pentagonal identities will hold $\Gamma$-wise, which is how we we compute these identities. We begin by defining the monoidal functor $\otimes:F_G(X) \times F_G(X) \to F_G(X)$. Define the functor on objects by setting, for $(A,T_A)$ and $(B,T_B)$ objects in $F_G(X)$,
	\[
	A \otimes B := \left\lbrace \quot{A}{\Gamma} \os{\Gamma}{\otimes} \quot{B}{\Gamma} \; \big| \; \Gamma \in \Sf(G)_0, \quot{A}{\Gamma} \in A, \quot{B}{\Gamma} \in B \right\rbrace.
	\]
	To define the transition isomorphisms on $A \otimes B$ first fix an $f \in \Sf(G)_1$ and write $\Dom f = \Gamma$, $\Codom f = \Gamma^{\prime}$. Then we define the morphism $\tau_f^{A \otimes B}$ via the composite
	\[
	\xymatrix{
		F(\overline{f})\left(\quot{A}{\Gamma^{\prime}} \os{\Gamma^{\prime}}{\otimes} \quot{B}{\Gamma^{\prime}} \right) \ar[rr]^-{\theta_f^{A,B}} \ar[drr]_{\tau_f^{A \otimes B}} & & F(\overline{f})(\quot{A}{\Gamma^{\prime}}) \os{\Gamma}{\otimes} F(\overline{f})(\quot{B}{\Gamma^{\prime}}) \ar[d]^{\tau_f^A \otimes \tau_f^B} \\
		& & \quot{A}{\Gamma} \os{\Gamma}{\otimes} \quot{B}{\Gamma}
	}
	\]
	where $\theta_f^{A,B}$ is used as an abuse of notation for
	\[
	\theta_f^{A,B} := \theta_f^{\quot{A}{\Gamma^{\prime}},\quot{B}{\Gamma^{\prime}}}.
	\]
	We then define the transition isomorphisms via the definition
	\[
	T_{A \otimes B} := \left\lbrace \tau_f^{A \otimes B} \; | \; f \in \Sf(G)_1 \right\rbrace. 
	\]
	
	To verify that the transition isomorphisms satisfy the cocycle condition, let $\Gamma \xrightarrow{f} \Gamma^{\prime} \xrightarrow{g} \Gamma^{\prime\prime}$ be a pair of composable morphisms in $\Sf(G)$. Now, using the pseudofunctoriality of the $\otimes$ functors together with the naturality of the $\theta_f$'s gives that
	\begin{align*}
	\tau_{g \circ f}^{A \otimes B} &= \left(\tau_{g \circ f}^{A} \os{\Gamma}{\otimes} \tau_{g \circ f}^{B}\right) \circ \theta_{g \circ f}^{A,B} \\
	&= \left(\big(\tau_f^A \circ F(\overline{f})\tau_g^A\big) \os{\Gamma}{\otimes}\big(\tau_f^B \circ F(\overline{f})\tau_g^B \big)\right) \circ \phi_{f,g}^{-1} \circ \theta_{g \circ f}^{A,B} \\
	&= \left(\big(\tau_f^A \circ F(\overline{f})\tau_g^A\big) \os{\Gamma}{\otimes}\big(\tau_f^B \circ F(\overline{f})\tau_g^B \big)\right) \\
	&\circ \phi_{f,g}^{-1} \circ \phi_{f,g} \circ \theta_{f}^{A,B} \circ F(\overline{f})\theta_g^{A,B} \circ \phi_{f,g}^{-1} \\
	&=\left(\big(\tau_f^A \circ F(\overline{f})\tau_g^A\big) \os{\Gamma}{\otimes}\big(\tau_f^B \circ F(\overline{f})\tau_g^B \big)\right) \circ \theta_{f}^{A,B} \circ F(\overline{f})\theta_g^{A,B} \circ \phi_{f,g}^{-1} \\
	&= \left(\tau_f^A \os{\Gamma}{\otimes} \tau_f^B\right) \circ \left(F(\overline{f})\tau_g^A \os{\Gamma}{\otimes} F(\overline{f})\tau_g^B\right) \circ \theta_{f}^{A,B} \circ F(\overline{f})\theta_{g}^{A,B} \circ \phi_{f,g}^{-1} \\
	&= \left(\tau_f^A \os{\Gamma}{\otimes} \tau_f^B\right) \circ \theta_f^{A,B} \circ F(\overline{f})\left(\tau_g^A \os{\Gamma^{\prime}}{\otimes} \tau_g^B\right) \circ F(\overline{f})\theta_{g}^{A,B} \circ \phi_{f,g}^{-1} \\
	&= \tau_f^{A \otimes B} \circ F(\overline{f})\tau_{g}^{A \otimes B} \circ \phi_{f,g}^{-1}.
	\end{align*}
	Thus we have that $\tau_{g \circ f}^{A \otimes B} \circ \phi_{f,g} = \tau_f^{A\otimes B} \circ F(\overline{f})\tau_g^{A \otimes B}$ and so 
	\[
	\big((A,T_A),(B,T_B)\big) \mapsto (A \otimes B, T_{A \otimes B})
	\] 
	determines the object assignment of the $\otimes$ functor.
	
	To define $\otimes$ on morphisms, fix a morphism $\Phi \in F_G(X)(A,B)$ and a morphism $\Psi \in F_G(X)(A^{\prime},B^{\prime})$. We then define the morphism \[
	\Phi \otimes \Psi:A \otimes A^{\prime} \to B \otimes B^{\prime}
	\] 
	by
	\[
	\Phi \otimes \Psi := \left\lbrace \quot{\varphi}{\Gamma} \os{\Gamma}{\otimes} \quot{\psi}{\Gamma} \; \big| \; \Gamma \in \Sf(G)_0, \quot{\varphi}{\Gamma} \in \Phi, \quot{\psi}{\Gamma} \in \Psi \right\rbrace.
	\]
	This is easily seen to be a morphism, as for any $f \in \Sf(G)_1$, we have the equalities $\quot{\varphi}{\Gamma} \circ \tau_f^A = \tau_f^{A^{\prime}} \circ F(\overline{f})\quot{\varphi}{\Gamma^{\prime}}$ and $\quot{\psi}{\Gamma} \circ \tau_f^{B} = \tau_f^{B^{\prime}} \circ F(\overline{f})\quot{\psi}{\Gamma^{\prime}}$ so it is immediate that
	\begin{align*}
	\quot{(\Phi \otimes \Psi)}{\Gamma} \circ \tau_f^{A \otimes B} &= \left(\quot{\varphi}{\Gamma} \os{\Gamma}{\otimes} \quot{\psi}{\Gamma}\right) \circ \left(\tau_f^A \os{\Gamma}{\otimes} \tau_f^B\right) \circ \theta_f^{A,B} \\
	&= \left(\tau_f^{A^{\prime}} \os{\Gamma}{\otimes} \tau_f^{B^{\prime}}\right) \circ \left(F(\overline{f})\quot{\varphi}{\Gamma^{\prime}} \os{\Gamma}{\otimes} F(\overline{f})\quot{\psi}{\Gamma^{\prime}}\right) \\
	&= \left(\tau_f^{A^{\prime}} \os{\Gamma}{\otimes} \tau_f^{B^{\prime}}\right) \circ \theta_f^{A^{\prime},B^{\prime}} \circ F(\overline{f})\left(\quot{\varphi}{\Gamma^{\prime}} \os{\Gamma^{\prime}}{\otimes} \quot{\psi}{\Gamma^{\prime}}\right) \\
	&= \tau_{f}^{A^{\prime} \otimes B^{\prime}} \circ F(\overline{f})(\quot{(\Phi \otimes \Psi)}{\Gamma^{\prime}}).
	\end{align*}
	It also follows by the fact that each of the $\os{\Gamma}{\otimes}$ are functors that $\otimes$ preserves compositions and identities, giving that $\otimes$ is a functor as well.
	
	Let us show that the object $I = \lbrace \quot{I}{\Gamma} \; | \; \Gamma \in \Sf(G)_0 \rbrace$ determines an object in $F_G(X)$. Define $T_I := \lbrace \sigma_f:F(\overline{f})\quot{I}{\Gamma^{\prime}} \to \quot{I}{\Gamma} \; | \; f \in \Sf(G)_1 \rbrace$ and note that by assumption for any $\Gamma \xrightarrow{f} \Gamma^{\prime} \xrightarrow{g} \Gamma^{\prime}$ in $\Sf(G)$, we have that
	\[
	\tau_f^{I} \circ F(\overline{f})\tau_g^I = \sigma_f \circ F(\overline{f})\sigma_g = \sigma_{g \circ f} \circ \phi_{f,g}= \tau_{g \circ f}^{I} \circ \phi_{f,g}.
	\]
	Thus $(I,T_I)$ is an object in $F_G(X)$. 
	
	Next we prove that there are natural isomorphisms $\lambda:I \otimes \id_{F_G(X)} \to \id_{F_G(X)}$ and $\rho:\id_{F_G(X)} \otimes I \to \id_{F_G(X)}$;  we will only show that $\lambda$ exists, as the existence of $\rho$ will follow mutatis mutandis.
	
	Fix an object $A \in F_G(X)_0$ and define $\lambda_A: I \otimes A \to A$ by
	\[
	\lambda_A := \left\lbrace \quot{\lambda_{\quot{A}{\Gamma}}}{\Gamma}:\quot{I}{\Gamma} \os{\Gamma}{\otimes} \quot{A}{\Gamma} \xrightarrow{\cong} \quot{A}{\Gamma} \; | \; \Gamma \in \Sf(G)_0 \right\rbrace.
	\]
	It remains to verify that $\lambda_A$ is a morphism in $F_G(X)$, so let $f: \Gamma \to \Gamma^{\prime}$ be a morphism in $\Sf(G)$. Observe that since the tensor functors vary pseudofunctorially over $\Sf(G)$ through the preservation isomorphisms, we get that the invertible $2$-cell
	\[
	\begin{tikzcd}
	F(\overline{f})\big(\quot{I}{\Gamma^{\prime}} \os{\Gamma^{\prime}}{\otimes} \quot{A}{\Gamma^{\prime}}\big) \ar[r, ""{name = U}]{}{\tau_f^I \os{\Gamma}{\otimes} \tau_f^A} \ar[d,swap]{}{F(\overline{f})\quot{\lambda}{\Gamma^{\prime}}} & \quot{A}{\Gamma} \os{\Gamma}{\otimes} \quot{I}{\Gamma} \ar[d]{}{\quot{\lambda_{\quot{A}{\Gamma}}}{\Gamma}} \\
	F(\overline{f})(\quot{A}{\Gamma^{\prime}}) \ar[r, swap, ""{name = L}]{}{\tau_f^A} & \quot{A}{\Gamma} \ar[from = U, to = L, Rightarrow, shorten >= 4pt, shorten <= 4pt]{}{\theta_f}
	\end{tikzcd}
	\]
	produces the commuting diagram:
	\[
	\xymatrix{
		F(\overline{f})\big(\quot{I}{\Gamma^{\prime}} \os{\Gamma^{\prime}}{\otimes} \quot{A}{\Gamma^{\prime}}\big)\ar[d]_{F(\overline{f})\quot{\lambda_{\quot{A}{\Gamma^{\prime}}}}{\Gamma^{\prime}}} \ar[r]^-{\theta_{f}^{I,A}} & F(\overline{f})(\quot{I}{\Gamma^{\prime}}) \os{\Gamma}{\otimes} F(\overline{f})(\quot{A}{\Gamma^{\prime}}) \ar[r]^-{\tau_f^I \os{\Gamma}{\otimes} \tau_f^A} & \quot{A}{\Gamma} \os{\Gamma}{\otimes} \quot{I}{\Gamma} \ar[d]^{\quot{\lambda_{\quot{A}{\Gamma}}}{\Gamma}} \\
		F(\overline{f})(\quot{A}{\Gamma^{\prime}}) \ar[rr]_-{\tau_f^A} & & \quot{A}{\Gamma}
	}
	\]
	This gives that
	\[
	\tau_f^A \circ F(\overline{f})F(\overline{f})\quot{\lambda_{\quot{A}{\Gamma^{\prime}}}}{\Gamma^{\prime}} = \quot{\lambda_{\quot{A}{\Gamma}}}{\Gamma} \circ \left(\tau_f^I \os{\Gamma}{\otimes} \tau_f^A\right) \circ \theta_f^{I,A} = \quot{\lambda_{\quot{A}{\Gamma}}}{\Gamma} \circ \tau_f^{I \otimes A},
	\]
	which proves that $\lambda_A$ is an $F_G(X)$ morphism (and in fact an isomorphism because each $\quot{\lambda}{\Gamma}$ is an isomorphism). That $\lambda$ is natural is immediate from the natruality of each $\quot{\lambda}{\Gamma}$, so we indeed get that $\lambda:I \otimes \id_{F_G(X)} \to \id_{F_G(X)}$ is a natural isomorphism. Similarly, we also have that $\rho:\id_{F_G(X)} \otimes I \to \id_{F_G(X)}$ is a natural isomorphism.
	
	We now need only to verify the existence of the associator(s). For this fix three objects $A, B, C  \in F_G(X)_0$ and define $\alpha_{A,B,C}:(A \otimes B) \otimes C \to A \otimes (B \otimes C)$ via
	\begin{align*}
	&\alpha_{A,B,C} \\
	&:= \left\lbrace \quot{\alpha_{A,B,C}}{\Gamma}:\left(\quot{A}{\Gamma} \os{\Gamma}{\otimes} \quot{B}{\Gamma}\right) \os{\Gamma}{\otimes} \quot{C}{\Gamma} \xrightarrow{\cong} \quot{A}{\Gamma} \os{\Gamma}{\otimes} \left(\quot{B}{\Gamma} \os{\Gamma}{\otimes} \quot{C}{\Gamma}\right) \; \big| \; \Gamma \in \Sf(G)_0\right\rbrace.
	\end{align*}
	We now prove that this is an $F_G(X)$ morphism so fix an $f \in \Sf(G)(\Gamma, \Gamma^{\prime})$ and write $\Gamma = \Dom f$ and $\Gamma^{\prime} = \Codom f$. Then we compute that
	\begin{align*}
	\tau_f^{(A \otimes B) \otimes C} &= \left(\tau_f^{A \otimes B} \os{\Gamma}{\otimes} \tau_f^C\right) \circ \theta_f^{A \otimes B,C} \\
	&= \left(\left(\tau_f^A \os{\Gamma}{\otimes} \tau_f^B\right) \os{\Gamma}{\otimes} \tau_f^C\right)  \circ \left( \theta_f^{A,B}\os{\Gamma}{\otimes}\id_{F(\overline{f})\quot{C}{\Gamma^{\prime}}}\right) \circ \theta_f^{A \otimes B,C}
	\end{align*}
	while on the other hand
	\begin{align*}
	\tau_f^{A \otimes (B \otimes C)} &= \left(\tau_f^A \os{\Gamma}{\otimes} \tau_f^{B \otimes C}\right) \circ \theta_{f}^{A,B \otimes C} \\
	&= \left(\tau_f^A \os{\Gamma}{\otimes} \left(\tau_f^B \os{\Gamma}{\otimes} \tau_f^C\right)\right) \circ \left(\id_{F(\overline{f})\quot{A}{\Gamma^{\prime}}} \os{\Gamma}{\otimes} \theta_{f}^{B,C}\right) \circ \theta_{f}^{A,B \otimes C}.
	\end{align*}
	Furthermore, using the pseudofuncotriality of the tensor functors and the naturality of the associators in each $\Gamma$-fibre category gives that the diagram
	\[
	\begin{tikzcd}
	F(\of)\left(\left(\AGammap \os{\Gamma^{\prime}}{\otimes} \BGammap\right)\os{\Gamma^{\prime}}{\otimes} \quot{C}{\Gamma^{\prime}}\right) \ar[ddddd, bend left = 90]{}{F(\of)\quot{\alpha_{A,B,C}}{\Gamma^{\prime}}} \ar[d]{}{\theta_f^{A \otimes B, C}} \\
	F(\of)\left(\AGammap \os{\Gamma^{\prime}}{\otimes}\BGammap\right)\os{\Gamma}{\otimes} F(\of)(\quot{C}{\Gamma^{\prime}}) \ar[d]{}{\theta_f^{A,B} \os{\Gamma}{\otimes} \id_{F(\of)\quot{C}{\Gamma^{\prime}}}} \\
	\left(F(\of)(\AGammap) \os{\Gamma}{\otimes} F(\of)(\BGammap)\right) \os{\Gamma}{\otimes} F(\of)(\quot{C}{\Gamma^{\prime}}) \ar[d]{}{\quot{\alpha_{F(\of)A,F(\of)B,F(\of)C}}{\Gamma}} \\
	F(\of)(\AGammap) \os{\Gamma}{\otimes}\left( F(\of)(\BGammap) \os{\Gamma}{\otimes} F(\of)(\quot{C}{\Gamma^{\prime}}) \right) \ar[d]{}{\id_{F(\of)\AGammap} \os{\Gamma}{\otimes} \left(\theta_f^{B,C}\right)^{-1}} \\
	F(\of)(\AGammap) \os{\Gamma}{\otimes} F(\of)\left(\BGammap \os{\Gamma^{\prime}}{\otimes} \quot{C}{\Gamma^{\prime}}\right) \ar[d]{}{\left(\theta_f^{A,B \otimes C}\right)^{-1}} \\
	F(\of)\left(\AGammap \os{\Gamma^{\prime}}{\otimes} \left(\BGammap \os{\Gamma^{\prime}}{\otimes} \quot{C}{\Gamma^{\prime}}\right)\right)
	\end{tikzcd}
	\]
	commutes. Substituting the induced identities above yields that the composite morphism $\tau_f^{A \otimes (B \otimes C)} \circ F(\overline{f})\quot{\alpha_{A,B,C}}{\Gamma^{\prime}}$ is equal to 
	\begin{align*}
	& \left(\tau_f^A \os{\Gamma}{\otimes} \left(\tau_f^B \os{\Gamma}{\otimes} \tau_f^C\right)\right) \circ \left(\id_{F(\overline{f})\quot{A}{\Gamma^{\prime}}} \os{\Gamma}{\otimes} \theta_{f}^{B,C}\right) \circ \theta_{f}^{A,B \otimes C} \circ F(\overline{f})\quot{\alpha_{A,B,C}}{\Gamma^{\prime}} \\
	&= \left(\tau_f^A \os{\Gamma}{\otimes} \left(\tau_f^B \os{\Gamma}{\otimes} \tau_f^C\right)\right) \circ \left(\id_{F(\overline{f})\quot{A}{\Gamma^{\prime}}} \os{\Gamma}{\otimes} \theta_{f}^{B,C}\right) \circ \theta_{f}^{A,B \otimes C} \\
	&\circ \left(\theta_f^{A,B\otimes C}\right)^{-1} \circ\left(\id_{F(\overline{f})\quot{A}{\Gamma^{\prime}}} \os{\Gamma}{\otimes} \left(\theta_f^{B,C}\right)^{-1}\right) \\
	&\quad\circ \quot{\alpha_{F(\overline{f})A,F(\overline{f})B,F(\overline{f})C}}{\Gamma} \circ \left(\theta_f^{A,B} \os{\Gamma}{\otimes} \id_{F(\overline{f})\quot{C}{\Gamma^{\prime}}}\right) \circ \theta_f^{A \otimes B, C} \\
	&= \left(\tau_f^A \os{\Gamma}{\otimes} \left(\tau_f^B \os{\Gamma}{\otimes} \tau_f^C\right)\right) \circ \quot{\alpha_{F(\overline{f})A,F(\overline{f})B,F(\overline{f})C}}{\Gamma} \\
	&\circ \left(\theta_f^{A,B} \os{\Gamma}{\otimes} \id_{F(\overline{f})\quot{C}{\Gamma^{\prime}}}\right) \circ \theta_f^{A \otimes B, C} \\
	&=\quot{\alpha_{A,B,C}}{\Gamma} \circ \left(\left(\tau_f^A \os{\Gamma}{\otimes} \tau_f^B\right) \os{\Gamma}{\otimes} \tau_f^C\right) \circ \left(\theta_f^{A,B} \os{\Gamma}{\otimes} \id_{F(\overline{f})\quot{C}{\Gamma^{\prime}}}\right) \circ \theta_f^{A \otimes B, C} \\
	&= \quot{\alpha_{A,B,C}}{\Gamma} \circ \tau_f^{(A \otimes B) \otimes C},
	\end{align*}
	which in turn proves that $\tau_f^{A \otimes (B \otimes C)} \circ F(\overline{f})\quot{\alpha_{A,B,C}}{\Gamma^{\prime}} =\quot{\alpha_{A,B,C}}{\Gamma} \circ \tau_f^{(A \otimes B) \otimes C}.$ Thus $\alpha_{A,B,C}$ is a morphism in $F_G(X)$ (and in fact an isomorphism because each $\quot{\alpha}{\Gamma}$ is an isomorphism). Moreover, the naturality of $\alpha$ is immediate from the fact that each $\quot{\alpha}{\Gamma}$ is a natural isomorphism.
\end{proof}
\begin{proposition}\label{Prop: Section 2: Equivariant cat is symmetric monoidal}
	Let $F:\SfResl_G(X)^{\op} \to \fCat$ be a monoidal pre-equivariant pseudofunctor. If each fibre category $F(\quot{X}{\Gamma})$ is symmetric monoidal, then $F_G(X)$ is a symmetric monoidal category.
\end{proposition}
\begin{proof}
	We know from Theorem \ref{Theorem: Section 2: Monoidal preequivariant pseudofunctor gives monoidal equivariant cat} that $F_G(X)$ is monoidal; we now only need to furnish $F_G(X)$ with the symmetric structure. This is routine, however; for all $\Gamma \in \Sf(G)_0$ let $\quot{s}{\Gamma}:F(\quot{X}{\Gamma}) \times F(\quot{X}{\Gamma}) \to F(\quot{X}{\Gamma}) \times F(\quot{X}{\Gamma})$ be the switching functor (so $\quot{s}{\Gamma}(A,B) = (B,A)$ for objects and $\quot{s}{\Gamma}(\varphi,\psi) = (\psi, \varphi)$ for morphisms) let $\quot{\zeta}{\Gamma}$ be the natural isomorphism fitting into the $2$-cell:
	\[
	\begin{tikzcd}
	F(\quot{X}{\Gamma}) \times F(\quot{X}{\Gamma}) \ar[r, bend left = 30, ""{name = U}]{}{\os{\Gamma}{\otimes}} \ar[r, bend right = 30, ""{name = L}, swap]{}{\os{\Gamma}{\otimes} \circ \quot{s}{\Gamma}} & F(\quot{X}{\Gamma}) \ar[from = U, to = L, Rightarrow, shorten <= 4pt, shorten >= 4pt]{}{\quot{\zeta}{\Gamma}}
	\end{tikzcd}
	\]
	It is then trivial to check that the switching functor $s:F_G(X) \times F_G(X) \to F_G(X) \times F_G(X)$ acts via
	\[
	s(A,B) = \left\lbrace \quot{s}{\Gamma}(\quot{A}{\Gamma},\quot{B}{\Gamma}) \; | \; \Gamma \in \Sf(G)_0 \right\rbrace = (B,A)
	\]
	on objects (with a similar action $s(T_{(A,B)}) = s(T_A,T_B) = (T_B,T_A)$ on transition isomorphisms) and acts via
	\[
	s(\Phi,\Psi) = \left\lbrace \quot{s}{\Gamma}(\quot{\varphi}{\Gamma},\quot{\psi}{\Gamma}) \; | \; \Gamma \in \Sf(G)_0 \right\rbrace = (\Psi,\Phi)
	\]
	on morphisms. Then for any $A, B \in F_G(X)_0$ the assignment
	\[
	\zeta_{A,B} := \left\lbrace \quot{\zeta_{A,B}}{\Gamma}:\quot{A}{\Gamma} \os{\Gamma}{\otimes} \quot{B}{\Gamma} \xrightarrow{\cong} \quot{B}{\Gamma} \os{\Gamma}{\otimes} \quot{A}{\Gamma} \; \big| \; f \in \Sf(G)_1 \right\rbrace
	\]
	produces a morphism in $F_G(X)$ by an argument mutatis mutandis to proving that the left unitor $\lambda$ is an $F_G(X)$-morphism. This in determines an invertible $2$-cell:
	\[
	\begin{tikzcd}
	F_G(X) \times F_G(X) \ar[rr, bend left = 30, ""{name = U}]{}{\otimes} \ar[rr, bend right = 30, swap, ""{name = L}]{}{\otimes \circ s} & & F_G(X) \ar[from = U, to = L, Rightarrow, shorten >= 4pt, shorten <= 4pt]{}{\zeta}
	\end{tikzcd}
	\]
	The natural transformation $\zeta$ is trivially an isomorphism, and also satisfies the required braiding axioms because it satisfies them at each fibre $F(\quot{X}{\Gamma})$. Thus we have that $F_G(X)$ is a symmetric monooidal category.
\end{proof}
\begin{example}\label{Example: Tensor Functor on EDC}
	There are two standard examples of monoidal pre-equivariant pseudofunctors that we should keep in mind. The first is that of the pre-equivariant pseudofunctor $F:\SfResl_G(X)^{\op} \to \fCat$ given on objects by $F(\XGamma) = \DbQl{\XGamma}$ and on morphisms by $F(\of) = \of^{\ast}$. Then each category $F(\XGamma)$ is monoidal with tensor functor
	\[
	\otimes_{\Gamma} := (-) \overset{L}{\otimes}_{\Gamma} (-).
	\]
	Note that the functor $(-)\otimes_{\Gamma}^{L}(-):\DbQl{\XGamma} \times \DbQl{\XGamma} \to \DbQl{\XGamma}$ is the (left) derived tensor functor of $\ell$-adic sheaves on $\XGamma$.
	By Theorem \ref{Theorem: Section 2: Monoidal preequivariant pseudofunctor gives monoidal equivariant cat} and Proposition \ref{Prop: Section 2: Equivariant cat is symmetric monoidal} this produces a symmetric monoidal functor $(-)\otimes^{L}(-)$ on $\DbeqQl{X}$.
		
	The second example to keep in mind is the pre-equivariant pseudofunctor $F:\SfResl_G(X)^{\op} \to \fCat$ given on objects by
	\[
	F(\XGamma) = \Shv(\XGamma, \text{{\'e}t})
	\]
	and on morphisms by
	\[
	F(\of) = \of^{\ast}:\Shv(\XGammap,\text{{\'e}t}) \to \Shv(\XGamma,\text{{\'e}t}).
	\]
	The monoidal functor in this case is defined to be the one induced by
	\[
	\otimes_{\Gamma} := (-) \times_{\Gamma} (-)
	\]
	This then induces an alternative construction of the product functor on the category $\Shv_G(X,\text{{\'e}t})$ by emphasizing its monoidal properties.
\end{example}

We now move to studying some cases when equivariant categories are regular categories and when equivariant categories have subobject classifiers.(cf.\@ Definition \ref{Defn: Subobject Classifier} below). We will defer the study of adjunctions to the next section, as we will need to introduce functors between equivariant categories (which we will call equivariant functors, perhaps perversely), as we will need these functors to formalize the kind of adjoints we need.

\begin{proposition}\label{Prop: Regularity of Equivariant Category}
Let $F$ be a pre-equivariant pseudofunctor on $X$ such that each fibre category $F(\quot{X}{\Gamma})$ is regular and for which each fibre functor $F(\overline{f})$ is finitely complete and preserves regular epimorphisms. Then $F_G(X)$ is regular.
\end{proposition}
\begin{proof}
Recall that a category $\Cscr$ is regular if and only if $\Cscr$ is finitely complete, admits all coequalizers of kernel pairs, and has the property that regular epimorphisms are stable under base change, i.e., regular epimorphisms are stable under pullback. We now proceed to show that $F_G(X)$ is regular by first showing that it is finitely complete; however, for this we simply apply Corollary \ref{Cor: Equivariant cat is finitely complete or cocomplete}.

To see that $F_G(X)$ admits coequalizers of kernel pairs, let $P \in F_G(X)(A,B)$ with $P = \lbrace \quot{\rho}{\Gamma} \; | \; \Gamma \in \Sf(G)_0\rbrace$. Now take the pullback
\[
\xymatrix{
A \times_B A \pullbackcorner \ar[r]^-{\Pi_0} \ar[d]_{\Pi_1} & A \ar[d]^{P} \\
A \ar[r]_{P} & B
}
\]
and note that by Theorem \ref{Thm: Section 2: Equivariant Cat has lims} the pullback takes the form
\[
A \times_B A := \lbrace \quot{A}{\Gamma} \times_{\quot{B}{\Gamma}} \quot{A}{\Gamma} \; | \; \Gamma \in \Sf(G)_0 \rbrace
\]
with the induced transition isomorphisms, while for $k \in \lbrace 0, 1\rbrace$ the morphisms $\Pi_k$ take the form
\[
\Pi_k := \lbrace \quot{\pi_k}{\Gamma}:\quot{A}{\Gamma} \times_{\quot{B}{\Gamma}} \quot{A}{\Gamma} \to \quot{A}{\Gamma} \; | \; \Gamma \in \Sf(G)_0 \rbrace
\]
where each $\quot{\pi_k}{\Gamma}$ is the corresponding first or second pullback map in the fibre category $F(\quot{X}{\Gamma})$. Now define the $\Vscr$-set
\[
C := \lbrace \Coeq(\quot{\pi_0}{\Gamma},\quot{\pi_1}{\Gamma}) \; | \; \Gamma \in \Sf(G)_0 \rbrace.
\]
By assumption, each coequalizer exists in each category $F(\quot{X}{\Gamma})$ and each functor $F(\overline{f})$ preserves these coequalizers. Thus by Theorem \ref{Theorem: Section 2: Equivariant Cat has colims} $(C,T_C)$ is an object in $F_G(X)$ (with induced transition isomorphisms $T_C$) and satisfies $C \cong \Coeq(\Pi_0,\Pi_1)$ with $\coeq(\Pi_0,\Pi_1):A \to C$ given by the set $\coeq(\Pi_0, \Pi_1) = \lbrace \coeq(\quot{\pi_0}{\Gamma}, \quot{\pi_1}{\Gamma}) \; | \; \Gamma \in \Sf(G)_0 \rbrace$. Thus $F_G(X)$ admits coequalizers of kernel pairs.

We now verify that regular epimorphisms are stable under base change. Assume that $\Sigma:A \to C$ is a regular epimorphism in $F_G(X)$ and let $P:B \to C$ be an arbitrary morphism. Take the pullback
\[
\xymatrix{
A \times_C B \ar[r]^-{\Phi} \ar[d]_{Z} \pullbackcorner & A \ar[d]^{\Sigma} \\
B \ar[r]_{P} & C
}
\]
where $Z := \lbrace \quot{\zeta}{\Gamma} \; | \; \Gamma \in \Sf(G)_0\rbrace$. Now let $(\Pi_0, \Pi_1)$ be the kernel pair of $Z$ and let $\Coeq(\Pi_0, \Pi_1)$ be the coequalizer of the kernel pair; note that again by Theorem \ref{Theorem: Section 2: Equivariant Cat has colims}, we have that 
\[
\Coeq(\Pi_0, \Pi_1) = \lbrace \Coeq(\quot{\pi_0}{\Gamma},\quot{\pi_1}{\Gamma}) \; | \; \Gamma \in \Sf(G)_0 \rbrace
\] 
and 
\[
\coeq(\Pi_0, \Pi_1) = \lbrace \coeq(\quot{\pi_0}{\Gamma}, \quot{\pi_1}{\Gamma}) \; | \; \Gamma \in \Sf(G)_0 \rbrace.
\]
Because each category $F(\quot{X}{\Gamma})$ is regular, we have that $\quot{\zeta}{\Gamma} \cong \coeq(\quot{\pi_0}{\Gamma},\quot{\pi_1}{\Gamma})$ via a unique isomorphism $\quot{\theta}{\Gamma}$ fitting into the diagram
\[
\xymatrix{
\quot{A}{\Gamma} \times_{\quot{C}{\Gamma}} \quot{B}{\Gamma} \ar@{=}[d] \ar[rr]^-{\coeq(\quot{\pi_0}{\Gamma},\quot{\pi_1}{\Gamma})} & & \Coeq(\quot{\pi_0}{\Gamma}, \quot{\pi_1}{\Gamma}) \ar@{-->}[d]^{\exists!\quot{\theta}{\Gamma}}_{\cong} \\
\quot{A}{\Gamma} \times_{\quot{C}{\Gamma}} \ar[rr]_-{\quot{\zeta}{\Gamma}} & & \quot{B}{\Gamma}
}
\]
in $F(\quot{X}{\Gamma})$. We need to check now that $\Theta = \lbrace \quot{\theta}{\Gamma} \; | \; \Gamma \in \Sf(G)_0 \rbrace$ is an $F_G(X)$ morphism, as it is automatically an isomorphism if it is an $F_G(X)$-morphism. For this we must prove that for all $f \in \Sf(G)_1$ (say with $\Dom f = \Gamma$ and $\Codom f = \Gamma^{\prime}$),
\[
\tau_f^B \circ F(\overline{f})\quot{\theta}{\Gamma^{\prime}} \overset{?}{=} \quot{\theta}{\Gamma} \circ \tau_f^{\Coeq(\Pi_0, \Pi_1)}.
\]
For this we recall that $Z:A \times_B A \to B$ and that $\coeq(\Pi_0, \Pi_1):A \times_B A \to \Coeq(\Pi_0, \Pi_1)$ are morphisms in $F_G(X)$. We then calculate that on one hand
\begin{align*}
\quot{\zeta}{\Gamma} \circ \tau_f^{A\times_C B} &= \quot{\theta}{\Gamma} \circ \coeq(\quot{\pi_0}{\Gamma}, \quot{\pi_1}{\Gamma}) \circ \tau_f^{A \times_C B} \\
&= \quot{\theta}{\Gamma} \circ \tau_f^{\Coeq(\Pi_0,\Pi_1)} \circ F(\overline{f})\coeq(\quot{\pi_0}{\Gamma^{\prime}}, \quot{\pi_1}{\Gamma^{\prime}})
\end{align*}
while on the other hand
\begin{align*}
\quot{\zeta}{\Gamma} \circ \tau_f^{A\times_C B} &= \tau_f^B \circ F(\overline{f})\quot{\zeta}{\Gamma^{\prime}} = \tau_f^B \circ F(\overline{f})\big(\quot{\theta}{\Gamma^{\prime}} \circ \coeq(\quot{\pi_0}{\Gamma},\quot{\pi_1}{\Gamma})\big) \\
&= \tau_f^B \circ F(\overline{f})\quot{\theta}{\Gamma^{\prime}} \circ F(\overline{f})\coeq(\quot{\pi_1}{\Gamma^{\prime}},\quot{\pi_1}{\Gamma^{\prime}}).
\end{align*}
Because $F(\overline{f})$ preserves regular epimorphisms, it follows that $F(\overline{f})\coeq(\quot{\pi_0}{\Gamma^{\prime}},\quot{\pi_1}{\Gamma^{\prime}})$ is epic. Thus, after canceling $F(\overline{f})\coeq(\quot{\pi_0}{\Gamma^{\prime}},\quot{\pi_1}{\Gamma^{\prime}})$ we get that
\[
\tau_f^B \circ F(\overline{f})\quot{\theta}{\Gamma^{\prime}} = \quot{\theta}{\Gamma} \circ \tau_f^{\Coeq(\Pi_0,\Pi_1)}.
\]
This shows that $\Theta$ is a morphism in $F_G(X)$ and hence proves that $Z$ is a regular epimorphism. This completes the proof that $F_G(X)$ is regular.
\end{proof}

Our next and final proposition for this section lies further in this direction of ``equivariant logic.'' We will show below that under mild assumptions that if all the fibre categories have subobject classifiers and if the fibre functors preserve these subobject classifiers, terminal objects, and subobject classification pullbacks, then $F_G(X)$ has a subobject classifier as well. The conditions on these functors may be seen as saying that the fibre functors (and the resolutions of the variety $X$) make the notion of subobject classification over $X$ equivariant with respect to the group action. In particular, this tells us also that we can descend subobject classification through the group action in this case.
\begin{definition}[{\cite[Chapter 1.3]{MacLaneMoerdijk}}]\label{Defn: Subobject Classifier}
	A category $\Cscr$ has a subobject classifier\index{Subobject Classifier} $\Omega$ if $\Cscr$ has a terminal object $\top$ and there is an object $\Omega \in \Cscr_0$ with a morphism $\true:\top \to \Omega$ such that for any monic $\mu:A \to B$ in $\Cscr$, a unique morphism $\chi_{\mu}:B \to \Omega$ making the diagram
	\[
	\xymatrix{
	A \ar@{-->}[r]^{\exists!} \ar[d]_-{\mu} & \top \ar[d]^{\true} \\
	B \ar[r]_{\exists!\chi_{\mu}} & \Omega
	}
	\]
	into a pullback diagram.
\end{definition}
\begin{example}
Let $\Cscr = \Set$ be the category of sets. Then the set $\Omega = \lbrace 0, 1 \rbrace$ is the subobject classier. The map $\true:\lbrace \ast \rbrace \to \Omega$ is given by $\ast \mapsto 1$ and for any monic $\mu:A \to B$ in $\Set$ the morphism $\chi_{\mu}:B \to \Omega$ is given by
\[
\chi_{\mu}(b) := \begin{cases}
1 & \exists\,a \in A.\, \mu(a) = b; \\
0 & \text{else}.
\end{cases}
\]
\end{example}
\begin{remark}
Let $\Cscr$ and $\Dscr$ be categories with subobject classifiers $\Omega_{\Cscr}$ and $\Omega_{\Dscr}$, respectively. We say that a functor $G:\Cscr \to \Dscr$ {preserves subobject classifing pullbacks}\index{Preserves Subobject Classifying Pullbacks} when it preserves subobject classifiers, terminal objects, and if for a monic $\mu \in \Cscr(A,B)$ which is preserved by $G$ the diagram\index[notation]{true@$\true$}
\[
\xymatrix{
GA \ar[d]_{G\mu} \ar@{-->}[r]^{G!_{A}} \ar[d]_{G\mu} & G\top_{\Cscr} \ar[d]^{G\true_{\Cscr}} \\
GB \ar[r]_{G\chi(\mu)} & G\Omega_{\Cscr}
}
\]
is isomorphic to the pullback diagram
\[
\xymatrix{
GA \ar[d]_{G\mu} \pullbackcorner\ar[r]^{!_{GA}} & \top_{\Dscr} \ar[d]^{\true_{\Dscr}} \\
GB \ar[r]_{\chi(G\mu)} & \Omega_{\Dscr}
}
\]
in $\Dscr$. That is, there are isomorphisms $G\top_{\Cscr} \to \top_{\Dscr}$ and $G\Omega_{\Cscr} \to \Omega_{\Dscr}$ for which the diagram
\[
\begin{tikzcd}
 & GA \ar[rr]{}{!_{GA}} \ar[dd,swap, near end]{}{G\mu} & & \top_{\Dscr} \ar[dd]{}{\true_{\Dscr}} \\
GA \ar[equals, ur] \ar[dd, swap]{}{G\mu} \ar[rr, near end, crossing over]{}{G!_{A}} & & G\top_{\Cscr} \ar[ur]{}{\cong}  \\
 & GB  \ar[rr,swap, near start]{}{\chi(G\mu)} & & \Omega_{\Dscr} \\
GB \ar[ur, equals] \ar[rr,swap]{}{G\chi(\mu)} & & G\Omega_{\Cscr} \ar[ur, swap]{}{\cong} \ar[from = 2-3, to = 4-3, near start, crossing over]{}{G\true_{\Cscr}}
\end{tikzcd}
\]
commutes in $\Dscr$.
\end{remark}
\begin{proposition}\label{Prop: Section 2: Equivariant Cat has Subobject Classifiers}
Let $F$ be a pre-equivariant pseudofunctor such that each fibre category $F(\quot{X}{\Gamma})$ has a subobject classifier, $\quot{\Omega}{\Gamma}$, and for which each fibre functor $F(\overline{f})$ preserves subobject classifiers, terminal objects, and subobject classifying pullbacks. Then $F_G(X)$ admits a subobject classifier $\Omega$ whose object $\Vscr$-set has the form
\[
\Omega = \lbrace \quot{\Omega}{\Gamma} \; | \; \Gamma \in \Sf(G)_0 \rbrace.
\]
\end{proposition}
\begin{proof}
We begin by defining our subobject classifier; since statement of the proposition gave $\Omega = \lbrace \quot{\Omega}{\Gamma} \; | \; \Gamma \in \Sf(G)_0 \rbrace$, we only need to give the transition isomorphisms $T_{\Omega}$. For each $f \in \Sf(G)_1$, say with $\Dom f = \Gamma$ and $\Codom f = \Gamma^{\prime}$, let $\theta_f:F(\overline{f})(\quot{\Omega}{\Gamma^{\prime}}) \xrightarrow{\cong} \quot{\Omega}{\Gamma}$ be the witness isomorphism induced by the fact that each $F(\overline{f})$ preserves subobject classifiers. We then define the $\Vscr$-set of transition isomorphisms by declaring that for all $f \in \Sf(G)_1$, $\tau_f^{\Omega} := \theta_f.$

Let us now verify the cocycle condition on the transition isomorphisms. Fix a composable pair of arrows $\Gamma \xrightarrow{f} \Gamma^{\prime} \xrightarrow{g} \Gamma^{\prime\prime}$ in $\Sf(G)$. We now verify that
\[
\tau_f^{\Omega} \circ F(\overline{f})\tau_g^{\Omega} = \tau_{g \circ f}^{\Omega} \circ \phi_{f,g}.
\]
For this begin by noting that the map $\tau_{g \circ f}^{\Omega}$ is the unique isomorphism between $F(\overline{g} \circ \overline{f})(\quot{\Omega}{\Gamma^{\prime\prime}})$ and $\quot{\Omega}{\Gamma}$; similarly, since $F(\overline{f})\tau_g^{\Omega}$ is the unique isomorphism between $F(\overline{f})\big(F(\overline{g})(\quot{\Omega}{\Gamma^{\prime\prime}})\big)$ and  $F(\overline{f})(\quot{\Omega}{\Gamma^{\prime}})$, while $\tau_f^{\Omega}$ is the unique isomorphism between $F(\overline{f})(\quot{\Omega}{\Gamma^{\prime}})$ and $\quot{\Omega}{\Gamma}$, we will be done if we know that $\phi_{f,g}^{\quot{\Omega}{\Gamma^{\prime\prime}}}$ preserves subobject classifiers. However, this is immediate, as $\phi_{f,g}^{\quot{\Omega}{\Gamma^{\prime\prime}}}$ is a natural isomorphism and both $F(\overline{f})(F(\overline{g})(\quot{\Omega}{\Gamma^{\prime\prime}}))$ and $F(\overline{g} \circ \overline{f})(\quot{\Omega}{\Gamma^{\prime\prime}})$ are subobject classifiers. Thus the diagram
\[
\xymatrix{
F(\overline{f})\big(F(\overline{g})(\quot{\Omega}{\Gamma^{\prime\prime}})\big) \ar[rr]^-{F(\overline{f})\tau_g^{\Omega}} \ar[d]_{\phi_{f,g}^{\quot{\Omega}{\Gamma^{\prime\prime}}}} & & F(\overline{f})(\quot{\Omega}{\Gamma^{\prime}}) \ar[d]^{\tau_f^{\Omega}} \\
F(\overline{g} \circ \overline{f})(\quot{\Omega}{\Gamma^{\prime\prime}}) \ar[rr]_-{\tau_{g \circ f}^{\Omega}} & & \quot{\Omega}{\Gamma}
}
\]
commutes, proving that $(\Omega,T_{\Omega})$ is an object in $F_G(X)$.

We now construct the map $\true:\top \to \Omega$ in $F_G(X)$. Define
\[
\true := \lbrace \quot{\true}{\Gamma}:\quot{\top}{\Gamma} \to \quot{\Omega}{\Gamma} \; | \; \Gamma \in \Sf(G)_0\rbrace;
\]
to prove that this is an $F_G(X)$ morphism, we must verify that for all $f \in \Sf(G)_1$,
\[
\quot{\true}{\Gamma} \circ \tau_{f}^{\top} \overset{?}{=} \tau_{f}^{\Omega} \circ F(\overline{f})\quot{\true}{\Gamma^{\prime}}.
\]
For this, however, note that because each functor preserves object classification pullbacks, given any monic $\mu:C \to D$ in $F_G(X)$ we have that the cube
\[
\begin{tikzcd}
& F(\overline{f})(\quot{C}{\Gamma^{\prime}}) \ar[rr]{}{!_{F(\overline{f})C}} \ar[dd,swap, near end]{}{F(\overline{f})\mu} & & \quot{\top}{\Gamma} \ar[dd]{}{\quot{\true}{\Gamma}} \\
F(\overline{f})(\quot{C}{\Gamma^{\prime}}) \ar[equals, ur] \ar[dd, swap]{}{F(\overline{f})\mu} \ar[rr, near end, crossing over]{}{F(\overline{f})!_{A}} & & F(\overline{f})(\quot{\top}{\Gamma^{\prime}}) \ar[ur]{}{\cong} \\
& F(\overline{f})(\quot{D}{\Gamma^{\prime}})  \ar[rr,swap, near start]{}{\chi(F(\overline{f})\mu)} & & \quot{\Omega}{\Gamma} \\
F(\overline{f})(\quot{D}{\Gamma^{\prime}}) \ar[ur, equals] \ar[rr,swap]{}{F(\overline{f})\chi(\mu)} & & F(\overline{f})(\quot{\Omega}{\Gamma^{\prime}}) \ar[ur, swap]{}{\cong}  \ar[from = 2-3, to = 4-3, near start, crossing over]{}{F(\overline{f})\quot{\true}{\Gamma^{\prime}}}
\end{tikzcd}
\]
commutes in $F(\quot{X}{\Gamma})$. Taking the right-most face shows that the square
\[
\xymatrix{
F(\overline{f})(\quot{\top}{\Gamma^{\prime}}) \ar[d]_{F(\overline{f})\quot{\true}{\Gamma^{\prime}}} \ar[r]^-{\cong} & \quot{\top}{\Gamma} \ar[d]^{\quot{\true}{\Gamma}} \\
F(\overline{f})(\quot{\Omega}{\Gamma^{\prime}}) \ar[r]_-{\cong} & \quot{\Omega}{\Gamma}
}
\]
commutes; because the horizontal isomorphisms in the square above are the unique isomorphisms between the objects, i.e., the witnessing isomorphisms, we have that the diagram takes the form
\[
\xymatrix{
F(\overline{f})(\quot{\top}{\Gamma^{\prime}}) \ar[d]_{F(\overline{f})\quot{\true}{\Gamma^{\prime}}} \ar[r]^-{\tau_f^{\top}} & \quot{\top}{\Gamma} \ar[d]^{\quot{\true}{\Gamma}} \\
F(\overline{f})(\quot{\Omega}{\Gamma^{\prime}}) \ar[r]_-{\tau_f^{\Omega}} & \quot{\Omega}{\Gamma}
}
\]
instead. However, this is exactly the verification that $\true$ is an $F_G(X)$ morphism.

We now construct the classifying map of an arbitrary monic in $F_G(X)$. Begin by letting $M = \lbrace \quot{\mu}{\Gamma} \; | \; \Gamma \in \Sf(G)_0\rbrace$ be a monic $M \in F_G(X)(A,B)$; note that each map $\quot{\mu}{\Gamma}$ is monic in $F(\quot{X}{\Gamma})$ by Theorem \ref{Thm: Section 2: Equivariant Cat has lims}. Then in each fibre category $F(\quot{X}{\Gamma})$ we have that there are (pullback) diagrams
\[
\xymatrix{
\quot{A}{\Gamma} \ar[r] \ar[d]_{\quot{\mu}{\Gamma}} & \quot{\top}{\Gamma} \ar[d]^{\quot{\true}{\Gamma}} \\
\quot{B}{\Gamma} \ar[r]_-{\chi(\quot{\mu}{\Gamma})} & \quot{\Omega}{\Gamma}
}
\]
We thus define $\chi(M):B \to \Omega$ via
\[
\chi(M) := \lbrace \chi(\quot{\mu}{\Gamma}) \; | \; \Gamma \in \Sf(G)_0, \quot{\mu}{\Gamma} \in M \rbrace;
\]
note that we now must show that for all $f \in \Sf(G)_1$, say with $\Dom f = \Gamma$ and $\Codom f = \Gamma^{\prime}$,
\[
\tau_f^{\Omega} \circ F(\overline{f})\chi(\quot{\mu}{\Gamma^{\prime}}) \overset{?}{=} \chi(\quot{\mu}{\Gamma}) \circ \tau_f^{B}.
\]
For this we note that since classifying maps are uniquely determined, it suffices to prove that the diagram
\[
\xymatrix{
\AGamma \ar[d]_{\quot{\mu}{\Gamma}} \ar@{-->}[rrrr]^-{\exists!} & & & & \quot{\top}{\Gamma} \ar[d]^{\quot{\true}{\Gamma}} \\
\BGamma \ar[rrrr]_-{\tau_f^{\Omega} \circ F(\of)\big(\chi(\quot{\mu}{\Gamma^{\prime}})\big) \circ \big(\tau_f^B\big)^{-1}} & & & & \quot{\Omega}{\Gamma}
}
\]
is a pullback diagram. Because the diagram is isomorphic to the standard pullback with $\chi(\quot{\mu}{\Gamma})$ morphism, the pullback routine is a straightforward fact to verify once we know that the diagram commutes. For this we calculate that
\begin{align*}
\tau_f^{\Omega} \circ F(\of)\left(\chi(\quot{\mu}{\Gamma^{\prime}})\right) \circ \left(\tau_f^B\right)^{-1} \circ \quot{\mu}{\Gamma} &= \tau_f^{\Omega} \circ F(\of)\left(\chi(\quot{\mu}{\Gamma^{\prime}})\right) \circ F(\of)\left(\quot{\mu}{\Gamma^{\prime}}\right) \circ \left(\tau_f^A\right)^{-1} \\
&= \tau_f^{\Omega} \circ F(\of)\left(\chi(\quot{\mu}{\Gamma^{\prime}}) \circ \quot{\mu}{\Gamma^{\prime}}\right) \circ \left(\tau_f^A\right)^{-1} \\
&= \tau_f^{\Omega} \circ F(\of)\left(\quot{\true}{\Gamma^{\prime}} \circ !_{\AGammap}\right) \circ \left(\tau_f^A\right)^{-1} \\
&= \tau_f^{\Omega} \circ F(\of)\left(\quot{\true}{\Gamma^{\prime}}\right) \circ F(\of)\left(!_{\AGammap}\right) \circ \left(\tau_f^A\right)^{-1} \\
&=\quot{\true}{\Gamma} \circ \tau_f^{\top} \circ F(\of)\left(!_{\AGammap}\right) \circ \left(\tau_f^A\right)^{-1} \\
&= \quot{\true}{\Gamma} \circ !_{\AGamma} \circ \tau_f^A \circ \left(\tau_f^A\right)^{-1} \\
&= \quot{\true}{\Gamma} \circ !_{\AGamma},
\end{align*}
which proves the diagram commutes. Using the pullback observation we made above it follows that
\[
\chi(\quot{\mu}{\Gamma}) = \tau_f^{\Omega} \circ F(\of)\left(\chi(\quot{\mu}{\Gamma^{\prime}})\right) \circ \left(\tau_f^B\right)^{-1}.
\]
so $\chi(M)$ is indeed an  $F_G(X)$-morphism. Furthermore, we note immediately that the diagram
\[
\xymatrix{
A \ar[d]_{M} \ar[r]^{!_A} & \top \ar[d]^{\true} \\
B \ar[r]_{\chi(M)} & \Omega
}
\]
commutes in $F_G(X)$, as for all $\Gamma \in \Sf(G)_0$ the diagrams
\[
\xymatrix{
\quot{A}{\Gamma} \pullbackcorner \ar[r]^{!_{\quot{A}{\Gamma}}} \ar[d]_{\quot{\mu}{\Gamma}} & \quot{!}{\Gamma} \ar[d]^{\quot{\true}{\Gamma}} \\
\quot{B}{\Gamma} \ar[r]_{\chi(\quot{\mu}{\Gamma})} & \quot{\Omega}{\Gamma}
}
\]
commute. We also note that because all these diagrams are pullbacks, we get that
\[
\xymatrix{
A \pullbackcorner \ar[d]_{M} \ar[r]^{!_A} & \top \ar[d]^{\true} \\
B \ar[r]_{\chi(M)} & \Omega
}
\]
is a pullback in $F_G(X)$. It also follows that the map $\chi(M)$ is unique, as any other map will be equal to $\chi(M)$ as it will be $\Gamma$-wise equal to $\chi(\quot{\mu}{\Gamma})$. This proves that $\Omega$ is a subobject classifier.
\end{proof}

\chapter{Introducing Equivariant Functors and Adjunctions}\label{Chapter 4}

In this section we begin a study of functors that arise between equivariant categories. While it would be of interest to give a comprehensive description of how functors $\Phi:F_G(X) \to E_H(Y)$ behave, in this paper we restrict ourselves primarily to three families of functors between equivariant categories that are induced from changes in equivariant descent data. These three families together form a class of functors that we will call equivariant functors, and the families themselves arise from step-wise change-of-structure or morphism-of-structure type constructions:
\begin{itemize}
	\item Change of Fibre/Pseudofunctor: These are functors of the form $F_G(X) \to K_G(X)$ induced by changes in descent data for some pre-equivariant pseudofunctors $F,K:\SfResl_G(X)^{\op} \to \fCat$;
	\item Change of Space: These are functors of the form $F_G(X) \to F_G(Y)$ or $F_G(Y) \to F_G(X)$, for some $G$-varieties $X$ and $Y$, based on translating equivariant descent data along or against equivariant morphisms $h:X \to Y$ or $h:Y \to X$.	
	\item Change of Group: These are functors of the form $F_G(X) \to F_H(X)$ induced by translating equivariant descent data through morphisms $\varphi:H \to G$, for some smooth algebraic groups $G,H$ over $\Spec K$.
\end{itemize}
\begin{remark}
The imprecision of what is meant by ``change in descent data'' will be cleared up in Sections \ref{Section: Section 3: Change of Fibre}, \ref{Section: Section 3: Change of Space}, and \ref{Section: Section 3: Change of Groups} (as we define explicitly and carefully how to change fibres, spaces, and groups, respectively); however, the precise definition of how to translate descent data in each of these directions are quite technical. We have opted for a lack of precision while giving the big picture of equivariant functors before diving into the technical weeds to study such objects.
\end{remark}
\begin{definition}\index{Equivariant Functor}
An equivariant functor is a functor $\underline{\alpha}:F_G(X) \to K_H(Y)$ which preserves equivariant data induced by any of the three following pieces of information and their compositions:
\begin{itemize}
	\item A pseudonatural transformation $\alpha:F \to K$ (a Change of Fibres; cf.\@ Definition \ref{Defn: Change of Fibre Functor});
	\item An equivariant morphism of varieties $h:X \to Y$ or $h:Y \to X$ (a Change of Spaces; cf.\@ Definition \ref{Defn: Change of Space Functor});
	\item A morphism of smooth algebraic groups $\varphi:H \to G$ (a Change of Groups).
\end{itemize}
\end{definition}

For what follows it will be helpful to have a category of equivariant categories over a fixed group $G$ and $G$-variety $X$. We define this below for later use in this section of the paper.
\begin{definition}
We define $\GEqCat_{/X}$\index[notation]{GEqCat@$\GEqCat_{/X}$} to be the (meta)category of equivariant categories over a fixed algebraic group $G$ and fixed left $G$-variety $X$:
\begin{itemize}
	\item Objects: Categories $F_G(X)$ for a pre-equivariant pseudofunctor $F:\SfResl_G(X)^{\op} \to \fCat$;
	\item Morphisms: Change of fibre equivariant functors $F_G(X) \to K_G(X)$ as induced by changes in descent data (cf.\@ Theorem \ref{Thm: Section 3: Psuedonatural trans lift to equivariant functors} and Definition \ref{Defn: Change of Fibre Functor} for a precise explanation of what this means);
	\item Composition and Identities: As in $\Cat$.
\end{itemize}
Similarly, we define the $2$-category $\fGEqCat_{/X}$
to be the (meta)2-category of $G$-equivariant categories over $X$ by defining:
\begin{itemize}
	\item Objects: $(\fGEqCat_{/X})_0 = (\GEqCat_{/X})_0$;
	\item $1$-Morphisms: $(\fGEqCat_{/X})_1 = (\GEqCat_{/X})_1$;
	\item $2$-Morphisms: Equivariant natural transformations (cf.\@ Lemma \ref{Lemma: Modifications give equivariant natural transformations} for a precise description of what these are).
	\item Compositions: As in $\fCat$.
\end{itemize} 
\end{definition}
The only things that need to be checked from these definitions is that the identity functor on an equivariant category is a Change of Fibre equivariant functor and that the composoition of Change of Firbe functors is again a Change of Fibre functor. However, this is shown in Corollary \ref{Cor: Section 3: Identity functors are change of fibre functors} (immediately following Theorem \ref{Thm: Section 3: Psuedonatural trans lift to equivariant functors}), as that is where it most naturally fits and where the proof of this fact makes the most sense.

\section{Equivariant Functors: Change of Fibre}\label{Section: Section 3: Change of Fibre}

We begin this section by studying perhaps the most straightforward type of equivariant functor: The change of fibre functors. These will allow us to compare different equivariant categories over the same space (such as give embeddings of categories $F_G(X) \to K_G(X)$ whenever we have embeddings $F(\quot{X}{\Gamma}) \to K(\quot{X}{\Gamma})$\ --- this is an immediate corollary of Theorem \ref{Thm: Section 3: Psuedonatural trans lift to equivariant functors}) and lift adjoints between equivariant categories (cf.\@ Theorem \ref{Thm: Section 3: Gamma-wise adjoints lift to equivariant adjoints}).

We will need pseudonatural transformations between pseudofunctors in order to have a unified theory of change of fibre equivariant functors; in fact, a careful reading of the proof of Theorem \ref{Thm: Section 3: Psuedonatural trans lift to equivariant functors} shows that the technical conditions necessary to move between the fibres categories $F(\quot{X}{\Gamma})$ and $K(\quot{X}{\Gamma})$ in a way compatible with the pseudofunctorial behaviour of both $F$ and $K$ exactly give a pseudonatural transformation between the pseudofunctors $F$ and $K$. Below we recall the definition of a pseudonatural transformation (as presented in \cite{Corner2019}); note that unlike the case of pseudofunctors, however, we present a less general version of pseudonatural transformations than what is given in the literature (cf.\@ \cite{Corner2019}, \cite{JohnstoneElephant1}, \cite{Lucas}, \cite{RiehlVertityElementss}). In particular, because we are working with strict $2$-categories (and not, in particular, bicategories) and the domain of the pseudofunctors, $\SfResl_G(X)^{\op}$, has only identity $2$-morphisms, the diagrams describing compatibility of pseudonatural transformations with the $2$-morphisms in $\SfResl_G(X)^{\op}$ and the unitors in $\fCat$ degenerate (for the unitors this is because our pseudofunctors are rigidified). As such we will omit these diagrams in our definition.
\begin{definition}
Let $F,E:\SfResl_G(X)^{\op} \to \fCat$ be pseudofunctors. A pseudonatural transformation\index{Pseudonatural Transformation} $\alpha:F \implies E$ is given by the following data:
\begin{itemize}
	\item For all $\Gamma \in \Sf(G)_0$ a functor $\quot{\alpha}{\Gamma}:F(\quot{X}{\Gamma}) \to E(\quot{X}{\Gamma})$;
	\item For all $f:\Gamma \to \Gamma^{\prime}$ in $\Sf(G)_1$, an invertible natural transformation $\quot{\alpha}{f}$:
	\[
	\begin{tikzcd}
	F(\quot{X}{\Gamma^{\prime}}) \ar[d,swap]{}{\quot{\alpha}{\Gamma^{\prime}}} \ar[r, ""{name = U}]{}{F(\overline{f})} & F(\quot{X}{\Gamma}) \ar[d]{}{\quot{\alpha}{\Gamma}} \ar[from = U, to =L, Rightarrow, shorten <= 5pt, shorten >= 5pt]{}{\quot{\alpha}{f}} \\
	E(\quot{X}{\Gamma^{\prime}}) \ar[r,swap, ""{name = L}]{}{E(\overline{f})} & E(\quot{X}{\Gamma})
	\end{tikzcd}
	\]
	such that for any composable pair of morphisms $\Gamma \xrightarrow{f} \Gamma^{\prime} \xrightarrow{g} \Gamma^{\prime\prime}$ in $\Sf(G)$, the pasting diagram
	\[
	\begin{tikzcd}
	F(\quot{X}{\Gamma^{\prime\prime}}) \ar[r,""{name = UL}]{}{F(\overline{g})} \ar[d,swap]{}{\quot{\alpha}{\Gamma^{\prime\prime}}} & F(\quot{X}{\Gamma^{\prime}}) \ar[d]{}{\quot{\alpha}{\Gamma^{\prime}}} \ar[r,""{name = UR}]{}{F(\overline{f})} & F(\quot{X}{\Gamma}) \ar[d]{}{\quot{\alpha}{\Gamma}} \\
	E(\quot{X}{\Gamma^{\prime\prime}}) \ar[r, swap,""{name = BL}]{}{E(\overline{g})} \ar[rr, swap, bend right = 40, ""{name = BB}]{}{E(\overline{g} \circ \overline{f})} & E(\quot{X}{\Gamma^{\prime}}) \ar[r,swap,""{name = BR}]{}{E(\overline{f})} & E(\quot{X}{\Gamma}) \ar[from = UL, to = BL, swap, Rightarrow, shorten >= 5pt, shorten <= 5pt]{}{\quot{\alpha}{g}} \ar[from = UR, to = BR, Rightarrow, shorten <= 5pt, shorten >= 5pt]{}{\quot{\alpha}{f}} \ar[from = 2-2, to = BB, shorten >= 5pt, shorten <= 5pt, Rightarrow]{}{\quot{\phi_{f,g}}{E}}
	\end{tikzcd}
	\]
	is equal to the pasting diagram:
	\[
	\begin{tikzcd} 
	& F(\quot{X}{\Gamma^{\prime}}) \ar[dr]{}{F(\overline{f})} \\
	F(\quot{X}{\Gamma^{\prime\prime}}) \ar[ur]{}{F(\overline{g})} \ar[rr, ""{name = M}]{}[description]{F(\overline{g} \circ \overline{f})} \ar[d, swap]{}{\quot{\alpha}{\Gamma^{\prime\prime}}} & & F(\quot{X}{\Gamma}) \ar[d]{}{\quot{\alpha}{\Gamma}} \\
	E(\quot{X}{\Gamma^{\prime\prime}}) \ar[rr, swap, ""{name = B}]{}{E(\overline{g}\circ \overline{f})} & & E(\quot{X}{\Gamma}) \ar[from = 1-2, to = M, Rightarrow, shorten <= 5pt, shorten >= 5pt]{}{\quot{\phi_{f,g}}{F}} \ar[from = M, to = B, Rightarrow, shorten >= 5pt, shorten <= 5pt]{}{\quot{\alpha}{g \circ f}}
	\end{tikzcd}
	\]
\end{itemize}
\end{definition}
\begin{remark}
Pseudonatural transformations form the $2$-morphisms in the $3$-category of $2$-categories. While we will not explicitly need $3$-categories in this paper, this remark helps us to realize that there is a (meta) $3$-category $\twoCat$ of $2$-categories and that the hom-object $\twoCat(\SfResl_G(X)^{\op},\fCat)$ is itself a $2$-category whose $0$-cells are pseudofunctors, $1$-cells are pseudonatural transformations, and $2$-cells are modifications. We will in fact use that $\twoCat(\SfResl_G(X)^{\op},\fCat)$ is a $2$-category later in this paper (cf.\@ Lemma \ref{Lemma: Section 3: 2-functors preserve adjoints} and Theorem \ref{Thm: Section 3: Gamma-wise adjoints lift to equivariant adjoints}, for example).
\end{remark}
\begin{example}
Let $F:\SfResl_G(X)^{\op} \to \fCat$ be a pre-equivariant pseudofunctor. The identity pseudonatural transformation $\iota_F:F \Rightarrow F$ is defined by setting first for all $\Gamma \in \Sf(G)_0$ $\quot{\iota_F}{\Gamma}$ to be the functor
\[
\quot{\iota_F}{\Gamma} := \id_{F(\XGamma)}:F(\XGamma) \to F(\XGamma)
\]
and setting, for all $f:\Gamma \to \Gamma^{\prime}$ in $\Sf(G)_1$, $\quot{\iota_F}{f}$ to be the identity natural transformation:
\[
\begin{tikzcd}
F(\XGammap) \ar[r, ""{name = U}]{}{F(\of)} \ar[d, equals, swap]{}{\quot{\iota_F}{\Gamma^{\prime}}} & F(\XGamma) \ar[d, equals]{}{\quot{\iota_F}{\Gamma}} \\
F(\XGammap) \ar[r, swap, ""{name = L}]{}{F(\of)} & F(\XGamma) \ar[from = U, to = L, Rightarrow, shorten >= 4pt, shorten <= 4pt]{}{\quot{\iota_F}{f}} \ar[from = U, to = L, Rightarrow, swap, shorten >= 4pt, shorten <= 4pt]{}{\id}
\end{tikzcd}
\]
Then $\iota_F:F \Rightarrow F$ is the identity pseudonatural transformation on $F$.
\end{example}
\begin{example}\label{Example: Shift functors}
The reason we have chosen the language of pseudonatural transformations is motivated by the following key example: Define the pre-equivariant pseudofunctor $D^b:\SfResl(X)^{\op} \to \fCat$ by $
D^{b}(\Gamma \times X) := D^{b}(\quot{X}{\Gamma})$
on objects and
$
D^{b}(f \times \id_X):= \of^{\ast}$
on morphisms. Then we can define a shift pseudonatural transformation $[1]:D^b \Rightarrow D^b$\index[notation]{Shift@$[1]:D_G^b(X) \to D_G^b(X)$} by setting
\[
\quot{[1]}{\Gamma}:D^b(\quot{X}{\Gamma}) \to D^{b}(\quot{X}{\Gamma})
\]
to be the degree $+1$ shift functor for all $\Gamma \in \Sf(G)_0$ and defining $\quot{\alpha}{f}$ to be the isomorphism witnessing the commutativity
\[
\of^{\ast} \circ \quot{[1]}{\Gamma} \cong \quot{[1]}{\Gamma} \circ \of^{\ast}.
\]
The functor $[1]:D^b_G(X) \to D^b_G(X)$ assembled out of $[1]:D^b \Rightarrow D^b$ (cf.\@ Theorem \ref{Thm: Section 3: Psuedonatural trans lift to equivariant functors} and Proposition \ref{Prop: Section Triangle: Equivariant suspensions}) then gives us our suspension functor on $D^b_G(X)$.
\end{example}
We now define a category that we will use later in this section: The category of pre-equivariant pseudofunctors.
\begin{definition}
Fix a smooth algebraic group $G$ over $\Spec K$ and a left $G$-variety $X$. The category $\PreEq(\SfResl_G(X)^{\op},\fCat)$\index[notation]{PreEq@$\PreEq(\SfResl_G(X)^{\op},\fCat)$} is the category of pre-equivariant pseudofunctors on $X$. It is defined by:
\begin{itemize}
	\item Objects: An object of $\PreEq(\SfResl_G(X)^{\op},\fCat)$ is a pre-equivariant pseudofunctor $F:\SfResl_G(X)^{\op} \to \fCat$;
	\item Morphisms: A morphism of $\PreEq(\SfResl_G(X)^{\op},\fCat)$ is a pseudonatural transformation;
	\item Composition: Vertical composition of pseudonatural transformations in $\twoCat$\index[notation]{TwoCat@$\twoCat$}, the strict $3$-category of $2$-categories;
	\item Identities: As in $\twoCat$ with respect to vertical composition.
\end{itemize}
\end{definition}

Our first use of pseudonatural transformations will be to translate between equivariant categories defined with respect to a fixed variety $X$ and algebraic group $G$. On one hand, we can see the invertible $2$-cells $\quot{\alpha}{f}$ as witnessing isomorphisms for the commutativity of the descent data over the fibre categories induced via the pseudofunctors in the form
\[
\quot{\alpha}{f}:\quot{\alpha}{\Gamma} \circ F(\overline{f}) \xrightarrow{\cong} E(\overline{f}) \circ \quot{\alpha}{\Gamma^{\prime}}. 
\]
Using this we want to reassemble a functor $\underline{\alpha}:F_G(X) \to E_G(X)$ given by simply applying the pseudonatural transformation to morphisms or objects na{\"i}evly, as the pseudonatural relations should be exactly what we need in order to move between the equivariant (descent) categories. While we can always make a straightforward functorial assignment on the object set $A$ of the pair $(A,T_A)$ and on morphisms $P \in F_G(X)_1$ (namely by applying the functors $\quot{\alpha}{\Gamma}$ degreewise to the objects $\AGamma$ and $\quot{\rho}{\Gamma}$, respectively), the assignment of the proposed functor on the transition isomorphisms is more subtle. Because the transition isomorphisms of an object in $E_G(X)$ needs to be of the form 
\[
\tau_f^B:E(\overline{f})\quot{B}{\Gamma^{\prime}} \to \quot{B}{\Gamma},
\] 
when building the transition isomorphisms of the objects $\alpha(\quot{X}{\Gamma})(F(\overline{f})\quot{A}{\Gamma^{\prime}})$, we need to pre-compose by the inverse of the transformation $\quot{\alpha}{f}$ at $\quot{A}{\Gamma^{\prime}}$. Explicitly, we will have to define our transition isomorphisms via the diagram:
\[
\xymatrix{
E(\overline{f})(\alpha(\quot{X}{\Gamma^{\prime}})\quot{A}{\Gamma^{\prime}}) \ar[rr]^-{\left(\quot{\alpha_{\quot{A}{\Gamma^{\prime}}}}{f}\right)^{-1}} \ar[drr] & & \alpha(\quot{X}{\Gamma})(F(\overline{f})\quot{A}{\Gamma^{\prime}}) \ar[d]^{\alpha(\quot{X}{\Gamma})\tau_f^A} \\
 & & \alpha(\quot{X}{\Gamma})\quot{A}{\Gamma}
}
\]
This observation leads us to the construction of functors between equivariant categories induced by pseudonatural transformations, which is proved below. Essentially, each pseudonatural transformation between pre-equivariant pseudofunctor yields a functor between equivariant categories induced by the transformation $\alpha$.

\begin{Theorem}\label{Thm: Section 3: Psuedonatural trans lift to equivariant functors}
Let $F,E:\SfResl_G(X)^{\op} \to \fCat$ be pre-equivariant pseudofunctors and let $\alpha:F \implies E$ be a pseudonatural transformation with fibre functors $\quot{\alpha}{\Gamma}:F(\quot{X}{\Gamma}) \to E(\quot{X}{\Gamma})$ and fibre transformations, for $f \in \Sf(G)_1$ and with $\Dom f = \Gamma, \Codom f = \Gamma^{\prime}$:
\[
\begin{tikzcd}
F(\quot{X}{\Gamma^{\prime}}) \ar[r, bend left = 30, ""{name = U}]{}{\quot{\alpha}{\Gamma} \circ F(\overline{f})} \ar[r, swap, bend right = 30, ""{name = L}]{}{E(\overline{f}) \circ \quot{\alpha}{\Gamma}^{\prime}} & E(\quot{X}{\Gamma}) \ar[from = U, to = L, Rightarrow, shorten >= 5pt, shorten <= 5pt]{}{\quot{\alpha}{f}}
\end{tikzcd}
\]
There then is a functor $\underline{\alpha}:F_G(X) \to E_G(X)$ given by sending morphisms $P = \lbrace \quot{\rho}{\Gamma} \; | \; \Gamma \in \Sf(G)_0\rbrace$ to
\[
\underline{\alpha}P := \lbrace \quot{\alpha}{\Gamma}(\quot{\rho}{\Gamma}) \; | \; \Gamma \in \Sf(G)_0, \quot{\rho}{\Gamma} \in P\rbrace,
\]
sending the set $A = \lbrace \quot{A}{\Gamma} \; | \; \Gamma \in \Sf(G)_0\rbrace$ of an object pair $(A,T_A)$ to
\[
\underline{\alpha}A := \lbrace \quot{\alpha}{\Gamma}(\quot{A}{\Gamma}) \; | \; \Gamma \in \Sf(G)_0, \quot{A}{\Gamma} \in A\rbrace
\]
and finally sending $T_A$ to the set of transition isomorphisms $T_{\underline{\alpha}A} = \lbrace \tau_f^{\underline{\alpha}A}\; | \; f \in \Sf(G)_1 \rbrace$ where each $\tau_f^{\underline{\alpha}A}$ is given by
\[
\tau_f^{\underline{\alpha}A} := \quot{\alpha}{\Gamma}(\tau_f^A) \circ \big(\quot{\alpha_{\quot{A}{\Gamma^{\prime}}}}{f}\big)^{-1}.
\]
\end{Theorem}
\begin{proof}
We first have to show that $(\underline{\alpha}A,T_{\underline{\alpha}A})$ is an object in $E_G(X)$. To show this it suffices to prove that the transition isomorphisms $\tau_f^{\underline{\alpha}A}$ satisfy the cocycle identity on composable arrows $f$ and $g$:
\[
\tau_f^{\underline{\alpha}A} \circ E(\overline{f})\tau_g^{\underline{\alpha}A} = \tau_{g \circ f}^{\underline{\alpha}A} \circ \quot{\phi_{f,g}}{E}.
\]
With this in mind, fix a composable pair of morphisms $\Gamma \xrightarrow{f} \Gamma^{\prime} \xrightarrow{g} \Gamma^{\prime\prime}$ in $\Sf(G)$. We begin by calculating that
\begin{align*}
\tau_{f}^{\underline{\alpha}A} \circ E(\overline{f})\tau_g^{\underline{\alpha}A} &= \quot{\alpha}{\Gamma}(\tau_f^A) \circ \left(\quot{\alpha_{\quot{A}{\Gamma^{\prime}}}}{f}\right)^{-1} \circ E(\overline{f})\big(\quot{\alpha}{\Gamma^{\prime}}(\tau_g^A) \circ \quot{\alpha_{\quot{A}{\Gamma^{\prime\prime}}}}{g}^{-1}\big) \\
&=\quot{\alpha}{\Gamma}(\tau_f^A) \circ \left(\quot{\alpha_{\quot{A}{\Gamma^{\prime}}}}{f}\right)^{-1} \circ E(\overline{f})\big(\quot{\alpha}{\Gamma^{\prime}}(\tau_g^A)\big) \circ E(\overline{f})\big(\quot{\alpha_{\quot{A}{\Gamma^{\prime\prime}}}}{g}\big)^{-1}.
\end{align*}
A routine calculation using the pseudonaturality of $\alpha$, and in particular the invertibility and naturality of the $\quot{\alpha}{f}$, gives that
\[
E(\overline{f})\big(\quot{\alpha}{\Gamma^{\prime}}(\tau_g^A)\big) = \left(\quot{\alpha_{\quot{A}{\Gamma^{\prime}}}}{f}\right) \circ \quot{\alpha}{\Gamma}\big(F(\overline{f})\tau_g^A\big) \circ \left(\quot{\alpha_{F(\overline{g})\quot{A}{\Gamma^{\prime\prime}}}}{f}\right)^{-1}.
\]
Substituting this into the above equation yields
\begin{align*}
&\quot{\alpha}{\Gamma}(\tau_f^A) \circ \left(\quot{\alpha_{\quot{A}{\Gamma^{\prime}}}}{f}\right)^{-1} \circ E(\overline{f})\big(\quot{\alpha}{\Gamma^{\prime}}(\tau_g^A)\big) \circ E(\overline{f})\big(\quot{\alpha_{\quot{A}{\Gamma^{\prime\prime}}}}{g}\big) \\
&= \quot{\alpha}{\Gamma}(\tau_f^A) \circ \left(\quot{\alpha_{\quot{A}{\Gamma^{\prime}}}}{f}\right)^{-1} \circ \left(\quot{\alpha_{\quot{A}{\Gamma^{\prime}}}}{f}\right) \circ \quot{\alpha}{\Gamma}\big(F(\overline{f})\tau_g^A\big) \\
&\circ \left(\quot{\alpha_{F(\overline{g})\quot{A}{\Gamma^{\prime\prime}}}}{f}\right)^{-1} \circ E(\overline{f})\big(\quot{\alpha_{\quot{A}{\Gamma^{\prime\prime}}}}{g}\big)^{-1} \\
&= \quot{\alpha}{\Gamma}(\tau_f^A) \circ \quot{\alpha}{\Gamma}(F(\overline{f})\tau_g^{A}) \circ \left(\quot{\alpha_{F(\overline{g})\quot{A}{\Gamma^{\prime\prime}}}}{f}\right)^{-1} \circ E(\overline{f})\big(\quot{\alpha_{\quot{A}{\Gamma^{\prime\prime}}}}{g}\big)^{-1} \\
&= \quot{\alpha}{\Gamma}\big(\tau_f^A \circ F(\overline{f})\tau_g^A\big) \circ \left(\quot{\alpha_{F(\overline{g})\quot{A}{\Gamma^{\prime\prime}}}}{f}\right)^{-1} \circ E(\overline{f})\big(\quot{\alpha_{\quot{A}{\Gamma^{\prime\prime}}}}{g}\big)^{-1} \\
&= \quot{\alpha}{\Gamma}(\tau_{g \circ f}^A \circ \quot{\phi_{f,g}}{F})\circ \left(\quot{\alpha_{F(\overline{g})\quot{A}{\Gamma^{\prime\prime}}}}{f}\right)^{-1} \circ E(\overline{f})\big(\quot{\alpha_{\quot{A}{\Gamma^{\prime\prime}}}}{g}\big)^{-1} \\
&= \quot{\alpha}{\Gamma}(\tau_{g \circ f}^{A}) \circ \quot{\alpha}{\Gamma}(\quot{\phi_{f,g}}{F})\circ \left(\quot{\alpha_{F(\overline{g})\quot{A}{\Gamma^{\prime\prime}}}}{f}\right)^{-1} \circ E(\overline{f})\big(\quot{\alpha_{\quot{A}{\Gamma^{\prime\prime}}}}{g}\big)^{-1}.
\end{align*}
From here we once again use the pseudonaturality of $\alpha$ to derive that the diagram
\[
\begin{tikzcd}
E(\overline{f})\big(E(\overline{g})\big(\quot{\alpha}{\Gamma^{\prime\prime}}(\quot{A}{\Gamma^{\prime\prime}})\big)\big) \ar[rr]{}{E(\overline{f})\big(\quot{\alpha_{\quot{A}{\Gamma^{\prime\prime}}}}{g}\big)^{-1}} \ar[dd, swap]{}{\quot{\phi_{f,g}}{E}} & & E(\overline{f})\big(\quot{\alpha}{\Gamma^{\prime}}(F(\overline{g})(\quot{A}{\Gamma^{\prime\prime}}))\big) \ar[d]{}{\big(\quot{\alpha_{F(\overline{g})\quot{A}{\Gamma^{\prime\prime}}}}{f}\big)^{-1}} \\
 & & \ar[d]{}{\quot{\alpha}{\Gamma}(\quot{\phi_{f,g}}{\Gamma})} \quot{\alpha}{\Gamma}\big(F(\overline{f})(F(\overline{g})(\quot{A}{\Gamma^{\prime\prime}}))\big) \\
E(\overline{g} \circ \overline{f})\big(\quot{\alpha}{\Gamma^{\prime\prime}}(\quot{A}{\Gamma^{\prime\prime}})\big)\ar[rr, swap]{}{\quot{\alpha_{\quot{A}{\Gamma^{\prime\prime}}}}{g \circ f}^{-1}} & & \quot{\alpha}{\Gamma}\big(F(\overline{g} \circ \overline{f})(\quot{A}{\Gamma^{\prime\prime}})\big)
\end{tikzcd}
\]
commutes. Using the implied identity
\[
\left(\quot{\alpha_{\quot{A}{\Gamma^{\prime\prime}}}}{g \circ f}\right)^{-1} \circ \quot{\phi_{f,g}}{E} = \quot{\alpha}{\Gamma}(\quot{\phi_{f,g}}{F})\circ \left(\quot{\alpha_{F(\overline{g})\quot{A}{\Gamma^{\prime\prime}}}}{f}\right)^{-1} \circ E(\overline{f})\big(\quot{\alpha_{\quot{A}{\Gamma^{\prime\prime}}}}{g}\big)^{-1}
\]
in the equation above then gives that
\begin{align*}
&\quot{\alpha}{\Gamma}(\tau_{g \circ f}^{A}) \circ \quot{\alpha}{\Gamma}(\quot{\phi_{f,g}}{F})\circ \left(\quot{\alpha_{F(\overline{g})\quot{A}{\Gamma^{\prime\prime}}}}{f}\right)^{-1} \circ E(\overline{f})\big(\quot{\alpha_{\quot{A}{\Gamma^{\prime\prime}}}}{g}\big)^{-1} \\&= \quot{\alpha}{\Gamma}(\tau_{g \circ f}^A) \circ \big(\quot{\alpha_{\quot{A}{\Gamma^{\prime\prime}}}}{g \circ f}\big)^{-1} \circ \quot{\phi_{f,g}}{E} = \tau_{g \circ f}^{\underline{\alpha}A} \circ \quot{\phi_{f,g}}{E}.
\end{align*}
Therefore we have that
\[
\tau_f^{\underline{\alpha}A} \circ E(\overline{f})\tau_g^{\underline{\alpha}A} = \tau_{g \circ f}^{\underline{\alpha}A} \circ \quot{\phi_{f,g}}{E},
\]
proving that $T_{\underline{\alpha}A}$ satisfies the cocycle condition and hence that $(\underline{\alpha}A,T_{\underline{\alpha}A})$ is an object in $E_G(X)$.

We now must prove that $\underline{\alpha}P$ is a morphism whenever $P \in F_G(X)(A,B)$; note that if we do this we are done, as from the definition 
\[
\underline{\alpha}P = \lbrace \quot{\alpha}{\Gamma}(\quot{\rho}{\Gamma}) \; | \; \Gamma \in \Sf(G)_0, \quot{\rho}{\Gamma} \rbrace
\]
it follows from the fact that the $\quot{\alpha}{\Gamma}$ are all functors that $\underline{\alpha}$ preserves compositions and identity morphisms. Thus we need only prove that for all $f \in \Sf(G)_1$ the identity
\[
\tau_f^{\underline{\alpha}B} \circ E(\overline{f})\quot{\underline{\alpha}\rho}{\Gamma^{\prime}} = \quot{\underline{\alpha}\rho}{\Gamma} \circ \tau_f^{\underline{\alpha}A}
\]
holds where $\Dom f = \Gamma$ and $\Codom f = \Gamma^{\prime}$. For this we first note that
\[
\quot{(\underline{\alpha}P)}{\Gamma^{\prime}} = \quot{\alpha}{\Gamma^{\prime}}(\quot{\rho}{\Gamma^{\prime}})
\]
and similarly $\quot{(\underline{\alpha}P)}{\Gamma} = \quot{\alpha}{\Gamma^{\prime}}(\quot{\rho}{\Gamma})$. We now, as in the proof of the cocycle identity, derive that the diagram
\[
\xymatrix{
E(\overline{f})\big(\quot{\alpha}{\Gamma^{\prime}}(\quot{A}{\Gamma^{\prime}})\big) \ar[d]_{\left(\quot{\alpha_{\quot{A}{\Gamma^{\prime}}}}{f}\right)^{-1}} \ar[rr]^-{E(\overline{f})\big(\quot{\alpha}{\Gamma^{\prime}}(\quot{\rho}{\Gamma^{\prime}})\big)} & & E(\overline{f})\big(\quot{\alpha}{\Gamma^{\prime}}(\quot{B}{\Gamma^{\prime}})\big) \\
\quot{\alpha}{\Gamma}\big(F(\overline{f})(\quot{A}{\Gamma^{\prime}})\big) \ar[rr]_-{\quot{\alpha}{\Gamma}\big(F(\overline{f})\quot{\rho}{\Gamma^{\prime}}\big)} & & \quot{\alpha}{\Gamma}\big(F(\overline{f})(\quot{B}{\Gamma^{\prime}})\big) \ar[u]_{\quot{\alpha_{\quot{B}{\Gamma^{\prime}}}}{f}}
}
\]
commutes in $E(\quot{X}{\Gamma})$. This allows us to compute that
\begin{align*}
&\tau_f^{\underline{\alpha}B} \circ E(\overline{f})\big(\quot{\alpha}{\Gamma^{\prime}}(\quot{\rho}{\Gamma^{\prime}})\big)\\
 &= \tau_f^{\underline{\alpha}B} \circ \quot{\alpha_{\quot{B}{\Gamma^{\prime}}}}{f} \circ \quot{\alpha}{\Gamma}\big(F(\overline{f})\quot{\rho}{\Gamma^{\prime}}\big)  \circ \left(\quot{\alpha_{\quot{A}{\Gamma^{\prime}}}}{f}\right)^{-1} \\
&= \quot{\alpha}{\Gamma}(\tau_f^B) \circ \left(\quot{\alpha_{\quot{B}{\Gamma^{\prime}}}}{f}\right)^{-1} \circ \quot{\alpha_{\quot{B}{\Gamma^{\prime}}}}{f} \circ \quot{\alpha}{\Gamma}\big(F(\overline{f})\quot{\rho}{\Gamma^{\prime}}\big)  \circ \left(\quot{\alpha_{\quot{A}{\Gamma^{\prime}}}}{f}\right)^{-1} \\
&= \quot{\alpha}{\Gamma}(\tau_f^B) \circ \quot{\alpha}{\Gamma}\big(F(\overline{f})\quot{\rho}{\Gamma^{\prime}}\big) \circ \left(\quot{\alpha_{\quot{A}{\Gamma^{\prime}}}}{f}\right)^{-1} \\
&= \quot{\alpha}{\Gamma}\big(\tau_f^B \circ F(\overline{f})\quot{\rho}{\Gamma^{\prime}}\big) \circ \left(\quot{\alpha_{\quot{A}{\Gamma^{\prime}}}}{f}\right)^{-1} \\
&= \quot{\alpha}{\Gamma}(\quot{\rho}{\Gamma} \circ \tau_f^A) \circ \left(\quot{\alpha_{\quot{A}{\Gamma^{\prime}}}}{f}\right)^{-1} = \quot{\alpha}{\Gamma}(\quot{\rho}{\Gamma}) \circ \quot{\alpha}{\Gamma}(\tau_f^A) \circ \left(\quot{\alpha_{\quot{A}{\Gamma^{\prime}}}}{f}\right)^{-1} \\
&= \quot{\alpha}{\Gamma}(\quot{\rho}{\Gamma}) \circ \tau_{f}^{\underline{\alpha}A}.
\end{align*}
Thus $\underline{\alpha}P$ is an $F_G(X)$-morphism.
\end{proof}
\begin{corollary}\label{Cor: Section 3: Identity functors are change of fibre functors}
The identity pseudonatural transformation $\iota_F:F \Rightarrow F$ gives rise to the identity functor $\id_{F_G(X)}:F_G(X) \to F_G(X)$. In particular $\id_{F_G(X)}$ is a change of fibre equivariant functor with $\id_{F_G(X)} = \underline{\iota_F}.$
\end{corollary}
\begin{proof}
Begin by noting that the pseudonatural transformation $\iota_F:F \Rightarrow F$ is defined by setting
\[
\quot{\iota_F}{\Gamma} := \id_{F(\quot{X}{\Gamma})}:F(\quot{X}{\Gamma}) \to F(\quot{X}{\Gamma})
\]
for all objects $\Gamma \in \Sf(G)_0$ and $\quot{\iota_F}{f} := \iota_{F(\overline{f})}$ for all $f \in \Sf(G)_1$, i.e., by setting $\quot{\iota_F}{f}$ to be the identity $2$-cell:
\[
\begin{tikzcd}
F(\quot{X}{\Gamma^{\prime}}) \ar[rr, bend left = 30, ""{name = U}]{}{F(\overline{f})} \ar[rr, bend right = 30, swap, ""{name = L}]{}{F(\overline{f})}  & & F(\quot{X}{\Gamma}) \ar[from = U, to = L, equals, shorten >= 5pt, shorten <= 5pt]{}{\iota_{F(\overline{f})}}
\end{tikzcd}
\]
It is then straightforward to check that $\iota_F$ is a pseudonatural transformation and that $\underline{\iota_F} = \id_{F_G(X)}$.
\end{proof}
\begin{corollary}
There is a functor \[\PreEq(\SfResl_G(X)^{\op},\fCat) \to \GEqCat_{/X}\] which sends a pre-equivariant pseudofunctor $F$ to the category $F_G(X)$ and a pseudonatural transformation $\alpha$ to the functor $\underline{\alpha}$.
\end{corollary}
\begin{proof}
Definition \ref{Defn: Equivariant Cats} and Theorem \ref{Thm: Section 3: Psuedonatural trans lift to equivariant functors} give the functions of the functor while Corollary \ref{Cor: Section 3: Identity functors are change of fibre functors} shows that this functor preserves identity morphisms. All that remains is to show that this functor preserves composition; however, this is a routine and tedious verification which is omitted here.
\end{proof}
\begin{definition}\label{Defn: Change of Fibre Functor}\index{Equivariant Functor! Change of Fibre Functor}
A Change of Fibre functor is a functor $R:F_G(X) \to E_G(X)$ where $R = \underline{\alpha}$ for some pseudonatural transformation $\alpha:F \Rightarrow E$.
\end{definition}

We now record an immediate result of the proof of Theorem \ref{Thm: Section 3: Psuedonatural trans lift to equivariant functors} above: The functor from the category of pre-equivariant pseudofunctors on $X$ to the category of $G$-equivariant categories on $X$ is faithful.
%

We now move on to prove the first key technical result about certain cases when change of fibre equivariant functors preserve limits of fixed shape. These will let us prove that under natural conditions (and in particular conditions we will be considering in this paper) the functors induced by pseudofunctors either preserve all limits of a specific shape (cf.\@ Corollary \ref{Cor: Section 3: (Co)limit preserving pseudonaturals give (co)limit preserving functors}) or simply preserve all limits or colimits (cf.\@ Corollary \ref{Cor: Section 3: (Co)Continuous pseudonaturals give (co)continuous functors}). These in turn are used to establish when certain equivariant functors are additive (cf.\@ Corollary \ref{Cor: Section 3: Additive equivraint change of fibre functors}), which may be used to begin establishing the properties of the equivariant categories we care about.

\begin{proposition}\label{Prop: Section 3: Change of fibre equivariant functor preserves limit of a specific shape}
Let $F$ and $E$ be pre-equivariant pseudofunctors on $X$ such that for a fixed index category $I$ the hypotheses of Theorem \ref{Thm: Section 2: Equivariant Cat has lims} are satisfied for $F$ and $E$ for diagram functors $d_F:I \to F_G(X), d_E:I \to E_G(X)$, and for all $\Gamma \in \Sf(G)_0$
\[
\xymatrix{
I \ar[drr]_{\quot{d_F}{\Gamma}} \ar[r]^-{d_F} & F_G(X) \ar[r]^-{\hat{\i}_F} & \prod\limits_{\Gamma \in \Sf(G)_0} F(\XGamma) \ar[d]^{\pi_{\Gamma}^{F}} \\
 & & F(\XGamma)
}
\] 
and:
\[
\xymatrix{
I \ar[drr]_{\quot{d_E}{\Gamma}} \ar[r]^-{d_E} & E_G(X) \ar[r]^-{\hat{\i}_E} & \prod\limits_{\Gamma \in \Sf(G)_0} E(\XGamma) \ar[d]^{\pi_{\Gamma}^{E}} \\
& & E(\XGamma)
}
\]
Assume furthermore that $\alpha:F \Rightarrow E$ is a pseudonatural transformation such that each functor $\quot{\alpha}{\Gamma}$ preserves the limit over each diagram $\quot{d_F}{\Gamma}$, i.e., the limit of $\quot{\alpha}{\Gamma} \circ \quot{d_F}{\Gamma}$ exists in $E(\quot{X}{\Gamma})$ and is isomorphic to $\quot{\alpha}{\Gamma}$ applied to the limit of $\quot{d_F}{\Gamma}$. Then $\underline{\alpha}$ preserves the limit of the induced diagram $d_F$ in $E_G(X)$.
\end{proposition}
\begin{proof}
Appealing to Theorem \ref{Thm: Section 2: Equivariant Cat has lims} gives that $\lim(d_F)$ exists in $F_G(X)$ and is given by the $\Gamma$-wise limit of the diagrams $\quot{d_F}{\Gamma}$; appealing to Theorem \ref{Thm: Section 2: Equivariant Cat has lims} again gives that the limit of $\underline{\alpha} \circ d_F$ exists in $E_G(X)$ and is given $\Gamma$-wise as
\[
\lim(\underline{\alpha} \circ d_F) = \lbrace \lim(\quot{\alpha}{\Gamma} \circ \quot{d_F}{\Gamma}) \; | \; \Gamma \in \Sf(G)_0 \rbrace.
\]
As such, to prove that $\underline{\alpha}$ preserves the limit of $d_F$, we only need to verify that there is an isomorphism
\[
\underline{\alpha}\left(\lim_{\substack{\longleftarrow \\ i \in I_0}} d_F(i)\right) \cong \lim_{\substack{\longleftarrow \\ i \in I_0}}\left(\underline{\alpha} \circ d_F\right)(i)
\]
in $E_G(X)$.

To prove the existence of the isomorphism claimed above, for all $i \in I_0$ let
\[
A_i := d_F(i) = \lbrace \quot{A_i}{\Gamma} \; | \; \Gamma \in \Sf(G)_0\rbrace = \lbrace \quot{d_F(i)}{\Gamma} \; | \; \Gamma \in \Sf(G)_0\rbrace
\]
and let
\[
A := \lim_{\substack{\longleftarrow \\ i \in I_0}} A_i = \lim_{\substack{\longleftarrow \\ i \in I_0}} d_F(i) = \left\lbrace \lim_{\substack{\longleftarrow \\ i \in I_0}} \quot{A_i}{\Gamma} \; | \; \Gamma \in \Sf(G)_0\right\rbrace.
\]
To construct the isomorphism
\[
\underline{\alpha}\left(\lim_{\substack{\longleftarrow \\ i \in I_0}} d_F(i)\right) \cong \lim_{\substack{\longleftarrow \\ i \in I_0}}\left(\underline{\alpha} \circ d_F\right)(i)
\]
we begin by defining the witness isomorphisms induced by the limit preservation of the $\quot{\alpha}{\Gamma}$ be given by
\[
\quot{\theta_{\quot{A_i}{\Gamma}}}{\Gamma}:\quot{\alpha}{\Gamma}\left(\lim_{\substack{\longleftarrow \\ i \in I_0}}\quot{A_i}{\Gamma}\right) \to \lim_{\substack{\longleftarrow \\ i \in I_0}}\quot{\alpha}{\Gamma}\left(\quot{A_i}{\Gamma}\right).
\]
and then define the $\Vscr$-set
\[
\Theta :=\lbrace \quot{\theta_{\quot{A_i}{\Gamma}}}{\Gamma} \; | \; \Gamma \in \Sf(G)_0\rbrace.
\]
Because each of the $\quot{\theta}{\Gamma}$ and $\tau_f^{-}$ maps are isomorphisms, we only need to verify that $\Theta$ is a $E_G(X)$ morphism to prove the proposition, i.e., we need to check that the diagram
\[
\xymatrix{
E(\overline{f})\big(\quot{\alpha}{\Gamma^{\prime}}(\lim(\quot{A_i}{\Gamma^{\prime}}))\big) \ar[d]_{\tau_f^{\underline{\alpha}A}} \ar[rr]^-{E(\overline{f})\quot{\theta_{\quot{A_i}{\Gamma^{\prime}}}}{\Gamma^{\prime}}} & & E(\overline{f})\big(\lim (\quot{\alpha}{\Gamma^{\prime}}(\quot{A_i}{\Gamma^{\prime}}))\big) \ar[d]^{\tau_f^{\lim \underline{\alpha}A_i}}  \\
\quot{\alpha}{\Gamma}\big(\lim(\quot{A_i}{\Gamma})\big) \ar[rr]_-{\quot{\theta_{\quot{A_i}{\Gamma}}}{\Gamma}} & & \lim(\quot{\alpha}{\Gamma}(\quot{A_i}{\Gamma}))
}
\]
commutes.

Recall first that $\quot{\theta_f}{E}$ is the witness isomorphism
\[
\quot{\theta_f}{E}:E(\overline{f})\big(\lim(\quot{\alpha}{\Gamma^{\prime}} \circ \quot{d_F}{\Gamma^{\prime}})\big) \xrightarrow{\cong} \lim\big(E(\overline{f}) \circ \quot{\alpha}{\Gamma^{\prime}} \circ \quot{d_F}{\Gamma^{\prime}}\big),
\]
and similarly for $\quot{\theta_f}{F}:F(\overline{f})\big(\lim (\quot{d_F}{\Gamma})\big) \xrightarrow{\cong} \lim\big(F(\overline{f})\circ \quot{d_F}{\Gamma^{\prime}}\big)$. By Theorem \ref{Thm: Section 3: Psuedonatural trans lift to equivariant functors} we have that
\[
\quot{\theta_{\quot{A_i}{\Gamma^{\prime}}}}{\Gamma} \circ \tau_f^{\underline{\alpha}A} = \quot{\theta_{\quot{A_i}{\Gamma^{\prime}}}}{\Gamma} \circ \quot{\alpha}{\Gamma}(\tau_f^{A}) \circ \left(\quot{\alpha_{\quot{A}{\Gamma^{\prime}}}}{f}\right)^{-1}.
\]

Note that by Theorem \ref{Thm: Section 2: Equivariant Cat has lims}, we have
\[
\tau_{f}^{\lim \overline{\alpha}A_i} = \lim\big(\tau_f^{\underline{\alpha}A_i}\big) \circ \quot{\theta_f}{E}
\]
so it follows that
\begin{align*}
\tau_{f}^{\lim \underline{\alpha}A} \circ E(\overline{f})\quot{\theta_{\quot{A_i}{\Gamma^{\prime}}}}{\Gamma^{\prime}} &= \lim\big(\quot{\alpha}{\Gamma}(\tau_f^{A_i})\big) \circ \quot{\theta_f}{E} \circ E(\overline{f})\quot{\theta_{\quot{A_i}{\Gamma^{\prime}}}}{\Gamma^{\prime}}.
\end{align*}
Now observe that from the uniqueness of morphisms (isomorphisms in this case) induced by the universal property of a limit, we derive that the diagram
\[
\xymatrix{
E(\overline{f})\big(\quot{\alpha}{\Gamma^{\prime}}(\lim (\quot{A_i}{\Gamma^{\prime}}))\big) \ar[d]_{E(\overline{f})\quot{\theta_{\quot{A_i}{\Gamma^{\prime}}}}{\Gamma^{\prime}}}  \ar[rr]^-{\left(\quot{\alpha_{\quot{A}{\Gamma^{\prime}}}}{f}\right)^{-1}} & & \quot{\alpha}{\Gamma}\big(F(\overline{f})(\lim(\quot{A_i}{\Gamma}))\big) \ar[d]^{\quot{\alpha}{\Gamma}\big(\quot{\theta_f}{F}\big)} \\
E(\overline{f})\big(\lim(\quot{\alpha}{\Gamma^{\prime}}(\quot{A_i}{\Gamma^{\prime}}))\big)  & & \quot{\alpha}{\Gamma}\big(\lim (F(\overline{f})(\quot{A_i}{\Gamma^{\prime}}))\big) \ar[d]^{\quot{\theta_{F(\overline{f})\quot{A_i}{\Gamma^{\prime}}}}{\Gamma}}\\
\lim E(\overline{f})\big(\quot{\alpha}{\Gamma^{\prime}}(\quot{A_i}{\Gamma^{\prime}})\big)     \ar[u]^{\left(\quot{\theta_f}{E}\right)^{-1}} & & \lim \left(\quot{\alpha}{\Gamma}\big(F(\overline{f})(\quot{A_i}{\Gamma^{\prime}})\big)\right) \ar[ll]^-{\lim\left(\quot{\alpha_{\quot{A_i}{\Gamma^{\prime}}}}{f}\right)}
}
\]
commutes in $E(\quot{X}{\Gamma})$. Substituting the induced identity into the equation above gives that
\begin{align*}
&\tau_{f}^{\lim \underline{\alpha}A} \circ E(\overline{f})\quot{\theta_{\quot{A_i}{\Gamma^{\prime}}}}{\Gamma^{\prime}} \\
&= \lim\big(\tau_f^{\underline{\alpha}A_i}\big) \circ \quot{\theta_f}{E} \circ E(\overline{f})\quot{\theta_{\quot{A}{\Gamma^{\prime}}}}{\Gamma^{\prime}} \\
&=\lim(\tau_f^{\underline{\alpha}A_i}) \circ  \quot{\theta_f}{E} \circ \left(\quot{\theta_f}{E}\right)^{-1} \circ \lim\left(\quot{\alpha_{\quot{A_i}{\Gamma^{\prime}}}}{f} \right) \circ \quot{\theta_{F(\overline{f})\quot{A_i}{\Gamma^{\prime}}}}{\Gamma} \\
&\circ \quot{\alpha}{\Gamma}\big(\quot{\theta_f}{F}\big) \circ \left(\quot{\alpha_{\quot{A}{\Gamma^{\prime}}}}{f}\right)^{-1} \\
&= \lim(\tau_f^{\underline{\alpha}A_i}) \circ \lim\left(\quot{\alpha_{\quot{A_i}{\Gamma^{\prime}}}}{f} \right) \circ \quot{\theta_{F(\overline{f})\quot{A_i}{\Gamma^{\prime}}}}{\Gamma} \circ \quot{\alpha}{\Gamma}\big(\quot{\theta_f}{F}\big) \circ \left(\quot{\alpha_{\quot{A}{\Gamma^{\prime}}}}{f}\right)^{-1} \\
&= \lim\left(\tau_f^{\underline{\alpha}A_i} \circ \quot{\alpha_{\quot{A_i}{\Gamma^{\prime}}}}{f} \right) \circ \quot{\theta_{F(\overline{f})\quot{A_i}{\Gamma^{\prime}}}}{\Gamma} \circ \quot{\alpha}{\Gamma}\big(\quot{\theta_f}{F}\big) \circ \left(\quot{\alpha_{\quot{A}{\Gamma^{\prime}}}}{f}\right)^{-1} \\
&= \lim\left(\quot{\alpha}{\Gamma}(\tau_f^{A_i}) \circ \quot{\alpha_{\quot{A_i}{\Gamma^{\prime}}}}{f} \circ \left(\alpha_{\quot{A_i}{\Gamma^{\prime}}}{f}\right)^{-1}\right) \circ \quot{\theta_{F(\overline{f})\quot{A_i}{\Gamma^{\prime}}}}{\Gamma} \\
&\circ \quot{\alpha}{\Gamma}\big(\quot{\theta_f}{F}\big) \circ \left(\quot{\alpha_{\quot{A}{\Gamma^{\prime}}}}{f}\right)^{-1} \\
&= \lim\left(\quot{\alpha}{\Gamma}(\tau_f^{A_i})\right) \circ \quot{\theta_{F(\overline{f})\quot{A_i}{\Gamma^{\prime}}}}{\Gamma} \circ \quot{\alpha}{\Gamma}\big(\quot{\theta_f}{F}\big) \circ \left(\quot{\alpha_{\quot{A}{\Gamma^{\prime}}}}{f}\right)^{-1}.
\end{align*}
We now use the uniqueness of limit cones to give the commutativity of the diagram
\[
\xymatrix{
\quot{\alpha}{\Gamma}\big(\lim(F(\overline{f})(\quot{A_i}{\Gamma}))\big) \ar[rr]^-{\quot{\theta_{F(\overline{f})\quot{A_i}{\Gamma^{\prime}}}}{\Gamma}} \ar[d]_{\quot{\alpha}{\Gamma}(\lim\tau_f^{A_i})} & &\lim \quot{\alpha}{\Gamma}\big(F(\overline{f})(\quot{A_i}{\Gamma^{\prime}})\big) \ar[d]^{\lim \quot{\alpha}{\Gamma}(\tau_f^{A_i})} \\
\quot{\alpha}{\Gamma}\big(\lim (\quot{A_i}{\Gamma})\big) \ar[rr]_-{\quot{\theta_{\quot{A_i}{\Gamma}}}{\Gamma}} & & \lim\big( \quot{\alpha}{\Gamma}\big(\quot{A_i}{\Gamma}\big)\big)
}
\]
in $E(\quot{X}{\Gamma})$. Using this allows us to further derive that
\begin{align*}
&\lim\left(\quot{\alpha}{\Gamma}(\tau_f^{A_i})\right) \circ \quot{\theta_{F(\overline{f})\quot{A_i}{\Gamma^{\prime}}}}{\Gamma} \circ \quot{\alpha}{\Gamma}\big(\quot{\theta_f}{F}\big) \circ \left(\quot{\alpha_{\quot{A}{\Gamma^{\prime}}}}{f}\right)^{-1}\\
&=\quot{\theta_{\quot{A_i}{\Gamma}}}{\Gamma} \circ \quot{\alpha}{\Gamma}\left(\lim \tau_f^{A_i}\right) \circ \quot{\alpha}{\Gamma}\big(\quot{\theta_f}{F}\big) \circ \left(\quot{\alpha_{\quot{A}{\Gamma^{\prime}}}}{f}\right)^{-1} \\
&= \quot{\theta_{\quot{A_i}{\Gamma}}}{\Gamma} \circ \quot{\alpha}{\Gamma}\left(\lim\big(\tau_f^{A_i}\big) \circ \quot{\theta_f}{F}\right) \circ \left(\quot{\alpha_{\quot{A}{\Gamma^{\prime}}}}{f}\right)^{-1} \\
&= \quot{\theta_{\quot{A_i}{\Gamma}}}{\Gamma} \circ \quot{\alpha}{\Gamma}\left(\tau_f^A\right) \circ \left(\quot{\alpha_{\quot{A}{\Gamma^{\prime}}}}{f}\right)^{-1} \\
&= \quot{\theta_{\quot{A_i}{\Gamma}}}{\Gamma} \circ \tau_f^{\underline{\alpha}A}.
\end{align*}
Thus it follows that $\tau_{f}^{\lim \underline{\alpha}A} \circ E(\overline{f})\quot{\theta_{\quot{A_i}{\Gamma^{\prime}}}}{\Gamma^{\prime}} = \quot{\theta_{\quot{A_i}{\Gamma}}}{\Gamma} \circ \tau_f^{\underline{\alpha}A}$, which in turn proves that $\Theta$ is a morphism in $E_G(X)$.
\end{proof}

\begin{corollary}\label{Cor: Section 3: (Co)limit preserving pseudonaturals give (co)limit preserving functors}
Assume that $F$ and $E$ are pre-equivariant pseudofunctors such that each fibre category $F(\quot{X}{\Gamma})$ and $E(\quot{X}{\Gamma})$ admit all limits of shape $I$ and that each fibre functor $F(\overline{f})$ and $E(\overline{f})$ preserves these limits. Assume furthermore that $\alpha:F \Rightarrow E$ is a pseudonatural transformation whose fibre-transition functors preserve limits of shape $I$. Then $\underline{\alpha}:F_G(X) \to E_G(X)$ preserves all limits of shape $I$ in $F_G(X)$. Dually, if $F$ and $E$ are pre-equivariant pseudofunctors such that each fibre category $F(\quot{X}{\Gamma})$ and $E(\quot{X}{\Gamma})$ admit all colimits of shape $I$ and that each fibre functor $F(\overline{f})$ and $E(\overline{f})$ preserves these colimits. Assume furthermore that $\alpha$ is a pseudonatural transformation whose fibre-transition functors preserve colimits of shape $I$. Then $\underline{\alpha}:F_G(X) \to E_G(X)$ preserves all colimits of shape $I$ in $F_G(X)$.
\end{corollary}
\begin{proof}
Simply use Proposition \ref{Prop: Section 3: Change of fibre equivariant functor preserves limit of a specific shape} for every limit to every diagram over the index category $I$ to give the result. To prove this for colimits, dualize Proposition \ref{Prop: Section 3: Change of fibre equivariant functor preserves limit of a specific shape} and proceed mutatis mutandis.
\end{proof}
\begin{corollary}\label{Cor: Section 3: (Co)Continuous pseudonaturals give (co)continuous functors}
Assume that $F$ and $E$ are pre-equivariant pseudofunctors such that each fibre category $F(\quot{X}{\Gamma})$ and $E(\quot{X}{\Gamma})$ is complete and for which each fibre functor $F(\overline{f})$ and $E(\overline{f})$ is continuous, i.e., preserves limits. If $\alpha:F \Rightarrow E$ is a pseudonatural transformation whose fibre-transition functors are all continuous, then $\underline{\alpha}:F_G(X) \to E_G(X)$ is continuous. Dually, if $F$ and $E$ are pre-equivariant pseudofunctors such that each fibre category $F(\quot{X}{\Gamma})$ and $E(\quot{X}{\Gamma})$ is cocomplete and for which each fibre functor $F(\overline{f})$ and $E(\overline{f})$ is cocontinuous and if $\alpha:F \Rightarrow E$ is a pseudonatural transformation whose fibre-transition functors are all cocontinuous, then $\underline{\alpha}:F_G(X) \to E_G(X)$ is cocontinuous
\end{corollary}
\begin{corollary}\label{Cor: Section 3: Additive equivraint change of fibre functors}
Assume that $F$ and $E$ are additive pre-equivariant pseudofunctors and that each fibre-transition functor $\quot{\alpha}{\Gamma}:F(\quot{X}{\Gamma}) \to E(\quot{X}{\Gamma})$ is an additive functor. Then the functor $\underline{\alpha}:F_G(X) \to E_G(X)$ is an additive functor of additive categories.
\end{corollary}

We proceed with a lonely proposition that shows that monoidal pseudonatural transformations give rise to monoidal equivariant functors. However, we first recall what it means to have a natural transformation be monoidal.
\begin{definition}\index{Monoidal Natural Transformation}
Let $(\Cscr,\otimes,I)$ and $(\Dscr, \boxtimes,J)$ be monoidal categories with monoidal functors $F,G:\Cscr \to \Dscr$. A natural transformation $\alpha:F \Rightarrow G$ is \textit{monoidal} if and only if the following hold: 
\begin{itemize}
	\item Given the unit preservation isomorphisms $\gamma:F(I) \xrightarrow{\cong} J$ and $\widetilde{\gamma}:F(I) \to J$, the diagram
	\[
	\xymatrix{
	F(I) \ar[dr]_{\gamma} \ar[r]^{\alpha_I} & G(I) \ar[d]^{\widetilde{\gamma}} \\
	 & J
	}
	\]
	commutes in $\Dscr$;
	\item For any objects $A,B \in \Cscr_0$, given the tensor preservation isomorphisms
	\[
	\theta_{F}^{A,B}:F(A \otimes B) \xrightarrow{\cong} FA \boxtimes FB
	\]
	and
	\[
	\theta_{G}^{A,B}:G(A \otimes B) \xrightarrow{\cong} GA \boxtimes GB,
	\]
	the diagram
	\[
	\xymatrix{
	F(A \otimes B) \ar[d]_{\alpha_{A \otimes B}} \ar[r]^{\theta_F^{A,B}} & FA \boxtimes FB \ar[d]^{\alpha_A \boxtimes \alpha_B} \\
	G(A \otimes B) \ar[r]_{\theta_{G}^{A,B}} & GA \boxtimes GB
	}
	\]
	commutes in $\Dscr$.
\end{itemize} 
\end{definition}
%

\begin{proposition}\label{Prop: Section 3: Monoidal pseudonats give rise to monoidal equivariant functors}
Let $F$ and $E$ be monoidal pre-equivariant pseudofunctors over $X$ and let $\alpha:F \Rightarrow E$ be a pseudonatural transformation such that each functor $\quot{\alpha}{\Gamma}$ is monoidal and each natural transformation $\quot{\alpha}{f}$ is a monoidal transformation. Then the functor $\underline{\alpha}:F_G(X) \to E_G(X)$ is monoidal.
\end{proposition}
\begin{proof}
For any $\Gamma \in \Sf(G)_0$ we let $\quot{I_F}{\Gamma}, \quot{I_E}{\Gamma}$ be the tensorial units in the categories $F(\XGamma)$ and $E(\XGamma)$, respectively; write
\[
\quot{\rho}{\Gamma}:\quot{\alpha}{\Gamma}(\quot{I_F}{\Gamma}) \xrightarrow{\cong} \quot{I_E}{\Gamma}
\]
for the unit preservation isomorphisms induced by the monoidal functors $\quot{\alpha}{\Gamma}$;  let $f \in \Sf(G)_1$ be arbitrary; and as in Theorem \ref{Theorem: Section 2: Monoidal preequivariant pseudofunctor gives monoidal equivariant cat}, let
\[
\sigma_f:F(\overline{f})(\quot{I_F}{\Gamma^{\prime}}) \to \quot{I_F}{\Gamma}
\]
and
\[
\sigma_f^{\prime}:E(\overline{f})(\quot{I_E}{\Gamma^{\prime}}) \to \quot{I_E}{\Gamma}
\]
be the unit preservation isomorphisms induced by the monoidal functors $F(\overline{f})$ and $E(\overline{f})$. Similarly, let
\[
\quot{\theta_f}{F}^{A,B}:F(\overline{f})\left(\AGammap \ds{\Gamma^{\prime}}{F}{\otimes} \BGammap \right) \xrightarrow{\cong} \left(F(\overline{f})(\AGammap)\right) \ds{\Gamma}{F}{\otimes} \left(F(\overline{f})(\BGammap)\right)
\]
and
\[
\quot{\theta_f}{E}^{A,B}:\quot{\theta_f}{E}^{A,B}:E(\overline{f})\left(\AGammap \ds{\Gamma^{\prime}}{E}{\otimes} \BGammap \right) \xrightarrow{\cong} \left(E(\overline{f})(\AGammap)\right) \ds{\Gamma}{E}{\otimes} \left(E(\overline{f})(\BGammap)\right)
\]
be the tensor preservation isomorphisms induced by the monoidal functors $F(\overline{f})$ and $E(\overline{f})$, respectively, where we write $\otimes_E^\Gamma$ and $\otimes_F^\Gamma$ as the tensor functors of $E(\XGamma)$ and $F(\XGamma)$, respectively. Finally let
\[
\quot{\theta_{\Gamma}}{\alpha}^{A,B}:\quot{\alpha}{\Gamma}\left(\AGamma \ds{\Gamma}{F}{\otimes} \BGamma\right) \xrightarrow{\cong} \quot{\alpha}{\Gamma}\left(\AGamma\right) \ds{\Gamma}{E}{\otimes} \quot{\alpha}{\Gamma}\left(\BGamma\right)
\]
be the tensor preservation isomorphism induced by $\quot{\alpha}{\Gamma}$.

We first prove that $\underline{\alpha}$ preserves the unit object. Let $I_F := \lbrace \quot{I_F}{\Gamma} \; | \; \Gamma \in \Sf(G)_0 \rbrace$ be the object set of the unit object $(I_F, T_{I_F})$ in $F_G(X)$; note that by Theorem \ref{Theorem: Section 2: Monoidal preequivariant pseudofunctor gives monoidal equivariant cat}, $T_{I_F} = \lbrace \sigma_f \; | \; f \in \Sf(G)_1 \rbrace$. Similarly, let $I_E := \lbrace \quot{I_E}{\Gamma} \; | \; \Gamma \in \Sf(G)_0 \rbrace$ be the object set of the unit object $(I_E,T_{I_E})$ in $E_G(X)$; as before, $T_{I_E} = \lbrace \sigma_{f}^{\prime} \; | \; f \in \Sf(G)_1 \rbrace$. Define the collection of morphisms $P$ by
\[
P := \lbrace \quot{\rho}{\Gamma}:\quot{\alpha}{\Gamma}(\quot{I_F}{\Gamma}) \to \quot{I_E}{\Gamma} \; | \; \Gamma \in \Sf(G)_0 \rbrace.
\]
Because each of the $\quot{\rho}{\Gamma}$ is an isomorphism, if we can prove that $P$ is a $E_G(X)$-morphism, we will have proved that $\underline{\alpha}$ preserves unit objects up to isomorphism. To this end fix an $f \in \Sf(G)_1$ and write $\Dom f = \Gamma$ and $\Codom f = \Gamma^{\prime}$. We must show that the diagram
\[
\xymatrix{
E(\overline{f})\big(\quot{\alpha}{\Gamma^{\prime}}(\quot{I_F}{\Gamma^{\prime}})\big) \ar[rr]^-{E(\overline{f}\quot{\rho}{\Gamma^{\prime}})} \ar[d]_{\tau_f^{\underline{\alpha}I_F}} & & E(\overline{f})(\quot{I_E}{\Gamma^{\prime}}) \ar[d]^{\tau_f^{I_E}} & \\
\quot{\alpha}{\Gamma}(\quot{I_F}{\Gamma}) \ar[rr]_-{\quot{\rho}{\Gamma}} & & \quot{I_E}{\Gamma}
}
\]
commutes. For this observe that the unit preservation isomorphism 
\[
\gamma: \quot{\alpha}{\Gamma}\big(F(\overline{f})(\quot{I_F}{\Gamma^{\prime}})\big) \to \quot{I_E}{\Gamma}
\] 
has the form
\[
\gamma = \quot{\rho}{\Gamma} \circ \quot{\alpha}{\Gamma}(\sigma_f) = \quot{\rho}{\Gamma} \circ \quot{\alpha}{\Gamma}(\tau_f^{I_F})
\]
while the unit preservation isomorphism $\widetilde{\gamma}:E(\overline{f})\big(\quot{\alpha}{\Gamma^{\prime}}(\quot{I_F}{\Gamma^{\prime}})\big) \to \quot{I_E}{\Gamma}$ takes the form
\[
\widetilde{\gamma} = \sigma_f^{\prime} \circ E(\overline{f})\quot{\rho}{\Gamma^{\prime}} = \tau_f^{I_E} \circ E(\overline{f})\quot{\rho}{\Gamma^{\prime}}.
\]
Using that $\quot{\alpha}{f}$ is a monoidal natural transformation gives that the diagram
\[
\xymatrix{
\quot{\alpha}{\Gamma}\big(F(\overline{f})(\quot{I_F}{\Gamma^{\prime}})\big) \ar[rr]^-{\gamma} \ar[d]_{\quot{\alpha_{\quot{I_F}{\Gamma^{\prime}}}}{f}} & & \quot{I_E}{\Gamma} \\
E(\overline{f})\big(\quot{\alpha}{\Gamma^{\prime}}(\quot{I_F}{\Gamma^{\prime}})\big) \ar[urr]_{\widetilde{\gamma}}
}
\]
commutes. This gives, however, that
\begin{align*}
\quot{\rho}{\Gamma} \circ \tau_{f}^{\underline{\alpha}I_F} &= \quot{\rho}{\Gamma} \circ \quot{\alpha}{\Gamma}(\tau_f^{I_F}) \circ \left(\quot{\alpha_{\quot{I_F}{\Gamma^{\prime}}}}{f}\right)^{-1}= \gamma \circ \left(\quot{\alpha_{\quot{I_F}{\Gamma^{\prime}}}}{f}\right)^{-1} \\
&= \widetilde{\gamma} \circ \quot{\alpha_{\quot{I_F}{\Gamma^{\prime}}}}{f} \circ \left(\quot{\alpha_{\quot{I_F}{\Gamma^{\prime}}}}{f}\right)^{-1} = \widetilde{\gamma} 
= \tau_f^{I_E} \circ E(\overline{f})\quot{\rho}{\Gamma^{\prime}}.
\end{align*}
Thus $P$ is a morphism in $E_G(X)$ and hence $\underline{\alpha}$ preserves the unit object.

We now prove that $\underline{\alpha}$ preserves the tensor functors up to natural isomorphism. For this fix objects $A, B \in F_G(X)_0$ and consider the isomorphisms, for all $\Gamma \in \Sf(G)_0$,
\[
\quot{\theta_{\Gamma}}{\alpha}^{A,B}:\quot{\alpha}{\Gamma}\left(\AGamma \ds{\Gamma}{F}{\otimes} \BGamma\right) \xrightarrow{\cong} \quot{\alpha}{\Gamma}\left(\AGamma\right) \ds{\Gamma}{E}{\otimes} \quot{\alpha}{\Gamma}\left(\BGamma\right).
\]
Define the collection $\Theta_{A,B}:\underline{\alpha}(A \otimes_F B) \to \underline{\alpha}A  \otimes_E \underline{\alpha}B$ by
\[
\Theta_{A,B}:= \left\lbrace \quot{\theta_{\Gamma}^{A,B}}{\alpha} \; : \; \Gamma \in \Sf(G)_0 \right\rbrace.
\]
Note that because each of the $\quot{\theta_{\Gamma}}{\alpha}$ are natural isomorphisms, if we can prove that $\Theta_{A,B}$ is a $E_G(X)$-morphism, then it will automatically induce a natural isomorphism $\underline{\alpha} \circ \big((-)\otimes_F(-)\big) \Rightarrow \underline{\alpha}(-) \otimes_E \underline{\alpha}(-)$. Thus to complete the proof of the proposition it suffices to prove that the diagram
\[
\xymatrix{
E(\overline{f})\left(\quot{\alpha}{\Gamma^{\prime}}\left(\AGammap \ds{\Gamma^{\prime}}{F}{\otimes} \BGammap \right)\right) \ar[rr]^-{E(\overline{f})\quot{\theta_{\Gamma^{\prime}}}{\alpha}^{A,B}} \ar[d]_{\tau_{f}^{\underline{\alpha}(A \otimes_F B)}} & & E(\overline{f})\left(\quot{\alpha}{\Gamma^{\prime}}(\AGammap) \ds{\Gamma^{\prime}}{E}{\otimes} \quot{\alpha}{\Gamma^{\prime}}(\BGammap)\right) \ar[d]^{\tau_{f}^{\underline{\alpha}A \otimes_E \underline{\alpha}B}} \\
\quot{\alpha}{\Gamma}\left(\AGamma \ds{\Gamma}{F}{\otimes} \BGamma\right) \ar[rr]_-{\quot{\theta_{\Gamma}}{\alpha}^{A,B}} & & \quot{\alpha}{\Gamma}(\AGamma) \ds{\Gamma}{E}{\otimes} \quot{\alpha}{\Gamma}(\BGamma)
}
\]
commutes. For this we first compute that by Theorems \ref{Theorem: Section 2: Monoidal preequivariant pseudofunctor gives monoidal equivariant cat} and \ref{Thm: Section 3: Psuedonatural trans lift to equivariant functors} we have one hand
\begin{align*}
\tau_f^{\underline{\alpha}(A \otimes_F B)} &= \quot{\alpha}{\Gamma}(\tau_f^{A \otimes_F B}) \circ \left(\quot{\alpha_{\AGammap \otimes \BGammap}}{f}\right)^{-1} \\
&= \quot{\alpha}{\Gamma}\left(\tau_f^A \ds{\Gamma}{F}{\otimes} \tau_f^B \right) \circ \quot{\alpha}{\Gamma}(\quot{\theta_{f}}{F}^{A,B}) \circ \left(\quot{\alpha_{\AGammap \otimes \BGammap}}{f}\right)^{-1},
\end{align*}
while on the other hand
\begin{align*}
 \tau_f^{\underline{\alpha}A \otimes_E \underline{\alpha}B} &=  \left(\tau_f^{\underline{\alpha}A} \ds{\Gamma}{E}{\otimes} \tau_f^{\underline{\alpha}B}\right) \circ \quot{\theta_f}{E}^{\underline{\alpha}A,\underline{\alpha}B} \\ 
 &= \left(\left(\quot{\alpha}{\Gamma}(\tau_f^A) \circ \left(\quot{\alpha_{\AGammap}}{f}\right)^{-1}\right) \ds{\Gamma}{E}{\otimes} \left(\quot{\alpha}{\Gamma}(\tau_f^B) \circ \left(\quot{\alpha_{\BGammap}}{f}\right)^{-1}\right)\right) \circ \quot{\theta_f}{E}^{\underline{\alpha}A,\underline{\alpha}B} \\
 &= \left(\quot{\alpha}{\Gamma}(\tau_f^A) \ds{\Gamma}{E}{\otimes} \quot{\alpha}{\Gamma}(\tau_f^{B}) \right) \circ \left(\quot{\alpha_{\AGammap}}{f} \ds{\Gamma}{E}{\otimes} \quot{\alpha_{\BGammap}}{f}\right)^{-1} \circ \quot{\theta_f}{E}^{\underline{\alpha}A,\underline{\alpha}B}.
\end{align*}
Now observe that the tensor preservation isomorphism
\[
\quot{\alpha}{\Gamma}\left(F(\overline{f})\left(\AGammap \ds{\Gamma^{\prime}}{F}{\otimes} \BGammap\right)\right) \xrightarrow[\quot{\theta_f}{\alpha \star F}^{A,B}]{\cong} \left(\quot{\alpha}{\Gamma}\big(F(\overline{f})(\AGammap)\big) \right) \ds{\Gamma}{E}{\otimes} \left(\quot{\alpha}{\Gamma}\big(F(\overline{f})(\BGammap)\big)\right)
\]
has the form
\[
\quot{\theta_f}{\alpha \star F}^{F(\overline{f})A,F(\overline{f})B} = \quot{\theta_{\Gamma}}{\alpha}^{A,B} \circ \quot{\alpha}{\Gamma}\big(\quot{\theta_f}{F}^{A,B}\big)
\]
while the tensor preservation isomorphism
\[
E(\overline{f})\left(\quot{\alpha}{\Gamma^{\prime}}\left(\AGammap \ds{\Gamma^{\prime}}{F}{\otimes} \BGammap\right)\right) \xrightarrow[\quot{\theta_{f}}{E \star \alpha}^{A,B}]{\cong} \left(E(\overline{f})\big(\quot{\alpha}{\Gamma^{\prime}}(\AGammap)\big)\right) \ds{\Gamma}{E}{\otimes} \left(E(\overline{f})\big(\quot{\alpha}{\Gamma^{\prime}}(\BGammap)\big)\right)
\]
has the form
\[
\quot{\theta_{f}}{E \star \alpha}^{A,B} = \quot{\theta_f}{E}^{\underline{\alpha}A,\underline{\alpha}B} \circ E(\overline{f})\big(\quot{\theta_{\Gamma^{\prime}}}{\alpha}^{A,B}\big),
\]
where in both cases $\alpha \star F$ and $E \star \alpha$ denote the whiskering of $\alpha$ by $F$ and $E$  by $\alpha$, respectively. Note also that since $\quot{\alpha}{f}$ is a monoidal natural isomorphism, so is $\quot{\alpha}{f}^{-1}$. The monoidal transformation $\quot{\alpha}{f}^{-1}$ also satisfies the identity
\[
\left(\quot{\alpha_{\AGammap}}{f} \ds{\Gamma}{E}{\otimes} \quot{\alpha_{\BGammap}}{f}\right)^{-1} \circ \quot{\theta_f}{E \star \alpha}^{A,B} = \quot{\theta_f}{\alpha \star F}^{A,B} \circ \left(\quot{\alpha_{\AGammap \otimes \BGammap}}{f}\right)^{-1}.
\]
Using these computations and identities in various ways we compute that
\begin{align*}
&\tau_f^{\underline{\alpha}A \otimes_E \underline{\alpha}B} \circ E(\overline{f})\quot{\theta_{\Gamma^{\prime}}}{\alpha}^{A,B} \\
&= \left(\tau_f^{\underline{\alpha}A} \ds{\Gamma}{E}{\otimes} \tau_f^{\underline{\alpha}B}\right) \circ \quot{\theta_f}{E}^{\underline{\alpha}A,\underline{\alpha}B} \circ E(\overline{f})\quot{\theta_{\Gamma^{\prime}}}{\alpha}^{A,B} = \left(\tau_f^{\underline{\alpha}A} \ds{\Gamma}{E}{\otimes} \tau_f^{\underline{\alpha}B}\right) \circ \quot{\theta_f}{E \star \alpha}^{A,B} \\
&= \left(\quot{\alpha}{\Gamma}(\tau_f^A) \ds{\Gamma}{E}{\otimes} \quot{\alpha}{\Gamma}(\tau_f^{B}) \right) \circ \left(\quot{\alpha_{\AGammap}}{f} \ds{\Gamma}{E}{\otimes} \quot{\alpha_{\BGammap}}{f}\right)^{-1} \circ \quot{\theta_f}{E \star \alpha}^{A,B} \\
&= \left(\quot{\alpha}{\Gamma}(\tau_f^A) \ds{\Gamma}{E}{\otimes} \quot{\alpha}{\Gamma}(\tau_f^{B}) \right) \circ \quot{\theta_f}{\alpha \star F}^{A,B} \circ \left(\quot{\alpha_{\AGammap \otimes \BGammap}}{f}\right)^{-1} \\
&= \left(\quot{\alpha}{\Gamma}(\tau_f^A) \ds{\Gamma}{E}{\otimes} \quot{\alpha}{\Gamma}(\tau_f^{B}) \right)\circ \quot{\theta_{\Gamma}}{\alpha}^{F(\overline{f})A,F(\overline{f})B} \circ \quot{\alpha}{\Gamma}\big(\quot{\theta_f}{F}^{A,B}\big) \circ \left(\quot{\alpha_{\AGammap \otimes \BGammap}}{f}\right)^{-1} \\
&= \quot{\theta_{\Gamma}}{\alpha}^{A,B} \circ \quot{\alpha}{\Gamma}\left(\tau_f^A \ds{\Gamma}{F}{\otimes} \tau_f^B\right) \circ \quot{\alpha}{\Gamma}\big(\quot{\theta_f}{F}^{A,B}\big) \circ \left(\quot{\alpha_{\AGammap \otimes \BGammap}}{f}\right)^{-1} \\
&= \quot{\theta_{\Gamma}}{\alpha}^{A,B} \circ \tau_{f}^{\underline{\alpha}(A \otimes_F B)}.
\end{align*}
This gives that $\tau_f^{\underline{\alpha}A \otimes_E \underline{\alpha}B} \circ E(\overline{f})\quot{\theta_{\Gamma^{\prime}}}{\alpha}^{A,B} = \quot{\theta_{\Gamma}}{\alpha}^{A,B} \circ \tau_{f}^{\underline{\alpha}(A \otimes_F B)}$ which proves that $\Theta_{A,B}$ is a morphism in $E_G(X)$.
\end{proof}

We now move to discuss how modifications ($3$-morphisms in the strict $3$-category $\twoCat$ of $2$-categories, i.e., morphisms between pseudonatural transformations) give rise to natural transformations of equivariant functors. These will not only serve as our first main source of examples of equivariant natural transformations, but will also give us a uniform language and structural strategy for building adjoints between equivariant categories that come from adjunctions that live between every fibre category. We proceed by first recalling the definition of a modification and then provide a structural lemma showing that modifications between pseudonatural transformations of pre-equivariant pseudofunctors give rise to natural transformations of equivariant functors.
\begin{definition}\index{Modification}
Let $F,E:\Cscr^{\op} \to \fCat$ be pseudofunctors for some $1$-category $\Cscr$ and let $\alpha,\beta:F \Rightarrow E$ be pseudonatural transformations. A modification $\eta:\alpha \Rrightarrow \beta$ is an assignment $\eta:\Cscr_0 \to \fCat_2$ where each natural transformation $\eta_A$ fits into a $2$-cell
\[
\begin{tikzcd}
F(A) \ar[r,bend left = 30, ""{name = U}]{}{\alpha_A} \ar[r,swap, bend right = 30,""{name = L}]{}{\beta_A} & E(A) \ar[from = U, to = L, Rightarrow,shorten >= 5pt, shorten <= 5pt]{}{\eta_A}
\end{tikzcd}
\]
in $\fCat$ such that for all morphisms $g:A \to B$ in $\Cscr$, the diagram of functors/natural transformations
\[
\xymatrix{
\alpha_A \circ F(g) \ar[rr]^-{\eta_A \ast \iota_{F(g)}} \ar[d]_{\alpha_g} & & \beta_A \circ F(g) \ar[d]^{\beta_g} \\
E(g) \circ \alpha_B \ar[rr]_-{\iota_{E(g)}\ast\eta_B} & & E(g) \circ \beta_B
}
\] 
commutes (where $\ast$ denotes horizontal composition of natural transformations).
\end{definition}
\begin{lemma}\label{Lemma: Modifications give equivariant natural transformations}
Let $F,E:\SfResl_G(X)^{\op} \to \fCat$ be pre-equivariant natural transformations and let $\alpha,\beta:F \Rightarrow E$ be pseudonatural transformations. If $\eta:\alpha \Rrightarrow \beta$ is a modification, then there exists a natural transformation $\underline{\eta}$ fitting into the $2$-cell
\[
\begin{tikzcd}
F_G(X) \ar[r,bend left = 30, ""{name = U}]{}{\underline{\alpha}} \ar[r,swap, bend right = 30,""{name = L}]{}{\underline{\beta}} & E_G(X) \ar[from = U, to = L, Rightarrow,shorten >= 5pt, shorten <= 5pt]{}{\underline{\eta}}
\end{tikzcd}
\]
given on objects $C = \lbrace \quot{C}{\Gamma} \; | \; \Gamma \in \Sf(G)_0\rbrace$ of $F_G(X)$ by
\[
\underline{\eta}_C := \lbrace \quot{\eta}{\Gamma}_{\quot{C}{\Gamma}}:\quot{\alpha}{\Gamma}(\quot{C}{\Gamma}) \to \quot{\beta}{\Gamma}(\quot{C}{\Gamma}) \; | \; \Gamma \in \Sf(G)_0\rbrace.
\]
\end{lemma}
\begin{proof}
Before we prove the naturality of $\underline{\eta}$, we must prove that for all objects $A$ of $F_G(X)$, the membership $\underline{\eta}_A \in E_G(X)(\underline{\alpha}A, \underline{\beta}A)$. That is, for all $f \in \Sf(G)_1$ (with $\Dom f = \Gamma$ and $\Codom f = \Gamma^{\prime}$) we have that
\[
\tau_f^{\underline{\beta}A} \circ E(\overline{f})\quot{\eta_{\quot{A}{\Gamma^{\prime}}}}{\Gamma^{\prime}} = \quot{\eta_{\quot{A}{\Gamma}}}{\Gamma} \circ \tau_f^{\underline{\alpha}A}.
\]
To this end recall that by Theorem \ref{Thm: Section 3: Psuedonatural trans lift to equivariant functors}
\[
\tau_{f}^{\underline{\alpha}A} = \quot{\alpha}{\Gamma}(\tau_f^A) \circ \left(\quot{\alpha_{\quot{A}{\Gamma^{\prime}}}}{f}\right)^{-1}, \qquad \tau_f^{\underline{\beta}A} = \quot{\beta}{\Gamma}(\tau_f^A) \circ \left(\quot{\beta_{\quot{A}{\Gamma^{\prime}}}}{f}\right)^{-1},
\]
and as usual we derive that the diagram
\[
\xymatrix{
F(\overline{f})\big(\quot{\alpha}{\Gamma^{\prime}}(\quot{A}{\Gamma^{\prime}})\big) \ar[rr]^-{E(\overline{f})\quot{\eta_{\quot{A}{\Gamma^{\prime}}}}{\Gamma^{\prime}}} \ar[d]_{\left(\quot{\alpha_{\quot{A}{\Gamma^{\prime}}}}{f}\right)^{-1}} & & E(\overline{f})\big(\quot{\beta}{\Gamma^{\prime}}(\quot{A}{\Gamma^{\prime}})\big) \\
\quot{\alpha}{\Gamma}\big(F(\overline{f})(\quot{A}{\Gamma^{\prime}})\big) \ar[rr]_-{\quot{\eta_{F(\overline{f})\quot{A}{\Gamma^{\prime}}}}{\Gamma}} & & \quot{\beta}{\Gamma}\big(F(\overline{f})(\quot{A}{\Gamma^{\prime}})\big) \ar[u]_{\quot{\beta_{\quot{A}{\Gamma^{\prime}}}}{f}}
}
\]
commutes. Using this and the naturality of $\quot{\eta}{\Gamma}$ gives us that
\begin{align*}
\tau_f^{\underline{\beta}f} \circ E(\overline{f})\quot{\eta_{\quot{A}{\Gamma^{\prime}}}}{\Gamma^{\prime}} &=\quot{\beta}{\Gamma}(\tau_f^A) \circ \left(\quot{\beta_{\AGammap}}{f}\right)^{-1} \circ E(\overline{f})\quot{\eta_{\AGammap}}{\Gamma^{\prime}} \\
&= \quot{\beta}{\Gamma}(\tau_f^A) \circ \left(\quot{\beta_{\AGammap}}{f}\right)^{-1} \circ \quot{\beta_{\AGammap}}{f} \circ \quot{\eta_{F(\overline{f})\AGammap}}{\Gamma} \circ \left(\quot{\alpha_{\AGammap}}{f}\right)^{-1} \\
&= \quot{\beta}{\Gamma}(\tau_f^A) \circ \quot{\eta_{F(\overline{f}\AGammap}}{\Gamma} \circ \left(\quot{\alpha_{\AGammap}}{f}\right)^{-1} \\
&= \quot{\eta_{\AGamma}}{\Gamma} \circ \quot{\alpha}{\Gamma}(\tau_f^{A}) \circ \left(\quot{\alpha_{\AGammap}}{f}\right)^{-1} = \quot{\eta_{\AGamma}}{\Gamma} \circ \tau_f^{\underline{\alpha}A},
\end{align*}
which proves that $\underline{\eta}_A:\underline{\alpha}A \to \underline{\beta}A$ is a $E_G(X)$-morphism.

We now verify the naturality of $\underline{\eta}$. For this let $P = \lbrace \quot{\rho}{\Gamma} \; | \; \Gamma \in \Sf(G)_0\rbrace \in F_G(X)(A,B)$. We must prove that $\underline{\eta}_B \circ \underline{\alpha}P = \underline{\beta}P \circ \underline{\eta}_A$.  To see this note for all $\Gamma \in \Sf(G)_0$ we get that
\[
\quot{\beta}{\Gamma}(\quot{\rho}{\Gamma}) \circ \quot{\eta_{\quot{A}{\Gamma}}}{\Gamma} = \quot{\eta_{\quot{B}{\Gamma}}}{\Gamma} \circ \quot{\alpha}{\Gamma}(\quot{\rho}{\Gamma})
\]
which shows the commutativity of
\[
\xymatrix{
\underline{\alpha}A \ar[r]^-{\underline{\alpha}P} \ar[d]_{\underline{\eta}_A} & \underline{\alpha}B \ar[d]^{\underline{\eta}_B} \\
\underline{\beta}A \ar[r]_-{\underline{\beta}P} & \underline{\beta}B
}
\]
in $E_G(X)$. Thus $\underline{\eta}$ is a natural transformation between $\underline{\alpha}$ and $\underline{\beta}$.
\end{proof}
\begin{definition}\index{Equivariant Natural Transformation}
Let $\underline{\alpha},\underline{\beta}:F_G(X) \to E_G(X)$ be equivariant functors. An {equivariant natural transformation} is a natural transformation $\underline{\eta}:\underline{\alpha} \to \underline{\beta}$ where $\underline{\eta}$ is induced by a modification $\eta:\alpha \Rrightarrow \beta$.
\end{definition}
\begin{example}
Consider the equivariant categories $F_G(X) = D^b_G(X) = E_G(X)$ be induced by the pseudofunctor $D^b(-)$ we have seen before and let $\alpha = [1] \circ [2]:D^b(-) \Rightarrow D^b(-)$ and $\beta = [2] \circ [1]:D^b(-) \Rightarrow D^b(-)$ be the degree $1 + 2$ and degree $2 + 1$ shift pseudonatural transformations. Then the natural isomorphisms
\[
\quot{[1]}{\Gamma} \circ \quot{[2]}{\Gamma} \xrightarrow[\eta_\Gamma]{\cong} \quot{[2]}{\Gamma} \circ \quot{1}{\Gamma}
\]
assemble to a modification $\eta:[2] \circ [1] \Rightarrow [1] \circ [2]$ which gives the equivariant natural transformation witnessing that $[2] \circ [1] \cong [1] \circ [2]$.
\end{example}

Going forward we will need to talk about a $2$-category of pre-equivariant pseudofunctors $\SfResl_G(X)^{\op} \to \fCat$. We denote this $2$-category by writing $\fPreEq(\SfResl_G(X)^{\op},\fCat)$\index[notation]{PreEqFrak@$\fPreEq(\SfResl_G(X)^{\op},\fCat)$} and define it as follows:
\begin{itemize}
	\item $0$-morphisms: Pre-equivariant pseudofunctors $F:\SfResl_G(X)^{\op} \to \fCat$.
	\item $1$-morphisms: Pseudonatural transformations $\alpha:F \Rightarrow E$.
	\item $2$-morphisms: Modifications $\eta:\alpha \to \beta$.
	\item $1$-morphism composition: Vertical composition of pseudonatural transformations.
	\item $2$-morphism vertical composition: As in $\twoCat$.
	\item $2$-morphism horizontal composition: As in $\twoCat$.
\end{itemize}

\begin{lemma}\label{Lemma: Section 3: The embedding of pre-equivariant functors to equivariant categories is a strict 2-functor}
There is a strict $2$-functor 
\[
\mathbf{i}:\fPreEq(\SfResl_G(X)^{\op},\fCat) \to \fGEqCat_{/X}
\] 
which sends pre-equivariant pseudofunctors $F$ to their equivariant categories $F_G(X)$, sends pseudonatural transformations $\alpha:F \to E$ to their induced equivariant functor $\underline{\alpha}$, and finally sends modifications $\eta$ to their induced equivariant natural transformations $\underline{\eta}$.
\end{lemma}
\begin{proof}
Because the $2$-functor is claimed to be strict, we have $\phi_{f,g} = \id_{\mathbf{i}(f,g)}.$ As such, we need only verify that $\mathbf{i}$ strictly preserves identity transformations, composition, and vertical composition of modifications.

We fist note that the preservation of identities is immediate by Corollary \ref{Cor: Section 3: Identity functors are change of fibre functors}. To see that (vertical) composition of pseudonatural transformations is preserved, let $\alpha:F \to E$ and $\beta:E \to L$ be pseudonatural transformations. Then the composite $\underline{\beta \circ\alpha}$ is given on objects $(A,T_A)$ as follows: First, the action of $\underline{\beta \circ \alpha}$ is given as
\begin{align*}
(\underline{\beta \circ \alpha})A &= \lbrace \quot{(\beta \circ \alpha)}{\Gamma}\quot{A}{\Gamma} \; | ; \Gamma \in \Sf(G)_0\rbrace = \lbrace (\quot{\beta}{\Gamma} \circ \quot{\alpha}{\Gamma})\quot{A}{\Gamma} \; | \; \Gamma \in \Sf(G)_0 \rbrace \\
&= (\underline{\beta} \circ \underline{\alpha})A.
\end{align*}
We also now need to prove that $\tau_{f}^{(\underline{\beta \circ \alpha})A} = \tau_f^{\underline{\beta}(\underline{\alpha}A)}$ for all $f \in \Sf(G)_1$. For this let $\Dom f = \Gamma$ and $\Codom f = \Gamma^{\prime}$. We now compute that on one hand
\[
\tau_{f}^{(\underline{\beta \circ \alpha})A} = \quot{\beta \circ \alpha}{\Gamma}(\tau_f^{A}) \circ \left(\quot{\alpha_{\quot{A}{\Gamma^{\prime}}}}{f}\right)^{-1} = \quot{\beta}{\Gamma}\left(\quot{\alpha}{\Gamma}\big(\tau_f^A\big)\right) \circ \left(\quot{\alpha_{\quot{A}{\Gamma^{\prime}}}}{f}\right)^{-1}.
\]
On the other hand we get from the pseudonaturality of $\alpha$ and $\beta$ that
\begin{align*}
\tau_f^{\underline{\beta}(\underline{\alpha}A)} &= \quot{\beta}{\Gamma}\left(\tau_f^{\underline{\alpha}A}\right) \circ \left(\quot{\beta_{\quot{\underline{\alpha}A}{\Gamma^{\prime}}}}{f}\right)^{-1}\\
& = \quot{\beta}{\Gamma}\left(\quot{\alpha}{\Gamma}\big(\tau_f^A\big) \circ \left(\quot{\alpha_{\quot{A}{\Gamma^{\prime}}}}{f}\right)^{-1}\right) \circ \left(\quot{\beta_{\quot{\alpha}{\Gamma^{\prime}}(\quot{A}{\Gamma^{\prime}})}}{f }\right)^{-1}  \\
&= \quot{\beta}{\Gamma}\left(\quot{\alpha}{\Gamma}\big(\tau_f^A\big)\right) \circ \quot{\beta}{\Gamma}\left(\quot{\alpha_{\quot{A}{\Gamma^{\prime}}}}{f}\right)^{-1} \circ \left(\quot{\beta_{\quot{\alpha}{\Gamma^{\prime}}(\quot{A}{\Gamma^{\prime}})}}{f }\right)^{-1} \\
&= \quot{\beta}{\Gamma}\left(\quot{\alpha}{\Gamma}\big(\tau_f^A\big)\right) \circ \left(\quot{\alpha_{\quot{A}{\Gamma^{\prime}}}}{f}\right)^{-1} \\
&= \tau_{f}^{(\underline{\beta \circ \alpha})A}
\end{align*}
which proves that $\tau_f^{\underline{\beta}(\underline{\alpha}A)} = \tau_{f}^{(\underline{\beta \circ \alpha})A}$. Thus it follows that
\begin{align*}
\mathbf{i}(\beta \circ \alpha)(A,T_A) &= (\underline{\beta \circ \alpha})(A,T_A) = ((\underline{\beta \circ \alpha})A, T_{\underline{\beta \circ \alpha}A}) = ((\underline{\beta} \circ \underline{\alpha})A, T_{(\underline{\beta} \circ \underline{\alpha})A}) \\
&= (\underline{\beta} \circ \underline{\alpha})(A,T_A) = (\mathbf{i}(\beta) \circ \mathbf{i}(\alpha))(A,T_A),
\end{align*}
and similarly $\mathbf{i}(\beta \circ \alpha)(P) = \mathbf{i}(\beta)\big(\mathbf{i}(\alpha)P\big)$ on morphisms $P$. Thus $\mathbf{i}$ preserves $1$-compositions strictly. 

To see that $2$-morphisms work out as expected, let $F,E:\SfResl_G(X)^{\op} \to \fCat$ be pre-equivariant pseudofunctors with pseudonatural transformations $\alpha,\beta,\gamma:F \Rightarrow E$ and modifications $\eta:\alpha \Rrightarrow \beta$ and $\epsilon:\beta \Rrightarrow \gamma$ as in the diagram
\[
\begin{tikzcd}
F \ar[rr, bend left = 60, ""{name = U}]{}{\alpha} \ar[rr, ""{name = M}]{}[description]{\beta} \ar[rr, swap, bend right = 60, ""{name = L}]{}{\gamma} & & E \ar[from = U, to = M, Rightarrow, shorten <= 4pt, shorten >= 4pt]{}{\eta} \ar[from = M, to = L, Rightarrow, shorten <= 4pt, shorten >= 4pt]{}{\epsilon}
\end{tikzcd}
\]
in the $2$-category $\twoCat(\SfResl_G(X)^{\op},\fCat)$. We now must check that
\begin{align*}
\mathbf{i}(\epsilon \circ \eta) = \mathbf{i}\epsilon \circ \mathbf{i}\eta.
\end{align*}
However, this is a routine check as  for all $\Gamma \in \Sf(G)_0$ and for all $A \in F_G(X)_0$
\begin{align*}
\quot{(\mathbf{i}(\epsilon \circ \eta)_A)}{\Gamma} &= \quot{(\underline{\epsilon \circ \eta})_A}{\Gamma} = \quot{(\epsilon \circ \eta)_{\AGamma}}{\Gamma} = \quot{\epsilon}{\Gamma}_{\AGamma} \circ \quot{\eta_{\AGamma}}{\Gamma} = \quot{\underline{\epsilon}_A}{\Gamma} \circ \quot{\underline{\eta}_A}{\Gamma} \\
&= \quot{\mathbf{i}(\epsilon)_A}{\Gamma} \circ \quot{\mathbf{i}(\eta)_A}{\Gamma}
\end{align*}
so we conclude that $\mathbf{i}(\epsilon \circ \eta) = \mathbf{i}(\epsilon) \circ \mathbf{i}(\eta)$, as was desired. The verification for horizontal composition is similar and omitted. Finally, the verification of the last remaining identity is trivial because the compositors $\phi_{\alpha,\beta}:\mathbf{i}(\beta) \circ \mathbf{i}(\alpha) \Rightarrow \mathbf{i}(\beta \circ \alpha)$ are all identity natural transformations. Thus $\mathbf{i}$ is a strict $2$-functor.
\end{proof}

We will now use modifications to describe adjoints that live in the $2$-category $\fPreEq(\SfResl_G(X)^{\op},\fCat)$. One benefit of this is that because adjoints in these $2$-categories are equational, they are preserved by any $2$-functor. In particular, the strict $2$-functor $\mathbf{i}$ of Lemma \ref{Lemma: Section 3: The embedding of pre-equivariant functors to equivariant categories is a strict 2-functor} lifts adjoints in $\fPreEq_{/X}(\SfResl_G(X)^{\op},\fCat)$ to adjoints in $\fGEqCat_{/X}$; because the $2$-category $\fGEqCat_{/X}$ is a subcategory of $\fCat$, the adjoints in $\fGEqCat_{/X}$ coincide with adjoints in $\fCat$. Before presenting this theorem, however, we will recall the definition of adjoints in a $2$-category as given in \cite{RiehlVertityElementss}; this definition is equivalent to the defnition of adjoints given in \cite{FreeAdj}, which defines a minimal $2$-category of the syntactic data required of an adjoint pair of functors and then gives an adjunction as a $2$-functor from this free adjunction category into a $2$-category $\Cfrak$.
\begin{definition}[{\cite[Defition 2.1.1]{RiehlVertityElementss}}]\index{Adjunctions!In a $2$-Category}
An adjunction in a $2$-category $\Cfrak$ is a pair of objects $A$ and $B$ together with $1$-morphisms $f:A \to B$ and $g:B \to A$ such that there exist $2$-morphisms
\[
\begin{tikzcd}
A \ar[r, bend left = 30, ""{name = UL}]{}{\id_A} \ar[bend right = 30, swap, r, ""{name = BL}]{}{g \circ f} & A \ar[from = UL, to = BL, Rightarrow, shorten >= 2pt, shorten <= 2pt]{}{\eta} & B \ar[r, bend left = 30, ""{name = UR}]{}{f \circ g} \ar[r, swap, bend right = 30, ""{name = BR}]{}{\id_B} & B \ar[from = UR, to = BR, Rightarrow, shorten <= 2pt, shorten >= 2pt]{}{\epsilon}
\end{tikzcd}
\]
which satisfy the triangle identities below. The pasting diagram
\[
\begin{tikzcd}
 & B \ar[rr, equals] \ar[dr]{}{g} & {} & B \\
A \ar[ur]{}{f} \ar[rr, equals] & {} \ar[u, Rightarrow, shorten >=5pt, shorten <= 5pt]{}{\eta} & A \ar[ur, swap]{}{f} \ar[u, swap, Rightarrow, shorten >=2pt, shorten <= 2pt]{}{\epsilon}
\end{tikzcd}
\]
is equal to the $2$-cell
\[
\begin{tikzcd}
A \ar[r, bend left = 30, ""{name = U}]{}{f} \ar[r, swap, bend right = 30, ""{name = B}]{}{f} & B \ar[from = U, to = B, Rightarrow, shorten >=2pt, shorten <= 2pt]{}{\iota_f}
\end{tikzcd}
\]
and the pasting diagram
\[
\begin{tikzcd}
 & A \ar[rr, equals] \ar[dr]{}{f} \ar[d, Rightarrow, shorten >=5pt, shorten <= 5pt]{}{\epsilon} & {} \ar[d, swap, Rightarrow, shorten >=2pt, shorten <= 2pt]{}{\eta} & A \\
B \ar[rr, equals] \ar[ur]{}{g} & {} & B \ar[ur, swap]{}{g} 
\end{tikzcd}
\]
is equal to the $2$-cell:
\[
\begin{tikzcd}
B \ar[r, bend left = 30, ""{name = U}]{}{g} \ar[r, swap, bend right = 30, ""{name = B}]{}{g} & A \ar[from = U, to = B, Rightarrow, shorten >=2pt, shorten <= 2pt]{}{\iota_g}
\end{tikzcd}
\]
We write $f \dashv g$ to show that the pair $(A,B, f, g, \eta, \epsilon)$ is an adjunction with $f$ the left adjoint of $g$ (and $g$ the right adjoint of $f$).
\end{definition}
\begin{lemma}[{\cite[Lemma 2.1.3]{RiehlVertityElementss}}]\label{Lemma: Section 3: 2-functors preserve adjoints}
If $F:\Cfrak \to \Dfrak$ is a strict $2$-functor, then $F$ preserves any adjoints that exist in $\Cfrak$ on the nose.
\end{lemma}
\begin{proof}
Because $F$ is strict, it preserves all equational data on the nose. Thus it sends the triangle identities in $\Cfrak$ on $(A,B,f,g,\eta,\epsilon)$ to the triangle identities in $\Dfrak$ on $(FA,FB,Ff,Fg,F\eta,F\epsilon)$.
\end{proof}
\begin{corollary}\label{Cor: Section 3: adjoints in PreEqSfResl get sent to adjoints in GEqCat}
The $2$-functor
 \[
\mathbf{i}:\fPreEq(\SfResl_G(X)^{\op},\fCat) \to \fGEqCat_{/X}
\]
sends adjunctions between pre-equivariant pseudofunctors to adjunctions between equivariant categories.
\end{corollary}
\begin{proof}
Apply Lemma \ref{Lemma: Section 3: 2-functors preserve adjoints} to the strict $2$-functor $\mathbf{i}$ of Lemma \ref{Lemma: Section 3: The embedding of pre-equivariant functors to equivariant categories is a strict 2-functor}.
\end{proof}

This gives us our first main strategy and source of examples for adjunctions between equivariant categories: Those that come from modifications. However, we would also like to know when we can take equivariant functors $\underline{\alpha}:F_G(X) \to E_G(X)$ and $\underline{\beta}:E_G(X) \to F_G(X)$ have adjoints $\quot{\alpha}{\Gamma} \dashv \quot{\beta}{\Gamma}:F(\quot{X}{\Gamma}) \to E(\quot{X}{\Gamma})$, is there an adjunction $\underline{\alpha} \dashv \underline{\beta}$? That is, is it the case that if we have functors and categories that already satisfy equivariant descent and if each functor determines a fibre-wise adjunction, does the adjunction descend? We prove that this is the case below.

\begin{Theorem}\label{Thm: Section 3: Gamma-wise adjoints lift to equivariant adjoints}
Let $F,E:\SfResl_G(X)^{\op} \to \fCat$ be pre-equivariant pseudofunctors and let $\alpha:F \Rightarrow E$ and $\beta:E \to F$ be pseudonatural transformations. Then if $\quot{\alpha}{\Gamma} \dashv \quot{\beta}{\Gamma}$ in $\fCat$ for all $\Gamma \in \Sf(G)_0$, there is an adjunction
\[
\begin{tikzcd}
F \ar[r, bend left = 30, ""{name = U}]{}{\alpha} & E \ar[l, bend left = 30, ""{name = L}]{}{{\beta}} \ar[from = U, to = L, symbol = \dashv]
\end{tikzcd}
\]
in the $2$-category $\fPreEq(\SfResl_G(X)^{\op},\fCat)$. In particular, this gives rise to an adjunction of categories:
\[
\begin{tikzcd}
F_G(X) \ar[r, bend left = 30, ""{name = U}]{}{\underline{\alpha}} & E_G(X) \ar[l, bend left = 30, ""{name = L}]{}{\underline{\beta}} \ar[from = U, to = L, symbol = \dashv]
\end{tikzcd}
\]
by applying the strict $2$-functor $\mathbf{i}$ to the adjunction $\alpha\dashv\beta$.
\end{Theorem}
\begin{proof}
For all $\Gamma \in \Sf(G)_0$, let $\quot{\eta}{\Gamma}$ be the unit of adjunction
\[
\quot{\eta}{\Gamma}:\id_{F(\quot{X}{\Gamma})} \Rightarrow \quot{\beta}{\Gamma} \circ \quot{\alpha}{\Gamma}
\]
and let $\quot{\epsilon}{\Gamma}$ be the counit of adjuction
\[
\quot{\epsilon}{\Gamma}:\quot{\alpha}{\Gamma} \circ \quot{\beta}{\Gamma} \Rightarrow \id_{E(\quot{X}{\Gamma})}
\]
We begin by showing that the collection $\eta = \lbrace \quot{\eta}{\Gamma} \; | \; \Gamma \in \Sf(G)_0\rbrace$ forms a modification $\eta:\iota_F \to \beta \circ \alpha$; proving that $\epsilon := \lbrace \alpha \circ \beta \to \iota_E \; | \; \Gamma \in \Sf(G)_0 \rbrace$ is a modification is dual and will be omitted.

Let $f \in \Sf(G)_1$ be arbitrary and write $\Dom f = \Gamma$, $\Codom f = \Gamma^{\prime}$. We now must show that the diagram
\[
\xymatrix{
F(\overline{f})\quot{A}{\Gamma^{\prime}} \ar[rr]^-{\quot{\eta_{F(\overline{f})\quot{A}{\Gamma^{\prime}}}}{\Gamma}} & & \big(\quot{\beta}{\Gamma} \circ \quot{\alpha}{\Gamma}\big)F(\overline{f})\quot{A}{\Gamma^{\prime}} \ar[d]^{\quot{(\beta \circ \alpha)_{\quot{A}{\Gamma^{\prime}}}}{f}} \\
F(\overline{f})\quot{A}{\Gamma^{\prime}} \ar@{=}[u] \ar[rr]_-{F(\overline{f})\quot{\eta_{\quot{A}{\Gamma^{\prime}}}}{\Gamma^{\prime}}} & & F(\overline{f})\big((\quot{\beta}{\Gamma^{\prime}} \circ \quot{\alpha}{\Gamma^{\prime}})\quot{A}{\Gamma^{\prime}}\big)
}
\]
commutes for all $\quot{A}{\Gamma^{\prime}}$ in $F(\quot{X}{\Gamma^{\prime}})_0$. Consider the diagram:
\[
\xymatrix{
F(\overline{f}) \ar[drr] \ar[d]_{F(\overline{f})\quot{\eta_{\quot{A}{\Gamma^{\prime}}}}{\Gamma^{\prime}}} \ar[rr]^-{\quot{\eta_{F(\overline{f})\quot{A}{\Gamma^{\prime}}}}{\Gamma}} & &  \big(\quot{\beta}{\Gamma} \circ \quot{\alpha}{\Gamma}\big)F(\overline{f})\quot{A}{\Gamma^{\prime}} \ar@{-->}[d]^{\exists!\big(\quot{\beta}{\Gamma} \circ \quot{\alpha}{\Gamma}\big)(\rho)} \\  
F(\overline{f})\big((\quot{\beta}{\Gamma^{\prime}} \circ \quot{\alpha}{\Gamma^{\prime}})\quot{A}{\Gamma^{\prime}}\big)\ar[rr]_-{\big(\quot{(\beta \circ \alpha)_{\quot{A}{\Gamma^{\prime}}}}{f}\big)^{-1}} & & \big(\quot{\beta}{\Gamma} \circ \quot{\alpha}{\Gamma}\big)F(\overline{f})\quot{A}{\Gamma^{\prime}}
}
\]
By the adjunction $\quot{\alpha}{\Gamma} \dashv \quot{\beta}{\Gamma}$ there exists a unique map $\rho:F(\overline{f})\quot{A}{\Gamma^{\prime}} \to F(\overline{f})\quot{A}{\Gamma^{\prime}}$ making the above diagram commute. It is routine to check using the pseudonaturality of $\alpha$ and $\beta$, as well as the fact that units of adjunction are given uniquely up to unique isomorphism, that the morphism $\rho$ must be the identity on $F(\overline{f})\quot{A}{\Gamma^{\prime}}$, giving that the diagram
\[
\xymatrix{
	F(\overline{f}) \ar[drr] \ar[d]_{F(\overline{f})\quot{\eta_{\quot{A}{\Gamma^{\prime}}}}{\Gamma^{\prime}}} \ar[rr]^-{\quot{\eta_{F(\overline{f})\quot{A}{\Gamma^{\prime}}}}{\Gamma}} & &  \big(\quot{\beta}{\Gamma} \circ \quot{\alpha}{\Gamma}\big)F(\overline{f})\quot{A}{\Gamma^{\prime}} \ar@{=}[d] \\  
	F(\overline{f})\big((\quot{\beta}{\Gamma^{\prime}} \circ \quot{\alpha}{\Gamma^{\prime}})\quot{A}{\Gamma^{\prime}}\big)\ar[rr]_-{\big(\quot{(\beta \circ \alpha)_{\quot{A}{\Gamma^{\prime}}}}{f}\big)^{-1}} & & \big(\quot{\beta}{\Gamma} \circ \quot{\alpha}{\Gamma}\big)F(\overline{f})\quot{A}{\Gamma^{\prime}}
}
\]
commutes. By post-composing the induced equation by $\quot{(\beta \circ \alpha)}{f}$ at $F(\overline{f})\quot{A}{\Gamma^{\prime}}$ we get 
\begin{align*}
\quot{(\beta \circ \alpha)_{\quot{A}{\Gamma^{\prime}}}}{f} \circ \quot{\eta_{F(\overline{f})\quot{A}{\Gamma}}}{\Gamma} &= \quot{(\beta \circ \alpha)_{\quot{A}{\Gamma^{\prime}}}}{f} \circ \left(\quot{(\beta \circ \alpha)_{\quot{A}{\Gamma^{\prime}}}}{f}\right)^{-1} \circ F(\overline{f})\quot{\eta_{\quot{A}{\Gamma^{\prime}}}}{\Gamma^{\prime}} \\
&= F(\overline{f})\quot{\eta_{\quot{A}{\Gamma^{\prime}}}}{\Gamma^{\prime}},
\end{align*}
which shows that the diagram
\[
\xymatrix{
	F(\overline{f})\quot{A}{\Gamma^{\prime}} \ar[rr]^-{\quot{\eta_{F(\overline{f})\quot{A}{\Gamma^{\prime}}}}{\Gamma^{\prime}}} & & \big(\quot{\beta}{\Gamma} \circ \quot{\alpha}{\Gamma}\big)F(\overline{f})\quot{A}{\Gamma^{\prime}} \ar[d]^{\quot{(\beta \circ \alpha)_{\quot{A}{\Gamma^{\prime}}}}{f}} \\
	F(\overline{f})\quot{A}{\Gamma^{\prime}} \ar@{=}[u] \ar[rr]_-{F(\overline{f})\quot{\eta_{\quot{A}{\Gamma^{\prime}}}}{\Gamma^{\prime}}} & & F(\overline{f})\big((\quot{\beta}{\Gamma^{\prime}} \circ \quot{\alpha}{\Gamma^{\prime}})\quot{A}{\Gamma^{\prime}}\big)
}
\]
commutes, showing that $\eta$ is a modification. That $\epsilon$ is a modification is given dually.

For the triangle identities, consider the pasting diagram in the $2$-category \\ $\fPreEq(\SfResl_G(X)^{\op},\fCat)$:
\[
\begin{tikzcd}
& E \ar[rr, equals] \ar[dr]{}{\beta} & {} & B \\
F \ar[ur]{}{\alpha} \ar[rr, equals] & {} \ar[u, Rightarrow, shorten >=5pt, shorten <= 5pt]{}{\eta} & A \ar[ur, swap]{}{\alpha} \ar[u, swap, Rightarrow, shorten >=2pt, shorten <= 2pt]{}{\epsilon}
\end{tikzcd}
\]
To prove that this glues (or untwists, if you prefer) to the correct diagram, note that for all $\Gamma \in \Sf(G)_0$ there is an equivalence of pasting diagrams
\[
\begin{tikzcd}
& E(\quot{X}{\Gamma}) \ar[rr, equals] \ar[ddr]{}{\quot{\beta}{\Gamma}} & {} & E(\quot{X}{\Gamma}) & & E(\quot{X}{\Gamma})\\ 
{}& {}& {}& {} & = \\
F(\quot{X}{\Gamma}) \ar[uur]{}{\quot{\alpha}{\Gamma}} \ar[rr, equals] & {} \ar[uu, Rightarrow, shorten >=5pt, shorten <= 5pt]{}{\quot{\eta}{\Gamma}} & A \ar[uur, swap]{}{\quot{\alpha}{\Gamma}} \ar[uu, swap, Rightarrow, shorten >=2pt, shorten <= 2pt]{}{\quot{\epsilon}{\Gamma}} & & & F(\quot{X}{\Gamma}) \ar[uu, bend left = 30, ""{name = L}]{}{\quot{\alpha}{\Gamma}} \ar[uu, bend right = 30, swap, ""{name = R}]{}{\quot{\alpha}{\Gamma}} \ar[Rightarrow, from = L, to = R, shorten <= 2pt, shorten >=2pt]{}{\iota_{\quot{\alpha}{\Gamma}}}
\end{tikzcd}
\]
which in turn shows that the pasting diagram in $\fPreEq(\SfResl_G(X)^{\op},\fCat)$ is equal to the $2$-cell:
\[
\begin{tikzcd}
F \ar[r, bend left = 30, ""{name = U}]{}{\alpha} \ar[swap, bend right = 30, r, ""{name = L}]{}{\alpha} & E \ar[Rightarrow, from = U, to = L, shorten <= 2pt, shorten >= 2pt]{}{\iota_{\alpha}}
\end{tikzcd}
\]
The other triangle identity is given dually and omitted. Thus we get that $(F,E,\alpha,\beta,\eta,\epsilon)$ describes an adjunction in $\fPreEq(\SfResl_G(X)^{\op},\fCat)$, and hence by Corollary \ref{Cor: Section 3: adjoints in PreEqSfResl get sent to adjoints in GEqCat} in $\fGEqCat_{/X}$ as well.
\end{proof}
\begin{corollary}\label{Cor: Section 3: Fibre-wise equivalences are equivalences of equivariant categories}
Let $\underline{\alpha} \dashv \underline{\beta}:F_G(X) \dashv E_G(X)$ be equivariant functors such that the adjoints are induced as in Theorem \ref{Thm: Section 3: Gamma-wise adjoints lift to equivariant adjoints}. Then if the adjoints $\quot{\alpha}{\Gamma} \dashv \quot{\beta}{\Gamma}$ are all equivalences of categories, so is the adjoint $\underline{\alpha} \dashv \underline{\beta}$.
\end{corollary}
\begin{proof}
By Theorem \ref{Thm: Section 3: Gamma-wise adjoints lift to equivariant adjoints}, we have that the adjoint equivalences $\quot{\alpha}{\Gamma} \dashv \quot{\beta}{\Gamma}$ lift/descend to and adjunction $\underline{\alpha} \dashv \underline{\beta}:F_G(X) \to E_G(X)$. From the $\Gamma$-wise natural isomorphisms $\quot{\alpha}{\Gamma} \circ \quot{\beta}{\Gamma} \cong \id_{E(\quot{X}{\Gamma})}$ and $\quot{\beta}{\Gamma} \circ \quot{\alpha}{\Gamma} \cong \id_{F(\quot{X}{\Gamma})}$ induced by the units and counits of adjunction, we get isomorphisms $\underline{\beta} \circ \underline{\alpha} \cong \id_{F_G(X)}$ and $\underline{\alpha} \circ \underline{\beta} \cong \id_{E_G(X)}$ induced by the unit and counit of adjunction.
\end{proof}
\begin{definition}\index{Equivariant Adjunction}\index{Adjunctions!Between Equivariant Categories|see {Equivariant Adjunction}}
An equivariant adjunction is an adjunction in $\fCat$ between equivariant categories of the form
\[
\begin{tikzcd}
F_G(X) \ar[r, bend left = 30, ""{name = 30}]{}{\underline{\alpha}} & E_G(X) \ar[l, bend left = 30, ""{name = L}]{}{\underline{\beta}} \ar[from = U, to = L, symbol = \dashv]
\end{tikzcd}
\]
where $\alpha:F \Rightarrow E$ and $\beta:E \Rightarrow F$ are pseudonatural transformations for which $\quot{\alpha}{\Gamma} \dashv \quot{\beta}{\Gamma}$ for all $\Gamma \in \Sf(G)_0$. If, in addition, each of $\quot{\alpha}{\Gamma}$ and $\quot{\beta}{\Gamma}$ are inverse equivalences, we say that $\alpha$ and $\beta$ are equivariant equivalences.\index{Equivariant Equivalence}
\end{definition}

Let us now proceed to examine the fact that Theorem \ref{Thm: Section 3: Gamma-wise adjoints lift to equivariant adjoints} allows us a clean way to discuss equivariant monoidal closed categories, equivariant toposes, and other such things because it allows us to deduce the existence of adjoints based on the existence of $\Gamma$-local adjoints. In particular, this tells us that under mild conditions the category of equivriant objects in a pseudofunctor fibred in toposes is again a topos which models the $G$-equivariant objects of the topos $F(X)$.

\begin{definition}
Let $F:(\Cscr, \otimes, I) \to (\Dscr, \boxtimes, J)$ be a functor between monoidal categories. We say that 
\end{definition}
\begin{lemma}\label{Lemma: Section 3: Equivariant Cat is monoidal closed}
	Let $F:\SfResl_G(X)^{\op} \to \fCat$ be a monodial pre-equivariant pseudofunctor for which each functor $F(\overline{f})$ is monoidal closed. Then if each category $F(\XGamma)$ is monoidal closed or symmetric monoidal closed, so is $F_G(X)$.
\end{lemma}
\begin{proof}
	By Theorem \ref{Theorem: Section 2: Monoidal preequivariant pseudofunctor gives monoidal equivariant cat} and Proposition \ref{Prop: Section 2: Equivariant cat is symmetric monoidal} it follows that in either case we have that $F_G(X)$ is monoidal or symmetric monoidal, repsectively. We now prove this lemma by providing the case in which $F_G(X)$ is only know to be monoidal; the symmetric monoidal case follows mutatis mutandis. From Theorem \ref{Thm: Section 3: Gamma-wise adjoints lift to equivariant adjoints} and the fact that we know there are adjoints
	\[
	[\quot{A}{\Gamma},-]_{\Gamma} \dashv \quot{A}{\Gamma} \os{\Gamma}{\otimes} (-)
	\]
	for all $\Gamma \in \Sf(G)_0$, it suffices to prove that for any object $A \in F_G(X)_0$, the functors $A \otimes (-):F_G(X) \to F_G(X)$ and $[A,-]:F_G(X) \to F_G(X)$ arise from pseudonatural transformations as in Theorem \ref{Thm: Section 3: Psuedonatural trans lift to equivariant functors}. To this end let
	\[
	\theta_f^{A,B}:F(\overline{f})\left(\quot{A}{\Gamma^{\prime}} \os{\Gamma^{\prime}}{\otimes} \quot{B}{\Gamma^{\prime}}\right) \xrightarrow{\cong} F(\overline{f})(\quot{A}{\Gamma^{\prime}}) \os{\Gamma}{\otimes} F(\overline{f})(\quot{B}{\Gamma^{\prime}})
	\]
	and
	\[
	\xi_{f}^{A,B}:F(\overline{f})[\quot{A}{\Gamma^{\prime}},\quot{B}{\Gamma^{\prime}}]_{\Gamma^{\prime}} \xrightarrow{\cong} [F(\overline{f})(\quot{A}{\Gamma^{\prime}}),F(\overline{f})(\quot{B}{\Gamma^{\prime}})]_{\Gamma}
	\]
	be the preservation isomorphisms. 
	
	Fix an object $A \in F_G(X)_0$. We then define pseudonatrual transformation $A \otimes (-):F \to F$ as follows: For all $\Gamma \in \Sf(G)_0$, define the functor
	\[
	\quot{(A \otimes (-))}{\Gamma}:F(\quot{X}{\Gamma}) \to F(\quot{X}{\Gamma})
	\]
	by
	\[
	\quot{(A \otimes (-))}{\Gamma} := \quot{A}{\Gamma} \os{\Gamma}{\otimes} (-).
	\]
	Fix a morphism $f \in \Sf(G)_1$ and write $\Dom f = \Gamma$ and $\Codom f = \Gamma^{\prime}$. Now observe that on one hand we have that for all $\quot{B}{\Gamma^{\prime}} \in F(\quot{X}{\Gamma^{\prime}})_0$,
	\[
	\big(\quot{(A \otimes (-))}{\Gamma} \circ F(\overline{f})\big)(\quot{B}{\Gamma^{\prime}}) = \quot{A}{\Gamma} \os{\Gamma}{\otimes} F(\overline{f})(\quot{B}{\Gamma^{\prime}})
	\]
	while on the other hand we have
	\[
	\big(F(\overline{f}) \circ \quot{(A \otimes (-))}{\Gamma^{\prime}}\big)(\quot{B}{\Gamma^{\prime}}) = F(\overline{f})\left(\quot{A}{\Gamma^{\prime}} \os{\Gamma^{\prime}}{\otimes} \quot{B}{\Gamma^{\prime}}\right).
	\]
	There is then an isomorphism
	\[
	\xymatrix{
	\quot{A}{\Gamma} \os{\Gamma}{\otimes} F(\overline{f})(\quot{B}{\Gamma^{\prime}}) \ar[rrr]^-{\left(\tau_f^A\right)^{-1}\os{\Gamma}{\otimes} \id_{F(\overline{f})\quot{B}{\Gamma^{\prime}}}} \ar[drrr] & & & F(\overline{f})(\AGammap) \os{\Gamma}{\otimes} F(\overline{f})(\quot{B}{\Gamma^{\prime}}) \ar[d]^{\left(\theta_f^{A,B}\right)^{-1}} \\
&	 & & F(\overline{f}) \big(\AGammap \os{\Gamma^{\prime}}{\otimes} \BGammap\big)
	}
	\]
	which is natural in $\quot{B}{\Gamma^{\prime}}$ by the bifunctoriality of $\otimes$ and the natruality of $\theta_f$ (and hence its inverse). This determines our natural transformation $\quot{(A \otimes (-))}{f}$ component of the pseudonatural transformation $A \otimes (-)$; its object assignment is given by
	\[
	\quot{(A \otimes (-))_{\quot{B}{\Gamma^{\prime}}}}{f} := \left(\theta_f^{A,B}\right)^{-1} \circ \left(\tau_f^A\right)^{-1} \os{\Gamma}{\otimes} \id_{F(\overline{f})(\quot{B}{\Gamma^{\prime}})}.
	\]
	A routine check shows that the collection 
	\[
	A \otimes (-) = \lbrace \quot{(A \otimes (-))}{\Gamma}, \quot{(A \otimes(-))}{f} \; | \; \Gamma \in \Sf(G)_0, f \in \Sf(G)_1\rbrace 
	\]
	determines a pseudonatrual transformation which we identify (through Lemma \ref{Lemma: Section 3: The embedding of pre-equivariant functors to equivariant categories is a strict 2-functor}) with the functor $A \otimes (-)$.
	
	The fact that $[A,-]$ gives rise to a pseudonatural transformation is dual to the argument above. More explicitly, the object functors $\quot{[A,-]}{\Gamma}$ are given by
	\[
	\quot{[A,-]}{\Gamma}(\quot{B}{\Gamma}) := [\quot{A}{\Gamma}, \quot{B}{\Gamma}]_{\Gamma}
	\]
	and the natural isomorphism $\quot{[A,-]}{f}:\quot{[A,-]}{\Gamma} \circ F(\overline{f}) \Rightarrow F(\overline{f}) \circ \quot{[A,-]}{\Gamma^{\prime}}$ are given by
	\[
	\quot{[A,-]_{-}}{f} := \left(\xi_f^{A,-}\right)^{-1} \circ [\tau_f^A,\id_{F(\overline{f})(-)}].
	\]
	As before, we identify the pseudonatural transformation $[A,-]$ with the functor $[A,-]:F_G(X) \to F_G(X)$. Now appealing to Theorem \ref{Thm: Section 3: Gamma-wise adjoints lift to equivariant adjoints} gives the adjunction $A \otimes (-) \dashv [A,-]$, as desired. 
\end{proof}
\begin{corollary}\label{Cor: Section 3: Equivariant Cartesian Closed Cats}
Let $F:\SfResl_G(X)^{\op} \to \fCat$ be a pre-equivariant pseudofunctor such that each category $F(\quot{X}{\Gamma})$ has finite products and each functor $F(\overline{f})$ preserves these products. Assume further that each category $F(\quot{X}{\Gamma})$ is Cartesian closed and that each fibre functor $F(\overline{f})$ is a Cartesian closed functor. Then $F_G(X)$ is Cartesian closed as well.
\end{corollary}
\begin{proof}
This is simply Lemma \ref{Lemma: Section 3: Equivariant Cat is monoidal closed} applied to a Cartesian monoidal closed structure.
\end{proof}

This leads us to a class of equivariant toposes that arise through equivariant descent valued in pre-equivariant pseudofunctors fibred over toposes.
\begin{proposition}\label{Proposition: Section 3: Equivariant Toposes}
Let $F:\SfResl_G(X)^{\op} \to \fCat$ be a pre-equivariant pseudofunctor such that each category $F(\quot{X}{\Gamma})$ is a topos and each morphism $F(\overline{f})$ is the pullback of a geometric morphism which preserves the subobject classifier and internal hom functors. Then $F_G(X)$ is a topos.
\end{proposition}
\begin{proof}
Recall that a topos is a category $\Cscr$ which has finite limits, is Cartesian closed, and has a subobject classifier (cf.\@ \cite{MacLaneMoerdijk}). By Proposition \ref{Prop: Section 2: Equivariant Cat has Subobject Classifiers}, $F_G(X)$ has a subobject classifier and by Corollary \ref{Cor: Equivariant cat is finitely complete or cocomplete} $F_G(X)$ is finitely complete. Finally by Corollary \ref{Cor: Section 3: Equivariant Cartesian Closed Cats} we have that $F_G(X)$ is Cartesian closed. Therefore $F_G(X)$ is a topos.
\end{proof}
\begin{remark}
An example of such a topos is one whose fibre functors $F(\overline{f})$ are all the inverse image functors of an atomic geometric morphism. In particular, if each functor $F(\overline{f}) = f^{\ast}$ is an inverse image of an {\'e}tale geometric morphism $f:F(\quot{X}{\Gamma}) \to F(\quot{X}{\Gamma^{\prime}})$, then $F_G(X)$ is a topos.
\end{remark}
\begin{example}
Consider the pre-equivariant pseudofunctor \\$F:\SfResl_G(X)^{\op} \to \fCat$ given by setting $F(\quot{X}{\Gamma}) := \Shv(X,\text{{\'e}t})$, the (small) topos of {\'e}tale sheaves of sets on $X$ and setting $F(\overline{f}) := \id_{\Shv(X,\text{{\'e}t})}$. Then each functor $F(\overline{f})$ is an inverse image of an {\'e}tale geometric morphism and $F_G(X)$ is a topos. This realizes the topos of {\'e}tale sheaves of sets on $X$ as an equivariant topos with respect to the trivial action of $G$ on $X$.
\end{example}

We now close this section with a topic of interest for equivariant homotopy theory (and in particular for a consideration of equivariant homological algebra): Equivariant localization. Our motivation for this is to see in which sense the equivariant derived category $D^b_G(X)$ arises as a localization of $\Shv^b_G(X)$ at some suitable class of quasi-isomorphisms. We will give some moderate technical conditions that show when we can essentially localize the fibre categories of a pre-equivariant pseudofunctor over the equivariant data specified by $\SfResl_G(X)$. These technical limitations should not prove too much of a hindrance in practice.

Begin by letting $F:\SfResl_G(X)^{\op} \to \fCat$ be a pre-equivariant pseudofunctor and, for all $\Gamma \in \Sf(G)_0$ let $\quot{S}{\Gamma} \subseteq F_G(X)_1$ (as a $\Vscr$-class, of course). Our goal now is to study when we can descend the data of the localizations $\quot{S}{\Gamma}^{-1}F(\XGamma)$ through the fibre functors $F(\overline{f})$, i.e., when there are functors $S^{-1}F(\overline{f}):\quot{S}{\Gamma^{\prime}}^{-1}F(\XGammap) \to \quot{S}{\Gamma}^{-1}F(\XGamma)$ which allow us to define a suitable localization $\lambda_S:F_G(X) \to S^{-1}F_G(X)$ (which we will call an equivariant localization). The evident condition we use to produce the pseudofunctor $S^{-1}F$ is that we need $F(\of)$ to carry every $\quot{S}{\Gamma^{\prime}}$-morphism to a $\quot{S}{\Gamma}$-morphism; it turns out this is the only structure constraint we need to derive the following proposition, which shows that taking localizations under this constraint gives rise to a pre-equivariant pseudofunctor $S^{-1}F$ whose fibre categories are all the localizations of the fibre categories of $F$. Subsequently, we show that there is a pseudonatural transformation $\lambda:F \implies S^{-1}F$ which realizes the category $S^{-1}F_G(X)$ as an equivariant localization. Finally, after discussing equivariant localization and the property it has (cf.\@ Propositions \ref{Prop: Section 3: Equivariant Localization Pseudofunctor}, \ref{Prop: Section 3: Equivariant localization pseudonat} and Theorem \ref{Theorem: Section 3: Equivariant localization is a 2-colimit}), we will compare this with the na{\"i}ve localization. In particular, we show that the equivariant localization is not a literal localization (cf.\@ Example \ref{Example: Section 3: Equivariant localizaiton is not the localization of equivariant category}) but does have an analogous universal property when restricted to equivariant categories and functors.

\begin{proposition}\label{Prop: Section 3: Equivariant Localization Pseudofunctor}
Let $F:\SfResl_G(X)^{\op} \to \fCat$ be a pre-equivariant pseudofunctor and for all $\Gamma \in \Sf(G)_0$ assume that there is a subclass $\quot{S}{\Gamma} \subseteq F(\XGamma)_1$ such that for all $f \in \Sf(G)_1$, $F(\overline{f})(\quot{S}{\Gamma^{\prime}}) \subseteq \quot{S}{\Gamma}$ where $f:\Gamma \to \Gamma^{\prime}$. Then if $S \subseteq F_G(X)_1$ is the subclass of morphisms
\[
S := \lbrace P \in F_G(X)_1 \; | \; \forall\,\Gamma \in \Sf(G)_0.\,\forall\,\quot{\rho}{\Gamma} \in P.\, \quot{\rho}{\Gamma} \in \quot{S}{\Gamma} \rbrace
\]
there is a pre-equivariant pseudofunctor $S^{-1}F:\SfResl_G(X)^{\op} \to \fCat$\index[notation]{SinverseF@$S^{-1}F$} which sends objects $\Gamma \times X$ to $(\quot{S}{\Gamma})^{-1}F(\XGamma)$ and sends morphisms $f \times \id_X:\Gamma \times X \to \Gamma^{\prime} \times X$ to the unique functor $(S^{-1}F)(\overline{f})$ fitting in the commuting square:
\[
\xymatrix{
F(\XGammap) \ar[rr]^{F(\overline{f})} \ar[d]_{\lambda_{\quot{S}{\Gamma^{\prime}}}} & & F(\XGamma) \ar[d]^{\lambda_{\quot{S}{\Gamma}}} \\
\quot{S}{\Gamma^{\prime}}^{-1}F(\XGammap) \ar@{-->}[rr]_{\exists!S^{-1}F(\overline{f})} & & \quot{S}{\Gamma}^{-1}F(\XGamma)
}
\]
\end{proposition}
\begin{proof}
We first note that the categories $(\quot{S}{\Gamma})^{-1}F(\XGamma)$ all exist by \cite[Section 1.1]{GabrielZisman}. The functor $(S^{-1}F)(\of)$ exists by first observing that for any $\varphi \in \quot{S}{\Gamma^{\prime}}$, since $F(\of)(\varphi) \in \quot{S}{\Gamma}$, $(\lambda_{\quot{S}{\Gamma}} \circ F(\of))(\varphi)$ is an isomorphism in $(\quot{S}{\Gamma})^{-1}F(\XGamma)$; applying the universal property of the localization to the diagram
\[
\xymatrix{
F(\XGammap) \ar[rr]^-{\lambda_{\quot{S}{\Gamma}} \circ F(\of)} \ar[dr]_{\lambda_{\quot{S}{\Gamma^{\prime}}}} & & (\quot{S}{\Gamma})^{-1}F(\XGamma) \\
 & (\quot{S}{\Gamma^{\prime}})^{-1}F(\XGammap) \ar@{-->}[ur]_{\exists!(SF)^{-1}(\of)}
}
\]
gives the existence of $(S^{-1}F)(\of)$.

Note that since the diagram
\[
\xymatrix{
F(\XGamma) \ar[d] \ar@{=}[r] & F(\XGamma) \ar[d] \\
\quot{S}{\Gamma}^{-1}F(\XGamma) \ar@{=}[r] & \quot{S}{\Gamma}^{-1}F(\XGamma)
}
\]
commutes, we have that $S^{-1}F(\overline{\id}_{\Gamma}) = \id_{\quot{S}{\Gamma}^{-1}F(\XGamma)}$. Thus we only need to establish the existence of the compositor natural isomorphisms, for $\Gamma \xrightarrow{f} \Gamma^{\prime} \xrightarrow{g} \Gamma^{\prime\prime}$ in $\Sf(G)$,
\[
\phi_{f,g}:S^{-1}F(\overline{f}) \circ S^{-1}F(\overline{g}) \Rightarrow S^{-1}F(\overline{g} \circ \overline{f})
\]
and their identities to prove that $S^{-1}F$ eixsts. For this fix $f$ and $g$ as above and consider that because $F$ is a pseudofunctor we have natural transformations $\phi_{f,g}$ fitting into a pasting diagram:
\[
\begin{tikzcd}
 & F(\XGammap) \ar[d, Rightarrow, shorten >= 4pt, shorten <= 4pt]{}{\phi_{f,g}} \ar[dr]{}{F(\overline{f})} \\
F(\XGammapp) \ar[ur]{}{F(\overline{g})} \ar[rr, swap]{}{F(\overline{g} \circ \overline{f})} & {} & F(\XGamma)
\end{tikzcd}
\]
We now recall (cf.\@ \cite{BorceuxCatAlg1}, \cite{GabrielZisman}, or \cite{MacLaneCWM}) that from the construction of the categories $(\quot{S}{\Gamma})^{-1}F(\XGamma)$ we have that for all $Y \in \Sf(G)_0$,
\[
(\quot{S}{Y}^{-1}F(\quot{X}{Y})_0 = F(\quot{X}{Y})_0.
\]
In particular, it follows from the universal property of the localizations and the constructions of the $\lambda$ functors that the object assignments of each functor $(S^{-1}F)(\of)$, $(S^{-1}F)(\overline{g})$, and $(S^{-1}F)(\overline{g} \circ \of)$ must be the same as the object assignments of $F(\of)$, $F(\overline{g})$, and $F(\overline{g} \circ \of)$, respectively. Then for any $A \in (\quot{S}{\Gamma^{\prime\prime}})^{-1}F(\XGammapp)_0$, we have
\[
\big(S^{-1}F(\overline{f}) \circ S^{-1}F(\overline{g})\big)A = \big(F(\overline{f}) \circ F(\overline{g})\big)A
\] 
(and similarly for $F(\overline{g} \circ \overline{f})$). Thus any natural transformation 
\[
\alpha:(S^{-1}F)(\overline{f}) \circ (S^{-1}F)(\overline{g}) \Rightarrow (S^{-1}F)(\overline{g} \circ \overline{f})
\]
must be a collection of $(\quot{S}{\Gamma})^{-1}F(\XGamma)$-morphisms
\[
\alpha = \lbrace \alpha_A:(F(\overline{f}) \circ F(\overline{g}))A \to F(\overline{g} \circ \overline{f})A \; | \; A \in (\quot{S}{\Gamma^{\prime\prime}})^{-1}F(\XGammapp) \rbrace.
\]

We now prove that $\phi_{f,g}$ extends to a compositor natural isomorphism 
\[
S^{-1}\phi_{f,g}:(S^{-1}F)(\overline{f}) \circ (S^{-1}F)(\overline{g}) \Rightarrow (S^{-1}F)(\overline{g} \circ \overline{f})
\] 
where
\[
S^{-1}\phi_{f,g}:= \lbrace \phi_{f,g}^{A}:(F(\overline{f}) \circ F(\overline{g}))A \to F(\overline{g} \circ \overline{f})A\; | \; A \in (\quot{S}{\Gamma^{\prime\prime}})^{-1}F(\XGammapp)_0 \rbrace.
\]
That this is an isomorphism for each $A$ is trivial from the fact that $\phi_{f,g}$ is; as such, we only need to prove the naturality holds, i.e., that for any zig-zag describing a morphism $\alpha:A \to B$ in $\quot{S}{\Gamma^{\prime\prime}}^{-1}F(\XGammapp)$, the diagram
\[
\xymatrix{
\big(S^{-1}F(\overline{f}) \circ S^{-1}F(\overline{g})\big)A \ar[rr]^-{\phi_{f,g}^{A}} \ar[d]_{\big(S^{-1}F(\overline{f}) \circ S^{-1}F(\overline{g})\big)\alpha} & & S^{-1}F(\overline{g} \circ \overline{f})A \ar[d]^{S^{-1}F(\overline{g} \circ \overline{f})\alpha} \\
\big(S^{-1}F(\overline{f}) \circ S^{-1}F(\overline{g})\big)B \ar[rr]_-{\phi_{f,g}^{B}} & & S^{-1}F(\overline{g} \circ \overline{f})B
}
\]
commutes.

To prove naturality of the square above, let $\alpha$ be represented by some zig-zag of the form
\[
\xymatrix{
A = A_0  & & A_1 &\cdots & A_n  & & A_{n+1} = B\\
 & C_0 \ar[ur]_{\rho_0} \ar[ul]^{\psi_0} &  & \cdots \ar[ur] \ar[ul] & & C_n \ar[ur]_{\rho_n} \ar[ul]^{\psi_n}
}
\]
where $\rho_i \in F(\XGammapp)_1$ and where $\psi_i \in \quot{S}{\Gamma^{\prime\prime}}$ for all $ 0 \leq i \leq n$ and for $n \in \N$. Note that not only is the argument we will present independent of the representation of $\alpha$, but proving the other cases (where $\alpha$ starts with a $\rho_i$, where $\alpha$ ends with a $\psi_i^{-1}$, etc.) follows mutatis mutandis and will be omitted.

Begin by observing that by construction, since
\[
\alpha = \prod_{i=0}^{n} \rho_i \circ \psi_i^{-1}
\]
from the universal properties defining all the objects we have that
\begin{align*}
&\big(S^{-1}F(\overline{f}) \circ S^{-1}F(\overline{g})\big)\alpha \\
&= \big(S^{-1}F(\overline{f}) \circ S^{-1}F(\overline{g})\big)\left(\prod_{i=0}^{n} \rho_i \circ \psi_i^{-1}\right) \\
& = \prod_{i=0}^{n} \big(S^{-1}F(\overline{f}) \circ S^{-1}F(\overline{g})\big)\rho_i \circ \big(S^{-1}F(\overline{f}) \circ S^{-1}F(\overline{g})\big)\psi_i^{-1} \\
&= \prod_{i=0}^{n} \big(F(\overline{f}) \circ F(\overline{g})\big)\rho_i \circ \big(F(\overline{f}) \circ F(\overline{g})\big)(\psi_i)^{-1}
\end{align*}
and similarly
\[
S^{-1}F(\overline{g} \circ \overline{f})\alpha = \prod_{i=0}^{n} F(\overline{g} \circ \overline{f})\rho_i \circ F(\overline{g} \circ \overline{f})(\psi_i)^{-1}.
\]
Now note that from the naturality of $\phi_{f,g}$ the diagrams
\[
\xymatrix{
\big(F(\overline{f}) \circ F(\overline{g})\big)C_i \ar[rr]^-{\phi_{f,g}^{C_i}} \ar[d]_{\big(F(\overline{f}) \circ F(\overline{g})\big)\rho_i} & & F(\overline{g} \circ \overline{f})C_i \ar[d]^{F(\overline{g} \circ \overline{f})\rho_i} \\
\big(F(\overline{f}) \circ F(\overline{g})\big)A_{i+1} \ar[rr]_-{\phi_{f,g}^{A_{i+1}}} & & F(\overline{g} \circ \overline{f})A_{i+1}
}
\]
commute in both $F(\XGamma)$ and $\quot{S}{\Gamma}^{-1}F(\XGamma)$ for all $0 \leq i \leq n$. Similarly, the diagrams
\[
\xymatrix{
	\big(F(\overline{f}) \circ F(\overline{g})\big)C_i \ar[rr]^-{\phi_{f,g}^{C_i}} \ar[d]_{\big(F(\overline{f}) \circ F(\overline{g})\big)\psi_i} & & F(\overline{g} \circ \overline{f})C_i \ar[d]^{F(\overline{g} \circ \overline{f})\psi_i} \\
	\big(F(\overline{f}) \circ F(\overline{g})\big)A_{i} \ar[rr]_-{\phi_{f,g}^{A_{i}}} & & F(\overline{g} \circ \overline{f})A_{i}
}
\]
commute in both $F(\XGamma)$ and $\quot{S}{\Gamma}^{-1}F(\XGamma)$ for all $0 \leq i \leq n$; however, upon inverting both vertical morphisms, we find that the diagrams
\[
\xymatrix{
	\big(F(\overline{f}) \circ F(\overline{g})\big)C_i \ar[rr]^-{\phi_{f,g}^{C_i}} & & F(\overline{g} \circ \overline{f})C_i  \\
	\big(F(\overline{f}) \circ F(\overline{g})\big)A_{i} \ar[rr]_-{\phi_{f,g}^{A_{i}}} \ar[u]^{\big(F(\overline{f}) \circ F(\overline{g})\big)(\psi_i)^{-1}} & & F(\overline{g} \circ \overline{f})A_{i} \ar[u]_{F(\overline{g} \circ \overline{f})(\psi_i)^{-1}}
}
\]
must also commute in $\quot{S}{\Gamma}^{-1}F(\XGamma)$ for all $0 \leq i \leq n$. Now observe that
\begin{align*}
&S^{-1}F(\overline{g} \circ \overline{f})\alpha \circ \phi_{f,g}^{A} \\
&= \left(\prod_{i=0}^{n} F(\overline{g} \circ \overline{f})\rho_i \circ F(\overline{g} \circ \overline{f})(\psi_i)^{-1}\right) \circ \phi_{f,g}^{A} \\
&= \left(\prod_{i=1}^{n} F(\overline{g} \circ \overline{f})\rho_i \circ F(\overline{g} \circ \overline{f})(\psi_i)^{-1}\right) \circ F(\overline{g} \circ \overline{f})\rho_0 \circ F(\overline{g} \circ \overline{f})(\psi_0)^{-1} \circ \phi_{f,g}^{A_0} \\
&= \left(\prod_{i=1}^{n} F(\overline{g} \circ \overline{f})\rho_i \circ F(\overline{g} \circ \overline{f})(\psi_i)^{-1}\right) \circ F(\overline{g} \circ \overline{f})\rho_0 \circ \phi_{f,g}^{C_0} \\
&\circ \big(F(\overline{f}) \circ F(\overline{g})\big)(\psi_0)^{-1} \\
&= \left(\prod_{i=1}^{n} F(\overline{g} \circ \overline{f})\rho_i \circ F(\overline{g} \circ \overline{f})(\psi_i)^{-1}\right) \circ \phi_{f,g}^{A_1}\circ \big(F(\overline{f}) \circ F(\overline{g})\big)\rho_0 \\
&\circ \big(F(\overline{f} \circ \overline{g})\big)(\psi_0)^{-1} \\ 
&= \phi_{f,g}^{A_{n+1}} \circ \prod_{i=0}^{n}\big(F(\overline{f}) \circ F(\overline{g})\big)\rho_i \circ \big(F(\overline{f}) \circ F(\overline{g})\big)(\psi_i)^{-1} \\
&= \phi_{f,g}^{A_{n+1}} \circ \big(F(\overline{f}) \circ F(\overline{g})\big)\alpha = \phi_{f,g}^{B} \circ \big(S^{-1}F(\overline{f}) \circ S^{-1}F(\overline{g})\big)\alpha.
\end{align*}
Thus $S^{-1}\phi_{f,g}$ is a natural isomorphism. That it satisfies all the required relations of a compositor follows in a straightforward but tedious calculation using the fact that $\phi_{f,g}$ is a compositor, proving that $S^{-1}F$ is a pre-equivariant pseudofunctor.
\end{proof}
\begin{proposition}\label{Prop: Section 3: Equivariant localization pseudonat}
Assume $F$ and $S$ are as in the statement of Proposition \ref{Prop: Section 3: Equivariant Localization Pseudofunctor}. There is a pseudonatural transformation $\lambda_S:F \Rightarrow S^{-1}F$\index[notation]{LambdaS@$\lambda_S$} whose object functor $\quot{\lambda}{\Gamma}$ is given by
\[
\quot{\lambda}{\Gamma} := \lambda_{\quot{S}{\Gamma}}:F(\XGamma) \to \quot{S}{\Gamma}^{-1}F(\XGamma).
\]
\end{proposition}
\begin{proof}
Since we have already defined $\quot{\lambda}{\Gamma}$, it suffices to define the natural isomorphisms $\quot{\lambda}{f}$. For this recall that the morphism $S^{-1}F(\overline{f})$ is defined via the universal property of the localization making the triangle
\[
\xymatrix{
F(\XGammap) \ar[rr]^{\quot{\lambda}{\Gamma} \circ F(\overline{f})} \ar[dr]_{\quot{\lambda}{\Gamma^{\prime}}} & & \quot{S}{\Gamma}^{-1}F(\XGamma) \\
 & \quot{S}{\Gamma^{\prime}}^{-1}F(\XGammap) \ar@{-->}[ur]_{\exists!S^{-1}F(\overline{f})}
}
\]
commute. However, this says that the diagram
\[
\xymatrix{
F(\XGammap) \ar[rr]^-{F(\overline{f})} \ar[d]_{\quot{\lambda}{\Gamma^{\prime}}} & & F(\XGamma) \ar[d]^{\quot{\lambda}{\Gamma}} \\
\quot{S}{\Gamma^{\prime}}^{-1}F(\XGammap) \ar[rr]_-{S^{-1}F(\overline{f})} & & \quot{s}{\Gamma}^{-1}F(\XGamma)
}
\] 
commutes on the nose, so we define $\quot{\lambda}{f}$ to simply be the equality natural transformation witnessing the equality $\quot{\lambda}{\Gamma} \circ F(\overline{f}) = S^{-1}F(\overline{f}) \circ \quot{\lambda}{\Gamma^{\prime}}$. The required equality on pasting diagrams then holds immediately based on the fact that the $\quot{\lambda}{f}$ are equalities and from how we defined our compositors in the proof of Proposition \ref{Prop: Section 3: Equivariant Localization Pseudofunctor}.
\end{proof}
\begin{Theorem}\label{Theorem: Section 3: Equivariant localization is a 2-colimit}
Let $F$ and $E$ be pre-equivariant pseudofunctors on $X$ and let 
\[
S = \lbrace \quot{S}{\Gamma} \; | \; \quot{S}{\Gamma} \subseteq F(\XGamma)_1, \Gamma \in \Sf(G)_0 \rbrace
\] 
be a collection of morphism classes as in Proposition \ref{Prop: Section 3: Equivariant Localization Pseudofunctor}. Assume there is a pseudonatural transformation $\alpha:F \to E$ such that for all $\Gamma \in \Sf(G)$ the functors $\quot{\alpha}{\Gamma}$ make 
\[
\quot{\alpha}{\Gamma}(\quot{S}{\Gamma}) \subseteq \Iso(E(\XGamma)) 
\]
hold. Then there exists a unique pseudonatural transformation $\eta:S^{-1}F \to E$ for which the induced diagram
\[
\xymatrix{
F_G(X) \ar[rr]^{\underline{\alpha}} \ar[dr]_{\underline{\lambda}_S} & & E_G(X) \\
 & S^{-1}F_G(X) \ar@{-->}[ur]_{\exists!\,\underline{\eta}}
}
\]
\index[notation]{SinverseFGX@$S^{-1}F_G(X)$}commutes in $\Cat$.
\end{Theorem}
\begin{proof}
We first construct the functor component of the pseudonatural transformation $\eta$. For this fix a $\Gamma \in \Sf(G)_0$ and define the functor $\quot{\eta}{\Gamma}:S^{-1}F(\XGamma) \to E(\XGamma)$ as the unique functor filling the commuting diagram:
\[
\xymatrix{
F(\XGamma) \ar[rr]^{\quot{\alpha}{\Gamma}} \ar[dr]_{\quot{\lambda}{\Gamma}} & & E(\XGamma) \\
 & S^{-1}F(\XGamma) \ar@{-->}[ur]_{\exists!\quot{\eta}{\Gamma}}
}
\]
To define the natural isomorphisms $\quot{\eta}{f}$, fix a morpshism $f \in \Sf(G)_1$ and write $\Dom f = \Gamma, \Codom f = \Gamma^{\prime}$. We observe that by assumption we have a $2$-cell
\[
\begin{tikzcd}
F(\XGammap) \ar[r, bend left = 30, ""{name = U}]{}{\quot{\alpha}{\Gamma} \circ F(\overline{f})} \ar[r, bend right = 30, swap, ""{name = L}]{}{E(\overline{f}) \circ \quot{\alpha}{\Gamma^{\prime}}} & E(\XGamma) \ar[Rightarrow, from = U, to = L, shorten <= 4pt, shorten >= 4pt]{}{\quot{\alpha}{f}}
\end{tikzcd}
\]
On one hand, note that the upper arrow in the $2$-cell factors in the strictly commuting square
\[
\xymatrix{
F(\XGamma) \ar[r]^-{\quot{\lambda}{\Gamma}} & S^{-1}F(\XGamma) \ar[d]^{\quot{\eta}{\Gamma}} \\
F(\XGammap) \ar[u]^{F(\overline{f})} \ar[r]_-{\quot{\alpha}{\Gamma} \circ F(\overline{f})} & E(\XGamma)
}
\]
while on the other hand the bottom arrow in the $2$-cell factors in the strictly commuting square:
\[
\xymatrix{
F(\XGammap) \ar[r]^-{E(\overline{f}) \circ \quot{\lambda}{\Gamma^{\prime}}} \ar[d]_{\quot{\lambda}{\Gamma^{\prime}}} & E(\XGamma) \\
S^{-1}F(\XGammap) \ar[r]_-{\quot{\eta}{\Gamma^{\prime}}} & E(\XGammap) \ar[u]_{E(\overline{f})}
}
\]
Combining these gives a pasting diagram:
\[
\begin{tikzcd}
F(\XGamma) \ar[r]{}{\quot{\lambda}{\Gamma}} & S^{-1}F(\XGamma) \ar[d]{}{\quot{\eta}{\Gamma}^{\prime}} \\
F(\XGammap) \ar[u]{}{F(\overline{f})} \ar[r, bend left = 30, ""{name = U}]{}{\quot{\alpha}{\Gamma} \circ F(\overline{f})} \ar[r, swap, bend right = 30, ""{name = L}]{}{E(\overline{f}) \circ \quot{\lambda}{\Gamma^{\prime}}} \ar[d, swap]{}{\quot{\lambda}{\Gamma^{\prime}}} & E(\XGamma) \\
S^{-1}F(\XGammap) \ar[r, swap]{}{\quot{\eta}{\Gamma^{\prime}}} & E(\XGammap) \ar[u, swap]{}{E(\overline{f})} \ar[from = U, to = L, Rightarrow, shorten <= 4pt, shorten >= 4pt]{}{\quot{\alpha}{f}}
\end{tikzcd}
\]
It is routine to check that since the square
\[
\xymatrix{
F(\XGammap) \ar[r]^-{\quot{\lambda}{\Gamma}} \ar[d]_{F(\overline{f})} & S^{-1}F(\XGammap) \ar[d]^{S^{-1}F(\overline{f})} \\
F(\XGamma) \ar[r]_-{\quot{\lambda}{\Gamma}} & S^{-1}F(\XGamma)
}
\]
commutes on the nose and the categories $S^{-1}F(\XGamma)$ and $S^{-1}F(\XGammap)$ arise as localizations, there is a unique (from the fact that the localization arises as a $(2,1)$-colimit) natural isomorphism 
\[
\quot{\eta}{f}:\quot{\eta}{\Gamma} \circ S^{-1}F(\overline{f}) \to E(\overline{f}) \circ \quot{\eta}{\Gamma^{\prime}}
\]
which fits into and factorizes the above pasting diagram as:
\[
\begin{tikzcd}
F(\XGamma) \ar[r]{}{\quot{\lambda}{\Gamma}} & S^{-1}F(\XGamma) \ar[d]{}{\quot{\eta}{\Gamma}^{\prime}} \\
F(\XGammap) \ar[u]{}{F(\overline{f})}  \ar[d, swap]{}{\quot{\lambda}{\Gamma^{\prime}}} & E(\XGamma) \\
S^{-1}F(\XGammap) \ar[uur, ""{name = M}]{}[description]{S^{-1}F(\overline{f})} \ar[r, swap]{}{\quot{\eta}{\Gamma^{\prime}}} & E(\XGammap) \ar[u, swap]{}{E(\overline{f})} \ar[from = M, to = 3-2, Rightarrow, swap, shorten <= 12pt, shorten >= 4pt]{}{\quot{\eta}{f}}
\end{tikzcd}
\]
However, this is exactly the $2$-cell
\[
\begin{tikzcd}
S^{-1}F(\XGammap) \ar[r, bend left = 30, ""{name = U}]{}{\quot{\eta}{\Gamma} \circ S^{-1}F(\overline{f})} \ar[r, swap, bend right = 30, ""{name = L}]{}{E(\overline{f}) \circ \quot{\eta}{\Gamma^{\prime}}} & E(\XGamma) \ar[from = U, to = L, Rightarrow, shorten <= 4pt, shorten >= 4pt]{}{\quot{\eta}{f}}
\end{tikzcd}
\]
that we will prove gives the morphism assignment of a pseudonatural transformation.

We now prove the pseudonaturality of $\eta$. For this fix a composable pair $\Gamma \xrightarrow{f} \Gamma^{\prime} \xrightarrow{g} \Gamma^{\prime\prime}$ in $\Sf(G)_0$. We begin by considering the pasting diagram:
\[
\begin{tikzcd}
 & S^{-1}F(\XGammap) \ar[d, Rightarrow, shorten <= 4pt, shorten >= 4pt]{}{\quot{\phi_{f,g}}{S^{-1}F}} \ar[dr]{}{S^{-1}F(\overline{f})} \\
S^{-1}F(\XGammapp) \ar[ur]{}{S^{-1}F(\overline{g})} \ar[d, swap]{}{\quot{\eta}{\Gamma^{\prime\prime}}} \ar[rr, ""{name = U}]{}[description]{S^{-1}F(\overline{g} \circ \overline{f})} & {} & S^{-1}F(\XGamma) \ar[d]{}{\quot{\eta}{\Gamma}} \\
E(\XGammapp) \ar[rr, swap, ""{name = L}]{}{E(\overline{g} \circ \overline{f})} & & E(\XGamma) \ar[from = U, to = L, Rightarrow, shorten <= 4pt, shorten >= 4pt]{}{\quot{\eta}{g \circ f}}
\end{tikzcd}
\]
A routine, but extremely tedious, check using that $\quot{\phi_{f,g}}{S^{1}F}^{A} = \quot{\phi_{f,g}}{F}^A$ for all $A \in S^{-1}F(\XGammapp)_0 = F(\XGammapp)_0$, the fact that $\quot{\lambda}{f}$ is an identity transformation, and the universal property that induced the $\quot{\eta}{f}$ allows us to translate the above diagram through the pseudonatruality of $\alpha$ (much like how we proved $\quot{\phi_{f,g}}{F}^{A} = \quot{\phi_{f,g}}{S^{-1}F}^{A}$), and then move back to the localizations to derive that the above pasting diagram is equivalent to the pasting diagram below:
\[
\begin{tikzcd}
S^{1}F(\XGammapp)  \ar[rr, ""{name = LeftT}]{}{S^{-1}F(\overline{g})} \ar[d, swap]{}{\quot{\eta}{\Gamma^{\prime\prime}}} & & S^{-1}F(\XGammap) \ar[rr, ""{name = MidT}]{}{S^{-1}F(\of)} \ar[d]{}{\quot{\eta}{\Gamma^{\prime}}} & & S^{-1}F(\XGamma) \ar[d]{}{\quot{\eta}{\Gamma}} \\
E(\XGammapp) \ar[rr, swap, ""{name = LeftB}]{}{E(\overline{g})}, \ar[rrrr, bend right = 30, swap, ""{name = BB}]{}{E(\overline{g} \circ \of)} & &  E(\XGammap) \ar[rr, swap, ""{name = MidB}]{}{E(\of)} & & E(\XGamma) \ar[from = LeftT, to = LeftB, Rightarrow, shorten >= 4pt, shorten <= 4pt]{}{\quot{\eta}{g}} \ar[from = MidT, to = MidB, Rightarrow, shorten >= 4pt, shorten <= 4pt]{}{\quot{\eta}{f}} \ar[from = 2-3, to = BB, Rightarrow, shorten <= 4pt, shorten >= 4pt]{}{\quot{\phi_{f,g}}{E}}
\end{tikzcd}
\]
Thus we have that $\eta$ is a pseudonatural transformation. Finally, we note that $\eta$ is even the unique pseudonatural transformation making
\[
\xymatrix{
F \ar[rr]^{\alpha} \ar[dr]_{\lambda} & & E \\
 & S^{-1}F \ar@{-->}[ur]_{\exists!\eta}
}
\] 
commute, as both $\quot{\eta}{\Gamma}$ and $\quot{\eta}{f}$ are induced by the universal properites of $1$-categories. This concludes the proof of the theorem.
\end{proof}
\begin{definition}\label{Defn: Equivariant Localization}\index{Equivariant Localization}
An equivariant localization of an equivariant category is a pair $(S^{-1}F_G(X),\underline{\lambda}_S)$ where $S \subseteq F_G(X)_1$ is a collection of morphisms as in Theorem \ref{Theorem: Section 3: Equivariant localization is a 2-colimit}.
\end{definition}

\begin{example}\label{Example: Section 3: Equivariant localizaiton is not the localization of equivariant category}
This example shows that if we na{\"i}vely localize an equivariant category at a class of morphisms, we can lose the equivariant nature of the category. This means in particular that equivariant localization is a delicate process; when localizing and preserving equivariant data, we need to make sure that we localize in a way that respects the descent data we have at hand.

Fix a variety $X$ over $\Spec K$, let $G = \Spec K$ be the trivial group, and let $R$ be a ring for which $a \in R$ is not a right zero divisor and $\afrak = Ra$ is a left ideal for which $\RMod(R/\afrak, R) = 0$ and $R/\afrak$ is a nonzero simple $R$-module (a specific example is $R = \Z_p$ and $a = p$). Now consider the pre-equivariant pseudofunctor $F:\SfResl_G(X)^{\op} \to \fCat$ given by $F(\XGamma) = \Ch(\RMod)$ and $F(\overline{f}) = \id_{\Ch(\RMod)}$. Consider now the complexes $A^{\bullet}$ and $B^{\bullet}$ in $\Ch(\RMod)$, where $A^{\bullet}$ is defined as
\[
\xymatrix{
\cdots \ar[r] & 0 \ar[r] & R \ar[r]^{\rho_a} & R \ar[r] & 0 \ar[r] & \cdots
}
\]
where $\rho_a$ is the right multiplication by $a$ map and $R$ appears in degrees $0$ and $1$, and $B^{\bullet}$ is defined by
\[
\xymatrix{
\cdots \ar[r] & 0 \ar[r] & R/\afrak \ar[r] & 0 \ar[r] & \cdots
}
\]
where $R/\afrak$ appears in degree $0$. It is routine to check that
\[
H^{n}(A^{\bullet}) = \begin{cases}
0 & \text{if}\, n \ne 0; \\
R/\afrak & \text{if}\, n = 0;
\end{cases}
\]
and that
\[
H^{n}(B^{\bullet}) = \begin{cases}
0 & \text{if}\, n \ne 0; \\
R/\afrak & \text{if}\, n = 0.
\end{cases}
\]
Because $R/\afrak$ is nonzero and simple, any nonzero morphism $\varphi^{\bullet}:A^{\bullet} \to B^{\bullet}$ extends to a quasi-isomorphism $H^{\ast}(\phi):H^{\ast}(A^{\bullet}) \xrightarrow{\cong} H^{\ast}(B^{\bullet})$. We take the morphism $\varphi^{\bullet}:A^{\bullet} \to B^{\bullet}$
\[
\varphi^{n} = \begin{cases}
0 & \text{if}\, n \ne 0; \\
\pi_{\afrak} & \text{if}\, n = 0
\end{cases}
\]
as our nonzero map, which is a quasi-isomorphism $A^{\bullet} \to B^{\bullet}$.

Now define $\quot{Q}{\Gamma} \subseteq F(\XGamma)_1$ by $\quot{Q}{\Gamma} := \lbrace f \in F(\XGamma)_1 \; | \; f \, \text{is\, a\, quasi-isomorphism}\rbrace$ and set $Q := \lbrace P \in F_G(X) \; | \; \forall\, \Gamma \in \Sf(G)_0 \forall\, \quot{\rho}{\Gamma} \in P, \quot{\rho}{\Gamma} \in \rbrace$. Write $(Q^{-1}F)_G(X)$ for the equivariant localization implied by Theorem \ref{Theorem: Section 3: Equivariant localization is a 2-colimit} and Definition \ref{Defn: Equivariant Localization} and write $Q^{-1}(F_G(X))$ for the localization of the category $F_G(X)$ at $Q$.

Our goal now is to show that $(Q^{-1}F)_G(X) \ne Q^{-1}(F_G(X))$. To show this we need only show that $(Q^{-1}F)_G(X)_0 \ne Q^{-1}(F_G(X))_0$. To this end consider the object $(A,T_A) \in Q^{-1}F_G(X)$ given by saying that the $\AGamma \in A$ take the form
\[
\AGamma = \begin{cases}
A^{\bullet} & \text{if}\, \Gamma \ne G; \\
B^{\bullet} & \text{if}\, \Gamma = G;
\end{cases}
\]
and defining the $\tau_f^A \in T_A$ by
\[
\tau_f^A = \begin{cases}
\varphi^{\bullet}:A^{\bullet} \to B^{\bullet} & \text{if}\,\exists\,\Gamma \in \Sf(G)_0. f = \nu_{\Gamma}:\Gamma \to \Spec K; \\
\id_{A^{\bullet}}:A^{\bullet} \to A^{\bullet} & \text{else}.
\end{cases}
\]
It is a trivial verification that $(A,T_A) \in (Q^{-1}F)_G(X)_0$. However, note that $(A,T_A) \notin F_G(X)_0$ because $\varphi^{\bullet}$ is a quasi-isomorphism (and hence an isomorphism in $Q^{-1}\Ch(\RMod)$ which is not an isomorphism in $\Ch(\RMod)$). As such we have $(Q^{-1}F)_G(X)_0 \ne Q^{-1}(F_G(X))_0$ and hence $(Q^{-1}F)_G(X) \ne Q^{-1}(F_G(X))$.
\end{example}

\section{Equivariant Functors: Change of Space}\label{Section: Section 3: Change of Space}
We now move to discuss when we can give functors which change the space of an equivariant category, i.e., equivariant functors of the form $F_G(X) \to E_G(Y)$ or $E_G(Y) \to F_G(X)$. We will start by studying the following situation: Given two left $G$-varieties $X$ and $Y$, we assume that we are given pre-equivariant pseudofunctors $F_X:\SfResl_G(X)^{\op} \to \fCat$ and $F_Y:\SfResl_G(Y)^{\op} \to \fCat$ which fit into the strictly commuting diagram:
\[
\xymatrix{
\SfResl_G(X)^{\op} \ar[dr]_{F_X} \ar[r]^-{\quo_X^{\op}} & \Var_{/\Spec K}^{\op} \ar[d]_{F} & \SfResl_G(Y)^{\op} \ar[l]_-{\quo_Y^{\op}} \ar[dl]^{F_Y} \\
 & \fCat &
}
\]
\begin{definition}
We call a pre-equivariant pseudofunctor \\$F:\Var_{/\Spec K}^{\op} \to \fCat$ a {simultaneous pre-equivariant pseudofunctor over $X$ and $Y$}\index{Pre-equivariant Pseudofunctor! Simultaneous} if there are pre-equivariant pseudofunctors $F_X$ on $X$ and $F_Y$ on $Y$ such that $F \circ \quo_Y^{\op} = F_Y$ and $F \circ \quo_X^{\op} = F_X$.
\end{definition}
\begin{example}
Let $G$ be a smooth algebraic group and let $X$ and $Y$ be any two $G$-varieties. The pseudofunctors $D^b(-)$, $\DbQl{-}, \Per(-), \Shv(-,\text{{\'e}t})$ are all simultaneous pre-equivariant pseudofunctors over $X$ and $Y$.
\end{example}

Now let $f \in \Sf(G)(\Gamma,\Gamma^{\prime})$ be given. Observe that since $h$ is a morphism of left $G$-varieties, we can produce a commuting cube in the category $\Var_{/\Spec K}$
\[
\begin{tikzcd}
 & \Gamma \times Y \ar[dd, swap, near end]{}{\quo_{\Gamma \times Y}} \ar[rr]{}{f \times \id_Y} & & \Gamma^{\prime} \times Y \ar[dd]{}{\quo_{\Gamma^{\prime} \times Y}} \\
\Gamma \times X \ar[dd, swap]{}{\quo_{\Gamma \times X}} \ar[rr, crossing over, near end]{}{f \times \id_X} \ar[ur]{}{\id_{\Gamma} \times h} & & \Gamma^{\prime} \times X \ar[ur,]{}{\id_{\Gamma^{\prime}}\times h} \\
 & \quot{Y}{\Gamma} \ar[rr, near end]{}{\quot{\overline{f}}{Y}} & & \quot{Y}{\Gamma^{\prime}} \\
\quot{X}{\Gamma} \ar[ur]{}{\quot{\overline{h}}{\Gamma}} \ar[rr, swap]{}{\quot{\overline{f}}{X}} & & \quot{X}{\Gamma^{\prime}} \ar[ur, swap]{}{\quot{\overline{h}}{\Gamma^{\prime}}} \ar[from = 2-3, to = 4-3, crossing over, near end]{}{\quo_{\Gamma^{\prime} \times X}}
\end{tikzcd}
\]\index[notation]{fX@$\fX$}\index[notation]{fY@$\fY$}\index[notation]{hGamma@$\hGamma$}\index[notation]{hGammaprime@$\hGammap$}where $\quot{\overline{f}}{X}:\quot{X}{\Gamma} \to \quot{X}{\Gamma^{\prime}}, \quot{\overline{f}}{Y}:\quot{Y}{\Gamma} \to \quot{Y}{\Gamma^{\prime}}, \quot{\overline{h}}{\Gamma}:\quot{X}{\Gamma} \to \quot{Y}{\Gamma}$, and $\quot{\overline{h}}{\Gamma^{\prime}}:\quot{X}{\Gamma^{\prime}} \to \quot{Y}{\Gamma^{\prime}}$ are the morphisms induced by the universal property of the quotient schemes. In particular, this gives rise to the commuting square
\[
\xymatrix{
\quot{X}{\Gamma} \ar[r]^-{\quot{\overline{f}}{X}} \ar[d]_{\quot{\overline{h}}{\Gamma}} & \quot{X}{\Gamma^{\prime}} \ar[d]^{\quot{\overline{h}}{\Gamma^{\prime}}} \\
\quot{Y}{\Gamma} \ar[r]_-{\quot{\overline{f}}{Y}} & \quot{Y}{\Gamma^{\prime}}
}
\]
of schemes. Define, for later use, $\quot{\overline{k}}{f} := \quot{\overline{f}}{Y} \circ \quot{\overline{h}}{\Gamma} = \quot{\overline{h}}{\Gamma^{\prime}} \circ \quot{\overline{f}}{X}$ for all morphisms $f \in \Sf(G)_1$\index[notation]{kf@$\quot{\overline{k}}{f}$}. Using the pseudofunctoriality of $F$ on $\Var_{/\Spec K}^{\op}$, we get not only that the compositors take values
\[
\phi_{\quot{f}{X}, \quot{h}{\Gamma^{\prime}}}:F(\quot{\overline{f}}{X}) \circ F(\quot{\overline{h}}{\Gamma^{\prime}}) \Rightarrow F(\quot{\overline{k}}{f})
\]
and
\[
\phi_{\quot{h}{\Gamma},\quot{f}{Y}}:F(\quot{\overline{h}}{\Gamma}) \circ F(\quot{\overline{f}}{Y}) \Rightarrow F(\quot{\overline{k}}{f}),
\]
we derive the existence of the commuting $2$-cells:
\[
\begin{tikzcd}
F(\quot{Y}{\Gamma^{\prime}}) \ar[rr, ""{name = U}]{}{F(\quot{\overline{h}}{\Gamma^{\prime}})} \ar[dd, swap]{}{F(\quot{\overline{f}}{Y})} \ar[ddrr, ""{name = M}]{}[description]{F(\quot{\overline{k}}{f})} & & F(\quot{X}{\Gamma^{\prime}}) \ar[dd]{}{F(\quot{\overline{f}}{X})} \\
\\
F(\quot{Y}{\Gamma}) \ar[rr, swap, ""{name = L}]{}{F(\quot{\overline{h}}{\Gamma})} & &  F(\quot{X}{\Gamma}) \ar[from = U, to = M, Rightarrow, shorten >= 5pt, shorten <= 5pt]{}{\phi_{\quot{f}{X}, \quot{h}{\Gamma^{\prime}}}} \ar[from = M, to = L, swap, Rightarrow, shorten >= 5pt, shorten <= 5pt]{}{\phi_{\quot{h}{\Gamma},\quot{f}{Y}}}
\end{tikzcd}
\]
A routine but tedious verification using the pseudofunctoriality of $F$ gives us also that for any composable arrows $\Gamma \xrightarrow{f} \Gamma^{\prime} \xrightarrow{g} \Gamma^{\prime\prime}$ in $\Sf(G)$, the pasting diagram
\[
\begin{tikzcd}
F(\quot{Y}{\Gamma^{\prime\prime}}) \ar[dd]{}{F(\quot{\overline{g}}{Y})} \ar[dddd, bend right = 80, swap, ""{name = LL}]{}{F(\quot{\overline{g}}{Y} \circ \quot{\overline{f}}{Y})} \ar[rrr, ""{name = UU}]{}{F(\quot{\overline{h}}{\Gamma^{\prime\prime}})} \ar[ddrrr, ""{name = UM}]{}[description]{F(\quot{\overline{k}}{g})} & & & F(\quot{X}{\Gamma^{\prime\prime}}) \ar[dddd, bend left = 80, ""{name = RR}]{}{F(\quot{\overline{g}}{X} \circ \quot{\overline{f}}{X})} \ar[dd]{}{F(\quot{\overline{g}}{X})}\\
\\
F(\quot{Y}{\Gamma^{\prime}}) \ar[dd]{}{F(\quot{\overline{f}}{Y})} \ar[rrr, ""{name = M}]{}[description]{F(\quot{\overline{h}}{\Gamma^{\prime}})} \ar[ddrrr, ""{name = LM}]{}[description]{F(\quot{\overline{k}}{f})} & & & F(\quot{X}{\Gamma^{\prime}}) \ar[dd]{}{F(\quot{\overline{f}}{X})} \\
\\
F(\quot{Y}{\Gamma}) \ar[rrr, swap, ""{name = BB}]{}{F(\quot{\overline{h}}{\Gamma})} & & & F(\quot{X}{\Gamma})
\ar[from = 3-4, to = RR, Rightarrow, shorten <=2pt, shorten >= 2pt]{}{}
\ar[from = LL, to = 3-1, Rightarrow, shorten <= 2pt, shorten >= 2pt]{}{}
\ar[from = UU, to = UM, Rightarrow, shorten <= 6pt, shorten >= 6pt]{}{}
\ar[from = UM, to = M, Rightarrow, shorten <= 6pt, shorten >= 6pt]
\ar[from = M, to = LM, Rightarrow, shorten <= 6pt, shorten >= 6pt]{}{}
\ar[from = LM, to = BB, Rightarrow, shorten <= 6pt, shorten >= 6pt]
\end{tikzcd}
\]
where each $2$-cell is a compositor or its inverse, is equivalent to the pasting diagram
\[
\begin{tikzcd}
F(\quot{Y}{\Gamma^{\prime\prime}}) \ar[dd, swap]{}{F(\quot{\overline{g}}{Y} \circ \quot{\overline{f}}{Y})} \ar[rr, ""{name = U}]{}{F(\quot{\overline{g}}{Y} \circ \quot{\overline{f}}{Y})} \ar[ddrr, ""{name = M}]{}[description]{F(\quot{\overline{k}}{g \circ f})} & & F(\quot{X}{\Gamma^{\prime\prime}}) \ar[dd]{}{F(\quot{\overline{g}}{X} \circ \quot{\overline{f}}{X})}\\
\\
F(\quot{Y}{\Gamma}) \ar[rr, swap, ""{name = B}]{}{F(\quot{\overline{h}}{\Gamma})} & & F(\quot{X}{\Gamma}) \ar[from = U, to = M, Rightarrow, shorten <= 6 pt, shorten >=6pt]{}{} \ar[from = M, to = B, Rightarrow, shorten <= 6pt, shorten >= 6pt]{}{}
\end{tikzcd}
\]
where the $2$-cells are again given first by the compositor $\phi_{\quot{g}{X} \circ \quot{f}{X}, \quot{h}{\Gamma^{\prime\prime}}}$ and second by $\phi_{\quot{h}{\Gamma}, \quot{g}{Y} \circ \quot{f}{Y}}^{-1}$. A routine analysis involving the use of the cocycle condition and the coherences of the natural isomorphisms then allows us to deduce the following technical lemma, which we only sketch. It is an extremely tedious but completely routine check using the pasting diagrams above and the tools to which we just referred.
\begin{lemma}\label{Lemma: Section 3: Compatible Family of things}
Let $\Gamma \xrightarrow{f} \Gamma^{\prime} \xrightarrow{g} \Gamma^{\prime\prime}$ be composable arrows in $\Sf(G)$ and let $h \in \GVar(X,Y)$ with $F:\Var_{/\Spec K}^{\op} \to \fCat$ a simultaneous pre-equivariant pseudofunctor over $X$ and $Y$. Then for any $\quot{A}{\Gamma^{\prime\prime}}$ there is an equality of morphisms
\begin{align*}
&\left(\phi_{\quot{h}{\Gamma}, \quot{g}{Y} \circ \quot{f}{Y}}^{\quot{A}{\Gamma^{\prime\prime}}}\right)^{-1} \phi_{\quot{g}{X} \circ \quot{f}{X}, \quot{h}{\Gamma^{\prime\prime}}}^{\quot{A}{\Gamma^{\prime\prime}}} \circ \phi_{\quot{f}{X}, \quot{g}{X} }^{F(\quot{\overline{h}}{\Gamma^{\prime\prime }})\quot{A}{\Gamma^{\prime\prime}}} \\ &= F(\quot{\overline{h}}{\Gamma})\phi_{\quot{f}{Y},\quot{g}{Y}}^{\quot{A}{\Gamma^{\prime\prime}}} \circ \left(\phi_{\quot{h}{\Gamma}, \quot{f}{Y}}^{F(\quot{\overline{g}}{Y})\quot{A}{\Gamma^{\prime\prime}}}\right)^{-1} \circ \phi_{\quot{f}{X}, \quot{h}{\Gamma^{\prime}}}^{F(\quot{\overline{g}}{Y})\quot{A}{\Gamma^{\prime\prime}}}  \circ F(\quot{\overline{f}}{X})\left(\phi_{\quot{h}{\Gamma^{\prime}},\quot{g}{Y}}^{\quot{A}{\Gamma^{\prime\prime}}}\right)^{-1} \\
&\circ F(\quot{\overline{f}}{X})\phi_{\quot{g}{X},\quot{h}{\Gamma^{\prime\prime}}}^{\quot{A}{\Gamma^{\prime\prime}}}
\end{align*}
\end{lemma}
\begin{proof}
We show that the morphisms type correctly, as that is the  the most difficult part of the lemma. For this note that on one hand the shorter of the two compositions is given by the diagram
\[
\begin{tikzcd}[cramped]
F(\quot{\overline{f}}{X})\big(F(\quot{\overline{g}}{X})\big(F(\quot{\overline{h}}{\Gamma^{\prime\prime}})(\quot{A}{\Gamma^{\prime\prime}})\big)\big) \ar[rr]{}{\phi_{\quot{f}{X},\quot{g}{X}}^{F(\quot{\overline{h}}{\Gamma^{\prime\prime}})\quot{A}{\Gamma^{\prime\prime}}}} & & F(\quot{\overline{g}}{X} \circ \quot{\overline{f}}{X})\big(F(\quot{\overline{h}}{\Gamma^{\prime\prime}})(\quot{A}{\Gamma^{\prime\prime}})\big) \ar[d]{}{\phi_{\quot{g}{X} \circ \quot{f}{X}, \quot{h}{\Gamma^{\prime\prime}}}^{\quot{A}{\Gamma^{\prime\prime}}}} \\
F(\quot{\overline{h}}{\Gamma}\big(F(\quot{\overline{g}}{Y} \circ \quot{\overline{f}}{Y})(\quot{A}{\Gamma^{\prime\prime}})\big) & & F(\quot{\overline{k}}{g \circ f})(\quot{A}{\Gamma^{\prime\prime}}) \ar[ll]{}{\left(\phi_{\quot{h}{\Gamma}, \quot{g}{Y}\circ\quot{f}{Y}}^{\quot{A}{\Gamma^{\prime\prime}}}\right)^{-1}}
\end{tikzcd}
\]
while the other composite is given by:
\[
\begin{tikzcd}[cramped]
F(\quot{\overline{f}}{X})\big(F(\quot{\overline{g}}{X})\big(F(\quot{\overline{h}}{\Gamma^{\prime\prime}})(\quot{A}{\Gamma^{\prime\prime}})\big)\big) \ar[rr]{}{F(\quot{\overline{f}}{X})\phi_{\quot{h}{\Gamma^{\prime}},\quot{g}{Y}}^{\quot{A}{\Gamma^{\prime\prime}}}} & & F(\quot{\overline{f}}{X})\big(F(\quot{\overline{h}}{\Gamma^{\prime\prime}} \circ \quot{\overline{g}}{X})(\quot{A}{\Gamma^{\prime\prime}})\big) \ar[d]{}{F(\quot{\overline{f}}{X})\left(\phi_{\quot{g}{X},\quot{h}{\Gamma^{\prime\prime}}}^{\quot{A}{\Gamma^{\prime\prime}}}\right)^{-1}} \\
F(\quot{\overline{k}}{f})\big(F(\quot{\overline{g}}{Y})(\quot{A}{\Gamma^{\prime\prime}})\big) \ar[d]{}{\left(\phi_{\quot{h}{\Gamma},\quot{f}{Y}}^{F(\quot{\overline{g}}{Y})\quot{A}{\Gamma^{\prime\prime}}}\right)^{-1}} & & F(\quot{\overline{f}}{X})\big(F(\quot{\overline{h}}{\Gamma^{\prime}})\big(F(\quot{\overline{g}}{Y})(\quot{A}{\Gamma^{\prime\prime}})\big)\big) \ar[ll]{}{\phi_{\quot{f}{X},\quot{h}{\Gamma^{\prime}}}^{F(\quot{\overline{g}}{Y})\quot{A}{\Gamma^{\prime\prime}}}} \\
F(\quot{\overline{h}}{\Gamma})\big(F(\quot{\overline{f}}{Y})\big(F(\quot{\overline{g}}{Y})(\quot{A}{\Gamma^{\prime\prime}})\big)\big) \ar[rr, swap]{}{F(\quot{\overline{h}}{\Gamma})\phi_{\quot{f}{Y},\quot{g}{Y}}^{\quot{A}{\Gamma^{\prime\prime}}}} & & F(\quot{\overline{h}}{\Gamma})\big(F(\quot{\overline{g}}{Y} \circ \quot{\overline{f}}{Y})(\quot{A}{\Gamma^{\prime\prime}})\big)
\end{tikzcd}
\]
After proceeding with the mechanincal proof invoking the consistent use of the pasting diagrams and pseudofunctoriality, we get our lemma.
\end{proof}
 
The technical lemma above is needed to give a pullback functor $h^{\ast}:F_G(Y) \to F_G(X)$ induced by the pseudofunctor $F$ whenever there is a morphism $h:X \to Y$. We show below that in all circumstances we can provide such a functor and that it looks more or less how we would expect at this point.

\begin{Theorem}\label{Thm: Section 3: Existence of pullback functors}
For any morphism $h \in \GVar(X,Y)$, there is a pullback functor $h^{\ast}:F_G(Y) \to F_G(X)$, where $F$ is a simultaneous pre-equivariant pseudofunctor over $X$ and $Y$, induced on object sets by
\[
\lbrace \quot{A}{\Gamma} \; | \; \Gamma \in \Sf(G)_0 \rbrace \mapsto \lbrace F(\quot{\overline{h}}{\Gamma})\quot{A}{\Gamma} \; | \; \Gamma \in \Sf(G)_0 \rbrace
\]
and on morphism sets by
\[
\lbrace \quot{\rho}{\Gamma} \; | \; \Gamma \in \Sf(G)_0 \rbrace \mapsto \lbrace F(\quot{\overline{h}}{\Gamma})\quot{\rho}{\Gamma} \; | \; G \in \Sf(G)_0 \rbrace.
\]
\end{Theorem}
\begin{proof}
Fix an object $(A,T_A)$ in $F_G(Y)_0$. We have already seen how to define $h^{\ast}A$, so we need only define $T_{h^{\ast}A}$. For this, fix an $f \in \Sf(G)_1$ and write $\Dom f = \Gamma$ and $\Codom f = \Gamma^{\prime}$. We now consider that $\tau_f^{h^{\ast}A}$ has to map from $F(\quot{\overline{f}}{X})\big(F(\quot{\overline{h}}{\Gamma^{\prime}})\quot{A}{\Gamma^{\prime}}\big) \to F(\quot{\overline{h}}{\Gamma})\quot{A}{\Gamma}$; however, $F(\quot{\overline{h}}{\Gamma})\tau_f^A:F(\quot{\overline{h}}{\Gamma})\big(F(\overline{f}_Y)\quot{A}{\Gamma^{\prime}}\big) \to F(\quot{\overline{h}}{\Gamma})\quot{A}{\Gamma}$ does not type correctly; to fix this we define $\tau_f^{h^{\ast}A}$ to be given by $F(\quot{\overline{h}}{\Gamma})\tau_f^A$ pre-composed by the compositors, i.e., we define $\tau_f^{h^{\ast}A}$ as in the diagram below:
\[
\xymatrix{
F(\quot{\overline{f}}{X})\big(F(\quot{\overline{h}}{\Gamma^{\prime}})(\quot{A}{\Gamma^{\prime}})\big) \ar[d]_{\tau_f^{h^{\ast}A}} \ar[rr]^-{\phi_{\quot{f}{X},\quot{h}{\Gamma^{\prime}}}^{\quot{A}{\Gamma^{\prime}}}} & & F(\quot{\overline{k}}{f})(\quot{A}{\Gamma^{\prime}}) \ar[d]^-{\left(\phi_{\quot{h}{\Gamma}, \quot{f}{Y}}^{\quot{A}{\Gamma^{\prime}}}\right)^{-1}} \\
 F(\quot{\overline{h}}{\Gamma})(\quot{A}{\Gamma}) & & F(\quot{\overline{h}}{\Gamma})\big(F(\quot{\overline{f}}{Y})(\quot{A}{\Gamma^{\prime}})\big) \ar[ll]^-{F(\quot{\overline{h}}{\Gamma})\tau_f^A} 
}
\]
We then set
\[
T_{h^{\ast}A} := \lbrace \tau_f^{h^{\ast}A} \; | \; f \in \Sf(G)_1\rbrace.
\]
Let us now verify that the $\tau_f^{h^{\ast}A}$ satisfy the cocycle condition. For this fix a composable pair of morphisms $\Gamma \xrightarrow{f} \Gamma^{\prime} \xrightarrow{g} \Gamma^{\prime\prime}$ in $\Sf(G)$ and consider that $\tau_f^{h^{\ast}A} \circ F(\fX)\tau_g^{h^{\ast}A}$ is equal to the expression
\begin{align*}
&F(\hGamma)\tau_f^A \circ \left(\phi_{\quot{h}{\Gamma},\quot{f}{Y}}^{\quot{A}{\Gamma^{\prime}}}\right)^{-1} \circ \phi_{\quot{f}{X},\quot{h}{\Gamma^{\prime}}}^{\quot{A}{\Gamma^{\prime}}} \circ F(\fX)\tau_g^{h^{\ast}A} \\
&=F(\hGamma)\tau_f^A \circ \left(\phi_{\quot{h}{\Gamma},\quot{f}{Y}}^{\quot{A}{\Gamma^{\prime}}}\right)^{-1} \circ \phi_{\quot{f}{X},\quot{h}{\Gamma^{\prime}}}^{\quot{A}{\Gamma^{\prime}}} \\
&\circ F(\fX)\left(F(\hGammap)\tau_g^A \circ \left(\phi_{\quot{h}{\Gamma^{\prime\prime}}, \quot{g}{Y}}^{ \quot{A}{\Gamma^{\prime\prime}}}\right)^{-1} \circ \phi_{\quot{g}{X},\quot{h}{\Gamma^{\prime\prime}}}^{\quot{A}{\Gamma^{\prime\prime}}}\right) \\
&= F(\hGamma)\tau_f^A \circ \left(\phi_{\quot{h}{\Gamma},\quot{f}{Y}}^{\quot{A}{\Gamma^{\prime}}}\right)^{-1} \circ \phi_{\quot{f}{X},\quot{h}{\Gamma^{\prime}}}^{\quot{A}{\Gamma^{\prime}}} \circ F(\fX)\left(F(\hGammap)\tau_g^A\right) \\
&\circ F(\fX)\left(\left(\phi_{\quot{h}{\Gamma^{\prime\prime}}, \quot{g}{Y}}^{ \quot{A}{\Gamma^{\prime\prime}}}\right)^{-1} \circ \phi_{\quot{g}{X},\quot{h}{\Gamma^{\prime\prime}}}^{\quot{A}{\Gamma^{\prime\prime}}}\right) \\
&= F(\hGamma)\tau_f^A \circ \left(\phi_{\quot{h}{\Gamma},\quot{f}{Y}}^{\quot{A}{\Gamma^{\prime}}}\right)^{-1} \circ F(\quot{\overline{k}}{f})\tau_g^A \circ \phi_{\quot{f}{X}, \quot{h}{\Gamma^{\prime}}}^{F(\gY)\quot{A}{\Gamma^{\prime\prime}}} \\
&\circ  F(\fX)\left(\left(\phi_{\quot{h}{\Gamma^{\prime\prime}}, \quot{g}{Y}}^{ \quot{A}{\Gamma^{\prime\prime}}}\right)^{-1} \circ \phi_{\quot{g}{X},\quot{h}{\Gamma^{\prime\prime}}}^{\quot{A}{\Gamma^{\prime\prime}}}\right) \\
&= F(\hGamma)\tau_f^A \circ F(\hGamma)\left(F(\fY)\tau_g^A\right) \circ \left(\phi_{\quot{h}{\Gamma},\quot{f}{Y}}^{F(\gY)\quot{A}{\Gamma^{\prime\prime}}}\right)^{-1} \circ \phi_{\quot{f}{X}, \quot{h}{\Gamma^{\prime}}}^{F(\gY)\quot{A}{\Gamma^{\prime\prime}}} \\
&\circ  F(\fX)\left(\left(\phi_{\quot{h}{\Gamma^{\prime\prime}}, \quot{g}{Y}}^{ \quot{A}{\Gamma^{\prime\prime}}}\right)^{-1} \circ \phi_{\quot{g}{X},\quot{h}{\Gamma^{\prime\prime}}}^{\quot{A}{\Gamma^{\prime\prime}}}\right) \\
&= F(\hGamma)\left(\tau_f^A \circ F(\fY)\tau_g^A\right) \circ \left(\phi_{\quot{h}{\Gamma},\quot{f}{Y}}^{F(\gY)\quot{A}{\Gamma^{\prime\prime}}}\right)^{-1} \circ \phi_{\quot{f}{X}, \quot{h}{\Gamma^{\prime}}}^{F(\gY)\quot{A}{\Gamma^{\prime\prime}}} \\
&\circ  F(\fX)\left(\left(\phi_{\quot{h}{\Gamma^{\prime\prime}}, \quot{g}{Y}}^{ \quot{A}{\Gamma^{\prime\prime}}}\right)^{-1} \circ \phi_{\quot{g}{X},\quot{h}{\Gamma^{\prime\prime}}}^{\quot{A}{\Gamma^{\prime\prime}}}\right) \\
&= F(\hGamma)\left(\tau_{g \circ f}^{A} \circ \phi_{\quot{f}{Y},\quot{g}{Y}}^{\quot{A}{\Gamma^{\prime\prime}}}\right) \circ \left(\phi_{\quot{h}{\Gamma},\quot{f}{Y}}^{F(\gY)\quot{A}{\Gamma^{\prime\prime}}}\right)^{-1} \circ \phi_{\quot{f}{X}, \quot{h}{\Gamma^{\prime}}}^{F(\gY)\quot{A}{\Gamma^{\prime\prime}}} \\
&\circ  F(\fX)\left(\left(\phi_{\quot{h}{\Gamma^{\prime\prime}}, \quot{g}{Y}}^{ \quot{A}{\Gamma^{\prime\prime}}}\right)^{-1} \circ \phi_{\quot{g}{X},\quot{h}{\Gamma^{\prime\prime}}}^{\quot{A}{\Gamma^{\prime\prime}}}\right) \\
& =F(\hGamma)\tau_{g \circ f}^{A} \circ F(\hGamma)\phi_{\quot{f}{Y},\quot{g}{Y}}^{\quot{A}{\Gamma^{\prime\prime}}} \circ \left(\phi_{\quot{h}{\Gamma},\quot{f}{Y}}^{F(\gY)\quot{A}{\Gamma^{\prime\prime}}}\right)^{-1} \circ \phi_{\quot{f}{X}, \quot{h}{\Gamma^{\prime}}}^{F(\gY)\quot{A}{\Gamma^{\prime\prime}}} \\
&\circ  F(\fX)\left(\phi_{\quot{h}{\Gamma^{\prime\prime}}, \quot{g}{Y}}^{ \quot{A}{\Gamma^{\prime\prime}}}\right)^{-1} \circ F(\fX)\phi_{\quot{g}{X},\quot{h}{\Gamma^{\prime\prime}}}^{\quot{A}{\Gamma^{\prime\prime}}}.
\end{align*}
From here invoking Lemma \ref{Lemma: Section 3: Compatible Family of things} gives that
\begin{align*}
&F(\hGamma)\tau_{g \circ f}^{A} \circ F(\hGamma)\phi_{\quot{f}{Y},\quot{g}{Y}}^{\quot{A}{\Gamma^{\prime\prime}}} \circ \left(\phi_{\quot{h}{\Gamma},\quot{f}{Y}}^{F(\gY)\quot{A}{\Gamma^{\prime\prime}}}\right)^{-1} \circ \phi_{\quot{f}{X}, \quot{h}{\Gamma^{\prime}}}^{F(\gY)\quot{A}{\Gamma^{\prime\prime}}} \\
&\circ  F(\fX)\left(\phi_{\quot{h}{\Gamma^{\prime\prime}}, \quot{g}{Y}}^{ \quot{A}{\Gamma^{\prime\prime}}}\right)^{-1} \circ F(\fX)\phi_{\quot{g}{X},\quot{h}{\Gamma^{\prime\prime}}}^{\quot{A}{\Gamma^{\prime\prime}}} \\
&= F(\hGamma)\tau_{g \circ f}^{A} \circ \left(\phi_{\quot{h}{\Gamma},\quot{g}{Y} \circ \quot{f}{Y}}^{\quot{A}{\Gamma^{\prime\prime}}}\right)^{-1} \circ \phi_{\quot{g}{X} \circ \quot{f}{X}, \quot{h}{\Gamma^{\prime\prime}}}^{\quot{A}{\Gamma^{\prime\prime}}} \circ \phi_{\quot{f}{X}, \quot{g}{X}}^{F(\hGammapp)\quot{A}{\Gamma^{\prime\prime}}} \\
&= \tau_{g \circ f}^{h^{\ast}A} \circ \phi_{\quot{f}{X}, \quot{g}{X}}.
\end{align*}
Therefore it follows that $\tau_{g \circ f}^{h^{\ast}A} \circ \phi_{\quot{f}{X}, \quot{g}{X}} = \tau_f^{h^{\ast}A} \circ F(\fX)\tau_g^{h^{\ast}A}$, which verifies that the maps in $T_{h^{\ast}A}$ satisfy the cocycle condition. Thus $(h^{\ast}A,T_{h^{\ast}A})$ is an $F_G(X)$ object.

We now prove that $h^{\ast}$ sends morphisms to morphisms; note that the pointwise application of the $F(\hGamma)$ in the definition of $h^{\ast}$ on morphisms implies that as soon as we know that $h^{\ast}$ sends morphisms to morphisms, it immediately preserves composition and identities. Let $P = \lbrace \quot{\rho}{\Gamma} \; | \; \Gamma \in \Sf(G)_0 \rbrace \in F_G(Y)(A,B)$ and fix a morphism $f \in \Sf(G)_1$. Write, as usual, $\Dom f = \Gamma, \Codom f = \Gamma^{\prime}$. We then see that the diagram
\[
\xymatrix{
F(\fX)\big(F(\hGammap)(\quot{A}{\Gamma^{\prime}})\big) \ar[d]_{\tau_f^{h^{\ast}A}} \ar[rrr]^-{F(\fX)\big(F(\hGammap)\quot{\rho}{\Gamma^{\prime}}\big)} & & & F(\fX)\big(F(\hGammap)(\quot{B}{\Gamma^{\prime}})\big) \ar[d]^{\tau_f^{h^{\ast}B}} \\
F(\hGamma)(\quot{A}{\Gamma}) \ar[rrr]_-{F(\hGamma)\quot{\rho}{\Gamma}} & & & F(\hGamma)(\quot{B}{\Gamma})
}
\]
is equivalent to the diagram:
\[
\xymatrix{
F(\fX)\big(F(\hGammap)(\quot{A}{\Gamma^{\prime}})\big) \ar[d]_{\phi_{\quot{f}{X},\quot{h}{\Gamma^{\prime}}}^{\quot{A}{\Gamma^{\prime}}}} \ar[rrr]^-{F(\fX)\big(F(\hGammap)\quot{\rho}{\Gamma^{\prime}}\big)} & & & F(\fX)\big(F(\hGammap)(\quot{B}{\Gamma^{\prime}})\big) \ar[d]^{\phi_{\quot{f}{X},\quot{h}{\Gamma^{\prime}}}^{\quot{A}{\Gamma^{\prime}}}} \\
F(\quot{\overline{k}}{f})(\quot{A}{\Gamma^{\prime}}) \ar[rrr]^-{F(\quot{\overline{k}}{f})} \ar[d]_{\left(\phi_{\quot{h}{\Gamma},\quot{f}{Y}}^{\quot{A}{\Gamma^{\prime}}}\right)^{-1}} & & & F(\quot{\overline{k}}{f})(\quot{B}{\Gamma^{\prime}}) \ar[d]^{\left(\phi_{\quot{h}{\Gamma},\quot{f}{Y}}^{\quot{B}{\Gamma^{\prime}}}\right)^{-1}} \\
F(\hGamma)\big(F(\fY)(\quot{A}{\Gamma^{\prime}})\big) \ar[rrr]^-{F(\hGamma)\big(F(\fY)\quot{\rho}{\Gamma^{\prime}}\big)} \ar[d]_{F(\hGamma)\tau_f^A} & & & F(\hGamma)\big(F(\fY)(\quot{B}{\Gamma^{\prime}})\big) \ar[d]^{F(\hGamma)\tau_f^B} \\
F(\hGamma)(\quot{A}{\Gamma}) \ar[rrr]_-{F(\hGamma)\quot{\rho}{\Gamma}} & & & F(\hGamma)(\quot{B}{\Gamma})
}
\]
Note that each rectangle in the above diagram commutes either by the naturality of the compositors (and their inverses) or by the fact that $P$ is an $F_G(Y)$ morphism, respectively. However, this implies that the original rectangle commutes, giving that $h^{\ast}P$ is a morphism.
\end{proof}

\begin{proposition}\label{Prop: Section 3: Pullbacks are lex/rex/additive whenever the F(h) are}
Let $F$ be a simultaneous pre-equivariant pseudofunctor over $X$ and $Y$ and let $h \in \GVar(X,Y)$. Assume $I$ is an index category such that for all $\Gamma \in \Sf(G)_0$ there are diagram functors $\quot{d}{\Gamma}:I \to F(\YGamma)$, and assume $d:I \to F_G(Y)$ is a diagram functor for which $d(i) = \lbrace \quot{d}{\Gamma}(i) \; | \; \Gamma \in \Sf(G)_0 \rbrace$. Assume also that each category $F(\YGamma)$ admits limits and colimits of shape $I$. Then if each functor $F(\quot{\overline{h}}{\Gamma})$ satisfies any of the following properties for all $\Gamma \in \Sf(G)_0$, so does $h^{\ast}$:
\begin{enumerate}
	\item $F(\quot{\overline{h}}{\Gamma})$ preserves all limits of specified shape $I$;
	\item $F(\quot{\overline{h}}{\Gamma})$ preserves all colimits of sepcified shape $I$;
	\item The pre-equivariant pseudofunctor $F$ is monoidal as well, and $F(\quot{\overline{h}}{\Gamma})$ is a monoidal functor.
\end{enumerate}
\end{proposition}
\begin{proof}
Note that $(1)$ and $(2)$ are dual, so we only prove $(1)$. Begin by writing $\quot{d(i)}{\Gamma} := \quot{A_i}{\Gamma}$ for all $i \in I_0$ and letting
\[
A_i := d(i) = \lbrace \quot{A_i}{\Gamma} \; | \; \Gamma \in \Sf(G) \rbrace.
\]
Now let
\[
\quot{A}{\Gamma} := \lim_{\substack{\longleftarrow \\ i \in I}} \quot{A_i}{\Gamma}
\]
and write
\[
A := \lim_{\substack{\longleftarrow \\ i \in I_0}} A_i;
\]
note that this limit exists in $F_G(Y)$ by Theorem \ref{Thm: Section 2: Equivariant Cat has lims}. Moreover, let $P_i \in F_G(X)(A,A_i)$ be the limit morphisms and write
\[
P_i := \lbrace \quot{\rho_i}{\Gamma} \; | \; \Gamma \in \Sf(G) \rbrace.
\]

We now show that $h^{\ast}A$ is a universal cone to the $h^{\ast}A_i$ in $F_G(X)$. To this end let $B \in F_G(X)_0$ be a cone over the $h^{\ast}A_i$; that is, $B \in F_G(X)_0$ such that there are morphisms $\Psi_i \in F_G(X)(B,h^{\ast}A_i)$ for which given any $i,j \in I_0$ and any $\alpha \in I(i,j),$ we have the commuting diagram
\[
\xymatrix{
 & B \ar[dl]_{\Psi_i} \ar[dr]^{\Psi_j} & \\
h^{\ast}A_i \ar[rr]_-{d(\alpha)} & & h^{\ast}A_j
}
\]
in $F_G(X)$. Write $\Psi_i = \lbrace \quot{\psi_i}{\Gamma} \; | \; \Gamma \in \Sf(G)_0\rbrace, \Psi_j = \lbrace \quot{\psi_j}{\Gamma} \; | \;  \Gamma \in \Sf(G)_0 \rbrace, d(\alpha) = \lbrace \quot{d(\alpha)}{\Gamma} \; | \; \Gamma \in \Sf(G)_0$. Now observe that since each functor $F(\hGamma)$ preserves limits, the object $F(\hGamma)\AGamma$ is a limit of the $F(\hGamma)\AGamma_i$. Thus there is a unique morphism $\quot{\zeta}{\Gamma} \in F_G(X)(\BGamma,F(\hGamma)\AGamma)$ making the diagram
\[
\xymatrix{
	& \quot{B}{\Gamma} \ar@{-->}[d]^{\quot{\zeta}{\Gamma}} \ar@/^/[ddr]^{\quot{\psi_j}{\Gamma}} \ar@/_/[ddl]_{\quot{\psi_i}{\Gamma}} & \\
	& F(\hGamma)(\quot{A}{\Gamma}) \ar[dr]_{\quot{F(\hGamma)\rho_j}{\Gamma}} \ar[dl]^{\quot{F(\hGamma)\rho_i}{\Gamma}} \\
	F(\hGamma)(\quot{A_i}{\Gamma}) \ar[rr]_-{\quot{d(\alpha)}{\Gamma}} & & F(\hGamma)(\quot{A_j}{\Gamma})
}
\]
commute in $F(\XGamma)$. We now define $Z:B \to h^{\ast}A$ via
\[
Z := \lbrace \quot{\zeta}{\Gamma} \; | \; \Gamma \in \Sf(G)_0\rbrace.
\]
Note that if $Z$ is an $F_G(X)$-morphism, then it is automatically unique, which in turn will prove that $h^{\ast}A$ is a limit of the $h^{\ast}A_i$. Thus it suffices to prove that $Z$ is a morphism in $F_G(X)$.

To prove that $Z$ is an $F_G(X)$-morphism, let $f \in \Sf(G)_1$ and write $\Dom f = \Gamma$, $\Codom f = \Gamma^{\prime}$. We now need to verify that the diagram
\[
\xymatrix{
F(\fX)(\BGammap) \ar[rr]^-{F(\fX)\quot{\zeta}{\Gamma^{\prime}}} \ar[d]_{\tau_f^{B}} & & F(\fX)\big(F(\hGamma)(\AGammap)\big) \ar[d]^{\tau_f^{h^{\ast}A}} \\
\BGamma \ar[rr]_-{\quot{\zeta}{\Gamma}} & & F(\hGamma)(\AGamma)
}
\]
commutes. However, this follows mutatis mutandis from the last part of the proof of Theorem \ref{Thm: Section 2: Equivariant Cat has lims}. To see how this is the case, note that the diagram
\[
\xymatrix{
 & F(\fX)(\BGammap) \ar[dl]_{F(\fX)\quot{\psi_i}{\Gamma^{\prime}}} \ar[dr]^{F(\fX)\quot{\psi_j}{\Gamma^{\prime}}} \\ 
F(\fX)\big(F(\hGammap)(\AGammap_i)\big) \ar[rr]_-{F(\fX)\big(F(\hGammap)\quot{d(\alpha)}{\Gamma^{\prime}}\big)} \ar[d]_{\tau_f^{h^{\ast}A_i}} & & F(\fX)\big(F(\hGammap)(\AGammap_j)\big) \ar[d]^{\tau_f^{h^{\ast}A_j}} \\
F(\hGamma)(\AGamma_i) \ar[rr]_-{F(\hGamma)\quot{d(\alpha)}{\Gamma}} & & F(\hGamma)(\AGamma_j)
}
\]
commutes because $B$ is a cone over the $h^{\ast}A_i$ and because $h^{\ast}d(\alpha)$ is an $F_G(X)$-morphism.  From here  proceed as in Theorem \ref{Thm: Section 2: Equivariant Cat has lims}.

We now prove $(3)$: In particular, we show that $h^{\ast}$ is monoidal. For this let us begin by codifying some notation. We will write $\otimes_X$ and $\otimes_Y$ for the tensor functors on $F_G(X)$ and $F_G(Y)$, respectively, and also write
\[
\OtimesXGamma:F(\XGamma) \times F(\XGamma) \to F(\XGamma), \qquad \OtimesYGamma:F(\YGamma) \times F(\YGamma) \to F(\YGamma)
\]
for the tensor functors on $F(\XGamma)$ and $F(\YGamma)$, respecitvely. Now, because each functor $F(\hGamma)$ is monoidal there is a natural isomorphism
\[
\theta_{\quot{h}{\Gamma}}^{A,B}:F(\hGamma)\left(A \OtimesYGamma B\right) \xrightarrow{\cong} F(\hGamma)A \OtimesXGamma F(\hGamma)B
\]
for all $A, B \in F(\YGamma)_0$ and a series of isomorphisms
\[
\sigma_{\quot{h}{\Gamma}}:F(\hGamma)\quot{I_Y}{\Gamma} \xrightarrow{\cong} \quot{I_X}{\Gamma}
\]
for each $\Gamma \in \Sf(G)_0$ which satisfy the cocycle condition for all composites of morphisms. Note that this trivially gives an isomorphism $h^{\ast}I_Y \cong I_X$; we omit these details and instead show that $h^{\ast}$ preserves tensors up to isomorphism. Now fix two objects $A,B \in F_G(X)_0$; we define the map $\Theta:h^{\ast}(A \otimes_Y B) h^{\ast} \to A \otimes_X h^{\ast}B$ by setting
\[
\Theta := \left\lbrace \theta_{\quot{h}{\Gamma}}^{A,B}:F(\hGamma)\left(A \OtimesYGamma B\right) \xrightarrow{\cong} F(\hGamma)A \OtimesXGamma F(\hGamma)B \; | \; \Gamma \in \Sf(G)_0 \right\rbrace
\]
recall that we are writing $\theta_{\quot{h}{\Gamma}}^{A,B}$ for each $\Gamma$ as an abuse of notation to reduce notational clutter in the superscript whenever possible. Note also that if we can show that $\Theta$ is an $F_G(X)$-morphism, it is both an isomorphism and the object assignment of a natural transformation immediately by the fact that all the $\theta_{\quot{h}{\Gamma}}$ are natrual isomorphisms as well.

We now prove that $\Theta$ is an $F_G(X)$-morphism. Let $f \in \Sf(G)_1$ and write $\Dom f = \Gamma$ and $\Codom f = \Gamma^{\prime}$. Note that we must prove that the diagram
\[
\theta_{\quot{h}{\Gamma}}^{A,B} \circ \tau_f^{h^{\ast}(A \otimes B)} = \tau_{f}^{h^{\ast}A \otimes_X h^{\ast}B} \circ F(\fX)\theta_{\quot{h}{\Gamma^{\prime}}}^{A,B}
\]
commutes. As a first step to this, we compute using Theorems \ref{Theorem: Section 2: Monoidal preequivariant pseudofunctor gives monoidal equivariant cat} and \ref{Thm: Section 3: Existence of pullback functors} that we have
\begin{align*}
\tau_f^{h^{\ast}(A \otimes_Y B)} &= F(\hGamma)\tau_f^{A \otimes_Y B} \circ \left(\phi_{\quot{h}{\Gamma}, \quot{f}{Y}}^{A \otimes B}\right)^{-1} \circ \phi_{\quot{f}{X}, \quot{h}{\Gamma^{\prime}}}^{A \otimes B} \\
&= F(\hGamma)\left(\left(\tau_f^A \OtimesYGamma \tau_f^B\right) \circ \theta_{\quot{f}{Y}}^{A,B}\right) \circ \left(\phi_{\quot{h}{\Gamma}, \quot{f}{Y}}^{A \otimes B}\right)^{-1} \circ \phi_{\quot{f}{X}, \quot{h}{\Gamma^{\prime}}}^{A \otimes B}
\end{align*}
while
\begin{align*}
&\tau_f^{h^{\ast}A \otimes_X h^{\ast}B} = \left(\tau_f^{h^{\ast}A} \OtimesXGamma \tau_f^{h^{\ast}B}\right) \circ \theta_{\quot{f}{X}}^{A,B} \\
&= \left(F(\hGamma)\tau_f^A \circ \left(\phi_{\quot{h}{\Gamma},\quot{f}{Y}}^{A}\right)^{-1} \circ \phi_{\quot{f}{X},\quot{h}{\Gamma^{\prime}}}^{A}\right) \OtimesXGamma \left(F(\hGamma)\tau_f^B \circ \left(\phi_{\quot{h}{\Gamma},\quot{f}{Y}}^{B}\right)^{-1} \circ \phi_{\quot{f}{X},\quot{h}{\Gamma^{\prime}}}^{B}\right) \\
&\circ \theta_{\quot{f}{X}}^{A,B} \\
 &= \left(F(\hGamma)\tau_f^A \OtimesXGamma F(\hGamma)\tau_f^B\right) \circ \left(\left(\phi_{\quot{h}{\Gamma},\quot{f}{Y}}^{A}\right)^{-1} \OtimesXGamma \left(\phi_{\quot{h}{\Gamma},\quot{f}{Y}}^{B}\right)^{-1}\right) \\
 &\circ \left(\phi_{\quot{f}{X},\quot{h}{\Gamma^{\prime}}}^{A} \OtimesXGamma \phi_{\quot{f}{X},\quot{h}{\Gamma^{\prime}}}^{B}\right) \circ \theta_{\quot{f}{X}}^{A,B}.
\end{align*}
Using these identities and the pseudonaturality of the $\theta$'s (cf.\@ Definition \ref{Defn: Section 2: Monoidal Preequivariant Pseudofunctor}) we compute that
\begin{align*}
&\tau_f^{h^{\ast}A \otimes_X h^{\ast}B} \circ F(\fX) \theta_{\quot{h}{\Gamma^{\prime}}}^{A,B} = \left(\tau_f^{h^{\ast}A} \OtimesXGamma \tau_f^{h^{\ast}B}\right) \circ \theta_{\quot{f}{X}}^{A,B} \circ F(\fX)\theta_{\quot{h}{\Gamma^{\prime}}}^{A,B} \\
&= \left(\tau_f^{h^{\ast}A} \OtimesXGamma \tau_f^{h^{\ast}B}\right) \circ \left(\left(\phi_{\quot{f}{X},\quot{h}{\Gamma^{\prime}}}^{A}\right)^{-1} \OtimesXGamma \left(\phi_{\quot{f}{X},\quot{h}{\Gamma^{\prime}}}^{B}\right)^{-1}\right) \circ \theta_{\quot{h}{\Gamma^{\prime}} \circ \quot{f}{X}}^{A,B} \circ \phi_{\quot{f}{X},\quot{h}{\Gamma^{\prime}}}^{A \otimes B} \\
&=  \left(F(\hGamma)\tau_f^A \OtimesXGamma F(\hGamma)\tau_f^B\right) \circ \left(\left(\phi_{\quot{h}{\Gamma},\quot{f}{Y}}^{A}\right)^{-1} \OtimesXGamma \left(\phi_{\quot{h}{\Gamma},\quot{f}{Y}}^{B}\right)^{-1}\right) \circ \theta_{\quot{h}{\Gamma^{\prime}} \circ \quot{f}{X}}^{A,B} \\
&\circ \phi_{\quot{f}{X},\quot{h}{\Gamma^{\prime}}}^{A \otimes B} \\
&= \left(F(\hGamma)\tau_f^A \OtimesXGamma F(\hGamma)\tau_f^B\right) \circ \left(\phi_{\quot{h}{\Gamma},\quot{f}{Y}}^{A} \OtimesXGamma \phi_{\quot{h}{\Gamma},\quot{f}{Y}}^{B}\right)^{-1} \circ \theta_{\quot{h}{\Gamma^{\prime}} \circ \quot{f}{X}}^{A,B} \circ \phi_{\quot{f}{X},\quot{h}{\Gamma^{\prime}}}^{A \otimes B}. 
\end{align*}
Expanding further gives that $\tau_f^{h^{\ast}A \otimes_X h^{\ast}B}$ is equal to
\begin{align*}
& \left(F(\hGamma)\tau_f^A \OtimesXGamma F(\hGamma)\tau_f^B\right) \circ \theta_{\quot{f}{X}}^{F(\hGamma)A,F(\hGamma)B} \circ F(\hGamma)\theta_{\quot{f}{Y}}^{A,B} \circ \left(\phi_{\quot{h}{\Gamma},\quot{f}{Y}}^{A \otimes B}\right)^{-1} \\
&\circ \phi_{\quot{f}{X},\quot{h}{\Gamma^{\prime}}}^{A \otimes B} \\
&= \theta_{\quot{h}{\Gamma}}^{A,B} \circ F(\hGamma)\left(\tau_f^{A}\OtimesYGamma \tau_f^B\right) \circ F(\hGamma)\theta_{\quot{f}{Y}}^{A,B} \circ \left(\phi_{\quot{h}{\Gamma},\quot{f}{Y}}^{A \otimes B}\right)^{-1} \circ \phi_{\quot{f}{X},\quot{h}{\Gamma^{\prime}}}^{A \otimes B} \\
&= \theta_{\quot{h}{\Gamma}}^{A,B} \circ F(\hGamma)\left(\left(\tau_f^{A} \OtimesYGamma \circ \theta_{\quot{f}{Y}}^{A,B}\right)\right) \circ \left(\phi_{\quot{h}{\Gamma},\quot{f}{Y}}^{A \otimes B}\right)^{-1} \circ \phi_{\quot{f}{X},\quot{h}{\Gamma^{\prime}}}^{A \otimes B} \\
&= \theta_{\quot{h}{\Gamma}}^{A,B} \circ F(\hGamma)\tau_{f}^{A \otimes_Y B} \circ \left(\phi_{\quot{h}{\Gamma},\quot{f}{Y}}^{A \otimes B}\right)^{-1} \circ \phi_{\quot{f}{X},\quot{h}{\Gamma^{\prime}}}^{A \otimes B}  \\
&= \theta_{\quot{h}{\Gamma}}^{A,B} \circ \tau_{f}^{h^{\ast}(A \otimes_Y B)}.
\end{align*}
Thus it follows that $\theta_{\quot{h}{\Gamma}}^{A,B} \circ \tau_{f}^{h^{\ast}(A \otimes_Y B)} = \tau_f^{h^{\ast}A \otimes_X h^{\ast}B}$ and so $\Theta$ is an $F_G(X)$-morphism. With the remarks we made earlier about the naturality of $\Theta$ and the fact that $\Theta$ is an isomorphism, we see that $h^{\ast}$ is monoidal.
\end{proof}
\begin{corollary}\label{Cor: Section 3: Additive pullbacks}
If $F$ is a simultaneous pre-equivariant pseudofunctor over $X$ and $Y$ and if $F$ is additive, then the pullback $h^{\ast}:F_G(Y) \to F_G(X)$ is additive as well.
\end{corollary}

We now discuss the existence of pushforward functors $h_{\ast}:F_G(X) \to F_G(Y)$. In particular, we will give a formal condition for when we can get pushforwards based on the existence of pushforwards $\left(\hGamma\right)_{\ast}:F(\quot{X}{\Gamma}) \to F(\quot{Y}{\Gamma})$. As a consequence, we will give the existence of pushforwards $Rh_{\ast}:\DbeqQl{X} \to \DbeqQl{Y}, Rh_{\ast}:D_G^{b}(\QCoh(X)) \to D^{b}(\QCoh(Y))$ and $Rh_{\ast}:D^{b}_G(\TShv(X,\text{{\'e}t})) \to D^{b}_G(\TShv(X,\text{{\'e}t})$, as well as proper pushforwards $Rh_{!}:D^{b}_G(\QCoh(X)) \to D^b_G(\QCoh(Y)), Rh_!:\DbeqQl{X} \to \DbeqQl{Y}$, where $\TShv(X,\text{{\'e}t})$ is the equivariant derived category of torsion {\'e}tale sheaves on $X$ (with torsion coprime to the characteristic of $X$, of course).

\begin{definition}\label{Defn: Section 3.2: Descent Pushforwards}
Let $F$ be a simultaneous pre-equivariant pseudofunctor over $X$ and $Y$ and let $h \in \GVar(X,Y)$. Then we say that $F$ admits {descent pushforwards along $h$}\index{Descent Pushforwards} if there is a collection of functors $\left(\hGamma\right)_{\ast}:F(\XGamma) \to F(\YGamma)$ such that for all $f \in \Sf(G)_1$, say with $f:\Gamma \to \Gamma^{\prime}$, there is a natural isomorphism
\[
\begin{tikzcd}
F(\XGammap) \ar[r, bend left = 30, ""{name = U}]{}{(\hGamma)_{\ast} \circ F(\fX)} \ar[r, bend right = 30, swap, ""{name = L}]{}{F(\fY) \circ (\hGammap)_{\ast}} & F(\YGamma) \ar[from = U, to = L, Rightarrow, shorten <= 4pt, shorten >= 4pt]{}{\quot{h}{f}_{\ast}}
\end{tikzcd}
\]
for which we have a commuting $2$-cell
\[
\begin{tikzcd}
F(\XGammap) \ar[r, ""{name = U}]{}{F(\fX)} \ar[d, swap]{}{(\hGammap)_{\ast}} & F(\XGamma) \ar[d]{}{(\hGamma)_{\ast}} \\
F(\YGammap) \ar[r, swap, ""{name = L}]{}{F(\fY)} & F(\YGamma) \ar[from = U, to = L, Rightarrow, shorten <= 4pt, shorten >=4pt]{}{\quot{h}{f}_{\ast}}
\end{tikzcd}
\]
and for which the pasting diagram, for any composable arrows $\Gamma \xrightarrow{f} \Gamma^{\prime} \xrightarrow{g} \Gamma^{\prime\prime}$ in $\Sf(G)$,
\[
\begin{tikzcd}
F(\XGammapp) \ar[r, ""{name = LU}]{}{F(\gX)} \ar[d, swap]{}{(\hGammapp)_{\ast}} & F(\XGammap) \ar[r, ""{name = RU}]{}{F(\fX)}  \ar[d]{}[description]{(\hGammap)_{\ast}} & F(\XGamma) \ar[d]{}{(\hGamma)_{\ast}} \\
F(\YGammapp) \ar[r, swap, ""{name = LB}]{}{F(\gY)} \ar[rr, swap, bend right = 50, ""{name = BB}]{}{F(\gY \circ \fY)} & F(\YGammap) \ar[r, swap, ""{name = RB}]{}{F(\fY)} & F(\YGamma) \ar[ from = LU, to = LB, Rightarrow, shorten <= 4pt, shorten >=4pt]{}{\quot{h}{g}_{\ast}} \ar[from = RU, to = RB, Rightarrow, shorten <= 4pt, shorten >=4pt]{}{\quot{h}{f}_{\ast}} \ar[from = 2-2, to = BB, Rightarrow, shorten <= 4pt, shorten >=4pt]{}{\phi_{\quot{f}{Y}, \quot{g}{Y}}}
\end{tikzcd}
\]
is equivalent to:
\[
\begin{tikzcd}
& F(\XGammap) \ar[dr]{}{F(\fX)}\\
F(\XGammapp) \ar[d, swap]{}{(\hGammapp)_{\ast}} \ar[ur]{}{F(\gX)} \ar[rr, ""{name = U}]{}[description]{F(\gX \circ \fX)} & & F(\XGamma) \ar[d]{}{(\hGamma)_{\ast}} \\
F(\YGammapp) \ar[rr, swap, ""{name = B}]{}{F(\gY \circ \fY)} & & F(\YGamma) \ar[from = 1-2, to = U, Rightarrow, shorten <= 4pt, shorten >=4pt]{}{\phi_{\quot{f}{X}, \quot{g}{X}}} \ar[from = U, to = B, Rightarrow, shorten <= 4pt, shorten >=4pt]{}{\quot{h}{g \circ f}_{\ast}}
\end{tikzcd}
\]

\end{definition}
\begin{remark}
The above definition essentially says that the collection of pushforwards $\left(\hGamma\right)_{\ast}$ and isomorphisms $\quot{h}{f}_{\ast}$ form a pseudonatural transformation on $F$; however, because this is only defined on the pieces of the pseudofunctor $F$ that come from the essential image of $\quo_X$ and $\quo_Y$, and because $\SfResl_G(X)$ is not generically equal to $\SfResl_G(Y)$, we cannot na{\"i}vely use the language of pseudofunctors. However, proofs involving these pushforwards carry over mutatis mutandis from the pseudonatural cases (cf.\@ Theorem \ref{Thm: Section 3: Psuedonatural trans lift to equivariant functors}, for example).
\end{remark}
\begin{Theorem}\label{Thm: Section 3: Existence of pushforwards}
Let $F$ be a pre-equivariant pseudofunctor over $X$ and $Y$, let $h \in \GVar(X,Y)$, and assume that $F$ admits descent pushforwards along $h$. Then there is a pushforward functor $h_{\ast}:F_G(X) \to F_G(Y)$ induced by the assignments
\[
\lbrace \quot{A}{\Gamma} \; | \; \Gamma \in \Sf(G)_0\rbrace \mapsto \lbrace (\hGamma)_{\ast}\quot{A}{\Gamma} \; | \; \Gamma \in \Sf(G)_0\rbrace
\]
on object sets with transition isomorphisms
\[
T_{h_{\ast}A} := \lbrace \quot{h}{\Gamma}_{\ast}(\tau_f^A) \circ (\quot{(h_{\ast})}{f}_{\quot{A}{\Gamma^{\prime}}})^{-1} \; | \; f \in \Sf(G)_0 \rbrace
\]
and on morphisms by the assignment
\[
\lbrace \quot{\rho}{\Gamma} \; | \; \Gamma \in \Sf(G)_0\rbrace \mapsto \lbrace (\hGamma)_{\ast}\quot{\rho}{\Gamma} \; | \; \Gamma \in \Sf(G)_0\rbrace.
\]
\end{Theorem}
\begin{proof}
While formally different, the conditions on $F$ admitting descent pushforwards is essentially the same as if the collections of $(\hGamma)_{\ast}$ and $\quot{h}{f}_{\ast}$ form a partial pseudonatural transformation 
\[
F|_{\quo^{\op}_X(\SfResl_G(X)^{\op})} \to F|_{\quo_Y^{\op}(\SfResl_G(Y)^{\op})}.
\] 
From here we can apply the proof of Theorem \ref{Thm: Section 3: Psuedonatural trans lift to equivariant functors} mutatis mutands to give the existence of the functor $h_{\ast}:F_G(X) \to F_G(Y)$.
\end{proof}
\begin{example}
Let $G$ be a smooth algebraic group and let $X$ and $Y$ be $G$-varieites with $h \in \GVar(X,Y)$. Then if $F = \Ab\Shv(-,\text{{\'e}t})$ is the Abelian {\'e}tale sheaf simultaneous pre-equivariant pseudofunctor, $F$ admits descent pushforwards along $h$ with pusforward functors the usual {\'e}tale direct image functors (because each direct image $\hGamma_{\ast}$ commutes with smooth pullback). Similarly, if $F = \DbQl{-}$ then $F$ admits descent pushforwards along both the functors $\hGamma_{\ast} := R(\hGamma)_{\ast}$ (derived direct image) and $\hGamma_{\ast} := R(\hGamma)_!$ (derived proper direct image).
\end{example}
\begin{remark}\label{Remark: Section 3.2: Equivariant proper pushforward}
	It is worth noting that while Definition \ref{Defn: Section 3.2: Descent Pushforwards} suggestively writes the family of functors as $(\hGamma)_{\ast}$ to bring to mind the pushforwards that come from a geometric morphisms. However, the formalism is significantly more flexible than just providing that particular example. For instance, we can use this formalism (together with Theorem \ref{Thm: Section 3: Existence of pushforwards}) to provide an equivariant proper pushforward, when given a morphism $h \in \GVar(X,Y)$,
	\[
	Rh_{!}:\DbeqQl{X} \to \DbeqQl{Y}
	\]
	as in the example above.
\end{remark}
\begin{proposition}\label{Prop: Section 3: Pushforwards preserve stuff}
Let $F$ be a simultaneous pre-equivariant pseudofunctor over $X$ and $Y$ and let $h \in \GVar(X,Y)$. Assume also that $F$ admits descent pushforwards along $h$. Assume that $I$ is an index category such that for all $\Gamma \in \Sf(G)_0$ there are diagram functors $\quot{d}{\Gamma}:I \to F(\YGamma)$ and assume that $d:I \to F_G(Y)$ is a diagram functor for which $d(i) = \lbrace \quot{d}{\Gamma}(i) \; | \; \Gamma \in \Sf(G)_0 \rbrace$. Assume also that each category $F(\YGamma)$ admits limits and colimits of shape $I$. Then if each functor $(\hGamma)_{\ast}$ satisfies any of the following properties for all $\Gamma \in \Sf(G)_0$, so does $h_{\ast}$:
\begin{enumerate}
	\item $(\hGamma)_{\ast}$ preserves all limits of specified shape $I$;
	\item $(\hGamma)_{\ast}$ preserves all colimits of sepcified shape $I$;
	\item The pre-equivariant pseudofunctor $F$ is monoidal and each $(\hGamma)_{\ast}$ and $\quot{h}{f}_{\ast}$ are monoidal functors and transformations.
\end{enumerate}
\end{proposition}
\begin{proof}
As in Proposition \ref{Prop: Section 3: Pullbacks are lex/rex/additive whenever the F(h) are}, $(2)$ is dual to $(1)$, so it suffices to prove $(1)$. For $(1)$, however, this follows mutatis mutandis to Proposition \ref{Prop: Section 3: Change of fibre equivariant functor preserves limit of a specific shape} by using that descent pushfowards along $h$ determine pseudonaturally varying families of functors. For $(3)$ we note that this follows mutatis mutandis to Propoisition \ref{Prop: Section 3: Monoidal pseudonats give rise to monoidal equivariant functors}.
\end{proof}
\begin{corollary}\label{Cor: Section 3: Additive pushforwards}
If each functor $(\hGamma)_{\ast}$ is additive, then so is $h_{\ast}$.
\end{corollary}

We now give the Change of Space analogue of Theorem \ref{Thm: Section 3: Gamma-wise adjoints lift to equivariant adjoints}. In particular, if a simulatneous pre-equivariant pseudofunctor has descent pushforwards along a morphism $h \in \GVar(X,Y)$ and if ther are adjoints $F(\hGamma) \dashv (\hGamma)_{\ast}$ for all $\Gamma \in \Sf(G)_0$, then we get an adjunction $h^{\ast} \dashv h_{\ast}:F_G(Y) \to F_G(X)$. This should be seen as simply saying that if we have fibre-wise adjunctions between $G$-equivariant categories that are determined by $G$-equivariant descent information, then these determine a $G$-equivariant adjunction. The proof, analogously to the case of Theorem \ref{Thm: Section 3: Existence of pushforwards}, can be deduced mutatis mutandis from the proof of Theorem \ref{Thm: Section 3: Gamma-wise adjoints lift to equivariant adjoints} using the same argument, so we only prove how to adapt the Change of Space situation to what was done in the Change of Fibre situation.

\begin{Theorem}\label{Theorem: Section 3: Change of space adjunctions}
Assume that $F$ is a simultaneous pre-equivariant pseudofunctor over $X$ and $Y$, that $h \in \GVar(X,Y)$, and that $F$ admits descent pushforwards along $h$. If there are adjoints $F(\hGamma) \dashv (\hGamma)_{\ast}$ for all $\Gamma \in \Sf(G)_0$, then there is an adjunction:
\[
\begin{tikzcd}
F_G(Y) \ar[r, bend left = 30, ""{name = U}]{}{h^{\ast}} & F_G(X) \ar[l, bend left = 30, ""{name = L}]{}{h_{\ast}} \ar[from = U, to = L, symbol = \dashv]
\end{tikzcd}
\]
If instead there are adjoints $(\hGamma)_{\ast} \dashv F(\hGamma)$ for all $\Gamma \in \Sf(G)_0$, then there is an adjunction
\[
\begin{tikzcd}
F_G(Y) \ar[r, bend left = 30, ""{name = U}]{}{h_{\ast}} & F_G(X) \ar[l, bend left = 30, ""{name = L}]{}{h^{\ast}} \ar[from = U, to = L, symbol = \dashv]
\end{tikzcd}
\]
\end{Theorem}
\begin{proof}
The second statement is dual to the first, so it suffices to prove only the first statement of the theorem. We show how to adapt this situation to the formalism described in Theorem \ref{Thm: Section 3: Gamma-wise adjoints lift to equivariant adjoints}. For this, let $\quot{\eta}{\Gamma}:\id_{F(\YGamma)} \Rightarrow (\hGamma)_{\ast} \circ F(\hGamma)$ and $\quot{\epsilon}{\Gamma}:F(\hGamma) \circ (\hGamma)_{\ast} \Rightarrow \id_{F(\YGamma)}$ be the unit and counit, respectively, of the adjunction $F(\hGamma) \dashv (\hGamma)_{\ast}$. Begin by observing that for any $f \in \Sf(G)_1$ with $f \in \Sf(G)(\Gamma, \Gamma^{\prime})$, we have that the composite
\[
\begin{tikzcd}[cramped]
\big((\hGamma)_{\ast} \circ F(\hGamma) \circ F(\fY)\big)(\quot{A}{\Gamma^{\prime}}) \ar[rr]{}{(\hGamma)_{\ast}\phi_{\quot{h}{\Gamma},\quot{f}{Y}}^{\quot{A}{\Gamma^{\prime}}}} & & \big((\hGamma)_{\ast} \circ F(\quot{\overline{k}}{f})\big)(\quot{A}{\Gamma^{\prime}}) \ar[d]{}{(\hGamma)_{\ast}\left(\phi_{\quot{f}{X}, \quot{h}{\Gamma^{\prime}}}^{\quot{A}{\Gamma^{\prime}}}\right)^{-1}} \\
\big(F(\fY) \circ (\hGammap)_{\ast} \circ F(\hGammap)\big)(\quot{A}{\Gamma^{\prime}}) & & \big((\hGamma)_{\ast} \circ F(\fX) \circ F(\hGammap)\big)(\quot{A}{\Gamma^{\prime}}) \ar[ll]{}{\quot{h}{f}_{\ast}^{F(\hGammap)\quot{A}{\Gamma^{\prime}}}}
\end{tikzcd}
\]
plays the Change of Space role of the natural isomorphism $\quot{(\beta \circ \alpha)}{f}$ at $\quot{A}{\Gamma^{\prime}}$ in Theorem \ref{Thm: Section 3: Gamma-wise adjoints lift to equivariant adjoints}. For the sake of brevity, define
\[
\quot{\alpha}{f} := \quot{h}{f}_{\ast} \circ (\hGamma)_{\ast}\left(\phi_{\quot{f}{X},\quot{h}{\Gamma^{\prime}}}\right)^{-1} \circ (\hGamma)_{\ast}\left(\phi_{\quot{h}{\Gamma},\quot{f}{Y}}\right).
\]
Proceeding now as in the proof of Theorem \ref{Thm: Section 3: Gamma-wise adjoints lift to equivariant adjoints} we get that the diagram
\[
\xymatrix{
F(\fY)(\quot{A}{\Gamma^{\prime}}) \ar[rr]^-{\quot{\eta_{F(\fY)\quot{A}{\Gamma^{\prime}}}}{\Gamma}} \ar@{=}[d] & & \big((\hGamma)_{\ast} \circ F(\hGamma) \circ F(\fY)\big)(\quot{A}{\Gamma^{\prime}}) \ar[d]^{\quot{\alpha_{\quot{A}{\Gamma^{\prime}}}}{f}} \\
F(\fY)(\quot{A}{\Gamma^{\prime}}) \ar[rr]_-{F(\fY)\quot{\eta_{\quot{A}{\Gamma^{\prime}}}}{\Gamma^{\prime}}} & & \big(F(\fY) \circ (\hGammap)_{\ast} \circ F(\hGammap)\big)(\quot{A}{\Gamma^{\prime}})
}
\]
commutes for any $\quot{A}{\Gamma^{\prime}} \in F(\YGammap)_0$, and dually for the diagram induced  by the counits. We can then deduce that the triangle identities on $F_G(Y) \to F_G(X)$ hold by using the piecewise triangle identities induced by $F(\hGamma) \dashv (\hGamma)_{\ast}$ and then using the compatibliity induced by the unit/counit diagrams to descend through the transition isomorphisms $\tau_f$ as necessary.
\end{proof}
\begin{remark}
This last step in the above proof is induced by Lemma \ref{Lemma: Section 3: 2-functors preserve adjoints} in the case of Theorem \ref{Thm: Section 3: Gamma-wise adjoints lift to equivariant adjoints} and does not need to be verified explicitly.
\end{remark}

Before moving on from these change of space functors, it will be useful to have functors $F_G(Y) \to F_G(X)$ which are not induced by the (simultaneous pre-equivariant) pseudofunctor $F$. A motivating example for this is the construction of an equivariant version of the functor $h^{!}:\DbeqQl{Y} \to \DbeqQl{X}$ when given a morphism $h \in \GVar(X,Y)$. It may be tempting to try to deduce this from Theorem \ref{Thm: Section 3: Existence of pullback functors} by taking our pseudofunctor $F$ to be given by $\Gamma \times X \mapsto \DbQl{\quot{X}{\Gamma}}$ and our morphisms by $f \times \id_X\mapsto \overline{f}^{!}$. Unfortunately, this does not give the same equivariant category. Indeed, even though the morphism $\overline{f}:\quot{X}{\Gamma} \to \quot{X}{\Gamma^{\prime}}$ is smooth, we only have that
\[
\overline{f}^{!} = \overline{f}^{\ast}[2d](d),
\]
where $d$ is the relative fibre dimension of $\overline{f}$. Because our $\Sf(G)$-varieties have morphisms with nontrivial fibre dimension, generically we have $\overline{f}^{\ast}[2d](d) \ne \overline{f}^{\ast}$, which in turn implies that $F_G(X) \ne \DbeqQl{X}$ (we call this category $F_G(X) := {}^{!}\DbQl{X}$). Instead, to develop $h^{!}$ we must take a different approach which mirrors how we built our equivariant pushforwards: We need to have a notion of descent pullbacks along a morphism $h$ to deduce the existence of a functor $F_G(Y) \to F_G(X)$.
\begin{definition}\label{Defn: Section 3.2: Descent Pullbacks}
Let $h \in \GVar(X,Y)$ and let $F$ be a simultaneous pre-equivariant pseudofunctor over $X$ and $Y$. We then say that $F$ admits {descent pullbacks along $h$}\index{Descent Pullbacks} if there is a collection of functors $(\hGamma)^{\ast}:F(\YGamma) \to F(\XGamma)$ such that for all $f \in \Sf(G)_1$ there is a natural isomorphism
\[
\begin{tikzcd}
F(\YGammap) \ar[r, bend left = 30, ""{name = U}]{}{(\hGamma)^{\ast} \circ F(\fY)} \ar[r, bend right = 30, swap, ""{name = L}]{}{F(\fX) \circ (\hGammap)^{\ast}} & F(\XGamma) \ar[from = U, to = L, Rightarrow, shorten <= 4pt, shorten >= 4pt]{}{\quot{h}{f}^{\ast}}
\end{tikzcd}
\]
for which we have a commuting $2$-cell
\[
\begin{tikzcd}
F(\YGammap) \ar[r, ""{name = U}]{}{F(\fY)} \ar[d, swap]{}{(\hGammap)^{\ast}} & F(\YGamma) \ar[d]{}{(\hGamma)^{\ast}} \\
F(\XGammap) \ar[r, swap, ""{name = L}]{}{F(\fX)} & F(\XGamma) \ar[from = U, to = L, Rightarrow, shorten <= 4pt, shorten >=4pt]{}{\quot{h}{f}^{\ast}}
\end{tikzcd}
\]
and for which the pasting diagram
\[
\begin{tikzcd}
F(\YGammapp) \ar[r, ""{name = LU}]{}{F(\gY)} \ar[d, swap]{}{(\hGammapp)^{\ast}} & F(\YGammap) \ar[r, ""{name = RU}]{}{F(\fY)}  \ar[d]{}[description]{(\hGammap)^{\ast}} & F(\YGamma) \ar[d]{}{(\hGamma)^{\ast}} \\
F(\XGammapp) \ar[r, swap, ""{name = LB}]{}{F(\gX)} \ar[rr, swap, bend right = 50, ""{name = BB}]{}{F(\gX \circ \fX)} & F(\XGammap) \ar[r, swap, ""{name = RB}]{}{F(\fY)} & F(\XGamma) \ar[ from = LU, to = LB, Rightarrow, shorten <= 4pt, shorten >=4pt]{}{\quot{h}{g}^{\ast}} \ar[from = RU, to = RB, Rightarrow, shorten <= 4pt, shorten >=4pt]{}{\quot{h}{f}^{\ast}} \ar[from = 2-2, to = BB, Rightarrow, shorten <= 4pt, shorten >=4pt]{}{\phi_{\quot{f}{X}, \quot{g}{X}}}
\end{tikzcd}
\]
is equivalent to:
\[
\begin{tikzcd}
& F(\YGammap) \ar[dr]{}{F(\fY)}\\
F(\YGammapp) \ar[d, swap]{}{(\hGammapp)^{\ast}} \ar[ur]{}{F(\gY)} \ar[rr, ""{name = U}]{}[description]{F(\gY \circ \fY)} & & F(\YGamma) \ar[d]{}{(\hGamma)^{\ast}} \\
F(\XGammapp) \ar[rr, swap, ""{name = B}]{}{F(\gX \circ \fX)} & & F(\XGamma) \ar[from = 1-2, to = U, Rightarrow, shorten <= 4pt, shorten >=4pt]{}{\phi_{\quot{f}{Y}, \quot{g}{Y}}} \ar[from = U, to = B, Rightarrow, shorten <= 4pt, shorten >=4pt]{}{\quot{h}{g \circ f}^{\ast}}
\end{tikzcd}
\]
\end{definition}
From here we can proceed to prove that these descent pullbacks give rise to equivariant pullback functors by following the proof of Theorem \ref{Thm: Section 3: Existence of pushforwards} and proceeding mutatis mutandis.
\begin{Theorem}\label{Thm: Section 3.2: Existence of descent pullbacks}
Let $F$ be a pre-equivariant pseodufunctor over $X$ and $Y$, let $h \in \GVar(X,Y)$, and assume that $F$ admits descent pullbacks along $h$. Then there is a pullback functor $h^{\ast}:F_G(X) \to F_G(Y)$ induced by the assignments
\[
\lbrace \AGamma \; | \; \Gamma \in \Sf(G)_0 \rbrace \mapsto \lbrace (\hGamma)^{\ast}\AGamma \; | \; \Gamma \in \Sf(G)_0 \rbrace
\]
on objects with transition isomorphisms given by
\[
T_{h^{\ast}A} \lbrace \quot{h}{\Gamma}^{\ast}((\tau_f^{A})) \circ (\quot{(h^{\ast})}{f}_{\AGammap})^{-1} \; | \; f \in \Sf(G)_0\rbrace,
\] 
while the assignment is given by
\[
\lbrace \quot{\rho}{\Gamma} \; | \; \Gamma \in \Sf(G)_0 \rbrace \mapsto \lbrace (\hGamma)^{\ast}\quot{\rho}{\Gamma} \; | \; \Gamma \in \Sf(G)_0 \rbrace
\]
on morphisms.
\end{Theorem}
\begin{proof}
As in the proof of Theorem \ref{Thm: Section 3: Existence of pushforwards}, we find that descent pullbacks essentially define a partial pseudonatural transformation 
\[
F|_{\quo_{X}^{\op}(\SfResl_G(Y)^{\op})} \to F|_{\quo_Y^{\op}(\SfResl_G(X)^{\op})}.
\]
We proceed then as in the proof of Theorem \ref{Thm: Section 3: Existence of pushforwards} by following the proof of Theorem \ref{Thm: Section 3: Psuedonatural trans lift to equivariant functors} mutatis mutandis to give rise to $h^{\ast}:F_G(Y) \to F_G(X)$.
\end{proof}
\begin{definition}\label{Defn: Change of Space Functor}\index{Equivariant Functor! Change of Space}
An equivariant functor $R:F_G(X) \to F_G(Y)$ is a Change of Space functor if:
\begin{itemize}
	\item There is an equivariant morphism $h \in \GVar(X,Y)$, $F$ admits descent pushforwards along $h$, and $R = h_{\ast}$;
	\item There is an equivariant morphism $h \in \GVar(Y,X)$, $F$ admits descent pullbacks along $h$, and $R = h^{\ast}$.
\end{itemize}
\end{definition}
\begin{remark}
Theorem \ref{Thm: Section 3: Existence of pullback functors} can be seen as a special case of Theorem \ref{Thm: Section 3.2: Existence of descent pullbacks} above. In fact, the significant amount of work in the construction can essentially be viewed as showing that for any simultaneous pre-equivariant pseudofunctor $F$ over $X$ and $Y$, the morphisms $F(\hGamma)$ which bridge the data over $Y$ to the data over $X$ give us a family of descent pullbacks.
\end{remark}
\begin{example}
In light of the remark above, we see that for any simultaneous pre-equivariant pseudofunctor $F$ over $X,Y$ and for any $h \in \GVar(X,Y)$, there are descent pullbacks along $h$ given by $\hGamma^{\ast} = F(\hGamma)$. For less pathological examples, take $F = \DbQl{-}$ and set $\hGamma^{\ast} := \hGamma^{!}:\DbQl{\YGamma} \to \DbQl{\XGamma}$. Then $F$ admits descent pullbacks along $h$ with pullback functors $\hGamma^{!}$.
\end{example}
\begin{proposition}\label{Prop: Section 3.2: Pullbacks preserve stuff}
	Let $F$ be a simultaneous pre-equivariant pseudofunctor over $X$ and $Y$ and let $h \in \GVar(X,Y)$. Assume also that $F$ admits descent pullbacks along $h$. Assume $I$ is an index category such that for all $\Gamma \in \Sf(G)_0$ there are diagram functors $\quot{d}{\Gamma}:I \to F(\YGamma)$, and assume $d:I \to F_G(Y)$ is a diagram functor for which $d(i) = \lbrace \quot{d}{\Gamma}(i) \; | \; \Gamma \in \Sf(G)_0 \rbrace$. Assume also that each category $F(\YGamma)$ admits limits and colimits of shape $I$. Then if each functor $(\hGamma)^{\ast}$ satisfies any of the following properties for all $\Gamma \in \Sf(G)_0$, so does $h^{\ast}$:
	\begin{enumerate}
		\item $(\hGamma)^{\ast}$ preserves all limits of specified shape $I$;
		\item $(\hGamma)^{\ast}$ preserves all colimits of sepcified shape $I$;
		\item The pre-equivariant pseudofunctor $F$ is monoidal, each of $(\hGamma)^{\ast}$ and $\quot{h}{f}^{\ast}$ are monoidal functors and transformations.
	\end{enumerate}
\end{proposition}
\begin{proof}
Proceed mutatis mutandis as in the prooof of Proposition \ref{Prop: Section 3: Pushforwards preserve stuff}.
\end{proof}
\begin{corollary}\label{Cor: Section 3.2: Additive descent pullbacks}
Let $F$ be a simultaneous pre-equivariant pseudofunctor that admits descent pullbacks along $h \in \GVar(X,Y)$. Then if each functor $(\hGamma)^{\ast}$ is additive, so is $h^{\ast}$.
\end{corollary}

Before moving on to discuss how Change of Space functors interact with Change of Fibre functors, we would like to provide an adjunction lifting theorem for descent pullbacks and pushforwards. In particular, we would like to lift Theorem \ref{Theorem: Section 3: Change of space adjunctions} to the case where the pullback $h^{\ast}$ is induced not just by the pseudofunctor $F$, but instead by a descent pullback along $F$.

\begin{Theorem}\label{Theorem: Section 3.2: Descent adjunctions lift to equivariant ajdunction}
Let $F$ be a  simultaneous pre-equivariant pseudofunctor over $X$ and $Y$ and assume that $F$ admits both descent pullbacks and pushforwards along $h \in \GVar(X,Y)$. If there are adjoints
\[
\begin{tikzcd}
F(\YGamma) \ar[r, bend left = 30, ""{name = U}]{}{(\hGamma)^{\ast}} & F(\XGamma) \ar[l, bend left = 30, ""{name = L}]{}{(\hGamma)_{\ast}} \ar[from = U, to = L, symbol = \dashv]
\end{tikzcd}
\]
or
\[
\begin{tikzcd}
F(\XGamma) \ar[r, bend left = 30, ""{name = U}]{}{(\hGamma)_{\ast}} & F(\YGamma) \ar[l, bend left = 30, ""{name = L}]{}{(\hGamma)^{\ast}} \ar[from = U, to = L, symbol = \dashv]
\end{tikzcd}
\]
for all $\Gamma \in \Sf(G)_0$, then there is an adjunction
\[
\begin{tikzcd}
F_G(Y) \ar[r, bend left = 30, ""{name = U}]{}{h^{\ast}} & F_G(X) \ar[l, bend left = 30, ""{name = L}]{}{h_{\ast}} \ar[from = U, to = L, symbol = \dashv]
\end{tikzcd}
\]
or
\[
\begin{tikzcd}
F_G(X) \ar[r, bend left = 30, ""{name = U}]{}{h_{\ast}} & F_G(Y) \ar[l, bend left = 30, ""{name = L}]{}{h^{\ast}} \ar[from = U, to = L, symbol = \dashv]
\end{tikzcd}
\]
\end{Theorem}
\begin{proof}
The case where $(\hGamma)_{\ast} \dashv (\hGamma)^{\ast}$ is dual to the other, so it suffices to prove the theorem only for the case in which $(\hGamma)^{\ast} \dashv (\hGamma)_{\ast}$ for all $\Gamma \in \Sf(G)_0$. As in the proof of Theorem \ref{Theorem: Section 3: Change of space adjunctions}, our strategy is to lift our current situation to the formalism of Theorem \ref{Thm: Section 3: Gamma-wise adjoints lift to equivariant adjoints}. Following the proof of Theorem \ref{Theorem: Section 3: Change of space adjunctions} and letting $\quot{\eta}{\Gamma}:\id_{F(\YGamma)} \Rightarrow (\hGamma)_{\ast} \circ (\hGamma)^{\ast}$ and $\quot{\epsilon}{\Gamma}:(\hGamma)^{\ast} \circ (\hGamma)_{\ast} \Rightarrow \id_{F(\XGamma)}$ be the units and counits of adjunction, respectively, we get that in this case the role of $\quot{\alpha}{f}$ in the proof of Theorem \ref{Thm: Section 3: Gamma-wise adjoints lift to equivariant adjoints} is played by the composite:
\[
\xymatrix{
\left((\hGamma)_{\ast} \circ (\hGamma)^{\ast} \circ F(\fY)\right)(\AGammap) \ar[rr]^-{(\hGamma)_{\ast}(\quot{\overline{h}}{f})^{\ast}_{\AGammap}} & & \left((\hGamma)_{\ast} \circ F(\fX) \circ (\hGammap)^{\ast}\right)(\AGammap) \ar[d]^-{(\quot{h}{f}_{\ast})_{(\hGammap)^{\ast}\AGammap}} \\
& & \left(F(\fX) \circ (\hGammap)_{\ast} \circ (\hGammap)^{\ast}\right)(\AGammap)
}
\]
From here setting
\[
\quot{\alpha}{f} := (\quot{h}{f}_{\ast})_{(\hGammap)^{\ast}(-)} \circ (\hGamma)_{\ast}(\quot{\overline{h}}{f})^{\ast} = \left(\quot{h}{f} \star (\hGammap)^{\ast}\right) \circ \left((\hGamma)_{\ast} \star \quot{h}{g}\right),
\]
where $\star$ denotes horizontal composition, and proceeding as in the proof of Theorem \ref{Theorem: Section 3: Change of space adjunctions} allows us to deduce that the desired adjunctions hold.
\end{proof}


We now move to discuss how to get equivariant functors which simultaneously change fibres and spaces. That is, we discuss equivariant functors of the form $E_G(Y) \to F_G(X)$ and $F_G(X) \to E_G(Y)$ where $F$ and $E$ are simultaneous pre-equivariant pseudofunctors over $X$ and $Y$, where $h \in \GVar(X,Y)$, and where we have have a pseudonatrual transformation $\alpha:F \to E$. Unfortunately, there is not much that we can say in complete generality, as to ask for a diagram
\[
\xymatrix{
F_G(X) \ar[r]^-{\underline{\alpha}_X} \ar[d]_{h^{\ast}_F} & E_G(X) \ar[d]^{h^{\ast}_E} \\
F_G(Y) \ar[r]_-{\underline{\alpha}_Y} & E_G(Y)
}
\]
to commute, even up to isomorphism, is to ask for natural isomorphisms
\[
\theta_{\Gamma}:(\quot{\alpha_Y}{\Gamma} \circ {\hGamma_F}^{\ast}) \xRightarrow{\cong} (\hGamma_E^{\ast} \circ \quot{\alpha_X}{\Gamma})
\]
for all $\Gamma \in \Sf(G)_0$ which behave sufficiently like a modification to commute between the coherences of both $\alpha$ and the descent pullbacks along $h$ (or pushforwards along $h$). However, when these conditions do arise, we can derive the existence of a natural isomorphism $(\underline{\alpha}_Y \circ h_F^{\ast}) \xRightarrow{\cong} (h_E^{\ast} \circ \underline{\alpha}_X).$ The conditions that describe how these commutativity isomorphisms must arise are quite technical (we need to explain what ``sufficiently like a modification'' means precisely) so we need to examine them in detail.

Let us describe how and why the technical conditions to which we alluded arise. Assume that we have two simultaneous pre-equivariant pseduofunctors $E,F:\SfResl_{G}(X)^{\op} \to\fCat$ over $X$ and $Y$ which have a pseudonatural transformation $\alpha:F \Rightarrow E$ between them as in the $2$-cell:
\[
\begin{tikzcd}
\Var_{/\Spec K}^{\op} \ar[rr, bend left = 30, ""{name = U}]{}{F} \ar[rr, bend right = 30, swap, ""{name = L}]{}{E} & & \fCat \ar[from = U, to = L, Rightarrow, shorten <= 4pt, shorten >= 4pt]{}{\alpha}
\end{tikzcd}
\]
Assume also that $h \in \GVar(X,Y)$ and that $E$ and $F$ admit descent pullbacks along $h$ (say with pullback families $\hGamma_F^{\ast}$ and $\hGamma_E^{\ast}$, respectively). Our technical conditions begin by asking for a natural isomorphism
\[
\begin{tikzcd}
F(\YGamma) \ar[r, ""{name = U}]{}{\hGamma_F^{\ast}} \ar[d, swap]{}{\quot{\alpha_Y}{\Gamma}} & F(\XGamma) \ar[d]{}{\quot{\alpha_X}{\Gamma}} \\
E(\YGamma) \ar[r, swap, ""{name = L}]{}{\hGamma_E^{\ast}} & E(\XGamma) \ar[from = U, to = L, Rightarrow, shorten >= 4pt, shorten <= 4pt]{}{\theta_{\Gamma}}
\end{tikzcd}
\]
for all $\Gamma \in \Sf(G)_0$. Unfortunately, this is not where the story ends. We need to know that these isomorphisms $\quot{\theta}{\Gamma}$ act sufficiently like modifications in the sense that they need to interact with the transition isomorphisms $\quot{\overline{h}^{\ast}}{f}$ and $\quot{\alpha}{f}$ in ways that allow us to produce an isomorphism in $E_G(X)$ which is $\Gamma$-locally built from the $\theta_{\Gamma}$. Writing down the condition to make
\[
\Theta_A = \lbrace \theta_{\Gamma}^{\AGamma} \; | \; \Gamma \in \Sf(G)_0 \rbrace
\]
an isomorphism in $E_G(X)$ whenever $A \in F_G(Y)_0$, we get that
\[
\xymatrix{
E(\fX)\big((\quot{\alpha_X}{\Gamma^{\prime}} \circ \hGammap_F^{\ast})(\AGammap)\big) \ar[rr]^-{E(\fX)\theta_{\Gamma^{\prime}}^{\AGammap}} \ar[d]_{\tau_f^{(\underline{\alpha} \circ h^{\ast})A}} & & E(\fX)\big((\hGammap_E^{\ast} \circ \quot{\alpha_Y}{\Gamma^{\prime}})(\AGammap)\big) \ar[d]^{\tau_{f}^{(h^{\ast} \circ \alpha)A}} \\
(\quot{\alpha_X}{\Gamma} \circ \hGamma_F^{\ast})(\AGamma) \ar[rr]_-{\theta_{\Gamma}^{\AGamma}} & & (\hGamma_E^{\ast} \circ \quot{\alpha_Y}{\Gamma})(\AGamma)
}
\]
must commute whenever $f \in \Sf(G)_1$ (with $\Dom f = \Gamma$ and $\Codom f = \Gamma^{\prime}$). Unwinding the definitions of the $\tau_f$ using Theorems \ref{Thm: Section 3: Psuedonatural trans lift to equivariant functors} and \ref{Thm: Section 3.2: Existence of descent pullbacks} gives that on one hand
\[
\theta_{\Gamma}^{\AGamma} \circ \tau_f^{(\underline{\alpha} \circ h^{\ast})A} = \theta_{\Gamma}^{A} \circ \quot{\alpha_X}{\Gamma}\left(\quot{h_F}{\Gamma}(\tau_f^A)\right) \circ \quot{\alpha_X}{\Gamma}\left(\quot{h_A}{f}\right)^{-1} \circ \left(\quot{\alpha_{h^{\ast}A}}{f}\right)^{-1}
\]
while on the other hand
\begin{align*}
&\tau_{f}^{(h^{\ast} \circ \alpha)A} \circ E(\fX)\theta_{\Gamma^{\prime}}^{\AGammap} \\
&= (\hGamma_E^{\ast} \circ \quot{\alpha_Y}{\Gamma})(\tau_f^A) \circ \hGamma_E^{\ast}\left(\quot{\alpha_A}{f}\right)^{-1} \circ \left(\quot{h_{\quot{\alpha(A)}{\Gamma^{\prime}}}}{f}\right)^{-1} \circ E(\fX)\theta_{\Gamma^{\prime}}^{\AGammap}.
\end{align*}
Because $\theta_{\Gamma}$ is a natural isomorphism $(\quot{\alpha_X}{\Gamma} \circ \hGamma_F^{\ast}) \Rightarrow (\hGamma_E^{\ast} \circ \quot{\alpha_Y}{\Gamma})$ we always have that
\[
\theta_{\Gamma}^{\AGamma} \circ \big(\quot{\alpha_X}{\Gamma} \circ \hGamma_F^{\ast}\big)(\tau_f^A) = \big(\hGamma_E^{\ast} \circ \quot{\alpha_Y}{\Gamma}\big)(\tau_f^{A}) \circ \theta_{\Gamma}^{F(\fY)A}
\]
which allows us to deduce that
\[
\theta_{\Gamma}^{\AGamma} \circ \tau_f^{(\underline{\alpha} \circ h^{\ast})A} = \big(\hGamma_E^{\ast} \circ \quot{\alpha_Y}{\Gamma}\big)(\tau_f^{A}) \circ \theta_{\Gamma}^{F(\fY)A} \circ \quot{\alpha_X}{\Gamma}\left(\quot{h_A}{f}\right)^{-1} \circ \left(\quot{\alpha_{h^{\ast}A}}{f}\right)^{-1}
\]
This tells us that it is both necessary and sufficient that the equality
\begin{align*}
&\theta_{\Gamma}^{F(\fY)A} \circ \quot{\alpha_X}{\Gamma}\left(\quot{h_A}{f}\right)^{-1} \circ \left(\quot{\alpha_{h^{\ast}A}}{f}\right)^{-1} \\
&= \hGamma_E^{\ast}\left(\quot{\alpha_A}{f}\right)^{-1} \circ \left(\quot{h_{\quot{\alpha(A)}{\Gamma^{\prime}}}}{f}\right)^{-1} \circ E(\fX)\theta_{\Gamma^{\prime}}^{\AGammap}
\end{align*}
hold in order for $\Theta_A = \lbrace \theta_{\Gamma}^{\AGamma} \; | \; \Gamma \in \Sf(G)_0\rbrace$ to be the natural isomorphism $(\underline{\alpha}_X \circ h_F^{\ast}) \xRightarrow{\cong} (h^{\ast}_E \circ \underline{\alpha}_Y)$. However, this equality has a more familiar incarnation. A straightforward manipulation allows us to express it as
\[
\quot{h_{\quot{\alpha_Y(A)}{\Gamma^{\prime}}}}{f} \circ \hGamma_E^{\ast}\left(\quot{\alpha_A}{f}\right) \circ \theta_{\Gamma}^{F(\fY)\AGamma} = E(\fX)\theta_{\Gamma^{\prime}}^{A} \circ \quot{\alpha_{h^{\ast}A}}{f} \circ \quot{\alpha_X}{\Gamma}\left(\quot{h_A}{f}\right),
\]
which is equivalent to saying that the diagram
\[
\xymatrix{
\left(\quot{\alpha_X}{\Gamma} \circ \hGamma_F^{\ast} \circ F(\fY)\right)(\AGammap) \ar[rr]^-{\theta_{\Gamma}^{F(\fY)A}} \ar[d]_{\quot{\alpha_X}{\Gamma}(\quot{h_A}{f})} & & \left(\hGamma_E^{\ast} \circ \quot{\alpha_Y}{\Gamma} \circ F(\fY)\right)(\AGammap) \ar[d]^{\hGamma_E^{\ast}(\quot{\alpha_{A}}{f})} \\
\left(\quot{\alpha_X}{\Gamma} \circ F(\fX) \circ \hGammap_F^{\ast}\right)(\AGammap) \ar[d]_{\quot{\alpha_{h^{\ast}A}}{f}} & & \left(\hGamma_E^{\ast} \circ E(\fY) \circ \quot{\alpha_Y}{\Gamma^{\prime}}\right)(\AGammap) \ar[d]^{\quot{h_{\quot{\alpha(A)}{\Gamma^{\prime}}}}{f}} \\
\left(E(\fX) \circ \quot{\alpha_X}{\Gamma^{\prime}} \circ \hGammap_F^{\ast}\right)(\AGammap) \ar[rr]_-{E(\fX)\theta_{\Gamma^{\prime}}^{A}} & & \left(E(\fX) \circ \hGammap_E^{\ast} \circ \quot{\alpha_Y}{\Gamma^{\prime}}\right)(\AGammap)
}
\]
commutes in $E(\XGamma)$. A comparison with this and the definition of what it means to be a modification allows us to conclude that this is what we mean when we say that the witness isomorphisms
\[
\theta_{\Gamma}^{A}:(\quot{\alpha_X}{\Gamma} \circ \hGamma^{\ast})A \xrightarrow{\cong} (\hGamma_E^{\ast} \circ \quot{\alpha_Y}{\Gamma})A
\]
need to be sufficiently like a modification. The situation follows mutatis mutandis for descent pushforwards along $h$ as well. These observations together provide a proof for the proposition below.

\begin{Theorem}\label{Prop: Necessary and Sufficeint Conditions for Change of Fibre commute with Change of Space}
Let $F$ and $E$ be simultaneous pre-equivariant pseudofunctors over $X$ and $Y$ and let $h \in \GVar(X,Y)$ and assume that there is a pseudonatural transformation $\alpha:F \Rightarrow E$ at the level
\[
\begin{tikzcd}
\Var_{/\Spec K}^{\op} \ar[rr, bend left = 30, ""{name = U}]{}{F} \ar[rr, bend right = 30, swap, ""{name = L}]{}{E} & & \fCat \ar[from = U, to = L, Rightarrow, shorten <= 4pt, shorten >= 4pt]{}{\alpha}
\end{tikzcd}
\]
Now assume that $E$ and $F$ both admit descent pullbacks along $h$. Then in order for the diagram
\[
\begin{tikzcd}
F_G(Y) \ar[r, ""{name = U}]{}{h_F^{\ast}} \ar[d, swap]{}{\underline{\alpha}_Y} & F_G(X) \ar[d]{}{\underline{\alpha}_X} \\
E_G(Y) \ar[r, swap, ""{name = L}]{}{h_E^{\ast}} & E_G(X) \ar[from = U, to = L, Rightarrow, shorten >= 4pt, shorten <= 4pt]{}{\Theta}
\end{tikzcd}
\]
commute up to natural isomorphism $\Theta$, it is necessary and sufficient to give natural isomorphisms 
\[
\theta_{\Gamma}:(\quot{\alpha_X}{\Gamma} \circ \hGamma_F^{\ast}) \xRightarrow{\cong} (\hGamma_E^{\ast} \circ \quot{\alpha_Y}{\Gamma})
\]
for which the identities
\[
\quot{h_{\quot{\alpha_Y(A)}{\Gamma^{\prime}}}}{f} \circ \hGamma_E^{\ast}\left(\quot{\alpha_A}{f}\right) \circ \theta_{\Gamma}^{F(\fY)\AGamma} = E(\fX)\theta_{\Gamma^{\prime}}^{A} \circ \quot{\alpha_{h^{\ast}A}}{f} \circ \quot{\alpha_X}{\Gamma}\left(\quot{h_A}{f}\right)
\]
hold for any $f \in \Sf(G)_1$. Dually, if $F$ and $E$ admit descent pushforwards, it is necessary and sufficent for the diagram
\[
\begin{tikzcd}
F_G(X) \ar[r, ""{name = U}]{}{h^F_{\ast}} \ar[d, swap]{}{\underline{\alpha}_X} & F_G(Y) \ar[d]{}{\underline{\alpha}_Y} \\
E_G(X) \ar[r, swap, ""{name = L}]{}{h^E_{\ast}} & E_G(X) \ar[from = U, to = L, Rightarrow, shorten >= 4pt, shorten <= 4pt]{}{\Theta}
\end{tikzcd}
\]
to commute up to the natural isomorphism $\Theta$ to give natural isomorphisms $\theta_{\Gamma}:(\quot{\alpha_Y}{\Gamma} \circ \hGamma_{\ast}^{F}) \xRightarrow{\cong} (\hGamma_{\ast}^{E} \circ \quot{\alpha_X}{\Gamma})$ for which the identities
\[
\quot{h_{\quot{\alpha_X(A)}{\Gamma^{\prime}}}}{f} \circ \hGamma^E_{\ast}\left(\quot{\alpha_A}{f}\right) \circ \theta_{\Gamma}^{F(\fY)\AGamma} = E(\fX)\theta_{\Gamma^{\prime}}^{A} \circ \quot{\alpha_{h_{\ast}A}}{f} \circ \quot{\alpha_Y}{\Gamma}\left(\quot{h_A}{f}\right)
\]
hold for any $f \in \Sf(G)_1$.
\end{Theorem}

The biggest reason for presenting the theorem in this way is that it allows us tools to deduce when we can provide the existence of such natural isomorphism based on how we move through the equivariant descent data that the pseudofunctors and pseudonatural transformations record. In view of Theorem \ref{Thm: Section 3: Psuedonatural trans lift to equivariant functors} and Lemma \ref{Lemma: Modifications give equivariant natural transformations}, it is perhaps unsurprising that we need the natural isomorphism $\Theta$ to behave sufficiently like modifications in order to assert that the diagram is natural in the descent data. However, this allows us to give a straightforward proof of Proposition \ref{Prop: Section 3: Change of  fibre then space is iso to changing space then fibre}, which says that when $F$ and $E$ are simultaneous pre-equivariant pseudofunctors over $X$ and $Y$ with a pseudonatural transformation $\alpha:F \Rightarrow E$, for any $h \in \GVar(X,Y)$ the descent pullbacks induced by the actions of the pseudofunctors give such a commuting diagram.

\begin{proposition}\label{Prop: Section 3: Change of  fibre then space is iso to changing space then fibre}
Let $F$ and $E$ be simultaneous pre-equivariant pseudofunctors over $X$ and $Y$ with a pseudonatrual transformation $\alpha:F \Rightarrow E$. Assume that $h \in \GVar(X,Y)$ and give both $F$ and $E$ the descent pullbacks induced by the functors $F(\hGamma)$ and $E(\hGamma)$. Then the diagram
\[
\begin{tikzcd}
F_G(Y) \ar[r, ""{name = U}]{}{h_F^{\ast}} \ar[d, swap]{}{\underline{\alpha}_Y} & F_G(X) \ar[d]{}{\underline{\alpha}_X} \\
E_G(Y) \ar[r, swap, ""{name = L}]{}{h_E^{\ast}} & E_G(X) \ar[from = U, to = L, Rightarrow, shorten <= 4pt, shorten >= 4pt]
\end{tikzcd}
\]
commutes up to natural isomorphism.
\end{proposition}
\begin{proof}
In light of Theorem \ref{Prop: Necessary and Sufficeint Conditions for Change of Fibre commute with Change of Space}, it suffices to prove that there are modification-like pseudonatural transformations. However, for any $\Gamma \in \Sf(G)_0$ the diagram
\[
\begin{tikzcd}
F(\YGamma) \ar[r, ""{name = U}]{}{F(\hGamma)} \ar[d, swap]{}{\quot{\alpha_Y}{\Gamma}} & F(\XGamma) \ar[d]{}{\quot{\alpha_X}{\Gamma}} \\
E(\YGamma) \ar[r, swap, ""{name = L}]{}{E(\hGamma)} & E(\XGamma) \ar[from = U, to = L, Rightarrow, shorten <= 4pt, shorten >= 4pt]{}{\quot{\alpha}{\hGamma}}
\end{tikzcd}
\]
is an invertible $2$-cell because $\alpha$ is a pseudonatural transformation. From here it is a routine but tedious calculation to verify that the desired identities hold for any $f \in \Sf(G)_1$ when $\theta_{\Gamma} := \quot{\alpha}{\hGamma}.$
\end{proof}

Before studying Change of Groups functors, we will illustrate the power of Change of Space functors and the descent pullback language to give an important result: The equivariant category $F_G(X)$ is formally na{\"i}evely equivariant in the sense that there are functors $\pi_2^{\ast}:F_G(X) \to F_G(G \times X)$ induced by the projection map $\pi_2:G \times X \to X$, $\alpha_X^{\ast}:F_G(X) \to F_G(G \times X)$ induced by the action map $\alpha_X:G \times X \to X$, and for any object $A \in F_G(X)$ there is an isomorphism
\[
\Theta:\alpha_X^{\ast}A \to \pi_2^{\ast}A
\]
which satisfies the cocycle condition of \cite[Definition 1.6]{GIT} (which we spell out later). Our first step here, however, is to construct the functors $\pi_2^{\ast}$ and $\alpha_X^{\ast}$.
\begin{lemma}
Any pre-equivariant pseudofunctor on $X$ is a simultaneous pre-equivariant pseudofunctor over $X$ and $G \times X$.
\end{lemma}
\begin{proof}
Because each variety $\Gamma \times (G \times X) \in \SfResl_G(G \times X)_0$ also satisfies
\[
\Gamma \times (G \times X) \cong (\Gamma \times G) \times X \in \SfResl_G(X)_0
\]
because $\Gamma \times G \in \Sf(G)_0$, the result follows by virtue of $F$ being defined on $\SfResl_G(X)^{\op}$.
\end{proof}
\begin{lemma}\label{Lemma: Section Change of Space: pi2 pullback}
For any pre-equivariant pseudofunctor $F$ over $X$, there is a functor
\[
\pi_2^{\ast}:F_G(X) \to F_G(G \times X).
\]
\end{lemma}
\begin{proof}
Since $\pi_2:G \times X \to X$ is $G$-equivariant, apply Theorem \ref{Thm: Section 3: Existence of pullback functors} to produce $\pi_w^{\ast}.$
\end{proof}
Constructing the functor $\alpha_X$, however, is more difficult. Because the action morphism $\alpha_X:G \times X \to X$ is rarely $G$-equivariant; the corresponding diagram takes the form (where $\alpha_{G \times X}$ uses the diagonal $G$-action on $G \times X$)
\[
\xymatrix{
G \times G \times X \ar[rr]^-{\alpha_{G \times X}} \ar[d]_{\id_G\times \alpha_X} & & G \times X \ar[d]^{\alpha_X} \\
G \times X \ar[rr]_-{\alpha_X} & & X
}
\]
which commutes if and only if 
\[
\alpha_{X}(\alpha_{G \times X}(h,g,x)) = hghx = hgx = \alpha_X((\id_G \times \alpha_X)(h,g,x))
\]
for every $x \in X$ and $g,h \in G$. Because this implies the action on $X$ is trivial, we cannot assume this. As such, to construct the pullback of the action $\alpha_X^{\ast}:F_G(X) \to F_G(G \times X)$ we need to use the language of descent pullbacks in order to overcome this obstruction.
\begin{lemma}\label{Lemma: Section Change of Space: Action map pullback}
For any pre-equivariant pseudofunctor $F$ on $X$ there is a functor
\[
\alpha_X^{\ast}:F_G(X) \to F_G(G \times X).
\]
\end{lemma}
\begin{proof}
Our strategy is to prove that $F$ admits descent pullbacks along $\alpha_X$. Begin letting $f \in \Sf(G)_1$ be arbitrary and write $\Dom f = \Gamma, \Codom f = \Gamma^{\prime}$. Then because $f$ is $G$-equivariant we get a commuting cube
\[
\begin{tikzcd}
 & \Gamma^{\prime} \times G \times X \ar[rr]{}{\id_{\Gamma^{\prime}} \times \alpha_{X}} \ar[dd, swap, near end]{}{\quo_{\Gamma^{\prime} \times G}} & & \Gamma^{\prime} \times X  \ar[dd]{}{\quo_{\Gamma^{\prime}}} \\
\Gamma \times G \times X \ar[dd, swap]{}{\quo_{\Gamma \times G}} \ar[ur]{}{f \times \id_{G \times X}} & & \Gamma \times X \ar[ur]{}{f \times \id_X} \\
 & \quot{X}{\Gamma^{\prime} \times G} \ar[rr, near start, swap]{}{\overline{\id_{\Gamma^{\prime}} \times \alpha_X}} & & \XGammap \\
\quot{X}{\Gamma \times G} \ar[ur]{}{\quot{\of}{G \times X}} \ar[rr, swap]{}{\overline{\id_{\Gamma} \times \alpha_X}} & & \XGamma \ar[ur, swap]{}{\fX} \ar[near end, from = 2-1, to = 2-3, crossing over]{}{\id_{\Gamma} \times \alpha_X} \ar[from = 2-3, to = 4-3, near start, crossing over]{}{\quo_{\Gamma}}
\end{tikzcd}
\]
which in turn induces the commuting square of quotient varieties
\[
\begin{tikzcd}
\quot{X}{\Gamma \times G} \ar[rr]{}{\overline{\id_{\Gamma} \times \alpha_X}} \ar[d, swap]{}{\quot{\of}{G \times X}} & & \XGamma \ar[d]{}{\fX} \\
\quot{X}{\Gamma^{\prime} \times G} \ar[rr, swap]{}{\overline{\id_{\Gamma^{\prime}} \times \alpha_X}} & & \XGammap
\end{tikzcd}
\]
so applying the functor $F$ to the diagram gives the pasting diagram:
\[
\begin{tikzcd}
F(\XGammap) \ar[dd, swap]{}{F(\fX)} \ar[rrr, ""{name = U}]{}{F(\overline{\id_{\Gamma^{\prime}} \times \alpha_X})}\ar[ddrrr, ""{name = M}]{}{} & & & F(\quot{X}{\Gamma^{\prime} \times G}) \ar[dd]{}{F(\quot{\of}{G \times X})}\ar[r, equals] & F(\quot{(G \times X)}{\Gamma^{\prime}}) \ar[dd]{}{F(\quot{\of}{G \times X})} \\
\\
F(\XGamma) \ar[rrr, swap, ""{name = L}]{}{F(\overline{\id_{\Gamma} \times \alpha_X})} & & & F(\quot{X}{\Gamma \times G}) \ar[r, equals] & F(\quot{(G \times X)}{\Gamma}) \ar[from = U, to = M, Rightarrow, near end, shorten <= 4pt, shorten >= 4pt]{}{\phi_{\quot{f}{G \times X},\id_{\Gamma^{\prime}} \times \alpha_X}} \ar[from = M, to = L, Rightarrow, shorten <= 4pt, shorten >= 4pt, swap]{}{\phi^{-1}_{\id_{\Gamma} \times \alpha_X, \quot{f}{X}}}
\end{tikzcd}
\]
From here we define our descent pullbacks by setting the $\Gamma$-local functors to be
\[
\quot{\alpha_X^{\ast}}{\Gamma} := F(\overline{\id_{\Gamma} \times \alpha_X})
\]
and defining our commutativity natural isomorphisms by
\[
\quot{\alpha_X^{\ast}}{f} := \phi_{\quot{f}{G \times X},\id_{\Gamma^{\prime}} \times \alpha_X}^{-1} \circ \phi_{\id_{\Gamma} \times \alpha_X, \quot{f}{X}}.
\]
From here it follows mutatis mutandis to the proof of Theorem \ref{Thm: Section 3: Existence of pullback functors} that this assignment satisfies what is required of descent pullbacks. Appealing from here to Theorem \ref{Thm: Section 3.2: Existence of descent pullbacks} gives the functor $\alpha_X^{\ast}$.
\end{proof}
\begin{remark}
In what follows, we will write
\[
\quot{\overline{\alpha}_X}{\Gamma} := \overline{\id_{\Gamma} \times \alpha_X}:\quot{X}{\Gamma \times G} \to \XGamma
\]\index[notation]{alphaoverline@$\quot{\overline{\alpha}_X}{\Gamma}$}
in order to reduce the already significant notational bloat.
\end{remark}
\begin{remark}\label{Remark: Section Change of Space: Explict pullback functors}
More explicitly, the functors of Lemmas \ref{Lemma: Section Change of Space: pi2 pullback} and \ref{Lemma: Section Change of Space: Action map pullback} are given as follows on objects. Let $(A,T_A) \in F_G(X)_0$. Then by Theorem \ref{Thm: Section 3: Existence of pullback functors}
\[
\pi_2^{\ast}A = \left\lbrace F(\quot{\opi_2}{\Gamma})\AGamma \; | \; \Gamma \in \Sf(G)_0 \right\rbrace,
\]
where $\quot{\pi_2}{\Gamma}$ is the projection map $\Gamma \times G \times X \to \Gamma \times X$, and the transition isomorphisms are given as
\[
T_{\pi_2^{\ast}A} = \left\lbrace F(\quot{\opi_2}{\Gamma})\tau_f^A \circ \left(\phi_{\quot{\pi_2}{\Gamma}, \quot{f}{X}}^{\AGammap}\right)^{-1} \circ \phi_{\quot{f}{G \times X}, \quot{\pi_2}{\Gamma^{\prime}}}^{\AGammap}  \; : \; f \in \Sf(G)_1 \right\rbrace.
\]
Similarly, by construction we also have that
\[
\alpha_X^{\ast}A = \left\lbrace F(\quot{\overline{\alpha}_X}{\Gamma})\AGamma \; | \; \Gamma \in \Sf(G)_0\right\rbrace
\]
while the transition isomorphisms are
\[
T_{\alpha_X^{\ast}A} = \left\lbrace F(\quot{\overline{\alpha}_X^{\ast}}{\Gamma})\tau_f^A \circ \left(\phi_{\quot{\alpha_X}{\Gamma}, \quot{f}{X}}^{\AGammap}\right)^{-1}  \circ \phi_{\quot{f}{G \times X}, \quot{\alpha_X}{\Gamma^{\prime}}}^{\AGammap} \; : \; f \in \Sf(G)_1 \right\rbrace.
\]
\end{remark}

In order to produce the isomorphism $\Theta$, we follow a process discovered by Clifton Cunningham during the BIRS focused research group ``the Voganish Project.'' The main idea here is to find, for each variety $\Gamma \in \Sf(G)_0$, an $\Sf(G)$-variety $\Gamma_c$ which induces an isomorphism of quotient varieties
\[
\quot{X}{\Gamma_c} \cong \quot{(G \times X)}{\Gamma}
\]
that is natural over $\Sf(G)$ with respect to the actions of $\quot{\overline{\alpha}_X}{\Gamma}$ and $\quot{\opi_2}{\Gamma}$ so that we can use the transition isomorphisms of an object of $F_G(X)$ to build an isomorphism with the objects $\pi_2^{\ast}A$ and $\alpha_X^{\ast}A$. We present this as the lemma below, which is the main ingredient in cooking up our desired isomorphism $\Theta$.
\begin{lemma}\label{Lemma: Clifton's Naivification Variety}\index[notation]{Gammac@$\Gamma_c$}
For all $\Gamma \in \Sf(G)_0$ there exists a variety $\Gamma_c \in \Sf(G)_0$ such that there is an isomorphism of quotient varieties
\[
\quot{(G \times X)}{\Gamma} \xrightarrow{\mu_{\Gamma}} \quot{X}{\Gamma_c}.
\]
Moreover, if $f \in \Sf(G)_1$ with $f \in \Sf(G)(\Gamma, \Gamma^{\prime})$ there is an $\Sf(G)$-morphism 
\[
f_c:\Gamma_c \to \Gamma^{\prime}_c
\] 
for which the diagram
\[
\xymatrix{
\quot{(G \times X)}{\Gamma} \ar[rr]^-{\mu_{\Gamma}} \ar[d]_{\quot{\of}{G \times X}} & & \quot{X}{\Gamma_c} \ar[d]^{\overline{f}_c} \\
\quot{(G \times X)}{\Gamma^{\prime}} \ar[rr]_-{\mu_{\Gamma^{\prime}}} & & \quot{X}{\Gamma^{\prime}_c}
}
\]
commutes. There are also $\Sf(G)$-morphisms $a_{\Gamma}, p_{\Gamma}:\Gamma_c \to \Gamma$ for which the diagrams
\[
\xymatrix{
\quot{(G \times X)}{\Gamma} \ar[r]^-{\mu_{\Gamma}} \ar[dr]_{\quot{\overline{\alpha}_X}{\Gamma}} & \quot{X}{\Gamma_c} \ar[d]^{\overline{a}_{\Gamma}} & \quot{(G \times X)}{\Gamma} \ar[r]^-{\mu_{\Gamma}} \ar[dr]_{\quot{\opi_2}{\Gamma}} & \quot{X}{\Gamma_c} \ar[d]^{\overline{p}_{\Gamma}}\\
 & \XGamma & & \XGamma
}
\]
commute. Finally, for any $f \in \Sf(G)_1$ the diagrams
\[
\xymatrix{
\Gamma_c \ar[d]_{f_c} \ar[r]^-{a_{\Gamma}} & \Gamma \ar[d]^{f} & \Gamma_c \ar[r]^-{p_{\Gamma}} \ar[d]_{f_c} & \Gamma \ar[d]^{f} \\
\Gamma_c^{\prime} \ar[r]_-{a_{\Gamma^{\prime}}} & \Gamma^{\prime} & \Gamma_c^{\prime} \ar[r]_-{p_{\Gamma^{\prime}}} & \Gamma^{\prime}
}
\]
commute as well.
\end{lemma}
\begin{proof}
We give this proof by doing a sketch using the stalks of the quotient varieties for readability. Let $\Gamma \in \Sf(G)_0$. Then define $\Gamma_c$ by setting
\[
\Gamma_c := \Gamma \times G
\]
with action map
\[
\alpha_{\Gamma_c} = \alpha_{\Gamma} \times \Ad
\]
where $\Ad$ is the adjoint action of $G$ on itself. More explicitly, $\Ad:G \times G \to G$ is the morphism given on sets by $\Ad(h,g) = hgh^{-1}$ and in categorical terms by
\[
\xymatrix{
G \times G \ar[d]_{\Delta \times \id_G} \ar[rrr]^-{\Ad} & & & G \\
G \times G \times G \ar[rr]_-{\id_G \times \inv \times \id_G} & & G \times G \times G \ar[r]_-{\id_G \times s} & G \times G \times G \ar[u]_{\mu \circ (\id_G \times \mu)}
}
\]
where $\inv:G \to G$ is the inversion map and $s$ is the switching isomorphism of the product. Now because $\Sf(G)$ admits products and both $\Gamma$ and $G$ are $\Sf(G)$-varieties, $\Gamma_c$ is a pure dimensional smooth variety. That it is a principal $G$-variety follows from the fact that $\Gamma \in \Sf(G)_0$ and Proposition \ref{Prop: Smooth quotient from principal G var}.

We define the morphism of quotient varieties $\mu_{\Gamma}:\quot{(G \times X)}{\Gamma} \to \quot{X}{\Gamma_c}$ by (in set-theoretic terms)
\[
[\gamma, (g,x)]_G \mapsto [(\gamma, g),x]_G;
\]
that this is an isomorphism is trivial to verify. Define $f_c:\Gamma_c \to \Gamma^{\prime}_c$ by $f_c := f \times \id_G$. Then we calculate that
\begin{align*}
(\of_c \circ \mu_{\Gamma})[\gamma, (g,x)]_G &= \of_c[(\gamma, g),x]_G = [(f(\gamma),g),x]_G = \mu_{\Gamma^{\prime}}[f(\gamma),(g,x)]_G \\
&= (\mu_{\Gamma^{\prime}} \circ \quot{\of}{G \times X})[\gamma, (g,x)]_G,
\end{align*}
which shows that the first diagram commutes.

We define the morphisms $a_{\Gamma}, p_{\Gamma}:\Gamma_c \to \Gamma$.  The map $a_{\Gamma}$ is defined by the diagram
\[
\xymatrix{
\Gamma_c \ar[d]_{a_{\Gamma}} \ar@{=}[r] & \Gamma \times G \ar[r]_-{\id_{\Gamma} \times \inv} & \Gamma \times G \ar[d]_{s}\\
\Gamma & & G \times \Gamma \ar[ll]^-{\alpha_{\Gamma}}
}
\]
(which says in set-theoretic terms that $a_{\Gamma}(\gamma,g) = g^{-1}\gamma$) while we define $p_{\Gamma}$ $=$ $\pi_2^{\Gamma G}:\Gamma \times G \to G$ to be the second projection. Because $\pi_2^{\Gamma G}$ is trivially verified to be an $\Sf(G)$-morphism, we only need to do this for $a_{\Gamma}$. For $a_{\Gamma}$ we note that it is smooth and of constant fibre dimension, so we only need to show that it is $G$-equivariant. We calculate that on one hand
\[
(a_{\Gamma} \circ \alpha_{\Gamma_c})(h,\gamma,g) = a_{\Gamma}(h\gamma, hgh^{-1}) = (hgh^{-1})^{-1}h\gamma = hg^{-1}hh^{-1}\gamma = hg^{-1}\gamma
\]
while on the other hand
\[
\big(\alpha_{\Gamma} \circ (\id_{G} \times a_{\Gamma})\big)(h,\gamma,g) = (h,g^{-1}\gamma) = hg^{-1}\gamma.
\]
Thus $a_{\Gamma}$ is $G$-equivariant and hence an $\Sf(G)$-morphism.

We now verify the commutativity of the triangles. For this note that the triangle involving $p_{\Gamma}$ commutes trivially, so we only verify the triangle for $a_{\Gamma}$. In this case we find that
\begin{align*}
(\overline{a}_{\Gamma} \circ \mu_{\Gamma})[\gamma, (g,x)]_G &= \overline{a}_{\Gamma}[(\gamma,g),x]_G = [g^{-1}\gamma,x]_G = [\gamma,gx]_G \\
&= {\overline{\alpha}_X}{\Gamma}[\gamma, (g, x)]_G,
\end{align*}
as was desired. Finally if $f \in \Sf(G)_1$, the last two squares commute trivially in the case of $p_{\Gamma}, p_{\Gamma^{\prime}}$ and exactly because $f$ is $G$-equivariant in the case of $a_{\Gamma}, a_{\Gamma^{\prime}}$.
\end{proof}

We now finally have the tools to provide our claimed isomorphism $\Theta:\alpha_X^{\ast}A \to \pi_2^{\ast}A$ for any object of $F_G(X)$. The construction of this isomorphism, in view of Lemma \ref{Lemma: Clifton's Naivification Variety} is straightforward, but the proof that it is an $F_G(X)$-morphism is quite difficult and involved. The techniques we use are essentially a generalization of an unpublished proof of C.\@ Cunningham and myself, as the proof we gave establishes (after generalizing) the case when $F$ is a strict pseudofunctor.

\begin{Theorem}\label{Theorem: Formal naive equivariance}
For any pre-equivariant pseudofunctor $F$ over $X$ and any $A \in F_G(X)_0$, there is an isomorphism
\[
\Theta_A:\alpha_X^{\ast}A \to \pi_2^{\ast}A
\]
such that for any $P \in F_G(X)_1$ the diagram
\[
\xymatrix{
\alpha_X^{\ast}A \ar[d]_{\alpha_X^{\ast}P} \ar[r]^-{\Theta_A} & \pi_2^{\ast}A \ar[d]^{\pi_2^{\ast}P} \\
\alpha_X^{\ast}B \ar[r]_-{\Theta_B} & \pi_2^{\ast}B
}
\]
commutes in $F_G(G \times X)$.
\end{Theorem}
\begin{proof}
We begin by fixing $A \in F_G(X)_0$ and defining the isomorphism $\quot{\theta}{\Gamma}$ as in the diagram below:
\[
\begin{tikzcd}
F(\quot{\overline{\alpha}_X}{\Gamma})(\AGamma) \ar[rrr]{}{\quot{\theta}{\Gamma}} \ar[d, equals]  && & F(\quot{\opi_2}{\Gamma}) \\
F(\aGamma \circ \muGamma)(\AGamma) \ar[d, swap]{}{\left(\phi_{\muGamma,a_{\Gamma}}^{\AGamma}\right)^{-1}} & & & F(\pGamma \circ \muGamma)(\AGamma) \ar[u, equals] \\
F(\muGamma)\left(F(\aGamma)(\AGamma)\right) \ar[r, swap]{}{F(\muGamma)\tau_{a_{\Gamma}}^{A}} & F(\muGamma)\quot{A}{\Gamma_c} \ar[rr, swap]{}{F(\muGamma)(\tau_{p_{\Gamma}}^{A})^{-1}} & & F(\muGamma)\left(F(\pGamma)(\AGamma)\right) \ar[u, swap]{}{\phi_{\muGamma, p_{\Gamma}}}
\end{tikzcd}
\]
We now begin our arduous task of showing that the equation 
\[
\tau_f^{\pi_2^{\ast}A} \circ F(\fGX)\quot{\theta}{\Gamma^{\prime}} = \quot{\theta}{\Gamma} \circ \tau_f^{\alpha_X^{\ast}A}
\] 
holds. For this we first find
\begin{align*}
&\tau_f^{\pi_2^{\ast}A} \circ F(\fGX)\quot{\theta}{\Gamma^{\prime}} \\
&= \tau_f^{\pi_2^{\ast}A} \circ F(\fGX)\left(\phi_{\muGammap,p_{\Gamma^{\prime}}}^{\AGammap} \circ F(\muGammap) \left(\left(\tau_{p_{\Gamma^{\prime}}}^{A}\right)^{-1} \circ \tau_{a_{\Gamma^{\prime}}}^{A}\right) \circ \left(\phi_{\muGammap, a_{\Gamma^{\prime}}}^{\AGammap}\right)^{-1}\right) \\
&=\tau_f^{\pi_2^{\ast}A} \circ F(\fGX)\left(\phi_{\muGammap,p_{\Gamma^{\prime}}}^{\AGammap}\right) \circ  \big(F(\fGX) \circ F(\muGammap)\big)\left(\left(\tau_{p_{\Gamma^{\prime}}}^{A}\right)^{-1} \circ \tau_{a_{\Gamma^{\prime}}}^{A}\right) \\
&\circ  F(\fGX)\left(\phi_{\muGammap, a_{\Gamma^{\prime}}}^{\AGammap}\right)^{-1}.
\end{align*} 
Now using the naturality of the compositor isomorphisms we calculate that 
\begin{align*}
& \big(F(\fGX) \circ F(\muGammap)\big)\left(\left(\tau_{p_{\Gamma^{\prime}}}^{A}\right)^{-1} \circ \tau_{a_{\Gamma^{\prime}}}^{A}\right) \\
&= \left(\phi_{\fGX, \muGammap}^{F(\pGammap)\AGammap}\right)^{-1} \circ  F(\muGammap \circ \fGX)\left(\left(\tau_{p_{\Gamma^{\prime}}}^{A}\right)^{-1} \circ\tau_{a_{\Gamma^{\prime}}}^{A}\right) \circ \phi_{\fGX, \muGammap}^{F(\aGammap)\AGammap}.
\end{align*}
Using that there is an equality of functors
\[
F(\muGammap \circ \fGX) = F(\of_c \circ \muGamma)
\]
we have
\begin{align*}
&F(\muGammap \circ \fGX)\left(\left(\tau_{p_{\Gamma^{\prime}}}^{A}\right)^{-1} \circ\tau_{a_{\Gamma^{\prime}}}^{A}\right) =F(\of_c \circ \muGamma)\left(\left(\tau_{p_{\Gamma^{\prime}}}^{A}\right)^{-1} \circ\tau_{a_{\Gamma^{\prime}}}^{A}\right)
\end{align*}
so using the same trick we write
\begin{align*}
&F(\of_c \circ \muGamma)\left(\left(\tau_{p_{\Gamma^{\prime}}}^{A}\right)^{-1} \circ\tau_{a_{\Gamma^{\prime}}}^{A}\right) \\
&= \phi_{\muGamma, \of_c}^{F(\aGammap)\AGammap} \circ \big(F(\muGamma) \circ F(\of_c)\big)\left(\left(\tau_{p_{\Gamma^{\prime}}}^{A}\right)^{-1} \circ\tau_{a_{\Gamma^{\prime}}}^{A}\right) \circ \left(\phi_{\muGamma, \of_c}^{F(\aGammap)\AGammap}\right)^{-1}
\end{align*}
and get that
\begin{align*}
&\tau_f^{\pi_2^{\ast}A} \circ F(\fGX)\quot{\theta}{\Gamma^{\prime}} \\
&= \tau_f^{\pi_2^{\ast}A} \circ F(\fGX)\left(\phi_{\muGammap,p_{\Gamma^{\prime}}}^{\AGammap} \circ F(\muGammap) \left(\left(\tau_{p_{\Gamma^{\prime}}}^{A}\right)^{-1} \circ \tau_{a_{\Gamma^{\prime}}}^{A}\right) \circ \left(\phi_{\muGammap, a_{\Gamma^{\prime}}}^{\AGammap}\right)^{-1}\right) \\
&=\tau_f^{\pi_2^{\ast}A} \circ F(\fGX)\left(\phi_{\muGammap,p_{\Gamma^{\prime}}}^{\AGammap}\right) \circ  \big(F(\fGX) \circ F(\muGammap)\big)\left(\left(\tau_{p_{\Gamma^{\prime}}}^{A}\right)^{-1} \circ \tau_{a_{\Gamma^{\prime}}}^{A}\right) \\
&\circ  F(\fGX)\left(\phi_{\muGammap, a_{\Gamma^{\prime}}}^{\AGammap}\right)^{-1} \\
&=\tau_f^{\pi_2^{\ast}A} \circ F(\fGX)\left(\phi_{\muGammap,p_{\Gamma^{\prime}}}^{\AGammap}\right) \circ  \left(\phi_{\fGX, \muGammap}^{F(\pGammap)\AGammap}\right)^{-1} \circ  \phi_{\muGamma, \of_c}^{F(\pGammap)\AGammap} \\
&\circ \big(F(\muGamma) \circ F(\of_c)\big)\left(\left(\tau_{p_{\Gamma^{\prime}}}^{A}\right)^{-1} \circ\tau_{a_{\Gamma^{\prime}}}^{A}\right) \circ \left(\phi_{\muGamma, \of_c}^{F(\aGammap)\AGammap}\right)^{-1} \circ \phi_{\fGX, \muGammap}^{F(\aGammap)\AGammap} \\
&\circ  F(\fGX)\left(\phi_{\muGammap, a_{\Gamma^{\prime}}}^{\AGammap}\right)^{-1}
\end{align*}
while
\begin{align*}
\tau_f^{\pi_2^{\ast}A} &= F(\opiGamma)\tau_f^A \circ \left(\phi_{\quot{\pi_2}{\Gamma}, \fX}^{\AGammap}\right)^{-1} \circ \phi_{\fGX, \quot{\pi_2}{\Gamma}}^{\AGammap} \\
&= F(\pGamma \circ \muGamma)\tau_f^A \circ \left(\phi_{\pGamma \circ \muGamma, \fX}^{\AGammap}\right)^{-1} \circ \phi_{\fGX, \pGammap \circ \muGammap}^{\AGammap} \\
&= \phi_{\muGamma, \pGamma}^{\AGamma} \circ \big(F(\muGamma) \circ F(\pGamma)\big)\tau_f^A \circ \left(\phi_{\muGamma, \pGamma}^{F(\fX)\AGammap}\right)^{-1} \circ \left(\phi_{\pGamma \circ \muGamma, \fX}^{\AGammap}\right)^{-1} \\
&\circ \phi_{\fGX, \pGammap \circ \muGammap}^{\AGammap}.
\end{align*}
In order to manipulate the (implicit) monster composition, we need to derive some identities with the compositor natural isomorphisms to simplify the expression. In particular, we will show that we can simplify the morphism
\begin{align*}
&\left(\phi_{\muGamma, \pGamma}^{F(\fX)\AGammap}\right)^{-1} \circ \left(\phi_{\pGamma \circ \muGamma, \fX}^{\AGammap}\right)^{-1} \circ \phi_{\fGX, \pGammap \circ \muGammap}^{\AGammap} \circ F(\fGX)\left(\phi_{\muGammap,p_{\Gamma^{\prime}}}^{\AGammap}\right) \\
&\circ  \left(\phi_{\fGX, \muGammap}^{F(\pGammap)\AGammap}\right)^{-1}
\end{align*}
by manipulating the coherences of the compositors.

For this we first calculate, using the pasting diagrams of pseudofunctors, that the diagram
\[
\begin{tikzcd}
F(\XGammap) \ar[r]{}{F(\fX)} \ar[rr, bend right = 60, swap, ""{name = L}]{}{F(\fX \circ \pGamma \circ \muGamma)} & F(\XGamma) \ar[r, bend left = 30, ""{name = U}]{}{F(\muGamma) \circ F(\pGamma)} \ar[r, bend right = 30, swap, ""{name = M}]{}{F(\pGamma \circ \muGamma)} & F(\quot{(G \times X)}{\Gamma}) \ar[from = U, to = M, Rightarrow, shorten <= 4pt, shorten >= 4pt] \ar[from = 1-2, to = L, Rightarrow, shorten <= 4pt, shorten >= 4pt]{}{}
\end{tikzcd}
\]
is equivalent to the pasting diagram:
\[
\begin{tikzcd}
F(\XGammap) \ar[r, bend left = 30, ""{name = U}]{}{F(\pGamma) \circ F(\fX)} \ar[r, bend right = 30, swap, ""{name = M}]{}{F(\fX \circ \pGamma)} \ar[rr, bend right = 70, swap, ""{name = L}]{}{F(\fX \circ \pGamma \circ \muGamma)} & F(\quot{X}{\Gamma_c}) \ar[r]{}{F(\muGamma)} & F(\quot{(G \times X)}{\Gamma}) \ar[from = U, to = M, Rightarrow, shorten <= 4pt, shorten >= 4pt] \ar[from = 1-2, to = L, Rightarrow, shorten <= 4pt, shorten >= 4pt]{}{}
\end{tikzcd}
\]
Note that in both cases, the $2$-cells are given by the corresponding compositors. This gives us that
\[
\left(\phi_{\mu_{\Gamma},\pGamma}^{F(\fX)\AGammap}\right)^{-1} \circ \left(\phi_{\pGamma \circ \muGamma, \fX}^{\AGammap}\right)^{-1} = F(\muGamma)\left(\phi_{\pGamma, \fX}^{\AGammap}\right)^{-1} \circ \left(\phi_{\muGamma, \fX \circ \pGamma}^{\AGammap}\right)^{-1}.
\]
Similarly, the pasting diagram
\[
\begin{tikzcd}
F(\XGammap) \ar[r, bend left = 30, ""{name = U}]{}{F(\muGammap) \circ F(\pGammap)} \ar[r, bend right = 30, swap, ""{name = M}]{}{F(\opiGammap)} \ar[rr, bend right = 60, swap, ""{name = L}]{}{F(\pGammap \circ \muGammap \circ \fGX)} & F(\quot{(G \times X)}{\Gamma^{\prime}}) \ar[r]{}{F(\fGX)} & F(\quot{(G \times X)}{\Gamma}) \ar[from = U, to = M, Rightarrow, shorten <= 4pt, shorten >= 4pt]{}{} \ar[from = 1-2, to = L, Rightarrow, shorten <= 4pt, shorten >= 4pt]{}{}
\end{tikzcd}
\]
is equivalent to the pasting diagram:
\[
\begin{tikzcd}
F(\XGammap) \ar[r]{}{F(\pGammap)} \ar[rr, bend right = 70, swap, ""{name = L}]{}{F(\pGammap \circ \muGammap \circ \fGX)} & F(\quot{X}{\Gamma_c^{\prime}}) \ar[r, bend left = 30, ""{name = U}]{}{F(\fGX) \circ F(\muGammap)} \ar[r, swap, bend right = 30, ""{name = M}]{}{F(\muGammap \circ \fGX)} &  F(\quot{(G \times X)}{\Gamma}) \ar[from = U, to = M, Rightarrow, shorten <= 4pt, shorten >= 4pt]{}{} \ar[from = 1-2, to = L, Rightarrow, shorten <= 4pt, shorten >= 4pt]{}{}
\end{tikzcd}
\]
As in the prior case that in both cases the $2$-cells are given by the corresponding compositors. This allows us to deduce that
\[
\phi_{\fGX, \pGammap \circ \muGammap}^{\AGammap} \circ F(\fGX)\phi_{\muGammap, \pGammap}^{\AGammap} = \phi_{\muGammap \circ \fGX, \pGammap}^{\AGammap} \circ \phi_{\fGX, \muGammap}^{F(\pGammap)\AGammap}.
\]
This allows us to compute that
\begin{align*}
&\left(\phi_{\muGamma, \pGamma}^{F(\fX)\AGammap}\right)^{-1} \circ \left(\phi_{\pGamma \circ \muGamma, \fX}^{\AGammap}\right)^{-1} \circ \phi_{\fGX, \pGammap \circ \muGammap}^{\AGammap} \circ F(\fGX)\left(\phi_{\muGammap,p_{\Gamma^{\prime}}}^{\AGammap}\right) \\
&\circ  \left(\phi_{\fGX, \muGammap}^{F(\pGammap)\AGammap}\right)^{-1} \\
&= F(\muGamma)\left(\phi_{\pGamma, \fX}^{\AGammap}\right)^{-1} \circ \left(\phi_{\muGamma, \fX \circ \pGamma}^{\AGammap}\right)^{-1} \circ \phi_{\muGammap \circ \fGX, \pGammap}^{\AGammap} \circ \phi_{\fGX, \muGammap}^{F(\pGammap)\AGammap} \\
&\circ \left(\phi_{\fGX, \muGammap}^{F(\pGammap)\AGammap}\right)^{-1} \\
&=F(\muGamma)\left(\phi_{\pGamma, \fX}^{\AGammap}\right)^{-1} \circ \left(\phi_{\muGamma, \fX \circ \pGamma}^{\AGammap}\right)^{-1} \circ \phi_{\muGammap \circ \fGX, \pGammap}^{\AGammap}.
\end{align*}
We will now use the same tricks to simplify the composition
\[
F(\muGamma)\left(\phi_{\pGamma, \fX}^{\AGammap}\right)^{-1} \circ \left(\phi_{\muGamma, \fX \circ \pGamma}^{\AGammap}\right)^{-1} \circ \phi_{\muGammap \circ \fGX, \pGammap}^{\AGammap} \circ  \phi_{\muGamma, \of_c}^{F(\pGammap)\AGammap}.
\]
First note that $\muGammap \circ \fGX = \of_c \circ \muGamma$ so
\[
\phi_{\muGammap \circ \fGX, \pGammap}^{\AGammap} \circ  \phi_{\muGamma, \of_c}^{F(\pGammap)\AGammap} = \phi_{\of_c \circ \muGamma, \pGammap}^{\AGammap} \circ  \phi_{\muGamma, \of_c}^{F(\pGammap)\AGammap}.
\]
Now because the pasting diagram
\[
\begin{tikzcd}
F(\XGammap) \ar[r]{}{F(\pGammap)} \ar[rr, bend right = 70, swap, ""{name = L}]{}{F(\pGammap \circ \of_c \circ \muGamma)} & F(\quot{X}{\Gamma_c^{\prime}}) \ar[r, bend left = 30, ""{name = U}]{}{F(\muGamma) \circ F(\of_c)} \ar[r, swap, bend right = 30, ""{name = M}]{}{F(\of_c \circ \muGamma)} &  F(\quot{(G \times X)}{\Gamma}) \ar[from = U, to = M, Rightarrow, shorten <= 4pt, shorten >= 4pt]{}{} \ar[from = 1-2, to = L, Rightarrow, shorten <= 4pt, shorten >= 4pt]{}{}
\end{tikzcd}
\]
is equivalent to
\[
\begin{tikzcd}
F(\XGammap) \ar[r, bend left = 30, ""{name = U}]{}{F(\of_c) \circ F(\pGammap)} \ar[r, bend right = 30, swap, ""{name = M}]{}{F(\pGammap \circ \of_c)} \ar[rr, bend right = 60, swap, ""{name = L}]{}{F(\pGammap \circ \of_c \circ \muGamma)} & F(\quot{(G \times X)}{\Gamma^{\prime}}) \ar[r]{}{F(\muGamma)} & F(\quot{(G \times X)}{\Gamma}) \ar[from = U, to = M, Rightarrow, shorten <= 4pt, shorten >= 4pt]{}{} \ar[from = 1-2, to = L, Rightarrow, shorten <= 4pt, shorten >= 4pt]{}{}
\end{tikzcd}
\]
we get that
\begin{align*}
 \phi_{\muGammap \circ \fGX, \pGammap}^{\AGammap} \circ  \phi_{\muGamma, \of_c}^{F(\pGammap)\AGammap} &= \phi_{\of_c \circ \muGamma, \pGammap}^{\AGammap} \circ  \phi_{\muGamma, \of_c}^{F(\pGammap)\AGammap} \\
 &=\phi_{\muGamma, \pGammap \circ \of_c}^{\AGammap} \circ F(\muGamma)\phi_{\of_c, \pGammap} \\
 &=\phi_{\muGamma, \fX \circ \pGamma}^{\AGammap} \circ F(\muGamma)\phi_{\of_c, \pGammap}^{\AGammap}.
\end{align*}
Using this we derive that
\begin{align*}
&F(\muGamma)\left(\phi_{\pGamma, \fX}^{\AGammap}\right)^{-1} \circ \left(\phi_{\muGamma, \fX \circ \pGamma}^{\AGammap}\right)^{-1} \circ \phi_{\muGammap \circ \fGX, \pGammap}^{\AGammap} \circ  \phi_{\muGamma, \of_c}^{F(\pGammap)\AGammap} \\
&= F(\muGamma)\left(\phi_{\pGamma, \fX}^{\AGammap}\right)^{-1} \circ \left(\phi_{\muGamma, \fX \circ \pGamma}^{\AGammap}\right)^{-1} \circ \phi_{\muGamma, \fX \circ \pGamma}^{\AGammap} \circ F(\muGamma)\phi_{\of_c, \pGammap}^{\AGammap} \\
&= F(\muGamma)\left(\phi_{\pGamma, \fX}^{\AGammap}\right)^{-1} \circ F(\muGamma)\phi_{\of_c, \pGammap}^{\AGammap}.
\end{align*}
In summary, we have shown so far that
\begin{align*}
&\left(\phi_{\muGamma,\pGamma}^{F(\fX)\AGammap}\right)^{-1} \circ \left(\phi_{\quot{\pi_2}{\Gamma},\fX}^{\AGammap}\right)^{-1} \circ \phi_{\fGX, \quot{\pi_2}{\Gamma^{\prime}}}^{\AGammap} 
\circ F(\fGX)\phi_{\muGammap, \pGammap}^{\AGammap} \\
&\circ \left(\phi_{\fGX, \muGammap}^{F(\pGammap)\AGammap}\right)^{-1} \circ  \phi_{\muGamma, \of_c}^{F(\pGammap)\AGammap} \\
&= F(\muGamma)\left(\phi_{\pGamma, \fX}^{\AGammap}\right)^{-1} \circ F(\muGamma)\phi_{\of_c, \pGammap}^{\AGammap}.
\end{align*}
Altogether this allows us to deduce that
\begin{align*}
&\tau^{\pi_2^{\ast}A} \circ F(\fGX)\quot{\theta}{\Gamma^{\prime}} \\
&= \phi_{\muGamma, \pGamma}^{\AGamma} \circ \big(F(\muGamma) \circ F(\pGamma)\big)\tau_f^A \circ \left(\phi_{\muGamma,\pGamma}^{F(\fX)\AGammap}\right)^{-1} \circ \left(\phi_{\quot{\pi_2}{\Gamma},\fX}^{\AGammap}\right)^{-1} \circ \phi_{\fGX, \quot{\pi_2}{\Gamma^{\prime}}}^{\AGammap} \\
&\circ F(\fGX)\phi_{\muGammap, \pGammap}^{\AGammap} \circ \left(\phi_{\fGX, \muGammap}^{F(\pGammap)\AGammap}\right)^{-1} \circ  \phi_{\muGamma, \of_c}^{F(\pGammap)\AGammap} \\
&\circ \big(F(\muGamma) \circ F(\of_c)\big)\left(\left(\tau_{p_{\Gamma^{\prime}}}^{A}\right)^{-1} \circ\tau_{a_{\Gamma^{\prime}}}^{A}\right) \circ \left(\phi_{\muGamma, \of_c}^{F(\aGammap)\AGammap}\right)^{-1} \circ \phi_{\fGX, \muGammap}^{F(\aGammap)\AGammap} \\
&\circ  F(\fGX)\left(\phi_{\muGammap, a_{\Gamma^{\prime}}}^{\AGammap}\right)^{-1} \\
&= \phi_{\muGamma, \pGamma}^{\AGamma} \circ \big(F(\muGamma) \circ F(\pGamma)\big)\tau_f^A \circ F(\muGamma)\left(\phi_{\pGamma, \fX}^{\AGammap}\right)^{-1} \circ F(\muGamma)\phi_{\of_c, \pGammap}^{\AGammap} \\
&\circ \big(F(\muGamma) \circ F(\of_c)\big)\left(\left(\tau_{p_{\Gamma^{\prime}}}^{A}\right)^{-1} \circ\tau_{a_{\Gamma^{\prime}}}^{A}\right) \circ \left(\phi_{\muGamma, \of_c}^{F(\aGammap)\AGammap}\right)^{-1} \circ \phi_{\fGX, \muGammap}^{F(\aGammap)\AGammap} \\
&\circ  F(\fGX)\left(\phi_{\muGammap, a_{\Gamma^{\prime}}}^{\AGammap}\right)^{-1}.
\end{align*}

To proceed, we first observe that
\begin{align*}
&\left(\phi_{\muGamma, \pGamma}^{\AGamma}\right)^{-1} \circ \tau_{f}^{\pi_2^{\ast}A} \circ F(\fGX)\quot{\theta}{\Gamma^{\prime}} \circ  F(\fGX)\left(\phi_{\muGammap, a_{\Gamma^{\prime}}}^{\AGammap}\right) \circ \left(\phi_{\fGX, \muGammap}^{F(\aGammap)\AGammap}\right)^{-1} \\
&\circ \left(\phi_{\muGamma, \of_c}^{F(\aGammap)\AGammap}\right) \\
&=\big(F(\muGamma) \circ F(\pGamma)\big)\tau_f^A \circ F(\muGamma)\left(\phi_{\pGamma, \fX}^{\AGammap}\right)^{-1} \circ F(\muGamma)\phi_{\of_c, \pGammap}^{\AGammap} \\
&\circ \big(F(\muGamma) \circ F(\of_c)\big)\left(\left(\tau_{p_{\Gamma^{\prime}}}^{A}\right)^{-1} \circ\tau_{a_{\Gamma^{\prime}}}^{A}\right)
\end{align*} 
We will now study what happens when we post-compose this expression by $F(\muGamma)\tau_{\pGamma}^A$. Explicitly,
\begin{align*}
&F(\muGamma)\tau_{\pGamma}^{A} \circ \big(F(\muGamma) \circ F(\pGamma)\big)\tau_f^A \circ F(\muGamma)\left(\phi_{\pGamma, \fX}^{\AGammap}\right)^{-1}\circ F(\muGamma)\phi_{\of_c, \pGammap}^{\AGammap} \\
& \circ \big(F(\muGamma) \circ F(\of_c)\big)\left(\left(\tau_{p_{\Gamma^{\prime}}}^{A}\right)^{-1} \circ\tau_{a_{\Gamma^{\prime}}}^{A}\right) \\
&= F(\muGamma)\left(\tau_{\pGamma}^{A} \circ F(\pGamma)\tau_f^A \circ \left(\phi_{\pGamma, \fX}^{\AGammap}\right)^{-1} \cdots \circ F(\of_c)\left(\left(\tau_{p_{\Gamma^{\prime}}}^{A}\right)^{-1} \circ\tau_{a_{\Gamma^{\prime}}}^{A}\right)\right) \\
&= F(\muGamma)\left(\tau_{\fX \circ \pGamma}^{A} \circ \phi_{\pGamma, \fX}^{\AGammap} \circ \left(\phi_{\pGamma, \fX}^{\AGammap}\right)^{-1} \circ \cdots \circ F(\of_c)\left(\left(\tau_{p_{\Gamma^{\prime}}}^{A}\right)^{-1} \circ\tau_{a_{\Gamma^{\prime}}}^{A}\right)\right) \\
&=F(\muGamma)\left(\tau_{\pGammap \circ \of_c}^{A} \circ \phi_{\pGamma, \of_c}^{\AGammap} F(\of_c)\left(\left(\tau_{p_{\Gamma^{\prime}}}^{A}\right)^{-1} \circ\tau_{a_{\Gamma^{\prime}}}^{A}\right)\right) \\
&=F(\muGamma)\left(\tau_{\of_c}^{A} \circ F(\of_c)\tau_{\pGammap}^A \circ \left(\phi_{\pGammap, \of_c}^{\AGammap}\right)^{-1} \phi_{\pGamma, \of_c}^{\AGammap} F(\of_c)\left(\left(\tau_{p_{\Gamma^{\prime}}}^{A}\right)^{-1} \circ\tau_{a_{\Gamma^{\prime}}}^{A}\right)\right) \\
&=F(\muGamma)\left(\tau_{\of_c}^{A} \circ F(\of_c)\tau_{\pGammap}^A \circ  F(\of_c)\left(\tau_{p_{\Gamma^{\prime}}}^{A}\right)^{-1} \circ F(\of_c)\tau_{a_{\Gamma^{\prime}}}^{A}\right) \\
&= F(\muGamma)\left(F(\tau_{f_c}^{A} \circ F(\of_c)\tau_{a_{\Gamma^{\prime}}}^{A}\right) = F(\muGamma)\left(\tau_{\aGammap \circ \of_c}^{A} \circ \phi_{f_c, \aGammap}^{\AGammap}\right) \\
&= F(\muGamma)\left(\tau_{\fX \circ \aGamma}^{A} \circ \phi_{f_c, \aGammap}^{\AGammap}\right) \\
&=F(\muGamma)\left(\tau_{\aGamma}^{A} \circ F(\aGamma)\tau_{f}^{A} \circ \left(\phi_{\aGamma, \fX}^{\AGammap}\right)^{-1}\circ\phi_{f_c, \aGammap}^{\AGammap}\right).
\end{align*}
This implies in turn that
\begin{align*}
&F(\muGamma)\tau_{\pGamma}^{A} \circ \left(\phi_{\muGamma, \pGamma}^{\AGamma}\right)^{-1} \circ \tau_{f}^{\pi_2^{\ast}A} \circ F(\fGX)\quot{\theta}{\Gamma^{\prime}} \circ  F(\fGX)\left(\phi_{\muGammap, a_{\Gamma^{\prime}}}^{\AGammap}\right) \\
&\circ \left(\phi_{\fGX, \muGammap}^{F(\aGammap)\AGammap}\right)^{-1} \\ 
&=F(\muGamma)\left(\tau_{\aGamma}^{A} \circ F(\aGamma)\tau_{f}^{A} \circ \left(\phi_{\aGamma, \fX}^{\AGammap}\right)^{-1}\circ\phi_{f_c, \aGammap}^{\AGammap}\right)
\end{align*}
which allows us to deduce that
\begin{align*}
&\tau_f^{\pi_2^{\ast}A} \circ F(\fGX)\quot{\theta}{\Gamma^{\prime}}\\
&= \phi_{\muGamma, \pGamma}^{\AGamma} \circ F(\muGamma)\left(\left(\tau_{\pGamma}^{A}\right)^{-1} \circ \tau_{\aGamma}^{A}\right) \circ \big(F(\muGamma) \circ F(\aGamma)\big)\tau_f^A \circ F(\muGamma)\left(\phi_{\aGamma, \fX}^{\AGammap}\right)^{-1}\\
&\circ F(\muGamma)\phi_{f_c, \aGammap}^{\AGammap} \circ \left(\phi_{\muGamma, \of_c}^{F(\aGammap)\AGammap}\right)^{-1} \circ \phi_{\fGX, \muGammap}^{F(\aGammap)\AGammap} \circ  F(\fGX)\left(\phi_{\muGammap, a_{\Gamma^{\prime}}}^{\AGammap}\right)^{-1}.
\end{align*}
To complete the proof that $\Theta$ is an $F_G(X)$-morphism, it suffices at this point to show that the equality
\begin{align*}
&\big(F(\muGamma) \circ F(\aGamma)\big)\tau_f^A \circ F(\muGamma)\left(\phi_{\aGamma, \fX}^{\AGammap}\right)^{-1} \circ F(\muGamma)\phi_{f_c, \aGammap}^{\AGammap}  \\
&\circ \left(\phi_{\muGamma, \of_c}^{F(\aGammap)\AGammap}\right)^{-1} \circ \phi_{\fGX, \muGammap}^{F(\aGammap)\AGammap} \circ  F(\fGX)\left(\phi_{\muGammap, a_{\Gamma^{\prime}}}^{\AGammap}\right)^{-1} \\
&= \left(\phi_{\muGamma, \aGamma}^{\AGamma}\right)^{-1} \circ F(\oalphaGamma)\tau_f^A \circ \left(\phi_{\oalphaGamma, \fX}^{\AGammap}\right)^{-1} \circ \phi_{\fGX, \oalphaGammap}^{\AGammap}
\end{align*}
holds. However, this follows from an extremely tedious manipulation of the naturality of the compositors and the relations on the compositors as induced by the pseudofunctoriality of $F$ and Lemma \ref{Lemma: Clifton's Naivification Variety}. From this it follows that
\begin{align*}
&\tau_f^{\pi_2^{\ast}A} \circ F(\fGX)\quot{\theta}{\Gamma^{\prime}}\\
&= \phi_{\muGamma, \pGamma}^{\AGamma} \circ F(\muGamma)\left(\left(\tau_{\pGamma}^{A}\right)^{-1} \circ \tau_{\aGamma}^{A}\right) \circ \big(F(\muGamma) \circ F(\aGamma)\big)\tau_f^A \circ F(\muGamma)\left(\phi_{\aGamma, \fX}^{\AGammap}\right)^{-1}\\
&\circ F(\muGamma)\phi_{f_c, \aGammap}^{\AGammap} \circ \left(\phi_{\muGamma, \of_c}^{F(\aGammap)\AGammap}\right)^{-1} \circ \phi_{\fGX, \muGammap}^{F(\aGammap)\AGammap} \circ  F(\fGX)\left(\phi_{\muGammap, a_{\Gamma^{\prime}}}^{\AGammap}\right)^{-1} \\
&=\phi_{\muGamma, \pGamma}^{\AGamma} \circ F(\muGamma)\left(\left(\tau_{\pGamma}^{A}\right)^{-1} \circ \tau_{\aGamma}^{A}\right) \circ \left(\phi_{\muGamma, \aGamma}^{\AGamma}\right)^{-1} \circ F(\oalphaGamma)\tau_f^A \circ \left(\phi_{\oalphaGamma, \fX}^{\AGammap}\right)^{-1} \\
&\circ \phi_{\fGX, \oalphaGammap}^{\AGammap} \\
&= \quot{\theta}{\Gamma} \circ \tau_f^{\alpha_X^{\ast}A}.
\end{align*}
This proves, finally, that $\tau_f^{\pi_2^{\ast}A} \circ F(\fGX)\quot{\theta}{\Gamma^{\prime}} = \quot{\theta}{\Gamma} \circ \tau_f^{\alpha_X^{\ast}A}$ and hence shows that $\Theta \in F_G(X)_1$. That $\Theta$ is an isomorphism is immediate from each $\quot{\theta}{\Gamma}$ being an isomorphism, which proves that there is an isomorphism
\[
\alpha_X^{\ast}A \xrightarrow[\Theta]{\cong} \pi_2^{\ast}A
\]
in $F_G(G \times X)$ for any $A \in F_G(X)_0$.

Having finally proved that $\Theta$ is an $F_G(X)$-isomorphism, we now show that the $\Theta$ are natural in $A$. For this fix a $P \in F_G(X)_1$ with $P \in F_G(X)(A,B).$ Fix a $\Gamma \in \Sf(G)_0$ and note that we must check that the diagram
\[
\xymatrix{
	\alpha_X^{\ast}A \ar[d]_{\alpha_X^{\ast}P} \ar[r]^-{\Theta_A} & \pi_2^{\ast}A \ar[d]^{\pi_2^{\ast}P} \\
	\alpha_X^{\ast}B \ar[r]_-{\Theta_B} & \pi_2^{\ast}B
}
\]
commutes. We calculate
\begin{align*}
&F(\opiGamma)\rhoGamma \circ \quot{\theta_A}{\Gamma} \\
&= F(\opiGamma)\rhoGamma \circ \phi_{\muGamma, p_{\Gamma}}^{\AGamma} \circ F(\muGamma)\left(\left(\tau_{p_{\Gamma}}^{A}\right)^{-1} \circ \tau_{a_{\Gamma}}^{A}\right) \circ \left(\phi_{\muGamma, a_{\Gamma}}^{\AGamma}\right)^{-1} \\
&= \phi_{\muGamma, p_{\Gamma}}^{\BGamma} \circ F(\muGamma)\left(F(\pGamma)\rhoGamma \circ \left(\tau_{p_{\Gamma}}^{A}\right)^{-1} \circ \tau_{a_{\Gamma}}^{A}\right) \circ \left(\phi_{\muGamma, a_{\Gamma}}^{\AGamma}\right)^{-1}  \\
&=\phi_{\muGamma, p_{\Gamma}}^{\BGamma} \circ F(\muGamma)\left(\left(\tau_{p_{\Gamma}}^{B}\right)^{-1} \circ \rhoGamma \circ \tau_{a_{\Gamma}}^{A}\right) \circ \left(\phi_{\muGamma, a_{\Gamma}}^{\AGamma}\right)^{-1}  \\
&= \phi_{\muGamma, p_{\Gamma}}^{\BGamma} \circ F(\muGamma)\left(\left(\tau_{p_{\Gamma}}^{B}\right)^{-1} \circ \tau_{a_{\Gamma}}^{B}\right) \circ F(\muGamma)F(\aGamma)\rhoGamma \circ \left(\phi_{\muGamma, a_{\Gamma}}^{\AGamma}\right)^{-1} \\
&= \phi_{\muGamma, p_{\Gamma}}^{\BGamma} \circ F(\muGamma)\left(\left(\tau_{p_{\Gamma}}^{B}\right)^{-1} \circ \tau_{a_{\Gamma}}^{B}\right)  \circ \left(\phi_{\muGamma, a_{\Gamma}}^{\BGamma}\right)^{-1}\circ F(\aGamma \circ \muGamma)\rhoGamma \\
&= \phi_{\muGamma, p_{\Gamma}}^{\BGamma} \circ F(\muGamma)\left(\left(\tau_{p_{\Gamma}}^{B}\right)^{-1} \circ \tau_{a_{\Gamma}}^{B}\right)  \circ \left(\phi_{\muGamma, a_{\Gamma}}^{\BGamma}\right)^{-1}\circ F(\oalphaGamma)\rhoGamma \\
&= \quot{\theta_B}{\Gamma} \circ F(\oalphaGamma)\rhoGamma,
\end{align*}
which verifies the commutativity of the diagram.
\end{proof}
Before proving the cocycle condition that $\Theta$ satisfies, we display it as a definition so that we know what it is we must prove.
\begin{definition}\label{Defn: GIT Cocycle}\index{GIT Cocycle Condition}
Let $F$ be a pre-equivariant pseudofunctor on $X$, let $A \in F_G(X)_0$, and let $\Theta:\alpha_X^{\ast}A \to \pi_2^{\ast}A$ be an isomorphism in $F_G(G \times X)$. We say that $\Theta$ satisfies the GIT-cocycle condition if, for morphisms
\begin{align*}
d_0^1 &:= \pi_{23}:G \times G \times X \to G \times X; \\
d_1^1 &:= \mu_G \times \id_X:G \times G \times X \to G \times X; \\
d_2^1 &:= \id_G \times \alpha_X:G \times G \times X \to G \times X
\end{align*}\index[notation]{dFaceMaps@$d_0^1, d_1^1, d_2^1$}
for all $\Gamma \in \Sf(G)_0$ the composite
\begin{align*}
\phi_{\quot{d_1^1}{\Gamma}, \opiGamma}^{\AGamma} \circ F(\quot{d_1^1}{\Gamma})\quot{\theta}{\Gamma} \circ \left(\phi_{\quot{d_1^1}{\Gamma}, \oalphaGamma}^{\AGamma}\right)
\end{align*}
is equal to the composite
\begin{align*}
&\phi_{\quot{d_0^1}{\Gamma}, \opiGamma}^{\AGamma} \circ F(\quot{d_0^1}{\Gamma})\quot{\theta}{\Gamma} \circ \left(\phi_{\quot{d_0^1}{\Gamma}, \oalphaGamma}^{\AGamma}\right)^{-1} \circ \phi_{\quot{d_2^1}{\Gamma},\opiGamma}^{\AGamma} \circ F(\quot{d_2^1}{\Gamma})\quot{\theta}{\Gamma} \circ \left(\phi_{\quot{d_2^1}{\Gamma},\oalphaGamma}^{\AGamma}\right)^{-1}.
\end{align*}
\end{definition}

Let us now motivate the GIT cocycle condition by studying how it arose in \cite{GIT}. In \cite{GIT} the authors first asked for a sheaf $\Fscr$ on $X$ with an isomorphism
\[
\theta:\alpha_X^{\ast}\Fscr \xrightarrow{\cong} \pi_2^{\ast} \Fscr
\]
in $\Shv(G \times X)$. This isomorphism $\theta$ was used to illustrate the how the action of $G$ on $X$ affects the space $X$ and the sheaf $\Fscr$, as well as how these conditions need to be suitably compatible. Let us describe this more explicitly. The action $\alpha_X$ needs to intertwine the multiplication of the group $G$ on itself together with the action on $X$ in the sense that the diagram
\[
\xymatrix{
G \times G \times X \ar[rr]^-{\id_G \times \alpha_X} \ar[d]_{\mu_G \times \id_X} & & G \times X \ar[d]^{\alpha_X} \\
G \times X \ar[rr]_-{\alpha_X} & & X
}
\]
commutes. When looking at this at the level of stalks of $\Fscr$, this requires an isomorphism
\[
\Fscr_x \cong \Fscr_{(gh)x} \cong \Fscr_{g(hx)}
\]
for any stalk $(g,h,x) \in G \times G \times X$. However, these isomorphisms should come from the equivariance $\theta$ lifted to $G \times G \times X$ by either changing the order of multiplication,  action, and projection more-or-less freely (provided we do so in a way that makes sense). There are three maps from $G \times G \times X$ to $G \times X$ that we use to lift $\theta$; these three maps coincide with the simplicial presentation of both the quotient variety $\underline{G\backslash X}_{\bullet}$ (cf.\@ Lemma \ref{Lemma: Simplicial Scheme of X}) and the simplicial information induced by the internal action groupoid $\mathbb{G\times X}$ (cf.\@ Definition \ref{Defn: Essentially alg theory internal gropuoid} and Appendix \ref{Appendix: Locally closed subtoposes}) and that both these presentations appear here implicitly is the shadow of something quite interesting. These three maps, together with the two maps from $G \times X \to X$ in the diagram
\[
\begin{tikzcd}
G \times G \times X \ar[rr, shift left = 3]{}{\pi_{23}} \ar[rr]{}[description]{\mu_G\times \id_X} \ar[rr, shift right = 3, swap]{}{\id_G \times \alpha_X} & & G \times X \ar[rr, shift left]{}{\pi_2} \ar[rr, shift right, swap]{}{\alpha_X} & & X
\end{tikzcd}
\]
satisfy the equations
\begin{align*}
 \alpha_X \circ (\mu_G\times \id_X) &= \alpha_X \circ (\id_G \times \alpha_X); \\
 \alpha_X \circ \pi_{23} &= \pi_2 \circ (\id_G \times \alpha_X); \\
 \pi_2 \circ \pi_{23} &= \pi_2 \circ (\mu_G \times \id_X).
\end{align*}
Define, in order to save space and illustrate the simplicial connection, the morphisms
\begin{align*}
\pi_{23} &:= d_0^1; \\
\mu_G \times \id_X &:= d_1^1; \\
\id_G \times \alpha_X &:= d_2^1.
\end{align*}
Using these maps and relations we then can produce a diagram (assuming that the underlying pseudofunctor $\Shv(-)$ is strict in this case)
\[
\begin{tikzcd}
(d_2^1)^{\ast}\big(\alpha_X^{\ast}\Fscr\big) \ar[r]{}{(d_2^1)^{\ast}\theta} \ar[d, equals] & (d_2^1)^{\ast}\big(\pi_{2}^{\ast}\Fscr\big) \ar[r, equals] & \left( \pi_2 \circ d_2^1\right)^{\ast}\Fscr \ar[d, equals] \\
(\alpha_X \circ d_2^1)^{\ast}\Fscr \ar[d, equals] & & (\alpha_X \circ d_0^1)^{\ast}\Fscr \ar[d, equals] \\
(\alpha_X \circ d_1^1)^{\ast}\Fscr \ar[d, equals] & & (d_0^1)^{\ast}\big(\alpha_X^{\ast}\Fscr\big) \ar[d]{}{(d_0^1)^{\ast}\theta} \\
(d_1^1)^{\ast}\alpha_X^{\ast}\Fscr \ar[dd, swap]{}{(d_1^1)^{\ast}\theta} &  & (d_0^1)^{\ast}\big(\pi_2^{\ast}\Fscr\big) \ar[d, equals] \\
&  & (\pi_2 \circ d_0^1)^{\ast}\Fscr \ar[d, equals] \\
(d_1^1)^{\ast}\big(\pi_2^{\ast}\Fscr\big)& & \ar[ll, equals] (\pi_2 \circ d_1^1)^{\ast}\Fscr
\end{tikzcd}
\]
The cocylce condition of \cite{GIT} asks exactly that the above diagram commute, i.e.,
\[
(d_1^1)^{\ast}\theta = (d_0^1)^{\ast}\theta \circ (d_2^1)^{\ast}\theta.
\] 
Our pseudofunctorial GIT-cocycle condition then asks for a pre-equivariant pseudofunctor $F$ over $X$ that for any $A \in F_G(X)_0$ and $\Gamma \in \Sf(G)_0$, the diagram
\[
\begin{tikzcd}
F(\oalphaGamma \circ \odGamma{2}{1})(\AGamma) \ar[rr]{}{\left(\phi_{\odGamma{2}{1}, \oalphaGamma}^{\AGamma}\right)^{-1}} \ar[d, equals] & & F(\odGamma{2}{1})\big(F(\oalphaGamma)(\AGamma)\big) \ar[d]{}{F(\odGamma{2}{1})\quot{\theta}{\Gamma}}  \\
F(\oalphaGamma \circ \odGamma{1}{1})(\AGamma) \ar[d, swap]{}{\left(\phi_{\odGamma{1}{1}, \oalphaGamma}^{\AGamma}\right)^{-1}} & & F(\odGamma{2}{1})\big(F(\opiGamma)(\AGamma)\big) \ar[d]{}{\phi_{\odGamma{2}{1}, \opiGamma}^{\AGamma}} \\
F(\odGamma{1}{1})\big(F(\oalphaGamma)(\AGamma)\big) \ar[d, swap]{}{F(\odGamma{1}{1})\quot{\theta}{\Gamma}} & & F(\opiGamma \circ \odGamma{2}{1}) \ar[d, equals] \\
F(\odGamma{1}{1})\big(F(\opiGamma)(\AGamma)\big) \ar[d, swap]{}{\phi_{\odGamma{1}{1}, \opiGamma}^{\AGamma}} & &  F(\oalphaGamma \circ \odGamma{0}{1})(\AGamma) \ar[d]{}{\left(\phi_{\odGamma{0}{1}, \oalphaGamma}^{\AGamma}\right)^{-1}} \\
F(\opiGamma \circ \odGamma{1}{1})(\AGamma) & & F(\odGamma{0}{1})\big(F(\oalphaGamma)(\AGamma)\big) \ar[d]{}{F(\odGamma{0}{1})\quot{\theta}{\Gamma}} \\
F(\opiGamma \circ \odGamma{0}{1})(\AGamma) \ar[u, equals] & & F(\odGamma{0}{1})\big(F(\opiGamma)(\AGamma)\big) \ar[ll]{}{\phi_{\odGamma{0}{1},\opiGamma}^{\AGamma}}
\end{tikzcd}
\]
commutes. Note that while we have not proved it explicitly, it follows mutatis mutandis to Lemmas \ref{Lemma: Section Change of Space: pi2 pullback} and \ref{Lemma: Section Change of Space: Action map pullback} that $F$ admits descent pullbacks against all of $d_0^1,$ $d_1^1,$ and $d_2^1$. This gives rise to induce the corresponding pullback functors
\[
(d_0^1)^{\ast}, (d_1^1)^{\ast}, (d_2^1)^{\ast}:F_G(G \times X) \to F_G(G \times G \times X).
\]

As we move to show that the isomorphism $\Theta$ of Theorem \ref{Theorem: Formal naive equivariance} satisfies the GIT cocycle condtion (cf.\@ Theorem \ref{Thm: GIT Cocycle Condition Formal Equivariance} below), we will need a short lemma on the interaction of the maps $\odGamma{2}{1}$ and $\odGamma{0}{1}$ with $\muGamma$.
\begin{lemma}\label{Lemma: MuGamma maps same with odgammas}
For any $\Gamma \in \Sf(G)_0$, if $\muGamma$ is the isomorphism of Lemma \ref{Lemma: Clifton's Naivification Variety}, then
\[
\muGamma \circ \odGamma{0}{1} = \muGamma \circ \odGamma{2}{1}.
\]
\end{lemma}
\begin{proof}
Observe that since
\[
\quot{d_0^1}{\Gamma} = \pi_{134}:\Gamma \times G \times G \times X \to \Gamma \times G \times X
\]
and
\[
\quot{d_2^1}{\Gamma} = \id_{\Gamma \times G} \times \alpha_{X}:\Gamma \times G \times G \times X \to \Gamma \times G \times X
\]
there are unique morphisms making the diagrams
\[
\xymatrix{
\Gamma \times G \times G \times X \ar@<.5ex>[r]^-{\quot{d_2^1}{\Gamma}} \ar@<-.5ex>[r]_-{\quot{d_0^1}{\Gamma}} \ar[d]_{\quo^{\Gamma\times G \times G}} & \Gamma \times G \times X \ar[d]^{\quo_{G \times \Gamma}} \\
\quot{(G \times G \times X)}{\Gamma} \ar@<.5ex>[r]^-{\odGamma{2}{1}} \ar@<-.5ex>[r]_-{\odGamma{0}{1}} & \quot{(G \times X)}{\Gamma}
}
\]
\[
\xymatrix{
\Gamma \times G \times G \times X \ar@<.5ex>[r]^-{\quot{d_2^1}{\Gamma}} \ar@<-.5ex>[r]_-{\quot{d_0^1}{\Gamma}} \ar[d]_{\quo_{(\Gamma \times G)_{c}}} & \Gamma \times G \times X \ar[d]^{\quo_{\Gamma_c}} \\
\quot{X}{(\Gamma \times G)_c} \ar@<.5ex>[r]^-{(\odGamma{2}{1})_c} \ar@<-.5ex>[r]_-{(\odGamma{0}{1})_c} & \quot{X}{\Gamma_c}
}
\]
commute. However, since the second diagram factors through the diagrams
\[
\xymatrix{
\Gamma \times G \times G \times X \ar@<.5ex>[r]^-{\quot{d_2^1}{\Gamma}} \ar@<-.5ex>[r]_-{\quot{d_0^1}{\Gamma}} \ar[d]_{\quo_{(\Gamma \times G)_{c}}} & \Gamma \times G \times X \ar[d]^{\quo_{\Gamma_c}} \\
\quot{(G \times G \times X)}{\Gamma} \ar@<.5ex>[r]^-{\odGamma{2}{1}} \ar@<-.5ex>[r]_-{\odGamma{0}{1}} \ar[d]_{\mu_{\Gamma \times G}} & \quot{(G \times X)}{\Gamma} \ar[d]^{\muGamma} \\
	\quot{X}{(\Gamma \times G)_c} \ar@<.5ex>[r]^-{(\odGamma{2}{1})_c} \ar@<-.5ex>[r]_-{(\odGamma{0}{1})_c} & \quot{X}{\Gamma_c}
}
\]
the result follows.
\end{proof}
\begin{corollary}\label{Cor: F of muG circ d01 is F of muG circ d21}
For any pre-equivariant pseudofunctor $F$ over $X$,
\[
F(\muGamma \circ \odGamma{0}{1}) = F(\muGamma \circ \odGamma{2}{1}).
\]
\end{corollary}
\begin{Theorem}\label{Thm: GIT Cocycle Condition Formal Equivariance}
Let $F$ be a pre-equivariant pseudofunctor over $X$ and let $A \in F_G(X)_0$. Then the equivariance $\Theta$ of Theorem \ref{Theorem: Formal naive equivariance} satisfies the GIT cocycle condition.
\end{Theorem}
\begin{proof}
Begin by noting that from Theorem \ref{Theorem: Formal naive equivariance} we must prove that
\begin{align*}
&\phi_{\quot{d_1^1}{\Gamma}, \opiGamma}^{\AGamma} \circ F(\quot{d_1^1}{\Gamma})\quot{\theta}{\Gamma} \circ \left(\phi_{\quot{d_1^1}{\Gamma}, \oalphaGamma}^{\AGamma}\right) \\
&= \phi_{\quot{d_1^1}{\Gamma}, \opiGamma}^{\AGamma} \circ F(\quot{d_1^1}{\Gamma})\left(\phi_{\muGamma, \pGamma}^{\AGamma}\circ F(\muGamma)\big(\big(\tau_{\pGamma}^{A}\big)^{-1} \circ \tau_{\aGamma}^{A}\big) \circ \left(\phi_{\muGamma, \aGamma}^{\AGamma}\right)^{-1}\right) \\
&\circ \left(\phi_{\quot{d_1^1}{\Gamma}, \oalphaGamma}^{\AGamma}\right)
\end{align*}
is equal to the composite
\begin{align*}
&\phi_{\quot{d_0^1}{\Gamma}, \opiGamma}^{\AGamma} \circ F(\quot{d_0^1}{\Gamma})\quot{\theta}{\Gamma} \circ \left(\phi_{\quot{d_0^1}{\Gamma}, \oalphaGamma}^{\AGamma}\right)^{-1} \circ \phi_{\quot{d_2^1}{\Gamma},\opiGamma}^{\AGamma} \circ F(\quot{d_2^1}{\Gamma})\quot{\theta}{\Gamma} \circ \left(\phi_{\quot{d_2^1}{\Gamma},\oalphaGamma}^{\AGamma}\right)^{-1} \\
&= \phi_{\quot{d_0^1}{\Gamma}, \opiGamma}^{\AGamma} \circ F(\odGamma{2}{1})\phi_{\muGamma, \pGamma}^{\AGamma} \circ \big(F(\odGamma{2}{1})\circ F(\muGamma)\big)\left(\big(\tau_{\pGamma}^{A}\big)^{-1} \circ \tau_{\aGamma}^A\right) \\
&\circ F(\odGamma{0}{1})\left(\phi_{\muGamma, \aGamma}^{\AGamma}\right)^{-1} \circ \left(\phi_{\quot{d_0^1}{\Gamma}, \oalphaGamma}^{\AGamma}\right)^{-1} \circ \phi_{\quot{d_2^1}{\Gamma},\opiGamma}^{\AGamma} \circ F(\odGamma{2}{1})\phi_{\muGamma, \pGamma}^{\AGamma} \\
&\circ \big(F(\odGamma{2}{1}) \circ F(\muGamma)\big)\left(\big(\tau_{\pGamma}^{A}\big)^{-1} \circ \tau_{\aGamma}^{A}\right) \circ F(\odGamma{2}{1})\left(\phi_{\muGamma, \aGamma}^{\AGamma}\right)^{-1} \circ \left(\phi_{\quot{d_2^1}{\Gamma},\oalphaGamma}^{\AGamma}\right)^{-1}.
\end{align*}
We will prove this by systematically attacking and simplifying the second expression. Begin by recalling that either by Lemma \ref{Lemma: Clifton's Naivification Variety} or the exposition prior to the Lemma, the identities
\begin{align*}
\opiGamma &= \pGamma \circ \muGamma; \\
\oalphaGamma &= \aGamma \circ \muGamma.
\end{align*}
both hold. Now, using these and the identities from before the statment of the theorem, we observe as in the proof of Thoerem \ref{Theorem: Formal naive equivariance} that the pasting diagram
\[
\begin{tikzcd}
F(\XGamma) \ar[r, bend left = 30, ""{name = U}]{}{F(\muGamma) \circ F(\pGammap)} \ar[r, bend right = 30, swap, ""{name = M}]{}{F(\opiGamma)} \ar[rr, bend right = 60, swap, ""{name = L}]{}{F(\pGamma \circ \muGamma \circ \odGamma{0}{1})} & F(\quot{(G \times X)}{\Gamma}) \ar[r]{}{F(\odGamma{0}{1})} & F(\quot{(G \times G \times X)}{\Gamma}) \ar[from = U, to = M, Rightarrow, shorten <= 4pt, shorten >= 4pt]{}{} \ar[from = 1-2, to = L, Rightarrow, shorten <= 4pt, shorten >= 4pt]{}{}
\end{tikzcd}
\]
is equal to the pasting diagram
\[
\begin{tikzcd}
F(\XGamma) \ar[r]{}{F(\pGamma)} \ar[rr, bend right = 60, swap, ""{name = L}]{}{F(\pGamma \circ \muGamma \circ \odGamma{0}{1})} & F(\quot{X}{\Gamma_c}) \ar[r, bend left = 30, ""{name = U}]{}{F(\odGamma{0}{1}) \circ F(\muGamma)} \ar[r, bend right = 30, swap, ""{name = M}]{}{F(\muGamma \circ \odGamma{0}{1})} & F(\quot{(G \times G \times X)}{\Gamma}) \ar[from = U, to = M, Rightarrow, shorten <= 4pt, shorten >= 4pt] \ar[from = 1-2, to = L, Rightarrow, shorten <= 4pt, shorten >= 4pt]{}{}
\end{tikzcd}
\]
where the the $2$-cells are again the compositor isomorphisms. This allows us to deduce that
\[
\phi_{\quot{d_0^1}{\Gamma}, \opiGamma}^{\AGamma} \circ F(\odGamma{0}{1})\phi_{\muGamma, \pGamma}^{\AGamma} = \phi_{\muGamma \circ \odGamma{0}{1}, \pGamma}^{\AGamma} \circ \phi_{\odGamma{0}{1}, \muGamma}^{F(\pGamma)\AGamma}.
\]
Similarly, using that
\[
F(\odGamma{0}{1})\left(\phi_{\muGamma, \aGamma}^{\AGamma}\right)^{-1} \circ \left(\phi_{\quot{d_0^1}{\Gamma}, \oalphaGamma}^{\AGamma}\right)^{-1} = \left(\phi_{\odGamma{0}{1}, \oalphaGamma}^{\AGamma} \circ F(\odGamma{0}{1})\phi_{\muGamma, \aGamma}^{\AGamma} \circ \right)^{-1}
\]
and proceeding to analyze the pasting diagrams as before, we get that
\[
F(\odGamma{0}{1})\left(\phi_{\muGamma, \aGamma}^{\AGamma}\right)^{-1} \circ \left(\phi_{\quot{d_0^1}{\Gamma}, \oalphaGamma}^{\AGamma}\right)^{-1} = \left(\phi_{\muGamma,\odGamma{0}{1}}^{F(\aGamma)}\right)^{-1} \circ \left(\phi_{\muGamma \circ \odGamma{0}{1}, \aGamma}^{\AGamma}\right)^{-1}.
\]
We also derive that
\[
\phi_{\quot{d_2^1}{\Gamma}, \opiGamma}^{\AGamma} \circ F(\odGamma{2}{1})\phi_{\muGamma, \pGamma}^{\AGamma} = \phi_{\muGamma \circ \odGamma{2}{1}, \pGamma}^{\AGamma} \circ \phi_{\odGamma{2}{1}, \muGamma}^{F(\pGamma)\AGamma}
\]
and
\[
F(\odGamma{2}{1})\left(\phi_{\muGamma, \aGamma}^{\AGamma}\right)^{-1} \circ \left(\phi_{\quot{d_2^1}{\Gamma}, \oalphaGamma}^{\AGamma}\right)^{-1} = \left(\phi_{\muGamma,\odGamma{2}{1}}^{F(\aGamma)}\right)^{-1} \circ \left(\phi_{\muGamma \circ \odGamma{2}{1}, \aGamma}^{\AGamma}\right)^{-1}
\]
through the same technique. Putting these together we get that
\begin{align*}
&\phi_{\quot{d_0^1}{\Gamma}, \opiGamma}^{\AGamma} \circ F(\odGamma{2}{1})\phi_{\muGamma, \pGamma}^{\AGamma} \circ \big(F(\odGamma{2}{1})\circ F(\muGamma)\big)\left(\big(\tau_{\pGamma}^{A}\big)^{-1} \circ \tau_{\aGamma}^A\right) \\
&\circ F(\odGamma{0}{1})\left(\phi_{\muGamma, \aGamma}^{\AGamma}\right)^{-1} \circ \left(\phi_{\quot{d_0^1}{\Gamma}, \oalphaGamma}^{\AGamma}\right)^{-1} \circ \phi_{\quot{d_2^1}{\Gamma},\opiGamma}^{\AGamma} \circ F(\odGamma{2}{1})\phi_{\muGamma, \pGamma}^{\AGamma} \\
&\circ \big(F(\odGamma{2}{1}) \circ F(\muGamma)\big)\left(\big(\tau_{\pGamma}^{A}\big)^{-1} \circ \tau_{\aGamma}^{A}\right) \circ F(\odGamma{2}{1})\left(\phi_{\muGamma, \aGamma}^{\AGamma}\right)^{-1} \circ \left(\phi_{\quot{d_2^1}{\Gamma},\oalphaGamma}^{\AGamma}\right)^{-1} \\
&= \phi_{\muGamma \circ \odGamma{0}{1}, \pGamma}^{\AGamma} \circ \phi_{\odGamma{0}{1}, \muGamma}^{F(\pGamma)\AGamma} \circ \big(F(\odGamma{2}{1})\circ F(\muGamma)\big)\left(\big(\tau_{\pGamma}^{A}\big)^{-1} \circ \tau_{\aGamma}^A\right) \\
& \circ \left(\phi_{\muGamma,\odGamma{0}{1}}^{F(\aGamma)}\right)^{-1} \circ \left(\phi_{\muGamma \circ \odGamma{0}{1}, \aGamma}^{\AGamma}\right)^{-1} \circ \phi_{\muGamma \circ \odGamma{2}{1}, \pGamma}^{\AGamma} \circ \phi_{\odGamma{2}{1}, \muGamma}^{F(\pGamma)\AGamma} \\
&\circ \big(F(\odGamma{2}{1}) \circ F(\muGamma)\big)\left(\big(\tau_{\pGamma}^{A}\big)^{-1} \circ \tau_{\aGamma}^{A}\right) \circ \left(\phi_{\muGamma,\odGamma{2}{1}}^{F(\aGamma)}\right)^{-1} \circ \left(\phi_{\muGamma \circ \odGamma{2}{1}, \aGamma}^{\AGamma}\right)^{-1}.
\end{align*}

To simplify the expression above, we now focus on manipulating the transition isomorphisms and their interactions with the compositor natural isomorphisms. To begin this process not that since the diagram
\[
\begin{tikzcd}
\big(F(\odGamma{0}{1})\circ F(\muGamma)\big)\big(F(\aGamma)(\AGamma)\big) \ar[rrr]{}{\big(F(\odGamma{0}{1}) \circ F(\muGamma)\big)\tau_{\aGamma}^{A}} & & & \big(F(\odGamma{0}{1}) \circ F(\muGamma)\big)(\quot{A}{\Gamma_c}) \\ F(\muGamma \circ \odGamma{0}{1})\big(F(\aGamma)(\AGamma)\big) \ar[u]{}{\left(\phi_{\odGamma{0}{1}, \muGamma}^{\AGamma}\right)^{-1}} \ar[rrr, swap]{}{F(\muGamma \circ \odGamma{0}{1})\tau_{\aGamma}^{A}} & & & F(\muGamma \circ \odGamma{0}{1})(\quot{A}{\Gamma_c}) \ar[u, swap]{}{\left(\phi_{\odGamma{0}{1}, \muGamma}^{\quot{A}{\Gamma_c}}\right)^{-1}}
\end{tikzcd}
\]
the equation
\[
\big(F(\odGamma{0}{1}) \circ F(\muGamma)\big)\tau_{\aGamma}^{A} \circ \left(\phi_{\odGamma{0}{1}, \muGamma}^{\AGamma}\right)^{-1} = \left(\phi_{\odGamma{0}{1}, \muGamma}^{\quot{A}{\Gamma_c}}\right)^{-1} \circ F(\muGamma \circ \odGamma{0}{1})\tau_{\aGamma}^{A}
\]
holds. Similarly, the naturality of the compositors and their inverses allows us to derive that
\[
\big(F(\odGamma{0}{1}) \circ F(\muGamma)\big)(\tau_{\pGamma}^{A})^{-1} \circ \left(\phi_{\odGamma{0}{1}, \muGamma}^{\quot{A}{\Gamma_c}}\right)^{-1} = \left(\phi_{\odGamma{0}{1}, \muGamma}^{F(\pGamma)\AGamma}\right)^{-1} \circ F(\muGamma \circ \odGamma{0}{1})(\tau_{\pGamma}^{A})^{-1}
\]
which allows us to further calculate that
\begin{align*}
&\phi_{\muGamma \circ \odGamma{0}{1}, \pGamma}^{\AGamma} \circ \phi_{\odGamma{0}{1}, \muGamma}^{F(\pGamma)\AGamma} \circ \big(F(\odGamma{2}{1})\circ F(\muGamma)\big)\left(\big(\tau_{\pGamma}^{A}\big)^{-1} \circ \tau_{\aGamma}^A\right)\circ \left(\phi_{\muGamma,\odGamma{0}{1}}^{F(\aGamma)}\right)^{-1} \\
&= \phi_{\muGamma \circ \odGamma{0}{1}, \pGamma}^{\AGamma} \circ \phi_{\odGamma{0}{1}, \muGamma}^{F(\pGamma)\AGamma} \circ \left(\phi_{\odGamma{0}{1}, \muGamma}^{F(\pGamma)\AGamma}\right)^{-1} \circ F(\muGamma \circ \odGamma{0}{1})\big((\tau_{\pGamma}^{A})^{-1} \circ \tau_{\aGamma}^{A}\big) \\
&= \phi_{\muGamma \circ \odGamma{0}{1}, \pGamma}^{\AGamma} \circ F(\muGamma \circ \odGamma{0}{1})\big((\tau_{\pGamma}^{A})^{-1} \circ \tau_{\aGamma}^{A}\big).
\end{align*}
Now by using Corollary \ref{Cor: F of muG circ d01 is F of muG circ d21} we get that 
\[
F(\muGamma \circ \odGamma{0}{1})(\quot{A}{\Gamma_c}) = F(\muGamma \circ \odGamma{2}{1})(\quot{A}{\Gamma_c}).
\]
Using this we derive that the diagram
\[
\begin{tikzcd}
F(\aGamma \circ \muGamma \circ \odGamma{0}{1})(\AGamma) \ar[rr, equals] \ar[d, swap]{}{\left(\phi_{\muGamma \circ \odGamma{0}{1}}^{\AGamma}\right)^{-1}} & & F(\pGamma \circ \muGamma \circ \odGamma{2}{1})(\AGamma) \ar[d]{}{\left(\phi_{\muGamma \circ \odGamma{2}{1}}^{\AGamma}\right)^{-1}} \\
F(\muGamma \circ \odGamma{0}{1})\big(F(\aGamma)(\AGamma)\big) \ar[d, swap]{}{F(\muGamma \circ \odGamma{0}{1})} & & F(\muGamma \circ \odGamma{2}{1})\big(F(\pGamma)(\AGamma)\big) \ar[d]{}{F(\muGamma \circ \odGamma{2}{1})\tau_{\pGamma}^{A}} \ar[ll]{}{\left(\phi_{\muGamma \circ \odGamma{0}{1}}^{\AGamma}\right)^{-1} \circ \phi_{\muGamma \circ \odGamma{2}{1}}^{\AGamma}} \\
F(\muGamma \circ \odGamma{0}{1})(\quot{A}{\Gamma_c}) \ar[rr, equals] & & F(\muGamma \circ \odGamma{2}{1})(\quot{A}{\Gamma_c})
\end{tikzcd}
\]
commutes which tells us that
\[
F(\muGamma \circ \odGamma{2}{1})\tau_{\pGamma}^{A} = F(\muGamma \circ \odGamma{0}{1})\tau_{\aGamma}^{A} \circ \left(\phi_{\muGamma \circ \odGamma{0}{1}}^{\AGamma}\right)^{-1} \circ \phi_{\muGamma \circ \odGamma{2}{1}}^{\AGamma}.
\]
This allows us to compute that
\begin{align*}
&F(\muGamma \circ \odGamma{0}{1})\tau_{\aGamma}^{A} \circ \left(\phi_{\muGamma \circ \odGamma{0}{1}}^{\AGamma}\right)^{-1} \circ \phi_{\muGamma \circ \odGamma{2}{1}}^{\AGamma} \circ \phi_{\odGamma{2}{1}, \muGamma}^{F(\pGamma)\AGamma} \\
&\circ \big(F(\odGamma{2}{1}) \circ F(\muGamma)\big)(\tau_{\pGamma}^{A})^{-1} \\
&=F(\muGamma \circ \odGamma{2}{1})\tau_{\pGamma}^{A} \circ \phi_{\muGamma \circ \odGamma{2}{1}}^{\AGamma} \circ \big(F(\odGamma{2}{1}) \circ F(\muGamma)\big)(\tau_{\pGamma}^{A})^{-1} \\
&=F(\muGamma \circ \odGamma{2}{1})\tau_{\pGamma}^{A} F(\muGamma \circ \odGamma{2}{1})(\tau_{\pGamma}^{A})^{-1} \circ \phi_{\muGamma \circ \odGamma{2}{1}}^{\quot{A}{\Gamma_c}} \\
&= \phi_{\muGamma \circ \odGamma{2}{1}}^{\quot{A}{\Gamma_c}}.
\end{align*}
Using these various calculations we find that
\begin{align*}
& \phi_{\muGamma \circ \odGamma{0}{1}, \pGamma}^{\AGamma} \circ \phi_{\odGamma{0}{1}, \muGamma}^{F(\pGamma)\AGamma} \circ \big(F(\odGamma{2}{1})\circ F(\muGamma)\big)\left(\big(\tau_{\pGamma}^{A}\big)^{-1} \circ \tau_{\aGamma}^A\right) \\
& \circ \left(\phi_{\muGamma,\odGamma{0}{1}}^{F(\aGamma)}\right)^{-1} \circ \left(\phi_{\muGamma \circ \odGamma{0}{1}, \aGamma}^{\AGamma}\right)^{-1} \circ \phi_{\muGamma \circ \odGamma{2}{1}, \pGamma}^{\AGamma} \circ \phi_{\odGamma{2}{1}, \muGamma}^{F(\pGamma)\AGamma} \\
&\circ \big(F(\odGamma{2}{1}) \circ F(\muGamma)\big)\left(\big(\tau_{\pGamma}^{A}\big)^{-1} \circ \tau_{\aGamma}^{A}\right) \circ \left(\phi_{\muGamma,\odGamma{2}{1}}^{F(\aGamma)}\right)^{-1} \circ \left(\phi_{\muGamma \circ \odGamma{2}{1}, \aGamma}^{\AGamma}\right)^{-1} \\
&= \phi_{\muGamma \circ \odGamma{0}{1}, \pGamma}^{\AGamma} \circ F(\muGamma \circ \odGamma{0}{1})\big((\tau_{\pGamma}^{A})^{-1} \circ \tau_{\aGamma}^{A}\big)  \circ \left(\phi_{\muGamma \circ \odGamma{0}{1}, \aGamma}^{\AGamma}\right)^{-1} \circ \phi_{\muGamma \circ \odGamma{2}{1}, \pGamma}^{\AGamma}  \\
&\circ \phi_{\odGamma{2}{1}, \muGamma}^{F(\pGamma)\AGamma} \circ \big(F(\odGamma{2}{1}) \circ F(\muGamma)\big)\left(\big(\tau_{\pGamma}^{A}\big)^{-1} \circ \tau_{\aGamma}^{A}\right) \circ \left(\phi_{\muGamma,\odGamma{2}{1}}^{F(\aGamma)}\right)^{-1} \\
&\circ \left(\phi_{\muGamma \circ \odGamma{2}{1}, \aGamma}^{\AGamma}\right)^{-1} \\
&= \phi_{\muGamma \circ \odGamma{0}{1}, \pGamma}^{\AGamma} \circ F(\muGamma \circ \odGamma{0}{1})(\tau_{\pGamma}^{A})^{-1} \circ \phi_{\odGamma{2}{1}, \muGamma}^{\quot{A}{\Gamma_c}} \circ \big(F(\odGamma{2}{1}) \circ F(\muGamma)\big)\tau_{\aGamma}^{A}  \\
&\circ \left(\phi_{\muGamma,\odGamma{2}{1}}^{F(\aGamma)}\right)^{-1}\circ \left(\phi_{\muGamma \circ \odGamma{2}{1}, \aGamma}^{\AGamma}\right)^{-1}  \\
&=\phi_{\muGamma \circ \odGamma{0}{1}, \pGamma}^{\AGamma} \circ F(\muGamma \circ \odGamma{0}{1})(\tau_{\pGamma}^{A})^{-1} \circ F(\muGamma \circ \odGamma{2}{1})\tau_{\aGamma}^{A} \circ \phi_{\odGamma{2}{1}, \muGamma}^{F(\aGamma)\AGamma} \\
&\circ \left(\phi_{\muGamma,\odGamma{2}{1}}^{F(\aGamma)}\right)^{-1}\circ \left(\phi_{\muGamma \circ \odGamma{2}{1}, \aGamma}^{\AGamma}\right)^{-1}  \\
&= \phi_{\muGamma \circ \odGamma{0}{1}, \pGamma}^{\AGamma} \circ F(\muGamma \circ \odGamma{0}{1})(\tau_{\pGamma}^{A})^{-1} \circ F(\muGamma \circ \odGamma{2}{1})\tau_{\aGamma}^{A} \circ \left(\phi_{\muGamma \circ \odGamma{2}{1}, \aGamma}^{\AGamma}\right)^{-1}.
\end{align*}

To further manipulate the expression above into our desired format, we will need some calculations to translate between the various compositions of $F(\pGamma \circ \muGamma \circ \odGamma{0}{1})$ and $F(\pGamma \circ \muGamma \circ \odGamma{1}{1})$ (and similarly with $F(\aGamma \circ \muGamma \circ \odGamma{2}{1})$ and $F(\aGamma \circ \muGamma \circ \odGamma{1}{1})$). Note that because
\[
F(\pGamma \circ \muGamma \circ \odGamma{0}{1}) = F(\pGamma \circ \muGamma \circ \odGamma{1}{1})
\]
we have the $2$-commutative diagram:
\[
\begin{tikzcd}
 &F(\quot{X}{\Gamma_c}) \ar[ddr, near start, ""{name = UpDiag}]{}{F(\muGamma \circ \odGamma{0}{1})} \ar[rr, ""{name = UpRight}]{}{F(\muGamma)} & & F(\quot{(G \times X)}{\Gamma}) \ar[ddl, ""{name = UpUp}]{}{F(\odGamma{0}{1})} \\
\\
F(\XGamma) \ar[ddr, swap]{}{F(\pGamma)} \ar[uur]{}{F(\pGamma)} \ar[rr, bend left = 20, ""{name = MidUp}]{}{F(\pGamma \circ \muGamma \circ \odGamma{0}{1})} \ar[rr, bend right = 20, swap, ""{name = MidLow}]{}{F(\pGamma \circ \muGamma \circ \odGamma{1}{1})} & & F(\quot{(G \times G \times X)}{\Gamma}) \\
\\
 & F(\quot{X}{\Gamma_c}) \ar[uur, swap, near start, ""{name = DownDiag}]{}{F(\muGamma \circ \odGamma{1}{1})} \ar[rr, swap, ""{name = DownRight}]{}{F(\muGamma)} & & F(\quot{(G \times X)}{\Gamma}) \ar[uul, swap]{}{F(\odGamma{1}{1})}
\ar[from = MidUp, to= MidLow, equals, shorten <= 4pt, shorten >= 4pt]
\ar[from = 1-2, to = MidUp, Rightarrow, shorten <= 4pt, shorten >= 8pt]
\ar[from = UpRight, to = 3-3, Rightarrow, shorten <= 4pt, shorten >= 4pt]
\ar[from = 5-2, to = MidLow, Rightarrow, shorten <= 4pt, shorten >= 8pt]
\ar[from = DownRight, to = 3-3, Rightarrow, shorten <= 4pt, shorten >= 4pt]
\end{tikzcd}
\]
These allow us to deduce that the (equivalent) pasting diagrams
\[
\begin{tikzcd}
F(\XGamma) \ar[r, bend left = 30, ""{name = U}]{}{F(\muGamma) \circ F(\pGammap)} \ar[r, bend right = 30, swap, ""{name = M}]{}{F(\opiGamma)} \ar[rr, bend right = 60, swap, ""{name = L}]{}{F(\pGamma \circ \muGamma \circ \odGamma{0}{1})} & F(\quot{(G \times X)}{\Gamma}) \ar[r]{}{F(\odGamma{0}{1})} & F(\quot{(G \times G \times X)}{\Gamma}) \ar[from = U, to = M, Rightarrow, shorten <= 4pt, shorten >= 4pt]{}{} \ar[from = 1-2, to = L, Rightarrow, shorten <= 4pt, shorten >= 4pt]{}{}
\end{tikzcd}
\]
\[
\begin{tikzcd}
F(\XGamma) \ar[r]{}{F(\pGamma)} \ar[rr, bend right = 60, swap, ""{name = L}]{}{F(\pGamma \circ \muGamma \circ \odGamma{0}{1})} & F(\quot{X}{\Gamma_c}) \ar[r, bend left = 30, ""{name = U}]{}{F(\odGamma{0}{1}) \circ F(\muGamma)} \ar[r, bend right = 30, swap, ""{name = M}]{}{F(\muGamma \circ \odGamma{0}{1})} & F(\quot{(G \times G \times X)}{\Gamma}) \ar[from = U, to = M, Rightarrow, shorten <= 4pt, shorten >= 4pt] \ar[from = 1-2, to = L, Rightarrow, shorten <= 4pt, shorten >= 4pt]{}{}
\end{tikzcd}
\]
are equal to the (also equivalent) pasting diagrams
\[
\begin{tikzcd}
F(\XGamma) \ar[r, bend left = 30, ""{name = U}]{}{F(\muGamma) \circ F(\pGammap)} \ar[r, bend right = 30, swap, ""{name = M}]{}{F(\opiGamma)} \ar[rr, bend right = 60, swap, ""{name = L}]{}{F(\pGamma \circ \muGamma \circ \odGamma{1}{1})} & F(\quot{(G \times X)}{\Gamma}) \ar[r]{}{F(\odGamma{1}{1})} & F(\quot{(G \times G \times X)}{\Gamma}) \ar[from = U, to = M, Rightarrow, shorten <= 4pt, shorten >= 4pt]{}{} \ar[from = 1-2, to = L, Rightarrow, shorten <= 4pt, shorten >= 4pt]{}{}
\end{tikzcd}
\]
\[
\begin{tikzcd}
F(\XGamma) \ar[r]{}{F(\pGamma)} \ar[rr, bend right = 60, swap, ""{name = L}]{}{F(\pGamma \circ \muGamma \circ \odGamma{1}{1})} & F(\quot{X}{\Gamma_c}) \ar[r, bend left = 30, ""{name = U}]{}{F(\odGamma{1}{1}) \circ F(\muGamma)} \ar[r, bend right = 30, swap, ""{name = M}]{}{F(\muGamma \circ \odGamma{1}{1})} & F(\quot{(G \times G \times X)}{\Gamma}) \ar[from = U, to = M, Rightarrow, shorten <= 4pt, shorten >= 4pt] \ar[from = 1-2, to = L, Rightarrow, shorten <= 4pt, shorten >= 4pt]{}{}
\end{tikzcd}
\]
and hence deduce that
\begin{align*}
\phi_{\opiGamma, \odGamma{1}{1}}^{\AGamma} \circ F(\odGamma{1}{1})\phi_{\muGamma, \pGamma}^{\AGamma} &= \phi_{\muGamma \circ \odGamma{1}{1}, \pGamma}^{\AGamma} \circ \phi_{\odGamma{1}{1}, \muGamma}^{F(\pGamma)\AGamma} = \phi_{\muGamma \circ \odGamma{0}{1}, \pGamma}^{\AGamma} \circ \phi_{\odGamma{0}{1}, \muGamma}^{F(\pGamma)\AGamma}.
\end{align*}
Similarly, from the equality
\[
F(\aGamma \circ \muGamma \circ \odGamma{2}{1}) = F(\aGamma \circ \muGamma \circ \odGamma{1}{1})
\]
we get the $2$-commutative diagram
\[
\begin{tikzcd}
&F(\quot{X}{\Gamma_c}) \ar[ddr, near start, ""{name = UpDiag}]{}{F(\muGamma \circ \odGamma{2}{1})} \ar[rr, ""{name = UpRight}]{}{F(\muGamma)} & & F(\quot{(G \times X)}{\Gamma}) \ar[ddl, ""{name = UpUp}]{}{F(\odGamma{2}{1})} \\
\\
F(\XGamma) \ar[ddr, swap]{}{F(\aGamma)} \ar[uur]{}{F(\aGamma)} \ar[rr, bend left = 20, ""{name = MidUp}]{}{F(\aGamma \circ \muGamma \circ \odGamma{2}{1})} \ar[rr, bend right = 20, swap, ""{name = MidLow}]{}{F(\aGamma \circ \muGamma \circ \odGamma{1}{1})} & & F(\quot{(G \times G \times X)}{\Gamma}) \\
\\
& F(\quot{X}{\Gamma_c}) \ar[uur, swap, near start, ""{name = DownDiag}]{}{F(\muGamma \circ \odGamma{1}{1})} \ar[rr, swap, ""{name = DownRight}]{}{F(\muGamma)} & & F(\quot{(G \times X)}{\Gamma}) \ar[uul, swap]{}{F(\odGamma{1}{1})}
\ar[from = MidUp, to= MidLow, equals, shorten <= 4pt, shorten >= 4pt]
\ar[from = 1-2, to = MidUp, Rightarrow, shorten <= 4pt, shorten >= 8pt]
\ar[from = UpRight, to = 3-3, Rightarrow, shorten <= 4pt, shorten >= 4pt]
\ar[from = 5-2, to = MidLow, Rightarrow, shorten <= 4pt, shorten >= 8pt]
\ar[from = DownRight, to = 3-3, Rightarrow, shorten <= 4pt, shorten >= 4pt]
\end{tikzcd}
\]
which allows us to deduce that
\begin{align*}
\left(\phi_{\odGamma{2}{1},\muGamma}^{F(\aGamma)\AGamma}\right)^{-1} \circ \left(\phi_{\muGamma \circ \odGamma{2}{1}, \aGamma}^{\AGamma}\right)^{-1} &= \left(\phi_{\odGamma{1}{1}, \muGamma}^{F(\aGamma)\AGamma}\right)^{-1} \circ \left(\phi_{\muGamma \circ \odGamma{1}{1}, \aGamma}^{\AGamma}\right)^{-1} \\
&= \left(F(\odGamma{1}{1})\phi_{\muGamma, \pGamma}^{\AGamma}\right)^{-1} \circ \left(\phi_{\oalphaGamma, \odGamma{1}{1}}^{\AGamma}\right)^{-1}.
\end{align*}

Before we can complete the proof of the theorem, we need three last observations. First note that by using the compositor isomorphisms we have that
\[
F(\muGamma \circ \odGamma{2}{1})\tau_{\aGamma}^{A} = \phi_{\odGamma{2}{1}, \muGamma}^{\quot{A}{\Gamma_c}} \circ \big(F(\odGamma{2}{1}) \circ F(\muGamma)\big)\tau_{\aGamma}^{A} \circ \left(\phi_{\odGamma{2}{1}, \muGamma}^{F(\aGamma)\AGamma}\right)^{-1}
\]
and similarly that
\[
F(\muGamma \circ \odGamma{0}{1})(\tau_{\pGamma}^{A})^{-1} = \phi_{\odGamma{0}{1}, \muGamma}^{F(\pGamma)\AGamma} \circ \big(F(\odGamma{0}{1}) \circ F(\muGamma)\big)(\tau_{\pGamma}^{A})^{-1} \circ \left(\phi_{\odGamma{0}{1}, \muGamma}^{\quot{A}{\Gamma_c}}\right)^{-1}.
\]
However, using that
\begin{align*}
&\big(F(\odGamma{0}{1}) \circ F(\muGamma)\big)(\tau_{\pGamma}^{A})^{-1} \circ \left(\phi_{\odGamma{0}{1}, \muGamma}^{\quot{A}{\Gamma_c}}\right)^{-1} \circ \phi_{\odGamma{2}{1}, \muGamma}^{\quot{A}{\Gamma_c}} \circ \big(F(\odGamma{2}{1}) \circ F(\muGamma)\big)\tau_{\aGamma}^{A} \\
&= \big(F(\odGamma{1}{1}) \circ F(\muGamma)\big)(\tau_{\pGamma}^{A})^{-1} \circ \left(\phi_{\odGamma{1}{1}, \muGamma}^{\quot{A}{\Gamma_c}}\right)^{-1} \circ \phi_{\odGamma{1}{1}, \muGamma}^{\quot{A}{\Gamma_c}} \circ \big(F(\odGamma{1}{1}) \circ F(\muGamma)\big)\tau_{\aGamma}^{A} \\
&= \big(F(\odGamma{1}{1}) \circ F(\muGamma)\big)(\tau_{\pGamma}^{A})^{-1} \circ \big(F(\odGamma{1}{1}) \circ F(\muGamma)\big)\tau_{\aGamma}^{A} \\
&= \big(F(\odGamma{1}{1}) \circ F(\muGamma)\big)\left((\tau_{\pGamma}^{A})^{-1} \circ \tau_{\aGamma}^{A}\right)
\end{align*}
we derive that
\begin{align*}
&\phi_{\quot{d_0^1}{\Gamma}, \opiGamma}^{\AGamma} \circ F(\quot{d_0^1}{\Gamma})\quot{\theta}{\Gamma} \circ \left(\phi_{\quot{d_0^1}{\Gamma}, \oalphaGamma}^{\AGamma}\right)^{-1} \circ \phi_{\quot{d_2^1}{\Gamma},\opiGamma}^{\AGamma} \circ F(\quot{d_2^1}{\Gamma})\quot{\theta}{\Gamma} \circ \left(\phi_{\quot{d_2^1}{\Gamma},\oalphaGamma}^{\AGamma}\right)^{-1} \\
&=\phi_{\muGamma \circ \odGamma{0}{1}, \pGamma}^{\AGamma} \circ F(\muGamma \circ \odGamma{0}{1})(\tau_{\pGamma}^{A})^{-1} \circ F(\muGamma \circ \odGamma{2}{1})\tau_{\aGamma}^{A} \circ \left(\phi_{\muGamma \circ \odGamma{2}{1}, \aGamma}^{\AGamma}\right)^{-1} \\
&= \phi_{\muGamma \circ \odGamma{0}{1}, \pGamma}^{\AGamma} \circ \phi_{\odGamma{0}{1}, \muGamma}^{F(\pGamma)\AGamma} \circ \big(F(\odGamma{1}{1}) \circ F(\muGamma)\big)\left((\tau_{\pGamma}^{A})^{-1} \circ \tau_{\aGamma}^{A}\right) \\
&\circ \left(\phi_{\odGamma{2}{1}, \muGamma}^{F(\aGamma)\AGamma}\right)^{-1} \circ \left(\phi_{\muGamma \circ \odGamma{2}{1}, \aGamma}^{\AGamma}\right)^{-1} \\
&= \phi_{\opiGamma, \odGamma{1}{1}}^{\AGamma} \circ F(\odGamma{1}{1})\phi_{\muGamma, \pGamma}^{\AGamma} \circ \big(F(\odGamma{1}{1}) \circ F(\muGamma)\big)\left((\tau_{\pGamma}^{A})^{-1} \circ \tau_{\aGamma}^{A}\right) \\ 
&\circ \left(F(\odGamma{1}{1})\phi_{\muGamma, \pGamma}^{\AGamma}\right)^{-1} \circ \left(\phi_{\oalphaGamma, \odGamma{1}{1}}^{\AGamma}\right)^{-1} \\
&= \phi_{\opiGamma, \odGamma{1}{1}}^{\AGamma} \circ F(\odGamma{1}{1})\quot{\theta}{\Gamma} \circ \left(\phi_{\oalphaGamma, \odGamma{1}{1}}^{\AGamma}\right)^{-1}.
\end{align*}
This establishes the GIT cocycle condition.
\end{proof}
\begin{corollary}
For any left $G$-variety $X$ and any pre-equivariant pseudofunctor $F$ on $X$, there is a natural isomorphism of functors $\pi_2^{\ast} \cong \alpha_X^{\ast}:F_G(X) \to F_G(G \times X)$.
\end{corollary}

\section{Equivariant Functors: Change of Groups}\label{Section: Section 3: Change of Groups}
We begin our discussion on the Change of Groups functors (cf.\@ Theorem \ref{Thm: Section 3: Change of groups functor}) by providing some key results that will allow us to descend through any morphism $\varphi:G \to H$ by using the induction spaces of Bernstein (cf.\@ Definition \ref{Defn: Induction Space} below) to translate from $\Sf(H)$ to $\Sf(G)$. The first such result (Proposition \ref{Prop: Section  3: Bernstein Quotient is Sf(H)}) shows how to realize the induction space as an $\Sf(H)$-object; the second is a corollary (cf.\@ Corollary \ref{Cor: Section 3: Functor from Sf(G) to Sf(H) by taking Bernstein quotients}) which shows a functorial assignment $\Sf(G) \to \Sf(H)$ by moving through the induction spaces. This allows us to describe the Change of Groups functors and then study how they interact with Change of Fibre, Change of Spaces, and even other Change of Groups functors.

In what follows fix smooth algebraic groups $G$ and $H$ over $\Spec K$ and let $\varphi \in \AlgGrp(G,H)$ be a morphism of algebraic groups, where $\AlgGrp$\index[notation]{AlgGrp@$\AlgGrp$} is the category of algebraic groups over $\Spec K$. Given the algebraic group $G$, we also write $\mu_G:G \times G \to G$\index[notation]{Mutipliciation@$\mu_G$} for the multiplication map, $(-1)_G:G \to G$\index[notation]{InverseMapGroup@$(-1)_G$} for the inversion map, and $1_G:\Spec K \to G$\index[notation]{Unit@$1_G$} for the unit morphism of $G$ (and write similar notation for $H$ by replacing every $G$ with an $H$).

\begin{lemma}\label{Lemma: Section 3: Obvious Pullback functor for restriction of action}
There is a pullback functor $\varphi^{\ast}:\HVar \to \GVar$ given by sending an $H$-variety $(X,\alpha_X:H \times X \to X)$ to the $G$-variety $(X,\alpha_X\circ (\varphi \times \id_X):G \times X \to X)$ and sending an $H$-equivariant morphism $f:X \to Y$ to itself.
\end{lemma}
\begin{proof}
The verification that the morphism
\[
\xymatrix{
G \times X \ar[r]^-{\varphi \times \id_X} & H \times X \ar[r]^-{\alpha_X} & X
}
\]
makes $X$ into a $G$-variety is a trivial application of the fact that $X$ is an $H$-variety and the fact that $\varphi$ is a morphism of algebraic groups. Similarly, that $f$ remains as a $G$-equivariant morphism follows because the diagram
\[
\xymatrix{
G \times X \ar[rrr]^-{\alpha_X \circ (\varphi \times \id_X)} \ar[d]_{\id_G \times f} & & & X \ar[d]^{f} \\
G \times Y \ar[rrr]_-{\alpha_Y \circ (\varphi \times \id_Y)} & & & Y
}
\]
factors as
\[
\xymatrix{
G \times X \ar[r]^-{\varphi \times \id_X} \ar[d]_{\id_G \times f} & H \times X \ar[r]^-{\alpha_X} \ar[d]^{\id_H \times f} & X \ar[d]^{f} \\
G \times Y \ar[r]_-{\varphi \times \id_Y} & H \times Y \ar[r]_-{\alpha_Y} & Y
}
\]
and both squares commute.
\end{proof}
\begin{remark}
We will write the action of $G$ on $\varphi^{\ast}X$ as $\varphi^{\ast}\alpha_X$ or $\alpha_{\varphi^{\ast}X}$ (depending on the situation).
\end{remark}
\begin{lemma}\label{Lemma: Section 3: PreEq Pseudo wrt Different GRoups}
Let $\varphi \in \AlgGrp(G,H)$ be a morphism between smooth algebraic groups $G$ and $H$ and let $X$ be an $H$-variety. If $F:\SfResl_H(X)^{\op} \to \fCat$ is a pre-equivariant pseudofunctor then there exists a pre-equivariant pseudofunctor $\quot{F}{G}:\SfResl_G(\varphi^{\ast}X) \to \fCat$ induced by $F$. Moreover, if there is a pseudonatural transformation between pseudofunctors
\[
\begin{tikzcd}
\Var_{/K}^{\op} \ar[r, bend left = 40, ""{name = U}]{}{F} \ar[r, bend right = 40, swap, ""{name = L}]{}{E} & \fCat \ar[from = U, to = L, Rightarrow, shorten <= 4pt, shorten >= 4pt]{}{\alpha}
\end{tikzcd}
\]
there is an induced pseudonatural transformation $\quot{\alpha}{G}:\quot{F}{G} \to \quot{E}{G}$.
\end{lemma}
\begin{proof}
Because $F$ is pre-equivariant, we have a commuting diagram
\[
\xymatrix{
\SfResl_H(X)^{\op} \ar[r]^-{\quo_H^{\op}} \ar[dr]_{F} & \Var_{/\Spec K}^{\op} \ar[d]^{\tilde{F}} \\
 & \fCat
}
\]
where $F:\Var_{/\Spec K}^{\op} \to \fCat$ is a pseudofunctor. We use this to define the pre-equivariant pseudofunctor $\quot{F}{G}$ via the diagram
\[
\xymatrix{
\SfResl_G(\varphi^{\ast}X)^{\op} \ar[r]^-{\quo_G^{\op}} \ar[dr]_{\quot{F}{G}} & \Var_{/\Spec K}^{\op} \ar[d]^{\tilde{F}} \\
 & \fCat
}
\]
that is, $\quot{F}{G} := \tilde{F} \circ \quo_{G}^{\op}$. This diagram evidently satisfies the definition of a pre-equivariant pseudofunctor $\quot{F}{G}:\SfResl_G(\varphi^{\ast}X) \to \fCat$ (cf.\@ Definition \ref{Defn: Prequivariant pseudofunctors}) and hence defines the desired pseudofunctor. The final claim of the lemma is induced by the horizontal composition
\[
\begin{tikzcd}
\SfResl_G(\varphi^{\ast}X)^{\op} \ar[r]{}{\quo^{\op}} & \Var_{/\Spec K}^{\op} \ar[r, bend left = 40, ""{name = U}]{}{F} \ar[r, bend right = 40, swap, ""{name = L}]{}{E} & \fCat \ar[from = U, to = L, Rightarrow, shorten <= 4pt, shorten >= 4pt]{}{\alpha}
\end{tikzcd}
\]
in $\fCAT$.
\end{proof}
\begin{remark}
When context matters and we have pre-equivariant pseudofunctors $F:\SfResl_H(X)^{\op} \to \fCat$ and $F:\SfResl_G(\varphi^{\ast}X)^{\op} \to \fCat$ induced by the same pseudofunctor $F:\Var_{/\Spec K}^{\op} \to \fCat$ in the sense that the diagram
\[
\xymatrix{
\SfResl_G(\varphi^{\ast}X)^{\op} \ar[drr]_{F} \ar[rr]^-{\quo_G^{\op}} & & \Var_{/\Spec K}^{\op} \ar[d]_{F} & & \SfResl_H(X)^{\op} \ar[dll]^{F} \ar[ll]_-{\quo_H^{\op}} \\
 & & \fCat
}
\]
commutes, we will write $\quot{F}{G}$ and $\quot{F}{H}$ for the corresponding pseudofunctors on $(X,\alpha_{\varphi^{\ast}X})$ and $(X,\alpha_X)$, respectively. 
\end{remark}
\begin{definition}\label{Remark: GH pre-equiv pseudofunctors}
Let $\varphi \in \AlgGrp(G,H)$ be a a morphism of algebraic groups and let $X$ be an $H$-variety. We will say that a pseudofunctor $F:\Sch_{/K}^{\op} \to \fCat$ is $(G,H)$-pre-equivariant (with components $\quot{F}{G}$ and $\quot{F}{H}$) when there is a strictly commuting diagram in $\fCAT$ of the form:
\[
\xymatrix{
\SfResl_H(X)^{\op} \ar[r]^-{\quo_H^{\op}} \ar[dr]_{\quot{F}{H}} & \Var_{/\Spec K}^{\op} \ar[d]^{F} & \SfResl_G(\varphi^{\ast}X) \ar[l]_-{\quo_G^{\op}} \ar[dl]^{\quot{F}{G}} \\
 & \fCat
}
\]
By abuse of notation we will call $F$ a $(G,H)$-pre-equivariant\index{Pre-equivariant Pseudofunctor! GH@ $(G,H)$-pre-equivariant Psuedofunctor} pseudofunctor on $X$ and leave all the above data implicit.
\end{definition}

We now give our pseudofunctorial generalization of the change of groups functor presented in \cite{LusztigCuspidal2}. In \cite{LusztigCuspidal2} this functor is given at heart by restricting a $\Vscr$-set $A = \lbrace \AGamma \; | \; \Gamma \in \Sf(G)_0 \rbrace$ to a specific sub-$\Vscr$-set which records the information coming from a morphism $\varphi:G \to H$ of algebraic groups; our construction will give the functor in essentially the same way, save that we distinguish when varieties are isomorphic and bang-on equal to each other.

For what follows we will need the induction space of Bernstein as described in \cite{PramodBook}, \cite{Bien}, and \cite{MirkovicVilonen}. To describe this, we first let $X$ be a $G$-variety and consider the variety $H \times X$. We now define a $G$-action $\alpha_{H \times X}$ on $H \times X$ via the diagram
\[
\begin{tikzcd}
G \times H \times X \ar[dd, swap]{}{\alpha_{H \times X}} \ar[rr]{}{\Delta_G \times \id_{H\times X}} & & G \times G \times H \times X \ar[r]{}{\cong} & H \times G \times G \times X \ar[d]{}{\id_H \times (-1)_G \times \id_{G \times X}}\\
 & & & H \times G \times G \times X \ar[d]{}{\id_H\times \varphi \times \id_{G \times X}} \\
H \times X & & & H \times H \times G \times X \ar[lll]{}{\mu_H\times \alpha_X}
\end{tikzcd}
\]
Note that in set-theoretic terms, this is given by 
\[
g\cdot(h,x) := (h\varphi(g)^{-1},gx).
\] 
Denote the quotient variety $G \backslash (H \times X)$ by $\quot{(H \times X)}{G}$; note that it exists because the action of $G$ on $H$ is image isomorphic to the action an algebraic subgroup of $H$ on the whole group, and because $H$ is a smooth free $H$-variety it follows that $H \times X$ admits a quotient by actions of algebraic subgroups as well.

We now define a further action of $H$ on $\quot{(H \times X)}{G}$. Write $[h,x]_G$ for a generic ``element'' of the scheme $\quot{(H \times X)}{G}$ and define the $H$-action on $\quot{(H \times X)}{G}$ by
\[
h\cdot[h^{\prime},x]_G := [hh^{\prime},x]_G.
\]
Note that this action is well-defined by a scalar restriction-style argument and the fact that any ``element'' $h^{\prime}$ coming from the algebraic subgroup of $H$ isomorphic to $G/\Ker \varphi$ gets absorbed/moved by the quotient action onto the $X$ component.
\begin{definition}\label{Defn: Induction Space}\index{Induction Space}
The scheme $\quot{(H \times X)}{G}$ together with the $H$-action
\[
h\cdot[h^{\prime},x]_G := [hh^{\prime},x]_G
\]
will be denoted $(H \times X)_G$ in this paper. Following \cite{PramodBook}, we call  this space the induction space of $X$. We will denote this space by $(H \times X)_{G}$.
\end{definition}
\begin{remark}
It is common in representation-theoretic literature (cf.\@ \cite{BernLun}, \cite{Bien}, \cite{LusztigCuspidal2}, \cite{MirkovicVilonen}, for instance) to see the above variety written as $H \times_G X$. We will not use this notation, as it conflicts with the pullback notation we have used throughout the paper up to this point. There should also be \textit{less} notational overlap and clash, if possible, so we propose this (further conflicting notation) instead. It is also worth remarking that in \cite{PramodBook} the space ${(H \times X)}_{G}$ is written as $H \times^{G} X$, although we do not adopt this convention here.
\end{remark}

We proceed with a key structural lemma that allows us to see that for any $\Gamma \in \Sf(G)_0$, the scheme $(H \times \Gamma)_G$ is an object in $\Sf(H)$. This works essentially because the actions that arise in this construction and quotients are all smooth, while pure dimension is automatic by construction and the fact that before taking quotients the action has fibres of constant pure dimension as well. To see that this is a principal $H$-variety on the surface seems less trivial than the pure dimension or smoothness of $(H \times \Gamma)_G$. However, because both $\Gamma$ and $H$ are {\'e}tale locally trivializable, it follows that $(H \times \Gamma)_G$ is as well. Thus the actual only nontrivial point is to check that the fibres take the correct form, i.e., for such a trivializing {\'e}tale cover $\lbrace U_i \to H \backslash (H \times \Gamma)_G \; | \; i \in I \rbrace$, we have $\quo^{-1}(U_i) \cong H \times U_i$.

\begin{proposition}\label{Prop: Section  3: Bernstein Quotient is Sf(H)}
Let $G$ and $H$ be smooth algebraic groups over $\Spec K$ with $\varphi \in \AlgGrp(G,H)$ and $\Gamma \in \Sf(G)$. Then $(H \times \Gamma)_G \in \Sf(H)_0$.
\end{proposition}
\begin{proof}
Using the comments prior to the statement of the proof, let $\lbrace \varphi_i:U_i \to G \backslash \Gamma \; | \; i \in I \rbrace$ be an {\'e}tale local trivialization of $\quo_{\Gamma}:\Gamma \to G \backslash \Gamma$, i.e., an {\'e}tale cover for which in the pullback
\[
\xymatrix{
\quo_{\Gamma}^{-1}(U_i) \ar[r] \ar[d] \pullbackcorner & \Gamma \ar[d]^{\quo_{\Gamma}} \\
U_i \ar[r]_-{\varphi_i} & G \backslash \Gamma
}
\]
we have $\quo_{\Gamma}^{-1}(U_i) \cong G \times U_i$ for all $i \in I$. Equip each scheme $U_i$ with the trivial $G$-action and note that the collection of maps
\[
\lbrace (\id_H \times \varphi_i)_G:H \backslash (H \times U_i)_G \to H \backslash (H \times \Gamma)_G \; | \; i \in I \rbrace,
\]
where each map $(\id_H \times \varphi_i)$ is determined via the assignments
\[
[h,u]_G \mapsto [h, \varphi_i(u)]_G,
\]
forms an {\'e}tale cover of $H \backslash(H \times \Gamma)_G$. It remains to check that the fibres have the correct format, i.e., that $\quo_{(H \times \Gamma)_G}^{-1}(H \backslash (H \times U_i)_G) \cong H \times U_i$. To check this we calculate that
\begin{align*}
\quo^{-1}\left(H \backslash (H \times U_i)_G\right) &\cong H \times G \backslash \quo^{-1}_{\Gamma}(U_i) \cong H \times G \backslash (G \times U_i) \cong H \times U_i;
\end{align*}
the first isomorphism holds because the $H$-action on $(H \times U_i)_{G}$ is of the form $h[h^{\prime},u]_G = [hh^{\prime},u]_G$ (and hence trivial on $U_i$) and because the elements $[h,u]_G$ are defined in terms of a $G$-quotient and can be absorbed to the second component of the product, the second isomorphism holds because $\quo_{\Gamma}^{-1}(U_i) \cong G \times U_i$, and the last isomorphism holds because the action of $G$ on $U_i$ is trivial.
\end{proof}

We now conclude that this produces a functor over $\Sf(G)$, which is a necessary step in concluding that we get a functor $F_H(X) \to F_G(X)$. Essentially this just says that the constructions above are functorial in $\Gamma$, which is a routine proof we present below. The crucial ingredient is that if we have a map $f \in \GVar(X,Y)$ then the induced morphism
\[
\id_H \times f:H \times X \to H \times Y
\]
is also $G$-equivariant because the actions $\alpha_{H \times X}$ and $\alpha_{H \times Y}$ do not intermingle the spaces $H$ and $X$ (and also $H$ and $Y$) in nontrivial ways; thus the equivariance of $f$ gets preserved upon taking products with $H$ and $\id_H$. In particular, this says that upon taking $G$-quotients, we get the commuting diagram
\[
\xymatrix{
H \times X \ar[rr]^-{\id_H\times f} \ar[d]_{\quo_{H \times X}} & & H \times Y \ar[d]^{\quo_{H\times Y}} \\
G\backslash{(H \times X)} \ar[rr]_-{(\id_H\times f)_G} & & G\backslash (H \times Y)
}
\]
of schemes. We claim the (suggestively notated) morphism $(\id_H\times f)_G$ is also $H$-equivariant, i.e., is a morphism $(H \times X)_G \to (H \times Y)_G$ of $H$-varieties. However, this is a routine argument; in set-theoretic terms, this morphism is given via
\[
[h,x]_G \mapsto [h,f(x)]_G
\]
and since the action is given by $h[h^{\prime},x]_G = [hh^{\prime},x]_G$, we get that
\[
(\id_H \times f)_G(h[h^{\prime},x]_G) = [hh^{\prime},f(x)]_G = h[h^{\prime},f(x)]_G = h\left((\id_H \times f)_G[h^{\prime},x]_G\right).
\]
It is also routine to verify using the argument above that this construction is functorial; identities get preserved immediately, while the fact that composition is preserved follows from the uniqueness of the induced maps between quotient schemes. Putting these observations together with Proposition \ref{Prop: Section  3: Bernstein Quotient is Sf(H)} then gives us the corollary below, which is crucial  our change of groups functors.

\begin{corollary}\label{Cor: Section 3: Functor from Sf(G) to Sf(H) by taking Bernstein quotients}
There is a functor $\Sf(G) \to \Sf(H)$ which sends a variety $\Gamma$ to the $H$-variety $(H \times \Gamma)_G$ and which sends a morphism $f:\Gamma \to \Gamma^{\prime}$ to $(\id_{H} \times f)_{G}:(H \times \Gamma)_G \to (H \times \Gamma^{\prime})_{G}$.
\end{corollary}
\begin{proof}
Given the remarks prior to the statement of this corollary and Proposition \ref{Prop: Section  3: Bernstein Quotient is Sf(H)}, all that we need do is check that the fibres of the smooth morphism $(\id_H \times f)_G$ are of constant pure dimension for all $f \in \Sf(G)_1$. However, this is immediate from the fact that $\id_H \times f$ has fibres of constant pure dimension ($\id_H$ trivially so and $f$ by assumption) and from the fact that the $G$-action and its quotient do not change fibre dimensions in a non-uniform fashion.
\end{proof}
We need one final lemma before we can present the change of groups functors. This lemma realizes an isomorphism of varieties
\[
\quot{X}{\Gamma} \cong \quot{X}{{(H \times\Gamma)}_{G}}
\]
which allows us to define an object $A \in F_G(X)$ by restricting an object of $F_H(X)$ to those pieces whose fibres are equivalent to fibres which come from the varieties $\XGamma$ for $\Gamma \in \Sf(G)_0$. It is worth remarking that this is an unproven side remark in \cite{LusztigCuspidal2}; however, it is absolutely fundamental to getting the change of groups functors running, so we sketch a proof of this fact here.
\begin{lemma}\label{Lemma: Section 3: Induction space quotient is iso to the Sf(G) quotient}
For any $\Gamma \in \Sf(G)_0$, there is an isomorphism of $K$-varieties
\[
G\backslash(\Gamma \times X)= \XGamma \cong \quot{X}{{(H \times \Gamma)}_{G}} = H \backslash \big({(H \times \Gamma)}_{G} \times X\big).
\]
\end{lemma}
\begin{proof}[Sketch]
We prove this by way of a set-theoretical group action sketch.

Begin by fixing a variety $\Gamma \in \Sf(G)_0$ and recall from Proposition \ref{Prop: Section  3: Bernstein Quotient is Sf(H)} that ${(H \times \Gamma)}_{G}$ is an object in $\Sf(H)$. Thus the quotient $H\backslash\big({(H \times \Gamma)}_{G} \times X\big) =: \quot{X}{{(H \times \Gamma)}_{G}}$ exists for all $H$-varieties $X$. Fix such an $X$. By taking the $H$ quotient of ${(H \times \Gamma)}_{G} \times X$ we get an $H$-variety which, upon pulling back to $\GVar$ as in Lemma \ref{Lemma: Section 3: Obvious Pullback functor for restriction of action} realizes $\quot{X}{{(H \times \Gamma)}_{G}}$ as a $G$-variety with a trivial $G$-action. In particular, it is routine to check that there is a trivial $G$-morphism $\Gamma \times X \to \quot{X}{{(H \times \Gamma)}_{G}}$ induced by the set-theoretic argument
\[
(\gamma,x) \mapsto \big[[1_H,\gamma]_G, x\big]_H;
\]
note that because the $H$-action on ${(H \times \Gamma)}_{G}$ is given by $h[h^{\prime},\gamma]_G = [hh^{\prime},\gamma]_G$ and $[h,g\gamma] = [h\varphi(g)^{-1},\gamma]_G$, we have that this well-defines a morphism. In particular, because $\XGamma$ is a categorical quotient, this implies that there is a unique morphism $\zeta:\XGamma \to \quot{X}{\quot{(H \times \Gamma)}{G}}$ making
\[
\xymatrix{
\Gamma \times X \ar[rr] \ar[dr]_{q} & & \quot{X}{{(H \times \Gamma)}_{G}} \\
 & \XGamma \ar@{-->}[ur]_{\exists!\zeta}
}
\]
commute.

We now construct a $\theta:\quot{X}{{(H \times \Gamma)}_{G}} \to \XGamma$ factoring the diagram:
\[
\xymatrix{
\Gamma \times X \ar[rr]^-{q} \ar[dr] & & \XGamma \\
 & \quot{X}{{(H \times \Gamma)}_{G}} \ar@{-->}[ur]_{\exists\,\theta}
}
\]
To construct this morphism, observe that the assignment
\[
\big[[h,\gamma]_G,x\big]_H \mapsto [\gamma,x]_G
\]
well-defines a morphism between $G$-varieties. Moreover, it is straightforward to check that the longer path of the diagram sends $(\gamma,x)$ to
\[
\theta\big[[1_H,\gamma]_G,x\big]_H = [\gamma,x]_G = q(\gamma,x)
\]
so the diagram does indeed commute. Moreover, it is also routine to check from this construction that from the universal property of $\XGamma$, we have that
\[
\theta \circ \zeta = \id_{\XGamma}.
\]
Thus to prove that $\XGamma \cong \quot{X}{\quot{(H \times \Gamma)}{G}},$ we only need to verify the other direction of the isomorphism.

To verify that $\zeta \circ \theta = \id_{\quot{X}{\quot{(H \times \Gamma)}{G}}}$, we consider that the functor $\varphi^{\ast}$ of Lemma \ref{Lemma: Section 3: Obvious Pullback functor for restriction of action} induces the commuting diagram
\[
\xymatrix{
\Gamma \times X \ar[dr] \ar[d] \\
{{(H \times \Gamma)}_{G}} \times X \ar[r] & \quot{X}{{(H \times \Gamma)}_{G}}
}
\]
as the map $(\gamma, x) \mapsto ([1_H,\gamma]_G,x)$ is $G$-equivariant from the action
\[
g ([1_H,\gamma]_G,x) = (g[1_H,\gamma]_G,gx) = ([1_H,\gamma]_G,\varphi(g)x).
\]
This allows us to note that upon taking $G$-orbits of the above action that $\zeta[\gamma,x] = \big[[h,\gamma]_G,x\big]_H$ is well-defined (as permuting the $g$-action and using how the quotient $[h,\gamma]_G$ is defined allows us manipulate the shift by an appropriate $G$-multiplication $x$ and then cancel via a quotient). Thus we have that
\[
(\zeta \circ \theta)\big[[h,\gamma]_G,x\big] = \zeta[\gamma,x]_G = \big[[h,\gamma]_G,x\big]
\]
so we get that $\zeta \circ \theta = \id_{\quot{X}{\quot{(H \times \Gamma)}{G}}}$. Thus $\XGamma \cong \quot{X}{{(H \times \Gamma)}_{G}}$, as desired.
\end{proof}

We finally have the tools at hand to give the change of fibre functors we will use in this paper. Essentially, given an $F_H(X)$-object $A$, we can restrict the $\Vscr$-set $\lbrace \AGamma \; | \; \Gamma \in \Sf(H)_0 \rbrace$ to the $\Vscr$-subset 
\[
\lbrace \quot{A}{(H \times \Gamma^{\prime})_G} \; | \; \Gamma^{\prime} \in \Sf(G)_0 \rbrace
\] 
varying through $\Sf(G)$ taking values in the quotient varieties $(H \times \Gamma^{\prime})_G$. However, to make this into an honest functor to $F_G(X)$, we need each of the objects to take values in $F(\XGamma)$ and not the induction spaces $\quot{X}{(H \times \Gamma)_G}$. For this, however, we use Lemma \ref{Lemma: Section 3: Induction space quotient is iso to the Sf(G) quotient} and apply the equivalence of categories
\[
F(h_{\Gamma}):F(\quot{X}{(H \times \Gamma)_G}) \xrightarrow{\simeq} F(\XGamma)
\] 
induced by the isomorphism $h_{\Gamma}:\XGamma \xrightarrow{\cong} \quot{X}{(H \times \Gamma)_G}$. That this does indeed define a functor is essentially due to Corollary \ref{Cor: Section 3: Functor from Sf(G) to Sf(H) by taking Bernstein quotients} and Lemma \ref{Lemma: Section 3: Induction space quotient is iso to the Sf(G) quotient}, as the two of these together say that we can pass through the equivalences $F(\XGamma) \simeq F(\quot{X}{{(H \times \Gamma)}_{G}})$ induced by the isomorphisms $\XGamma \cong \quot{X}{{(H \times \Gamma)}_{G}}$ to get a functor by following the proof of Theorem \ref{Thm: Section 3: Existence of pullback functors}.
\begin{Theorem}\label{Thm: Section 3: Change of groups functor}
Let $\varphi:G \to H$ be a morphism of algebraic groups and let $X$ be an $H$-variety. There is then a functor $\varphi^{\sharp}:F_H(X) \to F_G(X)$\index[notation]{PhiSharp@$\varphi^{\sharp}$} given by sending an object $(A,T_A)$ to the pair $(A^{\prime},T_{A^{\prime}})$ and a morphism $P$ to $P^{\prime}$, where the pair $(A^{\prime}, T_{A^{\prime}})$ is induced by the equation
\[
A^{\prime} := \lbrace F(h_{\Gamma})\quot{A}{(H \times \Gamma^{\prime})_G} \; | \; \Gamma^{\prime} \in \Sf(G)_0, \quot{A}{(H \times \Gamma^{\prime})_G} \in A \rbrace
\]
and on $P^{\prime}$ is induced by
\[
P^{\prime} := \lbrace F(h_{\Gamma})\quot{\rho}{(H \times \Gamma^{\prime})_G} \; | \; \Gamma^{\prime} \in \Sf(G)_0, \quot{\rho}{(H \times \Gamma^{\prime})_G} \in P  \rbrace,
\]
where $h_{\Gamma}$ is the isomorphism $\XGamma \xrightarrow{\cong} \quot{X}{{(H \times \Gamma)}_{G}}$ of Lemma \ref{Lemma: Section 3: Induction space quotient is iso to the Sf(G) quotient}.
\end{Theorem}
\begin{proof}
By Corollary \ref{Cor: Section 3: Functor from Sf(G) to Sf(H) by taking Bernstein quotients} and Lemma \ref{Lemma: Section 3: Induction space quotient is iso to the Sf(G) quotient} there are functorial isomorphisms
\[
\XGamma \to \quot{X}{{(H \times \Gamma)}_{G}}
\]
for all $\Gamma \in \Sf(G)_0$. Let $h_{\Gamma}:\XGamma \to \quot{X}{\quot{(H\times\Gamma)}{G}}$ be the induced isomorphism of varieties and consider the induced equivalence of categories
\[
F(h_{\Gamma}):F(\quot{X}{{(H \times \Gamma)}_{G}}) \xrightarrow{\simeq} F(\XGamma).
\]
Now observe that the diagram of schemes
\[
\xymatrix{
\XGamma \ar[dr]^{k_f} \ar[r]^-{h_{\Gamma}} \ar[d]_{\overline{f}} & \quot{X}{{(H \times \Gamma)}_{G}} \ar[d]^{\overline{(\id_H \times f)_G}} \\
\XGammap \ar[r]_-{h_{\Gamma^{\prime}}} & \quot{X}{{(H \times \Gamma^{\prime})}_{G}}
}
\]
commutes for all $f:\Gamma \to \Gamma^{\prime} \in \Sf(G)_1$. Set
\[
k_f := h_{\Gamma^{\prime}} \circ \overline{f}
\]
so that
\[
h_{\Gamma^{\prime}} \circ \overline{f} = k_f = \overline{(\id_H \times f)_G} \circ h_{\Gamma}.
\]
This in turn implies that if $(A,T_A) \in F_H(X)_0$ we must have that
\[
\tau_f^{\varphi^{\sharp}} = F(h_{\Gamma})\big(\tau_{\overline{(\id_H\times f)}_G}^{A}\big) \circ \left(\phi_{h_{\Gamma},\overline{(\id_H \times f)}_G}\right)^{-1} \circ \phi_{\overline{f},h_{\Gamma^{\prime}}}
\]
analogously to Theorem \ref{Thm: Section 3: Existence of pullback functors}.
We now check that $\varphi^{\sharp}$ defines a functor. However, this follows mutatis mutandis from the proof of Theorem \ref{Thm: Section 3: Existence of pullback functors}.
\end{proof}

\begin{definition}\index{Equivariant Functor! Change of Group}
	An equivariant functor $R:F_H(X) \to F_G(X)$ is a Change of Groups functor if there is a morphism of algebraic groups $\varphi:G \to H$ for which $R = \varphi^{\sharp}$.
\end{definition}

Let us now move to study how the Change of Groups functors interact when we change fibres and when we change spaces. In particular, we would ultimately like to show in this section how to interpret equivariant functors $F_G(X) \to E_H(Y)$ (when given a morphisms $\varphi:H \to G$ of smooth algebraic groups with $G$-varieties $X$ and $Y$) and in what sense an arbitrary equivariant functor depends on the order of composition, i.e., we want to study how any and all of the composite paths from $F_G(X) \to E_H(Y)$ in the cube
\[
\begin{tikzcd}
 & F_G(Y) & & F_H(Y) \\
F_G(X) & & F_H(X) \\
 & E_G(Y) & &E_H(Y) \\
E_G(X) & & E_H(X)
\ar[from = 1-2, to = 3-2] \ar[from = 3-2, to = 3-4] \ar[from = 2-1, to = 1-2]
\ar[from = 1-2, to = 1-4] \ar[from = 1-4, to = 3-4] \ar[from = 2-3, to = 1-4]
\ar[from = 4-1, to = 3-2] \ar[from = 4-1, to = 4-3] \ar[from = 4-3, to = 3-4]
\ar[from = 2-1, to = 4-1]
\ar[from = 2-1, to = 2-3, crossing over]
\ar[from = 2-3, to = 4-3, crossing over]
\end{tikzcd}
\]
compare when producing a functor $F_G(X) \to E_H(Y)$.

Let us begin by studying how Change of Fibre functors interact with Change of Group functors. Fix a morphism of smooth algebraic groups $\varphi:G \to H$ and let $X$ be an $H$-variety with $(G,H)$-pre-equivariant pseudofunctors $E$ and $F$ on $X$ (cf.\@ Definition \ref{Remark: GH pre-equiv pseudofunctors} for the definition of a $(G,H)$-pre-equivariant pseudofunctor). Now if $\alpha:F \to E$ is a pseudonatural transformation then there are induced functors
\[
\underline{\alpha}_H:F_H(X) \to E_H(X)
\] 
and
\[
\underline{\alpha}_G:F_G(X) \to E_G(X)
\]
from the whiskered transformations $\alpha_H:\quot{F}{H} \to \quot{E}{H}$ and $\alpha_G:\quot{F}{G} \to \quot{E}{G}$ (cf.\@ Lemma \ref{Lemma: Section 3: PreEq Pseudo wrt Different GRoups}). Our goal now is to study the relationship between the composites
\[
\xymatrix{
F_H(X) \ar[r]^{\varphi^{\sharp}_F}& F_G(X) \ar[r]^{\underline{\alpha}_G} & E_G(X)
}
\]
and
\[
\xymatrix{
F_H(X) \ar[r]^{\underline{\alpha}_H} & E_H(X) \ar[r]^{\varphi^{\sharp}_E} & E_G(X)
}
\]
which is done in Proposition \ref{Prop: Section 3: Change of groups interacts with Change of fibres} below.

In order to discuss the comparison of $\varphi_E^{\sharp} \circ \underline{\alpha}_H$ and $\underline{\alpha}_G \circ \varphi_F^{\sharp}$, we first calculate how they behave on objects $A = \lbrace \AGamma \; | \; \Gamma \in \Sf(G)_0\rbrace$. We calculate that, if $h_{\Gamma}:\XGamma \xrightarrow{\cong} \quot{X}{(H \times \Gamma)_G}$ is the isomorphism of varieties of Lemma \ref{Lemma: Section 3: Induction space quotient is iso to the Sf(G) quotient} for any $\Gamma \in \Sf(G)_0$, $\varphi_E^{\sharp} \circ \underline{\alpha}_H$ acts on an object $A$ as follows: For any $\Gamma \in \Sf(G)_0$,
\[
\quot{(\varphi_E^{\sharp} \circ \underline{\alpha}_H)(A)}{\Gamma} = E(h_{\Gamma})\left(\quot{\alpha}{(H \times \Gamma)_G}\big(\quot{A}{(H \times \Gamma)_G}\big)\right)
\]
so
\[
(\varphi_E^{\sharp} \circ \underline{\alpha}_H)(A) = \left\lbrace E(h_{\Gamma})\left(\quot{\alpha}{(H \times \Gamma)_G}\big(\quot{A}{(H \times \Gamma)_G}\big)\right) \; : \; \Gamma \in \Sf(G)_0\right\rbrace.
\]
Alternatively, we find that $\underline{\alpha}_G \circ \varphi_F^{\sharp}$ acts on $A$ as follows: For any $\Gamma \in \Sf(G)_0$,
\[
\quot{(\underline{\alpha}_G \circ \varphi_F^{\sharp})(A)}{\Gamma} = \quot{\alpha}{\Gamma}\left(F(h_{\Gamma})\quot{A}{(H \times \Gamma)_G}\right)
\]
so
\[
(\underline{\alpha}_G \circ \varphi_F^{\sharp})(A) = \left\lbrace \quot{\alpha}{\Gamma}\left(F(h_{\Gamma})\quot{A}{(H \times \Gamma)_G}\right) \; : \; \Gamma \in \Sf(G)_0 \right\rbrace.
\]
We now consider that since $F$ and $E$ are both $(H,G)$-pre-equivaraint and since $\alpha$ is a pseudnatrual transformation between $F$ and $E$, we have the invertible $2$-cell
\[
\begin{tikzcd}
F(\quot{X}{(H \times \Gamma)_G}) \ar[r, ""{name = U}]{}{F(h_{\Gamma})} \ar[d, swap]{}{\quot{\alpha}{(H \times \Gamma)_G}} & F(\XGamma) \ar[d]{}{\quot{\alpha}{\Gamma}} \\
E(\quot{X}{(H \times \Gamma)_{G}}) \ar[r, swap, ""{name = L}]{}{E(h_{\Gamma})} & E(\XGamma) \ar[from = U, to = L, Rightarrow, shorten <= 4pt, shorten >= 4pt]{}{\quot{\alpha}{h_{\Gamma}}}
\end{tikzcd}
\]
where $\quot{\alpha}{h_{\Gamma}}$ is a natural isomorphism. With this and the pseudonatural identities that the $\alpha$ satisfy, we can deduce via the techniques of Theorems \ref{Thm: Section 3: Psuedonatural trans lift to equivariant functors}, \ref{Thm: Section 3: Existence of pullback functors}, and Proposition \ref{Prop: Section 3: Change of  fibre then space is iso to changing space then fibre} that these isomorphisms assemble to form a natural isomorphism between the functors $(\underline{\alpha}_G \circ \varphi_F^{\sharp}) \xRightarrow{\cong} (\varphi^{\sharp}_E \circ \underline{\alpha}_H)$. In particular, the map
\[
\theta_A: (\underline{\alpha}_G \circ \varphi_F^{\sharp})(A) \to (\varphi^{\sharp}_E \circ \underline{\alpha}_H)A
\]
given by
\[
\theta_A:= \lbrace \quot{\alpha}{h_{\Gamma}}_{\quot{A}{(H \times \Gamma)_{G}}}:\quot{(\underline{\alpha}_G \circ \varphi_F^{\sharp})A}{\Gamma} \xrightarrow{\cong} \quot{(\varphi^{\sharp}_E \circ \underline{\alpha}_H)A}{\Gamma} \; | \; \Gamma \in \Sf(G)_0 \rbrace
\]
determines a morphism $\theta \in E_G(X)((\underline{\alpha}_G \circ \varphi_F^{\sharp})(A), (\varphi^{\sharp}_E \circ \underline{\alpha}_H)A)$. In particular, this is readily checked to be natural by using the naturality of the $\quot{\alpha}{h_{\Gamma}}$ together with the induced pseudonaturality between the fibre functors satisfying the conditions of Theorem \ref{Prop: Necessary and Sufficeint Conditions for Change of Fibre commute with Change of Space}. These observations give rise to the proposition below:

\begin{proposition}\label{Prop: Section 3: Change of groups interacts with Change of fibres}
There is a natural isomorphism of functors $\theta:\underline{\alpha}_G \circ \varphi^{\sharp}_F \xRightarrow{\cong} \varphi^{\sharp}_E \circ \underline{\alpha}_H$ which fits into the invertible $2$-cell of equivariant categories:
\[
\begin{tikzcd}
F_H(X) \ar[r, ""{name = U}]{}{\underline{\alpha}_H} \ar[d, swap]{}{\varphi^{\sharp}_F} & E_H(X) \ar[d]{}{\varphi^{\sharp}_E} \\
F_G(X) \ar[r, swap, ""{name = L}]{}{\underline{\alpha}_G} & E_G(X) \ar[from = L, to = U, Rightarrow, shorten >= 4pt, shorten <= 4pt]{}{\theta}
\end{tikzcd}
\]
\end{proposition}

We now move to discuss how Change of Group functors interact with Change of Space functors. For this, however, we will need to be in a remarkably technical situation; however, when considering this, just keep in mind that the technical jargon and adjectives that appear only mean that we are in a situation where we can talk about moving between groups and spaces at the same time. For this, we let $F$ be a simultaneous $(G,H)$-pre-equivariant pseudofunctor over $H$-varieties
$X$ and $Y$ with $g \in \HVar(X,Y)$. Assume furthermore that that the component functors $\quot{F}{G}$ and $\quot{F}{H}$ both admit descent pullbacks along $g$. To set some necessary notation, recall that for each $\Gamma \in \Sf(G)_0$ we have isomorphisms
\[
\quot{h_{\Gamma}}{X}:\XGamma \xrightarrow{\cong} \quot{X}{(H \times \Gamma)_G}, \qquad \quot{h_{\Gamma}}{Y}:\YGamma \xrightarrow{\cong} \quot{X}{(H \times \Gamma)_G}
\]
and morphisms induced by $g$ which fit into the following commuting diagram (where $f:\Gamma \to \Gamma^{\prime}$ is an $\Sf(G)$-morphism):
\[
\begin{tikzcd}
 & & & \quot{Y}{(H \times \Gamma)_G} \ar[rr]{}{\quot{\overline{(\id_H \times f)}_G}{Y}} & & \quot{Y}{(H \times \Gamma^{\prime})}  \\
 \quot{X}{(H \times \Gamma)_G}  \ar[urrr]{}{\quot{\overline{g}}{(H \times \Gamma)_G}} \ar[rr, swap]{}{\quot{\overline{(\id_H \times f)}_G}{X}} & & \quot{X}{(H \times \Gamma^{\prime})_G}   \\
 & & & \YGamma \ar[rr]{}{\fY}  \ar[uu]{}{\quot{h_{\Gamma}}{Y}}& & \YGammap\ar[uu, swap]{}{\quot{h_{\Gamma^{\prime}}}{Y}}\\
 \XGamma \ar[urrr]{}{\quot{\overline{g}}{\Gamma}} \ar[rr, swap]{}{\fX} \ar[uu]{}{\quot{h_{\Gamma}}{X}} & & \XGammap \ar[urrr, swap]{}{\quot{\overline{g}}{\Gamma^{\prime}}}
 \ar[from = 2-3, to = 1-6, crossing over, swap]{}{\quot{\overline{g}}{(H \times \Gamma^{\prime})_G}}
 \ar[from = 4-3, to = 2-3, crossing over, swap]{}{\quot{h_{\Gamma^{\prime}}}{X}}
\end{tikzcd}
\]
Because of this we have functors and composites
\[
\xymatrix{
F_H(Y) \ar[r]^-{g_H^{\ast}}& F_H(X) \ar[r]^-{\varphi_X^{\sharp}} & F_G(X)
}
\]
and
\[
\xymatrix{
F_H(Y) \ar[r]^-{\varphi_{Y}^{\sharp}} & F_G(Y) \ar[r]^-{g_G^{\ast}} & F_G(X)
}
\]
To see how these composites differ, we first calculate for an object $A \in F_H(Y)_0$ we have that on one hand for all $\Gamma \in \Sf(G)_0$, the $\Gamma$-local components of $(\varphi^{\sharp}_X \circ g_H^{\ast})A$ are given by
\[
\quot{(\varphi^{\sharp}_X \circ g_H^{\ast})A}{\Gamma} = F(\quot{h_{\Gamma}}{X})\left((\quot{\overline{g}}{(H \times \Gamma)_G})^{\ast}\quot{A}{(H \times \Gamma)_G}\right).
\]
On the other hand, the $\Gamma$-lcocal components of $(g_G^{\ast} \circ \varphi^{\sharp}_Y)A$ are given by
\[
\quot{(g_G^{\ast} \circ \varphi^{\sharp}_Y)A}{\Gamma} = \quot{\overline{g}}{\Gamma}^{\ast}\left(F(\quot{h_{\Gamma}}{Y})\quot{A}{(H \times \Gamma)_G}\right).
\]
Because each component of the $(G,H)$-pre-equivariant pseudofunctor $F$ admits descent pullbacks along $g$, we have an invertible $2$-cell in $\fCat$
\[
\begin{tikzcd}
F(\quot{Y}{(H \times \Gamma)_G}) \ar[r, ""{name = U}]{}{F(\quot{h_{\Gamma}}{Y})} \ar[d, swap]{}{\quot{\overline{g}}{(H \times \Gamma)_G}^{\ast}} & F(\YGamma) \ar[d]{}{\quot{\overline{g}}{\Gamma}^{\ast}} \\
F(\quot{X}{(H \times \Gamma)_G}) \ar[r, swap, ""{name = L}]{}{F(\quot{h_{\Gamma}}{X})} & F(\XGamma) \ar[from = U, to = L, Rightarrow, shorten >= 4pt, shorten <= 4pt]{}{\quot{g^{\ast}}{h_{\Gamma}}}
\end{tikzcd}
\]
which records the commutativity isomorphisms of $F$ with the pullbacks along $g$. We then obtain our natural isomorphism $\theta:(g_G^{\ast} \circ \varphi^{\sharp}_Y) \Rightarrow (\varphi^{\sharp}_X \circ g_H^{\ast})$ by setting the components, for $A \in F_H(Y)_0$, to be given by
\[
\theta_A := \left\lbrace\left(\quot{g^{\ast}}{h_{\Gamma}}\right)_{\quot{A}{(H \times \Gamma)}}:\quot{(g_G^{\ast} \circ \varphi^{\sharp}_Y)A}{\Gamma} \xrightarrow{\cong} \quot{(\varphi^{\sharp}_X \circ g_H^{\ast})A}{\Gamma} \; : \; \Gamma \in \Sf(G)_0 \right\rbrace.
\]
Proving that this is indeed a morphism in $F_G(X)$ is routine but tedious and follows the same techniques of those used in Proposition \ref{Prop: Section 3: Change of groups interacts with Change of fibres} together with those in Theorems \ref{Thm: Section 3: Existence of pullback functors}, \ref{Thm: Section 3: Existence of pushforwards}, \ref{Thm: Section 3.2: Existence of descent pullbacks}, and \ref{Prop: Necessary and Sufficeint Conditions for Change of Fibre commute with Change of Space}. In fact, by changing the local assumptions of having descent pullbacks to having descent pushforwards, we obtain the next two propositions:
\begin{proposition}\label{Prop: Section 3.2: Changer of groups and change of fibre by descent pullbacks}
Let $F$ be a simultaneous $(G,H)$-pre-equivariant-pseudofunctor over $H$-varieties $X$ and $Y$ and let $g \in \HVar(X,Y)$. Then if both components $\quot{F}{G}$ and $\quot{F}{H}$ admit descent pullbacks along $g$, the diagram below commutes up to invertible $2$-cell in $\fCat$:
\[
\begin{tikzcd}
F_H(Y) \ar[r, ""{name = U}]{}{\varphi_{Y}^{\sharp}} \ar[d, swap]{}{g_H^{\ast}} & F_G(Y) \ar[d]{}{g_G^{\ast}} \\
F_H(X) \ar[r, swap, ""{name = L}]{}{\varphi_X^{\sharp}} & F_G(X) \ar[from = U, to = L, Rightarrow, shorten >= 4pt, shorten <= 4pt]{}{\theta}
\end{tikzcd}
\]
\end{proposition}
\begin{proposition}\label{Prop: Section 3.2: Changer of groups and change of fibre by descent pushforwards}
	Let $F$ be a simultaneous $(G,H)$-pre-equivariant pseudofunctor over $H$-varieties $X$ and $Y$ and let $g \in \HVar(X,Y)$. Then if both components $\quot{F}{G}$ and $\quot{F}{H}$ admit descent pushforwards along $g$,  the diagram below commutes up to invertible $2$-cell in $\fCat$:
	\[
	\begin{tikzcd}
	F_H(X) \ar[r, ""{name = U}]{}{\varphi_{X}^{\sharp}} \ar[d, swap]{}{(g_H)_{\ast}} & F_G(X) \ar[d]{}{(g_G)_{\ast}} \\
	F_H(Y) \ar[r, swap, ""{name = L}]{}{\varphi_Y^{\sharp}} & F_G(Y) \ar[from = U, to = L, Rightarrow, shorten >= 4pt, shorten <= 4pt]{}{\theta}
	\end{tikzcd}
	\]
\end{proposition}
%

We now move to discuss how Change of Group functors interact with other Change of Group functors. In particular, if we have morphisms $\varphi \in \AlgGrp(G_0, G_1)$ and $\psi \in \AlgGrp(G_1, G_2)$ for smooth algebraic groups $G_0,$ $G_1,$ $G_2$, we would like to study how the composites $\varphi^{\sharp} \circ \psi^{\sharp}$ and $(\psi \circ \varphi)^{\sharp}$ compare. For this we will need a pre-equivariant pseudofunctor on a $G_2$-variety $X$ for both pairs $(G_0, G_1)$ and $(G_1, G_2)$; let us denote this by a $(G_0, G_1, G_2)$-pre-equivariant pseudofunctor on $X$. There are then two distinct Change of Groups functors from $F_{G_2}(X) \to F_{G_0}(X)$: First, the total composite
\[
\xymatrix{
F_{G_2}(X) \ar[r]^-{\psi^{\sharp}} & F_{G_1}(X) \ar[r]^-{\varphi^{\sharp}} & F_{G_0}(X)
}
\]
and second the morphism induced by $\psi \circ \varphi$:
\[
\xymatrix{
F_{G_2}(X) \ar[r]^-{(\psi \circ \varphi)^{\sharp}} & F_{G_0}(X)
}
\]
While in general these two functors should not be expected to be the same, in the proposition below we will prove that the are isomorphic to each other by a natural isomorphism induced simultaneously from the transition isomorphisms in $F_{G_2}(X)$ that an object must satisfy and the compositor isomorphisms that the pseudofunctor $F$ has at hand.

\begin{proposition}\label{Prop: Section 3.3: Chagnge of groups interacting with change of grups}
Let $\varphi \in \AlgGrp(G_0, G_1)$ and $\psi \in \AlgGrp(G_1, G_2)$ for smooth algebraic groups $G_0, G_1, G_2$ and let $F$ be a $(G_0, G_1, G_2)$-pre-equivariant pseudofunctor on a $G_2$-variety $X$. Then the diagram below commutes up to natural isomorphism in $\fCat$:
\[
\begin{tikzcd}
F_{G_2}(X) \ar[dr, swap]{}{(\psi \circ \varphi)^{\sharp}} \ar[rr]{}{\psi^{\sharp}} & {} &F_{G_1}(X) \ar[dl, swap]{}{\varphi^{\sharp}} \\
 & F_{G_0}(X) \ar[from = 2-2, to = 1-2, Rightarrow, shorten >= 4pt, shorten <= 4pt]{}{\alpha}
\end{tikzcd}
\]
\end{proposition}
\begin{proof}
Begin by observing that for any $\Gamma \in \Sf(G_0)_0$, it is straightforward but tedious to prove that there is a $G_2$-equivariant isomorphism of varieties
\[
(G_2 \times (G_1 \times \Gamma)_{G_0})_{G_1} \xrightarrow[\rho]{\cong} (G_2 \times \Gamma)_{G_0};
\]
in a set-theoretic sketch, this isomorphism is given by $[g_2, [g_1,\gamma]_{G_0}]_{G_1} \mapsto [g_2,\gamma]_{G_0}$ and that this is $G_2$-equivariant is immediate because the $G_2$-action occurs only in the $G_2$-component. Now, because each of the above varieties is an $\Sf(G_2)$-object by three applications of Proposition \ref{Prop: Section  3: Bernstein Quotient is Sf(H)} (one for $(G_2 \times\Gamma)_{G_0}$ directly, one for realizing $(G_1 \times \Gamma)_{G_0}$ as an $\Sf(G_1)$-variety, and then one for the case of realizing $(G_2 \times (G_1 \times \Gamma)_{G_0})_{G_1}$ as an $\Sf(G_2)$-variety) it follows from the fact that isomorphisms are smooth (and in fact {\'e}tale with constant fibres) that $\rho \in \Sf(G_2)_1$.

We calculate our two composites to begin our comparison. First fix an object $A \in F_{G_2}(X)$. Now observe that on one hand
\begin{align*}
\varphi^{\sharp}(\psi^{\sharp}A) &= \varphi^{\sharp}\left(\lbrace F(h^{G_1G_2}_{\Gamma^{\prime}})\quot{A}{(G_2 \times \Gamma^{\prime})_{G_1}} \; | \; \Gamma^{\prime} \in \Sf(G_1)_0 \rbrace\right) \\
&= \lbrace F(h_{\Gamma}^{G_0G_1})F(h^{G_1G_2}_{(G_1 \times \Gamma)_{G_0}}) \quot{A}{(G_2 \times (G_1 \times \Gamma)_{G_0})_{G_1}} \; | \; \Gamma \in \Sf(G_0)_0 \rbrace
\end{align*}
while on the other hand
\[
(\psi \circ \varphi)^{\sharp}A = \lbrace F(h_{\Gamma}^{G_0G_2})\quot{A}{(G_2 \times \Gamma)_{G_0}} \; | \; \Gamma \in \Sf(G)_0 \rbrace.
\]
To compare these two constructions, we note that it is routine to verify that the diagram
\[
\xymatrix{
\quot{X}{\Gamma} \ar[d]_{h_{\Gamma}^{G_0G_2}} \ar[rr]^-{h_{\Gamma}^{G_0G_1}} & & \quot{X}{(G_1 \times \Gamma)_{G_0}} \ar[d]^{h_{(G_1 \times \Gamma)_{G_0}}^{G_1G_2}} \\
\quot{X}{(G_2 \times \Gamma)_{G_0}} & & \quot{X}{(G_2 \times (G_1 \times \Gamma)_{G_0})_{G_1}} \ar[ll]^-{\overline{\rho}}
}
\]
commutes by using the universal properties of quotients. With this we produce the pasting diagram:
\[
\begin{tikzcd}
& F(\quot{X}{(G_1 \times \Gamma)_{G_0}}) \ar[dr, bend left = 30]{}{F(h_{\Gamma}^{G_0G_1})} & \\
F(\quot{X}{(G_2 \times (G_1 \times \Gamma)_{G_0})_{G_1}}) \ar[rr]{}[description]{F(h_{(G_1 \times \Gamma)_{G_0}} \circ h_{\Gamma})} \ar[ur, bend left = 30]{}{F(h_{(G_1 \times \Gamma)_{G_0}}^{G_1G_2})} & {} & F(\quot{X}{\Gamma}) \\
& F(\quot{X}{(G_2 \times \Gamma)_{G_0}}) \ar[ur, swap, bend right = 30]{}{F(h_{\Gamma}^{G_0G_2})} \ar[ul, bend left = 30]{}{F(\overline{\rho})} & \ar[from = 2-2, to = 3-2, swap, Rightarrow, shorten >= 4pt, shorten <= 4pt]{}{\phi_{h^{G_1G_2} \circ h^{G_0G_1}, \rho}} \ar[from = 1-2, to = 2-2, Rightarrow, shorten <= 4pt, shorten >=4pt, swap]{}{\phi_{h^{G_0G_1}, h^{G_1G_2}}}
\end{tikzcd}
\]
From the fact that $\rho \in \Sf(G_2)_1$, we have a transition isomorphism
\[
\tau_{\rho}^{A}:F(\overline{\rho})\quot{A}{(G_2 \times \Gamma)_{G_0}} \xrightarrow{\cong} \quot{A}{(G_2 \times (G_1 \times \Gamma)_{G_0})_{G_1}}.
\]
Applying the functor $F(h_{\Gamma}^{G_0G_1}) \circ F(h_{(G_1 \times \Gamma)_{G_0}}^{G_1G_2})$ to $\tau_{\rho}^{A}$ above gives an isomorphism:
\[
\xymatrix{
F(h_{\Gamma}^{G_0G_1})F(h_{(G_1 \times \Gamma)_{G_0}}^{G_1G_2})F(\overline{\rho})\quot{A}{(G_2 \times \Gamma)_{G_0}} \ar[d]^{\cong} \\
 F(h_{\Gamma}^{G_0G_1})F(h_{(G_1 \times \Gamma)_{G_0}}^{G_1G_2})\quot{A}{(G_2 \times (G_1 \times \Gamma)_{G_0})_{G_1}}
}
\]
Pre-composing this with the inverses of the two compositors in the pasting diagram then gives an isomorphism
\[
\quot{\alpha_A}{\Gamma} := F(h_{\Gamma}^{G_0G_1})F(h_{(G_1 \times \Gamma)_{G_0}}^{G_1G_2})\tau_{\rho}^{A} \circ \left(\phi_{h^{G_0G_1}, h^{G_1G_2}}^{F(\overline{\rho})A}\right)^{-1} \circ \left(\phi_{h^{G_1G_2} \circ h^{G_0G_1}, \rho}^{A}\right)^{-1}
\]
between the objects:
\[
\xymatrix{
\quot{(\psi \circ \varphi)^{\sharp}A}{\Gamma} \ar@{=}[r] \ar[d] & F(h_{\Gamma}^{G_0G_2})\quot{A}{(G_2 \times \Gamma)_{G_0}} \ar[d]^{\cong} \\ \quot{\varphi^{\sharp}(\psi^{\sharp}A)}{\Gamma} \ar@{=}[r] & F(h_{\Gamma}^{G_0G_1})F(h_{(G_1 \times \Gamma)_{G_0}}^{G_1G_2})\quot{A}{(G_2 \times (G_1 \times \Gamma)_{G_0})_{G_1}}
}
\]
Define the map $\alpha_A: (\psi \circ \varphi)^{\sharp}A \to (\varphi^{\sharp} \circ \psi^{\sharp})A$ by setting the $\Gamma$-local components as the morphisms $\quot{\alpha}{\Gamma}$ given above.

We now show that $\alpha$, as defined, is a morphism in $F_{G_0}(X)$ from $(\psi \circ \varphi)^{\sharp}A \to (\varphi^{\sharp} \circ \psi^{\sharp})A$. For this and for the rest of the proof we set some notation to reduce the already significant notational clutter. Make the following definitions of convenience:
\begin{align*}
\Gamma_1 &:= (G_1 \times \Gamma)_{G_0} & \Gamma_2 &:= (G_2 \times \Gamma_1)_{G_1} & \Gamma_3 &:= (G_2 \times \Gamma)_{G_0} \\
\Gamma^{\prime}_{1} &:= (G_1 \times \Gamma^{\prime})_{G_0} & \Gamma^{\prime}_2 &:= (G_2 \times \Gamma^{\prime}_1)_{G_1} & \Gamma_{3}^{\prime} &:= (G_2 \times \Gamma^{\prime})_{G_0} \\
f_1 &:= (\id_{G_1} \times f)_{G_0} & f_2 & := (\id_{G_2} \times f_1)_{G_1} & f_3 &:= (\id_{G_2} \times f)_{G_0} \\
h_{\Gamma}^{01} &:= h_{\Gamma}^{G_0G_1} & h_{\Gamma_1}^{12} &:= h_{\Gamma_1}^{G_1G_2} & h_{\Gamma}^{02} &:= h_{\Gamma}^{G_0G_2} \\
h_{\Gamma^{\prime}}^{01} &:= h_{\Gamma^{\prime}}^{G_0G_1} & h_{\Gamma_1^{\prime}}^{12} &:= h_{\Gamma_1^{\prime}}^{G_1G_2} & h_{\Gamma^{\prime}}^{02} &:= h_{\Gamma^{\prime}}^{G_0G_2}
\end{align*}
To prove that $\alpha_A$ is a morphism in $F_{G_0}(X)$, we must show that the diagram
\[
\xymatrix{
F(\of)\quot{(\psi \circ \varphi)^{\sharp}A}{\Gamma^{\prime}} \ar[rr]^-{F(\of)\quot{\alpha_A}{\Gamma^{\prime}}} \ar[d]_{\tau_f^{(\psi \circ \varphi)^{\sharp}A}} & & F(\of)\quot{\varphi^{\sharp}(\psi^{\sharp}A)}{\Gamma^{\prime}} \ar[d]^{\tau_f^{\varphi^{\sharp}(\psi^{\sharp}A)}} \\
\quot{(\psi \circ \varphi)^{\sharp}A}{\Gamma} \ar[rr]_-{\quot{\alpha}{\Gamma}} & & \quot{\varphi^{\sharp}(\psi^{\sharp}A)}{\Gamma}
}
\]
commutes for any $f \in \Sf(G)_1$. To establish this we first calculate using Theorem \ref{Thm: Section 3: Change of groups functor} that on one hand
\[
\tau_f^{(\psi \circ \varphi)^{\sharp}A} = F(h_{\Gamma}^{02})\tau_{f_3}^{A} \circ \left(\phi_{h_{\Gamma}^{02}, f_3}^{A}\right)^{-1} \circ \phi_{f,h_{\Gamma^{\prime}}^{02}}^{A}
\]
while on the other hand
\begin{align*}
\tau^{\varphi^{\sharp}(\psi^{\sharp}A)}_f &= F(h_{\Gamma}^{01})\tau_{f_1}^{\psi^{\sharp}A} \circ \left(\phi_{h_{\Gamma}^{01},f_1}^{\psi^{\sharp}A}\right)^{-1} \circ \phi_{f,h_{\Gamma^{\prime}}^{01}}^{\psi^{\sharp}A} \\
&= F(h_{\Gamma}^{01})\left(F(h_{\Gamma_1}^{12})\tau_{f_2}^{A} \circ \left(\phi_{h_{\Gamma_1},f_2}^{A}\right)^{-1} \circ \phi_{f_1, h_{\Gamma^{\prime}_1}}^{A}\right) \circ \left(\phi_{h_{\Gamma}^{01},f_1}^{\psi^{\sharp}A}\right)^{-1} \circ \phi_{f,h_{\Gamma^{\prime}}^{01}}^{\psi^{\sharp}A}.
\end{align*}
In order to prove the commutativity of the above diagram, we will break these morphisms into pieces and manipulate them piecewise by starting with $\tau_f^{\varphi^{\sharp}(\psi^{\sharp}A)} \circ F(\of)\quot{\alpha_A}{\Gamma^{\prime}}$ (or a piece of it, anyway). For this it will be particularly helpful to keep the commuting diagram of varieties
\[
\begin{tikzcd}
\quot{X}{\Gamma} \ar[rrr, bend left = 30]{}{h_{\Gamma}^{02}} \ar[r]{}{h_{\Gamma}^{01}} \ar[d, swap]{}{\of} & \quot{X}{\Gamma_1} \ar[d, swap]{}{\of_1} \ar[r]{}{h_{\Gamma_1}^{12}} & \quot{X}{\Gamma_2} \ar[d]{}{\of_2} \ar[r]{}{\rho} & \quot{X}{\Gamma_3} \ar[d]{}{\of_3} \\
\quot{X}{\Gamma^{\prime}} \ar[rrr, bend right = 30, swap]{}{h_{\Gamma^{\prime}}^{02}} \ar[r, swap]{}{h_{\Gamma^{\prime}}^{01}} & \quot{X}{\Gamma_1^{\prime}} \ar[r, swap]{}{h_{\Gamma_1^{\prime}}} & \quot{X}{\Gamma_2^{\prime}} \ar[r, swap]{}{\rho^{\prime}} & \quot{X}{\Gamma_3^{\prime}}
\end{tikzcd}
\] 
in mind. Begin by observing that
\begin{align*}
&\left(\phi_{h_{\Gamma}^{01},f_1}^{\psi^{\sharp}A}\right)^{-1} \circ \phi_{f,h_{\Gamma^{\prime}}^{01}}^{\psi^{\sharp}A} \circ F(\of)\quot{\alpha_A}{\Gamma^{\prime}} \\
&= \left(\phi_{h_{\Gamma}^{01},f_1}^{F(h_{\Gamma_1})\quot{A}{\Gamma_2}}\right)^{-1} \circ \phi_{f,h_{\Gamma^{\prime}}^{01}}^{F(h_{\Gamma_1})\quot{A}{\Gamma_2}}  \\
&\circ F(\of)\left(F(h_{\Gamma^{\prime}}^{01})F(h_{\Gamma_1^{\prime}}^{12})\tau_{\rho^{\prime}}^{A} \circ \left(\phi_{h_{\Gamma^{\prime}}, h_{\Gamma^{\prime}_1}}^{F(\overline{\rho}^{\prime})A}\right)^{-1} \circ \left(\phi_{h_{\Gamma_1^{\prime}} \circ h_{\Gamma^{\prime}}, \rho^{\prime}}^{A}\right)^{-1}\right) \\
&= \left(\phi_{h_{\Gamma}^{01},f_1}^{F(h_{\Gamma_1})\quot{A}{\Gamma_2}}\right)^{-1} \circ \phi_{f,h_{\Gamma^{\prime}}^{01}}^{F(h_{\Gamma_1})\quot{A}{\Gamma_2}} \circ F(\of)F(h_{\Gamma^{\prime}}^{01})F(h_{\Gamma_1^{\prime}}^{12})\tau_{\rho^{\prime}}^{A} \\
&\circ F(\of)\left(\left(\phi_{h_{\Gamma^{\prime}}, h_{\Gamma^{\prime}_1}}^{F(\overline{\rho}^{\prime})A}\right)^{-1} \circ \left(\phi_{h_{\Gamma_1^{\prime}} \circ h_{\Gamma^{\prime}}, \rho^{\prime}}^{A}\right)^{-1}\right) \\
&= \left(\phi_{h_{\Gamma}^{01},f_1}^{F(h_{\Gamma_1})\quot{A}{\Gamma_2}}\right)^{-1} \circ F(h_{\Gamma^{\prime}} \circ \of)F(h_{\Gamma^{\prime}_1}^{12})\tau_{{\rho}^{\prime}}^{A} \circ \phi_{f, h_{\Gamma^{\prime}}}^{F(h_{\Gamma_1^{\prime}})F(\overline{\rho}^{\prime})A} \\
&\circ F(\of)\left(\left(\phi_{h_{\Gamma^{\prime}}, h_{\Gamma^{\prime}_1}}^{F(\overline{\rho}^{\prime})A}\right)^{-1} \circ \left(\phi_{h_{\Gamma_1^{\prime}} \circ h_{\Gamma^{\prime}}, \rho^{\prime}}^{A}\right)^{-1}\right)  \\
&= F(h_{\Gamma}^{01})F(\of_1)F(h_{\Gamma_{1}^{\prime}}^{12})\tau_{\rho^{\prime}}^{A} \circ \left(\phi_{h_{\Gamma},f_1}^{F(h_{\Gamma_1^{\prime}})F(\overline{\rho}^{\prime})A}\right)^{-1} \circ \phi_{f, h_{\Gamma^{\prime}}}^{F(h_{\Gamma_1^{\prime}})F(\overline{\rho}^{\prime})A} \\
&\circ F(\of)\left(\left(\phi_{h_{\Gamma^{\prime}}, h_{\Gamma^{\prime}_1}}^{F(\overline{\rho}^{\prime})A}\right)^{-1} \circ \left(\phi_{h_{\Gamma_1^{\prime}} \circ h_{\Gamma^{\prime}}, \rho^{\prime}}^{A}\right)^{-1}\right).
\end{align*}
We now calculate that
\begin{align*}
&F(h_{\Gamma}^{01})\left(F(h_{\Gamma_1}^{12})\tau_{f_2}^{A} \circ \left(\phi_{h_{\Gamma_1},f_2}^{A}\right)^{-1} \circ \phi_{f_1, h_{\Gamma_1}^{\prime}}^{A}\right) \circ F(h_{\Gamma}^{01})F(\of_1)F(h_{\Gamma_1^{\prime}}^{12})\tau_{\rho^{\prime}}^{A} \\
&= F(h_{\Gamma}^{01})\left(F(h_{\Gamma_1}^{12})\tau_{f_2}^{A} \circ \left(\phi_{h_{\Gamma_1},f_2}^{A}\right)^{-1} \circ \phi_{f_1, h_{\Gamma_1}^{\prime}}^{A} \circ F(\of_1)F(h_{\Gamma_1^{\prime}}^{12})\tau_{\rho^{\prime}}^{A} \right) \\
&= F(h_{\Gamma}^{01})\left(F(h_{\Gamma_1}^{12})\tau_{f_2}^{A} \circ \left(\phi_{h_{\Gamma_1},f_2}^{A}\right)^{-1} \circ F(h_{\Gamma_1^{\prime}}^{12} \circ \of_1)\tau_{\rho^{\prime}}^{A} \circ \phi_{f_1, h_{\Gamma_1^{\prime}}}^{F(\overline{\rho}^{\prime})A} \right) \\
&= F(h_{\Gamma}^{01})\left(F(h_{\Gamma_1}^{12})\tau_{f_2}^{A} \circ F(h_{\Gamma_1})F(\of_2)\tau_{\rho^{\prime}}^{A} \circ \left(\phi_{h_{\Gamma_1}, f_2}^{F(\overline{\rho}^{\prime})A}\right)^{-1} \circ \phi_{f_1, h_{\Gamma_1^{\prime}}}^{F(\overline{\rho}^{\prime})A}\right) \\
&= F(h_{\Gamma}^{01})F(h_{\Gamma_{1}}^{12})\left(\tau_{f_2}^{A} \circ F(\of_2)\tau_{\rho^{\prime}}^{A}\right) \circ F(h_{\Gamma}^{01})\left(\left(\phi_{h_{\Gamma_1}, f_2}^{F(\overline{\rho}^{\prime})A}\right)^{-1} \circ \phi_{f_1, h_{\Gamma_1^{\prime}}}^{F(\overline{\rho}^{\prime})A}\right) \\
&= F(h_{\Gamma}^{01})F(h_{\Gamma_1}^{12})\left(\tau_{\rho^{\prime} \circ f_2}^{A} \circ \phi_{f_2,\rho^{\prime}}^{A}\right) \circ F(h_{\Gamma}^{01})\left(\left(\phi_{h_{\Gamma_1}, f_2}^{F(\overline{\rho}^{\prime})A}\right)^{-1} \circ \phi_{f_1, h_{\Gamma_1^{\prime}}}^{F(\overline{\rho}^{\prime})A}\right) \\
&= F(h_{\Gamma}^{01})F(h_{\Gamma_1}^{12})\left(\tau_{f_3 \circ \rho}^{A} \circ \phi_{f_2,\rho^{\prime}}^{A}\right) \circ F(h_{\Gamma}^{01})\left(\left(\phi_{h_{\Gamma_1}, f_2}^{F(\overline{\rho}^{\prime})A}\right)^{-1} \circ \phi_{f_1, h_{\Gamma_1^{\prime}}}^{F(\overline{\rho}^{\prime})A}\right) \\
&= F(h_{\Gamma}^{01})F(h_{\Gamma_1}^{12})\left(\tau_{\rho}^{A} \circ F(\rho)\tau_{f_3}^{A} \circ \left(\phi_{\rho,f_3}^{A}\right)^{-1} \circ \phi_{f_2,\rho^{\prime}}^{A}\right) \\
&\circ F(h_{\Gamma}^{01})\left(\left(\phi_{h_{\Gamma_1}, f_2}^{F(\overline{\rho}^{\prime})A}\right)^{-1} \circ \phi_{f_1, h_{\Gamma_1^{\prime}}}^{F(\overline{\rho}^{\prime})A}\right) \\
&= F(h_{\Gamma}^{01})F(h_{\Gamma_1}^{12})\tau_{\rho}^{A} \circ  F(h_{\Gamma}^{01})F(h_{\Gamma_1}^{12})\left(F(\rho)\tau_{f_3}^{A} \circ \left(\phi_{\rho,f_3}^{A}\right)^{-1} \circ \phi_{f_2,\rho^{\prime}}^{A}\right) \\
&\circ F(h_{\Gamma}^{01})\left(\left(\phi_{h_{\Gamma_1}, f_2}^{F(\overline{\rho}^{\prime})A}\right)^{-1} \circ \phi_{f_1, h_{\Gamma_1^{\prime}}}^{F(\overline{\rho}^{\prime})A}\right).
\end{align*}
Calculating that
\begin{align*}
&\quot{\alpha_A}{\Gamma} \circ \tau_f^{(\psi \circ \varphi)^{\sharp}A} \\
&= F(h_{\Gamma}^{01})F(h_{\Gamma_1}^{12})\tau_{\rho}^{A} \circ \left(\phi_{h^{01}_{\Gamma}, h_{\Gamma_1}}^{F(\overline{\rho})A}\right)^{-1} \circ \left(\phi_{h^{12}_{\Gamma_1} \circ h^{01}_{\Gamma}, \rho}\right)^{-1} \circ  F(h_{\Gamma}^{02})\tau_{f_3}^{A} \\
&\circ \left(\phi_{h_{\Gamma}^{02}, f_3}^{A}\right)^{-1} \circ \phi_{f,h_{\Gamma^{\prime}}^{02}}^{A}
\end{align*}
and realizing that by combining our derivations above we have calculated and manipulated the morphism $\tau_f^{\varphi^{\sharp}(\psi^{\sharp}A)} \circ F(\of)\quot{\alpha_A}{\Gamma^{\prime}}$, to prove that the two composites coincide it suffices to verify that the morphism
\begin{align*}
&F(h_{\Gamma}^{01})F(h_{\Gamma_1}^{12})\left(F(\rho)\tau_{f_3}^{A} \circ \left(\phi_{\rho,f_3}^{A}\right)^{-1} \circ \phi_{f_2,\rho^{\prime}}^{A}\right) \circ F(h_{\Gamma}^{01})\left(\left(\phi_{h_{\Gamma_1}, f_2}^{F(\overline{\rho}^{\prime})A}\right)^{-1} \circ \phi_{f_1, h_{\Gamma_1^{\prime}}}^{F(\overline{\rho}^{\prime})A}\right) \\
&\circ F(\of)\left(\left(\phi_{h_{\Gamma^{\prime}}, h_{\Gamma^{\prime}_1}}^{F(\overline{\rho}^{\prime})A}\right)^{-1} \circ \left(\phi_{h_{\Gamma_1^{\prime}} \circ h_{\Gamma^{\prime}}, \rho^{\prime}}^{A}\right)^{-1}\right)
\end{align*}
is equal to the morphism
\[
\left(\phi_{h^{01}_{\Gamma}, h_{\Gamma_1}}^{F(\overline{\rho})A}\right)^{-1} \circ \left(\phi_{h^{12}_{\Gamma_1} \circ h^{01}_{\Gamma}, \rho}\right)^{-1} \circ  F(h_{\Gamma}^{02})\tau_{f_3}^{A} \circ \left(\phi_{h_{\Gamma}^{02}, f_3}^{A}\right)^{-1} \circ \phi_{f,h_{\Gamma^{\prime}}^{02}}^{A}.
\]
Proving this is an extremely tedious but routine verification that follows from the pseudofunctoriality of $F$ and the induced relations provided by the pasting diagram:
\[
\begin{tikzcd}
F(\quot{X}{\Gamma_3^{\prime}}) \ar[dddrr, ""{name = LeftMid}] \ar[rrrrrr, bend left = 50, ""{name = UU}]{}{F(h_{\Gamma^{\prime}}^{02})} \ar[rr, ""{name = UpLeft}]{}[description]{F(\overline{\rho}^{\prime})} \ar[ddd, swap]{}{F(\of_3)} & & F(\quot{X}{\Gamma_2^{\prime}}) \ar[dddrr, ""{name = MM}]\ar[rr, ""{name = UL}]{}[description]{F(h_{\Gamma^{\prime}}^{12})} \ar[ddd, swap]{}{F(\of_2)} & & F(\quot{X}{\Gamma_1^{\prime}}) \ar[dddrr, ""{name = RightMid}]{}{} \ar[ddd, swap]{}{F(\of_1)} \ar[rr, ""{name = UpRight}]{}[description]{F(h_{\Gamma^{\prime}}^{01})} & & F(\quot{X}{\Gamma}) \ar[ddd]{}{F(\of)} \\
\\
\\
F(\quot{X}{\Gamma_3}) \ar[rrrrrr, swap, bend right = 50, ""{name = LL}]{}{F(h_{\Gamma}^{02})} \ar[rr, ""{name = LowLeft}]{}[description]{F(\overline{\rho})} & & F(\quot{X}{\Gamma_2}) \ar[rr, swap, ""{name = LU}]{}[description]{F(h_{\Gamma_1}^{12})} & & F(\quot{X}{\Gamma_1}) \ar[rr, ""{name = LowRight}]{}[description]{F(h_{\Gamma}^{01})} & & F(\quot{X}{\Gamma}) 
\ar[from = UU, to = UL,  Rightarrow, shorten <= 4pt, shorten >= 4pt]{}{(\phi_{h_{\Gamma_1^{\prime} \circ h_{\Gamma^{\prime}, \rho}}} \ast \phi_{h_{\Gamma^{\prime}}, h_{\Gamma_1^{\prime}}})^{-1}}
\ar[from = LU, to = LL, Rightarrow, shorten <= 4pt, shorten >= 4pt]{}{\phi_{h_{\Gamma_1} \circ h_{\Gamma}, \rho} \ast \phi_{h_{\Gamma}, h_{\Gamma_1}}}
\ar[from = UpLeft, to = LeftMid, Rightarrow, shorten >= 4pt, shorten <= 4pt, near end]{}{\phi_{ f_2,\rho^{\prime}}}
\ar[from = LeftMid, to = LowLeft, Rightarrow, shorten >= 4pt, shorten <= 4pt, swap]{}{(\phi_{ \rho,f_3})^{-1}}
\ar[from = UL, to = MM, Rightarrow, shorten <= 4pt, shorten >= 4pt, near end]{}{\phi_{f_1,h_{\Gamma^{\prime}}}}
\ar[from = MM, to = LU, Rightarrow, swap, shorten >=4pt, shorten <= 4pt]{}{(\phi_{h_{\Gamma_1}, f_2})^{-1}}
\ar[from = UpRight, to = RightMid, Rightarrow, near end, shorten <= 4pt, shorten >= 4pt]{}{\phi_{f,h_{\Gamma^{\prime}}}}
\ar[from = RightMid, to = LowRight, Rightarrow, swap, shorten <= 4pt, shorten >= 4pt]{}{(\phi_{h_{\Gamma},f_1})^{-1}}
\end{tikzcd}
\]
and the fact that it is equal to the invertible $2$-cell:
\[
\begin{tikzcd}
F(\quot{X}{\Gamma_3^{\prime}}) \ar[rrr, ""{name = U}]{}{F(h_{\Gamma^{\prime}}^{02})} \ar[d, swap]{}{F(\of_3)} & & & F(\XGammap) \ar[d]{}{F(\of)} \\
F(\quot{X}{\Gamma_3}) \ar[rrr, swap, ""{name = L}]{}{F(h_{\Gamma}^{02})} & & & F(\XGamma) \ar[from = U, to = L, Rightarrow, shorten <= 4pt, shorten >= 4pt]{}{\phi_{h_{\Gamma}^{02},f_3}^{-1} \circ \phi_{f,h_{\Gamma^{\prime}}^{02}}}
\end{tikzcd}
\]
Doing this shows that the two morphisms are indeed equivalent and establishes that
\[
\tau_{f}^{\varphi^{\sharp}(\psi^{\sharp}A)} \circ F(\overline{f})\quot{\alpha_A}{\Gamma^{\prime}} = \quot{\alpha_A}{\Gamma} \circ \tau_f^{(\psi \circ \varphi)^{\sharp}A},
\]
which in turn completes the verification that $\alpha_A:(\psi \circ \varphi)^{\sharp}A \to \varphi^{\sharp}(\psi^{\sharp}A)$ is a morphism in $F_{G_0}(X)$.

We finish by verifying the naturality of $\alpha$, i.e., that for any morphism
\[
\Sigma = \lbrace \quot{\sigma}{\Gamma} \; | \; \Gamma \in \Sf(G_2)_0 \rbrace \in F_{G_2}(X)(A,B)
\]
the diagram
\[
\xymatrix{
(\psi \circ \varphi)^{\sharp}A \ar[r]^-{\alpha_A} \ar[d]_{(\psi \circ \varphi)^{\sharp}\Sigma} & \varphi^{\sharp}(\psi^{\sharp}A) \ar[d]^{\varphi^{\sharp}(\psi^{\sharp}\Sigma)} \\
(\psi \circ \varphi)^{\sharp}B \ar[r]_-{\alpha_B} & \varphi^{\sharp}(\psi^{\sharp}B)
}
\]
commutes. Now observe that since $\Sigma$ is an $\Sf(G_2)$-morphism, for any $\Gamma \in \Sf(G_0)_0$ the diagram
\[
\xymatrix{
F(\overline{\rho})\quot{A}{\Gamma_3} \ar[rr]^-{F(\overline{\rho})\quot{\sigma}{\Gamma_3}} \ar[d]_{\tau_{\rho}^{A}} & & F(\overline{\rho})\quot{B}{\Gamma_3} \ar[d]^{\tau_{\rho}^{B}} \\
\quot{A}{\Gamma_2} \ar[rr]_-{\sigma_{\Gamma_2}} & & \quot{B}{\Gamma_2}
}
\]
commutes. Using this we verify that for an arbirtrary $\Gamma \in \Sf(G_0)_0$,
\begin{align*}
&\quot{\alpha_B}{\Gamma} \circ \quot{(\psi \circ \varphi)^{\sharp}\Sigma}{\Gamma} \\
&= F(h_{\Gamma}^{01})F(h_{\Gamma_1}^{12})\tau_{\rho}^{B} \circ \left(\phi_{h_{\Gamma},h_{\Gamma_1}}^{F(\overline{\rho})B}\right)^{-1} \circ \left(\phi_{h_{\Gamma_1} \circ h_{\Gamma}, \rho}^{B}\right)^{-1} \circ F(h_{\Gamma}^{02})\quot{\sigma}{\Gamma_3} \\
&=F(h_{\Gamma}^{01})F(h_{\Gamma_1}^{12})\tau_{\rho}^{B} \circ \left(\phi_{h_{\Gamma},h_{\Gamma_1}}^{F(\overline{\rho})B}\right)^{-1} \circ \left(\phi_{h_{\Gamma_1} \circ h_{\Gamma}, \rho}^{B}\right)^{-1} \circ F(\overline{\rho} \circ h_{\Gamma_1}^{12} \circ h_{\Gamma}^{01})\quot{\sigma}{\Gamma_3} \\
&= F(h_{\Gamma}^{01})F(h_{\Gamma_1}^{12})\tau_{\rho}^{B} \circ \left(\phi_{h_{\Gamma},h_{\Gamma_1}}^{F(\overline{\rho})B}\right)^{-1} \circ F(h_{\Gamma_1}^{12} \circ h_{\Gamma}^{01})F(\overline{\rho})\quot{\sigma}{\Gamma_3} \circ \left(\phi_{h_{\Gamma_1} \circ h_{\Gamma}, \rho}^{A}\right)^{-1} \\
&= F(h_{\Gamma}^{01})F(h_{\Gamma_1}^{12})\tau_{\rho}^{B} \circ F(h_{\Gamma}^{01})F(h_{\Gamma_1}^{12})F(\overline{\rho})\quot{\sigma}{\Gamma_3} \circ \left(\phi_{h_{\Gamma},h_{\Gamma_1}}^{F(\overline{\rho})A}\right)^{-1} \circ \left(\phi_{h_{\Gamma_1} \circ h_{\Gamma}, \rho}^{A}\right)^{-1} \\
&= F(h_{\Gamma}^{01})F(h_{\Gamma_1}^{12})\left(\tau_{\rho}^{B} \circ F(\overline{\rho})\quot{\sigma}{\Gamma_3}\right) \circ \left(\phi_{h_{\Gamma},h_{\Gamma_1}}^{F(\overline{\rho})A}\right)^{-1} \circ \left(\phi_{h_{\Gamma_1} \circ h_{\Gamma}, \rho}^{A}\right)^{-1} \\
&= F(h_{\Gamma}^{01})F(h_{\Gamma_1}^{12})\left(\quot{\sigma}{\Gamma_2} \circ \tau_{\rho}^{A}\right) \circ \left(\phi_{h_{\Gamma},h_{\Gamma_1}}^{F(\overline{\rho})A}\right)^{-1} \circ \left(\phi_{h_{\Gamma_1} \circ h_{\Gamma}, \rho}^{A}\right)^{-1} \\
&= F(h_{\Gamma}^{01})F(h_{\Gamma_1}^{12})\quot{\sigma}{\Gamma_2} \circ F(h_{\Gamma}^{01})F(h_{\Gamma_1}^{12})\tau_{\rho}^{A} \circ \left(\phi_{h_{\Gamma},h_{\Gamma_1}}^{F(\overline{\rho})A}\right)^{-1} \circ \left(\phi_{h_{\Gamma_1} \circ h_{\Gamma}, \rho}^{A}\right)^{-1} \\
&= \quot{\varphi^{\sharp}(\psi^{\sharp}\Sigma)}{\Gamma} \circ \quot{\alpha_A}{\Gamma},
\end{align*}
which establishes the naturality of $\alpha$.
\end{proof}

As a final topic for this section, we discuss the de-equivariantification functor. This is a useful consequence of the change of groups formalism; it gives rise to a functor $F_G(X) \to F(X)$ for any group $G$ and any space $X$. Since we always have a morphism of algebraic groups $1_G:\Spec K \to G$ for any algebraic group $G$, from Theorem \ref{Thm: Section 3: Change of groups functor} we get a functor $1_G^{\sharp}:F_G(X) \to F_{\Spec K}(X)$ that restricts along this inclusion of the identity of $G$. Recall that there is an equivalence of categories $F_{\Spec K}(X) \simeq F(X)$ (cf.\@ Proposition \ref{Prop: Section 2: Equivalence of K equivariant cat with F X}); let $\overline{F}:F_{\Spec K}(X) \to F(X)$ be one such functor witnessing the equivalence. We then obtain a functor
\[
F_G(X) \xrightarrow{1_G^{\sharp}}  F_{\Spec K}(X) \xrightarrow{\overline{F}} F(X)
\]
which suitably ``de-equivariantifies'' the category $F_G(X)$. We call the corresponding functor $\overline{F} \circ 1_G^{\ast}:F_G(X) \to F(X)$ the de-equivariantification functor. In fact, this can be seen as an explanation using the equivariant language exactly what the forgetful functor $\Forget:F_G(X) \to F(X)$ is doing. While the proof we present below is functionally very similar to that of Proposition \ref{Prop: Section 3.3: Chagnge of groups interacting with change of grups} above, we present it in detail below to illustrate how the forgetful functor is really a Change of Groups in disguise.
\begin{proposition}\label{Prop: Section 3.3: Forgetful functor is de-equivariantification functor}
Let $G$ be a smooth algebraic group and let $X$ be a left $G$-variety with unit map $1_G:\Spec K \to G$. Then the diagram below commutes up to invertible $2$-cell:
\[
\begin{tikzcd}
F_G(X) \ar[rr]{}{1_G^{\sharp}} \ar[dr, swap]{}{\Forget} & {} & F_{\Spec K}(X) \ar[dl]{}{\overline{F}} \\
 & F(X) \ar[from = 2-2, to = 1-2, Rightarrow, shorten >= 4pt, shorten <= 4pt]{}{\alpha}
\end{tikzcd}
\]
\end{proposition}
\begin{proof}
We begin by observing that there is a $G$-equivariant isomorphism of varieites in $\Sf(G)$ of the form
\[
(G \times \Spec K)_{\Spec K} \cong G
\] 
(this can be seen in set-theoretic terms via the mappings $[g,1]_{1} \mapsto g, g \mapsto [g,1]_{1}$), and this isomorphism induces the further isomorphism
\[
\quot{X}{(G \times \Spec K)_{\Spec K}} \xrightarrow[\rho]{\cong} \quot{X}{G}.
\]
Now consider the isomorphism
\[
h_{\Spec K}:\quot{X}{\Spec K} \xrightarrow{\cong} \quot{X}{(G \times \Spec K)_{\Spec K}}
\] 
and write the isomorphism from $X$ to $\quot{X}{\Spec K}$ of Proposition \ref{Prop: Section 2: Equivalence of K equivariant cat with F X} as
\[
X \xrightarrow[\varphi]{\cong} \quot{X}{\Spec K}.
\]
Finally, let $\psi:X \to \quot{X}{G}$ be the standard isomorphism. Let us now observe that on one hand we have that for any object $A \in F_G(X)_0$
\begin{align*}
\overline{F}(1_G^{\sharp}(A)) &= \overline{F}\left(\lbrace F(h_{\Gamma})\quot{A}{(G \times \Gamma)_{\Spec K}} \; | \; \Gamma \in \Sf(G)_0 \rbrace\right) \\
&= F(\varphi)F(h_{\Spec K})\quot{A}{(G \times \Spec K)_{\Spec K}}
\end{align*}
and similarly for morphisms. On the other hand we have that
\[
\Forget(A) = F(\psi)\quot{A}{G}.
\]
Now since the diagram
\[
\xymatrix{
\quot{X}{(G \times \Spec K)_{\Spec K}} \ar[d]_{\overline{\rho}} & & \quot{X}{\Spec K} \ar[ll]_-{h_{\Spec K}} \\
\quot{X}{G} & & X \ar[ll]^-{\psi} \ar[u]_{\varphi}
}
\]
of varieties commutes by the uniqueness of isomorphisms between quotient schemes, we have the following pasting diagram
\[
\begin{tikzcd}
 & F(\quot{X}{\Spec K}) \ar[dr]{}{F(\varphi)} & \\
F(\quot{X}{(G \times \Spec K)_{\Spec K}}) \ar[rr]{}[description]{F( h_{\Spec K} \circ \varphi)} \ar[ur]{}{F(h_{\Spec K})} & {} & F(X) \\
 & F(\quot{X}{G}) \ar[ur, swap]{}{F(\psi)} \ar[ul]{}{F(\overline{\rho})} & \ar[from = 2-2, to = 3-2, swap, Rightarrow, shorten >= 4pt, shorten <= 4pt]{}{\phi_{h \circ \varphi, \rho}} \ar[from = 1-2, to = 2-2, Rightarrow, shorten <= 4pt, shorten >=4pt, swap]{}{\phi_{\varphi, h}}
\end{tikzcd}
\]
in $\fCat$. Notice that because $\rho \in \Sf(G)_1$, we have an isomorphism
\[
\tau_{\rho}^{A}:F(\overline{\rho})\quot{A}{G} \xrightarrow{\cong} \quot{A}{(G \times \Spec K)_{\Spec K}}
\]
in $T_A$. Applying $F(\varphi) \circ F(h_{\Spec K})$ to this isomorphism gives us an isomorphism
\[
\xymatrix{
F(\varphi)F(h_{\Spec K})F(\overline{\rho})\quot{A}{G} \ar[drr] \ar[rr]^-{\cong}_-{F(\varphi)F(h_{\Spec K})\tau_{\rho}^{A}} & & F(\varphi)F(h_{\Spec K})\quot{A}{(G \times \Spec K)_{\Spec K}} \ar@{=}[d] \\
 & & (\overline{F}\circ 1_G^{\sharp})(A)
}
\]
To get our isomorphism to have the desired domain of $\Forget(A) = F(\psi)\quot{A}{G} = F(\overline{\rho} \circ h_{\Spec K} \circ \varphi)\quot{A}{G}$, we pre-compose by the inverses of the two natural isomorphisms to produce the isomorphism
\[
\alpha_A:= F(\varphi)F(h_{\Spec K})\tau_{\rho}^{A} \circ \left(\phi_{\varphi,h}^{F(\overline{\rho})\quot{A}{G}}\right)^{-1} \circ \left(\phi_{h \circ \varphi, \rho}^{\quot{A}{G}}\right)^{-1}.
\]
Then $\alpha_A:\Forget(A) \xrightarrow{\cong} (\overline{F} \circ 1_G^{\sharp})(A)$, as desired.

We end by showing that $\alpha$ is natural in $F_G(X)$, i.e.,that for a morphism
\[
\Sigma = \lbrace \quot{\sigma}{\Gamma} \; | \; \Gamma \in \Sf(G)_0 \rbrace \in F_G(X)(A,B)
\]
that the diagram
\[
\xymatrix{
\Forget(A) \ar[r]^-{\alpha_A} \ar[d]_{\Forget(\Sigma)} & (\overline{F} \circ 1_G^{\sharp})(A) \ar[d]^{(\overline{F}\circ 1_G^{\sharp})(\Sigma)} \\
\Forget(B) \ar[r]_-{\alpha_B} & (\overline{F} \circ 1_G^{\sharp})(B)
}
\]
commutes. First recall that since $\Sigma$ is an $F_G(X)$-morphisms, the diagram
\[
\xymatrix{
F(\overline{\rho})\quot{A}{G} \ar[d]_{\tau_{\rho}^{A}} \ar[rrr]^-{F(\overline{\rho})\quot{\sigma}{G}} & & & F(\overline{\rho})\quot{B}{G} \ar[d]^{\tau_{\rho}^{B}} \\
\quot{A}{(G \times \Spec K)_{\Spec K}} \ar[rrr]_-{\quot{\sigma}{(G \times \Spec K)_{\Spec K}}} & & & \quot{B}{(G \times \Spec K)_{\Spec K}}
}
\]
commutes in $F(\quot{X}{(G \times \Spec K)_{\Spec K}}).$ We then conclude the proof by calculating
\begin{align*}
&(\overline{F}\circ 1_G^{\sharp})(\Sigma) \circ \alpha_A \\
&= F(\varphi)F(h_{\Spec K})\quot{\sigma}{(G \times \Spec K)_{\Spec K}} \circ F(\varphi)F(h_{\Spec K})\tau_{\rho}^{A} \circ \left(\phi_{\varphi,h}^{F(\overline{\rho})\quot{A}{G}}\right)^{-1} \\
&\circ \left(\phi_{h \circ \varphi, \rho}^{\quot{A}{G}}\right)^{-1} \\
&= F(\varphi)F(h_{\Spec K})\left(\quot{\sigma}{(G \times \Spec K)_{\Spec K}} \circ \tau_{\rho}^{A}\right) \circ \left(\phi_{\varphi,h}^{F(\overline{\rho})\quot{A}{G}}\right)^{-1} \circ \left(\phi_{h \circ \varphi, \rho}^{\quot{A}{G}}\right)^{-1} \\
&= F(\varphi)F(h_{\Spec K})F(\tau_{\rho}^{B} \circ F(\overline{\rho})\quot{\sigma}{G}) \circ \left(\phi_{\varphi,h}^{F(\overline{\rho})\quot{A}{G}}\right)^{-1} \circ \left(\phi_{h \circ \varphi, \rho}^{\quot{A}{G}}\right)^{-1} \\
&= F(\varphi)F(h_{\Spec K})(\tau_{\rho}^{B}) \circ \left(\phi_{\varphi,h}^{F(\overline{\rho})\quot{B}{G}}\right)^{-1} \circ \left(\phi_{h \circ \varphi, \rho}^{\quot{B}{G}}\right)^{-1} \circ F(\psi)\quot{\sigma}{G} \\
&= \alpha_B \circ \Forget(\Sigma).
\end{align*}
\end{proof}

\chapter{Equivariant Triangulated Categories and Equivariant $t$-Structures}\label{Chapter 5}
We now move to the final topic of Part \ref{Chapter: EQCats} of the thesis: That of equivariant triangulated categories and $t$-structures on equivariant triangulated categories. Much of the structure that we describe and study in this section can be seen as an application to the various results proved about equivariant categories and equivariant functors (as well as a more in-depth study of triangulated pre-equivariant pseudofunctors; cf.\@ Definition \ref{Defn: Additive, triangulated, and locally small pseudofunctors}). However, it will be important for applications (such as when studying how the equivariant derived category relates to equivariant perverse sheaves, equivariant sheaves, or equivariant local systems) to have developed the language necessary to discuss such objects in a reasonable level of generality. Of particular importance is the theorem below (cf.\@ Theorem \ref{Theorem: Heart of equivariant t-structure is equivariant heart of t-structure}) which says that if $F$ is a triangulated pre-equivariant pseudofunctor on $X$ such that each functor $F(\of)$ restricts to an appropriate functor $F(\of)^{\heartsuit}$ between the hearts of the fibre categories $F(\XGamma)^{\heartsuit}$, then there is a canonical $t$-structure on $F_G(X)$ such that
\[
F_G(X)^{\heartsuit} = (F^{\heartsuit})_{G}(X).
\]
I like to think of this as the ``Change of Heart'' theorem.

\section{Triangulated Equivariant Categories}\label{Section: Triangulated Cats}
In the journey of describing the triangulations on equivariant categories, we will take our first steps by showing how to take an additive pre-equivariant pseudofunctor with suspensions on each fibre category and turn it into a suspended additive equivariant category. A benefit of this is that it explains why the shift functor on the equivariant derived category $[1]:D_G^{b}(X) \to D_G^{b}(X)$ (cf.\@ Example \ref{Example: Shift functors}) is an autoequivalence.

\begin{definition}
We say that an additive category $\Cscr$ is a category with suspension if there exists an auto-equivalence
\end{definition}
\begin{proposition}\label{Prop: Section Triangle: Equivariant suspensions}
Let $F:\SfResl_G(X)^{\op} \to \fCat$ be an additive pre-equivariant functor such that each category $F(\quot{X}{\Gamma})$ is a category with additive autoequivalence $\quot{\Sigma}{\Gamma}:F(\XGamma) \to F(\XGamma)$. Then $F_G(X)$ admits an additive autoequivalence $\Sigma:F_G(X) \to F_G(X)$ as well.
\end{proposition}
\begin{proof}
By Corollary \ref{Cor: Section 2: Additive pseudofunctor gives additive equivariat cat} we have that $F_G(X)$ is additive, while Theorem \ref{Thm: Section 3: Gamma-wise adjoints lift to equivariant adjoints} gives the adjoint $\Omega \dashv \Sigma$, where $\Omega:F_G(X) \to F_G(X)$ is the loop functor given by $\Omega = \lbrace \quot{\Omega}{\Gamma} \; | \; \Gamma \in \Sf(G)_0 \rbrace$ and $\Sigma:F_G(X) \to F_G(X)$ is the suspension functor given by $\Sigma = \lbrace \quot{\Sigma}{\Gamma} \; | \; \Gamma \in \Sf(G)_0 \rbrace$. By Corollary \ref{Cor: Section 3: Fibre-wise equivalences are equivalences of equivariant categories} this is an equivalence of categories. Finally, Corollary \ref{Cor: Section 3: Additive equivraint change of fibre functors} implies because each of the functors $\quot{\Sigma}{\Gamma}$ and $\quot{\Omega}{\Gamma}$ are additive, so are $\Omega$ and $\Sigma$.
\end{proof}
\begin{corollary}
The shift functors $[1]:D_G^b(X) \to D_G^b(X)$ and \\$[1]:\DbeqQl{X} \to \DbeqQl{X}$ are equivalences of categories.
\end{corollary}
Before preceeding, we recall the axioms of a triangulated category as spelled out in \cite{MayTrace}.
\begin{definition}
Let $\Ascr$ be an additive category with an additive autoequivalence $[1]:\Ascr \to \Ascr$. A triangle in $\Ascr$ is a triple of morphisms in $\Ascr$ of the form:
\[
\xymatrix{
	A \ar[r]^-{f} & B \ar[r]^-{g} & C \ar[r]^-{h} & A[1]
}
\]
Two triangle are said to be isomorphic if there is a commuting diagram of the form
\[
\xymatrix{
	A \ar[r]^-{f} \ar[d]_{\rho} & B \ar[r]^-{g} \ar[d]_{\varphi} & C \ar[r]^-{h} \ar[d]^{\psi} & A[1] \ar[d]^{\rho[1]} \\
	X \ar[r]_-{f^{\prime}} & Y \ar[r]_-{g^{\prime}} & Z \ar[r]_-{h^{\prime}} & X[1]
}
\]
where $\rho, \varphi, \psi$ are all isomorphisms.
\end{definition}
\begin{definition}{\cite{MayTrace}}
A triangulated category is an additive category $\Ascr$ with an autoequivalence $[1]:\Ascr \to \Ascr$ equipped with a collection of triangles called distinguished triangles which satisfy the following axioms:
\begin{itemize}
	\item[T1.] Let $A \in \Ascr_0$ and let $f:A \to B$ be in $\Ascr_1$. Then:
	\begin{itemize}
		\item The triangle $A \xrightarrow{\id} A \xrightarrow{0} 0 \to A[1]$ is distinguished;
		\item There exists a distinguished triangle which starts at $f$, i.e., there exist maps $g:B \to C, h:C \to A[1]$ such that the triangle 
		\[
		\xymatrix{
			A \ar[r]^-{f} & B \ar[r]^-{g} & C \ar[r]^-{h} & A[1]
		}
		\]
		is distinguished in $\Ascr$;
		\item Any triangle isomorphic to a distinguished triangle is itself distinguished.
	\end{itemize}
	\item[T2.] If the triangle
	\[
	\xymatrix{
		A \ar[r]^-{f} & B \ar[r]^-{g} & C \ar[r]^-{h} & A[1]
	}
	\]
	is distinguished, then so is:
	\[
	\xymatrix{
	B \ar[r]^-{g} & C \ar[r]^-{h} & A[1] \ar[r]^-{-f[1]} & B[1]
	}
	\]
	\item[T3.] Consider a diagram of the form
	\[
	\begin{tikzcd}
	A \ar[dr, swap]{}{f_0} \ar[rr, bend left = 30]{}{h_0} & & C \ar[rr, bend left = 30]{}{g_1} & & N \ar[rr, bend left = 30]{}{k_2} \ar[dr]{}{g_2} & & L[1] \\
	& B \ar[ur]{}{g_0} \ar[dr, swap]{}{f_1} & & M & & B[1] \ar[ur, swap]{}{f_1[1]} \\
	& & L \ar[rr, bend right = 30, swap]{}{f_2} & & A[1] \ar[ur, swap]{}{f_0[1]}
	\end{tikzcd}
	\]
	where the triangles $(f_0,f_1,f_2)$ and $(g_0, g_1, g_2)$ are distinguished, $h_0 = g_0 \circ f_0$, and $k_2 =f_1[1] \circ  g_2$. If there are morphisms $h_1:C \to M$ and $h_2:M \to A[1]$ which make the triangle $(h_0, h_1, h_2)$ distinguished, then there exist maps $k_0:L \to M$ and $k_1:M \to N$ which make the triangle $(k_0, k_1, k_2)$ distinguished and the diagram
	\[
		\begin{tikzcd}
	A \ar[dr, swap]{}{f_0} \ar[rr, bend left = 30]{}{h_0} & & C \ar[dr]{}{h_1} \ar[rr, bend left = 30]{}{g_1} & & N \ar[rr, bend left = 30]{}{k_2} \ar[dr]{}{g_2} & & L[1] \\
	& B \ar[ur]{}{g_0} \ar[dr, swap]{}{f_1} & & M  \ar[ur]{}{k_1} \ar[dr]{}{h_2} & & B[1] \ar[ur, swap]{}{f_1[1]} \\
	& & L \ar[ur]{}{k_0} \ar[rr, bend right = 30, swap]{}{f_2} & & A[1] \ar[ur, swap]{}{f_0[1]}
	\end{tikzcd}
	\]
	commute in $\Ascr$.
\end{itemize}
\end{definition}

The above construction is all we need to be able to prove that if the pre-equivariant pseudofunctor $F$ is triangulated (cf.\@ \ref{Defn: Additive, triangulated, and locally small pseudofunctors}) then so is $F_G(X)$. However, when proving that $F_G(X)$ is triangulated we will follow the refinement of the triangulation axioms as developed by May and presented in \cite{MayTrace} as opposed to the original axiomatization of Verdier in \cite{VerdierDerivedCat}. This refinement has the benefit of having fewer technical conditions to check: May's refinement only asks us to check Verdier's axiom TR1, a weaker form of Verdier's axiom TR2, and Verdier's axiom TR4; TR3 is entirely redundant. We also follow the naming and numbering of the triangulation axioms of \cite{MayTrace} instead of \cite{VerdierDerivedCat}.

\begin{Theorem}\label{Theorem: Section Triangle: Equivariant triangulation}
Let $F:\SfResl_G(X)^{\op} \to \fCat$ be a triangulated pre-equivariant pseudofunctor. Then the category $F_G(X)$ is a triangulated category with triangulation induced by saying that a triangle
\[
\xymatrix{
A \ar[r]^-{P} & B \ar[r]^-{\Sigma} & C \ar[r]^-{\Phi} & A[1]
}
\]
is distinguished if and only if for all $\Gamma \in \Sf(G)_0$ the triangle
\[
\xymatrix{
\AGamma \ar[r]^-{\quot{\rho}{\Gamma}} & \BGamma \ar[r]^-{\quot{\sigma}{\Gamma}} & \quot{C}{\Gamma} \ar[r]^-{\quot{\varphi}{\Gamma}} & \AGamma[1]
}
\]
is distinguished in $F(\XGamma)$.
\end{Theorem}
\begin{proof}
We begin by examining the first verifying Axiom T1 of \cite{MayTrace}. With this in mind let $A \in F_G(X)_0$. Then because for each $\Gamma \in \Sf(G)_0$ the triangle
\[
\xymatrix{
\AGamma \ar@{=}[r] & \AGamma \ar[r] & \quot{0}{\Gamma} \ar[r] & \AGamma[1]
}
\]
is distinguished in $F(\XGamma)$, it follows that the triangle
\[
\xymatrix{
A \ar@{=}[r] & A \ar[r] & 0 \ar[r] & A[1]
}
\]
is as well.

We now verify the next aspect of Axiom T1. Assume that the triangle
\[
\xymatrix{
A \ar[r]^-{P} & B \ar[r]^-{\Sigma} & C \ar[r]^-{\Phi} & A[1]
}
\]
is distinguished in $F_G(X)$ and that the triangle
\[
\xymatrix{
A^{\prime} \ar[r]^-{P^{\prime}} & B^{\prime} \ar[r]^-{\Sigma^{\prime}} & C^{\prime} \ar[r]^-{\Phi^{\prime}} & A^{\prime}[1]
}
\]
is isomorphic to $(A,B,C,P,\Sigma,\Phi)$. Then for all $\Gamma \in \Sf(G)_0$ the triangle $(\AGamma, \BGamma, \quot{C}{\Gamma}, \quot{\rho}{\Gamma}, \quot{\sigma}{\Gamma}, \quot{\varphi}{\Gamma})$ is distinguished and $(\AGamma^{\prime},\BGamma^{\prime},\quot{C^{\prime}}{\Gamma},\quot{\rho^{\prime}}{\Gamma},\quot{\sigma^{\prime}}{\Gamma},\quot{\varphi^{\prime}}{\Gamma})$ is isomorphic to this triangle as well. Thus the triangles \\$(\AGamma^{\prime},\BGamma^{\prime},\quot{C^{\prime}}{\Gamma},\quot{\rho^{\prime}}{\Gamma},\quot{\sigma^{\prime}}{\Gamma},\quot{\varphi^{\prime}}{\Gamma})$ are distinguished for all $\Gamma \in \Sf(G)_0$ so we conclude that the triangle
\[
\xymatrix{
A^{\prime} \ar[r]^-{P^{\prime}} & B^{\prime} \ar[r]^-{\Sigma^{\prime}} & C^{\prime} \ar[r]^-{\Phi^{\prime}} & A^{\prime}[1]
}
\]
is distinguished as well.

We now verify the last piece of Axiom T1. Let $P:A \to B$ be a morphism in $F_G(X)$ and note that we now must find a distinguished triangle in $F_G(X)$ which begins at $P$. For this note that for all $\Gamma \in \Sf(G)_0$ we can find distinguished triangles
\[
\xymatrix{
\AGamma \ar[r]^-{\quot{\rho}{\Gamma}} & \BGamma \ar[r]^-{\quot{\sigma}{\Gamma}} & \quot{C}{\Gamma} \ar[r]^-{\quot{\varphi}{\Gamma}} & \AGamma[1]
}
\]
in each category $F(\XGamma)$. Now using that each functor $F(\overline{f})$ is triangulated allows us to produce the commuting diagram
\[
\xymatrix{
F(\of)\AGammap \ar[rr]^-{F(\of)\quot{\rho}{\Gamma^{\prime}}} \ar[d]_{\tau_f^A} & & F(\of)\BGammap \ar[rr]^-{F(\of)\quot{\sigma}{\Gamma^{\prime}}} \ar[d]^{\tau_f^B} & & F(\of)\quot{C}{\Gamma^{\prime}} \ar[rr]^-{F(\of)\quot{\varphi}{\Gamma^{\prime}}} & & F(\of)\AGammap \ar[d]^{\tau_f^{A[1]}} \\
\AGamma \ar[rr]_-{\quot{\rho}{\Gamma}} & & \BGamma \ar[rr]_-{\quot{\sigma}{\Gamma}} & & \quot{C}{\Gamma} \ar[rr]_-{\quot{\varphi}{\Gamma}} & & \AGamma[1]
}
\]
in $F(\XGamma)$ where $f \in \Sf(G)_1$ with $\Dom f = \Gamma$ and $\Codom f = \Gamma^{\prime}$. Applying \cite[Lemmas 2.2, 2.3]{MayTrace} then allows us to deduce the existence of an isomorphism $\alpha:F(\of)\quot{C}{\Gamma^{\prime}} \to \quot{C}{\Gamma}$ which make the diagram
\[
\xymatrix{
F(\of)\AGammap \ar[rr]^-{F(\of)\quot{\rho}{\Gamma^{\prime}}} \ar[d]_{\tau_f^A} & & F(\of)\BGammap \ar[rr]^-{F(\of)\quot{\sigma}{\Gamma^{\prime}}} \ar[d]^{\tau_f^B} & & F(\of)\quot{C}{\Gamma^{\prime}} \ar[d]^{\alpha_f}_{\cong} \ar[rr]^-{F(\of)\quot{\varphi}{\Gamma^{\prime}}} & & F(\of)\AGammap \ar[d]^{\tau_f^{A[1]}} \\
\AGamma \ar[rr]_-{\quot{\rho}{\Gamma}} & & \BGamma \ar[rr]_-{\quot{\sigma}{\Gamma}} & & \quot{C}{\Gamma} \ar[rr]_-{\quot{\varphi}{\Gamma}} & & \AGamma[1]
}
\]
commute. Note that if we can prove that the isomorphisms $\alpha_f$ satisfy the cocycle condition they automatically make 
\[
\Sigma := \lbrace \quot{\sigma}{\Gamma} \; | \; \Gamma \in \Sf(G)_0\rbrace
\] 
and
\[
\Phi := \lbrace \quot{\varphi}{\Gamma} \; | \; \Gamma \in \Sf(G)_0 \rbrace
\]
morphisms in $F_G(X)$. These morphisms can be shown to make the diagrams, for any composable $f,g \in \Sf(G)_1$,
\[
\begin{tikzcd}
F(\of)F\quot{C}{\Gamma^{\prime}} \ar[rr, ""{name = U}]{}{\alpha_f} & & \quot{C}{\Gamma} \\
F(\of)F(\overline{g})\quot{C}{\Gamma^{\prime\prime}} \ar[u]{}{F(\of)\alpha_{g}} \ar[rr, swap, ""{name = L}]{}{\phi_{f,g}^{\quot{C}{\Gamma^{\prime\prime}}}} & & F(\overline{g}\circ \of)\quot{C}{\Gamma^{\prime\prime}} \ar[u, swap]{}{\alpha_{g \circ f}} \ar[from = U, to = L, Rightarrow, shorten >= 4pt, shorten <= 4pt]{}{\cong}
\end{tikzcd}
\]
commute up to a unique invertible $2$-cell in $F(\XGamma)$ because each morphism $\alpha_f$ is induced from a weak limit in $F(\XGamma)$, which exists because $F(\XGamma)$ is a triangulated category and $\alpha_f$ is induced by this weak universal property. As such we invoke the Axiom of Choice to select our isomorphisms $\alpha_f$ such that for all composable arrows $f,g \in \Sf(G)_1$ the identity $\alpha_{g \circ f} \circ \varphi_{f,g} = \alpha_f \circ F(\of)\alpha_g$. With this choice of isomorphisms we find that the pair
\[
(\lbrace \quot{C}{\Gamma} \; | \; \Gamma \in \Sf(G)_0 \rbrace, \lbrace \alpha_f \; | \; f \in \Sf(G)_1\rbrace)
\]
determines an object with morphisms $\Sigma:B \to C$ and $\Phi:C \to A[1]$ which make the triangle
\[
\xymatrix{
A \ar[r]^-{P} & B \ar[r]^-{\Sigma} & C \ar[r]^-{\Phi} & A[1]
}
\]
distinguished in $F_G(X)$. This completes the verification of Axiom T1.

We now proceed to verify Axiom T2. Assume that the triangle
\[
\xymatrix{
A \ar[r]^-{P} & B \ar[r]^-{\Sigma} & C \ar[r]^-{\Phi} & A[1]
}
\]
is distinguished in $F_G(X)$. Then for each $\Gamma \in \Sf(G)_0$ the triangle
\[
\xymatrix{
\AGamma \ar[r]^-{\quot{\rho}{\Gamma}} & \quot{B}{\Gamma} \ar[r]^-{\quot{\sigma}{\Gamma}} & \quot{C}{\Gamma} \ar[r]^-{\quot{\varphi}{\Gamma}} & \quot{A[1]}{\Gamma}
}
\]
is distinguished in $F(\XGamma)$, so applying Axiom T2 in $F(\XGamma)$ tells us that the triangle
\[
\xymatrix{
\BGamma \ar[r]^-{\quot{\sigma}{\Gamma}} & \quot{C}{\Gamma} \ar[r]^-{\quot{\varphi}{\Gamma}} & \AGamma[1] \ar[r]^-{-\quot{\rho[1]}{\Gamma}} & \BGamma[1]
}
\]
is distinguished as well. Because this happens for each $\Gamma \in \Sf(G)_0$ and these give the corresponding $\Gamma$-local descriptions of the objects in $F_G(X)$, we conclude that the triangle
\[
\xymatrix{
B \ar[r]^-{\Sigma} & C \ar[r]^-{\Phi} & A[1] \ar[r]^-{-P[1]} & B[1]
}
\]
is distinguished in $F_G(X)$. This establishes Axiom T2.

Finally we verify that Axiom T3 (Verdier's Axiom) holds in $F_G(X)$. Assume that we have a diagram in $F_G(X)$ of the form
\[
\begin{tikzcd}
A \ar[dr, swap]{}{P_0} \ar[rr, bend left = 30]{}{\Phi_0} & & C \ar[rr, bend left = 30]{}{\Sigma_1} & & N \ar[rr, bend left = 30]{}{\Psi_2} \ar[dr]{}{\Sigma_2} & & L[1] \\
 & B \ar[ur]{}{\Sigma_0} \ar[dr, swap]{}{P_1} & & M & & B[1] \ar[ur, swap]{}{P_1[1]} \\
 & & L \ar[rr, bend right = 30, swap]{}{P_2} & & A[1] \ar[ur, swap]{}{P_0[1]}
\end{tikzcd}
\]
in $F_G(X)$ in which the triangles $(P_0,P_1,P_2)$ and $(\Sigma_0,\Sigma_1,\Sigma_2)$ are both distinguished. Now assume that there are morphisms $\Phi_1:C \to M$ and $\Phi_2:M \to A[1]$ such that the triangle $(\Phi_0, \Phi_1, \Phi_2)$ is distinguished. We now need to construct morphisms $\Psi_0:L \to M$ and $\Psi_1:M \to N$ which make the triangle $(\Psi_0, \Psi_1, \Psi_2)$ and render the diagram commutative. To go about this, let $\Gamma \in \Sf(G)_0$ and invoke Axiom T3 for this $\Gamma$. We then produce morphisms $\quot{\psi_0}{\Gamma}:\quot{L}{\Gamma} \to \quot{M}{\Gamma}$ and $\quot{\psi_1}{\Gamma}:\quot{M}{\Gamma} \to \quot{N}{\Gamma}$ which make the diagram
\[
\begin{tikzcd}
\AGamma \ar[dr, swap]{}{\quot{\rho_0}{\Gamma}} \ar[rr, bend left = 30]{}{\quot{\varphi_0}{\Gamma}} & & \quot{C}{\Gamma} \ar[dr]{}{\quot{\varphi_1}{\Gamma}} \ar[rr, bend left = 30]{}{\quot{\sigma_1}{\Gamma}} & & \quot{N}{\Gamma} \ar[rr, bend left = 30]{}{\quot{\psi_2}{\Gamma}} \ar[dr]{}{\quot{\sigma_2}{\Gamma}} & & \quot{L[1]}{\Gamma} \\
& \BGamma \ar[ur]{}{\quot{\sigma_0}{\Gamma}} \ar[dr, swap]{}{\quot{\rho_1}{\Gamma}} & & \quot{M}{\Gamma} \ar[ur]{}{\quot{\psi_1}{\Gamma}} \ar[dr]{}{\quot{\varphi_2}{\Gamma}} & & \BGamma[1] \ar[ur, swap]{}{\quot{\rho_1[1]}{\Gamma}} \\
& & \quot{L}{\Gamma} \ar[rr, bend right = 30, swap]{}{\quot{\rho_2}{\Gamma}} \ar[ur]{}{\quot{\psi_0}{\Gamma}} & & \AGamma[1] \ar[ur, swap]{}{\quot{\rho_0[1]}{\Gamma}}
\end{tikzcd}
\]
commute with the triangle $(\quot{\psi_0}{\Gamma},\quot{\psi_1}{\Gamma}, \quot{\psi_2}{\Gamma})$ distinguished in addition to the others induced from the distinguished triangles given by assumption. Now using that for any $f \in \Sf(G)_1$ the functor $F(\of)$ is triangulated, together with the fact that each of $P_1, P_2, P_0[1], \Sigma_1, \Sigma_2 \Phi_1,$ and $\Phi_2$ are $F_G(X)$-morphisms allows us to deduce via tedious diagram chases that the diagrams
\[
\begin{tikzcd}
F(\of)\quot{L}{\Gamma^{\prime}} \ar[rr, ""{name = LU}]{}{F(\of)\quot{\psi_0}{\Gamma^{\prime}}} \ar[d]{}{\tau_f^L} & & F(\of)\quot{M}{\Gamma^{\prime}} \ar[d]{}{\tau_f^M} & F(\of)\quot{M}{\Gamma^{\prime}} \ar[rr, ""{name = RU}]{}{F(\of)\quot{\psi_1}{\Gamma^{\prime}}} \ar[d]{}{\tau_f^M} & & F(\of)\quot{N}{\Gamma^{\prime}} \ar[d]{}{\tau_f^N} \\
\quot{L}{\Gamma} \ar[rr, swap, ""{name = LL}]{}{\quot{\psi_0}{\Gamma}} & & \quot{M}{\Gamma} & \quot{M}{\Gamma} \ar[rr, swap, ""{name = RL}]{}{\quot{\psi_1}{\Gamma}} & & \quot{N}{\Gamma}
\ar[from = LU, to = LL, Rightarrow, shorten >= 4pt, shorten <= 4pt]{}{\cong} \ar[from = RU, to = RL, Rightarrow, shorten >= 4pt, shorten <= 4pt]{}{\cong}
\end{tikzcd}
\]
both commute up to a natural isomorphism induced from the $2$-cells that compare any two morphisms which both provide braids of triangles. We once again invoke the Axiom of Choice to redefine the $\quot{\psi_i}{\Gamma}$ if any of these natural isomorphisms are not strict identities and replace them with their strict counterparts. This allows us to deduce that $\Psi_0 := \lbrace \quot{\psi_0}{\Gamma} \; | \; \Gamma \in \Sf(G)_0 \rbrace$ and $\Psi_1 := \lbrace \quot{\psi_1}{\Gamma} \; | \; \Gamma \in \Sf(G)_0 \rbrace$ ar morphisms which makes the triangle $(\Psi_0,\Psi_1, \Psi_2)$ distinguished and renders the briad
\[
\begin{tikzcd}
A \ar[dr, swap]{}{P_0} \ar[rr, bend left = 30]{}{\Phi_0} & & C \ar[rr, bend left = 30]{}{\Sigma_1} \ar[dr]{}{\Phi_1} & & N \ar[rr, bend left = 30]{}{\Psi_2} \ar[dr]{}{\Sigma_2} & & L[1] \\
& B \ar[ur]{}{\Sigma_0} \ar[dr, swap]{}{P_1} & & M \ar[ur]{}{\Psi_1} \ar[dr]{}{\Phi_2} & & B[1] \ar[ur, swap]{}{P_1[1]} \\
& & L \ar[rr, bend right = 30, swap]{}{P_2} \ar[ur]{}{\Psi_0} & & A[1] \ar[ur, swap]{}{P_0[1]}
\end{tikzcd}
\]
commutative. This completes the verification of the Axiom T3 and finishes the proof of the theorem.
\end{proof}
\begin{corollary}
The categories $\DbeqQl{X}$ and $D_G^b(X)$ are triangulated.
\end{corollary}
\begin{remark}
In the above proof we have had to use the Axiom of Choice to artificially select morphisms to satisfy strictness conditions of our transition isomorphisms. This is largely because we are considering $1$-categorical information (strictly commuting diagrams) of what are inherently $2$-categorical constructions: The unique isomorphisms provided by weak limits and colimits. These issues disapper, however, in the $\infty$-categorical formalism, as the lack of strict equality becomes a feature rather than a bug; we record all such coherences and ways of moving between such isomorphisms as various higher morphisms and higher cells. In fact, by working with a notion of equivariant $\infty$-category (an $\infty$-pseudofunctor ``$F:\SfResl_G(X)^{\op} \to \mathfrak{QCat}$''\index[notation]{QCat@$\mathfrak{QCat}$} into the $(\infty,2)$-category of quasi-categories or ``$F:\SfResl_G(X)^{\op} \to \mathscr{K}$'' where $\mathscr{K}$ is an $\infty$-cosmos with a notion of stable $\infty$-category) whose fibre $\infty$-categories are all stable with triangulated $\infty$-fibre functors, we can sidestep these issues entirely. We do not explore these notations explicitly here, and instead defer our conjections to the Futre Work section of this thesis.
\end{remark}
For what follows we will need $t$-structures; as such, we recall their definition here for later use in the thesis.
\begin{definition}[{\cite[Definition 1.3.1]{BBD}}]
Let $\Tscr$ be a triangulated category. A $t$-structure on $\Tscr$ is a pair of full subcategories $(\Tscr^{\leq 0}, \Tscr^{\geq 0})$ which are both stable under equivalence such that the following are satisfied (where for any $n \in \Z$, we define $\Tscr^{\leq n} := \Tscr^{\leq 0}[-n]$ and $\Tscr^{\geq n} := \Tscr^{\geq 0}[-n]$):
\begin{itemize}
	\item  If $A$ is an object of $\Tscr^{\geq 1}$ and $B$ is an object of $\Tscr^{\leq 0}$, then $\Tscr(A,B) = 0$;
	\item We have inclusions of categories $\Tscr^{\leq 0} \hookrightarrow \Tscr^{\leq 1}$ and $\Tscr^{\geq 1} \hookrightarrow \Tscr^{\geq 0}$;
	\item For any object $A \in \Tscr_0$ there exists a truncation distinguished triangle
	\[
	\xymatrix{
	X \ar[r] & A \ar[r] & Y \ar[r] & X[1]
}
	\]
	where $X$ is an object of $\Tscr^{\leq 0}$ and $Y$ is an object of $\Tscr^{\geq 1}$.
\end{itemize}
Furthermore, given a $t$-structure $(\Tscr^{\leq 0}, \Tscr^{\geq 0})$, wee define the heart of $\Tscr$ to be the category
\[
\Tscr^{\heartsuit} := \Tscr^{\leq 0} \cap \Tscr^{\geq 0}.
\]
\end{definition}
\begin{example}\label{Example: Standard tstructure}
Let $\Ascr$ be an Abelian category and let $D(\Ascr)$ be the derived category of chain complexes of $\Ascr$. Define the subcategory $D(\Ascr)^{\leq 0}$ (respectivly $D(\Ascr)^{\geq 0}$) to be the subcategory of $D(\Ascr)$ generated by the complexes $X^{\bullet}$ for which $H^{n}(X^{\bullet}) = 0$ for all $n > 0$ (respectively generated by the complexes $X^{\bullet}$ for which $H^n(X^{\bullet}) = 0$ for all $n < 0$). The pair $(D(\Ascr)^{\leq 0}, D(\Ascr)^{\geq 0})$ is a $t$-structure on $D(\Ascr)$ and its heart is the subcategory of objects whose cohomology is concentrated in degree $0$. In particular, there is an equivalence
\[
D(\Ascr)^{\heartsuit} \simeq \Ascr.
\]
\end{example}
\begin{remark}
For any $t$-structure on a triangulated category $\Tscr$, the heart $\Tscr^{\heartsuit}$ is an Abelian category by \cite[Th{\'e}or{\'e}me 1.3.6]{BBD}. If there are multiple $t$-structures on a category $\Tscr$, we will denote their corresponding subcategories as $({}^{t}\Tscr^{\leq 0}, {}^{t}\Tscr^{\geq 0})$ and $({}^{t^{\prime}}\Tscr^{\leq 0}, {}^{t^{\prime}}\Tscr^{\geq 0})$ while writing the corresponding hearts as $\Tscr^{\heartsuit_t}$ and $\Tscr^{\heartsuit_{t^{\prime}}}$, respectively.
\end{remark}
\begin{definition}\label{Defn: Standard tstructure}
Let $\Dscr$ be a (bounded) derived category of some Abelian category or the bounded derived category of $\ell$-adic sheaves on some variety. We define the standard $t$-structure on $\Dscr$ to be the $t$-structure generated as in Example \ref{Example: Standard tstructure}.
\end{definition}

We now show how to build $t$-structures on $F_G(X)$ based on $t$-structures that come from those on $F(X)$. After proving the theorem (cf.\@ Theorem \ref{Theorem: Section Triangle: t-structure on our friend the ECat}) that says these are indeed $t$-structures on $F_G(X)$, we will then give some explicit examples of these $t$-structures as they may arise in nature. However, before proving this theorem, we will define truncated pre-equivariant pseudofunctors, which essentially are triangulated pre-equivariant pseudofunctors for which each fibre functor admits a $t$-structure and these $t$-structures vary pseudofunctorially in $F$. With this we will prove how to produce $t$-structures from pre-equivariant pseudofunctors and from there how to produce $t$-structures from the category $F(X)$.

\begin{definition}\label{Defn: Truncated preq pseudofunctor}\index{Pre-equivariant Pseudofunctor!Truncated}
Let $F:\SfResl_G(X)^{\op} \to \fCat$ be a triangulated pre-equivariant pseudofunctor. Assume each fibre category is has $t$-structure with truncation functors $(\quot{\tau^{\leq 0}}{\Gamma},\quot{\tau^{\geq 0}}{\Gamma})$ such that these give rise to pre-equivariant pseudofunctors $F^{\leq 0}:\SfResl_G(X)^{\op} \to \fCat$\index[notation]{Fleq0@$F^{\leq 0}$} and $F^{\geq 0}:\SfResl_{G}(X)^{\op} \to \fCat$\index[notation]{Fgeq0@$F^{\geq 0}$} for which the functors $F^{\geq 0}(\of)$ and $F^{\leq 0}(\of)$ are $t$-exact functors which coincide on all common subcategories. Then if for all $f \in \Sf(G)_1$ there are natural isomorphisms
\[
\begin{tikzcd}
F(\XGammap) \ar[r, ""{name = U}]{}{F(\of)} \ar[d, swap]{}{\quot{\tau^{\leq 0}}{\Gamma^{\prime}}} & F(\XGamma) \ar[d]{}{\quot{\tau^{\leq 0}}{\Gamma}} \\
F(\XGammap)^{\leq 0} \ar[r, swap, ""{name = L}]{}{F^{\leq 0}(\of)} & F(\XGamma)^{\leq 0} \ar[from = U, to = L, Rightarrow, shorten <= 4pt, shorten >= 4pt]{}{\theta_{f}^{\leq 0}}
\end{tikzcd}
\]
and
\[
\begin{tikzcd}
F(\XGammap) \ar[r, ""{name = U}]{}{F(\of)} \ar[d, swap]{}{\quot{\tau^{\geq 0}}{\Gamma^{\prime}}} & F(\XGamma) \ar[d]{}{\quot{\tau^{\geq 0}}{\Gamma}} \\
	F(\XGammap)^{\geq 0} \ar[r, swap, ""{name = L}]{}{F^{\geq 0}(\of)} & F(\XGamma)^{\geq 0} \ar[from = U, to = L, Rightarrow, shorten <= 4pt, shorten >= 4pt]{}{\theta_{f}^{\geq 0}}
	\end{tikzcd}
\]
which satisfy the coherence diagrams presented below, we say that the pre-equivariant pseudofunctor is truncated. The pasting diagram
\[
\begin{tikzcd}
F(\XGammapp) \ar[rr, ""{name = UL}]{}{F(\overline{g})} \ar[d, swap]{}{\quot{\tau^{\leq 0}}{\Gamma^{\prime\prime}}} \ar[rrrr, bend left = 30, ""{name = U}]{}{F(\overline{g} \circ \of)} &  & F(\XGammap) \ar[d]{}{\quot{\tau^{\leq 0}}{\Gamma^{\prime}}} \ar[rr, ""{name = UR}]{}{F(\of)} & & F(\XGamma) \ar[d]{}{\quot{\tau^{\leq 0}}{\Gamma}} \\
F(\XGammapp)^{\leq 0} \ar[rr, swap, ""{name = LL}]{}{F^{\leq 0}(\overline{g})} \ar[rrrr, bend right = 30, swap, ""{name = L}]{}{F^{\leq 0}(\overline{g} \circ \of)} & & F(\XGammap)^{\leq 0} \ar[rr, swap, ""{name = LR}]{}{F^{\leq 0}(\of)} & & F(\XGamma)^{\leq 0} \ar[from = U, to = 1-3, Rightarrow, shorten <= 4pt, shorten <= 4pt]{}{\quot{\phi_{f,g}}{F}^{-1}} \ar[from = UR, to = LR, Rightarrow, shorten <= 4pt, shorten >= 4pt]{}{\theta_{f}^{\leq 0}} \ar[from = UL,, to = LL, Rightarrow, shorten <= 4pt, shorten >= 4pt]{}{\theta_g^{\leq 0}} \ar[from = 2-3, to = L, shorten <= 4pt, shorten >= 4pt, Rightarrow]{}{\quot{\phi_{f,g}}{F^{\leq 0}}}
\end{tikzcd}
\]
is equal to the $2$-cell
\[
\begin{tikzcd}
F(\XGammapp) \ar[rr, ""{name = U}]{}{F(\overline{g} \circ \of)} \ar[d, swap]{}{\quot{\tau^{\leq 0}}{\Gamma^{\prime\prime}}} & & F(\XGamma) \ar[d]{}{\quot{\tau^{\leq 0}}{\Gamma}} \\
F(\XGammapp)^{\leq 0} \ar[rr, swap, ""{name = L}]{}{F^{\leq 0}( \overline{g} \circ \of)} & & F(\XGamma)^{\leq 0} \ar[from = U, to = L, Rightarrow, shorten <= 4pt, shorten >= 4pt]{}{\theta_{g \circ f}^{\leq 0}}
\end{tikzcd}
\]
and dually the pasting diagram
\[
\begin{tikzcd}
F(\XGammapp) \ar[rr, ""{name = UL}]{}{F(\overline{g})} \ar[d, swap]{}{\quot{\tau^{\geq 0}}{\Gamma^{\prime\prime}}} \ar[rrrr, bend left = 30, ""{name = U}]{}{F(\overline{g} \circ \of)} &  & F(\XGammap) \ar[d]{}{\quot{\tau^{\geq 0}}{\Gamma^{\prime}}} \ar[rr, ""{name = UR}]{}{F(\of)} & & F(\XGamma) \ar[d]{}{\quot{\tau^{\geq 0}}{\Gamma}} \\
F(\XGammapp)^{\geq 0} \ar[rr, swap, ""{name = LL}]{}{F^{\geq 0}(\overline{g})} \ar[rrrr, bend right = 30, swap, ""{name = L}]{}{F^{\geq 0}(\overline{g} \circ \of)} & & F(\XGammap)^{\geq 0} \ar[rr, swap, ""{name = LR}]{}{F^{\geq 0}(\of)} & & F(\XGamma)^{\geq 0} \ar[from = U, to = 1-3, Rightarrow, shorten <= 4pt, shorten <= 4pt]{}{\quot{\phi_{f,g}}{F}^{-1}} \ar[from = UR, to = LR, Rightarrow, shorten <= 4pt, shorten >= 4pt]{}{\theta_{f}^{\geq 0}} \ar[from = UL,, to = LL, Rightarrow, shorten <= 4pt, shorten >= 4pt]{}{\theta_g^{\geq 0}} \ar[from = 2-3, to = L, shorten <= 4pt, shorten >= 4pt, Rightarrow]{}{\quot{\phi_{f,g}}{F^{\geq 0}}}
\end{tikzcd}
\]
is equal to the $2$-cell:
\[
\begin{tikzcd}
F(\XGammapp) \ar[rr, ""{name = U}]{}{F(\overline{g} \circ \of)} \ar[d, swap]{}{\quot{\tau^{\geq 0}}{\Gamma^{\prime\prime}}} & & F(\XGamma) \ar[d]{}{\quot{\tau^{\geq 0}}{\Gamma}} \\
F(\XGammapp)^{\geq 0} \ar[rr, swap, ""{name = L}]{}{F^{\geq 0}( \overline{g} \circ \of)} & & F(\XGamma)^{\geq 0} \ar[from = U, to = L, Rightarrow, shorten <= 4pt, shorten >= 4pt]{}{\theta_{g \circ f}^{\geq 0}}
\end{tikzcd}
\]
\end{definition}
\begin{proposition}\label{Prop: Section Triangle: Truncated preeq gives truncation functors}
Let $F$ be a truncated pre-equivariant pseudofunctor on $X$. Then there are truncation functors
\[
\tau^{\leq 0}:F_G(X) \to F_G(X), \qquad \tau^{\geq 0}:F_G(X) \to F_G(X).
\]
The truncation $\tau^{\leq 0}$ is defined on objects by
\[
\tau^{\leq 0}(\lbrace \AGamma\; | \; \Gamma \in \Sf(G)_0 \rbrace) = \lbrace \quot{\tau^{\leq 0}}{\Gamma}(\AGamma) \; | \; \Gamma \in \Sf(G)_0 \rbrace,
\]
with transition isomorphisms $\tau_f^{A\leq 0}$ given by the diagram
\[
\xymatrix{
F(\of)\left(\quot{\tau^{\leq 0}}{\Gamma^{\prime}}\AGammap\right) \ar[rr]^-{\theta_f^{\leq 0}} \ar[drr]_{\tau_f^{A \leq 0}} & & \quot{\tau^{\leq 0}}{\Gamma}\left(F(\of)\AGammap\right) \ar[d]^{\quot{\tau^{\leq 0}}{\Gamma}\left(\tau_f^A\right)} \\
 & & \quot{\tau^{\leq 0}}{\Gamma}\AGamma
}
\]
and defined on morphisms by
\[
\tau^{\leq 0}\left(\lbrace \quot{\rho}{\Gamma} \; | \; \Gamma \in \Sf(G)_0 \rbrace\right) = \lbrace \quot{\tau^{\leq 0}}{\Gamma}\quot{\rho}{\Gamma} \; | \; \Gamma \in \Sf(G)_0 \rbrace.
\]
Similarly, the other truncation $\tau^{\geq 0}$ is defined on objects by
\[
\tau^{\geq 0}(\lbrace \AGamma\; | \; \Gamma \in \Sf(G)_0 \rbrace) = \lbrace \quot{\tau^{\geq 0}}{\Gamma}(\AGamma) \; | \; \Gamma \in \Sf(G)_0 \rbrace,
\]
with transition isomorphisms $\tau_f^{A\geq 0}$ given by the diagram
\[
\xymatrix{
	F(\of)\left(\quot{\tau^{\geq 0}}{\Gamma^{\prime}}\AGammap\right) \ar[rr]^-{\theta_f^{\geq 0}} \ar[drr]_{\tau_f^{\A leq 0}} & & \quot{\tau^{\geq 0}}{\Gamma}\left(F(\of)\AGammap\right) \ar[d]^{\quot{\tau^{\geq 0}}{\Gamma}\left(\tau_f^A\right)} \\
	& & \quot{\tau^{\geq 0}}{\Gamma}\AGamma
}
\]
and defined on morphisms by
\[
\tau^{\geq 0}\left(\lbrace \quot{\rho}{\Gamma} \; | \; \Gamma \in \Sf(G)_0 \rbrace\right) = \lbrace \quot{\tau^{\geq 0}}{\Gamma}\quot{\rho}{\Gamma} \; | \; \Gamma \in \Sf(G)_0 \rbrace.
\]
\end{proposition}
\begin{proof}
The technical conditions of this proposition are formally similar to the tensor functors in Theorem \ref{Theorem: Section 2: Monoidal preequivariant pseudofunctor gives monoidal equivariant cat}. As such, proceeding mutatis mutandis as in the proof of Theorem \ref{Theorem: Section 2: Monoidal preequivariant pseudofunctor gives monoidal equivariant cat} proves the proposition.
\end{proof}
By taking appropriate shifts and setting $F_G^{\leq n}(X) = F_G^{\leq 0}(X)[-n], F_G^{\geq n} = F_G^{\geq 0}(X)[-n]$ for $n \in \Z$, we obtain the following corollary.
\begin{corollary}\label{Cor: Section Triangle: Truncation functors for all $n$}
Let $F$ be a truncated pre-equivariant pseudofunctor on $X$. Then for all $n \in \Z$ there are categories $F_G^{\geq n}(X)$ and $F_G^{\leq n}$ equipped with truncation functors $\tau^{\geq n}$ and $\tau^{\leq n}$.
\end{corollary}

\begin{Theorem}\label{Theorem: Section Triangle: t-structure on our friend the ECat}
Let $F$ be a truncated pre-equivariant pseudofunctor. The truncation functors, when shifted for all $n \in \Z$, give rise to a $t$-structure on $F_G(X)$ for which the category of bounded above by degree $0$ objects of $F_G(X)$ is the equivariant category associated to $F^{\leq 0}$ while the category of objects of $F_G(X)$ bounded below by degree $0$ is the equivariant category associated to $F^{\geq 0}$. Notationally this is written as
\[
F^{\leq 0}_G(X) = F_G(X)^{\leq 0}, F_G^{\geq 0}(X) = F_G(X)^{\geq 0}.
\]
\end{Theorem}
\begin{proof}
We first prove that we can realize $F^{\leq 0}_G(X)$ and $F^{\geq 0}_G(X)$ as full subcategories of $F_G(X)$. For this we observe first that since each category $F(\XGamma)^{\leq 0}$ and $F(\XGamma)^{\geq 0}$ is a strictly full subcategory of each fibre category $F(\XGamma)$. Now let $A \in F_G^{\geq 0}(X)$; the case for $F_G^{\leq 0}(X)$ follows dually and so is omitted here. Because $A \in F_G^{\geq 0}(X)$, each fibre $\AGamma \in F{\geq 0}(\XGamma)_0 = (F(\XGamma)^{\geq 0})_{0}$ so in turn each $A \in F_G(X)_0$. Now any morphism $P:A \to B$ in $F_G(X)$ with $A,B \in F_G^{\geq 0}(X)_0$ has $\Gamma$-local components $\quot{\rho}{\Gamma}:\AGamma \to \BGamma$. Because $F(\XGamma)^{\geq 0}$ is full, we get that $\quot{\rho}{\Gamma} \in (F(\XGamma)^{\geq 0})_1$, and because this happens for each $\Gamma \in \Sf(G)_0$, it follows that $P \in F_G^{\geq 0}(X)_1$ as well. This shows the fullness of $F_G^{\geq 0}(X)$; the strictness (i.e., the fact that the subcategory is closed with respect to isomorphisms) follows similarly. We thus conclude that $F_G^{\geq 0}(X)$ and $F_G^{\leq 0}(X)$ are both strictly full subcategories of $F_G(X)$.

We now prove that $F_G(X)^{\geq 0}$ and $F_G(X)^{\leq 0}$ do indeed induce a $t$-structure on $F_G(X)$. For this we recall that the category $F_G(X)^{\geq n}$, for any $n \in \Z$, is induced by the construction on objects
\[
(F_G(X)^{\geq n})_0 = F_G(X)^{\geq 0}[-n] = \lbrace A[-n] \; | \; A \in F_G^{\geq 0}(X)_0 \rbrace.
\]
We now establish that for any $A \in F_G^{\leq 0}(X)$ and any $B \in F_G(X)^{\geq 1}$, the hom-set $F_G(X)(A,B) = 0$. For this we note that since $A \in (F_G(X)^{\leq 0})_0$ and $B \in (F_G(X)^{\geq 1})_0$, for every $\Gamma \in \Sf(G)$ we have that $\AGamma \in (F(\XGamma)^{\leq 0})_0$ and $\BGamma \in (F(\XGamma)^{\geq 1})_0$. Thus we get that
\[
F(\XGamma)(\AGamma,\BGamma) = \quot{0}{\Gamma}
\]
so we calculate that
\begin{align*}
0  &= \lbrace \quot{0}{\Gamma} \; | \; \Gamma \in \Sf(G)_0 \rbrace \subseteq F_G(X)(A,B) \\
&\subseteq \bigcup_{\Gamma \in \Sf(G)_0}F(\XGamma)(\AGamma,\BGamma) = \bigcup_{\Gamma \in \Sf(G)_0} \lbrace \quot{0}{\Gamma} \rbrace = \lbrace \quot{0}{\Gamma} \; | \; \Gamma \in \Sf(G)_0 \rbrace = 0,
\end{align*}
which establishes the claim.

We now prove that the categories $F_G^{\leq 0}(X)$ and $F_G^{\geq 0}(X)$ are stable under shifts in the sense that $F_G^{\leq -1}(X)$ is a subcategory of $F_G^{\leq 0}(X)$ and $F_G^{\geq 0}(X)$ is a subcategory of $F_G^{\geq -1}(X)$. However, because these identities hold $\Gamma$-locally and the shift functors are defined by shifting $\Gamma$-locally, this is a trivial observation.

We now prove that for any object $A$ that there is a distinguished triangle in $F_G(X)$ of the form
\[
\xymatrix{
	Z \ar[r] & A \ar[r] & Y \ar[r] & Z[1]
}
\]
with $Z \in F_G^{\leq 0}(X)_0$ and $Y \in F_G^{\geq 1}(X)_0$. For this, however, we define $Z := \tau^{\leq 0}A$, $Y := \tau^{\geq 1}A$, and let the morphisms be induced from the $\Gamma$-local distinguished triangles
\[
\xymatrix{
	\quot{\tau^{\leq 0}}{\Gamma}\AGamma \ar[r] & \AGamma \ar[r] & \quot{\tau^{\geq 1}}{\Gamma}\AGamma \ar[r] & \quot{\tau^{\leq 0}\AGamma[1]}{\Gamma}
}
\]
in each fibre category $F(\XGamma)$. However, because each of these triangles are distinguished we get that the triangle
\[
\xymatrix{
	\tau^{\leq 0}A \ar[r] & A \ar[r] & \tau^{\geq 1}A \ar[r] & \tau^{\leq  0}A[1]
}
\]
is distinguished in $F_G(X)$, which in turn completes the proof of the claim and hence the theorem.
\end{proof}
\begin{example}
We now give two different examples of truncated pseudofunctors on $X$ by using the pseudofunctor $F = \DbQl{-}$; in principle we can give four by using the same $t$-structures on $D^b(-)$ and proceeding mutatis mutandis.

The first truncated pseudofunctor we describe fully is the standard $t$-structure on $\DbQl{-}$ (cf.\@ Definition \ref{Defn: Standard tstructure}). For each $\Gamma \in \Sf(G)_0$ we say that
\[
F^{\leq 0}(\XGamma) := {}^{\text{stand}}\DbQl{\XGamma}^{\leq 0}
\]
and
\[
F^{\geq 0}(\XGamma) := {}^{\text{stand}}\DbQl{\XGamma}^{\geq 0}
\]
where stand denotes the standard $t$-structure on $\DbQl{-}$ and the truncations $\quot{\tau^{\leq 0}}{\Gamma}$ and $\quot{\tau^{\geq 0}}{\Gamma}$ are given from the standard $t$-structure. Then because the morphism $\of$ is smooth, the functor $\of^{\ast}$ is $t$-exact for the standard $t$-structure and so we set $F^{\leq 0}(\of) = F^{\geq 0}(\of) = \of^{\ast}$. Then the diagrams
\[
\xymatrix{
\DbQl{\XGammap} \ar[d]_{\quot{\tau^{\leq 0}}{\Gamma^{\prime}}} \ar[r]^-{\of^{\ast}} & \DbQl{\XGamma} \ar[d]^{\quot{\tau^{\leq 0}}{\Gamma}} \\
{}^{\text{stand}}\DbQl{\XGamma}^{\leq 0} \ar[r]_-{\of^{\ast}} & {}^{\text{stand}}\DbQl{\XGamma}
}
\]
\[
\xymatrix{
\DbQl{\XGammap} \ar[d]_{\quot{\tau^{\geq 0}}{\Gamma^{\prime}}} \ar[r]^-{\of^{\ast}} & \DbQl{\XGamma} \ar[d]^{\quot{\tau^{\geq 0}}{\Gamma}} \\
 {}^{\text{stand}}\DbQl{\XGammap}^{\geq 0} \ar[r]_-{\of^{\ast}} & {}^{\text{stand}}\DbQl{\XGamma}^{\geq 0}
}
\]
commute strictly, which implies that in this case the trucantion coherences are satisfied automatically. This is the standard $t$-structure on $\DbeqQl{X}$ and it is routine to see that $\DbeqQl{X}^{\heartsuit_{\text{stand}}} \simeq \Shv_G(X;\overline{\Q}_{\ell})$.\index[notation]{DbGQLStand@${}^{\text{stand}}\DbeqQl{X}^{\heartsuit}$}\index{Standard $t$-structure!On DbG@$\DbeqQl{X}$}

The second truncated pre-equivariant pseudofunctor we construct on $\DbeqQl{X}$ is induced from the perverse $t$-structure. For each $\Gamma \in \Sf(G)_0$ we define
\[
F^{\leq 0}(\XGamma) = {}^{p}\DbQl{\XGamma}^{\leq 0}
\]
and
\[
F^{\geq 0}(\XGamma) = {}^{p}\DbQl{\XGamma}^{\geq 0}.
\]
Now for each morphism $f \in \Sf(G)_1$ let $d_f \in \Z$ denote the relative (pure) dimension of $\of$. Then because $\of$ is smooth, we have from \cite[Page 108]{BBD} that $\of^{\ast}[d_f]$ is a $t$-exact functor for the pervese $t$-structure, so we set for all $f \in \Sf(G)_1$
\[
F^{\geq 0}(\of) = F^{\leq 0}(\of) = \of^{\ast}[d_f].
\]
Let each truncation functor $\quot{\tau^{\geq 0}}{\Gamma}, \quot{\tau^{\leq 0}}{\Gamma}$ be truncations for the perverse $t$-structure; we will omit the ${(-)}^{p}$ left-handed superscript to reduce notation clutter. Then a routine calculation shows that for any $f \in \Sf(G)_1$ the diagram
\[
\begin{tikzcd}
\DbQl{\XGammap} \ar[r, ""{name = LeftU}]{}{\of^{\ast}} \ar[d, swap]{}{\quot{\tau^{\leq 0}}{\Gamma^{\prime}}} &\DbQl{\XGamma} \ar[d]{}{\quot{\tau^{\leq 0}}{\Gamma}}  \\
{}^{p}\DbQl{\XGammap}^{\leq 0} \ar[r, swap, ""{name= LeftL}]{}{\of^{\ast}[d_f]} & {}^{p}\DbQl{\XGamma}^{\leq 0} \ar[from = LeftU, to = LeftL, Rightarrow, shorten <= 4pt, shorten >= 4pt] 
\end{tikzcd}
\]
\[
\begin{tikzcd}
\DbQl{\XGammap} \ar[r, ""{name = RightU}]{}{\of^{\ast}} \ar[d, swap]{}{\quot{\tau^{\geq 0}}{\Gamma^{\prime}}} & \DbQl{\XGamma} \ar[d]{}{\quot{\tau^{\geq 0}}{\Gamma}} \\
{}^{p}\DbQl{\XGammap}^{\geq 0} \ar[r, swap, ""{name = RightL}]{}{\of^{\ast}[d_f]} & {}^{p}\DbQl{\XGamma}^{\geq 0} \ar[from = RightU, to = RightL, Rightarrow, shorten <= 4pt, shorten >= 4pt]
\end{tikzcd}
\]\index[notation]{DbGQLPer@$\DbeqQl{X}^{\heartsuit_{\text{per}}}$}
commutes up to an invertible $2$-cell because of the truncations and the exactness of $\of^{\ast}[d_f]$, where the commutativity isomorphisms are induced based on how we identify and move between pullbacks. This defines the equivariant perverse $t$-structure on $\DbeqQl{X}$ and Theorem \ref{Theorem: Heart of equivariant t-structure is equivariant heart of t-structure} below shows that 
\[
(\DbeqQl{X})^{\heartsuit_{\text{per}}} = \Per_G(X).
\]
\end{example}

\begin{Theorem}\label{Theorem: Heart of equivariant t-structure is equivariant heart of t-structure}
	Let $F$ be a truncated pre-equivariant pseudofunctor on $X$ such that $F^{\geq 0}(\of) = F^{\leq 0}(\of)$ is $t$-exact for all $f \in \Sf(G)_1$. Define the Abelian pre-equivariant pseudofunctor $F^{\heartsuit}$ by setting
	\[
	F(\XGamma) := F^{\geq 0}(\XGamma) \cap F^{\leq 0}(\XGamma) = F(\XGamma)^{\heartsuit}
	\]
	and setting the fibre functors $F^{\leq 0}(\of) =: F^{\heartsuit}(\of)$. Then \index[notation]{FGXHeart@$F_G(X)^{\heartsuit}$}\index[notation]{FHeartGX@$(F^{\heartsuit})_G(X)$}
	\[
	(F^{\heartsuit})_G(X) = (F_G(X))^{\heartsuit}.
	\]
\end{Theorem}
\begin{proof}
Begin by observing that an object $A \in (F^{\heartsuit})_G(X)_0$ if and only if for all $\Gamma \in \Sf(G)_0$ the object $\AGamma \in F(\XGamma)^{\heartsuit}$ and if and only if each transition isomorphism $\tau_f^A:F(\of)\AGammap \xrightarrow{\cong} \AGamma$ is a transition isomorphism in $F(\XGamma)^{\heartsuit}$. On the other hand, by construction of the $t$-structure on $F_G(X)$ given in Theorem \ref{Thm: Section Triangle: t-structures on Equivariant cat} an object $A$ in $F_G(X)$ is in $(F_G(X))^{\heartsuit}_0$ as well if and only if $A \in (F_G^{\geq 0}(X) \cap F^{\leq 0}_G(X))_0$. However, unwinding this definition we find that $A$ lies in the heart of this $t$-structure if and only if $\Gamma$-locally we have that
\[
\AGamma \in \left(F(\XGamma)^{\leq 0} \cap F(\XGamma)^{\geq 0}\right)_0.
\]
However, this in turn implies that for all $\Gamma \in \Sf(G)_0$ we must have that $\AGamma \in F(\XGamma)^{\heartsuit}_0$, so the only aspect of morphisms that may differ is that of transition isomorphisms. However, this implies that because $\tau_f^{A}:F(\of)\AGammap \xrightarrow{\cong} \AGamma$ is a morphism in both $F(\XGamma)^{\geq 0}$ and $F(\XGamma)^{\leq 0}$ that the isomorphism $\tau_f^A$ is an isomorphism in $F(\XGamma)^{\geq 0} \cap F(\XGamma)^{\leq 0}$ and hence in the heart. As such, we find that $(F^{\heartsuit})_{G}(X)_0 = (F_G(X))^{\heartsuit}_0$. The verification that $(F^{\heartsuit})_{G}(X)_1 = (F_G(X))^{\heartsuit}_1$ follows a similar analysis as this using that $F^{\heartsuit}(\of) = F^{\geq 0}(\of) = F^{\leq 0}(\of)$ is $t$-exact, so we conclude that $F_G(X)^{\heartsuit} = (F^{\heartsuit})_G(X)$.
\end{proof}

\begin{remark}
For what follows it is worth introducing the notion of a replete subcategory $\Dscr$ of a category $\Cscr$. Such a subcategory is replete if for any $X \in \Dscr_0$ and any isomorphism $f:X \xrightarrow{\cong} Y$ in $\Cscr$, then both $Y$ and $f$ are objects and morphisms of $\Dscr$. In particualr, the subcategories $(\Dscr^{\leq 0}, \Dscr^{\geq 0})$ defining a $t$-structure on a category $\Dscr$ give two (potentially) separate examples of replete subcategories.
\end{remark}

\begin{Theorem}\label{Thm: Section Triangle: t-structures on Equivariant cat}
Let $F:\SfResl_G(X)^{\op} \to \fCat$ be a truncated pre-equivariant pseudofunctor and let $\Ascr := F(X)^{\heartsuit_t}$ be the heart of some $t$-structure on the category $F(X)$ for which there are invertible two-cells
\[
\begin{tikzcd}
F(\quot{X}{G}) \ar[r, ""{name = U}]{}{F(\psi)} \ar[d, swap]{}{\quot{\tau^{\leq 0}}{\Gamma}} & F(X) \ar[d]{}{{}^{t}\tau^{\leq 0}} \\
F(\quot{X}{G})^{\leq 0} \ar[r, swap, ""{name = L}]{}{F^{\leq 0}(\psi)} & {}^{t}F(X)^{\leq 0} \ar[from = U, to = L, Rightarrow, shorten <= 4pt, shorten >= 4pt]{}{\theta_{\psi}^{\leq 0}}
\end{tikzcd}
\]
and:
\[
\begin{tikzcd}
F(\quot{X}{G}) \ar[r, ""{name = U}]{}{F(\psi)} \ar[d, swap]{}{\quot{\tau^{\geq 0}}{\Gamma}} & F(X) \ar[d]{}{{}^{t}\tau^{\geq 0}} \\
F(\quot{X}{G})^{\geq 0} \ar[r, swap, ""{name = L}]{}{F^{\geq 0}(\psi)} & {}^{t}F(X)^{\geq 0} \ar[from = U, to = L, Rightarrow, shorten <= 4pt, shorten >= 4pt]{}{\theta_{\psi}^{\geq 0}}
\end{tikzcd}
\]
Then the subcategories generated by the objects
\[
{}^{t}F_G(X)^{\leq 0}_0 := \lbrace A \in F_G(X)_0 \; | \; F(\psi)\quot{A}{G} \in {}^{t}F(X)^{\leq 0}_0 \rbrace
\]
and
\[
{}^{t}F_G(X)^{\geq 0}_0 := \lbrace A \in F_G(X)_0 \; | \; F(\psi)\quot{A}{G} \in {}^{t}F(X)^{\leq 0}_0 \rbrace
\]
determine a $t$-structure on $F_G(X)$.
\end{Theorem}
\begin{proof}
We must check first that ${}^{t}F_G(X)^{\leq 0}$ and ${}^{t}F_G(X)^{\geq 0}$ are both strictly full subcategories of $F_G(X)$. For this we note that if we can prove this for ${}^{t}F_G(X)^{\geq 0}$ then the case for ${}^{t}F_G(X)^{\leq 0}$ follows mutatis mutandis, so it suffices to only prove it for ${}^{t}F_G(X)^{\geq 0}$. To this end we first show the fullness of ${}^{t}F_G(X)^{\geq 0}$. For this let $A, B \in ({}^{t}F_G(X)^{\geq 0})_0$ and assume that $P \in F_G(X)(A,B)$. Now since $P:A \to B$ is a morphism  $F(\psi)\quot{A}{G}, F(\psi)\quot{B}{G} \in {}^{t}F(X)^{\geq 0}_0$ and the category ${}^{t}F(X)^{\geq 0}$ is full, it follows that $F(\psi)\quot{\rho}{G} \in {}^{t}F(X)^{\geq 0}_1$ and hence $P \in {}^{t}F_G(X)^{\geq 0}_1$ as well. This shows that ${}^{t}F_G(X)^{\geq 0}$ is full.

To see that ${}^{t}F_G(X)^{\geq 0}$ is closed with respect to isomorphism, let $A$ be an object in ${}^{t}F_G(X)^{\geq 0}$ and let $B \in F_G(X)_0$ with $A \cong B$ through the morphism $P$. Now note that this implies that $F(\psi)\quot{A}{G} \cong F(\psi)\quot{B}{G}$ in $F(X)$. Because $F(\psi)\quot{A}{G}$ is an object in ${}^{t}F(X)^{\geq 0}$ and the category ${}^{t}F(X)^{\geq 0}$ is replete, it follows that $F(\psi)\quot{\rho}{G}$ is in ${}^{t}F(X)^{\geq 0}$ as well. Thus $B \in {}^{t}F_G(X)^{\geq 0}$ as well.

We now show that the categories ${}^{t}F_G(X)^{\geq 0}$ and ${}^{t}F_G(X)^{\leq 0}$ are appropriately translation invariant, i.e., that ${}^{t}F_G(X)^{\leq -1}$ is contained in ${}^{t}F_G(X)^{\leq 0}$ and ${}^{t}F_G(X)^{\geq 0}$ is contained in ${}^{t}F_G(X)^{\geq -1}$. However, both of these results are immediate from the structure of ${}^{t}F(X)^{\geq 0}$ and ${}^{t}F(X)^{\leq 0}$ together with the fact that the functor $F(\psi)$ commutes up to isomorphism with the shift functors $[n]$.

We now prove that if $A$ is an object of ${}^{t}F_G(X)^{\leq 0}$ and if $B$ is an object of ${}^{t}F_G(X)^{\geq 1}$, then $F_G(X)(A,B) = 0$. For this we note that first by the fact that $F(\psi)$ and $F(\psi^{-1})$ are inverse equivalences and hence both fully faithful that
\begin{align*}
F(\quot{X}{G})(A,B) &\cong F(\quot{X}{G})(F(\psi^{-1}F(\psi)\quot{A}{G},F(\psi^{-1})F(\psi)\quot{B}{G}))  \\
&\cong F(X)(F(\psi)\quot{A}{G},F(\psi)\quot{B}{G}) \cong 0
\end{align*}
where the last isomorphism holds because $F(\psi)\quot{B}{G} \in {}^{t}F(X)^{\geq 1}_0$ and $F(\psi)\quot{A}{G} \in F(X)^{\leq 0}_0.$ Now let $\Gamma \in \Sf(G)_0$ with an $\Sf(G)$-morphism $f:\Gamma \to G$. Then from the triangulated nature of the functors $F(\of)$ and the descent data these objects satisfy we get that
\[
F(\XGamma)(F(\of)\quot{A}{G},F(\of)\quot{B}{G}) \cong F(\quot{X}{G})(\quot{A}{G},\quot{B}{G}) \cong 0
\]
so it follows that
\[
F(\XGamma)(\AGamma,\BGamma) \cong F(\XGamma)(F(\of)\quot{A}{G},F(\of)\quot{B}{G}) \cong 0.
\]
Finally, because this class of varieties is sufficiently cofinal in $\Sf(G)_0$ (and in fact contains every variety $G^n$ for $n \geq 1$ by the fact that the projections  and multiplications $G^{n+1} \to G^{n}$ are all smooth $G$-equivariant maps of constant fibre dimension) by equivariant descent we get that
\[
F_G(X)(A,B) \cong 0,
\]
as was desired.

Finally, we verify that there is a distinguished triangle
\[
\xymatrix{
Z \ar[r] & A \ar[r] & Y \ar[r] & Z[1]
}
\]
for any $A \in F_G(X)_0$ with $Z \in {}^{t}F_G(X)^{\leq 0}_0$ and $Y \in {}^{t}F_G(X)^{\geq 1}_0$. For this we note that if we can show that, for the truncation functors $\tau^{\leq 0}$ and $\tau^{\geq 1}$, that $\tau^{\leq 0}A \in {}^{t}F_G(X)^{\leq 0}_0$ and $\tau^{\geq 1}A \in {}^{t}F_G(X)^{\geq 1}_0$ then we will be done by the proof of Theorem \ref{Thm: Section Triangle: t-structures on Equivariant cat}. For this we calculate that
\begin{align*}
F(\psi)\left(\quot{\left(\tau^{\leq 0}A\right)}{G}\right) & = F(\psi)\left(\quot{\tau^{\leq 0}}{G}\quot{A}{G}\right) \cong {}^{t}\tau^{\leq 0}\left(F(\psi)\quot{A}{G}\right) \in {}^{t}F(X)^{\leq 0}_0
\end{align*}
and similarly we have that 
\[
F(\psi)\left(\quot{\tau^{\geq 1}}{\Gamma}\right) \cong {}^{t}\tau^{\geq 1}\left(F(\psi)\quot{A}{G}\right) \in {}^{t}F(X)^{\geq 1}_0.
\]
Thus we conclude that the our distinguished triangle of the desired form is given by
\[
\xymatrix{
\tau^{\leq 0}A \ar[r] & A \ar[r] & \tau^{\geq 1}A \ar[r] & \tau^{\leq 0}A[1]
}
\]
which in turn proves the theorem.
\end{proof}

\part{The Four Flavors of Equivariant Derived Category}\label{Chapter: EDC Comp}
\chapter{Introduction}\label{Chapter Intro Comparison}
When studying a group object $G$ in a category $\Cscr$, we have seen historically that studying the objects which carry a $G$-action and the equivariant $G$-morphisms between such objects is a fruitful way to learn about the group $G$ itself. For instance, this approach gives rise to the myriad and important equivalences and situations described below:
\begin{itemize}
	\item For a triply connected topological space $X$, the equivalence
	\[
	\Loc(X;\C) \simeq \pi_1(X,x)\mathbf{Rep}
	\]\index[notation]{Pi1Xx@$\pi_1(X,x)\mathbf{Rep}$}
	of complex local systems on $X$ with complex $\pi_1(X,x)$-representations.
	\item For a group $G$ in $\Set$, the topos $G-\Set$ of left $G$-sets and $G$-equivariant morphisms. This gives a natural location to perform group cohomology, as the topos cohomology of $G\Set$ is group cohomology $H^1(G,-)$.
	\item If $W_K$ is the Weil-Deligne group of a local field $K$, the realization of $W_K$-representations in terms of the local class field theory of $K$ is foundational for the Langlands Programme.
	\item The theory of complex representations of finite and compact groups $G$ is important in the theory of harmonic analysis and in Lie theory.
\end{itemize}
In the examples above, we have been able to avoid describing how it is that actions of groups on objects, say a $G$ action on a space $X$, gives rise to an equivariance of sheaves on $X$. While it is straightforward to say when a function $f:X \to Y$ between spaces $X$ and $Y$ that both carry $G$-actions is equivariant (we simply say that $f(gx) = gf(x)$ for all $x \in X, g \in G$), it is not as straightforward to say when a sheaf $\Fscr$ is equivariant with respect to the $G$-action on $X$. This  problem, however, was eventually solved by Mumford and Fogarty in \cite{GIT}: The solution is to say that a sheaf $\Fscr$ is $G$-equivariant on $X$ if there is an isomorphism of sheaves
\[
\theta:\pi_2^{\ast}\Fscr \xrightarrow{\cong} \alpha_X^{\ast}\Fscr
\]
on the space $G \times X$ with projection map $\pi_2:G \times X \to X$ and action map $\alpha_X:G \times X \to X$ such that this isomorphism satisfies the cocycle condition
\[
\big(\alpha_X \circ (\id_G\times \alpha_X)\big)^{\ast}\theta \circ (\alpha_X \circ (\pi_{23})^{\ast})\theta = (\mu \times \id_X)^{\ast}\theta.
\] 
In this paper we call the map $\theta$ an equivariance of $\Fscr$.\index{Equivariance! Of a Sheaf $\Fscr$} That this is the right notion comes from the fact that on every stalk of the form $(g,x)$ in $G \times X$ we end up with isomorphisms
\[
\Fscr_x = \Fscr_{\pi_2(g,x)} = \pi_2^{\ast}\Fscr_{(g,x)} \cong \alpha_X^{\ast}\Fscr_{(g,x)} = \Fscr_{\alpha_X(g,x)} = \Fscr_{gx}.
\]
This same condition is straightforward to extend to other various Abelian categories of (complexes) of sheaves and to scheme-theoretic settings: For a variety $X$ carrying an action of an algebraic group $G$, an object $A$ of $\Cscr(X)$ is equivariant if and only if there is an equivariance
\[
\theta:\pi_2^{\ast}A \xrightarrow{\cong} \alpha_X^{\ast}A
\]
between objects in $\Cscr(G \times X)$ where $\Cscr$ is any of the following categorical constructions:
\begin{itemize}
	\item $\Per(X)$ and $\Per(G \times X)$, the categories of perverse sheaves on $X$ and $G \times X$;
	\item $\Loc(X)$ and $\Loc(G \times X)$, the categories of local systems on $X$ and $G \times X$;
	\item $\Shv(X;\overline{\Q}_{\ell})$ and $\Shv(G \times X;\overline{\Q}_{\ell})$, the categories of $\ell$-adic sheaves on $X$ and $G \times X$;
	\item $\Per(X;\overline{\Q}_{\ell})$ and $\Per(G \times X;\overline{\Q}_{\ell})$, the categories of $\ell$-adic perverse sheaves on $X$ and $G \times X$;
	\item $\Loc(X;\overline{\Q}_{\ell})$ and $\Loc(G \times X;\overline{\Q}_{\ell})$, the categories of $\ell$-adic local systems on $X$ and $G \times X$.
\end{itemize}
However, when trying to extend this to the full level of $D^b_c(X)$ or $\DbQl{X}$ through the na{\"i}ve approach of asking for such an isomorphism, things fail quickly as the category of such ``na{\"i}vely equivariant''\index{Naive@Na{\"i}ve Equivariance} complexes (complexes with an equivariance $\theta:\pi_2^{\ast}A \to \alpha_X^{\ast}A$) need not even be triangulated. Following an example spelled out Torsten Ekedahl in a MathOverflow answer \cite{Torsten}, let us explain what happens: Let $G$ be a group, let $n \in \N$, and let $p \in \Z$ be a nonzero prime. Find an additive $G$-action on the group $(\Z/p\Z)^n$ which does not lift to an additive action on $(\Z/p^2\Z)^n$ and find the Bockstein morphism $(\Z/p\Z)^n \to (\Z/p\Z)^n[1]$ in $D^b(\ast)$ which has mapping cone $(\Z/p^2\Z)^n$. Then this morphism gives a $G$-equivariant morphism in $D^b(\ast)$ which has no compatible $G$-action on the mapping cone, and hence cannot fit into a distinguished triangle in $D_G(\ast)$. In particular, while there is an equivariance for $\pi_2^{\ast}(\Z/p\Z)^n$ and $\alpha^{\ast}(\Z/p\Z)^n$, there is no mapping cone for this action in $D^b(\ast)$ so the category $D_G^{\text{na{\"i}ve}}(X)$ of na{\"i}vely equivariant complexes cannot be triangulated.

There have been various approaches to solving this ``lack of triangulation'' problem in the literature by various different authors. Perhaps the most well-known approach was that of Bernstein-Lunts in the 90's (cf.\@ \cite{BernLun}) where the equivariant bounded derived category was defined essentially to be given as effective descent data taken over certain acyclic resolutions of the $G$-action on $X$ by free $G$-spaces that sit over $X$. This approach was later extended to the scheme theoretic case, at least in principle, by Pramod Achar in \cite{PramodBook}. This extension restricts the resolutions, however, from all free $G$-spaces to instead acyclic $G$-resolutions of the scheme/variety $X$. This approach to the equivariant bound derived category has seen use in \cite{CFMMX}, \cite{Juteauetal}, \cite{ladic}, \cite{LaszloOlsson}, \cite{MirkVilDuality}, and \cite{Schnurer}, for instance.

An alternative to the Bernstein-Lunts approach to the equivariant derived category is that of the simplicial equivariant derived category. This uses a simplicial space $\underline{G \backslash X}_{\bullet}$ to approximate the quotient space $G \backslash X$ and takes a certain category of simplicial sheaves on this simplicial space to provide the setting for an equivariant derived category. This was the category used, for instance, in Deligne's paper \cite{DeligneHodge3} and is known only in the topological case to be equivalent to the bounded derived category of Bernstein-Lunts.

A very formally similar, but alternate, approach to simplicial equivariant derived category is the equivariant derived category of $\ell$-adic sheaves on the quotient stack $[G \backslash X]$ of Behrend introduced in \cite{Behrend} and seen later in \cite{LasOls}, \cite{ArtinStacks2}, and \cite{LaszloOlsson}. This category uses a pair of fibred Grothendieck toposes over the category $\mathbf{\Delta}^{\op}$ to produce constructible derived $\ell$-adic sheaves overtop an algebraic stack $\Xfrak$. By setting $\Xfrak = [G \backslash X]$ we obtain another competing perspective on the equivariant derived category which is formally similar but technically distinct from the others.

The fourth and final main approach to the equivariant derived category we consider in this paper is the one pioneered by Lusztig in \cite{LusztigCuspidal2}. The equivariant derived category of Lusztig is a given as a construction which takes a $G$-scheme $X$ and resolves the action of $G$ on $X$ by a certain class of smooth free principal $G$-varieties which nicely control the singularities of the action and then takes descent data over these resolutions in order to produce the equivariance that equivaraint derived complexes need to have. This approach to equivariance is the one primarily used in \cite{CFMMX} to define equivariant vanishing and nearby cycles. The equivariant derived category of Lusztig is a formally powerful construction (and may be readily generalized to equivariant categories indexed by pseudofunctors as in Section \ref{Section: Equivariant Categories} above), but it also is intimately related to equivariant cohomology and homology of \cite{LusztigCuspidal1, LusztigCuspidal2} and is used in \cite{LusztigCuspidal2} to construct the graded affine Hecke algebra. 

Four distinct formulations of the equivariant derived category is, frankly, too many different formalisms to carry around at all moments. When studying equivariant derived complexes, one needs to keep in mind what exact situation they are in, the tools they have at hand, and how to build them. However, it turns out that worrying about this is unnecessary: In this inquiry, we establish that all four equivariant derived categories are equivalent to each other when $G$ is a smooth aglebraic group and $X$ is a variety. In fact this is the main result of Part \ref{Chapter: EDC Comp} present below; cf.\@ Theorem \ref{Theorem: Section 4.4: Main Theorem}).
\begin{Theorem}[{cf.\@ Theorem \ref{Theorem: Section 4.4: Main Theorem}}]
	Let $G$ be a smooth affine algebraic group and let $X$ be a variety. Then if $\DbeqQl{X}$ is the equivariant derived category of Lusztig, $\quot{\DbeqQl{X}}{ABL}$ is the Achar-Bernstein-Lunts equivariant derived category, $\Dbeqsimp{X}$ is the simplical equivariant derived category, and $\DbQl{[G \backslash X]}$ is the stacky equivariant derived category there is a four-way equivalence
	\[
	\DbeqQl{X} \simeq \quot{\DbeqQl{X}}{ABL} \simeq \Dbeqsimp{X} \simeq \DbQl{[G \backslash X]}.
	\]
	Furthermore, these equivalences are all compatible with the standard $t$-structure on each triangulated category and give a mutual equivalence of the heart of each $t$-structure with the category $\Shv_G^{\text{na{\"i}ve}}(X;\overline{\Q}_{\ell})$.
\end{Theorem}

What this essentially means that all one needs to do when studying problems involving the equivariant derived category is keep in mind what situation they are in, what tools they have available, and what their ultimate goal is. For instance, if one needs to compute something explicitly, the simplicial model is the most straightforward setting; if one wants to make descent theoretic arguments, the ABL-model is quite helpful by building effectiveness into the story innately; if one is in a stack-theoretic setting, the category $\Dbeqstack{X}$ is quite helpful to have at hand; and if one is studying the interaction between various other equivariant categories and needs to know the general structure of $\DbeqQl{X}$, or if one is in a setting adjacent to the graded affine Hecke algebra, the category of Lusztig is ideal. 

\section{Summary of Main Results}

Let us now proceed to explain some of the key results of Part \ref{Chapter: EDC Comp}. The first main result we present is that the two different perspectives on a resolution-based approach to the equivariant derived category, those of Lusztig and Achar-Bernstein-Lunts, are equivalent by discussing acyclic descent:
\begin{Theorem}[{cf.\@ Corollary \ref{Cor: Section 4: EDCs are equiv}}]
	For any smooth affine algebraic group $G$ and left $G$-variety $X$, there is an equivalence of categories:
	\[
	\begin{tikzcd}
	\DbeqQl{X} \ar[r, bend left = 30, ""{name = U}] & \quot{\DbeqQl{X}}{ABL} \ar[l, bend left = 30, ""{name = L}] \ar[from = U, to = L, symbol = \simeq]
	\end{tikzcd}
	\]
	Furthermore, this equivalence of categories is compatible with the standard and pervese $t$-structures on each category.
\end{Theorem}

Once we have shown that the two categories are equivalent, we prove that the equivariant derived category $\DbeqQl{X}$ satisfies Desiderata \ref{Desiderata}. Before we do this, we fill in the remaining gaps for what we require and culminate in the theorem below.
\begin{Theorem}[cf.\@ Theorem \ref{Theorem: EDC sats desiderata}]
Let $G$ be a smooth algebraic group and let $X$ be a left $G$-variety. Then $\DbeqQl{X}$ satisfies Desiderata \ref{Desiderata}.
\end{Theorem}

After showing that $\DbeqQl{X} \simeq \quot{\DbeqQl{X}}{ABL}$ and that $\DbeqQl{X}$ satisfies Desiderata \ref{Desiderata}, we move on to compare Lusztig's equivariant derived category $\DbeqQl{X}$ with the simplicial category $\Dbeqsimp{X}$. In the process of performing this comparison we provide a ``sanity check'' theorem that states that when $G$ is a smooth group scheme and $U$ is a $G$-scheme that admits a $G$-quotient $G \backslash U$ for which the quotient morphism $\quo_U:U \to G\backslash U$ is smooth, then there is an equivalence of categories between the normal bounded derived category on $G \backslash U$ and the simplicial equivariant category on $\underline{G \backslash U}_{\bullet}$. We present this as the theorem below for detail:
\begin{Theorem}[{cf.\@ Proposition \ref{Prop: Section 4: Quot derived cat is equiv to simplicial derived cat}}]\label{Thm: Statment of results for Prop 5.5.19}
	Let $G$ be a smooth group scheme over a field $\Spec K$ and let $U$ be a left $G$-scheme which admits a $G$-quotient for which the morphism $\quo_{G}:U \to G \backslash U$ is smooth. Then there is an equivalence of categories
	\[
	\Dbeqsimp{U} \simeq D_c^b(G \backslash U; \overline{\Q}_{\ell}).
	\]
\end{Theorem}
As a straightforward consequence to this equivalence, we can provide a very concrete description of what the pullback against the quotient morphism is doing on the categorical level, i.e., we can give a very concrete factorization of $\quo_{U}^{\ast}$. While we only perform this for the case in which $U = G^n \times X$ for a smooth affine algebraic group $G$, $n \geq 1$, and a left $G$-variety $X$, it is enough of a de-mystification of what the quotient pullback is doing that it is, to me at least, a key lemma of Part \ref{Chapter: EDC Comp}.
\begin{lemma}[{cf.\@ Lemma \ref{Lemma: Section 4: Factorization of RGn}}]
	For any smooth affine algebraic group $G$, for any $n \in \N$ with $n \geq 1$, and for any left $G$-variety $X$, if $L_{G^n}$ is the equivalence of Theorem \ref{Thm: Statment of results for Prop 5.5.19} below
	\[
	L_{G^n}:\DbQl{G \backslash (G^n \times X)} \xrightarrow{\simeq} \Dbeqsimp{(G^n \times X)}
	\]
	then there is a factorization of functors
	\[
	\xymatrix{
		\DbQl{G \backslash (G^n \times X)} \ar[drr]_{\quo_{G^n \times X}^{\ast}} \ar[rr]^-{L_{G^n}} & & \Dbeqsimp{(G^n \times X)} \ar[d]^{\tr_0} \\
		& & \DbQl{G^n \times X}
	}
	\]
	where $\tr_0$ is the degree $0$ truncation of a sequence of complexes of simplicial sheaves on the simplicial scheme $\underline{G \backslash (G^n \times X)}_{\bullet}$.
\end{lemma}

After using the above lemma and theorem to perform very explicit computations and translations between $\DbeqQl{X}$ and $\Dbeqsimp{X}$, we complete the next section of Part \ref{Chapter: EDC Comp} by concluding finally that the two categories are equivalent.

\begin{Theorem}
	For any smooth affine algebraic group $G$ and any left $G$-variety $X$, there is an equivalence of categories:
	\[
	\begin{tikzcd}
	\DbeqQl{X} \ar[r, bend left = 30, ""{name = U}] & \Dbeqsimp{X} \ar[l, bend left = 30, ""{name = L}] \ar[from = U, to = L, symbol = \simeq]
	\end{tikzcd}
	\]
\end{Theorem}

After showing that the equivalence $\DbQl{X} \simeq \Dbeqsimp{X}$, we move on to compare the simplicial equivariant derived category $\Dbeqsimp{X}$ with the stacky equivariant derived category of Behrend $\DbQl{[G \backslash X]}$ constructed in \cite{Behrend}. Our strategy towards proving this equivalence is to prove that the simplicial scheme $\underline{G \backslash X}_{\bullet}$ is isomorphic in the category $[\mathbf{\Delta}^{\op},\Sch]$ of simplicial schemes to the simplicial presentation of the stack $\underline{[G \backslash X]}_{\bullet}$ induced from the fact that the quotient stack $[G \backslash X]$ is an algebraic stack with a presentation as an fppf quotient of a internal groupoid in the category $\Var_{/\Spec K}$ of varieties over a field $K$. In doing this we prove the next main proposition, which states the isomorphism of simplicial schemes (more accurately, simplicial varieties).
\begin{Theorem}[{cf.\@ Proposition \ref{Prop: Section 4: Two simplicial presentations are isomorphic}}]
	Let $G$ be a smooth affine algebraic group and let $X$ be a left $G$-variety. Then there is an isomorphism of simplicial varieties
	\[
	\underline{G \backslash X}_{\bullet} \cong \underline{[G \backslash X]}_{\bullet}.
	\]
\end{Theorem}

Proving the theorem above requires some introduction of coskeleta of simplicial schemes and the truncation functors of simplicial schemes. However, an important consequence of the theorem above is that this isomorphism lifts to an equivalence of simplicial equivariant derived categories, which provides us with the final tool needed to prove the simplicial and stacky equivalence.

\begin{corollary}[{cf.\@ Corollary \ref{Cor: Section 4: Bounded derived simplicial cat is equiv}}]
	There is an equivalence of categories
	\[
	\Dbeqsimp{X} \simeq \Dbeqstack{X}
	\]
	that can be lifted to an isomorphism by choosing compatible pullbacks appropriately.
\end{corollary}

After providing this corollary and reviewing the construction of the category $\DbQl{[G \backslash X]}$ of \cite{Behrend}, we finally prove that $\DbQl{[G \backslash X]}$ is equivalent to $\Dbeqsimp{X}$ are equivalent.
\begin{Theorem}[{cf.\@ Theorem \ref{Theorem: Section 4: Equivalence of Behrend's Derived Cat and Simplicial Derived Cat}}]
	For any smooth affine algebraic group $G$ and left $G$-variety $X$, there is an equivalence of categories
	\[
	\DbeqQl{[G \backslash X]} \simeq \Dbeqsimp{X}.
	\]
\end{Theorem}

Finally, these main results together allow us to conclude with the main result of Part \ref{Chapter: EDC Comp}. While we have presented it already in this introductory setting, we restate the theorem here for the sake of having a complete record of the main results of Part \ref{Chapter: EDC Comp}.
\begin{Theorem}[{cf.\@ Theorem \ref{Theorem: Section 4.4: Main Theorem}}]
	Let $G$ be a smooth affine algebraic group and let $X$ be a variety. Then if $\DbeqQl{X}$ is the equivariant derived category of Lusztig, $\quot{\DbeqQl{X}}{ABL}$ is the Achar-Bernstein-Lunts equivariant derived category, $\Dbeqsimp{X}$ is the simplical equivariant derived category, and $\DbQl{[G \backslash X]}$ is the stacky equivariant derived category there is a four-way equivalence
	\[
	\DbeqQl{X} \simeq \quot{\DbeqQl{X}}{ABL} \simeq \Dbeqsimp{X} \simeq \DbQl{[G \backslash X]}.
	\]
\end{Theorem}

\section{Structure of the Paper}
The paper is laid out with the following structure in mind:
\begin{itemize}
\item Chapter \ref{Section: Lusztig EDC is equiva to ABL EDC} proves the equivalence of the category $\DbeqQl{X}$ with the category $\ABLDbeqQl{X}$ and is largely devoted to the descent-theoretic techniques between these objects. This chapter largely involves more high level techniques; explicit versions of the results laid out in here are often given in the simplicial case for two main reasons: First, the simplicial case benefits much more than the general descent theoretic case from the more explicit detail, as it is the situation in which computations are easiest to perform; second, the proofs in either case largely follow one from the other, so we have deferred the main technical detail to where it is more beneficial to be seen in order to reduce the already significant length of Part \ref{Chapter: EDC Comp}. 
\item Chapter \ref{Section: Lusztig EDC is the Simplicial EDC} proves the equivalence of categories between the category $\DbeqQl{X}$ and $\Dbeqsimp{X}$. While this chapter is quite long, it is explicit because the simplicial equivariant derived category is well-suited to computations and I felt it is valuable to have illustrations of how to manipulate and calculate these objects explicitly and carefully. As a result of this approach to the equivalence, we have provided various results like showing how there are transition isomorphisms for degeneracy maps $G^n \times X \to G^{n+1} \times X$ when $A \in \DbeqQl{X}_0$, as well as more technical results of that type. This chapter is largely independent, but uses some results from Chapter \ref{Section: Lusztig EDC is equiva to ABL EDC}.
\item Chapter \ref{Section: Simplicial EDC is the stcky EDC} proves that the simplicial equivariant derived category is equivalent to the stacky equivariant derived category. Few results from prior sections are used here, but significant stack-theoretic background and Grothendieck topos theoretic background (in the style of \cite{SGA4}) is required in this section. The highlights of the theory that we use are spelled out in Appendices \ref{Appendix: Stackification via Torsors}, \ref{Appendix: Locally closed subtoposes}, and \ref{Appendix: Stratifications and c-structures} below.
\end{itemize}

\chapter{The Achar-Bernstein-Lunts Equivariant Derived Category}\label{Section: Lusztig EDC is equiva to ABL EDC}

\section{The Equivalence $\DbeqQl{X} \simeq \ABLDbeqQl{X}$}
Let us now recall the definition of an $n$-acyclic resolution of $X$, where $n \in \N$. Essentially this is a variety $U$ over $X$ which resolves any $n$-or-smaller dimensional singularities of $X$ in the sense that any perverse cohomology of a complex $\Fscr$ on $X$, when pulled back to $U$ should arise only from the nature of $\Fscr$ and $X$, and not be contributed to in any way by $U$ within those dimensions. Formally, we capture this in the following definition.
\begin{definition}[{\cite[Definition 6.2.4]{PramodBook}}]\index{nAcyclicMor@$n$-Acyclic Morphism}
	A smooth morphism $f:U \to X$ is $n$-acyclic if for any smooth map $g:V \to X$ the base change
	\[
	\xymatrix{
		U \times_X V \ar[r]^-{h} \ar[d] \pullbackcorner & V \ar[d]^{g}\\
		U \ar[r]_-{f} & X
	}
	\]
	satisfies, for any perverse sheaf $\Fscr$ on $V$, that the natural map
	\[
	\Fscr \to {}^{p}\tau^{\leq n}\left(\per Rh_{\ast}\left(\per h^{\ast}\Fscr\right)\right)
	\]
	is an isomorphism. A smooth morphism $f \in \GVar(U,X)$ is an $n$-acyclic resolution of $X$ if $U$ is a principal $G$-variety and $f$ is $n$-acyclic.
\end{definition}
\begin{remark}
We say that a morphism $f:Y \to X$ is $\infty$-acyclic if it is $n$-acyclic for all $n \in \N$.
\end{remark}
\begin{remark}
In the definition above, the truncation functor $\per\tau^{\leq n}$ is the truncation functor ${}^{t}\tau^{\leq n}$ for the perverese $t$-structure on $\DbQl{V}$ (cf.\@ \cite[D{\'e}finition 2.1.2]{BBD}) while $\per h^{\ast}$ is the perverse-exact pullback functor $\per h^{\ast} = h^{\ast}[-\dim V]$ and $\per Rh_{\ast}$ is the perverse pushforward functor.
\end{remark}

\begin{lemma}\label{Lemma: Section 5.5: Acyclic Base Change}
	Let $f:Y \to X$ be an $n$-acyclic morphism  and consider the base change diagram of schemes below:
	\[
	\xymatrix{
		Y \times_X Z \ar[r]^-{p_1} \ar[d]_{p_2} & Y \ar[d]^{f} \\
		Z \ar[r]_-{g} & X
	}
	\]
	Then the morphism $p_2$ is $n$-acyclic as well.
\end{lemma}

A convenient characterization of a class of $n$-acyclic morphisms and $n$-acyclic varieties is given below. We can give these as a class of varieties following our intuition of those which do not contribute any cohomology in dimensions $1 \leq q \leq n$.
\begin{proposition}[{\cite[Lemma 6.2.5]{PramodBook}}]\label{Prop: Cohomology Prop}
	Let $k$ be a field and let $X$ be a smooth connected $k$-variety such that $H^q(X;\Z_{\ell}) = 0$ for $1 \leq q \leq n$ and for which $H^{n+1}(X;\Z_{\ell})$ is a free $\Z_{\ell}$-module. Then the map $X \to \Spec k$ is $n$-acyclic.
\end{proposition}
\begin{proof}[Sketch]
We present here a sketch of the argument presented in \cite{AcharEtAl}. In particular, we will use the following fact: Given $m \in \Z$, if $A \in D^b(\Z_{\ell}\mathbf{Mod}_{\text{f.g.}})_0^{\geq m}$ with $H^{m+1}(A)$ a free $\Z_{\ell}$-module, then for any variety $Y$ over $K$ and for any $\Pscr \in \Per(Y;\overline{\Q}_{\ell})_0$,
\[
\nu_Y^{\ast}\left(\overline{\Q}_{\ell} \overset{L}{\underset{\Z_{\ell}}{\otimes}} A \right) \overset{L}{\otimes} \Pscr \in {}^{p}D^b_c(Y;\overline{\Q}_{\ell})^{\geq m}.
\]
where $\nu_{Y}:Y \to \Spec K$ is the structure map and we regard $\overline{\Q}_{\ell}$ as a sheaf of (non-finitely generated) $\Z_{\ell}$-modules on $\Spec K$.

Write $\Gamma$ for the {\'e}tale global sections functor. Observe that to prove that $\nu_X:X \to \Spec K$ is $n$-acylcic, we must consider a pullback of the form
\[
\xymatrix{
	X \times_{\Spec K} Y \ar[r]^-{\pi_2} \ar[d]_{\pi_1} \pullbackcorner & Y \ar[d]^{\nu_Y} \\
	X \ar[r]_-{\nu_X} & \Spec K
}
\]
with $\nu_Y$ smooth. However, because this pullback is the is the product of varieties we find that for any $\Pscr \in \Per(Y;\overline{\Q}_{\ell})_0$,
\begin{align*}
\left({}^{p}R(\pi_2)_{\ast} \circ {}^{p}(\pi_2)^{\ast}\right)(\Pscr) &\cong \left(R(\pi_2)_{\ast} \circ \pi_2^{\ast}\right)(\Pscr) \cong (\nu_X \times \id_Y)^{\ast}\left(\mathbbm{1}_X \boxtimes \Pscr\right) \\
&\cong \nu_Y^{\ast}(R\Gamma(\mathbbm{1}_X)) \overset{L}{\otimes} \Pscr \cong \nu_Y^{\ast}\left(\overline{\Q}_{\ell} \overset{L}{\underset{\Z_{\ell}}{\otimes}} R\Gamma(\mathbbm{1}^{\Z_{\ell}}_{X})\right) \overset{L}{\otimes} \Pscr.
\end{align*}
where$\mathbbm{1}_{X}$ is the constant $\overline{\Q}_{\ell}$-sheaf on $X$, $\mathbbm{1}^{\Z_{\ell}}_X$ is the constant $\ell$-adic sheaf on $X$, and $\boxtimes$ is the external tensor functor. Now, because the variety $X$ is connected we find that
\[
R^0\Gamma(\mathbbm{1}^{\Z_{\ell}}_{X}) \cong H^0(X;\Z_{\ell}) \cong \Z_{\ell}
\]
while our assumptions on $X$ give that for $1 \leq k \leq n$ we have
\[
R^k\Gamma(\mathbbm{1}_X^{\Z_{\ell}}) \cong H^k(X;\Z_{\ell}) \cong 0.
\]
Thus we find that $\tau^{\leq n}(R\Gamma(\mathbbm{1}_X^{\Z_{\ell}})) \cong \Z_{\ell}[0]$ (the complex with group $\Z_{\ell}$ in degree zero and $0$ in all other degrees). This gives rise to the distinguished triangle
\[
\xymatrix{
	\Z_{\ell}[0] \ar[r] & R\Gamma(\mathbbm{1}^{\Z_{\ell}}_X) \ar[r] & \tau^{\geq n+1} R\Gamma(\mathbbm{1}_X^{\Z_{\ell}}) \ar[r] & \Z_{\ell}[1]
}
\]
in $D^b(\Z_{\ell}\mathbf{Mod}_{\textbf{f.g.}}).$ Passing through the equivalence of categories \\ $D^b(\Z_{\ell}\mathbf{Mod}_{\text{f.g.}}) \simeq D^b_c(\Spec K;\Z_{\ell})$ and applying the functor
\[
\nu_{Y}^{\ast}\left(\overline{\Q}_{\ell} \overset{L}{\underset{\Z_{\ell}}{\otimes}} (-)\right) \overset{L}{\otimes} \Pscr
\] 
to our above distinguished triangle gives the distinguished triangle
\[
\xymatrix{
	\Pscr \ar[r] & ({}^{p}R(\pi_2)_{\ast} \circ {}^{p}(\pi_2)^{\ast})(\Pscr) \ar[r] & \nu_Y^{\ast}\left(\overline{\Q}_{\ell} \overset{L}{\underset{\Z_{\ell}}{\otimes}} \tau^{\geq n+ 1}(R\Gamma(\mathbbm{1}_X^{\Z_{\ell}}))\right) \overset{L}{\otimes} \Pscr \ar[d] \\
	& & \Pscr[1]
}
\]
in $\DbQl{Y}$. However, since
\[
R^{n+1}\Gamma(\mathbbm{1}_X^{\Z_{\ell}}) \cong H^{n+1}(X;\Z_{\ell})
\]
is a free $\Z_{\ell}$-module, our stated fact at the start of the sketch shows us that
\[
\nu_{Y}^{\ast}\left(\overline{\Q}_{\ell} \overset{L}{\underset{\Z_{\ell}}{\otimes}} \tau^{\geq n+ 1}(R\Gamma(\mathbbm{1}_X^{\Z_{\ell}}))\right) \overset{L}{\otimes} \Pscr \in {}^{p}\DbQl{Y}^{\geq n + 1}_0.
\]
This implies that the distinguished triangle
\[
\xymatrix{
	\Pscr \ar[r] & ({}^{p}R(\pi_2)_{\ast} \circ {}^{p}(\pi_2)^{\ast})(\Pscr) \ar[r] & \nu_Y^{\ast}\left(\overline{\Q}_{\ell} \overset{L}{\underset{\Z_{\ell}}{\otimes}} \tau^{\geq n+ 1}(R\Gamma(\mathbbm{1}_X^{\Z_{\ell}}))\right) \overset{L}{\otimes} \Pscr \ar[d] \\
	& & \Pscr[1]
}
\]
is isomorphic to a truncation distinguished triangle, which in turn implies
\[
\Pscr \cong {}^{p}\tau^{\geq n + 1}({}^{p}R(\pi_2)_{\ast} \circ {}^{p}(\pi_2))\Pscr,
\]
showing that $\nu_X$ is $n$-acyclic.
\end{proof}

We now need a result, which is in both \cite{LusztigCuspidal1} and \cite{PramodBook}, which says that for all smooth affine algebraic groups $G$ and all $G$-varieties $X$, there exists an $n$-acyclic resolution $j:U \to X$. While this is not the language that the lemma/proposition below uses, a careful reading of the argument on page 149 of \cite{LusztigCuspidal1} shows that the smooth free varities Lusztig constructs are exactly the affine case of the $n$-acyclic varieties that are constructed in \cite[Proposition 6.2.10]{PramodBook}. We will present the argument of \cite{PramodBook}, as it also makes the smooth freeness of the varieties we construct apparent. The existence of the (yet-to-be-defined) category $\ABLDbeqQl{X}$ is also heavily dependent on the existence of $n$-acyclic resolutions of $X$ for all $n \in \N$, so we provide a sketch of the proof here for that reason as well. We only sketch the full proof, however.
\begin{proposition}[\cite{PramodBook}, Proposition 6.2.10]\label{Prop: Section 4: Existence of n acyclic maps}
	Let $n \in \N$ and let $G$ be an affine algebraic group with $X$ a $G$-variety. Then there is an $n$-acyclic resolution $j:U \to X$.
\end{proposition}
\begin{proof}[Sketch]
Using \cite[Theorem 4.8]{milneiAG}, we have that $G$ is a closed subgroup of some $\GL_k$ for some $k \in \N$; moreover, because $G$ is smooth by assumption, $G$ is a smooth closed subgroup of $\GL_k$. Now fix $m \geq k$ and embed $\GL_k$ (call this embedding $i$) into the group $\GL_m$ by using the scheme-theoretic version of the matrix inclusion
	\[
	\begin{pmatrix}
	a_{11} & \cdots & a_{1k} \\
	\vdots & & \vdots \\
	a_{k1} & \cdots & a_{kk}
	\end{pmatrix} \mapsto \begin{pmatrix}
	a_{11} & \cdots & a_{1k} & 0 \\
	\vdots & & \vdots \\
	a_{k1} & \cdots & a_{kk} & 0 \\
	0 & 0 & 0 & I_{m-k}
	\end{pmatrix}
	\]
	where $I_{m-k}$ is the $(m-k) \times (m-k)$ identity matrix. We then can induce a closed subgroup $P$ of $\GL_m$ by the varietal equations describing the span of matrices of the form
	\[
	\begin{pmatrix}
	I_{k} & a_{1,k+1} & \cdots & a_{1m} \\
	0 & \ast & \ast & \ast \\
	0 & a_{m, k+1} & \cdots & a_{mm}
	\end{pmatrix}.
	\]
	When we conjugate $\GL_m$ by entries $\widetilde{\gamma} = i(\gamma)$ for $\gamma \in \GL_k$, we find that the group $P$ is preserved by this conjugation. Thus $G$ preserves $P$ via conjugation. We thus define the group $G \ltimes P$, which is also a closed linear subgroup of $\GL_m$.
	
	Define $W_{m,k}^K$ by the quotient $W_{n,k}^K := \GL_m/P$ (in the complex analytic case $(W_{n,k}^{\C})_{\text{an}}$ is a noncompact Stiefel manifold). The left action of $G$ on $\GL_m$ induced by the embeddings $G \to \GL_k \xrightarrow{i} \GL_m$ induces a left action of $G$ on $W_{m,k}^K$ which realizes $W_{m,k}$ as a smooth principal $G$-variety; in fact, there is an isomorphism of schemes
	\[
	G\backslash W_{m,k}^K \cong \GL_m/(G \ltimes P).
	\]
	It follows also from Theorem \ref{Thm: Comparison} that for any $q \in \Z$, we have isomorphisms of cohomology groups
	\[
	H^{q}(W_{m,k}^{K};\Z_{\ell}) \cong H_{\text{sing}}^{q}\left((W_{m,k}^{\C})_{\text{an}};\Z_{\ell}\right).
	\]
	Moreover, by Theorem \ref{Thm: Stiefel Computation} we have that $H^{q}_{an}(W_{m,k}^K;\Z_{\ell}) \cong 0$ for all $1 \leq q \leq 2m - 2$ and that for any $q \in \Z$, $H^{q}_{an}(W_{m,k}^{K};\Z_{\ell})$ is a free $\Z_{\ell}$-module. Thus it follows that the map $W_{m,k}^K \to \Spec K$ is $(2m-2)$-acyclic by Proposition \ref{Prop: Cohomology Prop}.
	
	We now introduce $X$ into the story. By a base change argument and from the fact that $W_{m,k}^K$ is a smooth free $G$-variety with $K$-morphism $(2m-2)$-acyclic, we have that the variety $W_{m,k}^K \times X$ is a principal $G$-variety as well and that the projection morphism
	\[
	\pi_2:W_{m,k}^K \times X \to X
	\]
	is itself a $(2m-2)$-acyclic resolution of $X$. 
\end{proof}
\begin{lemma}\label{Lemma: Section 5.5: Cofinal acyclic resls}
	Let $\Acal$ denote the class of varieties $\Gamma \in \Sf(G)_0$ for which there is an $n$-acyclic map $\Gamma \to \Gamma^{\prime}$ for some $n \in \N \cup \lbrace \infty \rbrace$ and some $\Gamma^{\prime} \in \Sf(G)_0$ with $\Gamma \ne \Gamma^{\prime}$. Regard $\Sf(G)_0$ as a preorder by saying that $\Gamma^{\prime} \leq \Gamma$ if there is an $n$-acyclic $\Sf(G)$-morphism $f:\Gamma \to \Gamma^{\prime}$ for some $n \in \N \cup \lbrace \infty \rbrace$. Then $\Acal$ is cofinal in $\Sf(G)$.
\end{lemma}
\begin{proof}
	Let $\Gamma \in \Sf(G)_0$ be arbitrary and let $\Gamma^{\prime} \in \Acal$ be an $n$-acyclic resolution of $\Spec K$ for some $n \in \N$, i.e., the structure map $\nu:\Gamma^{\prime} \to \Spec K$ is $n$-acyclic; that these exist is the content of Proposition \ref{Prop: Section 4: Existence of n acyclic maps}. Form the pullback square
	\[
	\xymatrix{
		\Gamma \times \Gamma^{\prime} \ar[r]^{p_2} \ar[d]_{p_1} \pullbackcorner & \Gamma^{\prime} \ar[d]^{\nu} \\
		\Gamma \ar[r] & \Spec K
	}
	\]
	and note that $p_1$ is an $n$-acyclic resolution of $\Gamma$. Then $\Gamma \leq \Gamma \times \Gamma^{\prime}$ so $\Acal$ is indeed cofinal.
\end{proof}

We will need to have an explicit notion of the ABL-cocycle condition we require  the descent data to satisfy before we can define the ABL-equivraiant derived category. This will look like a descent category taken over the acyclic resolutions (cf.\@ \cite{Vistoli} for the general case of descent categories over sites). Let us spell out this condition: 
\begin{definition}\label{Defn: ABL Cocycle cond}\index{ABL-Cocycle Condition}
	Fix three acyclic resolutions $j:U \to X$, $i:V \to X$, and $k:W \to X$ and note that the product $j\times i \times k:U \times V \times W\to X$ is also an acyclic resolution of $X$ as well. We then form the projections in the diagram
	\[
	\begin{tikzcd}
	& U \times V \times W \ar[dr]{}{\pi_{23}} \ar[d]{}[description]{\pi_{13}} \ar[dl, swap]{}{\pi_{12}} & \\
	U \times V \ar[d, swap]{}{\pi_1}  & U \times W \ar[dr, near end]{}[description]{\pi_2} \ar[dl, swap, near end]{}[description]{\pi_1} & V \times W  \ar[d]{}{\pi_2} \\
	U & V & W \ar[from = 2-1, to = 3-2, crossing over, near end]{}[description]{\pi_{2}} \ar[from = 2-3, to = 3-2, crossing over, near end]{}[description]{\pi_1}
	\end{tikzcd}
	\]
	Taking all $G$-quotients in the diagram gives the induced commuting diagrm:
	\[
	\begin{tikzcd}
	& G\backslash (U \times V \times W) \ar[dr]{}{\overline{\pi}_{23}} \ar[d]{}[description]{\overline{\pi}_{13}} \ar[dl, swap]{}{\overline{\pi}_{12}} & \\
	G\backslash (U \times V) \ar[d, swap]{}{\overline{\pi}_1}  & G \backslash(U \times W) \ar[dr, near end]{}[description]{\overline{\pi}_2} \ar[dl, swap, near end]{}[description]{\overline{\pi}_1} & G\backslash (V \times W)  \ar[d]{}{\overline{\pi}_2} \\
	G\backslash U & G \backslash V & G\backslash W \ar[from = 2-3, to = 3-2, crossing over, near end]{}[description]{\overline{\pi}_1} \ar[from = 2-1, to = 3-2, crossing over, near end]{}[description]{\overline{\pi}_{2}}
	\end{tikzcd}
	\]
	The ABL-cocycle condition is then the existence of isomorphisms
	\[
	\overline{\pi}_1^{\ast}\Fscr_{j:U \to X} \xrightarrow{\varphi_{12}} \overline{\pi}_2^{\ast}\Fscr_{i:V \to X}
	\]
	\[
	\overline{\pi}_1^{\ast}\Fscr_{j:U \to X} \xrightarrow{\varphi_{13}} \overline{\pi}_3^{\ast}\Fscr_{k:W \to X}
	\]
	\[
	\overline{\pi}_2^{\ast}\Fscr_{i:V \to X} \xrightarrow{\varphi_{23}} \Fscr_{k:W \to X}
	\]
	which make the diagram, for projection maps $\pi^{123}_{1}:U \times V \times W \to U, \pi^{123}_{2}:U \times V \times W \to V,$ and $\pi_{3}^{123}:U \times V \times W \to W$,
	\[
	\xymatrix{
		(\overline{\pi}^{123}_{1})^{\ast}\Fscr_{j:U \to X} \ar[dr]_{\overline{\pi}_{12}^{\ast}\varphi_{12}} \ar[rr]^-{\overline{\pi}_{13}\varphi_{13}} & & (\overline{\pi}^{123}_{3})^{\ast}\Fscr_{k:W \to X} \\
		& (\overline{\pi}^{123}_{2})\Fscr_{i:V \to X} \ar[ur]_{\overline{\pi}_{23}^{\ast}\varphi_{23}}
	}
	\]
	commute in the category $\DbQl{G \backslash(U \times V \times W)}.$
\end{definition}
\begin{remark}\index{Cocycle Condition}
	While the cocycle condition we present here was certainly not first discovered in \cite{PramodBook} or \cite{BernLun}, this is the most convenient label for this particular cocycle condition in this paper. More explicitly, any other cocycle condition we use will be in the presence of isomorphisms $\of^{\ast}(A) \xrightarrow{\cong} B$ for some objects $A$ and $B$ and some morphisms $\of$ subject to coherence conditions; cf.\@ the equivariant derived category $\DbeqQl{X}$ of Lusztig and the simplicial equivariant derived category $\Dbeqsimp{X}$.
\end{remark}

We now can define the algebraic version of the Bernstein-Lunts equivariant derived category as extended by Achar. Essentially, this equivariant derived category is given as the category of descent data with values in constructible bounded derived categories taken over all acyclic resolutions $j:U \to X$, together with an effective descent condition which states that the descent data is all locally isomorphic to the pullback of a single object $\Fscr \in D_{c}^{b}(X)_0$. 
\begin{definition}[\cite{PramodBook}]\index{Equivariant Derived Category! ABL $\ell$-adic Equivariant Derived Category}
	The category $\ABLDbeqQl{X}$\index[notation]{DbGXABL@$\ABLDbeqQl{X}$} is defined by:
	\begin{itemize}
		\item Objects: Sets 
		\[
		\lbrace \Fscr, \Fscr_{j:U \to X}, \theta_{j:U \to X} \; |n \in \N, \; j:U \to X \, \text{is\, an\,} n\text{-acyclic\, resolution\,} \rbrace,
		\] 
		where:
		\begin{itemize}
			\item $\Fscr \in D_c^{b}(X)_0$;
			\item For all acyclic resolutions $j:U \to X$, $\Fscr_{j:U \to X} \in D^b_c(G\backslash U)_0$;
			\item $\theta_{j:U \to X}$ is an isomorphism
			\[
			j^{\ast}\Fscr \xrightarrow{\theta_{j:U \to X}} q_U^{\ast}\Fscr_{j:U \to X}
			\]
			in $D^b_c(U)_0$ (where $q_U:U \to G\backslash U$ is the quotient map)
		\end{itemize}
		where the sheaves/complexes of sheaves $\Fscr$ are subject to the ABL-cocycle condition of definition \ref{Defn: ABL Cocycle cond}.
		\item Morphisms: A morphism 
		\[
		\rho:(\Fscr, \Fscr_{j:U \to X}, \theta_{j:U \to X}) \to (\Gscr, \Gscr_{j:U \to X}, \sigma_{j:U \to X})
		\] 
		is given by the following information:
		A morphism $\rho \in D_c^b(X)(\Fscr, \Gscr)$ and, for all acyclic resolutions $j:U \to X$, morphisms \[
		\rho_{j} \in \DbQl{G \backslash U}\left(\Fscr_{j:U \to X}, \Gscr_{j:U \to X}\right)
		\] 
		which make the diagrams
		\[
		\xymatrix{
			j^{\ast}\Fscr \ar[r]^{j^{\ast}\rho} \ar[d]_{\theta_{j:U \to X}} & j^{\ast}\Gscr \ar[d]^{\sigma_{j:U \to X}} \\
			q_U^{\ast}\Fscr_{j:U \to X} \ar[r]_{q_U^{\ast}\rho_j} & q_U^{\ast}\Gscr_{j:U \to X}
		}
		\]
		and
		\[
		\xymatrix{
			\overline{\pi}_1^{\ast}\Fscr_{j:U \to X} \ar[d]_{\varphi_{12}^{\Fscr}} \ar[r]^-{\overline{\pi}_1^{\ast}\rho_j} & \overline{\pi}_1 \Gscr_{j:U \to X} \ar[d]^{\varphi_{12}^{\Gscr}} \\
			\overline{\pi}_2^{\ast}\Fscr_{i:V \to X} \ar[r]_-{\overline{\pi}_2^{\ast}\rho_i} & \overline{\pi}_2^{\ast}\Gscr_{i:V \to X}
		}
		\]
		commute whenever $j:U \to X$ and $i:V \to X$ are acyclic resolutions of $X$ and $\overline{\pi}_{1}:G\backslash (U \times V) \to G\backslash U$ and $\overline{\pi}_2:G\backslash(U \times V) \to G\backslash V$ are the quotients of the projection morphisms.
		\item Composition: Term-wise.
		\item Identities: The identity on a triple $\Fscr = (\Fscr, \Fscr_{j:U \to X}, \theta_{j:U \to X})$ is the morphism $(\id_{\Fscr}, \id_{\Fscr_{j:U \to X}})$.
	\end{itemize}
\end{definition}

Instead of comparing $\DbeqQl{X}$ with $\ABLDbeqQl{X}$ directly, we will use an intermediate category, $\Glue(\SfResl_G(X)_{\text{acyc}})$, a  category equivalent to $\ABLDbeqQl{X}$ by Proposition 6.2.17 of \cite{PramodBook} (cf.\@ Proposition \ref{Prop: Section 4: ABL derived cat is SfResl Glue Cat} below as well). Before we can define $\Glue(\SfResl_{G}(X)_{\text{acyc}})$, we will need to make some technical observations for later use. If we have any two $\Sf(G)$ varieties $\Gamma, \Gamma^{\prime}$ we write the product of $\Gamma \times X$ and $\Gamma^{\prime} \times X$ in $\SfResl_G(X)$ as $\Gamma \times \Gamma^{\prime} \times X$; note that this makes sense because the product in $\SfResl_G(X)$ is the pullback over the structure maps $\pi_2^{\Gamma}:\Gamma \times X \to X$ and $\pi_2^{\Gamma^{\prime}}:\Gamma^{\prime} \times X \to X$ and there is an isomorphism of varieties
\[
(\Gamma \times X) \times_X (\Gamma^{\prime} \times X) \cong \Gamma \times \Gamma^{\prime} \times X.
\] 
We then write the projections out of $\Gamma \times \Gamma^{\prime} \times X$ as in the span below:
\[
\xymatrix{
	& \Gamma \times \Gamma^{\prime} \times X \ar[dr]^{\pi_2^{\Gamma\Gamma^{\prime}}} \ar[dl]_{\pi_1^{\Gamma\Gamma^{\prime}}} & \\
	\Gamma \times X & & \Gamma^{\prime} \times X
}
\]
This will be necessary when we consider the induced morphisms
\[
\xymatrix{
	& G\backslash(\Gamma \times \Gamma^{\prime} \times X) \ar[dr]^{\opi_{1}^{\Gamma\Gamma^{\prime}}} \ar[dl]_{\opi_2^{\Gamma\Gamma^{\prime}}}  \\
	G \backslash(\Gamma \times X) & & G \backslash(\Gamma^{\prime} \times X)
}
\]
when we define and work with $\Glue(\SfResl_G(X)_{\text{acyc}})$.

\begin{definition}
	We define the category
	\[
	\SfResl_G(X)_{\text{acyc}}
	\]\index[notation]{SfReslGAcyc@$\SfResl_G(X)_{\text{acyc}}$}
	to be the full subcategory of $\SfResl_G(X)$ generated by the objects $\Gamma \times X \in \SfResl_G(X)_0$ such that the projection morphism 
	\[
	\pi_2^{\Gamma}:\Gamma \times X \to X
	\]
	is an $n$-acyclic resolution of $X$ for some $n \in \N$. Similarly, $\Sf(G)_{\text{acyc}}$ is the full subcategory of $\Sf(G)$ generated by the objects $\Gamma$ for which the structure morphism $\Gamma \to \Spec K$ is $n$-acyclic for some $n \in \N$.
\end{definition}

\begin{definition}[\cite{PramodBook}]\index[notation]{GlueSfResl@$\Glue(\SfResl_G(X))$}\index{Acylcic Gluing Category}
	The category $\Glue(\SfResl_G(X))$ is defined by:
	\begin{itemize}
		\item Objects: Triples
		\[
		(\Fscr, \lbrace \quot{\Fscr}{\Gamma} \; | \; \Gamma \in (\Sf(G)_{\text{acyc}})_0 \rbrace, \lbrace \quot{\theta}{\Gamma} \; | \; \Gamma \in \Sf(G)_0\rbrace)
		\] 
		where $\Fscr \in\DbQl{X}_0$, for all $\Gamma \times X \in (\SfResl_G(X)_{\text{acyc}})_0$ the complex $\quot{\Fscr}{\Gamma} \in \DbQl{\XGamma}_0$, and $\quot{\theta}{\Gamma}$ is an isomorphism
		\[
		\quot{\theta}{\Gamma}:(\pi_2^{\Gamma})^{\ast}\Fscr \xrightarrow{\cong} q_{\Gamma}^{\ast}\quot{\Fscr}{\Gamma}
		\]
		in $\DbQl{\Gamma \times X}$. We also assert that the $\quot{\Fscr}{\Gamma}$ satisfy the ABL-cocycle condition.
		\item Morphisms: As in $\ABLDbeqQl{X}$, save now restricted those morphisms whose $\quot{\rho}{\Gamma}$-components take local values in $\DbQl{\XGamma}$.
		\item Composition: Term-wise.
		\item Identities: The identity on an object $\Fscr$ is $(\id_{\Fscr}, \id_{\quot{\Fscr}{\Gamma}})$.
	\end{itemize}
\end{definition}
\begin{remark}\label{Remark: Equivalence of All Resls with Acylcic Resls}
	This category can be defined for any category $\Cscr$ of $G$-resolutions of $X$, provided that $\Cscr$ has all products of resolutions and provided that we have at least one $n$-acyclic resolution for each $n \in \N$. In such a case we get an equivalence of categories $\Glue(\Cscr) \simeq \Glue(\Cscr_{\text{acyc}})$\index[notation]{GlueC@$\Glue(\Cscr)$} by acyclic descent, where $\Cscr_{\text{acyc}}$ is the subcategory of $\Cscr$ induced by the acyclic resolutions of $X$ that lie in $\Cscr$. We will implicitly use this as the definition of the category $\Glue(\lbrace j_{\alpha}:U_{\alpha} \to X \; | \; \alpha  \in I \rbrace)$ throughout Part \ref{Chapter: EDC Comp}, as we can then define $\Glue(\SfResl_G(X))$ on all varieties $\Gamma \times X$ instead of just those that arise from $\SfResl_G(X)_{\text{acyc}}$. In particular, we have
	\[
	\Glue(\SfResl_G(X)) \simeq \Glue(\SfResl_{G}(X)_{\text{acyc}})
	\]
	so when working up to equivalence our choice of category does not matter.
\end{remark}

\begin{proposition}[\cite{PramodBook}, Proposition 6.2.17]\label{Prop: Section 4: ABL derived cat is SfResl Glue Cat}
	Let $C$ be a collection $\lbrace j_{\alpha}:U_{\alpha} \to X \; | \; \alpha \in I \rbrace$ of acyclic resolutions $j_{\alpha}$ of $X$ with the property that for all $n \in \N$, there is an $n$-acyclic resolution of $X$. there is an equivalence of categories
	\[
	\Glue(C) \simeq \quot{D_G^b(X)}{ABL}.
	\]
\end{proposition}
\begin{corollary}
	There is an equivalence of categories
	\[
	\Glue(\SfResl_G(X)) \simeq \ABLDbeqQl{X} \simeq \Glue(\SfResl_G(X)_{\text{acyc}}).
	\]
\end{corollary}
\begin{proof}
	By Proposition \ref{Prop: Section 4: Existence of n acyclic maps}, we have that there are $n$-acyclic $\Sf(G)$ varieties for all $n \in \N$.  Using Remark \ref{Remark: Equivalence of All Resls with Acylcic Resls} and Proposition \ref{Prop: Section 4: ABL derived cat is SfResl Glue Cat} then gives the desired equivalence.
\end{proof}

We now move to compare the categories $\DbeqQl{X}$ and $\Glue(\SfResl_G(X))$ by proving some structural lemmas about the three categories $\DbeqQl{X}$, $\SfResl_G(X)$, and $\Glue(\SfResl_G(X))$ so that we know how to compare them. After this we will describe the comparison functors we use and finally prove that they give an equivalence of categories.
\begin{lemma}\label{Lemma: Section 4: Cocycle condition for equivariant derived cat}
	Let $(A,T_A)$ be an object in $\DbeqQl{X}$. Then for all $\Gamma, \Gamma^{\prime} \in \Sf(G)_0$, there are isomorphisms 
	\[
	(\overline{\pi}_1^{\Gamma\Gamma^{\prime}})^{\ast}\AGamma \cong (\overline{\pi}_2^{\Gamma\Gamma^{\prime}})^{\ast}\quot{A}{\Gamma^{\prime}}
	\]
	which satisfy the ABL-cocycle condition.
\end{lemma}
\begin{proof}
	Begin by noting that  the variety $\Gamma \times \Gamma^{\prime} \in \Sf(G)_0$ and the morphisms $\pi_1:\Gamma \times \Gamma^{\prime} \to \Gamma$ and $\pi_2:\Gamma \times \Gamma^{\prime} \to \Gamma^{\prime}$ are both $\Sf(G)$-morphism. Consequently the morphisms $\pi_{1}^{\Gamma\Gamma^{\prime}}:\Gamma \times \Gamma^{\prime} \times X \to \Gamma \times X$ and $\pi_2^{\Gamma\Gamma^{\prime}}:\Gamma \times \Gamma^{\prime} \times X \to \Gamma^{\prime} \times X$ are both $\SfResl_G(X)$-morphisms. Since $(A,T_A)$ is an object in $\DbeqQl{X}$, we have isomorphisms
	\[
	(\overline{\pi}_1^{\Gamma\Gamma^{\prime}})^{\ast}\AGamma \xrightarrow{\tau_{\pi_1}^{A}} \quot{A}{\Gamma \times \Gamma^{\prime}}
	\]
	and
	\[
	(\overline{\pi}_{2}^{\Gamma\Gamma^{\prime}})^{\ast}\AGammap \xrightarrow{\tau_{\pi_2}^{A}} \quot{A}{\Gamma \times \Gamma^{\prime}}
	\]
	which satisfy the cocycle condition of Definition \ref{Defn: Equivariant Cats}. We now define the isomorphism
	\[
	\varphi_{\Gamma\Gamma^{\prime}}:(\overline{\pi}_1^{\Gamma\Gamma^{\prime}})^{\ast}\AGamma \xrightarrow{\cong} (\overline{\pi}_2^{\Gamma\Gamma^{\prime}})^{\ast}\AGammap
	\]
	via the diagram:
	\[
	\xymatrix{
		\left(\opi_1^{\Gamma\Gamma^{\prime}}\right)^{\ast}\AGammap \ar[rr]^-{\varphi_{\Gamma\Gamma^{\prime}}} \ar[dr]_{\tau_{\pi_1}^A} & & \left(\opi_2^{\Gamma\Gamma^{\prime}}\right)^{\ast}\AGamma \\
		& \quot{A}{\Gamma \times \Gamma^{\prime}} \ar[ur]_{(\tau_{\pi_2}^{A})^{-1}}
	}
	\]
	
	To see the ABL-cocycle condition, note that by Lemmas \ref{Lemma: Section 5.5: Acyclic Base Change} if $\Gamma \in \Sf(G)_0$ is $n$-acyclic, the resolution $\pi_2^{\Gamma}:\Gamma \times X \to X$ is $n$-acyclic as well. From here carrying the proof of Lemma \ref{Lemma: Section 5.5: Cofinal acyclic resls} over to the $\SfResl_G(X)$ case allows us to deduce that the class $\Acal_X$ of acyclic resolutions in $\SfResl_{G}(X)$ is also cofinal. As such we can always arrange our three objects $\Gamma \times X, \Gamma^{\prime} \times X, \Gamma^{\prime\prime} \times X$ to be given so that each of there is a commuting diagram
	\[
	\begin{tikzcd}
	& \Gamma \times \Gamma^{\prime} \times \Gamma^{\prime\prime} \times X \ar[ddl, swap]{}{\pi_{12}^{\Gamma\Gamma^{\prime}\Gamma^{\prime\prime}}} \ar[dd]{}[description]{\pi_{13}^{\Gamma\Gamma^{\prime}\Gamma^{\prime\prime}}} \ar[ddr]{}{\pi^{\Gamma\Gamma^{\prime}\Gamma^{\prime\prime}}_{23}} \\
	\\
	\Gamma \times \Gamma^{\prime} \times X \ar[dd, swap]{}{\pi^{\Gamma\Gamma^{\prime}}_1} & \Gamma \times \Gamma^{\prime\prime} \times X \ar[ddl, near end, swap]{}{\pi^{\Gamma\Gamma^{\prime\prime}}_1} \ar[ddr, near end]{}{\pi_{2}^{\Gamma\Gamma^{\prime\prime}}} & \Gamma^{\prime} \times \Gamma^{\prime\prime} \times X \ar[dd]{}{\pi^{\Gamma^{\prime}\Gamma^{\prime\prime}}_2} \\
	\\
	\Gamma \times X \ar[r, swap]{}{f \times \id_X} \ar[rr, bend right = 30, swap]{}{(g \circ f) \times \id_X} & \Gamma^{\prime} \times X \ar[r, swap]{}{g \times \id_X} & \Gamma^{\prime\prime} \times X \ar[from = 3-1, to = 5-2, crossing over, near end]{}{\pi^{\Gamma\Gamma^{\prime}}_2} \ar[from = 3-3, to =  5-2, crossing over, near end, swap]{}{\pi_{1}^{\Gamma^{\prime}\Gamma^{\prime\prime}}}-
	\end{tikzcd}
	\]
	with morphisms in $f, g \in\Sf(G)_1$. We now note that on one hand we have
	\begin{align*}
	&\left(\opi_{23}^{\Gamma\Gamma^{\prime}\Gamma^{\prime\prime}}\right)^{\ast}\varphi_{\Gamma^{\prime}\Gamma^{\prime\prime}} \circ \left(\opi_{12}^{\Gamma\Gamma^{\prime\prime}}\right)^{\ast}\varphi_{\Gamma\Gamma^{\prime}}\\
	 &= \left(\opi_{23}^{\Gamma\Gamma^{\prime}\Gamma^{\prime\prime}}\right)^{\ast}\left(\big(\tau_{\pi_{2}^{\Gamma^{\prime}\Gamma^{\prime\prime}}}^{A}\big)^{-1} \circ \tau_{\pi_{1}^{\Gamma^{\prime}\Gamma^{\prime\prime}}}^{A}\right) \circ \left(\opi_{12}^{\Gamma\Gamma^{\prime\prime}}\right)^{\ast}\left(\big(\tau^{A}_{\pi_2^{\Gamma\Gamma^{\prime}}}\big)^{-1} \circ \tau_{\pi_1^{\Gamma\Gamma^{\prime}}}^{A}\right).
	\end{align*}
	Using the commutative diagram above to get that $\opi_2^{\Gamma\Gamma^{\prime}} = \of \circ \opi_1^{\Gamma\Gamma^{\prime}}$ and deduce that
	\[
	\tau_{\pi_2^{\Gamma\Gamma^{\prime}}}^{A} = \tau_{\pi_1^{\Gamma\Gamma^{\prime}}}^{A} \circ \left(\opi_1^{\Gamma\Gamma^{\prime}}\right)^{\ast}\tau_f^A.
	\]
	Substituting this for $\tau_{\pi_2^{\Gamma\Gamma^{\prime}}}^{A}$ in the equation above gives
	\[
	\big(\tau_{\pi_2^{\Gamma\Gamma^{\prime}}}^{A}\big)^{-1} \circ \tau_{\pi_1^{\Gamma\Gamma^{\prime}}}^{A} = \left(\tau_{\pi_1^{\Gamma\Gamma^{\prime}}}^{A} \circ \left(\opi_1^{\Gamma\Gamma^{\prime}}\right)^{\ast}\tau_{f}^{A}\right)^{-1} \circ \tau_{\pi_{1}^{\Gamma\Gamma^{\prime}}}^{A} = \left(\opi_1^{\Gamma\Gamma^{\prime}}\right)^{\ast}(\tau_{f}^A)^{-1};
	\]
	similarly using $\pi_2^{\Gamma^{\prime}\Gamma^{\prime\prime}} = (g \circ \id_X) \circ \pi_{1}^{\Gamma^{\prime}\Gamma^{\prime\prime}}$ allows us to deduce
	\[
	\big(\tau_{\pi_{2}^{\Gamma^{\prime}\Gamma^{\prime\prime}}}^{A}\big)^{-1} \circ \tau_{\pi_{1}^{\Gamma^{\prime}\Gamma^{\prime\prime}}}^{A} = \big(\opi_1^{\Gamma^{\prime}\Gamma^{\prime\prime}}\big)^{\ast}(\tau_g^A)^{-1}.
	\]
	as well. Plugging these into our equation gives
	\begin{align*}
	&\left(\opi_{23}^{\Gamma\Gamma^{\prime}\Gamma^{\prime\prime}}\right)^{\ast}\varphi_{\Gamma^{\prime}\Gamma^{\prime\prime}} \circ \left(\opi_{12}^{\Gamma\Gamma^{\prime\prime}}\right)^{\ast}\varphi_{\Gamma\Gamma^{\prime}} \\
	&= \left(\opi_{23}^{\Gamma\Gamma^{\prime}\Gamma^{\prime\prime}}\right)^{\ast}\left(\big(\opi_1^{\Gamma^{\prime}\Gamma^{\prime\prime}}\big)^{\ast}(\tau_g^A)^{-1}\right) \circ \left(\opi_{12}^{\Gamma\Gamma^{\prime\prime}}\right)^{\ast}\left(\big(\opi_1^{\Gamma\Gamma^{\prime}}\big)^{\ast}(\tau_{f}^A)^{-1}\right) \\
	&= \left(\opi_1^{\Gamma^{\prime}\Gamma^{\prime\prime}} \circ \opi_{23}^{\Gamma\Gamma^{\prime}\Gamma^{\prime\prime}}\right)^{\ast}(\tau_g^{A})^{-1} \circ \left(\opi_1^{\Gamma\Gamma^{\prime}} \circ \opi_{12}^{\Gamma\Gamma^{\prime}\Gamma^{\prime\prime}}\right)^{\ast}(\tau_f^A)^{-1} \\
	&= \left(\opi_2^{\Gamma\Gamma^{\prime}} \circ \opi_{12}^{\Gamma\Gamma^{\prime}\Gamma^{\prime\prime}}\right)^{\ast}(\tau_g^{A})^{-1} \circ \left(\opi_1^{\Gamma\Gamma^{\prime}} \circ \opi_{12}^{\Gamma\Gamma^{\prime}\Gamma^{\prime\prime}}\right)^{\ast}(\tau_f^A)^{-1} \\
	&= \left(\of \circ \opi_1^{\Gamma\Gamma^{\prime}} \circ \opi_{12}^{\Gamma\Gamma^{\prime}\Gamma^{\prime\prime}}\right)^{\ast}(\tau_g^{A})^{-1} \circ \left(\opi_1^{\Gamma\Gamma^{\prime}} \circ \opi_{12}^{\Gamma\Gamma^{\prime}\Gamma^{\prime\prime}}\right)^{\ast}(\tau_f^A)^{-1} \\
	&= \left(\opi_1^{\Gamma\Gamma^{\prime}} \circ \opi_{12}^{\Gamma\Gamma^{\prime}\Gamma^{\prime\prime}}\right)^{\ast}(\of^{\ast}\tau_g^{A})^{-1} \circ \left(\opi_1^{\Gamma\Gamma^{\prime}} \circ \opi_{12}^{\Gamma\Gamma^{\prime}\Gamma^{\prime\prime}}\right)^{\ast}(\tau_f^A)^{-1} \\
	&=\left(\opi_1^{\Gamma\Gamma^{\prime}} \circ \opi_{12}^{\Gamma\Gamma^{\prime}\Gamma^{\prime\prime}}\right)^{\ast}\left(\tau_f^A \circ \of^{\ast}\tau_g^A\right)^{-1} \\
	&=\left(\opi_1^{\Gamma\Gamma^{\prime}} \circ \opi_{12}^{\Gamma\Gamma^{\prime}\Gamma^{\prime\prime}}\right)^{\ast}(\tau_{g \circ f}^{A})^{-1} = \left(\opi_1^{\Gamma\Gamma^{\prime\prime}} \circ \opi_{13}^{\Gamma\Gamma^{\prime}\Gamma^{\prime\prime}}\right)^{\ast}(\tau_{g \circ f}^{A})^{-1} \\
	&= \left(\opi_{13}^{\Gamma\Gamma^{\prime}\Gamma^{\prime\prime}}\right)^{\ast}\left(\big(\opi_1^{\Gamma\Gamma^{\prime\prime}}\big)^{\ast}\tau_{g \circ f}^{A}\right)^{-1}
	\end{align*}
	We now observe that since $\pi_2^{\Gamma \Gamma^{\prime\prime}} = \big((g \circ f) \times \id_X\big) \circ \pi_1^{\Gamma\Gamma^{\prime\prime}}$ we get that
	\[
	\tau_{\pi_2^{\Gamma\Gamma^{\prime\prime}}}^{A} = \tau_{\pi_1^{\Gamma\Gamma^{\prime\prime}}}^{A} \circ \left(\opi_{1}^{\Gamma\Gamma^{\prime}}\right)^{\ast}\tau_{g \circ f}^{A}
	\]
	which may be rearranged to
	\[
	\left(\opi_{1}^{\Gamma\Gamma^{\prime}}\right)^{\ast}\tau_{g \circ f}^{A} = \left(\tau_{\pi_1^{\Gamma\Gamma^{\prime\prime}}}^{A}\right)^{-1} \circ \tau_{\pi_2^{\Gamma\Gamma^{\prime\prime}}}^{A}.
	\]
	Substituting this in place of $\tau_{g \circ f}^{A}$ then yields
	\begin{align*}
	\left(\opi_{13}^{\Gamma\Gamma^{\prime}\Gamma^{\prime\prime}}\right)^{\ast}\left(\big(\opi_1^{\Gamma\Gamma^{\prime\prime}}\big)^{\ast}\tau_{g \circ f}^{A}\right)^{-1} &= \left(\opi_{13}^{\Gamma\Gamma^{\prime}\Gamma^{\prime\prime}}\right)^{\ast}\left(\left(\tau_{\pi_1^{\Gamma\Gamma^{\prime\prime}}}^{A}\right)^{-1} \circ \tau_{\pi_2^{\Gamma\Gamma^{\prime\prime}}}^{A}\right)^{-1}\\ 
	&= \left(\opi_{13}^{\Gamma\Gamma^{\prime}\Gamma^{\prime\prime}}\right)^{\ast}\left(\left(\tau_{\pi_2^{\Gamma\Gamma^{\prime\prime}}}^{A}\right)^{-1}  \circ \tau_{\pi_1^{\Gamma\Gamma^{\prime\prime}}}^{A}\right) 	\\
	&= \left(\opi_{13}^{\Gamma\Gamma^{\prime}\Gamma^{\prime\prime}}\right)^{\ast}\varphi_{\Gamma\Gamma^{\prime\prime}}.
	\end{align*}
	With it follows that
	\[
	\left(\opi_{23}^{\Gamma\Gamma^{\prime}\Gamma^{\prime\prime}}\right)^{\ast}\varphi_{\Gamma^{\prime}\Gamma^{\prime\prime}} \circ \left(\opi_{12}^{\Gamma\Gamma^{\prime\prime}}\right)^{\ast}\varphi_{\Gamma\Gamma^{\prime}} = \left(\opi_{13}^{\Gamma\Gamma^{\prime}\Gamma^{\prime\prime}}\right)^{\ast}\varphi_{\Gamma\Gamma^{\prime\prime}}.
	\]
	Because this is the ABL-cocycle condition, this completes the proof of the lemma.
\end{proof}

\begin{lemma}\label{Lemma: Section 4: Cocycle on ABL cat is cocycle in our cat}
	Let $(\Fscr, \quot{\Fscr}{\Gamma},\quot{\theta}{\Gamma})$ be an object in $\Glue(\SfResl_G(X))$. Then for any $f \in \Sf(G)_1$ there is an isomorphism
	\[
	\tau_{f}^A:\overline{f}^{\ast}\AGammap \to \AGamma
	\]
	which satisfies the cocycle condition
	\[
	\tau_{g \circ f}^{A} = \tau_f^A \circ \overline{f}^{\ast}\tau_g^A
	\]
	for all composable arrows $\Gamma \xrightarrow{f} \Gamma^{\prime} \xrightarrow{g} \Gamma^{\prime\prime}$ in $\Sf(G)$.
\end{lemma}
\begin{proof}
	Begin by recalling that the product $\Gamma \times \Gamma^{\prime} \times X$ is $\Sf(G)$ arises as the pullback
	\[
	\xymatrix{
		\Gamma \times \Gamma^{\prime} \times X \ar[r] \ar[d] \pullbackcorner & \Gamma^{\prime} \times X \ar[d]^{\pi_2^{\Gamma^{\prime}}} \\
		\Gamma \times X \ar[r]_-{\pi_2^{\Gamma}} & X
	}
	\]
	in $\GVar$. In particular, if $f \in \Sf(G)(\Gamma, \Gamma^{\prime})$ the diagram
	\[
	\xymatrix{
		\Gamma \times \Gamma^{\prime} \times X \ar[r] \ar[d] \pullbackcorner & \Gamma^{\prime} \times X \ar[d]^{\pi_2^{\Gamma^{\prime}}} \\
		\Gamma \times X \ar[ur]_{f \times \id_X} \ar[r]_-{\pi_2^{\Gamma}} & X
	}
	\]
	commutes, which implies that the diagram of quotient varieites
	\[
	\xymatrix{
		& \quot{X}{\Gamma \times \Gamma^{\prime}} \ar[dr]^{\overline{\pi}_2^{\Gamma\Gamma^{\prime}}} \ar[dl]_{\overline{\pi}_1^{\Gamma\Gamma^{\prime}}} & \\
		\XGamma \ar[rr]_-{\overline{f}} & & \XGammap
	}
	\]
	commutes as well. We now calculate that
	\[
	(\overline{\pi}_2^{\Gamma\Gamma^{\prime}})^{\ast}\quot{\Fscr}{\Gamma^{\prime}} = (\overline{f} \circ \overline{\pi}_1^{\Gamma\Gamma^{\prime}})^{\ast}\quot{\Fscr}{\Gamma^{\prime}} = (\overline{\pi}_1^{\Gamma\Gamma^{\prime}})^{\ast}\left(\overline{f}^{\ast}\quot{\Fscr}{\Gamma^{\prime}}\right).
	\]
	From the isomorphism $(\overline{\pi}_1^{\Gamma\Gamma^{\prime}})^{\ast}\quot{\Fscr}{\Gamma} \cong (\overline{\pi}_2^{\Gamma\Gamma^{\prime}})^{\ast}\quot{\Fscr}{\Gamma^{\prime}}$ we induce an isomorphism
	\[
	(\overline{\pi}_1^{\Gamma\Gamma^{\prime}})^{\ast}\quot{\Fscr}{\Gamma}  \cong (\overline{\pi}_1^{\Gamma\Gamma^{\prime}})^{\ast}\left(\overline{f}^{\ast}\quot{\Fscr}{\Gamma^{\prime}}\right).
	\]
	To see that this produces our desired isomorphism of complexes, we give a sketch by checking on stalks; while the argument we present is not entirely rigorous, it sketches the idea in a much more clear fashion and the translation to general points may be made by looking at the points each covers by the corresponding projection maps, i.e., if $\alpha$ is a point of $\quot{X}{\Gamma \times \Gamma^{\prime}}$ the translation is given by looking at the stalks $\opi_1^{\Gamma\Gamma^{\prime}}(\alpha)$ of $\XGamma$ and $\opi_2^{\Gamma\Gamma^{\prime}}(\alpha)$ of $\XGammap$. At any rate, assume that $[\gamma, \gamma^{\prime},x]_{G}$ is a point of $\quot{X}{\Gamma \times \Gamma^{\prime}}$. Because of the isomorphism $\big(\opi_1^{\Gamma\Gamma^{\prime}}\big)^{\ast}(\quot{\Fscr}{\Gamma}) \cong \big(\opi_1^{\Gamma\Gamma^{\prime}}\big)^{\ast}(\of^{\ast}\quot{\Fscr}{\Gamma^{\prime}})$ we get an isomorphism of stalks
	\[
	\big(\opi_1^{\Gamma\Gamma^{\prime}}\big)^{\ast}(\quot{\Fscr}{\Gamma})_{[\gamma,\gamma^{\prime},x]_{G}} \cong \big(\opi_1^{\Gamma\Gamma^{\prime}}\big)^{\ast}(\of^{\ast}\quot{\Fscr}{\Gamma^{\prime}})_{[\gamma,\gamma^{\prime},x]_G}.
	\]
	However, using the isomorphism
	\[
	\big(\opi_1^{\Gamma\Gamma^{\prime}}\big)^{\ast}(\quot{\Fscr}{\Gamma})_{[\gamma,\gamma^{\prime},x]_{G}} \cong \quot{\Fscr}{\Gamma}_{\opi_1^{\Gamma\Gamma^{\prime}}[\gamma,\gamma^{\prime},x]_{G}} = \quot{\Fscr}{\Gamma}_{[\gamma,x]_{G}}
	\]
	(and similarly for $\big(\opi_1^{\Gamma\Gamma^{\prime}}\big)^{\ast}(\of^{\ast}\quot{\Fscr}{\Gamma^{\prime}})_{[\gamma,\gamma^{\prime},x]_G}$) we get that there are isomorphisms
	\[
	\quot{\Fscr}{\Gamma}_{[\gamma,x]_{G}} \cong (\of^{\ast}\quot{\Fscr}{\Gamma^{\prime}})_{[\gamma,x]_{G}}
	\]
	for each point $[\gamma,x]_{G}$ of $\XGamma$. Gluing these morphisms together (which may be done by construction) produces an induced morphism $\tau_f^{A}:\of^{\ast}(\quot{\Fscr}{\Gamma^{\prime}}) \to \quot{\Fscr}{\Gamma}$ with stalks the isomorphisms described above. However, since the {\'e}tale topology has enough points and stalks, we conclude from the isomorphisms at all stalks that $\tau_f^A$ is an isomorphism of complexes of $\ell$-adic sheaves on $\XGamma$.
	
	That these isomorphisms satisfy the desired cocycle condition is a tedious but routine check in the same vein as the one in Lemma \ref{Lemma: Section 4: Cocycle condition for equivariant derived cat}. In particular, assuming the existence of the composable pair of morphisms $\Gamma \xrightarrow{f} \Gamma^{\prime} \xrightarrow{g} \Gamma^{\prime\prime}$ in $\Sf(G)$, we can deduce the existence of the commuting diagram:
	\[
	\begin{tikzcd}
	& \quot{X}{\Gamma \times \Gamma^{\prime} \times \Gamma^{\prime\prime}} \ar[dd]{}[description]{\opi_{13}^{\Gamma\Gamma^{\prime}\Gamma^{\prime\prime}}} \ar[ddr]{}{\opi_{23}^{\Gamma\Gamma^{\prime}\Gamma^{\prime\prime}}} \ar[ddl, swap]{}{\opi_{12}^{\Gamma\Gamma^{\prime}\Gamma^{\prime\prime}}} \\
	\\
	\quot{X}{\Gamma \times \Gamma^{\prime}} \ar[dd, swap]{}{\opi_{1}^{\Gamma\Gamma^{\prime}}} & \quot{X}{\Gamma \times \Gamma^{\prime\prime}} \ar[ddr, near end]{}{\opi_2^{\Gamma\Gamma^{\prime\prime}}} \ar[ddl, near end, swap]{}{\opi_{1}^{\Gamma\Gamma^{\prime\prime}}} & \quot{X}{\Gamma^{\prime} \times \Gamma^{\prime\prime}} \ar[dd]{}{\opi_{2}^{\Gamma^{\prime}\Gamma^{\prime\prime}}} \\
	\\
	\XGamma \ar[r, swap]{}{\of} \ar[rr, swap, bend right = 30]{}{\overline{g} \circ \of} & \XGammap \ar[r, swap]{}{\overline{g}} & \XGammapp \ar[from = 3-1, to = 5-2, crossing over, near end]{}{\opi_{2}^{\Gamma\Gamma^{\prime}}} \ar[from = 3-3, to = 5-2, crossing over, near end, swap]{}{\opi_1^{\Gamma^{\prime}\Gamma^{\prime\prime}}}
	\end{tikzcd}
	\]
	From here doing the same sort of cancellation tricks with stalks that we did earlier, coupled with the ABL-cocycle identity and liberal use of the identities induced from the diagram on pullback functors  allows us to deduce that $\tau_{g \circ f}^{A} = \tau_f^A \circ \of^{\ast}\tau_g^A$.
\end{proof}

This lemma above gives rise to our first comparison functor, as it allows us to take an effective gluing descent object from $\Glue(\SfResl_G(X))$ and turn it into an object in $\DbeqQl{X}$. We will present the existence and definition of this functor below, and then move on with the remaining necessary constructions that allow us to produce the inverse comparison.
\begin{proposition}\label{Prop: Section 4: Comparison functor from ABL to equivariant cat}
	There is a functor $\beta:\Glue(\SfResl_G(X)) \to \DbeqQl{X}$ which is defined on objects by sending a triple
	\[
	(\Fscr, \Fscr_{\Gamma \times X \to X}, \theta_{\Gamma \times X \to X})
	\] 
	to the object
	\[
	\big(\lbrace \Fscr_{\Gamma \times X \to X} \; | \; \Gamma \in (\Sf(G)_0 \rbrace, \lbrace \tau_f^{\Fscr} \; | \; f \in \Sf(G)_1 \rbrace\big),
	\]
	where the $\tau_f^{\Fscr}$ are those isomorphisms given in Lemma \ref{Lemma: Section 4: Cocycle on ABL cat is cocycle in our cat}, and defined on morphisms by sending a morphism pair
	\[
	(\rho,\rho_{\Gamma \times X \to \Gamma} )
	\]
	to the morphism
	\[
	\lbrace \rho_{\Gamma \times X \to X} \; | \; \Gamma \in \Sf(G)_0 \rbrace.
	\]
\end{proposition}
\begin{proof}
	That this is well-defined is immediate from Lemma \ref{Lemma: Section 4: Cocycle on ABL cat is cocycle in our cat}, as it is only the cocycle condition that causes trouble. It follows from construction that this assignment is functorial as well.
\end{proof}

We now proceed to consrtruct the comparison functor $\alpha:\DbeqQl{X} \to \Glue(\SfResl_G(X))$: The $ABL$-ification functor. What we essentially need to do to define this functor is take descent data $(A,T_A)$ in $\DbeqQl{X}$ and then use it to produce a sheaf on $X$ which is effective with respect to the descent data $(A,T_A)$. Explicitly, given an object $(A,T_A)$ with
\[
A = \lbrace \AGamma \; | \; \Gamma \in \Sf(G)_0 \rbrace
\]
and
\[
T_A = \lbrace \tau_f^A \; | \; f \in \Sf(G)_1 \rbrace
\]
we need to be able to construct an object $(\Fscr, \Fscr_{\Gamma \times X \to X}, \theta_{\Gamma})$ in $\Glue(\SfResl_G(X))$. By Lemma \ref{Lemma: Section 4: Cocycle condition for equivariant derived cat} we can define, for all $\Gamma \in \Sf(G)_0$,
\[
\Fscr_{\Gamma \times X \to X} := \AGamma,
\]
as the collection of complexes of sheaves
\[
A = \lbrace \AGamma \; | \; \Gamma \in \Sf(G)_0 \rbrace
\] 
satisfies the appropriate cocycle conditions (they are induced by $T_A$ as in Lemma \ref{Lemma: Section 4: Cocycle condition for equivariant derived cat}). We now need to construct the objects $\Fscr$ and $\theta_{\Gamma}$. To see a natural object $\Fscr$ on $X$, consider first that there is a non-$G$-equivariant isomorphism of $K$-schemes
\[
X \xrightarrow[\cong]{\psi} G \backslash (G \times X) = \quot{X}{G};
\]
in particular, this isomorphism induces an equivalence of categories $\psi^{\ast}:\DbQl{\quot{X}{G}} \xrightarrow{\simeq} \DbQl{X}$. We thus define the object $\Fscr$ via
\[
\Fscr := \psi^{\ast}\quot{A}{G}.
\]
We now prove the following crucial lemma that will allow us to conclude that the assignment $(A, T_A) \mapsto (\Fscr, \Fscr_{\Gamma \times X \to X}, \theta_{\Gamma})$ is the object assignment of our yet-to-be-defined comparison functor $\alpha$.
\begin{lemma}\label{Lemma: Section 4: Descent on object in DGbX is effective}
	Let $(A,T_A) \in D_G^b(X)_0$ with $A = \lbrace \AGamma \; | \; \Gamma \in \Sf(G)_0\rbrace$ and $T_A = \lbrace \tau_f^A \; | \; f \in \Sf(G)_1 \rbrace$. Then for all $\Gamma \in \Sf(G)_0$ there is an isomorphism
	\[
	\theta_{\Gamma}:(\pi_2^{\Gamma})^{\ast}\psi^{\ast}\quot{A}{G} \xrightarrow{\cong} (\quo_{\Gamma})^{\ast}\AGamma.
	\]
\end{lemma}
\begin{proof}
	To prove this we begin by letting $\Gamma \in \Sf(G)_0$ and considering the diagram:
	\[
	\xymatrix{
		\Gamma \times X \ar[r]^-{\pi_2^{\Gamma}}\ar[d]_{\quo_{\Gamma}} & X \ar[d]^{\psi} \\
		\XGamma \ar[r]_-{\overline{\varphi}} & \quot{X}{G}
	}
	\]
	Because the class of acyclic resolutions in $\SfResl_G(X)$ is cofinal by Lemma \ref{Lemma: Section 5.5: Cofinal acyclic resls} and the proof of Lemma \ref{Lemma: Section 4: Cocycle condition for equivariant derived cat}, and because of  the presence of the isomorphism $\psi$, it suffices to show by acyclic descent that for an object in $\SfResl_{G}(X)$ with $n$-acyclic resolution $\pi_2^{\Gamma}:\Gamma \times X \to X$ and an $\Sf(G)$-morphism $f:\Gamma \to G$ that
	\[
	(\quo_{\Gamma})^{\ast}(\AGamma) \cong (\pi_2^{\Gamma})^{\ast}\big(\psi^{\ast}\quot{A}{G}\big).
	\]
	However, for this observe that the diagram
	\[
	\begin{tikzcd}
	\Gamma \times X \ar[rr]{}{\pi_2^{\Gamma}} \ar[dr, swap]{}{f \times \id_X} \ar[dd, swap]{}{\quo_{\Gamma}} & & X \ar[dd, equals] \\
	& G \times X \ar[ur, swap]{}{\pi_2^G}  \\
	\XGamma \ar[rr, near start]{}{\overline{\varphi}} \ar[dr, swap]{}{\of} & & X \ar[dl]{}{\psi} \\
	& \quot{X}{G} \ar[from = 2-2, to = 4-2, crossing over, near start]{}[description]{\quo_G}
	\end{tikzcd}
	\]
	commutes in $\Var_{/\Spec K}$. From this we derive that
	\begin{align*}
	\left(\pi_2^{\Gamma}\right)^{\ast}\left(\psi^{\ast}\quot{A}{G}\right) &= \left(\psi \circ \pi_2^{\Gamma}\right)^{\ast}(\quot{A}{G}) \cong \left(\of \circ \quo_{\Gamma}\right)^{\ast}(\quot{A}{G}) \cong \left(\quo_{\Gamma}\right)^{\ast}\left(\of^{\ast}\big(\quot{A}{G}\big)\right) \\
	&\cong \left(\quo_{\Gamma}\right)^{\ast}(\AGamma),
	\end{align*}
	where the first isomorphism is given by the pesudofunctoriality of the assignment $\DbQl{-}$ and the second isomorphism is $\quo_{\Gamma}^{\ast}(\tau_f^A)$.
\end{proof}

\begin{proposition}\label{Prop: Section 4: Comparison of EDC with ABL EDC}
	There is a functor $\alpha:\DbeqQl{X} \to \Glue(\SfResl_G(X))$ induced on objects by
	\[
	\lbrace \AGamma \; | \; \Gamma \in \Sf(G)_0\rbrace \mapsto (\psi^{\ast}\quot{A}{G}, \AGamma, \theta_{\Gamma})
	\]
	and on morphisms by
	\[
	\lbrace \quot{\rho}{\Gamma} \; | \; \Gamma \in \Sf(G)_0 \rbrace \mapsto (\psi^{\ast}\quot{\rho}{G}, \quot{\rho}{\Gamma}).
	\]
\end{proposition}
\begin{proof}
	That the object assignment is defined (so in particular takes values in $\Glue(\SfResl_G(X))_0$) follows from Lemma \ref{Lemma: Section 4: Descent on object in DGbX is effective} and that the morphism assignment is defined (takes values in $\Glue(\SfResl_G(X))_1$) is a routine but tedious check using the compatibility between cocycle conditions described by Lemma \ref{Lemma: Section 4: Cocycle condition for equivariant derived cat}. This assignment is functorial because $\psi^{\ast}$ is itself a functor, which proves the proposition.
\end{proof}

Propositions \ref{Prop: Section 4: Comparison functor from ABL to equivariant cat} and \ref{Prop: Section 4: Comparison of EDC with ABL EDC} give us comparison functors 
\[
\beta:\Glue(\SfResl_G(X)) \to \DbQl{X}
\]
and 
\[
\alpha:\DbeqQl{X} \to \Glue(\SfResl_G(X)).
\] 
We will now prove the central result of this section, which shows that these form an equivalence of categories, which yields another proof that there is a $t$-structure on $\DbeqQl{X}$ whose heart is equivalent to the category of na{\"i}vely equivariant perverse sheaves.

\begin{Theorem}\label{Thm: Section 4: EDC is ABL EDC}
	The functors $\alpha, \beta$ of Propositions \ref{Prop: Section 4: Comparison functor from ABL to equivariant cat} and \ref{Prop: Section 4: Comparison of EDC with ABL EDC} induce an adjoint equivalence of categories:
	\[
	\begin{tikzcd}
	\DbeqQl{X} \ar[rr, bend left = 30, ""{name = U}]{}{\alpha} & & \Glue(\SfResl_G(X)) \ar[ll, bend left = 30, ""{name = L}]{}{\beta} \ar[from = U, to = L, symbol = \simeq]
	\end{tikzcd}
	\]
\end{Theorem}
\begin{proof}
	We will prove that $\alpha \circ \beta \cong \id_{\Glue(\SfResl_G(X))}$ and that $\beta \circ \alpha \cong \id_{\DbeqQl{X}}$ and then invoke the classical categorical fact that any equivalence can be promoted to an adjoint equivalence (cf.\@ \cite[Proposition 4.4.5]{riehl2017category}). For this fix an object
	\[
	(A,T_A) = \bigg(\lbrace \AGamma \; | \; \Gamma \in \Sf(G)_0 \rbrace, \lbrace \tau_f^A \; | \; f \in \Sf(G)_1 \rbrace\bigg)
	\] 
	of $\DbeqQl{X}$ and an object $(\Fscr, \Fscr_{\Gamma}, \theta_{\Gamma})$ of $\Glue(\SfResl_G(X))$.
	
	We begin by showing that $\alpha \circ \beta \cong \id_{\Glue(\SfResl_G(X))}$. For this it suffices to prove that 
	\[
	(\alpha \circ \beta)(\Fscr, \Fscr_{\Gamma}, \theta_{\Gamma}) \cong (\Fscr, \Fscr_{\Gamma}, \theta_{\Gamma}).
	\] 
	Begin by calculating that
	\begin{align*}
	(\alpha \circ \beta)(\Fscr, \Fscr_{\Gamma}, \theta_{\Gamma}) &= \alpha\left(\lbrace \Fscr_{\Gamma} \; | \; \Gamma \in \Sf(G)_0 \rbrace, \lbrace \tau_f^{\Fscr} \; | \; f \in \Sf(G)_1 \rbrace\right) \\
	&= \left(\psi^{\ast}\Fscr_G, \Fscr_{\Gamma}, \theta_{\Gamma}^{\prime}\right).
	\end{align*}
	It thus suffices to prove that there is an isomorphism $\Fscr \cong \psi^{\ast}\Fscr_G$; however, because $\psi:X \to \quot{X}{G}$ is a $G$-invariant (non-equivariant) isomorphism the cocycle condition allows us to derive isomorphisms
	\[
	(\pi_2^{G})^{\ast}\Fscr \cong (\quo_G)^{\ast}\Fscr_G \cong (\pi_2^{G})^{\ast}(\psi^{\ast}\Fscr_G)
	\]
	which allow us to conclude using the invariance of $\psi$ that $\Fscr \cong \psi^{\ast}\Fscr_G$.
	
	We now prove the other equivalence, i.e., that $(\beta \circ \alpha)(A,T_A) \cong (A,T_A)$. For this we note that
	\[
	(\beta \circ \alpha)(A,T_A) = \left(\lbrace \AGamma \; | \; \Gamma \in \Sf(G)_0 \rbrace, \lbrace \sigma_f^A \; | \; f \in \Sf(G)_1 \rbrace\right);
	\]
	to prove that there is an isomorphism of $A$ with $\beta(\alpha A)$, we need to show that there is an isomorphism $P:A \to \beta(\alpha A)$, $P = \lbrace \quot{\rho}{\Gamma} \; | \; \Gamma \in \Sf(G)_0\rbrace$, such that for any $f \in \Sf(G)_1$, say with $\Dom f = \Gamma $ and $\Codom f = \Gamma^{\prime}$, the diagram
	\[
	\xymatrix{
		\overline{f}^{\ast}\AGammap \ar[r]^{\overline{f}^{\ast}\quot{\rho}{\Gamma^{\prime}}} \ar[d]_{\tau_f^A} & \overline{f}^{\ast}\AGammap \ar[d]^{\sigma_f^A} \\
		\AGamma \ar[r]_{\quot{\rho}{\Gamma}} & \AGamma
	}
	\]
	commutes whenever $\Gamma\times X$ has an $n$-acyclic resolution to some variety $\Gamma^{\prime} \times X$ and for some $n \in \N$. However, from Lemmas \ref{Lemma: Section 4: Cocycle on ABL cat is cocycle in our cat} and \ref{Lemma: Section 4: Cocycle condition for equivariant derived cat} we note that the transition isomorphisms $\sigma_{f}^{A}$  are induced from the ABL-cocycle data on $(\psi^{\ast}\quot{A}{G},\AGamma,\theta_{\Gamma})$, which is in turn induced from the isomorphisms
	\[
	\xymatrix{
		(\overline{\pi}_1^{\Gamma\Gamma^{\prime}})^{\ast}\AGamma \ar[rr]^-{\varphi_{\Gamma\Gamma^{\prime}}} \ar[dr]_{\tau_{\pi_1}^{A}} & & (\overline{\pi}_2^{\Gamma\Gamma^{\prime}})^{\ast}\AGammap \\
		& \quot{A}{\Gamma \times \Gamma^{\prime}} \ar[ur]_{(\tau_{\pi_2}^A)^{-1}}
	}
	\]
	valid for all acyclic $\Gamma, \Gamma^{\prime} \in \Sf(G)_0$. This implies that the descent data induced by the $\tau_f$ and $\sigma_f$ are the same, which in turn allows us to induce our desired isomorphism $P$ by doing the stalk tricks that we did in Lemma \ref{Lemma: Section 4: Cocycle on ABL cat is cocycle in our cat}. Chasing this detail out proves that $A \cong \beta(\alpha A)$, which in turn shows that $\alpha$ and $\beta$ are inverse equivalences.
\end{proof}
\begin{Corollary}\label{Cor: Section 4: EDCs are equiv}
	There is an equivalence of catgories
	\[
	\DbeqQl{X} \simeq \ABLDbeqQl{X}.
	\]
\end{Corollary}
\begin{proof}
	Compose the equivalence of Theorem \ref{Thm: Section 4: EDC is ABL EDC} with the equivalence of Proposition \ref{Prop: Section 4: ABL derived cat is SfResl Glue Cat}.
\end{proof}

We now need a few straightforward observations about how this equivalence interacts with the standard $t$-structures on $\DbeqQl{X}$ and $\ABLDbeqQl{X}$ as well as the perverse $t$-structure on both $\DbeqQl{X}$ and $\ABLDbeqQl{X}$. The interaction is much more straightforward on the standard $t$-structure, so let us begin there. 

The standard $t$-structure on $\DbeqQl{X}$ is induced from saying that a triangle
\[
\xymatrix{
	A \ar[r] & B \ar[r] & C \ar[r] & A[1]
}
\]
is distinguished if and only if for all $\Gamma \in \Sf(G)_0$, the triangle
\[
\xymatrix{
	\AGamma \ar[r] & \BGamma \ar[r] & \quot{C}{\Gamma} \ar[r] & \AGamma[1]
}
\]
is distinguished in $\DbQl{\XGamma}$. Similarly, the standard $t$-structure on the category $\ABLDbeqQl{X}$ is induced by saying that a triangle\index{Standard $t$-structure!On ABLDbG@$\ABLDbeqQl{X}$}
\[
\xymatrix{
	(\Fscr,\Fscr_{U}, \theta_{U}) \ar[r] & (\Gscr, \Gscr_{U}, \theta_{U}) \ar[r] & (\Hscr, \Hscr_U, \theta_U) \ar[r] & (\Fscr[1], \Fscr_U[1], \theta_U[1])
}
\]
is distinguished if and only if the triangles
\[
\xymatrix{
	\Fscr_U \ar[r] & \Gscr_U \ar[r] & \Hscr_U \ar[r] & \Fscr_U[1]
}
\]
are distinguished in $\DbQl{G\backslash U}$ for all acyclic resolutions $U$ of $X$. Translating through the equivalence of Corollary \ref{Cor: Section 4: EDCs are equiv} and using the exactness of the functors $\alpha$ and $\beta$ gives an equivalence of the standard hearts. However, we can do better: These hearts are also equivalent to the category of na{\"i}vely equivariant ($\ell$-adic) sheaves.
\begin{definition}\index{Na{\"i}ve Equivariant Sheaf}\index[notation]{ShvNaive@$\Shv_G^{\text{na{\"i}ve}}$}
The category $\Shv_G^{\text{na{\"i}ve}}(X;\overline{\Q}_{\ell})$ of na{\"i}vely equivariant $\ell$-adic sheaves on $X$ is defined as follows:
\begin{itemize}
	\item Objects: Pairs $(\Fscr, \theta)$ where $\Fscr \in \Shv(X;\overline{\Q}_{\ell})$ and $\theta$ is an isomorphism
	\[
	\theta:\alpha_X^{\ast}\Fscr \xrightarrow{\cong} \pi_2^{\ast}\Fscr
	\]
	in $\Shv(G \times X; \overline{\Q}_{\ell})$ which satisfies the GIT cocycle condition (cf.\@ Definition \ref{Defn: GIT Cocycle}).
	\item Morphisms: A morphism $\rho:(\Fscr,\theta) \to (\Gscr,\sigma)$ is a morphism $\rho:\Fscr\to \Gscr$ of $\ell$-adic sheaves for which the diagram
	\[
	\xymatrix{
	\alpha_X^{\ast}\Fscr \ar[r]^{\theta} \ar[d]_{\alpha_X^{\ast}\rho} & \pi_2^{\ast}\Fscr \ar[d]^{\pi_2^{\ast}\rho} \\
	\alpha_X^{\ast}\Gscr \ar[r]_-{\sigma} & \pi_2^{\ast}\Gscr
	}
	\]
	commutes in $\Shv(G \times X; \overline{\Q}_{\ell})$.
	\item Composition and Identities: As in $\Shv(X;\overline{\Q}_{\ell})$.
\end{itemize}
\end{definition}
\begin{proposition}\label{Prop: Equiv of Standard t-structures on EDCs}
	There is an equivalence of categories
	\[
	{}^{\text{stand}}\DbeqQl{X}^{\heartsuit} \simeq \ABLDbeqQl{X}^{\heartsuit_{\text{stand}}} \simeq \Shv_G^{\text{na{\"i}ve}}(X;\overline{\Q}_{\ell})
	\]
\end{proposition}
\begin{proof}
In light of the equivalence ${}^{\text{stand}}\DbeqQl{X}^{\heartsuit} \simeq \ABLDbeqQl{X}^{\heartsuit_{\text{stand}}}$, it suffices to show that the category ${}^{\text{stand}}\DbeqQl{X}$ is equivalent to the category $\Shv_G^{\text{na{\"i}ve}}(X;\overline{\Q}_{\ell}).$ Before proving this, however, observe that for any object $A \in {}^{\text{stand}}\DbeqQl{X}^{\heartsuit}_0$ we have $\Gamma$-locally $\AGamma \in A$ determines an object in ${}^{\text{stand}}\DbQl{\XGamma}^{\heartsuit}$. By the equivalence $\Shv(\XGamma;\overline{\Q}_{\ell}) \simeq {}^{\text{stand}}\DbQl{\XGamma}^{\heartsuit}$, we can regard each $\AGamma$ as an $\ell$-adic sheaf on $\XGamma$. Note that in particular, this equivalence in the direction ${}^{\text{stand}}\DbQl{\XGamma} \to \Shv(\XGamma;\overline{\Q}_{\ell})$ is given by the cohomology functor $\quot{H}{\Gamma}^{0}$. In the interest of keeping track of this, we write $\quot{H^0}{\Gamma}:\DbQl{\XGamma} \to \Shv(\XGamma;\overline{\Q}_{\ell})$ and $[0]:\Shv(\XGamma;\overline{\Q}_{\ell}) \to \DbQl{\XGamma}$ for the quasi-inverse of $\quot{H^0}{\Gamma}$ which embeds the sheaf $\Fscr$ as a complex concentrated only in degree $0$.
	
We begin by setting some notational conventions. Let $\psi_X:\quot{X}{G} \to X$ be the non-equivariant isomorphism of varieties described earlier and note that since the inverse map $\psi^{-1}_X:X \to \quot{X}{G}$ is $G$-trivial, for any $\Gamma \in \Sf(G)_0$ there is a unique morphism $\overline{\Gamma}_X:\XGamma \to \quot{X}{G}$ making the diagram
\[
\begin{tikzcd}
\Gamma \times X \ar[d, swap]{}{\quo_{\Gamma}} \ar[r]{}{\pi_2} & X \ar[d, shift left = 1]{}{\psi_X^{-1}} \\
\XGamma \ar[r, swap]{}{\overline{\Gamma}_X} & \quot{X}{G} \ar[u, shift left = 1]{}{\psi_X}
\end{tikzcd}
\]
commute. Similarly, there is a non-equivariant isomorphism of varieties $\psi_{G \times X}:\quot{(G \times X)}{G} \to G \times X$ such that for any $\Gamma \in \Sf(G)_0$ there is a unique morphism of varieties $\overline{\Gamma}_{G \times X}:\quot{(G \times X)}{\Gamma} \to \quot{(G \times X)}{G}$ making the diagram
\[
\begin{tikzcd}
\Gamma \times G \times X \ar[d, swap]{}{\quo_{\Gamma}} \ar[r]{}{\pi_2} & G \times X \ar[d, shift left = 1]{}{\psi_{G \times X}^{-1}} \\
\quot{(G \times X)}{\Gamma} \ar[r, swap]{}{\overline{\Gamma}_{G \times X}} & \quot{(G \times X)}{G} \ar[u, shift left = 1]{}{\psi_{G \times X}}
\end{tikzcd}
\]
commute as well. Note also that these unique morphisms are natural in the sense that for any morphism $f:\XGamma \to \XGammap$ the induced diagram
\[
\xymatrix{
\XGamma \ar[rr]^-{f} \ar[dr]_{\overline{\Gamma}_X} & & \XGammap \ar[dl]^{\overline{\Gamma}^{\prime}_X} \\
 & \quot{X}{G}
}
\]
commutes by the uniqueness of the morphisms $\overline{\Gamma}_X$ and $\overline{\Gamma}^{\prime}_X$. It is worth noting also that both the morphisms $\overline{G}_X$ and $\overline{G}_{G \times X}$ are the corresponding identity morphisms on $\quot{X}{G}$ and $\quot{(G \times X)}{G}$, respectively.

Let us now proceed by defining our functor 
\[
E:\Shv_G^{\text{na{\"i}ve}}(X;\overline{\Q}_{\ell}) \to {}^{\text{stand}}\DbeqQl{X}.
\]
To do this first fix an object $(\Fscr,\theta) \in \Shv_G(X;\overline{\Q}_{\ell})_0$ and define the object set $E(\Fscr)$ by
\[
E(\Fscr) := \lbrace \overline{\Gamma}_X^{\ast}\big(\psi^{\ast}_X\Fscr\big)[0] \; | \; \Gamma \in \Sf(G)_0 \rbrace
\]
so that $\quot{E(\Fscr)}{\Gamma} = \overline{\Gamma}^{\ast}_X\big(\psi_X^{\ast}\Fscr\big).$ We now define the transition isomorphisms $T_{E(\Fscr)}$ as follows: Fix a $\Gamma \in \Sf(G)_0$ of the form where either $\Gamma = G$ or $\Gamma \not\cong G^n \times \Gamma^{\prime}$ for any $\Gamma^{\prime} \in \Sf(G)_0$ and any $n \geq 1$. Now consider the diagram of varieties
\[
\begin{tikzcd}
&\quot{(G \times X)}{\Gamma} \ar[rr]{}{\overline{\Gamma}_{G \times X}} \ar[dd, shift right = 1, swap]{}{\oalphaGamma} \ar[dd, shift left = 1]{}{\opiGamma} \ar[dl, swap]{}{\muGamma} \ar[dr]{}{\muGamma} & & \quot{(G \times X)}{G} \ar[rr]{}{\psi_{G \times X}} \ar[dd, swap, shift right = 1]{}{\quot{\overline{\alpha}_X}{G}} \ar[dd, shift left = 1]{}{\quot{\opi_2}{G}} & & G \times X \ar[dd, swap, shift right = 1]{}{\alpha_X} \ar[dd, shift left = 1]{}{\pi_2} \\
\quot{X}{\Gamma_c} \ar[dr, swap]{}{\aGamma} & & \quot{X}{\Gamma_c} \ar[dl]{}{\pGamma} \\
&\XGamma \ar[rr, swap]{}{\overline{\Gamma}_X} & & \quot{X}{G} \ar[rr, swap]{}{\psi_X} & & X
\end{tikzcd}
\]
Now observe that
\begin{align*}
\oalphaGamma^{\ast}\left(\overline{\Gamma}_X^{\ast}\big(\psi_X^{\ast}\Fscr\big)\right)[0] &= \overline{\Gamma}^{\ast}_{G \times X}\left(\quot{\overline{\alpha}_X}{G}^{\ast}\big(\psi_X^{\ast}\Fscr\big)\right)[0] \\
&= \overline{\Gamma}_{G \times X}^{\ast}\left(\psi_{G \times X}^{\ast}\big(\alpha_X^{\ast} \Fscr\big)\right)[0]
\end{align*}
and similarly
\[
\opiGamma^{\ast}\left(\overline{\Gamma}_X^{\ast}\big(\psi_X^{\ast}\Fscr\big)\right)[0] = \overline{\Gamma}_{G \times X}^{\ast}\left(\psi_{G \times X}^{\ast}\big(\pi_2^{\ast}\Fscr\big)\right)[0].
\]
This allows us to give an isomorphism
\[
\overline{\Gamma}_{G \times X}^{\ast}\left(\psi_{G \times X}^{\ast}\big(\theta\big)\right)[0]:\oalphaGamma^{\ast}\left(\overline{\Gamma}_X^{\ast}\big(\psi_X^{\ast}\Fscr\big)\right)[0] \xrightarrow{\cong} \opiGamma^{\ast}\left(\overline{\Gamma}_X^{\ast}\big(\psi_X^{\ast}\Fscr\big)\right)[0]
\]
which essentially tells us how we have to build our transitions. In fact, we induce our transition isomorphisms over $\Gamma$ by asserting that they are given so that the equation
\[
\muGamma^{\ast}\big(\tau_{\pGamma}^{A}\big)^{-1} \circ \muGamma^{\ast}\tau_{\aGamma}^{A} = \overline{\Gamma}_{G \times X}^{\ast}\left(\psi_{G \times X}^{\ast}\big(\theta\big)\right)[0]
\]
holds and then freely generating the $\theta_f^A$ for all morphisms $f \in \Sf(G)_1$ by using the diagrams
\[
\xymatrix{
	\XGamma \ar[rr]^-{f} \ar[dr]_{\overline{\Gamma}_X} & & \XGammap \ar[dl]^{\overline{\Gamma}^{\prime}_X} \\
	& \quot{X}{G}
}
\]
and
\[
\begin{tikzcd}
 & \quot{(G \times X)}{G} \\
\quot{(G \times X)}{\Gamma} \ar[ur]{}{\overline{\Gamma}_{G \times X}} \ar[dd, shift left = 1]{}{\opiGamma} \ar[dd, shift right = 1, swap]{}{\oalphaGamma} \ar[rr]{}{\quot{\of}{G \times X}} \ar[dr]{}{\muGamma} & & \quot{(G \times X)}{\Gamma^{\prime}} \ar[ul, swap]{}{\overline{\Gamma}^{\prime}_{G \times X}}\ar[dr]{}{\muGamma} \\
 & \quot{X}{\Gamma_c} \ar[rr, swap, near start]{}{\of_c} \ar[dl, shift left = 1]{}{\pGamma} \ar[dl, shift right = 1, swap]{}{\aGamma} & & \quot{X}{\Gamma_c} \ar[dl, swap, shift right = 1]{}{\aGammap} \ar[dl, shift left = 1]{}{\pGammap} \\
\XGamma \ar[dr, swap]{}{\overline{\Gamma}_X} \ar[rr, swap]{}{\fX} & & \XGammap \ar[dl]{}{\overline{\Gamma}^{\prime}_X} \\
& \quot{X}{G}
 \ar[from = 2-3, to = 4-3, crossing over, shift left = 1, near start]{}{\opiGammap} \ar[from = 2-3, to = 4-3, crossing over, shift right = 1, near start, swap]{}{\oalphaGammap} 
\end{tikzcd}
\]
to produce the necessary diagrams. Note that the cocycle condition is satisfied by using the GIT cocycle condition to show it must hold.

To define $E$ on morphisms we send a map $\rho:(\Fscr,\theta) \to (\Gscr,\sigma)$ to the map
\[
E(\rho) = \lbrace \overline{\Gamma}_X^{\ast}\big(\psi_X^{\ast}(\rho)\big)[0] \; | \; \Gamma \in \Sf(G)_0 \rbrace.
\]
That this is indeed a morphism in ${}^{\text{stand}}\DbeqQl{X}^{\heartsuit}$ is routine by the conditions we use to define $T_{E(\Fscr)}$ and the fact that the diagram
\[
\xymatrix{
\alpha_X^{\ast}\Fscr \ar[d]_{\alpha_X^{\ast}\rho} \ar[r]^-{\theta} & \pi_2^{\ast}\Fscr \ar[d]^{\pi_2^{\ast}\rho} \\
\alpha_X^{\ast}\Gscr \ar[r]_-{\sigma} & \pi_2^{\ast}\Gscr
}
\]
commutes in $\Shv_G(G \times X;\overline{\Q}_{\ell})$. That this is a functor is immediate from the fact that the morphisms $E(\rho)$ are induced by functors as well.

We now define the functor
\[
F:{}^{\text{stand}}\DbeqQl{X}^{\heartsuit} \to \Shv_G^{\text{na{\"i}ve}}(X;\overline{\Q}_{\ell})
\]
which we will show is an inverse equivalence to $E$. Define the functor $F$ by sending an object $A$ of ${}^{\text{stand}}\DbeqQl{X}^{\heartsuit}$ to the $\ell$-adic sheaf
\[
F(A) = (\psi^{-1}_X)^{\ast}\left(\quot{H^0}{G}\big(\quot{A}{G}\big)\right).
\]
To define the equivariance $\theta$ on $FA$ we calculate that
\[
\alpha_X^{\ast}(FA) = \alpha_X^{\ast}\left(\big(\psi_X^{-1}\big)^{\ast}\big(\quot{H^0}{G}(\quot{A}{G})\big)\right) = (\psi_{G \times X}^{-1})^{\ast}\left(\quot{\overline{\alpha}_X}{G}^{\ast}\big(\quot{H^0}{G}(\quot{A}{G})\big)\right)
\]
while
\[
\pi_2^{\ast}(FA) = (\psi_{G \times X}^{-1})^{\ast}\left(\quot{\opi_2}{G}^{\ast}\big(\quot{H^0}{G}(\quot{A}{G})\big)\right),
\]
so we define the equivariance $\theta$ on $FA$ by
\[
\theta = (\psi_{G \times X}^{-1})^{\ast}\quot{\Theta_{H^0}(A)}{G}
\]
where $\quot{\Theta_{H^0}(A)}{G}$ is the $\quot{X}{G}$-component of the isomorphism $\Theta$ of Theorem \ref{Theorem: Formal naive equivariance} for the object $H^0(A)$; note that this satisfies the GIT cocycle condition by Theorem \ref{Thm: GIT Cocycle Condition Formal Equivariance}. We define $F$ on morphisms $P$ by
\[
FP = (\psi_X^{-1})^{\ast}\quot{\rho}{G}.
\]
That this morphism commutes with the equivariances $\theta_{FA}$ and $\theta_{FB}$ follows from the naturality condition at the end of Theorem \ref{Theorem: Formal naive equivariance}. finally, that this is a functor is immediately verified from the fact that $FP$ is the functorial image of a morphism in the collection $P$.

Observe that a careful reading of the proofs of Theorem \ref{Thm: Section 4: EDC is ABL EDC} and Corollary \ref{Cor: Section 4: EDCs are equiv} only involve exact functors, and hence standard-exact functors. Consequently, by using the inverse equivalences $\quot{H}{G}^{0}$ and $[0]$ we get that there is a natural isomorphism between $E \circ F$ and the identity functor by proceeding mutatis mutandis as in the proof of Theorem \ref{Thm: Section 4: EDC is ABL EDC}. We can recover the transition isomorphisms (up to equivalence) of the object $A$ from $(E \circ F)(A) $ by using the natural isomorphism
\begin{align*}
\muGamma^{\ast}(\tau_{\pGamma}^{E(F(A))})^{-1} \circ \muGamma^{\ast}\tau_{\aGamma}^{E(F(A))} &= \psi_{G \times X}^{\ast}\left((\psi_{G \times X}^{\ast})\big(\quot{H^0}{G}(\Theta_A)\big)\right)[0] \\
&\cong \quot{H^0}{G}(\Theta_A)[0] = \quot{\Theta_A}{G} \\
&= \muGamma^{\ast}(\pGamma^A)^{-1} \circ \muGamma^{\ast}\tau_{\aGamma}^{A}.
\end{align*}
This allows us to deduce that
\[
E \circ F \cong \id_{{}^{\text{stand}}\DbeqQl{X}^{\heartsuit}},
\]
which is one half of the desired equivalence.

We now only need to verify that for any $(\Fscr,\theta) \in \Shv_G^{\text{na{\"i}ve}}(X;\overline{\Q}_{\ell})_0$ $(\Fscr,\theta)$ is isomorphic to $(F \circ E)(\Fscr,\theta)$. For this, however, we compute that $F \circ E$ produces the underlying $\ell$-adic sheaf
\[
(F \circ E)(\Fscr) = F(E(\Fscr),T_{E(\Fscr)}) = (\psi_X^{-1})^{\ast}\big(\quot{H^0}{G}\big(\psi_X^{\ast}\Fscr\big)[0]\big) = (\psi_X^{-1})^{\ast}\big(\psi_X^{\ast}\Fscr\big)
\]
which is naturally isomorphic to $\Fscr$ via the map
\[
\nu_{\Fscr}:\Fscr \xrightarrow{\cong} (\psi_X^{-1})^{\ast}\big(\psi_X^{\ast}\Fscr\big)
\]
induced by the fact that $\Shv_G(-)$ is a pseudofunctor. That this isomorphism commutes with the equivariances comes from the fact that the diagram
\[
\xymatrix{
\alpha_X^{\ast}\Fscr \ar[r]^-{\alpha_X^{\ast}\nu_{\Fscr}}  \ar[d]_{\theta} & \alpha_X^{\ast}(F \circ E)\Fscr \ar[d]^{\theta_{(F \circ E)\Fscr}} \\
\pi_2^{\ast}\Fscr \ar[r]_-{\pi_2^{\ast}\nu_{\Fscr}} & \pi_2^{\ast}(F \circ E)\Fscr
}
\]
can be rewritten as the commuting diagram
\[
\begin{tikzcd}
\alpha_X^{\ast}\Fscr \ar[rr, bend left = 30]{}{\eta_{\alpha_X^{\ast}\Fscr}} \ar[r]{}{\alpha_X^{\ast}\nu_{\Fscr}} \ar[d, swap]{}{\theta} & \alpha_X^{\ast}(F \circ E)\Fscr \ar[d]{}{\theta_{(F \circ E)\Fscr}} \ar[r, equals] & (\psi_{G \times X}^{-1})^{\ast}\big(\psi_{G \times X}^{\ast}\big(\alpha_X^{\ast}\Fscr\big)\big) \ar[d]{}{(\psi_{G \times X}^{-1})^{\ast}\big(\psi_{G \times X}^{\ast}\theta\big)} \\
\pi_2^{\ast}\Fscr \ar[rr, bend right = 30, swap]{}{\eta_{\pi_2^{\ast}\Fscr}} \ar[r, swap]{}{\pi_2^{\ast}\nu_{\Fscr}} & \pi_2^{\ast}(F \circ E)\Fscr \ar[r, equals] & (\psi_{G \times X}^{-1})^{\ast}\big(\psi_{G \times X}^{\ast}\big(\pi_2^{\ast}\Fscr\big)\big)
\end{tikzcd}
\]
where $\eta$ is the induced natural isomorphism
\[
\eta:\id_{\Shv(G \times X;\overline{\Q}_{\ell})} \xRightarrow{\cong} (\psi_{G \times X}^{-1})^{\ast} \circ \psi_{G \times X}^{\ast}.
\]
Thus the map $\nu_{\Fscr}:(F \circ E)\Fscr \to \Fscr$ is an isomorphism
\[
(\Fscr,\theta) \cong \big((F \circ E)\Fscr, \theta_{(F \circ E)\Fscr}\big)
\]
and so we conclude that
\[
E \circ F \cong \id_{\Shv_G^{\text{na{\"i}ve}}(X;\overline{\Q}_{\ell})}.
\]
This provides the second half of our desired equivalence and completes the proof of the proposition.
\end{proof}
Because each of the functors in play here is induced by exact pullback functors, applying this proof to the case when each na{\"i}vely equivariant sheaf $(\Fscr,\theta)$ is in fact a na{\"i}vely equivariant local system $(\Lcal, \theta)$, we get the following equivalence of categories.
\begin{corollary}
There is an equivalence of categories
\[
\Loc_G(X;\overline{\Q}_{\ell}) \simeq \Loc_G^{\text{na{\"i}ve}}(X;\overline{\Q}_{\ell}).
\]
\end{corollary}
\begin{proof}
In light of Proposition \ref{Prop: Equiv of Standard t-structures on EDCs}, the only thing we need to check is that the functors $E$ and $F$ preserve the local trivialization aspects of local systems, i.e., that if $\Lcal \in \Loc_G(X;\overline{\Q}_{\ell})_0$, $F(\Lcal) \in \Loc_G^{\text{na{\"i}ve}}(X;\overline{\Q}_{\ell})$ and that if $(\Lcal, \theta) \in \Loc_G(X;\overline{\Q}_{\ell})_0$ then $E(\Lcal) \in \Loc_G(X;\overline{\Q}_{\ell})_0$. 

Doing this check for the functor $F$ is straightforward: If $\Lcal \in \Loc_G(X;\overline{\Q}_{\ell})_0$ then each object $\quot{\Lcal}{\Gamma} \in \Loc(\XGamma;\overline{\Q}_{\ell})_0$. Because the morphism $\psi_X$ is an isomorphism and hence {\'e}tale, the pullback functors $\psi_X^{\ast}$ and $(\psi_X^{-1})^{\ast}$ both preserve local systems. Thus 
\[
F(\Lcal) = (\psi_X^{-1})^{\ast}\quot{\Lcal}{G} \in \Loc_G^{\text{na{\"i}ve}}(X;\overline{\Q}_{\ell})_0,
\]
as was desired.

We now perform the check for $E$. Let $(\Lcal, \theta)$ be a na{\"i}vely equivariant local system. We must show that each sheaf
\[
\quot{E(\Lcal)}{\Gamma} = \overline{\Gamma}_X^{\ast}\big(\psi_X^{\ast}\Lcal\big)
\]
is a local system on $\XGamma$. However, we already know that $\psi_X^{\ast}\Lcal$ is a local system, so it suffices to show that $\overline{\Gamma}_X^{\ast}$ preserves local systems. However, since $\overline{\Gamma}_X$ is induced by the diagram
\[
\xymatrix{
\Gamma \times X \ar[r]^-{\pi_2} \ar[d]_{\quo_{\Gamma}} & X \ar[d]^{\psi_X^{-1}} \\
\XGamma \ar[r]_-{\overline{\Gamma}_X} & \quot{X}{G}
}
\]
where $\quo_{\Gamma}$ is surjective and smooth (by being a geometric quotient of a princiapl $G$-variety) and the map $\psi_X^{-1} \circ \pi_2$ is smooth by being a composition of smooth morphisms. Thus appealing to Proposition \ref{Prop: Section 2.1: Smooth quotient map} gives that $\overline{\Gamma}$ is smooth so $\overline{\Gamma}^{\ast}$ preserves local systems. Thus
\[
\quot{E(\Lcal)}{\Gamma} = \overline{\Gamma}_X^{\ast}\big(\psi_X^{\ast}\Lcal\big) \in \Loc(\XGamma;\overline{\Q}_{\ell})
\]
so it follows that $E(\Lcal) \in \Loc_G(X;\overline{\Q}_{\ell})$. This proves the corollary.
\end{proof}

A careful reading of the proof of Proposition \ref{Prop: Equiv of Standard t-structures on EDCs} shows that we only needed to use the truncation/embedding adjunction $\tau^{\geq 0}\circ\tau^{\leq 0} = H^0 \dashv [0]:\DbQl{X} \to \Shv(X;\overline{\Q}_{\ell})$ to make this proof work. As such if we replace the standard $t$-structure with the perverse $t$-structure and proceed mutatis mutandis together with the observation that isomorphisms are dimension $0$ morphisms so $\psi^{\ast} = {}^{p}\psi^{\ast}$, we derive the following theorem.
\begin{Theorem}\label{Theorem: Equivalence of hearts with Equivariant Perverse Sheaves}
There is an equivalence of categories
\[
\ABLDbeqQl{X}^{\heartsuit_{\text{per}}} \simeq \DbeqQl{X}^{\heartsuit_{\text{per}}} \simeq \Per_G(X;\overline{\Q}_{\ell}).
\]
\end{Theorem}
\begin{proof}[Sketch]
Let $\Gamma \in \Sf(G)_0$. Begin by writing
\[
\quot{{}^{p}H^0}{\Gamma}:\DbQl{\XGamma} \to \Per(\XGamma;\overline{\Q}_{\ell})
\]
for the degree $0$ perverse cohomology at $\Gamma$ and let
\[
{}^{p}[0]:\Per(\XGamma;\overline{\Q}_{\ell}) \to \DbQl{\XGamma}
\]
denote the degree $0$ inclusion of perverse sheaves on $\XGamma$ into the derived category; note that ${}^{p}[0]$ and $\quot{{}^{p}H^0}{\Gamma}$ are inverse equivalences when restricted to $\DbQl{\XGamma}^{\heartsuit_{\text{per}}}$. We then define the pervese analogue of the functor $E$ in Proposition \ref{Prop: Equiv of Standard t-structures on EDCs} by, for any $(\Pscr, \theta) \in \Per_G^{\text{na{\"i}ve}}(X;\overline{\Q}_{\ell})_0$,
\[
{}^{p}E(\Pscr) := \left\lbrace{}^{p}\overline{\Gamma}_X^{\ast}\left(\psi_X^{\ast} \Pscr\right){}^{p}[0] \; | \; \Gamma \in \Sf(G)_0 \right\rbrace
\]
and inducing the transition isomorphisms as was done in the proof of Proposition \ref{Prop: Equiv of Standard t-structures on EDCs}. Note that for all $\Gamma \in \Sf(G)_0$ each object is a perverse sheaf on $\XGamma$ because we are taking perverse pullbacks along $\overline{\Gamma}_X$ while because $\psi_X$ is an isomorphism (and hence a smooth morphism of relative dimension $0$), $\psi_X^{\ast} = \psi_X^{\ast}[0] = {}^{p}\psi_X^{\ast}$. Defining ${}^{p}E$ on morphisms is done similarly and gives rise to a functor
\[
{}^{p}E:\Per_G^{\text{na{\"i}ve}}(X;\overline{\Q}_{\ell}) \to \DbeqQl{X}^{\heartsuit_{\text{per}}}.
\]

Defining the perverse analogue of the  functor $F$ used in Proposition \ref{Prop: Equiv of Standard t-structures on EDCs} is straightforward: For any $A \in \DbeqQl{X}^{\heartsuit_{\text{per}}}_0$, we define 
\[
{}^{p}F(A) := (\psi_X^{-1})^{\ast}\left(\quot{{}^{p}H^0}{G}(\quot{A}{G})\right);
\]
that this is indeed a perverse sheaf follows from the fact that since $\psi_X^{-1}$ is an isomorphism (and hence a smooth morphism of relative dimension $0$), ${}^{p}(\psi^{-1}_X)^{\ast} = (\psi_X^{-1})^{\ast}[0] = (\psi_X^{-1})^{\ast}.$ The equivariance $\theta$ is built by taking
\[
\theta := (\psi_{G \times X}^{-1})^{\ast}\quot{\Theta_{{}^{p}H^0}(A)}{G}
\]
and proceeding as in Proposition \ref{Prop: Equiv of Standard t-structures on EDCs}. We omit the definition of ${}^{p}F$ on morphisms, as it is straightforward, and use that this gives rise to a functor 
\[
{}^{p}F:\DbeqQl{X}^{\heartsuit_{\text{per}}} \to \Per_G^{\text{na{\"i}ve}}(X;\overline{\Q}_{\ell}).
\]

Finally, to verify the isomorphisms 
\[
{}^{p}E \circ {}^{p}F \cong \id_{\DbeqQl{X}^{\heartsuit_{\text{per}}}}
\]
and 
\[
{}^{p}F \circ {}{^p}E \cong \id_{\Per_G^{\text{na{\"i}ve}}(X;\overline{\Q}_{\ell})},
\] 
simply proceed mutatis mutands as in the proof of Proposition \ref{Prop: Equiv of Standard t-structures on EDCs}.
\end{proof}

\section{The Category $\DbeqQl{X}$ and Desiderata \ref{Desiderata}}
We now close this section with a peripheral but important result: The equivariant derived category $\DbeqQl{X}$ satisfies Desiderata \ref{Desiderata}. Most of the various constructions are consequences of the constructions we performed in Sections \ref{Section 2: Cat theory of ECat}, \ref{Section: Section 3: Change of Fibre}, \ref{Section: Section 3: Change of Space}, and in Chapter \ref{Chapter 5}. However, because Desiderata relies on the existence of $t$-structures which have hearts equivalent to the categories of equivariant sheaves and equivariant perverse sheaves, which is the content of Proposition \ref{Prop: Equiv of Standard t-structures on EDCs} and Theorem \ref{Theorem: Equivalence of hearts with Equivariant Perverse Sheaves}, we have delayed the proof until here. Let us begin by collecting some of the scheme-theoretic results we need for proving $\DbeqQl{X}$ satisfies Desiderata \ref{Desiderata}.

\begin{lemma}[{\cite[\href{https://stacks.math.columbia.edu/tag/0F7L}{Tag 0F7L}]{stacks-project}}]\label{LEmma: Proper base change commutes iwth pullback}
	Consider a pullback diagram of quasi-compact and quasi-separated schemes
\[
\xymatrix{
	W \ar[r]^-{g} \pullbackcorner \ar[d]_{k} & Z \ar[d]^{f} \\
	Y \ar[r]_-{h} & X 
}
\]
where $f$ is a separated morphism of finite type. Then there is a natural isomorphism
\[
h^{\ast} \circ Rf_! \cong Rk_! \circ g^{\ast}.
\]
\end{lemma}

In order to apply this lemma to the varities and schemes we use, we need some short results saying that the diagrams we consider satisfy the conditions above.
\begin{lemma}\label{Lemma: Section ?: Obvious pullback square}
	For any $h \in \GVar(X,Y)$ and for any $f \in \Sf(G)(\Gamma,\Gamma^{\prime})$, the diagram
	\[
	\xymatrix{
		\Gamma \times X \ar[r]^{f \times \id_X} \ar[d]_{\id_{\Gamma} \times h} & \Gamma^{\prime} \times X \ar[d]^{\id_{\Gamma^{\prime}} \times h} \\
		\Gamma \times Y \ar[r]_{f \times \id_Y} & \Gamma^{\prime} \times Y
	}
	\]
	is a pullback.
\end{lemma}
\begin{proof}
	For any scheme $Z$ with morphisms $\psi = \langle \psi_\Gamma^{\prime}, \psi_X\rangle:Z \to \Gamma^{\prime} \times X$ and $\varphi = \langle \varphi_{\Gamma}, \varphi_{Y}\rangle:Z \to \Gamma \times Y$ making the induced outer square commute, simply set the unique map $\theta:Z \to \Gamma \times X$ to be given by $\theta := \langle \varphi_{\Gamma}, \psi_{X}\rangle$. From here the verification is routine.
\end{proof}
\begin{lemma}\label{Lemma: Section ?: We can use base change}
	For any $h \in \GVar(X,Y)$ and for any $f \in \Sf(G)(\Gamma,\Gamma^{\prime})$, the diagram
	\[
	\xymatrix{
		\XGamma \ar[r]^{\fX} \ar[d]_{\hGamma} & \XGammap \ar[d]^{\hGammap} \\
		\YGamma \ar[r]_{\fY} & \YGammap
	}
	\]
\end{lemma}
\begin{proof}
	Begin by considering the diagram of schemes
	\[
	\xymatrix{
		\XGamma \ar@/^/[drr]^{\fX} \ar@/_/[ddr]_{\hGamma} \ar@{-->}[dr]^{\exists!\theta}& & \\
		& \XGammap \times_{\YGammap} \YGamma \pullbackcorner \ar[d] \ar[r] & \XGammap \ar[d]^{\hGammap} \\
		& \YGamma \ar[r]_-{\fY} & \YGammap
	}
	\]
	produced by pulling back against $\fY$ and $\hGammap$. We claim that the morphism $\theta$ is an isomorphism. This is a tedious but routine result proved by realizing that $\XGammap \times_{\YGammap} \YGamma$ is the scheme which universally factors the quotient morphism $\overline{f \times h}$, where $f \times h$ is produced as the diagonal of the pullback square in Lemma \ref{Lemma: Section ?: Obvious pullback square}. However, because by construction we have that $\XGamma$ also universally factors this diagram from the fact that
	\[
	\hGammap \circ \fX = \overline{f \times h} = \fY \circ \hGamma,
	\]
	we get that $\XGamma$ is also the scheme which universally factors $\overline{f \times h}$. Thus it must be the case that $\theta$ is an isomorphism, as both objects must factor through each other.
\end{proof}

With these together we can finally prove the explicit existence of the six functors which make up the six functor yoga for the equivariant derived category.
\begin{proposition}\label{Prop: Existence of push/pull}
For any $G$-equivariant morphism $h:X \to Y$, there is an equivariant pullback/pushforward adjunction
\[
\begin{tikzcd}
	\DbeqQl{X} \ar[r, swap, bend right = 30, ""{name = L}]{}{Rh_{\ast}} & \DbeqQl{Y} \ar[l, bend right =30, swap, ""{name = U}]{}{h^{\ast}} \ar[from = U, to = L, symbol = \dashv]
\end{tikzcd}
\]
where $h^{\ast}$ is given on objects $A \in \DbeqQl{Y}_0$ by
\[
\lbrace \AGamma \; | \; \Gamma \in \Sf(G)_0 \rbrace \mapsto \lbrace \hGamma^{\ast}\AGamma \; | ; \Gamma \in \Sf(G)_0 \rbrace
\]
and $Rh_{\ast}$ is given on objects $B \in \DbeqQl{X}_0$ by
\[
\lbrace \BGamma\; | \; \Gamma \in \Sf(G)_0 \rbrace \mapsto \lbrace R(\hGamma)_{\ast}\BGamma \; | \; \Gamma \in \Sf(G)_0 \rbrace.
\]
\end{proposition}
\begin{proof}
The existence of $h^{\ast}$ follows directly from Theorem \ref{Thm: Section 3: Existence of pullback functors} applied to the simultaneous pre-equivariant pseudofunctor $\DbeqQl{-}$. A routine check involving the pre-equivariant pseudofunctor $\DbeqQl{-}$ and using the Smooth Base Change Theorem applied to the diagrams
\[
\xymatrix{
\XGamma \ar[r]^-{\fX} \ar[d]_{\hGamma} & \XGammap \ar[d]^{\hGammap} \\
\YGamma \ar[r]_-{\fY} & \YGammap
}
\]
for $f \in \Sf(G)(\Gamma, \Gamma^{\prime})$ (where we note that all maps in sight are between varieties and hence quasi-compact and separated, the $\fX$ and $\fY$ are smooth, and all maps in sight are finitely presented) shows that $\DbeqQl{-}$ admits descent pushforwards along $R(\hGamma)_{\ast}$; Theorem \ref{Thm: Section 3: Existence of pushforwards} thus proves the existence of $Rh_{\ast}$. Finally, using Theorem \ref{Theorem: Section 3.2: Descent adjunctions lift to equivariant ajdunction} together with the adjoints
\[
\begin{tikzcd}
	\DbQl{\XGamma} \ar[r, swap, bend right = 30, ""{name = L}]{}{R(\hGamma)_{\ast}} & \DbQl{\YGamma} \ar[l, bend right =30, swap, ""{name = U}]{}{(\hGamma)^{\ast}} \ar[from = U, to = L, symbol = \dashv]
\end{tikzcd}
\]
gives the desired adjunction:
\[
\begin{tikzcd}
	\DbeqQl{X} \ar[r, swap, bend right = 30, ""{name = L}]{}{Rh_{\ast}} & \DbeqQl{Y} \ar[l, bend right =30, swap, ""{name = U}]{}{h^{\ast}} \ar[from = U, to = L, symbol = \dashv]
\end{tikzcd}
\]
\end{proof}
\begin{proposition}\label{Prop: Existence of proper pushpull}
For any $G$-equivariant morphism $h:X \to Y$, there is an equivariant proper pushforward/pullback adjunction
\[
\begin{tikzcd}
	\DbeqQl{X} \ar[r, bend left = 30, ""{name = L}]{}{Rh_{!}} & \DbeqQl{Y} \ar[l, bend left =30,  ""{name = U}]{}{h^{\ast}} \ar[from = L, to = U, symbol = \dashv]
\end{tikzcd}
\]
where $h^{!}$ is given on objects $A \in \DbeqQl{Y}_0$ by
\[
\lbrace \AGamma \; | \; \Gamma \in \Sf(G)_0 \rbrace \mapsto \lbrace \hGamma^{!}\AGamma \; | ; \Gamma \in \Sf(G)_0 \rbrace
\]
and $Rh_{!}$ is given on objects $B \in \DbeqQl{X}_0$ by
\[
\lbrace \BGamma\; | \; \Gamma \in \Sf(G)_0 \rbrace \mapsto \lbrace R(\hGamma)_{!}\BGamma \; | \; \Gamma \in \Sf(G)_0 \rbrace.
\]
\end{proposition}
\begin{proof}
The commutativity $\hGamma^{!} \circ \fX^{\ast} \cong \fY^{\ast} \circ \hGammap^{!}$ for all $f \in \Sf(G)(\Gamma,\Gamma^{\prime})$ allows us to deduce that $\DbeqQl{-}$ admits descent pullbacks along the functors $\hGamma^{!}$. Thus Theorem \ref{Thm: Section 3.2: Existence of descent pullbacks} gives the existence of $h^{!}:\DbeqQl{Y} \to \DbeqQl{X}$.

To see that $Rh_!$ exists, note that for any $\Gamma \in \Sf(G)(\Gamma,\Gamma^{\prime})$ Lemma \ref{Lemma: Section ?: We can use base change} gives that the square
\[
\xymatrix{
\XGamma \ar[d]_{\fX} \ar[r]^-{\hGamma} & \YGamma \ar[d]^{\fY} \\
\XGammap \ar[r]_-{\hGammap} & \YGammap
}
\]
is a pullback. Moreover, each map is quasi-compact, separated, and finitely presented (as they are all morphisms between quasi-compact, finitely presented, and separated schemes) and the morphisms $\fX$ and $\fY$ are both smooth. Thus appealing to Lemma \ref{LEmma: Proper base change commutes iwth pullback} gives natural isomorphisms
\[
R(\hGamma)_! \circ \fX^{\ast} \cong \fY^{\ast} \circ R(\hGammap)_!
\]
that allow us to prove that $\DbeqQl{-}$ admits descent pushforwards along $Rh_!$. Applying Theorem \ref{Thm: Section 3: Existence of pushforwards} gives the existence of $Rh_!$. Finally, using the adjunction analogously to how it was used in the proof of Proposition \ref{Prop: Existence of push/pull} and appealing to Theorem \ref{Theorem: Section 3.2: Descent adjunctions lift to equivariant ajdunction} gives the desired adjunction:
	\[
\begin{tikzcd}
	\DbeqQl{X} \ar[r, bend left = 30, ""{name = L}]{}{Rh_{!}} & \DbeqQl{Y} \ar[l, bend left =30, ""{name = U}]{}{h^{!}} \ar[from = L, to = U, symbol = \dashv]
\end{tikzcd}
\]
\end{proof}
\begin{proposition}\label{Prop: EDC is symmetric monoidal closed}
The category $\DbeqQl{X}$ is symmetric moniodal closed with tensor functor given by
\[
A \os{L}{\otimes} B := \left\lbrace \AGamma \ds{L}{\XGamma}{\otimes} \BGamma \; : \; \Gamma \in \Sf(G)_0, \AGamma \in A, \BGamma \in B \right\rbrace
\]
and internal hom given by
\[
R[A,B] := \left\lbrace R[\AGamma, \BGamma]_{\Gamma} \; : \; \Gamma \in \Sf(G)_0, \AGamma \in A, \BGamma \in B \right\rbrace
\]
for all $A, B \in \DbeqQl{X}_0$.
\end{proposition}
\begin{proof}
We first define the symmetric monoidal functor $(-) \otimes^{L} (-):\DbeqQl{X} \times \DbeqQl{X} \to \DbeqQl{X}$ by, for $A, B \in \DbeqQl{X}_0$,
\[
A \os{L}{\otimes} B := \left\lbrace \AGamma \os{L}{\us{\Gamma}{\otimes}} \BGamma \; : \; \Gamma \in \Sf(G)_0, \AGamma \in A, \BGamma \in B \right\rbrace
\]
where the functor $(-)\os{L}{\us{\Gamma}{\otimes}}(-)$ is the (derived) tensor on $\DbQl{\XGamma}$. Note that this is the functor described in Example \ref{Example: Tensor Functor on EDC}, which is symmetric monoidal by Theorem \ref{Theorem: Section 2: Monoidal preequivariant pseudofunctor gives monoidal equivariant cat} and Proposition \ref{Prop: Section 2: Equivariant cat is symmetric monoidal}.

Let us now describe the internal hom functor 
\[
R[-,-]:\DbeqQl{X}^{\op} \times \DbeqQl{X} \to \DbeqQl{X}
\] 
by using Theorem \ref{Thm: Section 3: Psuedonatural trans lift to equivariant functors}, i.e., by describing a pseudonatural transformation $R[-,-]:\DbQl{-}^{\op} \times \DbQl{-} \Rightarrow \DbQl{-}$. We define the functor component
\[
R[-,-]_{\Gamma}:\DbQl{\XGamma}^{\op} \times \DbQl{\XGamma} \to \DbQl{\XGamma}
\]
to be the (derived) internal hom functor on $\DbQl{\XGamma}$. To define the coherence isomorphisms $\quot{R[-,-]}{f}$ for a morphism $f \in \Sf(G)(\Gamma, \Gamma^{\prime})$ we observe the following: By \cite[4.2.5.1, pg. 109]{BBD} for any smooth morphism $h:Z \to W$ of varieties, there is a natural isomorphism
\[
\nu_{(A,B)}:R[h^{\ast}A,h^{\ast}B]_{Z} \cong h^{\ast}\left(R[A,B]_{W}\right).
\]
In particular, using that the morphisms $\of:\XGamma \to \XGammap$ are smooth shows that each functor $\of^{\ast}$ is symmetric monoidal closed. In particular, it follows that
the pre-equivariant pseudofunctor $\DbQl{-}$ satisfies the hypotheses of Lemma \ref{Lemma: Section 3: Equivariant Cat is monoidal closed}; thus we conclude that $\DbeqQl{X}$ is symmetric monoidal closed.
\end{proof}

We now will prove the way in which the six equivariant functors interact with each other. In particular, we will show the relations between each of the functors and the existence of open/closed immersion distinguished triangles.

\begin{proposition}\label{Prop: Equivariant Grothendiekc Context}
For any morphism $h \in \GVar(X,Y)$ and any objects $A \in \DbeqQl{X}_0, B \in \DbeqQl{Y}$ there are natural isomorphisms
\begin{align*}
\big(Rh_!A\big) \os{L}{\otimes}_Y B &\cong Rh_!\left(A \os{L}{\otimes}_X h^{\ast}B\right), \\
R[Rh_!A,B]_Y &\cong Rh_{\ast}\left(R[A,h^{!}B]_X\right), \\
h^{!}R[A,B]_Y &\cong R[h^{\ast}A,h^{!}B]_X.
\end{align*}
\end{proposition}
\begin{proof}
This is a routine but tedious check involving the fact that for any $\Gamma \in \Sf(G)_0$ we have isomorphisms
\begin{align*}
\big(R(\hGamma)_!\AGamma\big) \ds{L}{\YGamma}{\otimes} \BGamma &\cong R(\hGamma)_!\left(\AGamma \ds{L}{\XGamma}{\otimes} (\hGamma)^{\ast}\BGamma\right), \\
R\left[R(\hGamma)_!\AGamma,\BGamma\right]_{\YGamma} &\cong R(\hGamma)_{\ast}\left(R[\AGamma,(\hGamma)^{!}\BGamma]_{\XGamma}\right), \\
(\hGamma)^{!}R[\AGamma,\BGamma]_{\YGamma} &\cong R[(\hGamma)^{\ast}\AGamma,h^{!}\BGamma]_{\XGamma}
\end{align*}
for any $\AGamma \in A$ and any $\BGamma \in B$ which commute suitably with the transition isomorphisms of each object (which itself follows from the construction of the functors).
\end{proof}

We now prove the open/closed exact triangles in $\DbeqQl{X}$. For this, however, we will need some basic results on how these triangles all exist on the $\Gamma$-local level. Our basic technique will be to show that immersions descend through $\Sf(G)$, i.e., that if we have an (open/closed) immersion of schemes $i:U \to V$ then the induced map $\quot{\overline{\imath}}{\Gamma}:\quot{U}{\Gamma} \to \quot{V}{\Gamma}$ is an (open/closed) immersion as well.
\begin{lemma}\label{Lemma: Section: Yet More Stuff: Closed immersions descend}
	Let $i:V \to X$ be a closed immersion of varieties. Then for any $\Gamma \in \Sf(G)_0$ the morphism
	\[
	\overline{(\id_{\Gamma} \times i)}:\quot{V}{\Gamma} \to \XGamma
	\]
	is a closed immersion.
\end{lemma}
\begin{proof}
	Recall that a map $\zeta:Y \to Z$ determines a closed immersion if and only if $\lvert \zeta \rvert(\lvert Y \rvert)$ is closed in $\lvert Z  \rvert$ and there is a short exact sequence of quasi-coherent sheaves
	\[
	\xymatrix{
		0 \ar[r] & \Iscr \ar[r] & \CalO_Z \ar[r]^-{\zeta^{\sharp}} &  \zeta_{\ast}\CalO_Y \ar[r] & 0
	}
	\]
	where $\Iscr$ is a sheaf of ideals on $Z$. We will first prove that we obtain the desired short exact sequence of quasi-coherent sheaves on $\XGamma$ and then show the set
	\[
	\lvert \overline{(\id_{\Gamma} \times i)} \rvert \left(\lvert \quot{V}{\Gamma} \rvert\right)
	\]
	is closed in $\lvert \XGamma \rvert$. Consider the commuting diagram
	\[
	\xymatrix{
		\Gamma \times V \ar[r]^-{\id_{\Gamma} \times i} \ar[d]_{\quo_{\Gamma \times V}} & \Gamma \times X \ar[d]^{\quo_{\Gamma \times X}} \\
		\quot{V}{\Gamma} \ar[r]_-{\overline{(\id_{\Gamma} \times i)}} & \XGamma
	}
	\]
	of varieties. Because the identity map is closed and $i$ is closed, $\id_{\Gamma} \times i$ is closed and so there is a short exact sequence of quasi-coherent sheaves
	\[
	\xymatrix{
		0 \ar[r] & \Kscr \ar[r] & \CalO_{\Gamma \times X} \ar[r] & (\id_{\Gamma} \times i)_{\ast}\CalO_{\Gamma \times V} \ar[r] & 0
	}
	\]
	on $\Gamma \times X$ where $\Kscr$ is a sheaf of ideals on $\Gamma \times X$. 
	We now apply the functor $(\quo_{\Gamma \times X})_{\ast}$ and use the fact that since $\quo_{\Gamma \times X}:\Gamma \times X \to \XGamma$ is a smooth morphism of varieties (and hence is a smooth separated quasi-compact map between Noetherian schemes) that the pushforward gives an exact functor
	\[
	(\quo_{\Gamma \times X})_{\ast}:\QCoh(\Gamma \times X) \to \QCoh(\XGamma).
	\]
	As such the sequence of quasi-coherent sheaves
	\[
	\xymatrix{
		0 \ar[r] & (\quo_{\Gamma \times X})_{\ast}\Iscr \ar[r] & (\quo_{\Gamma \times X})_{\ast}\CalO_{\Gamma \times X} \ar[r] & (\quo_{\Gamma \times X})_{\ast}\left((\id_{\Gamma} \times i)_{\ast}\CalO_{\Gamma \times V}\right) \ar[r] & 0
	}
	\]
	is short exact (and hence the sheaf $(\quo_{\Gamma \times X})_{\ast}\Iscr$ is a sheaf of ideals on $\XGamma$). Furthermore, because the quotient $\XGamma$ is a geometric quotient, the Zariski topology on $\XGamma$ coincides with the quotient topology. Thus there is an isomorphism
	\[
	(\quo_{\Gamma \times X})_{\ast}\CalO_{\Gamma \times X} \cong \CalO_{\XGamma}
	\]
	so we get that our short exact sequence above is isomorphic to the short exact sequence:
	\[
	\xymatrix{
		0 \ar[r] & (\quo_{\Gamma \times X})_{\ast}\Iscr \ar[r] & \CalO_{\XGamma} \ar[d] \\ 
		& 0 & \ar[l] (\quo_{\Gamma \times X})_{\ast}\left((\id_{\Gamma} \times i)_{\ast}\CalO_{\Gamma \times V}\right) 
	}
	\]
	Now from
	\begin{align*}
		(\quo_{\Gamma \times X})_{\ast}\left((\id_{\Gamma} \times i)_{\ast}\CalO_{\Gamma \times V}\right) &\cong \left(\quo_{\Gamma \times X} \circ (\id_\Gamma \times i)\right)_{\ast}\CalO_{\Gamma \times V} \\
		& \cong \left(\overline{(\id_{\Gamma} \times i)} \circ \quo_{\Gamma \times V}\right)_{\ast}\CalO_{\Gamma \times V} \\ 
		& \cong \left(\overline{(\id_{\Gamma} \times i)}\right)_{\ast}\left((\quo_{\Gamma \times V})_{\ast}\CalO_{\Gamma \times V}\right) \\
		&\cong \left(\overline{(\id_{\Gamma} \times i)}\right)_{\ast}\CalO_{\quot{V}{\Gamma}},
	\end{align*}
	where the last isomorphism follows mutatis mutandis to $(\quo_{\Gamma \times X})_{\ast}\CalO_{\Gamma \times X} \cong \CalO_{\XGamma}$, we find that the short exact sequence is isomorphic to the short exact sequence
	\[
	\xymatrix{
		0 \ar[r] & (\quo_{\Gamma \times X})_{\ast}\Iscr \ar[r] & \CalO_{\XGamma} \ar[r] & \left(\overline{(\id_{\Gamma} \times i)}\right)_{\ast}\CalO_{\quot{V}{\Gamma}} \ar[r] & 0
	}
	\]
	of quasi-coherent sheaves on $\XGamma$.
	
	We now argue that $\lvert \overline{(\id_{\Gamma} \times i)} \rvert \left(\lvert \quot{V}{\Gamma} \rvert\right)$ is closed in $\lvert \XGamma\rvert$. For this we recall that since the topology on $\XGamma$ is a quotient topology, for any set $W \subseteq \lvert \XGamma \rvert,$ $W$ is closed (open, respectively) if and only if $\lvert \quo_{\Gamma \times X}\rvert^{-1}(W)$ is closed (open, respectively) in $\lvert \Gamma \times X \rvert$. As such consider that
	\[
	\lvert \overline{(\id_{\Gamma} \times i)} \rvert \left(\lvert \quot{V}{\Gamma} \rvert\right) = \lbrace [\gamma, i(v)]_G \; : \; \gamma \in \lvert\Gamma\rvert, v \in \lvert V \rvert \rbrace
	\]
	and note
	\begin{align*}
		&\lvert \quo_{\Gamma \times X}\rvert^{-1}\left(\lvert \overline{(\id_{\Gamma} \times i)} \rvert \left(\lvert \quot{V}{\Gamma} \rvert\right)\right) \\
		&= \left\lbrace (\gamma, x) \in \lvert \Gamma \times X \rvert \; : \; (\gamma, x) \in \lvert \overline{(\id_{\Gamma} \times i)} \rvert \left(\lvert \quot{V}{\Gamma} \rvert\right)\right\rbrace \\
		&= \left\lbrace (\gamma, x) \; | \; \exists\, v\in \lvert V \rvert.\, [\gamma, x]_G = [\gamma,i(v)]_G \right\rbrace = \lvert \Gamma \times V \rvert \\
		&= \lvert\id_{\Gamma} \times i\rvert \left(\lvert\Gamma \times V \rvert\right).
	\end{align*}
	Because $\id_{\Gamma} \times i$ is closed, the above preimage is closed as well. This shows that $\lvert \overline{(id_{\Gamma} \times i)}\rvert(\lvert \quot{V}{\Gamma}\rvert)$ is closed in $\lvert \XGamma\rvert$. Thus $\overline{(\id_{\Gamma} \times i)}$ is a closed immersion.
\end{proof}
\begin{lemma}\label{Lemma: Section: Yet More Stuff: Open immerions descend}
	Let $j:U \to X$ be an open immersion of varieties and let $\Gamma \in \Sf(G)_0$. The map
	\[
	\overline{(\id_{\Gamma} \times i)}:\quot{U}{\Gamma} \to \XGamma
	\]
	is an open immersion of varieties.
\end{lemma}
\begin{proof}
	That the subset $\lvert\overline{(\id_{\Gamma} \times j)}\rvert(\lvert \quot{U}{\Gamma}\rvert)$ is open in $\lvert \XGamma \rvert$ is dual to the fact that the image $\lvert\overline{(\id_{\Gamma} \times i)}\rvert\left(\lvert \quot{V}{\Gamma}\rvert\right)$ in Lemma \ref{Lemma: Section: Yet More Stuff: Closed immersions descend} is closed and hence is omitted. By \cite[Proposition 4.2.2.a]{EGA1} we only need to show that for any point $[\gamma, u]_G\in \lvert \quot{U}{\Gamma} \rvert$ there is an isomorphism of local rings
	\[
	(\overline{(\id_{\Gamma} \times j)})^{\ast}_{[\gamma,u]_G}:\CalO_{\XGamma,[\gamma,j(u)]_G} \xrightarrow{\cong} \CalO_{\quot{U}{\Gamma}, [\gamma,u]_G}.
	\]
	For this we note that because $\id_{\Gamma} \times j$ is itself an open immersion, for each $(\gamma,u) \in \lvert \Gamma \times U\rvert$ there is an isomorphism of local rings
	\[
	(\id_{\Gamma} \times j)^{\ast}_{(\gamma,u)}:\CalO_{\Gamma \times X, (\gamma, j(u))} \xrightarrow{\cong} \CalO_{\Gamma \times U, (\gamma,u)}.
	\]
	Using the smooth surjectivity of the quotient morphisms, taking stalks, and the commuting diagram
	\[
	\xymatrix{
		\Gamma \times U \ar[d]_{\quo_{\Gamma \times U}} \ar[r]^-{\id_{\Gamma} \times j} & \Gamma \times X \ar[d]^{\quo_{\Gamma \times X}} \\
		\quot{U}{\Gamma} \ar[r]_-{\overline{(\id_{\Gamma} \times j)}} & \XGamma
	}
	\]
	we find that the isomorphism $(\id_{\Gamma} \times j)^{\ast}_{(j,u)}$ descends to an isomorphism of local rings based induced by the stalks at the point $\quo_{\Gamma \times U}(\gamma,u) = [\gamma,u]_G$:
	\[
	\rho_{[\gamma,u]_G}:(\quo_{\Gamma \times X})_{\ast}\left(\CalO_{\Gamma \times X, (\gamma,j(u))}\right) \xrightarrow{\cong} (\quo_{\Gamma \times U})_{\ast}\left(\CalO_{\Gamma \times U, (\gamma, u)}\right)
	\]
	which is equal to
	\[
	\rho_{[\gamma, u]_{G}}:\left((\quo_{\Gamma \times X})_{\ast}\CalO_{\Gamma \times X}\right)_{[\gamma, j(u)]_G} \xrightarrow{\cong} \left((\quo_{\Gamma \times U})_{\ast}\CalO_{\Gamma \times U}\right)_{[\gamma,u]_G}.
	\]
	Using the isomorphisms $(\quo_{\Gamma \times X})_{\ast}\CalO_{\Gamma \times X} \cong \CalO_{\XGamma}$ and $(\quo_{\Gamma \times U})_{\ast}\CalO_{\Gamma \times U} \cong \CalO_{\quot{U}{\Gamma}}$ we then find $\rho_{[\gamma, u]_G}$ fits into a commuting diagram 
	\[
	\xymatrix{
		\left((\quo_{\Gamma \times X})_{\ast}\CalO_{\Gamma \times X}\right)_{[\gamma, j(u)]_G} \ar[rr]^-{\rho_{[\gamma,u]_G}} \ar[d]_{\cong} & & \left((\quo_{\Gamma \times U})_{\ast}\CalO_{\Gamma \times U}\right)_{[\gamma,u]_G} \ar[d]^{\cong} \\
		\CalO_{\XGamma, [\gamma, j(u)]_G} \ar[rr]_-{\overline{(\id_{\Gamma} \times j)}^{\ast}_{[\gamma, j(u)]_G}} &  & \CalO_{\quot{U}{\Gamma}, [\gamma, u]_G}
	}
	\]
	where the bottom map is $\overline{(\id_{\Gamma} \times j)}_{[\gamma,u]_G}^{\ast}$ by construction of taking these particular stalks and passing through the various isomorphisms. From a diagram chase it follows that $\overline{(\id_{\Gamma} \times j)}_{[\gamma,u]_G}^{\ast}$ is an isomorphism and hence that $\overline{(\id_{\Gamma} \times j)}$ is an open immersion.
\end{proof}
\begin{lemma}\label{Lemma: Yet More Stuff: Immersions}
	Let $h:W \to X$ be an immersion of varieties and let $\Gamma \in \Sf(G)_0$. Then the map
	\[
	\overline{(\id_{\Gamma} \times h)}:\quot{W}{\Gamma} \to \XGamma
	\]
	is an immersion of schemes as well.
\end{lemma}
\begin{proof}
	If $h = i \circ j$ for $i$ closed and $j$ open, applying Lemmas \ref{Lemma: Section: Yet More Stuff: Open immerions descend} and \ref{Lemma: Section: Yet More Stuff: Closed immersions descend} together with the observation
	\[
	\overline{(\id_{\Gamma} \times h)} = \overline{(\id_{\Gamma} \times (i \circ j)} = \overline{(\id_{\Gamma} \times i)} \circ \overline{(\id_{\Gamma} \times j)}
	\]
	gives the result. In the other case, where $h = j \circ i$ for $j$ open and $i$ closed, we simply apply Lemmas \ref{Lemma: Section: Yet More Stuff: Closed immersions descend} and \ref{Lemma: Section: Yet More Stuff: Open immerions descend} in the opposite order of the prior case.
\end{proof}

\begin{proposition}\label{Prop: Exact open/closed triangle}
Let $j:U \to X$ and be a $G$-equivariant open immersion and let $i:V \to X$ be the corresponding closed immersion. Then for any object $A \in \DbeqQl{X}_0$ there is a distinguished triangle of the form
\[
\xymatrix{
Rj_!\left(j^{!}A\right) \ar[r]^-{\epsilon_A} & A \ar[r]^-{\eta_A} & Ri_{\ast}\left(i^{\ast}A\right) \ar[r] & Rj_!\left(j^!A\right)[1]
}
\]
where $\epsilon:Rj_{!} \circ j^{!} \Rightarrow \id_{\DbeqQl{X}}$ is the counit of adjunction and $\eta:\id_{\DbeqQl{X}} \Rightarrow Ri_{\ast} \circ i^{\ast}$ is the unit of adjunction.
\end{proposition}
\begin{proof}
By Theorem \ref{Theorem: Section Triangle: Equivariant triangulation} it suffices to prove that for any $\Gamma \in \Sf(G)_0$ there is a distinguished triangle of the form
\[
\xymatrix{
	R\left(\quot{\overline{\jmath}}{\Gamma}\right)_!\left(\left(\quot{\overline{\jmath}}{\Gamma}\right)^{!}\AGamma\right) \ar[r]^-{\epsilon_{\AGamma}} & \AGamma \ar[r]^-{\eta_{\AGamma}} & R\left(\quot{\overline{\imath}}{\Gamma}\right)_{\ast}\left(\left(\quot{\overline{\imath}}{\Gamma}\right)^{\ast}\AGamma\right) \ar[d] \\
	& &  R\left(\quot{\overline{\jmath}}{\Gamma}\right)_!\left(\left(\quot{\overline{\jmath}}{\Gamma}\right)^!\AGamma\right)[1]
}
\]
in each category $\DbQl{\XGamma}$. We observe first that such an distinguished triangle exists in each derived category $\DbQl{Z}$, for a variety $Z$, provided that the closed immersion $i$ is the closed complement of the open immersion $j$. Thus, since each category $\XGamma$ is a variety and $\quot{\overline{\jmath}}{\Gamma}$ is open by Lemma \ref{Lemma: Section: Yet More Stuff: Open immerions descend} (and similarly $\quot{\overline{\imath}}{\Gamma}$ is closed by Lemma \ref{Lemma: Section: Yet More Stuff: Closed immersions descend}), we only need to prove that the morphisms $\quot{\overline{\imath}}{\Gamma}$ and $\quot{\overline{\jmath}}{\Gamma}$ remain complementary to each other, i.e., $\quot{V}{\Gamma}$ is the closed complement to $\quot{U}{\Gamma}$.

To show that $\quot{U}{\Gamma}$ is the closed complement of $\quot{U}{\Gamma}$ it suffices to prove that the square
\[
\xymatrix{
\quot{U}{\Gamma} \times_{\XGamma} \quot{V}{\Gamma} \pullbackcorner \ar[rr]^-{p_2} \ar[d]_{p_1} & & \quot{V}{\Gamma} \ar[d]^{\quot{\overline{\imath}}{\Gamma}} \\
\quot{U}{\Gamma} \ar[rr]_-{\quot{\overline{\jmath}}{\Gamma}} & & \XGamma
}
\]
is both a pullback and pushout square, as since $\quot{\overline{\jmath}}{\Gamma}$ and $\quot{\overline{\imath}}{\Gamma}$ are immersions, the underlying space of the pushout taken over their pullback (intersection) is
\[
\left\lvert \quot{U}{\Gamma} \coprod_{\quot{U}{\Gamma} \times_{\XGamma} \quot{V}{\Gamma}} \quot{V}{\Gamma} \right\rvert = \left\lvert \quot{U}{\Gamma}\right\rvert \cup \left\lvert \quot{V}{\Gamma} \right\rvert,
\]
where we identify the underlying space of each variety with its image in $\XGamma$. To do this assume we have a diagram of the form
\[
\xymatrix{
\quot{U}{\Gamma} \times_{\XGamma} \quot{V}{\Gamma} \pullbackcorner \ar[rr]^-{p_2} \ar[d]_{p_1} & & \quot{V}{\Gamma} \ar[d]_{\quot{\overline{\imath}}{\Gamma}} \ar@/^/[ddr]^{k} \\
\quot{U}{\Gamma} \ar[rr]_-{\quot{\overline{\jmath}}{\Gamma}} \ar@/_/[drrr]_-{h} & & \XGamma \\
 & & & Z
}
\]
for some scheme $Z$. To construct a unique morphism $\alpha:\XGamma \to Z$, it suffices to construct a map along the underlying points of the space. Define a map $\alpha:\XGamma \to Z$ by
\[
\alpha[\gamma,x]_G := \begin{cases}
h[\gamma,x]_G & \text{if}\, x \in \lvert\quot{U}{\Gamma}\rvert \\
k[\gamma,x]_G & \text{if}\, x \in \lvert\quot{V}{\Gamma}\rvert.
\end{cases}
\]
To see that this is well-defined, it suffices to show that $[\gamma,u]_G \ne [\gamma^{\prime},v]_G$ for any $u \in \lvert \quot{U}{\Gamma}\rvert$ and any $v \in \lvert \quot{V}{\Gamma}\rvert$. So we calculate that because the fibre above the topological points of the quotients are orbit spaces of $\Gamma \times X$, $[\gamma,u]_G = [\gamma^{\prime},v]_G$ exactly when $g(\gamma,u) = (\gamma^{\prime},v)$ for some $g \in \lvert G \rvert$. But then, because the action of $G$ on $\Gamma \times X$ is diagonal, we have that $g\gamma = \gamma^{\prime}$ and $gu = v$. Using that the immersion $j:U \to X$ is a $G$-equivariant open immersion implies that $G\lvert U \rvert \subseteq \lvert U \rvert$; thus, there is no $g \in \lvert G \rvert$ for which $gu = v$ and so we cannot have that $[\gamma,u]_G = [\gamma,v]_G$. This shows that $\alpha$ is well-defined and hence that the diagram
\[
\xymatrix{
	\quot{U}{\Gamma} \times_{\XGamma} \quot{V}{\Gamma} \pullbackcorner \ar[rr]^-{p_2} \ar[d]_{p_1} & & \quot{V}{\Gamma} \ar[d]_{\quot{\overline{\imath}}{\Gamma}} \ar@/^/[ddr]^{k} \\
	\quot{U}{\Gamma} \ar[rr]_-{\quot{\overline{\jmath}}{\Gamma}} \ar@/_/[drrr]_-{h} & & \XGamma \ar[dr]^{\alpha} \\
	& & & Z
}
\]
commutes. That this map is unique follows from the fact that $\lvert \quot{U}{\Gamma} \rvert$ and $\lvert \quot{V}{\Gamma} \rvert$ are disjoint and together cover $\lvert \XGamma \rvert$, so any map defined on $\XGamma$ which must agree with $k$ on $\quot{V}{\Gamma}$ and $h$ on $\quot{U}{\Gamma}$ must agree everywhere and hence be the same map. Thus it follows that the diagram takes the form
\[
\xymatrix{
	\quot{U}{\Gamma} \times_{\XGamma} \quot{V}{\Gamma} \pullbackcorner \ar[rr]^-{p_2} \ar[d]_{p_1} & & \quot{V}{\Gamma} \ar[d]_{\quot{\overline{\imath}}{\Gamma}} \ar@/^/[ddr]^{k} \\
	\quot{U}{\Gamma} \ar[rr]_-{\quot{\overline{\jmath}}{\Gamma}} \ar@/_/[drrr]_-{h} & & \XGamma \ar@{-->}[dr]^{\exists!\alpha} \\
	& & & Z
}
\]
and hence the square
\[
\xymatrix{
	\quot{U}{\Gamma} \times_{\XGamma} \quot{V}{\Gamma} \pullbackcorner \ar[rr]^-{p_2} \ar[d]_{p_1} & & \quot{V}{\Gamma} \ar[d]^{\quot{\overline{\imath}}{\Gamma}} \\
	\quot{U}{\Gamma} \ar[rr]_-{\quot{\overline{\jmath}}{\Gamma}} & & \pushoutcorner \XGamma
}
\]
is a pushout-pullback square. This gives that $\quot{U}{\Gamma}$ is the open complement of $\quot{V}{\Gamma}$ and hence implies the existence of the distinguished triangle
\[
\xymatrix{
	R\left(\quot{\overline{\jmath}}{\Gamma}\right)_!\left(\left(\quot{\overline{\jmath}}{\Gamma}\right)^{!}\AGamma\right) \ar[r]^-{\epsilon_{\AGamma}} & \AGamma \ar[r]^-{\eta_{\AGamma}} & R\left(\quot{\overline{\imath}}{\Gamma}\right)_{\ast}\left(\left(\quot{\overline{\imath}}{\Gamma}\right)^{\ast}\AGamma\right) \ar[d] \\
	& &  R\left(\quot{\overline{\jmath}}{\Gamma}\right)_!\left(\left(\quot{\overline{\jmath}}{\Gamma}\right)^!\AGamma\right)[1]
}
\]
in $\DbQl{\XGamma}$. Because these triangles are distinguished for every $\Gamma$, we get that the triangle
\[
\xymatrix{
	Rj_!\left(j^{!}A\right) \ar[r]^-{\epsilon_A} & A \ar[r]^-{\eta_A} & Ri_{\ast}\left(i^{\ast}A\right) \ar[r] & Rj_!\left(j^!A\right)[1]
}
\]
is distinguished in $\DbeqQl{X}$, as was desired.
\end{proof}

We now prove that the forgetful functor of Example \ref{Example: Section 2: The forgetful functor} and Proposition \ref{Prop: Section 3.3: Forgetful functor is de-equivariantification functor} and the cohomological functors $H^0_G:\DbeqQl{X} \to \Shv_G(X;\overline{\Q}_{\ell})$, ${}^{p}H_G^0:\DbeqQl{X} \to \Per_G(X;\overline{\Q}_{\ell})$ commute up to isomorphism. Because of the equivalences of Proposition \ref{Prop: Equiv of Standard t-structures on EDCs} and Theorem \ref{Theorem: Equivalence of hearts with Equivariant Perverse Sheaves} it follows that we have invertible $2$-cells
\[
\begin{tikzcd}
\DbeqQl{X} \ar[rr, bend left = 30, ""{name = Up}]{}{\Forget} \ar[d, swap]{}{H_G^0} \ar[r, ""{name = LU}]{}{\DbQl{1_G^{\sharp}}} & D^b_{\Spec K}(X;\overline{\Q}_{\ell}) \ar[d]{}{H_{\Spec K}^{0}} \ar[r, ""{name = RU}]{}{\overline{D}^b} & \DbQl{X} \ar[d]{}{H^0} \\
\Shv_G(X;\overline{\Q}_{\ell}) \ar[r, swap, ""{name = LL}]{}{\Shv(1_G^{\sharp})} \ar[rr, swap, bend right = 30, ""{name = Down}]{}{\Forget} & \Shv_{\Spec K}(X;\overline{\Q}_{\ell}) \ar[r, swap, ""{name = RL}]{}{\overline{\Shv}} & \Shv(X;\overline{\Q}_{\ell}) \ar[from = LU, to = LL, Rightarrow, shorten <= 4pt, shorten >= 4pt]{}{\cong} \ar[from = RU, to = RL, Rightarrow, shorten <= 4pt, shorten >= 4pt]{}{\cong} \ar[from = Up, to = 1-2, Rightarrow, shorten >= 4pt, shorten <= 4pt]{}{\cong} \ar[from = 2-2, to = Down, Rightarrow, shorten <= 4pt, shorten >= 4pt]{}{\cong}
\end{tikzcd}
\]
where $\DbQl{1_G^{\sharp}}$ is the Change of Groups functor along the identity inclusion $1_G:\Spec K \to G$ against the pre-equivariant pseudofunctor $\DbQl{-}$ (and similarly for $\Shv(1_G^{\sharp})$; cf.\@ Theorem \ref{Thm: Section 3: Change of groups functor}) and where $\overline{D}^b$ and $\overline{\Shv}$ are the equivalence witnesses for
\[
D^b_{\Spec K}(X;\overline{\Q}_{\ell}) \simeq \DbQl{X}
\]
and
\[
\Shv_{\Spec K}(X;\overline{\Q}_{\ell}) \simeq \Shv(X;\overline{\Q}_{\ell})
\]
furnished by Proposition \ref{Prop: Section 2: Equivalence of K equivariant cat with F X}. Note that the inner squares form invertible $2$-cells via a routine use of Proposition \ref{Prop: Section 3: Change of groups interacts with Change of fibres} and the fact that the equivalences of Proposition \ref{Prop: Section 2: Equivalence of K equivariant cat with F X} are natural. Pasting all the cells together we get the invertible $2$-cell:
\[
\begin{tikzcd}
	\DbeqQl{X} \ar[rr, ""{name = Up}]{}{\Forget} \ar[d, swap]{}{H_G^0} & & \DbQl{X} \ar[d]{}{H^0} \\
	\Shv_G(X;\overline{\Q}_{\ell}) \ar[rr, swap, ""{name = Down}]{}{\Forget} & & \Shv(X;\overline{\Q}_{\ell}) \ar[from = Up, to = Down, Rightarrow, shorten >= 4pt, shorten <= 4pt]{}{\cong}
\end{tikzcd}
\]
Adding the invertible $2$-cell induced via the equivalence \\ $\Shv_G(X;\overline{\Q}_{\ell}) \simeq \Shv_G^{\text{na{\"i}ve}}(X;\overline{\Q}_{\ell})$ of Proposition \ref{Prop: Equiv of Standard t-structures on EDCs} fitting into the pasting diagram
\[
\begin{tikzcd}
	\DbeqQl{X} \ar[rr, ""{name = Up}]{}{\Forget} \ar[d, swap]{}{H_G^0} & & \DbQl{X} \ar[d]{}{H^0} \\
	\Shv_G(X;\overline{\Q}_{\ell}) \ar[rr, swap, ""{name = Down}]{}{\Forget} \ar[dr, swap]{}{\simeq} & & \Shv(X;\overline{\Q}_{\ell}) \ar[from = Up, to = Down, Rightarrow, shorten >= 4pt, shorten <= 4pt]{}{\cong} \\
	 & \Shv_G^{\text{na{\"i}ve}}(X;\overline{\Q}_{\ell}) \ar[ur, swap]{}{\Forget} \ar[from = Down, to = 3-2, Rightarrow, shorten <= 6pt, shorten >= 4pt]{}{\cong}
\end{tikzcd}
\]
produces the total diagram:
\[
\begin{tikzcd}
	\DbeqQl{X} \ar[rr, ""{name = Up}]{}{\Forget} \ar[d, swap]{}{H_G^0} & & \DbQl{X} \ar[d]{}{H^0} \\
	\Shv_G^{\text{na{\"i}ve}}(X;\overline{\Q}_{\ell}) \ar[rr, swap, ""{name = Down}]{}{\Forget} & & \Shv(X;\overline{\Q}_{\ell}) \ar[from = Up, to = Down, Rightarrow, shorten >= 4pt, shorten <= 4pt]{}{\cong}
\end{tikzcd}
\]
Note that the existence of the invertible $2$-cell
\[
\begin{tikzcd}
	\DbeqQl{X} \ar[rr, ""{name = Up}]{}{\Forget} \ar[d, swap]{}{{}^{p}H_G^0} & & \DbQl{X} \ar[d]{}{{}^{p}H^0} \\
	\Per_G^{\text{na{\"i}ve}}(X;\overline{\Q}_{\ell}) \ar[rr, swap, ""{name = Down}]{}{\Forget} & & \Per(X;\overline{\Q}_{\ell}) \ar[from = Up, to = Down, Rightarrow, shorten >= 4pt, shorten <= 4pt]{}{\cong}
	\end{tikzcd}
\]
follows mutatis mutandis to the diagram above, save from using Theorem \ref{Theorem: Equivalence of hearts with Equivariant Perverse Sheaves} instead of Proposition \ref{Prop: Equiv of Standard t-structures on EDCs}. Putting these together gives the proposition below.

\begin{proposition}\label{Prop: Forget restricts to naive forget}
For any $G$-variety $X$ there are invertible $2$-cells
\[
\begin{tikzcd}
	\DbeqQl{X} \ar[rr, ""{name = Up}]{}{\Forget} \ar[d, swap]{}{H_G^0} & & \DbQl{X} \ar[d]{}{H^0} \\
	\Shv_G^{\text{na{\"i}ve}}(X;\overline{\Q}_{\ell}) \ar[rr, swap, ""{name = Down}]{}{\Forget} & & \Shv(X;\overline{\Q}_{\ell}) \ar[from = Up, to = Down, Rightarrow, shorten >= 4pt, shorten <= 4pt]{}{\cong}
\end{tikzcd}
\]
and:
\[
\begin{tikzcd}
	\DbeqQl{X} \ar[rr, ""{name = Up}]{}{\Forget} \ar[d, swap]{}{{}^{p}H_G^0} & & \DbQl{X} \ar[d]{}{{}^{p}H^0} \\
	\Per_G^{\text{na{\"i}ve}}(X;\overline{\Q}_{\ell}) \ar[rr, swap, ""{name = Down}]{}{\Forget} & & \Per(X;\overline{\Q}_{\ell}) \ar[from = Up, to = Down, Rightarrow, shorten >= 4pt, shorten <= 4pt]{}{\cong}
\end{tikzcd}
\]
That is, the forgetful functors on the equivariant derived category restrict to the forgetful functors on equivariant perverse sheaves and equivarant $\ell$-adic sheaves.
\end{proposition}

Finally we show that the equivariant six functor formalism commutes with the forgetful functor up to isomorphism, which is the last ingredient we need to prove that the categories $\DbeqQl{X}$ satisfy Desiderata \ref{Desiderata}.

\begin{Theorem}\label{Thm: Six functor forgetful commute}
Let $h \in \GVar(X,Y)$. Then each of the six functors of Propositions \ref{Prop: Existence of push/pull}, \ref{Prop: Existence of proper pushpull}, \ref{Prop: EDC is symmetric monoidal closed} commute with the forgetful functors up to isomorphism.
\end{Theorem}
\begin{proof}
We begin by writing the non-equivariant isomorphisms $\psi_X:X \to \quot{X}{G}$ and $\psi_{Y}:Y\to \quot{Y}{G}$. Note that these are natural in the sense that the diagram of varieties
\[
\xymatrix{
X \ar[r]^-{h} \ar[d]_{\psi_{X}} & Y \ar[d]^{\psi_Y} \\
\quot{X}{G} \ar[r]_-{\quot{\overline{h}}{G}} & \quot{Y}{G}
}
\]
commutes. Similarly, write $\Forget_X:\DbeqQl{X} \to \DbQl{X}$ and $\Forget_Y:\DbeqQl{Y} \to \DbQl{Y}$.

Let us now show that the natural isomorphism
\[
\Forget_X \circ h^{\ast} \cong h^{\ast} \circ \Forget_Y
\]
exists. For this we calculate that on one hand, for any $A \in \DbeqQl{Y}_0$,
\begin{align*}
\Forget_X\left(h^{\ast}A\right) &= \Forget_X\left(\left\lbrace \hGamma^{\ast}\left(\AGamma\right) \; | \; \Gamma \in \Sf(G)_0, \AGamma \in A \right\rbrace\right) \\
&= \psi_X^{\ast}\left(\quot{\overline{h}}{G}^{\ast}\left(\quot{A}{G}\right)\right) = \left(\psi_X^{\ast} \circ \quot{\overline{h}}{G}^{\ast}\right)\quot{A}{G}.
\end{align*}
On the other hand
\[
h^{\ast}\left(\Forget_Y(A)\right) = h^{\ast}\left(\psi_Y^{\ast}\left(\quot{A}{G}\right)\right) = \left(h^{\ast} \circ \psi_Y^{\ast}\right)\left(\quot{A}{G}\right);
\]
using the natural isomorphism of functors $h^{\ast} \circ \psi_Y^{\ast} \cong \psi_X^{\ast} \circ \quot{\overline{h}}{G}^{\ast}$ thus gives the desired invertible $2$-cell:
\[
\begin{tikzcd}
	\DbeqQl{Y} \ar[r, ""{name = Up}]{}{h^{\ast}} \ar[d, swap]{}{\Forget_Y} & \DbeqQl{X} \ar[d]{}{\Forget_X} \\
	\DbQl{Y} \ar[r, swap, ""{name = Down}]{}{h^{\ast}} & \DbQl{X} \ar[from = Up, to = Down, Rightarrow, shorten >= 4pt, shorten <= 4pt]{}{\cong}
\end{tikzcd}
\]

We now show the existence of the natural isomorphism
\[
\Forget_X \circ h^{!} \cong h^{!} \circ \Forget_Y.
\]
Because the diagram
\[
\xymatrix{
	X \ar[r]^-{h} \ar[d]_{\psi_{X}} & Y \ar[d]^{\psi_Y} \\
	\quot{X}{G} \ar[r]_-{\quot{\overline{h}}{G}} & \quot{Y}{G}
}
\]
is a pullback diagram with $\psi_X$ and $\psi_Y$ smooth (by virtue of being isomorphisms) we have a natural isomorphism of functors $h^{!}\circ\psi_Y^{\ast}  \cong \psi_X^{\ast} \circ \quot{\overline{h}}{G}^{!}$. Using this we calculate for any $A \in \DbeqQl{Y}_0$
\begin{align*}
\Forget_X\left(h^{!}A\right) &= \Forget_X\left(\left\lbrace \hGamma^{!}\left(\AGamma\right) \; | \; \Gamma \in \Sf(G)_0, \AGamma \in A \right\rbrace\right) \\
&= \psi_X^{\ast}\left(\quot{\overline{h}}{G}^{!}\left(\quot{A}{G}\right)\right) = \left(\psi_X^{\ast} \circ \quot{\overline{h}}{G}^{!}\right)\quot{A}{G} \\
&\cong \left(h^{!}\circ\psi_Y^{\ast}\right)\quot{A}{G} = h^{!}\left(\psi_Y^{\ast}(\quot{A}{G})\right) \\
&= h^{!}\left(\Forget_Y\left\lbrace \AGamma \; | \; \Gamma \in \Sf(G)_0\right\rbrace\right) = h^{!}\left(\Forget_Y(A)\right)
\end{align*}
This establishes the invertible $2$-cell:
\[
\begin{tikzcd}
	\DbeqQl{Y} \ar[r, ""{name = Up}]{}{h^{!}} \ar[d, swap]{}{\Forget_Y} & \DbeqQl{X} \ar[d]{}{\Forget_X} \\
	\DbQl{Y} \ar[r, swap, ""{name = Down}]{}{h^{!}} & \DbQl{X} \ar[from = Up, to = Down, Rightarrow, shorten >= 4pt, shorten <= 4pt]{}{\cong}
\end{tikzcd}
\]

We now prove the existence of the natural isomorphism
\[
Rh_{\ast} \circ \Forget_X \cong \Forget_Y \circ Rh_{\ast}.
\]
As in the case above, we use that the diagram
\[
\xymatrix{
	X \ar[r]^-{h} \ar[d]_{\psi_{X}} & Y \ar[d]^{\psi_Y} \\
	\quot{X}{G} \ar[r]_-{\quot{\overline{h}}{G}} & \quot{Y}{G}
}
\]
is a pullback diagram of varieties with both $\psi_X$ and $\psi_Y$ smooth. Thus we get isomorphisms of functors $Rh_{\ast} \circ \psi_X^{\ast} \cong \psi_Y^{\ast} \circ R(\quot{\overline{h}}{G})_{\ast}$ by the Smooth Base Change Theorem. Proceeding as in the two prior cases we derive that for any $A \in \DbeqQl{X}_0$,
\begin{align*}
(Rh_{\ast} \circ \Forget_X)A &= \left(Rh_{\ast} \circ \psi_X^{\ast}\right)\quot{A}{G} \cong \left(\psi_Y^{\ast} \circ R(\quot{\overline{h}}{G})_{\ast}\right)\quot{A}{G} \\
&= (\Forget_Y \circ Rh_{\ast})A
\end{align*}
which gives the invertible $2$-cell:
\[
\begin{tikzcd}
	\DbeqQl{X} \ar[r, ""{name = Up}]{}{Rh_{\ast}} \ar[d, swap]{}{\Forget_Y} & \DbeqQl{Y} \ar[d]{}{\Forget_X} \\
	\DbQl{X} \ar[r, swap, ""{name = Down}]{}{Rh_{\ast}} & \DbQl{Y} \ar[from = Up, to = Down, Rightarrow, shorten >= 4pt, shorten <= 4pt]{}{\cong}
\end{tikzcd}
\]
In fact, using that the diagram
\[
\xymatrix{
	X \ar[r]^-{h} \ar[d]_{\psi_{X}} & Y \ar[d]^{\psi_Y} \\
	\quot{X}{G} \ar[r]_-{\quot{\overline{h}}{G}} & \quot{Y}{G}
}
\]
satisfies the hypotheses of Lemma \ref{LEmma: Proper base change commutes iwth pullback} (because every morphism and scheme in sight is finite type and separated) we get that $ Rh_{!} \circ \psi_X^{\ast}  \cong \psi_Y^{\ast} \circ R(\quot{\overline{h}}{G})_{!}$. Proceeding mutatis mutandis to the pusforward case gives the invertible $2$-cell:
\[
\begin{tikzcd}
	\DbeqQl{X} \ar[r, ""{name = Up}]{}{Rh_{!}} \ar[d, swap]{}{\Forget_Y} & \DbeqQl{Y} \ar[d]{}{\Forget_X} \\
	\DbQl{X} \ar[r, swap, ""{name = Down}]{}{Rh_{!}} & \DbQl{Y} \ar[from = Up, to = Down, Rightarrow, shorten >= 4pt, shorten <= 4pt]{}{\cong}
\end{tikzcd}
\]

We now show that the tensor products commute with the forgetful functor, i.e., there is a natural isomorphism
\[
\Forget_X\left(A \os{L}{\otimes} B\right) \cong \Forget_X(A) \os{L}{\otimes} \Forget_X(B)
\]
for any $A,B \in \DbeqQl{X}$. For this we calculate that on one hand
\begin{align*}
\Forget_X\left(A \os{L}{\otimes} B\right) &= \psi_X^{\ast}\left(\quot{A}{G} \ds{L}{\XGamma}{\otimes} \quot{B}{G}\right)
\end{align*}
while on the other hand
\[
\Forget_X(A) \os{L}{\otimes} \Forget_X(B) = \psi_X^{\ast}\left(\quot{A}{G}\right) \os{L}{\otimes} \psi_X^{\ast}\left(\quot{B}{G}\right).
\]
Because $\psi:X \to \quot{X}{G}$ is smooth, there is a a natural isomorphism
\[
\psi_X^{\ast}\left(\quot{A}{G} \ds{L}{\XGamma}{\otimes} \quot{B}{G}\right) \cong \psi_X^{\ast}\left(\quot{A}{G}\right) \os{L}{\otimes} \psi_X^{\ast}\left(\quot{B}{G}\right)
\]
which gives an invertible $2$-cell:
\[
\begin{tikzcd}
	\DbeqQl{X} \times \DbeqQl{X} \ar[rr, ""{name = Up}]{}{\os{L}{\otimes}} \ar[d, swap]{}{\Forget_X \times \Forget_X} & & \DbeqQl{X} \ar[d]{}{\Forget_X} \\
\DbQl{X} \times \DbQl{X} \ar[rr, swap, ""{name = Down}]{}{\os{L}{\otimes}} & & \DbQl{X} \ar[from = Up, to = Down, Rightarrow, shorten >= 4pt, shorten <= 4pt]{}{\cong}
\end{tikzcd}
\]
Following this argument mutatis mutandis for the functors $R[-,-]$ gives an invertible $2$-cell
\[
\begin{tikzcd}
	\DbeqQl{X}^{\op} \times \DbeqQl{X} \ar[rr, ""{name = Up}]{}{R[-,-]} \ar[d, swap]{}{\Forget_X^{\op} \times \Forget_X} & & \DbeqQl{X} \ar[d]{}{\Forget_X} \\
	\DbQl{X}^{\op} \times \DbQl{X} \ar[rr, swap, ""{name = Down}]{}{R[-,-]} & & \DbQl{X} \ar[from = Up, to = Down, Rightarrow, shorten >= 4pt, shorten <= 4pt]{}{\cong}
\end{tikzcd}
\]
and completes the proof of the proposition.
\end{proof}
We now have the ingredients necessary to prove that the equivariant derived category $\DbeqQl{X}$ satisfies Desiderata \ref{Desiderata}.

\begin{Theorem}\label{Theorem: EDC sats desiderata}
The equivariant derived category $\DbeqQl{X}$ satisfies Desiderata \ref{Desiderata}.
\end{Theorem}
\begin{proof}
We prove that $\DbeqQl{X}$ satisfies the conditions of Desiderata \ref{Desiderata} in the order the properties are listed there.
\begin{itemize}
	\item That the category $\DbeqQl{X}$ is triangulated is Theorem \ref{Theorem: Section Triangle: Equivariant triangulation}
	\item The existence of the standard and perverse equivariant $t$-structures with
	\[
	{}^{\text{stand}}\DbeqQl{X}^{\heartsuit} \simeq \Shv_G(X;\overline{\Q}_{\ell}) \simeq \Shv_G^{\text{na{\"i}ve}}(X;\overline{\Q}_{\ell})
	\]
	and
	\[
	{}^{p}\DbeqQl{X}^{\heartsuit} \simeq \Per_G(X;\overline{\Q}_{\ell}) \simeq \Per_G^{\text{na{\"i}ve}}(X;\overline{\Q}_{\ell})
	\]
	are Proposition \ref{Prop: Equiv of Standard t-structures on EDCs} and Theorem \ref{Theorem: Equivalence of hearts with Equivariant Perverse Sheaves}.
	\item The existence of the forgetful functors come from Example \ref{Example: Section 2: The forgetful functor} and the commutativity of the forgetful functors with restriction to the hearts is Proposition \ref{Prop: Forget restricts to naive forget}.
	\item The existence of the equivariant six functors is Propositions \ref{Prop: Existence of push/pull}, \ref{Prop: Existence of proper pushpull}, and \ref{Prop: EDC is symmetric monoidal closed}; the commutativity of these functors with the forgetful functors is Theorem \ref{Thm: Six functor forgetful commute}.
	\item The open-closed distinguished triangle is Proposition \ref{Prop: Exact open/closed triangle}.
	\item The relationships
		\begin{align*}
		\big(Rh_!A\big) \ds{L}{Y}{\otimes} B &\cong Rh_!\left(A \ds{L}{X}{\otimes} h^{\ast}B\right), \\
		R[Rh_!A,B]_Y &\cong Rh_{\ast}\left(R[A,h^{!}B]_X\right), \\
		h^{!}R[A,B]_Y &\cong R[h^{\ast}A,h^{!}B]_X.
	\end{align*}
	between the equivariant six functors is Proposition \ref{Prop: Equivariant Grothendiekc Context}.
\end{itemize}
\end{proof}

\chapter{The Simplicial Equivariant Derived Equivalence}\label{Section: Lusztig EDC is the Simplicial EDC}

Our next goal in this paper is to give a very explicit proof of the equivalence $\DbeqQl{X}$ with the simplicial equivariant derived category $\Dbeqsimp{X}$. The big point of this perspective is that the simplicial equivariant derived category is very straightforward to give computations with, and developing the various interaction functors that it has explicitly really highlights how this category can be used. It is also the bridge through which we prove the stacky equivalence $\DbeqQl{X} \simeq \DbQl{[G \backslash X]}$, so getting to know this category in detail now will pay off later on.

Another large benefit of working with $\Dbeqsimp{X}$ is that it allows us to give an ``adic'' formal metaphor to the notion of equivariant sheaves. Because we can regard an equivariant (derived) complex of $\ell$-adic sheaves as a sequence of derived complexes of $\ell$-adic sheaves on the varieties $G^n \times X$ which have isomorphisms against the pullbacks $h^{\ast}A_n \to A_m$ when $h$ is a map coming from a simplicial morphism (cf.\@ Lemma \ref{Lemma: Simplicial Scheme of X} and Definition \ref{Defn: Simplicial sheaves}). In this way we can see each complex $A_n$ in $\DbQl{G^n \times X}$ as reading out an $n$-infinitesimal neighborhood of how the group acts on $A = (A_n)_{n \in \N}$ while the pullback isomorphisms $h^{\ast}A_n \to A_m$ tell us how to move these infinitesimal neighborhoods around by using either the group multiplication/insertion of identities in the $G^n$-components or by using the $G$-action (in the $G\times X$ component). A handy consequence of this formal analogy is Lemma \ref{Lemma: Section 4: Quot pullback for An}, which may be rephrased as saying that for all $n \geq 1$ the pullback of the quotient functor
\[
\quo_{G^n}^{\ast}:\DbQl{\quot{X}{G^n}} \to \DbQl{G^n \times X}
\]
essentially acts by picking out the sheaves which are concentrated within only an infinitesimal neighborhood of degree $0$. That is, the (derived complexes of) sheaves which arise as pullbacks against the quotient morphism play the role of constant terms in a $p$-adic series expansion. I do want to point out, however, that it is unknown if there is anything more than a formal analogy going on or if this is what is happening in more concrete terms. It would be of interest to examine the equivariant simplicial cohomology of $\DbeqQl{G^n \times X}$ and see in this case if a sheaf can be approximated by all the cohomology sheaves and then if those which arise from pullbacks are concentrated in simplicial equivariant degree $0$ only, but at this point whether or not this happens is unknown to the author (as is the precise notion of cohomology one would need to use to give the ``simplicial equivariant'' cohomology whose $H^n$ records $n$-infinitesimal neighborhoods as $H_{dR}^n$ records $n$-infinitesimal neighborhoods of manifolds).

We now proceed to discuss the theory of simplicial schemes and simplicial sheaves on a simplicial scheme. However, before describing these objects we will need to review the topologist's simplex category $\mathbf{\Delta}$, as it is the category for which simplicial objects in a category $\Cscr$ arise as $\Cscr$-valued presheaves.
\begin{definition}\index{The Topologist's Simplex Category}\index[notation]{Delta@$\mathbf{\Delta}$}
The topologist's simplex category $\mathbf{\Delta}$ is defined as follows:
\begin{itemize}
	\item Objects: For all $n \in \N$ a linearly ordered set
	\[
	[n] = \lbrace 0, \cdots, n\rbrace.
	\]
	\item Morphisms: Monotonic increasing functions of sets, i.e., a morphism $f:[m] \to [n]$ is a function in $\mathbf{\Delta}_1$ if and only if for all $j,k \in [m]$ with $j \leq k $ we have that $f(j) \leq f(k)$.
	\item Composition: As in $\Set$.
	\item Identities: As in $\Set$.
\end{itemize}
\end{definition}

It is worth noting that the category $\mathbf{\Delta}$ is equivalent to the category of finite linearly ordered sets and order preserving maps. 

\begin{definition}\index{Simplicial! Scheme}
	Let $\mathbf{\Delta}$ be the topologist's simplex category. A simplicial scheme is a functor
	\[
	\underline{X}:\mathbf{\Delta}^{\op} \to \Sch.
	\]
	Similarly, for any scheme $S$, a simplicial $S$-scheme is a functor
	\[
	\underline{X}:\mathbf{\Delta}^{\op} \to \Sch_{/S}.
	\]
	A morphism of simplicial schemes, $\underline{f}:\underline{X}_{\bullet} \to \underline{Y}_{\bullet}$ is a natural transformation between functors. In more detail, a morphism of simplicial schemes is a sequence of morphisms of schemes $f_n:X_n \to Y_n$ such that for all $h \in \mathbf{\Delta}_1$, the diagram
	\[
	\xymatrix{
		X_m \ar[r]^-{f_m} \ar[d]_{\underline{X}(h)} & Y_m \ar[d]^{\underline{Y}(h)} \\
		X_n \ar[r]_{f_n} & Y_n
	}
	\]
	commutes. Similarly, a simplicial variety $\underline{X}_{\bullet}$\index{Simplicial! Variety}  is a simplicial $(\Spec K)$-scheme which factors through the category of $K$-varieties.
	
	We write $[\mathbf{\Delta}^{\op},\Sch]$\index[notation]{sSchemes@$[\mathbf{\Delta}^{\op},\Sch]$} for the category of simplicial schemes and we write $[\mathbf{\Delta}^{\op},\Var_{/\Spec K}]$\index[notation]{sVar@$[\mathbf{\Delta}^{\op},\Var_{/\Spec K}]$} for the category of simplicial varieties.
\end{definition}
\begin{remark}\label{Remark: (Co)Face Maps and (Co)Degeneracy Maps}
By \cite[Page 4]{GoerssJardine} it there is a class of maps in $\mathbf{\Delta}$ which generates all morphisms in $\mathbf{\Delta}_1$ under composition. These are the coface and codegeneracy maps, and they satisfy a collection of cosimplicial identities \cite[Page 4, Equation 1.2]{GoerssJardine}. We will write $\delta_k^n:[n ] \to [n+1]$\index[notation]{deltankn@$\delta_k^n$} and $\sigma_{k}^{n}:[n+1] \to [n]$\index[notation]{zetank@$\sigma_k^n$}\index{Simplicial! Coface Maps}\index{Simplicial! Codegeneracy Maps} for the coface and codegeneracy maps in $\mathbf{\Delta}$ whenever necessary. Moreover, we call the corresponding opposite morphisms of a simplicial object $\underline{X}_{\bullet}$ a face map $d_k^n:\underline{X}_{n+1} \to \underline{X}_n$ and a degeneracy map $s_k^n:\underline{X}_n \to \underline{X}_{n+1}$.

The upshot of these the existence of this maps is they tell us that to ensure that an object $\underline{X}_{\bullet}$ is a simplicial object in a category $\Cscr$, it suffices to define $X_n$ for all $n \in \N$ and then to define all face maps\index{Simplicial! Face Maps}\index{Simplicial! Degeneracy Maps} $d_i^n:X_{n+1} \to X_n$\index[notation]{dnk@$d_k^n$} and degeneracy maps $s_i^n:X_n \to X_{n+1}$\index[notation]{snk@$s_i^n$} in such a way that they satisfy the standard simplicial identities (cf.\@ \cite[Page 5, Equation 1.3]{GoerssJardine}). In such a case we can fully ensure that $\underline{X}_{\bullet}$ is defined by setting $\underline{X}(\delta_k^n) := d_k^n$ and $\underline{X}(\sigma_k^n) = s_k^n$.
\end{remark}
\begin{example}
	Let $X$ be a scheme. Then the functor $\const(X)_{\bullet}:\mathbf{\Delta}^{\op} \to \Sch$ given by $\const(X)([n]) := X$ for all $n \in \N$ and $\const(X)(h) = \id_X$ for all $h \in \mathbf{\Delta}_1$ is a simplicial scheme.
	
	A more interesting example is the following. Define the simplicial scheme $\underline{\Abb}^{\bullet}_{K}$ by:
	\begin{itemize}
		\item For all $n \in \N$, $\underline{\Abb}^{\bullet}_K([n]) = \Abb_K^{n}$.
		\item The morphisms are induced by setting the face maps $d_k^n:\Abb^{n+1}_K \to \Abb_K^n$ to be given by
		\[
		d_k^n = \begin{cases}
		\pi_{1,\cdots,n}^{\Abb_K^{n+1}} & \text{if}\, k = 0; \\
		\id_{\Abb^{k}_K} \times \alpha \times \id_{\Abb_K^{n+1-k}} & \text{if}\, 1 \leq k \leq n+1;
		\end{cases}
		\]
		where $\alpha:\Abb_K^1 \times \Abb_K^1 \to \Abb_K^1$ is the addition of $\Gbb_{a,K}$ and setting the degeneracy maps $s_k^n:\Abb_K^{n} \to \Abb_K^{n+1}$ to be given by the insertion of the identity $0_K:\Spec K \to \Abb_K^1$ in the $k$-th component of the product $\Abb_K^n \to \Abb_K^{n+1}$ for $0 \leq k \leq n$. Then $\underline{\Abb}_K^{\bullet}$ is a simplicial scheme which gives a simplicial presentation of the pro-variety $\Abb_K^{\infty}$.
	\end{itemize}
\end{example}
\begin{remark}\index[notation]{Xbullet@$\underline{X}_{\bullet}$}\index[notation]{Ybullet@$\underline{Y}_{\bullet}$}\index[notation]{funderline@$\underline{f}$}
When working with simplicial objects, we will write $\underline{X}_{\bullet}$ and $\underline{Y}_{\bullet}$ as simplicial schemes/varieties and write $\underline{X}_{n} := \underline{X}([n])$ and $\underline{Y}_{n} := \underline{Y}([n])$. Similarly, morphisms between simplicial objects $\underline{X}_{\bullet} \to \underline{Y}_{\bullet}$ are denoted by underlined morphisms $\underline{f}$ where each $\underline{f}_n = \underline{f}([n]):X_n \to Y_n$ is a morphism.
\end{remark}

There is a functor
\[
\const:\Sch \to [\mathbf{\Delta}^{\op},\Sch]
\]\index[notation]{const@$\const$}
which sends a scheme $X$ to the constant simplicial scheme $\const(X)_{\bullet}$ on $X$; this is the simplicial scheme $\const(X)_{\bullet}$ where $\const(X)_n = X$ for all $n \in \N$ and $\const(X)_{\bullet}(h) = \id_X$ for all $h \in \mathbf{\Delta}_1$. We use this to define augmentations of simplicial schemes (cf.\@ Definition \ref{Defn: Section 4: Augmentations of simplicial schemes}), which are in turn used to prove the equivalence $D^b_{eq}(\underline{G \backslash X};\overline{\Q}_{\ell}) \simeq \DbeqQl{X}$. For this we need two crucial observations which are well-known to experts; however, we write them down here to explicitly use later. Both these observations involve the language of simplicial sheaves, however, so we move to define those now. When working with a simplicial sheaf on a simplicial scheme $\underline{X}_{\bullet}$, we will follow the conventions of \cite{BernLun} and \cite{DeligneHodge3}.
\begin{definition}[{\cite[Section 5.1.6]{DeligneHodge3}}]\index{Simplicial! Sheaf}
	A simplicial ({\'e}tale) sheaf $\underline{\Fscr}$ on a simplicial scheme $\underline{X}_{\bullet}$ is a sequence $(\Fscr_n)_{n \in \N}$ of {\'e}tale sheaves $\Fscr_n$ on $X_n$ such that if $h \in \mathbf{\Delta}([m],[n])$, then there is a morphism of sheaves on $X_m$
	\[
	\alpha_h:h^{\ast}\Fscr_n \to \Fscr_m
	\]
	such that for all composable arrows $[m] \xrightarrow{h} [n] \xrightarrow{g} [s]$ in $\mathbf{\Delta}$, the cocycle condition
	\[
	\alpha_{g \circ h} = \alpha_h \circ h^{\ast}\alpha_g
	\]
	is satisfied. Furthermore, we say that $\underline{\Fscr}$ is a simplicial $\ell$-adic sheaf on $\underline{X}_{\bullet}$ if each sheaf $\Fscr_n$ is an $\ell$-adic sheaf on $X_n$.
	
	A morphism of simplicial sheaves $\underline{\Fscr} \to \underline{\Gscr}$ is a sequence of natural transformations
	\[
	\rho_n:\Fscr_n \to \Gscr_n
	\]
	for all $n \in \N$ such that for any $h \in \mathbf{\Delta}_1$, the diagram
	\[
	\xymatrix{
		h^{\ast}\Fscr_{m} \ar[d]_{\alpha_h} \ar[r]^{h^{\ast}\rho_m} & h^{\ast}\Gscr_m \ar[d]^{\beta_h} \\
		\Fscr_n \ar[r]_{\rho_n} & \Gscr_n
	}
	\]
	commutes. 
	
	Finally, for simplicial schemes $\underline{X}_{\bullet}$, we write $\Shv(\underline{X}_{\bullet};\overline{\Q}_{\ell})$\index{sShvG@$\Shv_G(\underline{X}_{\bullet};\overline{\Q}_{\ell})$} for the category of $\ell$-adic simplicial sheaves on the simplicial scheme $\underline{X}_{\bullet}$.
\end{definition}
\begin{proposition}\label{Prop: Section 4: ladic simplicial sheaves are Abelian}
	For any simplicial scheme $\underline{X}_{\bullet}$, the category $\Shv(\underline{X}_{\bullet};\overline{\Q}_{\ell})$ is Abelian.
\end{proposition}
\begin{proof}
	We begin by first defining the zero simplicial sheaf by $\underline{0} = (0_n)$, i.e., the simplicial sheaf which is degree-wise the zero sheaf. For any two $\ell$-adic simplicial sheaves $\underline{\Fscr}$ and $\underline{\Gscr}$, we define the direct sum $\underline{\Fscr} \oplus \underline{\Gscr}$ by
	\[
	(\underline{\Fscr} \oplus  \underline{\Gscr})_n := \Fscr_n \oplus \Gscr_n;
	\]
	similarly, for any two morphisms $\rho, \psi$ of simplicial $\ell$-adic sheaves with common domain and codomain, we define $\rho + \psi$ by setting
	\[
	(\rho +\psi)_n := \rho_n + \psi_n.
	\]
	That this gives rise to an additive category is trivially verified and omitted. 
	
	In a similar vein, for any morphism $\rho$ we find that $\underline{\Ker \rho}$ is given via the assignment
	\[
	(\underline{\Ker \rho})_n := \Ker \rho_n
	\]
	and dually for cokernels. That these are suitably regular, i.e., that every epic and monic are regular, follow from the fact that these statements hold on the level of $\ell$-adic sheaves for each $X_n$. It thus follows that $\Shv(\underline{X}_{\bullet};\overline{\Q}_{\ell})$ is an Abelian category.
\end{proof}
\begin{definition}\label{Remark: Section 4: ell-adic simplicial sheaves}\index{Simplicial! Sheaves! Derived Category of}
	If $\underline{X}_{\bullet}$ is a simplicial scheme, we will write $D^b(\underline{X}_{\bullet};\overline{\Q}_{\ell})$ for the bounded derived category of $\ell$-adic simplicial sheaves and $\DbQl{\underline{X}_{\bullet}}$ for the category of simplicial (complexes of) $\ell$-adic constructible sheaves.
\end{definition}
\begin{proposition}\label{Prop: Simplicial morphisms induce pullback functors of sheaf categories}
	For any morphisms of simplicial schemes $\underline{f}:\underline{X}_{\bullet} \to \underline{Y}_{\bullet}$, there is an induced exact (additive) functor
	\[
	\underline{f}^{\ast}:\Shv(\underline{Y}_{\bullet};\overline{\Q}_{\ell}) \to \Shv(\underline{X}_{\bullet};\overline{\Q}_{\ell}).
	\]
\end{proposition}
\begin{proof}
	We define $\underline{f}^{\ast}$ by sending a simplicial sheaf $(\Fscr_{n})_{n \in \N}$ to the sequence of sheaves $(f_n^{\ast}\Fscr_{n})_{n \in \N}$. To see that this is a simplicial sheaf on $\underline{X}_{\bullet}$ we must find structure maps $\alpha_h$ for each $h \in \mathbf{\Delta}_1$ which satisfy the cocycle condition. Fix an aribitrary $h \in \mathbf{\Delta}_1$ and write $h:[n] \to [m]$. We now need to construct a morphism
	\[
	\underline{X}(h)^{\ast}(f_n^{\ast}\Fscr_{n}) \xrightarrow{\alpha_h} f_{m}^{\ast}\Fscr_m
	\]
	For this we recall from the fact that $\underline{f}$ is a morphism of simplicial schemes that the diagram
	\[
	\xymatrix{
		X_m \ar[r]^-{f_m} \ar[d]_{\underline{X}(h)} & Y_m \ar[d]^{\underline{Y}(h)} \\
		X_n \ar[r]_{f_n} & Y_n
	}
	\]
	commutes. As such, because the pseudofunctor $\Shv(-;\overline{\Q}_{\ell})$ is strict over schemes, we find that 
	\[
	\underline{X}(h)^{\ast}(f_n^{\ast}{\Fscr_n}) = (f_n \circ \underline{X}(h))^{\ast}\Fscr_n = (\underline{Y}(h) \circ f_m)^{\ast}\Fscr_n =  f_m^{\ast}(\underline{Y}(h)^{\ast}\Fscr_n).
	\]
	From this we define our structure morphism $\alpha_h:\underline{X}(h)^{\ast}(f_n^{\ast}\Fscr_n) \to f_{m}^{\ast}\Fscr_m$ via the diagram:
	\[
	\xymatrix{
		\underline{X}(h)^{\ast}(f_n^{\ast}{\Fscr_n}) \ar[d]_{\alpha_h} \ar@{=}[r] & (f_n \circ \underline{X}(h))^{\ast}\Fscr_n \ar@{=}[d] \\
		f_m^{\ast}\Fscr_m & f_m^{\ast}(\underline{Y}(h)^{\ast}\Fscr_n) \ar[l]^{f_m^{\ast}\beta_h}
	}
	\]
	That these maps satisfy the cocycle condition follows mutatis mutandis to the proof of the cocycle condition in Theorem \ref{Thm: Section 3: Existence of pullback functors}, save with $\Sf(G)$ replaced by $\mathbf{\Delta}$ and the simplified condition that the transition isomorphisms are all identity maps. Thus $\underline{f}^{\ast}\underline{\Fscr}$ is a simplicial sheaf on $\underline{X}$. The assignment on morphisms is given by
	\[
	\underline{\rho} = (\rho_n) \mapsto (f_n^{\ast}\rho_n) = \underline{f}^{\ast}\underline{\rho},
	\]
	which is can be verified to be a morphism in $\Shv(\underline{X}_{\bullet};\overline{\Q}_{\ell})$ mutatis mutandis to the proof of Theorem \ref{Thm: Section 3: Existence of pullback functors} as well. This yields immediately that $\underline{f}^{\ast}$ is a functor; that it is exact follows immediately from the isomorphisms 
	\[
	f_n^{\ast}(\Fscr_n) = f_n^{\ast}\lim_{\substack{\longleftarrow \\ i \in I}} \Fscr_{i,n} \cong \lim_{\substack{\longleftarrow \\ i \in I}} f_n^{\ast}\Fscr_{i,n}
	\]
	for any finite limit of $\ell$-adic sheaves together with the fact that limits are computed degreewise in the category $\Shv(\underline{X}_{\bullet};\overline{\Q}_{\ell})$, i.e., 
	\[
	\lim(\Fscr_{i})_{n \in \N} = (\lim \Fscr_{i})_{n \in \N}.
	\]
\end{proof}
\begin{corollary}
	A map of simplicial schemes $\underline{f}^{\ast}:\underline{X}_{\bullet} \to \underline{Y}_{\bullet}$ induces a derived pullback functor
	\[
	f^{\ast}:D^b_c(\underline{Y}_{\bullet};\overline{\Q}_{\ell}) \to D^b_c(\underline{X}_{\bullet};\overline{\Q}_{\ell})
	\]
	which restricts to a constructible pullback functor
	\[
	f^{\ast}:\DbQl{\underline{Y}_{\bullet}} \to \DbQl{\underline{X}_{\bullet}}.
	\]
\end{corollary}
\begin{proposition}\label{Prop: Section 4: Equiv of categories for constant simplicial scheme}
	For any scheme $X$, there is a a canonical equivalence of categories
	\[
	D^b_c(X;\overline{\Q}_{\ell}) \simeq D^b_c(\const(X)_{\bullet};\overline{\Q}_{\ell}).
	\]
\end{proposition}
\begin{proof}
Recall that the simplicial scheme  $\const(X)_{\bullet}$ is defined by setting, for all $n$, $\const(X)_n = X$ and setting all the face and degeneracy maps to be the identity morphisms. We define a functor $F:\DbQl{X} \to \DbQl{\const(X)_{\bullet}}$ by defining $FA := (A)_{n \in \N}$ (the constant sequence on $A$) on objects and $F\rho := (\rho)_{n \in \N}$ (the constant sequence at $\rho$) on morphisms. In the other direction, define $H:\DbQl{\const(X)_{\bullet}} \to \DbQl{X}$ by $H(A_n)_{n \in \N} := A_0$ on objects and $H(\rho_n)_{n \in \N} := \rho_0$ on morphisms. That both these assignments are functors is straightforward to check.

We first show that $H \circ F \cong \id_{\DbQl{X}}$ For this we calculate, for $A \in \DbQl{X}_0$,
\[
(H \circ F)A = H(A)_{n \in \N} = A = \id_{\DbQl{X}}(A)
\]
which is equal to (and hence isomorphic) to the identity functor.

Let us now check the other composite and note for an object $(A_n)$ of $\DbQl{\const(X)_{\bullet}}_0$,
\[
(F \circ H)(A_n)_{n \in \N} = FA_0 = (A_0)_{n \in \N}.
\]
We now must show that $(A_n)_{n \in \N} \cong (A_0)_{n \in \N}$. This is immediate, however, because the structure isomorphisms, for any $h:[0] \to [n]$ in $\mathbf{\Delta}$,
\[
\alpha_h:h^{\ast}A_n \xrightarrow{\cong} A_0
\]
corresponds to an isomorphism $\id^{\ast}A_n \xrightarrow{\cong} A_0$. Because the pullback against the identity morphism is itself an identity, it follows that $\alpha_h$ is an isomorphism
\[
A_n = \id^{\ast}A_n \xrightarrow{\cong} A_0.
\]
For each $n \in \N$, let $\nu_n:[n] \to [0]$ denote the terminal map in $\mathbf{\Delta}$. We define our isomorphism $\alpha:(A_0)_{n \in \N} \to (A_n)_{n \in \N}$ by, for each $n \in \N$,
\[
\alpha_n:A_0 \xrightarrow{\alpha_{\nu_{n}}} \id^{\ast}A_n = A_n.
\]
To see this is indeed an isomorphism, it suffices to prove that $\alpha$ commutes with the structure isomorphisms. For this we note that $\alpha_{(A_0)} = \id$ for any simplicial map, while if $h:[m] \to [n]$ is a morphism in $\mathbf{\Delta}$, $h^{\ast} = \id^{\ast}$ and
\[
\nu_n = \nu_m \circ h.
\]
Thus applying the cocycle condition we find
\[
\alpha_m \circ \alpha_{A_0} = \alpha_m \circ \id_{A_0} = \alpha_m = \alpha_{\nu_m} = \alpha_{h \circ \mu_n} = \alpha_{\nu_n} \circ \nu_n^{\ast}\alpha_h = \alpha_{\nu_n} \circ \alpha_h
\]
so the diagram
\[
\xymatrix{
A_n \ar[r]^-{\alpha_h} & A_m \\
A_0 \ar[u]^{\alpha_{\nu_n}} \ar@{=}[r] & A_0 \ar[u]_{\alpha_{\nu_m}}
}
\]
commutes. This shows that $\alpha$ is an isomorphism and hence that $H \circ F \cong \id_{\DbQl{\const(X)_{\bullet}}}$.
\end{proof}
We now proceed to discuss the simplicial quotient scheme of a $G$-variety. Let $G$ be a smooth affine algebraic group over $K$ and let $X$ be a (left) $G$-variety. We now want to define a simplcial scheme $\underline{G \backslash X}_{\bullet}$ to play the role of a quotient stack and/or quotient scheme for the $G$-action on $X$. Let us proceed by using a schematic version of the simplical space denoted $[G \backslash X]_{\bullet}$ in \cite{JamesMracek}. 
\begin{definition}\index{Simplicial! Quotient Scheme}\index{Simplicial! Quotient Variety}\index[notation]{GslashXbullet@$\underline{G \backslash X}_{\bullet}$}
We define $\underline{G \backslash X}_{\bullet}$ by first setting
\[
\underline{G \backslash X}_n := G^n \times X
\]
for all $n \in \N$. We next define the face maps $d_i^n:\underline{G \backslash X}_{n+1} \to \underline{G \backslash X}_n$ by:
\[
d_i^n := \begin{cases}
\pi_{1, \cdots, n}^{G^{n+1}} \times \id_X & \text{if}\, i = 0; \\
\id_{G^{i-1}} \times \mu \times \id_{G^{n+1-(i+1)} \times X} & \text{if}\, 0 < i < n+1 \\
\id_{G^{n}} \times \alpha_X & \text{if}\, i = n+1.
\end{cases}
\]
The degeneracies $s_i^{n}:\underline{G \backslash X}_{n} \to \underline{G \times X}_{n+1}$ are given by inserting the identity of $G$ into the $i$-th component of the product $G^{n+1}$, i.e., for $0 \leq i \leq n$ we define $s_i^{n}:G^{n} \times X \to G^{n+1} \times X$ to be the morphism factoring as in the commuting diagram:
\[
\xymatrix{
	G^{n} \times X \ar[rrrr]^-{s_i^{n}} \ar[d]_{\cong} & & & & G^{n+1} \times X \\
	G^{i} \times \Spec K \times G^{n-i} \times X \ar[rrrr]_-{\id_{G^i} \times 1_G \times \id_{G^{n-i}} \times \id_X} & & & & G^{i} \times G \times G^{n-i} \times X \ar[u]_{\cong}
}
\]
\end{definition}
Let us explain these maps a little more, as while they are correct, they are perhaps written in an unintuitive way. In a set-theoretic abuse of notation, they are given by sending a tuple $(g_0, \cdots, g_{n},x)$ to
\[
d_i^{n}(g_0, \cdots, g_n, x) = \begin{cases}
(g_1, \cdots, g_n,x) & \text{if}\, i = 0; \\
(g_0, \cdots, g_{i-1}g_i, \cdots, g_n, x) & \text{if}\, 1 \leq  i \leq n; \\
(g_0, \cdots, g_{n-1}, g_nx) & \text{if}\, i = n+1.
\end{cases}
\]
The degeneracy maps are given by 
\[
s_k^n(g_0, \cdots, g_{n-1}) = (g_0, \cdots, g_{k-1}, 1_G, g_k, \cdots, g_{n-1})
\]
in this same notation.
\begin{lemma}\label{Lemma: Simplicial Scheme of X}
	The data $(\underline{G \backslash X}_n, d_i^{n}, s_i^{n})_{n \in \N, 0 \leq i \leq n+1}$ determines a simplicial scheme $\underline{G \backslash X}_{\bullet}$.
\end{lemma}
\begin{proof}
	A routine verification using the fact that $G$ is a group scheme and that $X$ is a left $G$-variety shows that the simplicial identities hold.
\end{proof}
\begin{lemma}[\cite{JamesMracek}]\label{Lemma: Section 4: equivariant maps give simplicial morphisms between simplicial quotients}
	For any $G$-equivariant morphism $f \in \GSch(X,Y)$, there is an induced simplicial morphism of simplicial shcemes
	\[
	\underline{f}:\underline{G \backslash X}_{\bullet} \to \underline{G \backslash Y}_{\bullet}
	\]
	which is given by $\underline{f}_n := \id_{G^n} \times f$.
\end{lemma}
\begin{proof}
	To prove the lemma we first must prove that for any $n \in \N$ and any $0 \leq i \leq n+1$, the diagram
	\[
	\xymatrix{
		G^{n+1} \times X \ar[rr]^-{\id_{G^{n+1}} \times f} \ar[d]_{d_i^{n+1}} & & G^{n+1} \times Y \ar[d]^{d_{i}^{n+1}} \\
		G^{n} \times X \ar[rr]_-{\id_{G^n} \times f} & & G^n \times Y
	}
	\]
	commutes. If $i = 0$ then $d_{i}^{n+1} = \pi_{1, \cdots, n}^{G^{n+1}} \times \id_X$ and $\pi_{1,\cdots, n}^{G^{n+1}} \times \id_Y$, so the commutativity follows trivially. Analogously, if $i = n+1$, then $d_i^{n+1} = \id_{G^n} \times \alpha_X$ or $\id_{G^n} \times \alpha_Y$, so in this case the commutativity follows from the equivariance of $f$, i.e., because the diagram
	\[
	\xymatrix{
		G \times X \ar[d]_{\id_G \times f} \ar[r]^-{\alpha_X} & X \ar[d]^{f} \\
		G \times Y \ar[r]_-{\alpha_Y} & Y
	}
	\]
	commutes, so does the diagram:
	\[
	\xymatrix{
		G^n \times G \times X \ar[rr]^-{\id_{G^n} \times \alpha_X} \ar[d]_{\id_{G^{n+1}} \times f} & & G^n \times X \ar[d]^{\id_{G^n} \times f} \\
		G^n \times G \times Y \ar[rr]_-{\id_{G^n}\times \alpha_Y} & & G^{n} \times Y
	}
	\]
	If $0 < i < n+1$, the map $d_{i}^{n+1}$ multiplies the $(i-1)$-th and $i$-th components of the product of groups together and leaves both the $X$ and $Y$ components alone; in this case because $f_{n} = \id_{G^{n}} \times f$, the commutativity of the diagram is trivial to verify. This establishes that the $f_n$ commute appropriately with the face maps.
	
	Finally, we must verify that the $f_n$ commute appropriately with the degeneracy maps $s_i^{n}$, i.e., we must verify that the diagram
	\[
	\xymatrix{
		G^n \times X \ar[rr]^-{\id_{G^{n}} \times f} \ar[d]_{s_i^{n}} & & G^n \times Y \ar[d]^{s_i^{n}} \\
		G^{n+1} \times X \ar[rr]_-{\id_{G^{n+1}} \times f} & & G^{n+1} \times Y
	}
	\]
	commutes for all $n \in \N$ and all $0 \leq i \leq n$. However, this is immediate from the fact that each of the $s_i^{n}$ simply insert the unit of $G$ into the $i$-th component of the product $G^{n+1}$. This proves that the $f_n$ commute appropriately with the degeneracy maps and hence proves the lemma.
\end{proof}
\begin{corollary}\label{Cor: Section 4: Functor from equivariant schemes to simplicial schemes}
	There is a functor $\underline{(\cdot)}:\GVar \to [\mathbf{\Delta}^{\op},\Sch]$ given by sending objects $X \mapsto \underline{G \backslash X}_{\bullet}$ and by sending morphisms to $f \mapsto \underline{f}$.
\end{corollary}
\begin{remark}\label{Remark: Section 4: face maps are equivariant}
	Let $G$ be a group scheme and $X$ a left $G$-scheme. For all $n \in \N$ and all $0 \leq k \leq n+1$ the morphism $d_k^n:G^{n+1} \times X \to G^{n} \times X$ is $G$-equivariant. This is mainly by construction and the fact that $G$ is a group and $X$ is a left $G$-scheme; in fact, it is routine to show that the diagrams
	\[
	\xymatrix{
		G \times G^{n+1} \times X \ar[rr]^-{\alpha_{G^{n + 1} \times X}} \ar[d]_{\id_G \times d_k^n} & & G^{n + 1} \times X \ar[d]^{d_k^n} \\
		G \times G^{n} \times X \ar[rr]_-{\alpha_{G^n \times X}} & & G^n \times X
	}
	\]
	always commute. Combining this observation with Corollary \ref{Cor: Section 4: Functor from equivariant schemes to simplicial schemes} above allows us to note that there are induced simplicial morphisms 
	\[
	\underline{d}_k^n:\underline{G \backslash (G^{n+1} \times X)}_{\bullet} \to \underline{G \backslash(G^n \times X)} 
	\]
	for all $n \in \N$ and all $0 \leq k \leq n+1$.
\end{remark}
\begin{definition}[\cite{BernLun}, \cite{DeligneHodge3}]\label{Defn: Simplicial sheaves}\index{Simplicial!Sheaf! Equivariant}\index{Equivariant Derived Category! Simplicial}
	Fix a $G$-variety $X$. We say that a simplicial sheaf $\underline{\Fscr}$ on $\underline{G\backslash X}_{\bullet}$ is equivariant if each structure morphism, for all $h \in \mathbf{\Delta}_1$,
	\[
	\alpha_h:h^{\ast}\Fscr_{n} \to \Fscr_m
	\]
	is an isomorphism. We write $\Shv_{\text{eq}}(\underline{G \backslash X}_{\bullet};\overline{\Q}_{\ell})$\index[notation]{Shveq@$\Shv_{\text{eq}}(\underline{G \backslash X}_{\bullet};\overline{\Q}_{\ell})$} for the full subcategory of equivariant simplicial $\ell$-adic sheaves on $\underline{G\backslash X}_{\bullet}$. Similarly, we write $D^{b}_{\text{eq}}(\underline{G \backslash U}_{\bullet}; \overline{\Q}_{\ell})$\index[notation]{Dbeqsimplicial@$\Dbeqsimp{X}$} for the full subcategory of $D^{b}_c(\underline{G \backslash U}_{\bullet}; \overline{\Q}_{\ell})$ generated by complexes of equivariant simplical $\ell$-adic sheaves.
\end{definition}
\begin{remark}
	It is a routine check to show that $\Shv_{\text{eq}}(\underline{G\backslash X}_{\bullet};\overline{\Q}_{\ell})$ is a full Abelian subcategory of $\Shv(\underline{G \backslash X}_{\bullet};\overline{\Q}_{\ell})$ which induces a full inclusion $\Dbeqsimp{X} \hookrightarrow D^b_c(\underline{G \backslash X}_{\bullet};\overline{\Q}_{\ell})$.
\end{remark}

We now need the notion of what it means to have an augmentation of a simplicial scheme $\underline{X}_{\bullet}$ by a scheme $U$. While we will only need to use augmentations to prove that for $G$-schemes $U$ which admit a quotient scheme $G \backslash U$, there is an augmentation $\epsilon$ of the simplicial scheme $\underline{G \backslash U}_{\bullet}$ by $G \backslash U$ (as we will prove that this augmentation $\epsilon$ induces an equivalence of categories $D^b_{\text{eq}}(\underline{G \backslash U}_{\bullet};\overline{\Q}_{\ell}) \simeq D^{b}(G \backslash U; \overline{\Q}_{\ell})$). We will now need to discuss some general nonsense on augmentations of simplicial schemes (and in doing so we follow \cite{BarrKennisonRaphel}) in order to show that there is an augmentation of the form we will use.

\begin{definition}[{\cite[Page 475]{BarrKennisonRaphel}}]\label{Defn: Section 4: Augmentations of simplicial schemes}\index{Simplicial! Scheme! Augmented}\index{{Augmented Simplicial Scheme}|see {Simplicial, Scheme, Augmented}}
	An augmentation of a simplicial scheme $\underline{X}_{\bullet}$ by a scheme $U$ is a morphism of simplicial schemes $\underline{\epsilon}:\underline{X}_{\bullet} \to \const(U)_{\bullet}$.\index[notation]{epsilonunderline@$\underline{\epsilon}$}
\end{definition}
\begin{lemma}\label{Lemma: Section 4: How to define an augmentation of simplicial schemes}
	Let $\underline{X}$ be a simplicial scheme and let $U$ be a scheme. To give an augmentation $\underline{\epsilon}:\underline{X} \to \const(U)_{\bullet}$ it suffices to give  amorphism $\epsilon:X_0 \to U$ such that $\epsilon$ coequalizes the face maps $d_0^0$ and $d_1^0$, i.e., the diagram
	\[
	\xymatrix{
		X_1 \ar@<.5ex>[r]^{d_1^1} \ar@<-.5ex>[r]_{d_0^1} & X_0 \ar[r]^{\epsilon} & U
	}
	\]
	commutes in $\Sch$.
\end{lemma}
\begin{proof}
$\implies$: If there is a morphism of simplicial schemes $\underline{\epsilon_0}:\underline{X}_{\bullet} \to \const(U)_{\bullet}$, then because the diagram
\[
\begin{tikzcd}
X_1 \ar[d, shift left = 5, swap]{}{d_0^0} \ar[d, shift right = 5]{}{d_0^1} \ar[r]{}{\epsilon_1} & U \ar[d, equals] \\
X_0 \ar[r, swap]{}{\epsilon_0} & U
\end{tikzcd}
\]
taking $\epsilon := \epsilon_0$ does the job.

$\impliedby:$ Assume we have a morphism $\epsilon:X_0 \to U$ of schemes for which the diagram
\[
\xymatrix{
X_1 \ar@<.5ex>[r]^-{d_0^0} \ar@-.5ex>[r]_-{d_1^0} & X_0 \ar[r]^-{\epsilon} & U
}
\]
commutes. We then define our morphism $\underline{\epsilon}:\underline{X}_{\bullet} \to \const(U)_{\bullet}$ by
\[
\underline{\epsilon}_n := \begin{cases}
\epsilon & \text{if}\, n = 0; \\
\epsilon \circ d_0^0 \circ \cdots \circ d_0^{n-1} & \text{if}\, n \geq 1.
\end{cases}
\]
It is then a routine check to showw that for any $0 \leq k \leq n +1$ the diagrams
\[
\xymatrix{
X_{n+1} \ar[d]_{d_{k}^{n}} \ar[r]^-{\epsilon_{n+1}} & U \ar@{=}[d] \\
X_{n} \ar[r]_-{\epsilon_n} & U
}
\]
and
\[
\xymatrix{
X_n \ar[r]^-{\epsilon_n} \ar[d]_{s_k^n} & U \ar@{=}[d] \\
X_{n+1} \ar[r]_-{\epsilon_{n+1}} & U
}
\]
commute by induction on $n$, which completes the proof.
\end{proof}
\begin{example}
	Let $\underline{X}_{\bullet}$ be a simplicial $S$-scheme. Then the structure map $X_0 \to S$ allows us to realize $\const(S)$ as an augmentation of $\underline{X}_{\bullet}$.
\end{example}
\begin{lemma}\label{Lemma: Section 4: Augmentation of quotient simplicial scheme}
	For any $G$-scheme $U$ admitting a quotient scheme $G \backslash U$, there is an augmentation of the simplicial scheme $\underline{G \backslash U}_{\bullet}$ by $G \backslash U$.
\end{lemma}
\begin{proof}
	Following Lemma \ref{Lemma: Section 4: How to define an augmentation of simplicial schemes} it suffices to define a morphism $\epsilon:U \to G \backslash U$ which coequalizes $d_0^0$ and $d_1^0$. However, because $d_0^0 = \pi_2:G \times X \to X$ and $d_1^0 = \alpha_U:G \times U \to U$, we define $\epsilon:U \to G \backslash U$ by setting $\epsilon := \quo:U \to G\backslash U$. By definition of quotient varieties, however, this is automatic.
\end{proof}

Before proceeding to prove the equivalence of $D^b(G\backslash U;\overline{\Q}_{\ell})$ and $D^b(\underline{G \backslash U};\overline{\Q}_{\ell})$ we recall a result of Deligne in Section 6.1.2.b of \cite{DeligneHodge3}. It states that there is an equivalence of categories
\[
\Shv_{\text{eq}}(\underline{G \backslash X}_{\bullet}) \simeq \Shv_G^{\text{na{\"i}ve}}(X)
\]
in the topological case. However, this extends without issue to the schematic case when $X$ is a $G$-scheme that admits a quotient scheme $G \backslash X$ and iterating the eequivalence of sheaf categories through the $\ell$-adic construction induces the equivalence of categories below.
\begin{proposition}[\cite{DeligneHodge3}]\label{Prop: Section 4: equiv sheaves equiv}
	For any $G$-scheme $U$ which has a $G$-quotient scheme $G \backslash U$ and smooth quotient map $\quo_U:U \to G\backslash U$ , there is an equivalence of categories
	\[
	\Shv_{\text{eq}}(\underline{G \backslash U}_{\bullet}; \overline{\Q}_{\ell}) \simeq \Shv_G^{\text{na{\"i}ve}}(U; \overline{\Q}_{\ell}) \simeq \Shv(G \backslash U; \overline{\Q}_{\ell}).
	\]
\end{proposition}
\begin{proof}
	The equivalence $\Shv_{\text{eq}}(\underline{G \backslash U}_{\bullet}; \overline{\Q}_{\ell}) \simeq \Shv_G^{\text{na{\"i}ve}}(U; \overline{\Q}_{\ell})$ is what is described in \cite{DeligneHodge3} (with the translation to schematic and $\ell$-adic formalisms sketched above). Recall from Proposition \ref{Prop: Equiv of Standard t-structures on EDCs} we have an equivalence of categories 
	\[
	{}^{\text{stand}}\DbeqQl{X}^{\heartsuit} \simeq \ABLDbeqQl{X}^{\heartsuit_{\text{stand}}} \simeq \Shv_G^{\text{na{\"i}ve}}(X;\overline{\Q}_{\ell}).
	\]
	We now use the equivalence
	\[
	\DbQl{G \backslash U} \simeq \ABLDbeqQl{U}
	\]
	of \cite[Corollary 6.2.13]{PramodBook} (with the observation that all that is needed to give the equivalence is the smoothness of the quotient morphism $\quo_{U}:U \to G\backslash U$) together with the induced equivalence 
	\[
	{}^{\text{stand}}D^b_c(G \backslash U;\overline{\Q}_{\ell})^{\heartsuit} \simeq \Shv(G\backslash U;\overline{\Q}_{\ell})
	\] to give the chain of equivalences
	\begin{align*}
	\Shv_G^{\text{na{\"i}ve}}(U;\overline{\Q}_{\ell}) &\simeq \quot{D^b_G(U;\overline{\Q}_{\ell})^{\heartsuit_{\text{stand}}}}{ABL} \simeq {}^{\text{stand}}D^b_c(G\backslash U;\overline{\Q}_{\ell})^{\heartsuit} \\
	&\simeq \Shv(G\backslash U; \overline{\Q}_{\ell}),
	\end{align*}
	as was to be shown.
\end{proof}

\begin{proposition}\label{Prop: Section 4: Quot derived cat is equiv to simplicial derived cat}
	For any $G$-scheme $U$ which admits a quotient scheme $G \backslash U$ with smooth quotient map $\quo:U \to G \backslash U$, there is an equivalence of categories
	\[
	\Dbeqsimp{U}  \simeq D^b(G\backslash U;\overline{\Q}_{\ell}).
	\]
\end{proposition}
\begin{proof}
	As above, let $\underline{\epsilon} = \underline{\quo}:\underline{G \backslash U}_{\bullet} \to \const(G \backslash U)$ be the augmentation of Lemma \ref{Lemma: Section 4: Augmentation of quotient simplicial scheme} and recall from Proposition \ref{Prop: Section 4: Equiv of categories for constant simplicial scheme} that
	\[
	D^{b}(G \backslash U; \overline{\Q}_{\ell}) \simeq D^{b}(\const(G \backslash U); \overline{\Q}_{\ell}).
	\]
	Now we observe that the image of the simplicial pullback functor $\underline{\quo}:D^{b}(\const(G \backslash U );\overline{\Q}_{\ell}) \to D^{b}(\underline{G \backslash U}_{\bullet};\overline{\Q}_{\ell})$ takes values in $D_{\text{eq}}^{b}(\underline{G \backslash U}_{\bullet};\overline{\Q}_{\ell})$. Thus from the equivalence $\Shv_{\text{eq}}(\underline{G \backslash U}_{\bullet}; \overline{\Q}_{\ell}) \simeq \Shv(G \backslash U; \overline{\Q}_{\ell})$ in Proposition \ref{Prop: Section 4: equiv sheaves equiv} above, it suffices to prove that for any complexes $\underline{A}, \underline{B} \in D^{b}(\const(G\backslash U);\overline{\Q}_{\ell})_0$ and without loss of generality for any $n \in \N$, there is an isomorphism
	\[
	\Ext_{D^{b}(G\backslash U; \overline{\Q}_{\ell})}^{n}(A_0, B_0) \cong \Ext_{D^{b}_{\text{eq}}(\underline{G \backslash U}_{\bullet};\overline{\Q}_{\ell})}^{n}(\underline{\quo}^{\ast}\underline{A},\underline{\quo}^{\ast}\underline{B}).
	\]
	We now consider the augmentation $\underline{\quo}:\underline{G \backslash U}_{\bullet} \to \const(G \backslash U)$. For $n = 0$ we have that $\quo_0$ is given by $\quo:U \to G \backslash U$. If $n \geq 1$, then $\quo_n$ is inductively defined by the top-most horizontal arrow in the commuting diagram
	\[
	\xymatrix{
		X_{n} \ar[rr]^-{\quo_{n-1} \circ d_{0}^{n}} \ar[d]_{d^{n}_{0}} & & G \backslash U \ar@{=}[d] \\
		X_{n-1} \ar[rr]_-{\quo_{n-1}} & & G \backslash U
	}
	\]
	of schemes (cf.\@ the proof of Lemma \ref{Lemma: Section 4: How to define an augmentation of simplicial schemes}). Because the complexes $\underline{A}$ and $\underline{B}$ are equivariant and each morphism $\quo_{n}$ is smooth (at $n = 0$ this follows by assumption while for $n \geq 1$ this follows from a routine induction and the fact that $d_0^{n}$ is always smooth by the fact that $G$ is smooth) it suffices to prove in fact that for any $n \in \N$ there is an isomorphism
	\[
	\Ext_{\DbQl{G \backslash U}}^{n}(A_0,B_0) \cong \Ext_{\DbQl{U}}^{n}(\quo^{\ast}A_0, \quo^{\ast}B_0).
	\]
	However, because $\quo:U \to G \backslash U$ is smooth, it is flat and hence the derived pullback functor $\quo^{\ast}:D^b(G \backslash U; \overline{\Q}_{\ell}) \to D^{b}(U;\overline{\Q}_{\ell})$ is exact. Thus the isomorphism
	\[
	\Ext_{\DbQl{G \backslash U}}^{n}(A_0,B_0) \cong \Ext_{\DbQl{U}}^{n}(\quo^{\ast}A_0, \quo^{\ast}B_0)
	\]
	holds, which proves the proposition.
\end{proof}

With Proposition \ref{Prop: Section 4: Quot derived cat is equiv to simplicial derived cat}, we can now proceed to use the induced equivalences to explicitly prove the equivalence $\DbeqQl{X}$ with $\Dbeqsimp{X}$. However, while on the surface these constructions are quite similar, because each has as many different moving pieces we need a significant number of technical lemmas and tools to show how to move between the descent theoretic approach of $\DbeqQl{X}$ and the simplicial approach of $\Dbeqsimp{X}$. 

Perhaps the most notable obstacle to overcome when moving from $\DbeqQl{X}$ to $\Dbeqsimp{X}$ is that while the transition isomorphisms $T_{A}$ of an object $A$ has isomorphisms for most face maps (we still need to construct the transition isomorphisms coming from the face maps of the form $d_{n+1}^{n}$), it has no transition isomorphisms for any of the degeneracy maps $s_k^n:\quot{X}{G^n} \to \quot{X}{G^{n+1}}$. Another issue is that we even have to construct a sheaf on $X$ together with the transition isomorphisms $\tau_{d_0^0}, \tau_{d_0^1}, \tau_{s_0^0}$ in order to move towards $\Dbeqsimp{X}$. We will first show how to provide the degeneracy isomorphism $\tau_{s_k^n}$ by recording a straightforward lemma about why we are missing the degeneracy maps from $\Sf(G)_1$.

\begin{lemma}\label{Lemma: Section 4: Simplicial group scheme maps are not all Sf(G)}
	For any $n \in \N$ and any $0 \leq k \leq n$, the face map $d_k^{n}:G^{n+1} \times X \to G^n \times X$ is an $\Sf(G)$-morphism, while if $G$ is nontrivial the degeneracy maps $s_{\ell}^n:G^n \times X \to G^{n+1} \times X$ need not be $\Sf(G)$-morphisms.
\end{lemma}
\begin{proof}
	The first claim follows from the fact that since $G$ is a smooth algebraic group, as this implies that the multiplication map $\mu:G \times G \to G$ and projection morphisms are smooth with constant fibre dimension as well. The second claim follows from the fact that if $G$ is a nontrivial group, the embedding $1_G:X \to G \times X$ does not have constant fibre dimensions; it generically has a jump in dimension at the identity of $G$ and drops in dimension away from the identity.
\end{proof}
\begin{remark}
The reason why we cannot say outright that $s_{k}^{n}$ is not in $\Sf(G)_1$ is because when $G = \Spec K$ is the trivial group, each morphism $s_k^n$ is a smooth $G$-isomorphism and hence lies in $\Sf(G)_1$.
\end{remark}
\begin{lemma}\label{Lemma: Section 4: Presence of equivariant object gives degeneracy cocycle isos}
	Let $(A,T_A) \in D_G^b(X;\overline{\Q}_{\ell})_0$. Then for every $n \geq 1$ and every $0 \leq k \leq n$ there is an isomorphism
	\[
	(\overline{s}_{k}^{n})^{\ast}\quot{A}{G^{n+1}} \xrightarrow[\tau_{s_{k}^{n}}^{A}]{\cong} \quot{A}{G^{n}}
	\]
	which satisfies the cocycle condition with respect to morphisms in $\mathbf{\Delta}_1$ which do not contain a $d_{m+1}^{m}$ face map for any $m \in \N$.
\end{lemma}
\begin{proof}
	We first construct the isomorphism. First fix $s_k^n$ as in the statment of the lemma and let $d_{k}^n:G^{n+1} \times X$ be a face map. We then note that the triangle
	\[
	\xymatrix{
		G^{n} \times X \ar[r] \ar@{=}[dr] & G^{n+1} \times X \ar[d]^{d_k^n} \\
		& G^n \times X
	}
	\]
	commutes because $d_k^n:G^{n+1} \times X \to G^{n} \times X$ is given by multiplying the two adjacent copies of the group $G$ starting at the $k$-th index together and $s_k^n = \id_{G^{k}} \times 1_G \times \id_{G^{n-k}}$. Now note that the universal property of the quotient groups induces a morphism $\overline{s}_k^n:\quot{X}{G^n} \to \quot{X}{G^{n+1}}$ which further induces a commuting diagram
	\[
	\xymatrix{
		G^n \times X \ar[r]^-{s_k^n} \ar[d]_{\quo_{G^{n}}} & G^{n+1} \times X \ar[d]^{\quo_{G^{n+1}}} \\
		\quot{X}{G^{n}} \ar[r]_-{\overline{s}_k^n} & \quot{X}{G^{n+1}}	
	}
	\]
	of varieties. It then follows that the diagram
	\[
	\begin{tikzcd}
	& & \quot{X}{G^{n+1}} \ar[ddr]{}{\overline{d}_k^n}\\
	& G^{n+1} \times X \ar[ur]{}{\quo_{G^{n+1}}}  & \\
	& \quot{X}{G^n} \ar[rr, equals] \ar[uur, swap]{}{\overline{s}_k^n} & & \quot{X}{G^n}\\
	G^n \times X \ar[rr, equals] \ar[uur]{}{s_k^n} \ar[ur, swap]{}{\quo_{G^n}} & & G^n \times X \ar[ur, swap]{}{\quo_{G^n}}
	\ar[from = 2-2, to = 4-3, crossing over, swap, near end]{}{d_k^n}
	\end{tikzcd}
	\]
	commutes as well. Now write $d := d_k^n$ and consider that since $d_k^n$ is an $\Sf(G)$-morphism by Lemma \ref{Lemma: Section 4: Simplicial group scheme maps are not all Sf(G)}, there is an isomorphism $\tau_d^A \in T_A$ of the form
	\[
	\tau_d^A:\overline{d}^{\ast}\quot{A}{G^n} \to \quot{A}{G^{n+1}}.
	\]
	Applying the functor $(\overline{s}_k^n)^{\ast}$ to the isomorphism above yields a map
	\[
	(\overline{s}_k^n)^{\ast}\big(\overline{d}^{\ast}\quot{A}{G^{n}}\big) \xrightarrow[(\overline{s}_k^n)^{\ast}\tau_d^A]{\cong} (\overline{s}_k^n)^{\ast}\quot{A}{G^{n+1}}
	\]
	which satisfies the cocycle condition with respect to face maps $d_a^m$ where $m \in \N$ and $m \geq 1$ while $0 \leq a \leq m$, as it is the functorial image of such an isomorphism. Now observe that
	\[
	(\overline{s}_k^n)^{\ast}\big(\overline{d}^{\ast}\quot{A}{G^{n}}\big) = (\overline{d} \circ \overline{s}_k^n)^{\ast}\quot{A}{G^n} = (\id_{\quot{X}{G^n}})^{\ast}\quot{A}{G^{n}} = \quot{A}{G^n}
	\]
	so our isomorphism has domain and codomain
	\[
	\quot{A}{G^n} \xrightarrow[(\overline{s}_k^n)^{\ast}\tau_d^A]{\cong} (\overline{s}_k^n)^{\ast} \quot{A}{G^{n+1}}.
	\]
	We thus define our isomorphism $\tau_{s_k^n}^A:(\overline{s}_k^n)^{\ast}\quot{A}{G^{n+1}} \to \quot{A}{G^n}$ by
	\[
	\tau_{s_k^n}^{A} := (\overline{s}_k^n)^{\ast}(\tau_d^A)^{-1}.
	\]
	
	We now establish the cocycle condition for degneracy maps. Note that for this we have the diagram of morphisms
	\[
	\xymatrix{
		G^n \times X \ar@<1ex>[r]^-{s_k^n} & G^{n+1} \times X  \ar@<1ex>[l]^-{d_{k}^{n}} \ar@<1ex>[r]^-{s_r^{n+1}} & G^{n+2} \times X \ar@<1ex>[l]^-{d_{r}^{n+1}}	
	}
	\]
	where $0 \leq k \leq n$ and $0 \leq r \leq n+1$. From here, we will omit the subscripts in each morphism for readability. Now we calculate that
	\begin{align*}
	\tau_{s^{n+1} \circ s^{n}}^{A} &= (\overline{s}^{n+1} \circ \overline{s}^{n})^{\ast}(\tau_{d^{n} \circ d^{n+1}}^A)^{-1} = (\overline{s}^{n+1} \circ \overline{s}^{n})^{\ast}\left(\tau_{d^{n+1}}^A \circ (\overline{d}^{n+1})^{\ast}\tau_{d^{n}}^{A}\right)^{-1} \\
	&= \left((\overline{s}^{n+1} \circ \overline{s}^{n})^{\ast}\tau_{d^{n+1}}^{A} \circ (\overline{s}^{n+1} \circ \overline{s}^{n})^{\ast}(\overline{d}^{n+1})^{\ast}\tau_{d^n}^{A}\right)^{-1} \\
	&= \left((\overline{s}^n)^{\ast}\big((\overline{s}^{n+1})^{\ast}\tau_{d^{n+1}}^{A}\big) \circ (\overline{d}^{n+1} \circ \overline{s}^{n+1} \circ \overline{s}^{n})\tau_{d^{n}}^{A} \right)^{-1} \\
	&= \left((\overline{s}^{n})^{\ast}(\tau_{s^{n+1}}^{A})^{-1} \circ (\overline{s}^n)^{\ast}(\tau_{d^{n}}^{A})^{-1}\right)^{-1} \\
	&= \left((\overline{s}^{n})^{\ast}(\tau_{d^{n+1}}^A)^{-1} \circ (\tau_{s^{n}}^{A})^{-1}\right)^{-1} = \tau_{s^{n}}^{A} \circ (\overline{s}^n)^{\ast}\tau_{s^{n+1}}^{A}.
	\end{align*}
	
	We now must construct the isomorphisms
	\[
	\tau_{s_k^n \circ d_{r}^n}, \quad \tau_{d_{r}^{n} \circ s_k^n}
	\]
	in order to prove that these satisfy the cocycle condition. For this we note that by using the simplicial identities, it suffices to define exactly one of the $\tau_{s \circ d}$ or $\tau_{d \circ s}$, as one can always be converted to the other by playing combinatorial games and shifting indices $n$ appropriately. For this we fix $0 \leq k, r \leq n$. We now define the transition isomorphism
	\[
	(\overline{d}_r^n \circ \overline{s}_k^n)^{\ast}\quot{A}{G^{n}} \xrightarrow[\tau_{d_r^n \circ s_k^n}^{A}]{\cong} \quot{A}{G^n}
	\]
	by
	\[
	\tau_{d_r^n \circ s_k^n}^{A} := \tau_{s_k^n}^{A} \circ (\overline{s}_k^n)^{\ast}\tau_{d_r^n}^{A}.
	\]
	This satisfies the cocycle condition immediately and completes the proof of the lemma.
\end{proof}

We will now proceed to start constructing our functor $L:\DbeqQl{X} \to \Dbeqsimp{X}$ which we will later prove (cf.\@ Theorem \ref{Theorem: Section 4.2: Simplicial equivalence with derived category}) is one functor in a pair of inverse equivalences. Our basic strategy in constructing $L$ will be, at least to give the sheaves on $G^n \times X$ for $n \geq 1$, to go through the equivalences of Proposition \ref{Prop: Section 4: Quot derived cat is equiv to simplicial derived cat} to move from a complex $\quot{A}{G^n}$ on $\quot{X}{G^n}$ to a complex on $G^n \times X$. Let us explain this process before we proceed to the technical lemmas and propositions that show it works. Let $(A,T_A) \in D_G^b(X;\overline{\Q}_{\ell})_0$ and note that
\[
A = \lbrace \AGamma \; | \; \Gamma \in \Sf(G)_0 \rbrace
\]
where each $\AGamma \in D_c^b(G\backslash(\Gamma \times X);\overline{\Q}_{\ell})_0.$ Recall also that the quotient map $\quo_{\Gamma}:\Gamma \times X \to \XGamma$ is smooth as well by Proposition \ref{Prop: Smooth quotient from principal G var}. Now for each $\Gamma \in \Sf(G)_0$, let the equivalence of Proposition \ref{Prop: Section 4: Quot derived cat is equiv to simplicial derived cat} be denoted by:
\[
\begin{tikzcd}
D_{c}^{b}(\XGamma;\overline{\Q}_{\ell}) \ar[rr, bend left = 30, ""{name = U}]{}{L_{\Gamma}} & & \Dbeqsimp{(X \times \Gamma)} \ar[ll, bend left = 30, ""{name = L}]{}{R_{\Gamma}} \ar[from = U, to = L, symbol = \simeq]
\end{tikzcd}
\]\index[notation]{LGamma@$L_{\Gamma}$}\index[notation]{RGamma@$R_{\Gamma}$}
Under this assignment, we then get for every complex $\AGamma$ a simplicial complex of sheaves 
\[
(L_{\Gamma}(\AGamma)_n) \in \Dbeqsimp{(G^n \times X)}_0.
\] 
We would now like to use these functors to produce a family of equivariant simplicial complexes of $\ell$-adic sheaves
\[
\lbrace \AGamma \; | \; \Gamma \in \Sf(G)_0 \rbrace \mapsto \lbrace \quot{A}{G^n} \; | \; n \in \N, n \geq 1 \rbrace \mapsto \lbrace L_{G^n}(\quot{A}{G^n}) \; | \; n \in \N, n \geq 1 \rbrace
\]\index[notation]{LGothen@$L{G^n}$}
on each of the simplicial schemes $\underline{G \backslash (G^n \times X)}_{\bullet}$ (for $n \geq 1$) and then use these to finally build an equivariant simplicial complex of $\ell$-adic sheaves on $\underline{G \backslash X}_{\bullet}$ by taking the appropriate degree $0$ truncations. However, for this we also need to be able to produce a sequence of structure isomorphisms
\[
h^{\ast}L_{G^{n}}(\quot{A}{G^{n}}) \xrightarrow[\cong]{\alpha_h} L_{G^{m}}(\quot{A}{G^{m}})
\]
for every $h \in \mathbf{\Delta}_1$ with $\Dom h = [m]$ and $\Codom h = [n]$. The set of transition isomorphisms $T_A$, by Lemmas \ref{Lemma: Section 4: Simplicial group scheme maps are not all Sf(G)} and \ref{Lemma: Section 4: Presence of equivariant object gives degeneracy cocycle isos}, allow us to produce such isomorphisms whenever the morphisms $h$ involve no face maps of the form $d_{n+1}^{n}$ and neither of the maps $s_0^0$ or $d_0^0$. Our first main step towards showing that $T_A$ gives us most of the structure isomorphisms is to prove Proposition \ref{Prop: Section 4: Quot derived cat is equiv to simplicial derived cat} is natural at least for $\Sf(G)$-morphisms, degeneracy morphisms $s_k^n:G^n \times X \to G^{n+1} \times X$, and for the action face maps $d_{n+1}^{n}$.
\begin{lemma}\label{Lemma: Section 4: Naturality of Prop Quot Simp Equiv}
	Let $\Gamma, \Gamma^{\prime} \in \Sf(G)_0$ and let $f \in \Sf(G)(\Gamma,\Gamma^{\prime})$. Then for any $G$-variety $X$, the diagram of categories
	\[
	\xymatrix{
		D^b_c(\XGammap;\overline{\Q}_{\ell}) \ar[d]_{L_{\Gamma^{\prime}}} \ar[rr]^-{\overline{f}^{\ast}} & & D_c^b(\XGamma;\overline{\Q}_{\ell}) \ar[d]^{L_{\Gamma}} \\
		\Dbeqsimp{(\Gamma^{\prime} \times X)} \ar[rr]_-{\underline{f \times \id_X}} & & \Dbeqsimp{(\Gamma \times X)}
	}
	\]
	commutes strictly.
\end{lemma}
\begin{proof}
	This is a routine check that ultimately comes down to the following facts: The diagram
	\[
	\xymatrix{
		X \times \Gamma \ar[r]^-{f \times \id_X} \ar[d]_{\quo_{\Gamma}} & X \times \Gamma^{\prime} \ar[d]^{\quo_{\Gamma^{\prime}}} \\
		\XGamma \ar[r]_-{\overline{f}} &  \XGammap
	}
	\]
	commutes and is $G$-equivariant, $D_{c}^{b}(\XGammap;\overline{\Q}_{\ell}) \simeq D_{\text{eq}}^{b}(\const(\XGammap);\overline{\Q}_{\ell})$ (by Proposition \ref{Prop: Section 4: Equiv of categories for constant simplicial scheme}), and in this paper the pseudofunctor $D_c^b(-;\overline{\Q}_{\ell})$ is assumed to be strict.
\end{proof}
\begin{lemma}\label{Lemma: Section 4: Commutativity of face map pullbacks}
	Let $n \in \N$ with $n \geq 1$ and let $0 \leq k \leq n+1$. Then the diagram
	\[
	\xymatrix{
		D_c^b(\quot{X}{G^n};\overline{\Q}_{\ell}) \ar[rr]^-{(\overline{d}_k^n)^{\ast}} \ar[d]_{L_{G^n}} & & D_c^b(\quot{X}{G^{n+1}};\overline{\Q}_{\ell}) \ar[d]^{L_{G^{n+1}}} \\
		\Dbeqsimp{(G^n \times X)} \ar[rr]_-{(\underline{d}_k^n)^{\ast}} & & \Dbeqsimp{(G^{n+1} \times X)}
	}
	\]
	commutes strictly.
\end{lemma}
\begin{proof}
	Let $n \in \N$ with $n \geq 1$ and note that if $0 \leq k \leq n$, this follows from Lemmas \ref{Lemma: Section 4: Simplicial group scheme maps are not all Sf(G)} and \ref{Lemma: Section 4: Naturality of Prop Quot Simp Equiv} above. Now set $k = n+1$ and recall that $d_{n+1}^{n}$ is the map
	\[
	\id_{G^n} \times \alpha_X:G^{n+1} \times X \to G^{n} \times X.
	\]
	Then the functor
	\[
	(\underline{\id_{G^n} \times \alpha_X})^{\ast}:\Dbeqsimp{(G^n \times X)} \to \Dbeqsimp{(G^{n+1} \times X)}
	\]
	acts by pulling back against the $\alpha_X$ component in each level of an object $(A_n)$. Now let 
	\[
	\epsilon_n:\underline{G \backslash( G^n \times X)}_{\bullet} \to \const(G\backslash(G^n \times X))
	\] 
	be the augmentation of simplicial schemes for $G^n \times X$ and let 
	\[
	\epsilon_{n+1}:\underline{G \backslash (G^{n+1} \times X)}_{\bullet} \to \const(G \backslash(G^{n+1} \times X))
	\] 
	be the augmentation of simplicial schemes for $G^{n+1} \times X$. From Proposition \ref{Prop: Section 4: Quot derived cat is equiv to simplicial derived cat} we have that the functor $L_{G^n}$ acts on a complex $A \in D_c^b(\quot{A}{G^n};\overline{\Q}_{\ell})$ by
	\[
	L_{G^n}(A) = ((\epsilon_n^m)^{\ast}A)_{m \in \N},
	\]
	where $\epsilon_n^m$ is the degree $m$-component of the simplicial map $\underline{\epsilon}_n$; similarly, the functor $L_{G^n+1}$ acts on a complex $B \in D_c^b(\quot{X}{G^{n+1}};\overline{\Q}_{\ell})$ by
	\[
	L_{G^{n+1}}(B) = ((\epsilon_{n+1}^{m})^{\ast}B)_{m \in \N},
	\]
	where $\epsilon_{n+1}^m$ is the degree $m$-component of the simplicial map $\epsilon_{n+1}$. Thus we compute that on one hand we have
	\[
	L_{G^{n+1}}\big((\overline{d}_{n+1}^{n})^{\ast}A\big) = L_{G^{n+1}}\big((\overline{\id_{G^n} \times \alpha_X})^{\ast}A\big) = \big((\epsilon^m_{n+1})^{\ast}\big((\overline{\id_{G^n} \times \alpha_X})^{\ast}A\big)\big)_{m \in \N}
	\]
	while on the other hand we have
	\begin{align*}
	(\underline{d}_{n+1}^n)^{\ast}\big(L_{G^n}(A)\big) &= (\underline{d}_{n+1}^{n})^{\ast}\big((\epsilon_n^m)^{\ast}A\big)_{m \in \N}\\ 
	&= (\id_{G^m} \times \id_{G^n} \times \alpha_X)^{\ast} \big((\epsilon_n^m)^{\ast}A\big)_{m \in \N} \\ 
	&= (\id_{G^{m+n}} \times \alpha_X)^{\ast}\big((\epsilon_n^m)^{\ast}A\big)_{m \in \N}.
	\end{align*}
	If we can show that for any $m \in \N$
	\[
	(\id_{G^{n+m}} \times \alpha_X)^{\ast}\big((\epsilon_n^m)^{\ast}A\big) = (\epsilon_{n+1}^{m})^{\ast}\big((\overline{\id_{G^n} \times \alpha_X})^{\ast}A\big)
	\]
	we will be done. For this recall that the morphism $\epsilon_n^m$ is defined as $\epsilon^m_n = d^{\prime} \circ \quo_{G^{n}}$ where $d^{\prime}$ is the composition of face maps $d_0^{t + n}$ where $0 \leq t \leq m$. Similarly, $\epsilon_{n+1}^m:G^{m+n+1} \times X \to G\backslash (G^{n+1} \times X)$ is defined by $\epsilon_{n + 1}^{m}$ as $d^{\prime\prime} \circ \quo_{G^{n+1}} = \epsilon_{n+1}^{m}$ where $d^{\prime\prime}$ is the composition of face maps $d_0^{t +n + 1}$ where $0 \leq t \leq m$. Now consider that the diagrams
	\[
	\xymatrix{
		G^{m + n + 1} \times X \ar[d]_{\id_{G^{m+n}} \times \alpha_X} \ar[r]^-{d^{\prime\prime}} & G^{n+1} \times X \ar[d]^{\id_{G^n} \times \alpha_X} \\
		G^{m+n} \ar[r]_-{d^{\prime}} & G^{n} \times X
	}
	\]
	and
	\[
	\xymatrix{
		G^{n+1} \ar[r]^{\quo_{G^{n+1}}} \ar[d]_{\id_{G^n} \times \alpha_X} & \quot{X}{G^{n+1}} \ar[d]^{\overline{\id_{G^n} \times \alpha_X}} \\
		G^n \times X \ar[r]_-{\quo_{G^n}} & \quot{X}{G^{n}}
	}
	\]
	both commute by construction. Thus we find that the diagram
	\[
	\xymatrix{
		G^{m + n + 1} \times X \ar[d]_{\id_{G^{m+n}} \times \alpha_X} \ar[r]^-{d^{\prime\prime}} & G^{n+1} \ar[r]^{\quo_{G^{n+1}}} \ar[d]_{\id_{G^n} \times \alpha_X} & \quot{X}{G^{n+1}} \ar[d]^{\overline{\id_{G^n} \times \alpha_X}} \\
		G^{m+n} \ar[r]_-{d^{\prime}} & G^n \times X \ar[r]_-{\quo_{G^n}} & \quot{X}{G^{n}}
	}
	\]
	commutes; however, this diagram factors the diagram
	\[
	\xymatrix{
		G^{m+n+1} \ar[r]^-{\epsilon_{n+1}^{m}} \ar[d]_{\id_{G^{m+n}} \times \alpha_X} & \quot{X}{G^{n+1}} \ar[d]^{\overline{\id_{G^{n}} \times \alpha_X}} \\
		G^{m+n} \times X \ar[r]_-{\epsilon_{n}^{m}} & \quot{X}{G^n}
	}
	\]
	so we derive that
	\[
	\epsilon_{n}^m \circ (\id_{G^{m+n}} \times \alpha_X) = \overline{\id_{G^{n}} \times \alpha_X} \circ \epsilon_{n+1}^{m}.
	\]
	From this and the strictness of the pseudofunctor $D^b_c(-;\overline{\Q}_{\ell})$ we find that
	\[
	(\id_{G^{n+m}} \times \alpha_X)^{\ast}\big((\epsilon_n^m)^{\ast}A\big) = (\epsilon_{n+1}^{m})^{\ast}\big((\overline{\id_{G^n} \times \alpha_X})^{\ast}A\big)
	\]
	for any $A \in D_c^b(\quot{X}{G^n};\overline{\Q}_{\ell})$, as was desired.
\end{proof}
\begin{lemma}\label{Lemma: Section 4: Simplicial degeneracies are natural with equiv}
	Let $n \in \N$ with $n \geq 1$ and let $0 \leq k \leq n$. Then the diagram
	\[
	\xymatrix{
		D^b_c(\quot{X}{G^{n+1}};\overline{\Q}_{\ell}) \ar[d]_{L_{G^{n+1}}} \ar[rr]^-{(\overline{s}_{k}^{n})^{\ast}} & & D^b_c(\quot{X}{G^n};\overline{\Q}_{\ell}) \ar[d]^{L_{G^n}} \\
		\Dbeqsimp{(G^{n+1} \times X)} \ar[rr]_-{(\underline{s}_{k}^{n})^{\ast}} & & \Dbeqsimp{(G^n \times X)}
	}
	\]
	commutes strictly.
\end{lemma}
\begin{proof}
	Begin by observing that the morphism
	\[
	\overline{s}_{k}^{n}:\quot{X}{G^n} \to \quot{X}{G^{n+1}}
	\]
	is defined via the commuting diagram:
	\[
	\xymatrix{
		G^{n} \times X \ar[rr]^-{\quo_{G^n}} \ar[d]_{s_k^n} & & \quot{X}{G^n} \ar[d]^{\overline{s}_k^n} \\
		G^{n+1} \times X \ar[rr]_-{\quo_{G^{n+1}}} & & \quot{X}{G^{n+1}}
	}
	\]
	As in Lemma \ref{Lemma: Section 4: Commutativity of face map pullbacks}, write $\epsilon_n:\underline{G \backslash(G^n \times X)}_{\bullet} \to \const(G \backslash (G^n \times X))$. We then observe that on one hand we have
	\[
	L_{G^n}\big((\overline{s}_k^n)^{\ast}\quot{A}{G^{n+1}}\big) = \left((\epsilon_n^m)^{\ast}\big((\overline{s}_k^n)^{\ast}\quot{A}{G^{n+1}}\big)\right))_{m \in \N}.
	\]
	while on the other we have
	\[
	(\underline{s}_k^n)^{\ast}L_{G^{n+1}}(\quot{A}{G^{n+1}}) = \left((\id_{G^m} \times s_{\ell}^{n})^{\ast}(\epsilon_{n+1}^{m})^{\ast}\quot{A}{G^{n+1}}\right)_{m \in \N}.
	\]
	Now for each $m \in \N$ we recall that $\epsilon_{n + 1}^m$ is the morphism
	\[
	\xymatrix{
		G^{m+n+1} \ar[r]^-{\pi_m} & G^{n+1} \ar[rr]^-{\quo_{G^{n+1}}} & & \quot{X}{G^{n+1}}
	}
	\]
	and the $m$-th component of the map $\underline{s}_k^n$ is the morphism
	\[
	\id_{G^{m}} \times s_k^n:G^{m+n} \times X \to G^{m+n+1} \times X.
	\]
	A routine check shows that the diagram
	\[
	\xymatrix{
		G^{m + n} \times X \ar[r]^-{\pi_m} \ar[d]_{\id_{G^m}\times s_k^n} & G^{n} \times X \ar[d]^{s_k^n} \\
		G^{m + n + 1} \times X \ar[r]_-{\pi_m} & G^{n+1} \times X	
	}
	\]
	commutes, which in turn implies that
	\[
	\xymatrix{
		G^{m + n} \times X \ar[r]^-{\pi_m} \ar[d]_{\id_{G^m}\times s_k^n}	& G^{n} \times X \ar[rr]^-{\quo_{G^n}} \ar[d]_{s_k^n} & & \quot{X}{G^n} \ar[d]^{\overline{s}_k^n} \\
		G^{m+n+1} \times X \ar[r]_-{\pi_m} &	G^{n+1} \times X \ar[rr]_-{\quo_{G^{n+1}}} & & \quot{X}{G^{n+1}}
	}
	\]
	commutes as well. However, the above diagram has outer square 
	\[
	\xymatrix{
		G^{m+n} \times X \ar[r]^-{\epsilon_{n}^{m}} \ar[d]_{\id_{G^m} \times s_k^n} & \quot{X}{G^{n}} \ar[d]^{\overline{s}_k^n} \\
		G^{m + n + 1} \times X \ar[r]_-{\epsilon_{n+1}^{m}} & \quot{X}{G^{n+1}}
	}
	\]
	so this implies that
	\[
	\overline{s}_k^n \circ \epsilon_n^m = \epsilon_{n+1}^{m} \circ (\id_{G^m} \times s_k^n)
	\]
	for all $m \in \N$. Thus we conclude that
	\begin{align*}
	(\underline{s}_k^n)^{\ast}L_{G^{n+1}}(\quot{A}{G^{n+1}}) &= \left((\id_{G^m} \times s_{\ell}^{n})^{\ast}(\epsilon_{n+1}^{m})^{\ast}\quot{A}{G^{n+1}}\right)_{m \in \N} \\
	&= \left((\epsilon_n^m)^{\ast}\big((\overline{s}_k^n)^{\ast}\quot{A}{G^{n+1}}\big)\right))_{m \in \N} \\
	&= L_{G^n}\big((\overline{s}_k^n)^{\ast}\quot{A}{G^{n+1}}\big),
	\end{align*}
	which proves the lemma.
\end{proof}

\begin{proposition}\label{Prop: Section 4: Simplicial cocycle structure isos}
	For any $n \in \N$ with $n \geq 1$ and any $0 \leq k \leq n$, given face map $d = d_k^n:G^{n+1} \times X \to G^{n} \times X$, then for any $(A,T_A) \in D_G^b(X;\overline{\Q}_{\ell})$ there is an isomorphism
	\[
	\underline{d}^{\ast}L_{G^{n}}(\quot{A}{G^n}) \xrightarrow[\alpha_d]{\cong} L_{G^{n+1}}(\quot{A}{G^{n+1}}).
	\]
	Furthermore, given the degeneracy maps $s = s_k^n$ there are also isomorphisms
	\[
	\underline{s}^{\ast}:L_{G^{n+1}}(\quot{A}{G^{n+1}}) \xrightarrow[\alpha_s]{\cong} L_{G^{n}}(\quot{A}{G^{n}}).
	\]
	Finally, the $\alpha_d$ and $\alpha_s$ satisfy the cocycle condition with respect to morphisms $\underline{h}$ that come from $\mathbf{\Delta}$ that do not contain a composition of a $d_{n+1}^{n}$, $d_k^0$, or $s_0^0$ map.
\end{proposition}
\begin{proof}
	Fix an $n, k \in \N$ with $n \geq 1$ and $0 \leq k \leq n$. In this case the face map $d_k^n$ is an $\Sf(G)$-morphism by Lemma \ref{Lemma: Section 4: Simplicial group scheme maps are not all Sf(G)}, so applying Lemma \ref{Lemma: Section 4: Commutativity of face map pullbacks} gives that for $d = d_k^n$,
	\[
	L_{G^{n+1}} \circ \overline{d}^{\ast} = \underline{d}^{\ast} \circ L_{G^{n}}.
	\]
	Now consider the transition isomorphism $\tau_d^A:\overline{d}^{\ast}\quot{A}{G^n} \to \quot{A}{G^{n+1}}$ in $T_A$. Applying the functor $L_{G^{n+1}}$ gives an isomorphism
	\[
	L_{G^{n+1}}(\tau_d^A):L_{G^{n+1}}\big(\overline{d}^{\ast}\quot{A}{G^{n}}\big) \xrightarrow{\cong} L_{G^{n+1}}(\quot{A}{G^{n+1}}),
	\]
	which by is in fact an isomorphism
	\[
	L_{G^{n+1}}(\tau_d^A):\underline{d}^{\ast}\big(L_{G^{n}}(\quot{A}{G^{n}})\big) \xrightarrow{\cong} L_{G^{n+1}}(\quot{A}{G^{n+1}}).
	\]
	We thus make the definition
	\[
	\alpha_d := L_{G^{n+1}}(\tau_d^A).
	\]
	Now let $h \in \mathbf{\Delta}_1$ be a composite of face maps $d^{m}_{\ell}$, for $m \geq 1$ and $0 \leq \ell \leq m$. We first assume that $d_{\ell}^{m} \circ d$ is defined and show the cocycle condition; in the case where $d \circ d_{\ell}^{m}$ is defined, the cocycle condition follows mutatis mutandis. Define the transition isomorphism
	\[
	\alpha_{d^m_{\ell} \circ d} := L_{G^{n+1}}(\tau_{d_{\ell}^{m} \circ d}^A)
	\] 
	and make a similar definition when $d \circ d_{\ell}^{m}$ is defined. In our case note that we have that the equation
	\[
	\tau_{d_{\ell}^m \circ d}^A = \tau_{d_{\ell}^{m}}^{A} \circ (\overline{d}_{\ell}^m)^{\ast}\tau_{d}^{A}
	\]
	holds by virtue of the transition isomorphisms being in $T_A$. Now applying the functor $L_{G^{m+1}}$ and using Lemma \ref{Lemma: Section 4: Commutativity of face map pullbacks} gives
	\begin{align*}
	L_{G^{m+1}}(\tau_{d_{\ell}^m \circ d}^{A}) &= L_{G^{m+1}}\left(\tau_{d_{\ell}^{m}}^{A} \circ (\overline{d}_{\ell}^m)^{\ast}\tau_{d}^{A}\right) = L_{G^{m+1}}(\tau_{d_{\ell}^{m}}) \circ L_{G^{m+1}}\left((\overline{d}_{\ell}^m)^{\ast}\tau_{d}^{A}\right) \\
	&= \alpha_{d_{\ell}^{m}} \circ (\underline{d}_{\ell}^{m})^{\ast}\left(L_{G^{n+1}}(\tau_{d}^{A})\right) = \alpha_{d_{\ell}^{m}} \circ (\underline{d}_{\ell}^{m})^{\ast}\alpha_{d},
	\end{align*}
	as was to be verified.
	
	To construct the structure isomorphisms for the degeneracy maps, we proceed as above. Fix an $n \in \N$ with $n \geq 1$ and let $0 \leq t \leq n$. Now define $s = s_t^n$ and note that by Lemma \ref{Lemma: Section 4: Presence of equivariant object gives degeneracy cocycle isos} there are isomorphisms
	\[
	\tau_{s}^{A}:(\overline{s})^{\ast}\quot{A}{G^{n+1}} \xrightarrow{\cong} \quot{A}{G^{n}}.
	\]
	Applying the functor $L_{G^{n}}$ gives an isomorphism
	\[
	L_{G^{n}}\big((\overline{s})^{\ast}\quot{A}{G^{n+1}}\big) \xrightarrow[L_{G^{n}}(\tau_s^A)]{\cong} L_{G^{n}}(\quot{A}{G{^n}}),
	\]
	which by Lemma \ref{Lemma: Section 4: Simplicial degeneracies are natural with equiv} is an isomorphism
	\[
	\underline{s}^{\ast}L_{G^{n+1}}(\quot{A}{G^{n+1}}) \xrightarrow[L_G{^{n}}(\tau_{s}^{A})]{\cong} L_{G^n}(\quot{A}{G^n}).
	\]
	We thus define the structure isomorphism $\alpha_{s}:\underline{s}^{\ast}L_{G^{n+1}}(\quot{A}{G^{n+1}}) \to L_{G^{n}}(\quot{A}{G^{n}})$ by setting
	\[
	\alpha_s := L_{G^{n}}(\tau_s^A).
	\]
	We now verify the cocycle condition of the $\alpha_{s}$ with respect to other degeneracy maps. For this consider $s_r^m$ for $m \in \N$ with $m \geq 1$ and $0 \leq r \leq m$ and assume that $s^m_r \circ s$ is defined; the case $s \circ s_r^m$ is defined is omitted. Define
	\[
	\alpha_{s_r^m \circ s} := L_{G^{n}}(\tau_{s_r^m \circ s})
	\]
	and make the similar definition for $\alpha_{s \circ s^{\ell}_{i}}$. We now compute that, using Lemma \ref{Lemma: Section 4: Presence of equivariant object gives degeneracy cocycle isos} again,
	\begin{align*}
	\alpha_{s_r^m \circ s} &= L_{G^n}(\tau_{s_r^m \circ s}^A) = L_{G^{n}}(\tau_{s}^A \circ \overline{s}^{\ast}\tau_{s_r^m}^A) = L_{G^n}(\tau_s^A) \circ L_{G^n}(\overline{s}^{\ast}\tau_{s^m_r}^{A}) \\
	&=\alpha_s \circ \underline{s}^{\ast}L_{G^m}(\tau_{s_r^m}^A) = \alpha_s \circ \underline{s}^{\ast}\alpha_{s_r^m},
	\end{align*}
	as was to be verified.
	
	We finally now verify the cocycle condition for mixed composites of face and degeneracy maps. For this we assume that $d \circ s$ is defined; in this case we have the diagram
	\[
	\xymatrix{
		G^{n+1} \times X \ar@<1ex>[r]^-{d} & G^{n} \times X \ar@<1ex>[l]^-{s}
	}
	\]
	for some $n \geq 1$. We then define
	\[
	\alpha_{d \circ s} := L_{G^n}(\tau_{d \circ s}^A), \qquad \alpha_{s \circ d} := L_{G^{n+1}}(\tau_{s \circ d}^A).
	\] 
	By Lemma \ref{Lemma: Section 4: Presence of equivariant object gives degeneracy cocycle isos} we have that
	\[
	\tau_{d \circ s}^A = \tau_s^A \circ \overline{s}^{\ast}\tau_d^A,
	\]
	so applying the functor $L_{G^n}$ and using Lemmas \ref{Lemma: Section 4: Presence of equivariant object gives degeneracy cocycle isos}, \ref{Lemma: Section 4: Commutativity of face map pullbacks}, and \ref{Lemma: Section 4: Simplicial degeneracies are natural with equiv} we calculate that
	\begin{align*}
	\alpha_{d \circ s} &= L_{G^n}(\tau_{d \circ s}^A) = L_{G^n}(\tau_{s}^A \circ \overline{s}^{\ast}\tau_d^A) = L_{G^n}(\tau_s^A) \circ L_{G^n}(\overline{s}^{\ast}\tau_d^A) \\
	&= \alpha_s \circ \underline{s}^{\ast}L_{G^{n+1}}(\tau_d^A) = \alpha_s \circ \underline{s}^{\ast}\alpha_{d}.
	\end{align*}
	Similarly,
	\begin{align*}
	\alpha_{s \circ d} &= L_{G^{n+1}}(\tau_{s \circ d}^A) = L_{G^{n+1}}(\tau_d^A \circ \overline{d}^{\ast}\tau_s^A) = L_{G^{n+1}}(\tau_d^A) \circ L_{G^{n+!}}(\overline{d}^{\ast}\tau_s^A) \\
	&= \alpha_d \circ \underline{d}^{\ast}L_{G^n}(\tau_s^A) = \underline{d}^{\ast}\alpha_s,
	\end{align*}
	as was to be shown. This verifies the cocycle condition and proves the proposition.
\end{proof}

We now begin our study of how to extend our structure isomorphisms $\alpha_h$ to include face maps of the form $d_{n+1}^{n}$. This ultimately comes down to an induction, so in this process we will also provide the structure isomorphisms for the face and degeneracy maps $d_0^0, d_0^1$, and $s_0^0$. The next lemma should be thought of as the base case for establishing the structure isomorphisms at the $d_{n+1}^{n}$ face maps in Proposition \ref{Prop: Section 4: Action map simplicial isos in deg 0}.
\begin{lemma}\label{Lemma: Section 4: Object in LG0}
	Let $A \in D_G^b(X;\overline{\Q}_{\ell})$. There is then an object $A_0$ of $D^b_c(X;\overline{\Q}_{\ell})$ such that if $L_G(\quot{A}{G})_0$ denotes the degree $0$ component of $L_G(\quot{A}{G})$,
	\[
	(d_0^0)^{\ast}A_0 = L_G(\quot{A}{G})_0 = (d_1^0)^{\ast}A_0, \quad (s_0^0)^{\ast}L_G(\quot{A}{G})_0 = A_0.
	\]
\end{lemma}
\begin{proof}
	Begin by recalling that there is a non-$G$-equivariant isomorphism of varieties
	\[
	\psi:X \xrightarrow{\cong} \quot{X}{G}
	\]
	which equips $X$ with the trivial $G$-action. Consequently, it follows that the diagrams
	\[
	\xymatrix{
		& G \times X  \ar[dl]_{\quo_{G}} \ar@<.5ex>[dr]^{d_0^0} \ar@<-.5ex>[dr]_{d_1^0} & \\
		\quot{X}{G}  	& & X \ar[ll]^-{\psi}
	}
	\]
	\[
	\xymatrix{
		& G \times X \ar[dl]_{\quo_{G}} & \\
		\quot{X}{G} & & X \ar[ll]^-{\psi} \ar[ul]_{s_0^0}	
	}
	\]
	both commute. Now define
	\[
	A_0 := \psi^{\ast}\quot{A}{G} \in D_c^b(X;\overline{\Q}_{\ell}).
	\]
	We then calculate that
	\begin{align*}
	(d_0^0)^{\ast}A_0 &= (d_0^0)^{\ast}\psi^{\ast}\quot{A}{G} = (\psi \circ d_0^0)^{\ast}\quot{A}{G} =\quo_{G}^{\ast}\quot{A}{G} = L_G(\quot{A}{G})_0 \\
	& = \quo_{G}^{\ast}\quot{A}{G} = (\psi \circ d_1^0)^{\ast}\quot{A}{G} = (d_1^0)^{\ast}\psi^{\ast}\quot{A}{G} = (d_1^0)^{\ast}A_0.
	\end{align*}
	In the other direction we get that
	\[
	(s_0)^{\ast}L_G(\quot{A}{G})_0 = (s_0^0)^{\ast}\quo_G^{\ast}\quot{A}{G} = (\quo_G \circ s_0^0)^{\ast}\quot{A}{G} = \psi^{\ast}\quot{A}{G} = A_0,
	\]
	which proves the lemma.
\end{proof}
\begin{proposition}\label{Prop: Section 4: Action map simplicial isos in deg 0}
	Fix an object $(A,T_A) \in D_G^b(X;\overline{\Q}_{\ell})_0$ and let $L_{G^{n}}(\quot{A}{G}^n)_0$ denote the $0$-th degree object of $L_{G^n}(\quot{A}{G^n})$. Then for any $n \in \N$, there are structure isomorphisms
	\[
	\alpha_{d_{n+1}^{n}}:(d_{n+1}^{n})^{\ast}L_{G^n}(\quot{A}{G^n})_0 \xrightarrow{\cong} L_{G^{n+1}}(\quot{A}{G^{n+1}})_0,
	\]
	where $L_{G^0}(\quot{A}{G^0})$ is defined to be the object $A_0$ of Lemma \ref{Lemma: Section 4: Object in LG0}. Furthermore, there are structure isomorphisms
	\[
	\begin{cases}
	(s_t^n \circ d_{n+1}^{n})^{\ast}L_{G^n+1}(\quot{A}{G^{n+1}})_0 \xrightarrow[\alpha_{s_t^n \circ d_{n+1}^{n}}]{\cong} L_{G^{n+1}}(\quot{A}{G^{n+1}})_0; & \text{for}\,\, 0 \leq t \leq n; \\
	(d_{n+1}^{n} \circ s_t^n)^{\ast}L_{G^{n}}(\quot{A}{G^n})_0 \xrightarrow[\alpha_{d_{n+1}^{n} \circ s_t^n}]{\cong} L_{G^{n}}(\quot{A}{G^n})_0 & \text{for}\,\, 0 \leq t \leq n;\\
	(d_{k}^{n-1} \circ d_{n+1}^{n})^{\ast}L_{G^{n-1}}(\quot{A}{G^{n-1}})_0 \xrightarrow[\alpha_{d_{n+1}^{n} \circ d_{k}^{n}}]{\cong} L_{G^{n+1}}(\quot{A}{G^{n+1}})_0 & \text{if}\,\, n \geq 1,\\
	& \text{and}\,\, 0 \leq k \leq n;\\
	(d_{n+1}^{n}\circ d_r^{n+1})^{\ast}L_{G^{n}}(\quot{A}{G^n})_0 \xrightarrow[\alpha_{d_r^{n+1} \circ d_{n+1}^n}]{\cong} L_{G^{n+2}}(\quot{A}{G^{n+2}})_0 & \text{if}\,\, 0 \leq r \leq n+1;
	\end{cases}
	\]
	which make it so that $\alpha_{d_{n+1}^{n}}$ satisfies the cocycle condition with respect to the morphisms from $\mathbf{\Delta}_1$.
\end{proposition}
\begin{proof}
	We prove this via induction on $n$. For this first note that the base case $n = 0$ is given in Lemma \ref{Lemma: Section 4: Object in LG0}. The only thing yet to be done is define thee isomorphisms $\alpha_{d_{1}^{0} \circ d_{k}^{1}}:(d_1^0 \circ d_k^1)^{\ast}A_0 \xrightarrow{\cong} \quot{A}{G^2}$ for $k = 0, 1$. These are defined via the diagram
	\[
	\xymatrix{
		(d_1^0 \circ d_k^1)^{\ast}A_0 \ar@{=}[r] \ar[dr]_{\alpha_{d_1^0 \circ d_k^1}} & (d_k^1)^{\ast}(d_1^0)^{\ast}A_0 \ar@{=}[r] & (d_k^1)^{\ast}(\epsilon_1^0)^{\ast}\quot{A}{G} \ar@{=}[d] \\
		& L_{G^{2}}(\quot{A}{G^2})_0 & (d_k^1)^{\ast}L_{G}(\quot{A}{G})_0 \ar[l]^{\alpha_{d_k^1}}
	}
	\]
	for $k = 0, 1$. That the satisfy the cocycle condition is immediate.
	
	We now proceed via induction. Let $k \in \N$ and assume that there is an isomorphism
	\[
	\alpha_{d_{m+1}^{m}}:L_{G^k}(\quot{A}{G^{m}})_0 \xrightarrow{\cong} L_{G^{m+1}}(\quot{A}{G^{m+1}})_0,
	\]
	as well as all other structure isomorphisms claimed in the statement of the lemma which satisfy the cocycle conditions as claimed for relevant morphisms in $\mathbf{\Delta}_1$ for all $0 \leq m \leq k$. Let us begin by finding the structure isomorphism from $(d_{k+2}^{k+1})^{\ast}L_{G^{k+1}}(\quot{A}{G^{k+1}})_0$ to $L_{G^{k+2}}(\quot{A}{G^{k+2}})_0$. For this we begin by observing that the cube
	\[
	\begin{tikzcd}
	& G^{k+1} \times X \ar[rr]{}{d_{k+1}^{k}} \ar[dd, near start]{}{\quo_{G^{k+1}}} & & G^{k} \times X \ar[dd]{}{\quo_{G^k}} \\
	G^{k+2} \times X \ar[ur]{}{d_{0}^{k+1}} \ar[dd, swap]{}{\quo_{G^{k+2}}}  & & G^{k+1} \times X \ar[ur, swap]{}{d_{0}^{k}} \\
	& \quot{X}{G^{k+1}} \ar[rr, near end]{}{\overline{d}_{k+1}^{k}} & & \quot{X}{G^{k}}\\
	\quot{X}{G^{k+2}} \ar[ur]{}{\overline{d}_0^{k+1}} \ar[rr, swap]{}{\overline{d}_{k+2}^{k+1}} & & \quot{X}{G^{k+1}} \ar[ur, swap]{}{\overline{d}_0^k}
	\ar[from = 2-1, to = 2-3, crossing over, near start]{}{d_{k+2}^{k+1}} \ar[from = 2-3, to = 4-3, crossing over, swap, near start]{}{\quo_{G^{k+1}}}
	\end{tikzcd}
	\]
	commutes whenever $k \geq 1$; if $k = 0$ then we use the map $\psi:X \to \quot{X}{G}$ described in Lemma \ref{Lemma: Section 4: Object in LG0} instead of $\quo_{G^k}$, as it is unlikely that $G \backslash X$ exists. Appealing to Lemma \ref{Lemma: Section 4: Commutativity of face map pullbacks} if $k \geq 1$ or giving a manipulation mutatis mutandis to Lemma \ref{Lemma: Section 4: Object in LG0} if $k = 0$, we deduce the existence of structure isomorphisms
	\[
	\alpha_{d^{m}_i}:(d_i^m)^{\ast}L_{G^{m}}(\quot{A}{G^m})_0 \to L_{G^{m+1}}(\quot{A}{G^{m+1}})_0
	\]
	for all $0 \leq m \leq k-1$ and for all $0 \leq i \leq m+1$. We now define $\alpha_{d_{k+2}^{k+1}}$ by the following diagram:
	\[
	\xymatrix{
		(d_{k+2}^{k+1})^{\ast}L_{G^{k+1}}(\quot{A}{G^{k+1}})_0 \ar[rr]^-{(d_{k+2}^{k+1})^{\ast}\alpha_{d_0^k}^{-1}}\ar[dd]_{\alpha_{k+2}^{k+1}} & & (d_{k+2}^{k+1})^{\ast}(d_0^{k})^{\ast}L_{G^{k}}(\quot{A}{G^k})_0 \ar@{=}[d] \\
		& & (d_0^k \circ d_{k+2}^{k+1})^{\ast}L_{G^{k}}(\quot{A}{G^k})_0 \ar@{=}[dd] \\
		L_{G^{k+2}}(\quot{A}{G^{k+2}}) & & \\
		& & (d_{k+1}^{k} \circ d_0^{k+1})^{\ast}L_{G^k}(\quot{A}{G^k})_0 \ar@{=}[d] \\
		(d_0^{k+1})^{\ast}L_{G^{k+1}}(\quot{A}{G^{k+1}})_0 \ar[uu]^{\alpha_{d_0^{k+1}}} & & (d_0^{k+1})^{\ast}(d_{k+1}^{k})^{\ast}L_{G^{k}}(\quot{A}{G^k})_0 \ar[ll]^-{(d_0^{k+1})^{\ast}\alpha_{d_{k+1}^{k}}}
	}
	\]
	More explicitly, we define
	\[
	\alpha_{d_{k+2}^{k+1}} := \alpha_{d_0^{k+1}} \circ (d_0^{k+1})^{\ast}\alpha_{d_{k+1}^{k}} \circ (d_{k+2}^{k+1})^{\ast}\alpha_{d_{0}^{k}}^{-1};
	\]
	note that we use the inductive hypothesis in using the isomorphism $\alpha_{d_{k+1}^{k}}$.
	
	We now must construct the isomorphisms that make $\alpha_{d_{k+2}^{k+1}}$ satisfy the coycle condition. For this we begin by considering the composite $d_{k+2}^{k+1} \circ s_t^{k+1}$ for any $0 \leq t \leq k+1$. Now observe that by the simplicial identities we have that
	\[
	d_{k+2}^{k+1} \circ s_t^{k+1} = \begin{cases}
	\id & \text{if}\, t = k+1; \\
	s^{k}_{t} \circ d_{k+1}^{k} & \text{if}\, 0 \leq t \leq k.
	\end{cases}
	\]
	In this case we define
	\[
	\alpha_{d_{k+2}^{k+1} \circ s_t^{k+1}} = \begin{cases}
	\id & \text{if}\, t = k+1; \\
	\alpha_{s^k_t \circ d_{k+1}^{k}} & \text{if}\, 0 \leq t \leq k;
	\end{cases}
	\]
	and note that in both cases $\alpha_{d_{k+2}^{k+1} \circ s_k^t}$ satisfies the cocycle condition (either by the fact that the identity does so trivially or by the inductive hypothesis).
	
	Once again fix $0 \leq t \leq k+1$ and consider the composite $s_t^{k+1} \circ d_{k+2}^{k+1}$. We define the isomorphism
	\[
	\alpha_{s_t^{k+1} \circ d_{k+2}^{k+1}} := \alpha_{d_{k+2}^{k+1}} \circ (d_{k+2}^{k+1})^{\ast}\alpha_{s_t^{k+1}}
	\]
	and note that this, together with the definition above, makes $\alpha_{d_{k+2}^{k+1}}$ satisfy the cocycle condition with respect to the transition isomorphisms $\alpha_{s^{k+1}_{t}}$.
	
	We now construct the transition isomorphisms with respect to the face maps $d_{r}^{k}$ and $d_{i}^{k+2}$, for $0 \leq r \leq k+1$ and $0 \leq i \leq k+2$; note that the construction of the transition isomorphism for $d_{k+3}^{k+2} \circ d_{k+2}^{k+1}$ is handled by the next case of the induction, i.e., the $k+2$ case. We set
	\[
	\alpha_{d_{k+1}^{k+2} \circ d_{i}^{k+2}} := \alpha_{d_{i}^{k+2}} \circ (d_{i}^{k+2})^{\ast}\alpha_{d_{k+1}^{k+2}},
	\]
	which makes $\alpha_{d_{k+2}^{k+1}}$ satisfy the cocycle condition with respect to pre-composition by face maps by construction. Similarly, we set
	\[
	\alpha_{d_{r}^{k} \circ d^{k+1}_{k+2}} := \alpha_{d_{k+2}^{k+1}} \circ (d_{k+2}^{k+1})^{\ast}\alpha_{d_{r}^{k}},
	\]
	which makes $\alpha_{d_{k+2}^{k+1}}$ satisfy the cocycle condition with respect to post-comoposition by face maps. Combining this all together we complete the inductive step of the proof and in turn complete the proof of the proposition.
\end{proof}

While the isomorphisms defined above are all given at the level of degree $0$-components, they all lift through the simplicial techniques by following the strategies used in Lemmas \ref{Lemma: Section 4: Commutativity of face map pullbacks} and \ref{Lemma: Section 4: Simplicial degeneracies are natural with equiv} and Proposition \ref{Prop: Section 4: Simplicial cocycle structure isos}. In particular, ignoring the $n = 0$ case allows us to build simplicial transition isomorphisms between every $L_{G^{n}}(\quot{A}{G^n})$ for all face and degeneracy maps $d_k^n$ and $s_t^n$, provided $n \geq 1$.
\begin{corollary}
	For all $n \in \N$ with $n \geq 1$, there are structure isomorphisms
	\[
	\alpha_{d_k^n}:(\underline{d}_k^n)^{\ast}L_{G^n}(\quot{A}{G^n}) \xrightarrow{\cong} L_{G^{n+1}}(\quot{A}{G^{n+1}}),
	\]
	\[
	\alpha_{s_t^n}:(\underline{s}_t^n)^{\ast}(L_{G^{n+1}}(\quot{A}{G^{n+1}})) \xrightarrow{\cong} L_{G^{n}}(\quot{A}{G^{n}})
	\]
	for all $0 \leq k \leq n+1, 0 \leq t \leq n,$ which satisfy the cocycle condition with respect to morphisms in $\mathbf{\Delta}$.
\end{corollary}

\begin{proposition}\label{Prop: Section 4: The functor L from derived equiv cat to simplicial equiv cat}
	For any smooth affine algebraic group $G$ and any left $G$-variety $X$ there is a functor
	\[
	L:D_G^b(X;\overline{\Q}_{\ell}) \to \Dbeqsimp{X}.
	\]\index[notation]{L@$L:\DbeqQl{X} \to \Dbeqsimp{X}$}
	The functor $L$ is given on objects by sending $(A,T_A)$ to the sequence $(LA_n)$ where
	\[
	LA_n := \begin{cases}
	\psi^{\ast}(\quot{A}{G}) & \text{if}\,\, n = 0; \\
	L_{G^n}(\quot{A}{G^n})_0 & \text{if}\,\, n \geq 1;
	\end{cases}
	\]
	and is given on morphisms by sending a morphism $P = \lbrace \quot{\rho}{\Gamma} \; | \; \Gamma \in \Sf(G)_0\rbrace$ to the sequence $(\rho_n)$ where
	\[
	\rho_n := \begin{cases}
	\psi^{\ast}(\quot{\rho}{G}) & \text{if}\,\,n = 0; \\
	L_{G^n}(\quot{\rho}{G^n})_0 & \text{if}\,\, n \geq 1.
	\end{cases}
	\] 
\end{proposition}
\begin{proof}
	We first check that the assignments do indeed give a functor. First we note that when $n = 0$, $LA_0 = \psi^{\ast}\quot{A}{G} \in D^b_c(X;\overline{\Q}_{\ell})_0$ by the constructions in Lemma \ref{Lemma: Section 4: Object in LG0}. Similarly, by construction each object $L_{G^n}(\quot{A}{G^{n}}) \in \Dbeqsimp{(G^n \times X)}_0$. As such taking the degree $0$ component gives
	\[
	LA_n = L_{G^n}(\quot{A}{G^n})_0 \in D^b_c(G^{n} \times X; \overline{\Q}_{\ell})_0.
	\]
	We now just need to know about the transition isomorphisms to conclude that $LA$ is indeed an object of $\Dbeqsimp{X}$. However, by Lemma \ref{Lemma: Section 4: Object in LG0} and Propositions \ref{Prop: Section 4: Simplicial cocycle structure isos} and \ref{Prop: Section 4: Action map simplicial isos in deg 0} we have that the structure isomorphisms exist and are given by either the degree $0$ component of the isomorphisms $\alpha_{h}:\underline{h}^{\ast}L_{G^m}(\quot{A}{G^m}) \xrightarrow{\cong} L_{G^n}(\quot{A}{G^n})$ for $m,n \geq 1$ or the structure maps of Lemma \ref{Lemma: Section 4: Object in LG0} when $n = 0$. Thus we get that $(LA_n)_{n \in \N}$ indeed is an object in $\Dbeqsimp{X}$.
	
	The check that, for $P \in D^b_G(X;\overline{\Q}_{\ell})$, the sequence $(L\rho_n)$ is a morphism in $\Dbeqsimp{X}$ is more straightforward. For $n = 0$ the check is trivial. For $n \geq 1$ we note that each of the diagrams
	\[
	\xymatrix{
		\overline{f}^{\ast}\AGammap \ar[r]^-{\overline{f}^{\ast}\quot{\rho}{\Gamma^{\prime}}} \ar[d]_{\tau_f^A} & \overline{f}^{\ast}\BGammap \ar[d]^{\tau_f^B} \\
		\AGamma \ar[r]_-{\quot{\rho}{\Gamma}} & \BGamma
	}
	\]
	commutes for all $f \in \Sf(G)_1$, it is routine to show that the $L\rho_n$ commute with all the structure isomorphisms $\alpha_d$, $\alpha_s$, $\alpha_{d \circ s}$, and $\alpha_{s \circ d}$, for $d$ a face map and $s$ a degeneracy, by using Lemmas \ref{Lemma: Section 4: Presence of equivariant object gives degeneracy cocycle isos}, \ref{Lemma: Section 4: Commutativity of face map pullbacks}, \ref{Lemma: Section 4: Simplicial degeneracies are natural with equiv}, \ref{Lemma: Section 4: Object in LG0}, and Propositions \ref{Prop: Section 4: Simplicial cocycle structure isos} and \ref{Prop: Section 4: Action map simplicial isos in deg 0}. Thus we conclude that $L$ is a functor with the declared assignments.
\end{proof}

We now proceed to construct the inverse equivalence of $L:\DbeqQl{X} \to \Dbeqsimp{X}$: The functor $R:\Dbeqsimp{X} \to D_G^b(X;\overline{\Q}_{\ell})$ whose construction will occupy the majority of the remainder of this section. For this we will pass through various functors involving the acyclic descent that is hidden in the simplicial data of $\Dbeqsimp{X}$. The construction we perform relies a quick lemma which sates that there are always acyclic resolutions of arbitrary degree above each variety $G^n$.

\begin{lemma}\label{Lemma: Section 4: The acyclic maps for all G^n}
	Let $G$ be a smooth affine algebraic group. For all $m,n \in \N$ there is an $\Sf(G)$-variety $\Gamma$ and an $m$-acyclic morphism $\Gamma \to G^n$.
\end{lemma}
\begin{proof}
	Let $m \in \N$ be arbitrary. We now prove the lemma by induction on $n$. Set $k = 0$ for the base case and note that $G^0 \cong \Spec K$. Now apply Proposition \ref{Prop: Section 4: Existence of n acyclic maps} with $X = G$ and $m$ to construct an $m$-acyclic variety $\Gamma \to \Spec K$; this $\Gamma$ is an $m$-acyclic resolution of $\Spec K \cong G^0$. Moreover, by the remarks preceeding the statement of Propsoition \ref{Prop: Section 4: Existence of n acyclic maps}, we also have $\Gamma \in \Sf(G)_0$ as well. This proves the base case.
	
	We proceed with the induction. Let $k \in \N$ and assume that there is an $m$-acyclic resolution $\gamma_{k,m}:\Gamma_{k,m} \to G^k$.\index[notation]{Gammakn@$\Gamma_{k,n}$}\index[notation]{gammakn@$\gamma_{k,n}$} Now define the variety
	\[
	\Gamma_{k+1,m} := G \times \Gamma_{k,m}
	\]
	and note that $\Gamma_{k+1,m} \in \Sf(G)_0$. We now define our morphism $\gamma_{k+1,m}:\Gamma_{k+1,m} \to G^{k+1}$ through the diagram below:
	\[
	\xymatrix{
		\Gamma_{k+1,m} \ar[rr]^-{\gamma_{k+1,m}} \ar@{=}[d] & & G^{k+1} \ar@{=}[d] \\
		G \times \Gamma_{k,m} \ar[rr]_-{\id_{G} \times \gamma_{k,m}} & & G \times G^{k}
	}
	\]
	Note first that since $G$ is a smooth principal $G$-variety, the morphism $\id_{G}:G \to G$ is $\infty$-acyclic. Moreover, because $\gamma_{k+1,m} = \id_G \times \gamma_{k,m}$ is a product of maps, and hence a base change of an $m$-acyclic and $\infty$-acyclic map over the terminal object, it follows that it is $m$-acyclic as well. This completes the inductive step.
\end{proof}
An important observation that we made in the proof above is recorderd here for later use.
\begin{lemma}\label{Lemma: Section 4: Infty Acyclic Resls for each Gn}
	For each variety $G^n$, there is an $\infty$-acyclic resolution over $G^n$.
\end{lemma}
\begin{proof}
	Because $G^n$ is a principal $G$ variety for any $n \in \N$ with $n \geq 1$, the map $\id_{G^n}:G^{n} \to G^{n}$ is $\infty$-acyclic.
\end{proof}

We can now proceed to develop the functor $R$; however, before this let us describe our strategy. To build our functor $R$, we will move in the stages suggested by the diagram below
\[
\xymatrix{
	\Dbeqsimp{X} \ar[d]_{\overline{R}} \ar[r]^{R} & D_G^b(X;\overline{\Q}_{\ell}) \\
	\Glue(\Lambda^X) \ar[r]_{\simeq} & \quot{D^b_G(X;\overline{\Q}_{\ell})}{ABL} \ar[u]_{\simeq}
}
\]
where we will define the functor $\overline{R}$ and the set $\Lambda^X$ have yet to be defined. The main idea here is to use the two lemmas above to define a set of acyclic resolutions above each variety $G^n$ through which to descend. In particular, we will show that $\Lambda^X$ is a set of acyclic resolutions which has an $n$-acyclic resolution of $X$ for every $n \in \N$ so that the bottom equivalence passes through Remark \ref{Remark: Equivalence of All Resls with Acylcic Resls} and Proposition \ref{Prop: Section 4: ABL derived cat is SfResl Glue Cat}. It is also worth observing that the rightmost equivalence passes through Corollary \ref{Cor: Section 4: EDCs are equiv}.

We begin our construction of $\Lambda^X$ by defining intermediate sets of resolutions. Fix an $n \in \N$ and define the set
\[
\Lambda_n := \begin{cases}
\lbrace \gamma_{m,n}:\Gamma_{m,n} \to G^n \; | \; m \in \N \rbrace & \text{if}\, n = 0 \\
\lbrace \gamma_{m,n}:\Gamma_{m,n} \to G^n \; | \; m \in \N \rbrace \cup \lbrace \id_{G^n}:G^n \to G^n \rbrace & \text{if}\, n \geq 1
\end{cases}
\]\index[notation]{Lambdan@$\Lambda_n$}
where each $\gamma_{m,n}:\Gamma_{m,n} \to G^n$ is one of the $m$-acyclic resolutions of $G^n$ constructed in Lemma \ref{Lemma: Section 4: The acyclic maps for all G^n}. Define the set
\[
\Lambda := \bigcup_{n \in \N} \Lambda_n
\]\index[notation]{Lambda@$\Lambda$}
and note that because $\Lambda_0 = \lbrace \gamma_{m,0}:\Gamma_{m,0} \to \Spec K \; | \; n \in \N \rbrace$, $\Lambda$ contains an $m$-acyclic resolution of $\Spec K$ for every $m \in \N$ and $\Lambda_0 \subseteq \Lambda$, we have that $\Lambda$ contains an $m$-acyclic resolution of $\Spec K$ for all $m \in \N$. Once again fixing $n \in \N$, define the $X$-relative version of $\Lambda_n$ via
\[
\Lambda^X_n := \begin{cases}
\lbrace \Gamma_{m,n} \times X \xrightarrow{\gamma_{m,n} \times \id_{X}} G^{n} \times X \; | \; m \in \N \rbrace & \text{if}\, n = 0; \\
\lbrace \Gamma_{m,n} \times X \xrightarrow{\gamma_{m,n} \times \id_{X}} G^{n} \times X \; | \; m \in \N \rbrace \cup \lbrace \id_{G^n \times X}\rbrace & \text{if}\, n \geq 1;
\end{cases}
\]\index[notation]{LambdanX@$\Lambda_n^X$}
and note that $\Lambda_0^X$ contains an $m$-acyclic resolution of $X$ for all $m \in \N$. This can be seen because the maps $\gamma_{m,0} \times \id_{X}:\Gamma_{m,0} \times X \to \Spec K \times X$ are isomorphic to the maps $\pi_2^{\Gamma_{m,0}}:\Gamma_{m,0} \times X \to X$, which itself is the base change of $\gamma_{m,0}$ in the diagram
\[
\xymatrix{
	\Gamma_{m,0} \times X \ar[r]^-{\pi_1^{\Gamma_{m,0}}} \ar[d]_{\pi_{2}^{\Gamma_{m,0}}} \pullbackcorner & \Gamma_{m,0} \ar[d]^{\gamma_{m,0}} \\
	X \ar[r] & \Spec K
}
\]
where the pullback is isomorphic to the product because it is taken over the terminal object $\Spec K$. We then define the $X$-relative version of $\Lambda$ to be
\[
\Lambda^X := \bigcup_{n \in \N} \Lambda_n^X
\]\index[notation]{LambdaX@$\Lambda^X$}
and note that $\Lambda^X$ has $m$-acyclic resolutions of $X$ for all $m \in \N$. Because of this we get the following helpful lemma.
\begin{lemma}\label{Lemma: Section 4: Equivalence of GlueLambda with ABL cat}
	There is an equivalence of categories
	\[
	\Glue(\Lambda^X) \simeq \quot{D^b_G(X;\overline{\Q}_{\ell})}{ABL},
	\]\index[notation]{GlueLambdaX@$\Glue(\Lambda^X)$}
\end{lemma}
\begin{proof}
	This is a routine application of Remark \ref{Remark: Equivalence of All Resls with Acylcic Resls} and Proposition \ref{Prop: Section 4: ABL derived cat is SfResl Glue Cat} together with the observations made before the statment of the lemma.
\end{proof}
\begin{remark}
	We will write $\id_{G^n}:G^n \to G^n$ as $\gamma_{\infty,n}:\Gamma_{\infty,n} \to G^n \times X$ in order to have notational consistency when defining constructions along the sets $\Lambda$ and $\Lambda^X$. In this case we want to think of $\id_{G^{n}}$ as an $\infty$-acyclic resolution of $G^n$, so the notation is not as bad as it otherwise could be.
\end{remark}
We now need to construct our functor $\overline{R}:\Dbeqsimp{X} \to \Glue(\Lambda^X)$. For this we recall that the objects of $\Glue(\Lambda^X)$ are triples $(A,A_{{m,n}},\theta_{\Gamma_{m,n}})$ where $n \in \N$, $m \in \N \cup \lbrace \infty \rbrace$, $A \in D_c^b(X;\overline{\Q}_{\ell})_0$, each $A_{{m,n}} \in D_c^b(\quot{X}{\Gamma_{m,n}};\overline{\Q}_{\ell})_0$, and 
\[
\theta_{m,n}:(\pi_2^{\Gamma_{m,n}})^{\ast}A \xrightarrow{\cong} (\quo_{\Gamma_{m,n}})^{\ast}A_{m,n}
\] 
is an isomorphism in $D_c^b(\Gamma_{m,n} \times X;\overline{\Q}_{\ell})$ satisfying the ABL-cocycle condition. Note that the morphisms $\pi_2^{\Gamma_{m,n}}$ and $\quo_{\Gamma_{m,n}}$ come from the diagram:
\[
\xymatrix{
	X & & \Gamma_{m,n} \times X \ar[ll]_-{\pi_2^{\Gamma_{m,n}}} \ar[rr]^-{\quo_{\Gamma_{m,n}}} & & \quot{X}{\Gamma_{m,n}}
}
\] 
We will soon see that as in our construction of the functor $L$ (cf.\@ Proposition \ref{Prop: Section 4: The functor L from derived equiv cat to simplicial equiv cat}), defining things for objects and morphisms in $\Dbeqsimp{X}$ when \\$n \geq 1$ is not too difficult, as we then find ourselves over some nontrivial product $G^n \times X$ and hence can use the machinery involving quotients as we did prior in this chapter; it is the case when we have no power of $G$ to temper the (potential) harshness of $X$ that we have to work to overcome.

Let us begin our construction of $\overline{R}$ by beginning with the more straightforward aspect of defining $\overline{R}$. Fix an $n \geq 1$ and let $m \in \N \cup \lbrace \infty \rbrace$. We now observe that since $\Gamma_{m,n} \in \Sf(G)_0$ by construction, we have a commuting diagram of schemes:
\[
\xymatrix{
	\Gamma_{m,n} \times X \ar[rr]^-{\gamma_{m,n} \times \id_{X}} \ar[d]_{\quo_{\Gamma_{m,n}}} & & G^n \times X \ar[d]^{\quo_{G^n}} \\
	\quot{X}{\Gamma_{m,n}} \ar[rr]_-{\overline{\gamma}_{m,n}} & & \quot{X}{G^n} 
}
\]
Our ultimate goal is to make an assignment as close to something like sending the complex $A_n$ on $G^n \times X$ to the pullback $\overline{\gamma}_{m,n}^{\ast}A_n$ as possible; however, this does not type correctly, as one can see: $A_n$ is a complex of $\ell$-adic sheaves on $G^n \times X$, not on $\quot{X}{G^n}$. Thus to get this working we need to do some trickery that essentially amounts to sending $A_n$ to a simplicial approximation in $\Dbeqsimp{(G^n \times X)}$ and then invoking Proposition \ref{Prop: Section 4: Quot derived cat is equiv to simplicial derived cat} to define something over $\quot{X}{G^n}$.

Let us examine this strategy in more detail. Begin by defining a functor
\[
\underline{(-)}:D^b_c(G^n \times X; \overline{\Q}_{\ell}) \to \Dbeqsimp{(G^n \times X)}
\]\index[notation]{Underline@$\underline{(-)}$}
as follows: We send an object $A \in D_c^b(X;\overline{\Q}_{\ell})$ to the sequence $(A_m) := (\pi_m^{\ast}A)$, where for all $m \in \N$, $\pi_m$ is the projection morphism
\[
\pi_m:G^m \times G^n \times X \to G^n \times X
\]
for all $m \in \N$. Similarly, this is defined on morphisms by sending a morphism $\rho$ to the sequence
$(\rho_m) := (\pi_m^{\ast}\rho)$. Note that this is essentially the result of what happens to an object in a category over the $G$-object $G^n \times X$ when applying the simplicialization functor of Corollary \ref{Cor: Section 4: Functor from equivariant schemes to simplicial schemes}. Perhaps unsurprisingly, it is routine to check that the image of $A$ under this functor is in fact an equivariant object in $D^b_c(\underline{G \backslash (G^n \times X)}_{\bullet};\overline{\Q}_{\ell})$.  We will write $\underline{A}$ for the image of an object $A$ under this functor, and similarly for morphisms. We now obtain a functor from $D_c^b(G^n \times X;\overline{\Q}_{\ell})$ to $D_c^b(\quot{G^n}{X};\overline{\Q}_{\ell})$ by post-composing the functor $\underline{(-)}$ with the inverse equivalence $R_{G^n}$ of $L_{G^n}$ asserted by Proposition \ref{Prop: Section 4: Quot derived cat is equiv to simplicial derived cat}; explicitly we have our functor $\overline{R}_{G^n}:D_c^b(G^n \times X;\overline{\Q}_{\ell}) \to D_c^b(\quot{X}{G^n})$ is given by the diagram:
\[
\xymatrix{
	D_c^b(G^n \times X; \overline{\Q}_{\ell}) \ar[r]^-{\underline{(-)}} & \Dbeqsimp{(G^n \times X)} \ar[r]^-{R_{G^n}}_-{\simeq} & D^b_c(\quot{X}{G^n};\overline{\Q}_{\ell})
}
\]
It is worth noting tha this functor $\overline{R}_{G^n}$ defined above plays our role of a (simplicial) equivariant pushforward.

Now that we have a functor that sends each complex $A$ of $D_c^b(G^n \times X;\overline{\Q}_{\ell})$ to a complex $R_{G^n}(\underline{A})$ in $D_c^b(\quot{X}{G^n};\overline{\Q}_{\ell})$ for $n \geq 1$, we need to construct complexes in $D_c^b(\quot{X}{\Gamma_{m,n}};\overline{\Q}_{\ell})$ for each $m \in \N \cup \lbrace \infty \rbrace$. However, this step of the construction is straightforward; from the commuting diagram
\[
\xymatrix{
	\Gamma_{m,n} \times X \ar[rr]^-{\gamma_{m,n} \times \id_{X}} \ar[d]_{\quo_{\Gamma_{m,n}}} & & G^n \times X \ar[d]^{\quo_{G^n}} \\
	\quot{X}{\Gamma_{m,n}} \ar[rr]_-{\overline{\gamma}_{m,n}} & & \quot{X}{G^n} 
}
\]
we get that the object
\[
\overline{\gamma}_{m,n}^{\ast}(R_{G^n}(\underline{A})) = \overline{\gamma}_{m,n}^{\ast}A_{\infty,n}
\]
does the job. We will show later that there are indeed ABL-cocycle isomorphisms between these objects as necessary, but for the moment we will instead move to consider what happens when $n = 0$.

When $n = 0$, our task of constructing a functor $D_c^b(X;\overline{\Q}_{\ell}) \to D_c^b(\quot{X}{\Gamma_{m,0}};\overline{\Q}_{\ell})$ is more subtle. However, we use Lemma \ref{Lemma: Section 4: Object in LG0} to guide us. Fix an $m \in \N \cup \lbrace \infty \rbrace$ and note that we have the span
\[
\xymatrix{
	\Gamma_{m,0} \times X \ar[r]^-{\pi_{2}^{\Gamma_{m,0}}} \ar[d]_{\quo_{\Gamma_{m,0}}} & X \\
	\quot{X}{\Gamma_{m,0}}
}
\]
from which we must deduce our functor. However, recall that we have an isomorphism $\psi:X \to \quot{X}{G}$ which equips $X$ with the trivial $G$-action. As such, from the universal property of the quotient, we find a unique map which renders the diagrams
\[
\xymatrix{
	\Gamma_{m,0} \times X \ar[r]^-{\pi_{2}^{\Gamma_{m,0}}} \ar[d]_{\quo_{\Gamma_{m,0}}} & X \ar[d]^{\psi} & \Gamma_{m,0} \times X \ar[r]^-{\pi_{2}^{\Gamma_{m,0}}} \ar[d]_{\quo_{\Gamma_{m,0}}} & X \\
	\quot{X}{\Gamma_{m,0}} \ar@{-->}[r]_-{\exists!\overline{\gamma}_{m,n}} & \quot{X}{G} & \quot{X}{\Gamma_{m,0}} \ar@{-->}[r]_-{\exists!\overline{\gamma}_{m,n}} & \quot{X}{G} \ar[u]_{\psi^{-1}}
}
\]
commutative. As such we can define a functor $D_c^b(X;\overline{\Q}_{\ell}) \to D_c^b(\quot{X}{\Gamma_{m,n}};\overline{\Q}_{\ell})$ via the assignments
\[
A \mapsto \overline{\gamma}_{m,0}^{\ast}(\psi^{-1})^{\ast}A
\]
on objects $A$ and
\[
\rho \mapsto \overline{\gamma}_{m,0}^{\ast}(\psi^{-1})^{\ast}\rho
\]
on morphisms $\rho$. 

We will put these assignments together now to build our functor $\overline{R}$ after first proving some technical lemmas that tell us we satisfy the structural conditions necessary to build such a functor. In the meantime, we will give the necessary definitions of the objects that will make up our $\overline{R}$ functor before moving on to prove some basic lemmas and propositions that allow us to deduce the existence of the isomorphisms $\theta_{m,n}$ for $n \geq 1$ and the ABL-cocycle condition on our objects $A_{m,n}$ whenever $n \geq 1$.

\begin{definition}
	For every $n \in \N$ with $n \geq 1$, the functor 
	\[
	\overline{R}_{G^n}:\DbQl{G^n \times X} \to \DbQl{\quot{X}{G^n}}
	\]\index[notation]{RoverlineGn@$\overline{R}_{G^n}$} 
	is given by
	\[
	\overline{R}_{G^n} := R_{G^n} \circ \underline{(-)}.
	\]
\end{definition}
\begin{definition}
	Fix an object $(A_n) \in \Dbeqsimp{X}$ and let $m,n \in \N$. If $n \geq 1$ we first define
	\[
	A_{\infty,n} := \overline{R}_{G^n}(A_n)
	\]\index[notation]{Asubinftyn@$A_{\infty,n}$}
	and then set
	\[
	A_{m,n} := \overline{\gamma}^{\ast}_{m,n}A_{\infty,n} = \overline{\gamma}^{\ast}\overline{R}_{G^n}(A_n).
	\]\index[notation]{Amn@$A_{m,n}$}
	If $n = 0$ we set
	\[
	A_{m,0} := \overline{\gamma}_{m,0}^{\ast}(\psi^{-1})^{\ast}A_0
	\]
	and note that in this case we need not define $A_{\infty,0}$.
\end{definition}
\begin{lemma}\label{Lemma: Section 4: Factorization of RGn}
	Let $n \geq 1$. Then the functor
	\[
	\quo_{G^{n}}^{\ast}:D_c^b(\quot{X}{G^n};\overline{\Q}_{\ell}) \to D_c^b(G^n \times X;\overline{\Q}_{\ell})
	\]
	factors as
	\[
	\xymatrix{
		D_c^b(\quot{X}{G^n};\overline{\Q}_{\ell}) \ar[r]^-{L_{G^n}} \ar[dr]_{\quo_{G^n}^{\ast}} & \Dbeqsimp{(G^n \times X)} \ar[d]^{[0]} \\
		& D_c^b(G^n \times X;\overline{\Q}_{\ell})
	}
	\]
	where $[0]:\Dbeqsimp{(G^n \times X)} \to D_c^b(G^n \times X;\overline{\Q}_{\ell})$ is the degree $0$-truncation, i.e., $(A_n)_{n \in \N} \mapsto A_0$.
\end{lemma}
\begin{proof}
	We begin by recalling that for a complex $A \in D_c^b(\quot{X}{G^n};\overline{\Q}_{\ell})$,
	\[
	L_{G^n}(A)[0] = \big(((\epsilon_{n}^{m})^{\ast}A)_{m \in \N}\big)[0] = (\epsilon_n^0)^{\ast}A = \quo_{G^n}^{\ast}A.
	\]
	The calculation is similar for morphisms, which proves the lemma.
\end{proof}
\begin{lemma}\label{Lemma: Section 4: Quot pullback for An}
	For any $n \geq 1$,
	\[
	\quo_{G^n}^{\ast} \circ \overline{R}_{G^n} \cong \id_{D_c^b(G^n \times X;\overline{\Q}_{\ell})}.
	\]
\end{lemma}
\begin{proof}
	We begin by noting that
	\[
	[0] \circ \underline{(-)} = \id_{D_c^b(G^n \times X;\overline{\Q}_{\ell})}
	\]
	and because $L_{G^n}$ and $R_{G^n}$ are inverse equivalences, 
	\[
	L_{G^n} \circ R_{G^n} \cong \id_{\Dbeqsimp{(G^n \times X)}}.
	\] 
	We thus calculate using Lemma \ref{Lemma: Section 4: Factorization of RGn} that
	\begin{align*}
	\quo_{G^n}^{\ast} \circ \overline{R}_{G^n} &= [0] \circ L_{G^n} \circ R_{G^n} \circ \underline{(-)} \cong [0] \circ \id_{\Dbeqsimp{(G^n \times X)}} \circ \underline{(-)}\\ &= [0] \circ \underline{(-)} 
	= \id_{D_c^b(G^n \times X;\overline{\Q}_{\ell})},
	\end{align*}
	as was to be shown.
\end{proof}

\begin{lemma}\label{Lemma: Section 4.2: Isomorphism at Gn}
	Let $(A_n) \in \Dbeqsimp{X}$. Then there is an isomorphism, for any $n \geq 1$ and for $A_{\infty,n} = \overline{\gamma}_{\infty,n}^{\ast}\overline{R}_{G^n}(A_n)$,
	\[
	\theta_{\infty,n}:\pi_2^{\ast} A_0 \xrightarrow{\cong} \quo_{G^{n}}^{\ast}A_{\infty,n}.
	\]
\end{lemma}
\begin{proof}
	Begin by noting that $\pi_2:G^n \times X \to X$ is a simplicial morphism; in particular,
	\[
	\pi_2 = d_0^{n-1} \circ \cdots \circ d_0^0 = d,
	\]
	so there is an isomorphism in $D_c^b(G \times X; \overline{\Q}_{\ell})$
	\[
	\alpha_d:d^{\ast}A_0 \xrightarrow{\cong} A_n
	\]
	satisfying the cocycle condition. We will now be done if we can prove that $A_n \cong \quo_{G^n}^{\ast}A_{\infty,n}$. For this we note, however, that
	\[
	\quo_{G^n}^{\ast}A_{\infty,n} = \quo_{G^n}^{\ast}\overline{R}_{G^n}(A_n) = \quo_{G^n}^{\ast}R_{G^n}(\underline{A}_n).
	\]
	Invoking Lemma \ref{Lemma: Section 4: Quot pullback for An} we find a natural isomorphism $\quo_{G^n}^{\ast}A_{\infty,n} \cong A_n$; as such, it follows that there is an isomorphism
	\[
	\xymatrix{
		\pi_n^{\ast}A_0 \ar[r]^{\alpha_d} \ar[dr]_{\theta_{\infty,n}} & A_n \ar[d]^{\cong} \\
		& \quo_{G^n}^{\ast}A_{\infty,n}
	}
	\]
	as claimed.
\end{proof}
\begin{lemma}\label{Lemma: Section 4.2: Structure Isos for n > 0}
	Let $n \geq 1$ and let $m \in \N \cup \lbrace \infty \rbrace$. Then there is an isomorphism
	\[
	\theta_{n,m}:(\pi_2^{\Gamma_{n,m}})^{\ast}A_0 \xrightarrow{\cong} \quo_{\Gamma_{m,n}}^{\ast}A_{m,n}
	\]\index[notation]{Thetamn@$\theta_{m,n}$}
	in $D_c^b(\Gamma_{m,n} \times X;\overline{\Q}_{\ell})$.
\end{lemma}
\begin{proof}
	The case where $m = \infty$ is given in Lemma \ref{Lemma: Section 4.2: Isomorphism at Gn}, so we can safely take $m \in \N$. Now observe that the diagram
	\[
	\begin{tikzcd}
	\Gamma_{n,m} \times X \ar[rrrr, bend left = 30]{}{\pi_{2}^{m,n}} \ar[rr]{}{\gamma_{m,n} \times \id_X} \ar[d, swap]{}{\quo_{\Gamma_{m,n}}} & & G^n \times X \ar[rr]{}{\pi_2^n} \ar[d]{}{\quo_{G^n}} & & X \\
	\quot{X}{\Gamma_{m,n}} \ar[rr, swap]{}{\overline{\gamma}_{m,n}} & & \quot{X}{G^n}
	\end{tikzcd}
	\]
	commutes in $\Sch$. With this we calculate that
	\begin{align*}
	\quo_{\Gamma_{m,n}}^{\ast}A_{m,n} &= \quo_{\Gamma_{m,n}}^{\ast}\overline{\gamma}_{m,n}^{\ast}A_{\infty,n} = (\overline{\gamma}_{m,n} \circ \quo_{\Gamma_{m,n}})^{\ast}A_{\infty,n} \\
	&= (\quo_{G^n} \circ (\gamma_{m,n} \times \id_{G^n}))^{\ast}A_{\infty,n} \\
&	= (\gamma_{m,n} \times \id_X)^{\ast}\quo_{G^n}^{\ast}A_{\infty,n}
	\end{align*}
	and similarly that
	\[
	(\pi_2^{m,n})^{\ast}A_0 = (\pi_2^n \circ (\gamma_{m,n} \circ \id_X))^{\ast}A_0 = (\gamma_{m,n} \times \id_x)^{\ast}(\pi_2^n)^{\ast}A_0.
	\]
	Now use Lemma \ref{Lemma: Section 4.2: Isomorphism at Gn} to give an isomorphism $\theta_{\infty,n}:(\pi_2^n)^{\ast}A_0 \xrightarrow{\cong} \quo_{G^{n}}^{\ast}A_{\infty,n}$. Applying the functor $(\gamma_{m,n} \times \id_X)^{\ast}$ to the isomorphism $\theta_{\infty,n}$ gives an isomorphism
	\[
	(\gamma_{m,n} \times \id_{X})^{\ast}\theta_{\infty,n}:(\gamma_{m,n} \times \id_X)^{\ast}(\pi_2^n)^{\ast}A_0 \xrightarrow{\cong} (\gamma_{m,n} \times \id_X)^{\ast}\quo_{G^n}^{\ast}A_{\infty,n}
	\]
	which is equal to the desired isomorphism
	\[
	(\gamma_{m,n} \times \id_{X})^{\ast}\theta_{\infty,n}: (\pi_2^{m,n})^{\ast}A_0 \xrightarrow{\cong} \quo_{\Gamma_{m,n}}^{\ast}A_{m,n}.
	\]
	Setting
	\[
	\theta_{m,n}:= (\gamma_{m,n} \times \id_{X})^{\ast}\theta_{\infty,n}
	\]
	proves the lemma.
\end{proof}

Now that we have shown how the $\theta_{m,n}$ arise whenever $n \geq 1$, we need to show that for $n \geq 1$ the $A_{m,n}$ satisfy the ABL-cocycle condition as well.

\begin{lemma}\label{Lemma: Section 4.2: Commutativity of d pullback with RGn}
	Fix any $m,n \in \N$ with $m,n \geq 1$. Then if $n \geq m$ and $h:G^n \times X \to G^m \times X$ is a simplicial map (a composition of face and degeneracy maps), the diagram
	\[
	\xymatrix{
		D_c^b(G^m \times X) \ar[r]^-{h^{\ast}} \ar[d]_{\overline{R}_{G^m}} & D_c^b(G^n \times X) \ar[d]^{\overline{R}_{G^n}} \\
		D_{c}^{b}(\quot{X}{G^m}) \ar[r]_-{\overline{h}^{\ast}} & \DbQl{\quot{X}{G^n}}
	}
	\]
\end{lemma}
\begin{proof}
	Begin by observing that the commutativity of the diagrams
	\[
	\xymatrix{
		\DbQl{G^m \times X} \ar[r]^-{h^{\ast}} \ar[d]_{\underline{(-)}} & \DbQl{G^n \times X} \ar[d]^{\underline{(-)}} \\
		\Dbeqsimp{(G^m \times X)} \ar[r]_-{\underline{h}^{\ast}} & \Dbeqsimp{(G^n \times X)}
	}
	\]
	is both routine and trivial to verify, and hence is omitted. However,  the diagram
	\[
	\xymatrix{
		\DbQl{G^m \times X} \ar[r]^-{h^{\ast}} \ar[d]_{\underline{(-)}} & \DbQl{G^n \times X} \ar[d]^{\underline{(-)}} \\
		\Dbeqsimp{(G^m \times X)} \ar[r]_-{\underline{h}^{\ast}} \ar[d]_{R_{G^m}} & \Dbeqsimp{(G^n \times X)} \ar[d]^{R_{G^n}} \\
		\DbQl{\quot{X}{G^m}} \ar[r]_-{\overline{h}^{\ast}} & \DbQl{\quot{X}{G^n}}
	}
	\]
	factorizes the diagram
	\[
	\xymatrix{
		\DbQl{G^m \times X} \ar[r]^-{h^{\ast}} \ar[d]_{\overline{R}_{G^m}} & \DbQl{G^n \times X} \ar[d]^{\overline{R}_{G^n}} \\
		\DbQl{\quot{X}{G^m}} \ar[r]_-{\overline{h}^{\ast}} & \DbQl{\quot{X}{G^n}}
	}
	\]
	which concludes the proof.
\end{proof}
\begin{proposition}\label{Prop: Section 4.2: Cocycle isos for n > 0 and m = infty}
	Let $k,n \in \N$ with $k,n \geq 1$ and $n \geq k$ and consider the diagram of schemes:
	\[
	\xymatrix{
		& G^{n} \times G^{k} \times X \ar[dr]^{\pi_{23}} \ar[dl]_{\pi_{13}} & \\
		G^n \times X \ar[dr]_{\pi_2^n} & & G^k \ar[dl]^{\pi_2^k} \\
		& X
	}
	\]
	For any $(A_n) \in \Dbeqsimp{X}_0$, there are then isomorphisms
	\[
	\varphi_{nk}:\opi_{13}^{\ast}A_{\infty,n} \xrightarrow{\cong} \opi_{23}^{\ast}A_{\infty,k}
	\]\index[notation]{phivarnk@$\phi_{nk}$}
	which satisfy the ABL-cocycle condition.
\end{proposition}
\begin{proof}
	We first construct our desired isomorphisms. For this begin by observing that since $n \geq s$, there is a simplicial map $G^n \times X \to G^s \times X$ making the diagram
	\[
	\xymatrix{
		& G^{n} \times G^{s} \times X \ar[dr]^{\pi_{s}^{ns}} \ar[dl]_{\pi_{n}^{ns}} & \\
		G^n \times X \ar[rr]^-{d} \ar[dr]_{\pi_2^n} & & G^s \ar[dl]^{\pi_2^k} \\
		& X
	}
	\]
	commute; $d = \id_{G^n \times X}$ if $n = s$ and otherwise $d$ is a composite of $d^{k}_0$ morphisms that project away coordinates. Because the above diagram commutes, so does the diagram:
	\[
	\xymatrix{
		& \quot{X}{G^{n+s}} \ar[dr]^{\opi_{s}^{ns}} \ar[dl]_{\opi_{n}^{ns}} & \\
		\quot{X}{G^n} \ar[rr]_-{\overline{d}} & & \quot{X}{G^s}
	}
	\]
	Now recall that since $(A_n) \in \Dbeqsimp{X}$, there is an isomorphism $\alpha_{d}:d^{\ast}A_s \xrightarrow{\cong} A_n$ in $\DbQl{G^n \times X}.$ Applying the functor $\overline{R}_{G^n}$ to this isomorphism, we obtain an isomorphism
	\[
	\overline{R}_{G^n}(\alpha_d):\overline{R}_{G^n}(d^{\ast}A_s) \xrightarrow{\cong} \overline{R}_{G^n}(A_n)
	\]
	which, by Lemma \ref{Lemma: Section 4.2: Commutativity of d pullback with RGn}, is equal to an isomorphism
	\[
	\overline{d}^{\ast}\overline{R}_{G^{s}}(A_s) \xrightarrow{\cong} \overline{R}_{G^n}(A_n).
	\]
	To get the desired isomorphism, we apply the functor $(\opi_{n}^{ns})^{\ast}$ to the above isomorphism, together with the observation that
	\[
	(\opi_{n}^{ns})^{\ast}\overline{d}^{\ast}R_{G^s}(A_s) = (\overline{d} \circ \opi_n^{ns})^{\ast}\overline{R}_{G^s}(A_s) = (\opi_{s}^{ns})^{\ast}\overline{R}_{G^s}(A_s).
	\]
	This gives an isomorphism
	\[
	(\opi_{n}^{ns})^{\ast}\overline{R}_{G^n}(\alpha_d):(\opi_s^{ns})^{\ast}\overline{R}_{G^s}(A_s) \xrightarrow{\cong} \overline{R}_{G^n}(A_n).
	\]
	Defining
	\[
	\varphi_{ns} := \left((\opi_{n}^{ns})^{\ast}\overline{R}_{G^n}(\alpha_d)\right)^{-1}
	\]
	and recalling that $\overline{R}_{G^n}(A_n) = A_{\infty,n}$ and $\overline{R}_{G^s}(A_s) = A_{\infty,s}$ then provides us with our desired isomorphism $\varphi_{ns}:(\opi_{n}^{ns})^{\ast}A_{\infty,n} \xrightarrow{\cong} (\opi_{s}^{ns})^{\ast}A_{\infty,s}$.
	
	We now prove the ABL-cocycle condition holds. For this fix $n,s,j \in N$ with $n,s,j \geq 1$ and assume that $n \geq s \geq j$. As above, there are simplicial maps $d_{ns}:G^n \times X \to G^s \times X, d_{sj}:G^s \times X \to G^j \times X$, and $d_{nj}:G^n \times X \to G^j \times X$ which factor the projection maps. Furthermore, these maps satisfy $d_{nj} = d_{sj} \circ d_{ns}$ and make the diagram
	\[
	\begin{tikzcd}
	& G^{n+s+j} \times X \ar[ddr]{}{\pi_{nj}^{nsj}} \ar[dd]{}{\pi_{sj}^{nsj}} \ar[ddl, swap]{}{\pi_{ns}^{nsj}} \\ 
	\\
	G^{n + s} \times X \ar[dd, swap]{}{\pi_{n}^{ns}} & G^{n + j} \ar[ddl, near start]{}{\pi_{n}^{nj}} \ar[ddr, swap, near start]{}{\pi_{j}^{nj}} \times X & G^{s + j} \times X \ar[dd]{}{\pi_j^{sj}} \\
	\\
	G^n \times X \ar[r]{}{d_{ns}} \ar[rr, bend right = 30, swap]{}{d_{nj}} & G^s \times X \ar[r]{}{d_{sj}} &  G^j \times X \\
	\ar[from = 3-1, to = 5-2, crossing over, near start]{}{\pi_s^{ns}} \ar[from = 3-3, to = 5-2, crossing over, swap, near start]{}{\pi^{sj}_{j}} \\
	\end{tikzcd}
	\]
	commute. We now must show that
	\[
	(\opi_{nj}^{nsj})^{\ast}\varphi_{nj} = (\opi_{sj}^{nsj})^{\ast}\varphi_{sj} \circ (\opi_{ns}^{nsj})^{\ast}\varphi_{ns}.
	\]
	However, because each isomorphism $\varphi_{nj}, \varphi_{ns}, \varphi_{sj}$ are defined as the inverse of some other morphism, the above cocycle condition holds if and only if
	\[
	(\opi_{nj}^{nsj})^{\ast}\varphi_{nj}^{-1} = (\opi_{ns}^{nsj})^{\ast}\varphi_{ns}^{-1} \circ (\opi_{sj}^{nsj})^{\ast}\varphi_{sj}^{-1},
	\]
	which is what we will prove instead. To spell things out more explicitly, we will show that the equality
	\[
	(\opi_{nj}^{nsj})^{\ast}(\opi_{n}^{nj})^{\ast}\overline{R}_{G^n}(\alpha_{d_{nj}}) = (\opi_{ns}^{nsj})^{\ast}(\opi_{n}^{ns})^{\ast}\overline{R}_{G^n}(\alpha_{d_{nj}}) \circ (\opi_{sj}^{nsj})^{\ast}(\opi_{j}^{js})^{\ast}\overline{R}_{G^s}(\alpha_{d_{js}})
	\]
	holds. For this we recall that since $d_{nj} = d_{sj} \circ d_{ns}$ we have that
	\[
	\alpha_{d_{nj}} = \alpha_{d_{sj} \circ d_{ns}} = \alpha_{d_{ns}} \circ d_{ns}^{\ast}\alpha_{d_{sj}}.
	\]
	This allows us to deduce that
	\begin{align*}
	&(\opi_{nj}^{nsj})^{\ast}(\opi_{n}^{nj})^{\ast}\overline{R}_{G^n}(\alpha_{d_{nj}}) \\
	&= (\opi_{nj}^{nsj})^{\ast}(\opi_{n}^{nj})^{\ast}\left(\overline{R}_{G^n}(\alpha_{d_{ns}}) \circ \overline{R}_{G^n}(d_{ns}^{\ast}\alpha_{d_{sj}})\right) \\
	&= \left(\opi_{n}^{nj} \circ \opi_{nj}^{nsj}\right)^{\ast}\left(\overline{R}_{G^n}(\alpha_{d_{ns}}) \circ \overline{R}_{G^n}(d_{ns}^{\ast}\alpha_{d_{sj}})\right) \\
	&= \left(\opi_{n}^{nj} \circ \opi_{nj}^{nsj}\right)^{\ast}\overline{R}_{G^n}(\alpha_{d_{ns}}) \circ \left(\opi_{n}^{nj} \circ \opi_{nj}^{nsj}\right)^{\ast}\overline{R}_{G^n}(d_{ns}^{\ast}\alpha_{d_{sj}}).
	\end{align*}
	We now consider that since $\pi_{n}^{nj} \circ \pi_{nj}^{nsj} = \pi_{n}^{ns} \circ \pi_{ns}^{nsj}$, we have that $(\opi_{n}^{nj} \circ \opi_{nj}^{nsj})^{\ast} = (\opi_{n}^{ns} \circ \opi_{ns}^{nsj})^{\ast} = (\opi_{ns}^{nsj})^{\ast}(\opi_{n}^{ns})^{\ast}$. This allows us to write
	\[
	\left(\opi_{n}^{nj} \circ \opi_{nj}^{nsj}\right)^{\ast}\overline{R}_{G^n}(\alpha_{d_{ns}}) = (\opi_{ns}^{nsj})^{\ast}(\opi_{n}^{ns})^{\ast}\overline{R}_{G^n}(\alpha_{d_{ns}}) = (\opi_{ns}^{nsj})\varphi_{ns}^{-1}.
	\]
	Now recall that from Lemma \ref{Lemma: Section 4.2: Commutativity of d pullback with RGn} we have that $\overline{R}_{G^n}(d_{ns}^{\ast}\alpha_{d_{sj}}) = \overline{d}_{ns}^{\ast}\overline{R}_{G^s}(\alpha_{d_{sj}})$. We also observe that
	\begin{align*}
	\left(\opi_{n}^{nj} \circ \opi_{nj}^{nsj}\right)^{\ast}\overline{d}_{ns}^{\ast} &= \left(\overline{d}_{ns} \circ \opi_{n}^{nj} \circ \opi_{nj}^{nsj} \right)^{\ast} = \left(\overline{d}_{ns} \circ \opi_{n}^{ns} \circ \opi_{ns}^{nsj}\right)^{\ast} = \left(\opi_{s}^{ns} \circ \opi_{ns}^{nsj}\right)^{\ast}\\ 
	&= \left(\opi_{s}^{sj} \circ \opi_{sj}^{nsj}\right)^{\ast} 
	= (\opi_{sj}^{nsj})^{\ast}(\opi_{sj}^{nsj})^{\ast}.
	\end{align*}
	Combining these observations gives that 
	\begin{align*}
	\left(\opi_{n}^{nj} \circ \opi_{nj}^{nsj}\right)^{\ast}\overline{R}_{G^n}(d_{ns}^{\ast}\alpha_{d_{sj}}) &= \left(\opi_{n}^{nj} \circ \opi_{nj}^{nsj}\right)^{\ast}\overline{d}_{ns}^{\ast}\overline{R}_{G^s}(\alpha_{d_{sj}})\\ 
	&= (\opi_{sj}^{nsj})^{\ast}(\opi_{sj}^{nsj})^{\ast}\overline{R}_{G^s}(\alpha_{d_{sj}}) 
	= (\opi_{sj}^{nsj})^{\ast}\varphi_{sj}^{-1}.
	\end{align*}
	Substituting these observations then gives
	\begin{align*}
	&(\opi_{nj}^{nsj})^{\ast}\varphi_{nj}^{-1}  =(\opi_{nj}^{nsj})^{\ast}(\opi_{n}^{nj})^{\ast}\overline{R}_{G^n}(\alpha_{d_{nj}}) \\
	&= \left(\opi_{n}^{nj} \circ \opi_{nj}^{nsj}\right)^{\ast}\overline{R}_{G^n}(\alpha_{d_{ns}}) \circ \left(\opi_{n}^{nj} \circ \opi_{nj}^{nsj}\right)^{\ast}\overline{R}_{G^n}(d_{ns}^{\ast}\alpha_{d_{sj}}) \\
	&=(\opi_{ns}^{nsj})^{\ast}(\opi_{n}^{ns})^{\ast}\overline{R}_{G^n}(\alpha_{d_{ns}}) \circ (\opi_{sj}^{nsj})^{\ast}(\opi_{sj}^{nsj})^{\ast}\overline{R}_{G^s}(\alpha_{d_{sj}}) \\
	&=(\opi_{ns}^{nsj})\varphi_{ns}^{-1} \circ (\opi_{sj}^{nsj})^{\ast}\varphi_{sj}^{-1}.
	\end{align*}
	This concludes the proof of the propostion, as taking inverses proves
	\[
	(\opi_{nj}^{nsj})^{\ast}\varphi_{nj} = (\opi_{sj}^{nsj})^{\ast}\varphi_{sj} \circ (\opi_{ns}^{nsj})\varphi_{ns},
	\]
	as was desired.
\end{proof}
Let us set the stage for the proposition below, as it will have significant (extra) notational bloat that makes proving the ABL-cocycle condition on our isomorphisms more confusing. However, it builds upon Proposition \ref{Prop: Section 4.2: Cocycle isos for n > 0 and m = infty} above and  uses largely the same technique. Fix some $n,t \in \N$ with $n,t \geq 1$ and let $m,s \in \N \cup \lbrace \infty \rbrace$. We then have $\Lambda^X$ morphisms $\gamma_{m,n} \times \id_{X}: \Gamma_{m,n} \times X \to G^n \times X$ and these resolutions factor the projections:
\[
\xymatrix{
	\Gamma_{m,n} \times X \ar[drr]_{\pi_2^{\Gamma_{m,n}}} \ar[rr]^-{\gamma_{m,n} \times \id_X} & & G^n \times X \ar[d]^{\pi_2^{G^n}} \\
	& & X
}
\]
When we are in the presence of two resolutions 
\[
\Gamma_{m,n} \times X \xrightarrow{\pi_2^{\Gamma_{m,n}}} X \xleftarrow{\pi_2^{\Gamma_{s,t}}} \Gamma_{s,t} \times X
\]
we then can proceed by performing manipulations analogous to those used in Proposition \ref{Prop: Section 4.2: Cocycle isos for n > 0 and m = infty} to produce our desired isomorphisms and prove that the satisfy the cocycle condition. Explicitly, if we assume that without loss of generality $n \geq t$ we get a commuting diagram
\[
\begin{tikzcd}
\Gamma_{m,n} \times X \ar[ddr, swap]{}{\gamma_{m,n} \times \id_X} & & \Gamma_{m,n}\times \Gamma_{s,t} \times X \ar[rr]{}{\pi_{s,t}^{\Gamma_{m,n}\Gamma_{s,t}}} \ar[ll, swap]{}{\pi_{m,n}^{\Gamma_{n,m}\Gamma_{s,t}}} \ar[d]{}{\gamma_{m,n} \times \gamma_{s,t} \times \id_X} & & \Gamma_{s,t} \times X \ar[ddl]{}{\gamma_{s,t} \times \id_X} \\
& & G^n \times G^t \times X \ar[dr]{}{\pi_t^{G^nG^t}} \ar[dl, swap]{}{\pi_{n}^{G^nG^t}} \\
& G^n \times X \ar[dr, swap]{}{\pi_2^{G^n}} \ar[rr]{}{d} & & G^t \times X \ar[dl]{}{\pi_2^{G^t}} \\
& & X
\end{tikzcd}
\]
factorizing the resolutions and projections out of the product $\Gamma_{m,n} \times \Gamma_{s,t} \times X$ together with its projections to $\Gamma_{m,n} \times X$ and $\Gamma_{s,t} \times X$. Note that $d$ is a simplicial map given by either the identity or successive projections depending on if $n = t$ or $n > t$, respectively. We will see that the magic happens upon taking quotients, as this gives rise to the further commuting diagram:
\[
\begin{tikzcd}
\quot{X}{\Gamma_{m,n}} \ar[ddr, swap]{}{\ogamma_{m,n}} & & \quot{X}{\Gamma_{m,n}\times\Gamma_{s,t}} \ar[rr]{}{\opi_{s,t}^{\Gamma_{m,n}\Gamma_{s,t}}} \ar[ll, swap]{}{\opi_{m,n}^{\Gamma_{n,m}\Gamma_{s,t}}} \ar[d]{}{\overline{\gamma_{m,n} \times \gamma_{s,t}}} & & \quot{X}{\Gamma_{s,t}} \ar[ddl]{}{\ogamma_{s,t}} \\
& & \quot{X}{G^n \times G^t} \ar[dr]{}{\opi_t^{G^nG^t}} \ar[dl, swap]{}{\opi_{n}^{G^nG^t}} \\
& \quot{X}{G^n}  \ar[rr, swap]{}{\overline{d}} & & \quot{X}{G^t}  \\
\end{tikzcd}
\]
These two diagrams, together with the notational conventions they set, will be of crucial use in the proposition below.
\begin{proposition}\label{Prop: Section 4.2: Cocycle maps for all n,t geq 1 and all m,s}
	Let $n,t \in \N$ with $n,t \geq 1$ and let $m,s \in \N \cup \lbrace \infty \rbrace$. Then there are isomorphisms
	\[
	\varphi^{m,n}_{s,t}:(\overline{\pi}_{m,n}^{\Gamma_{m,n}\Gamma_{s,t}})^{\ast}A_{m,n} \xrightarrow{\cong} (\opi_{s,t}^{\Gamma_{m,n}\Gamma_{s,t}})^{\ast}A_{s,t}
	\]\index[notation]{phivarmnst@$\varphi_{s,t}^{m,n}$}
	which satisfy the ABL-cocycle condition.
\end{proposition}
\begin{proof}
	Recall that Proposition \ref{Prop: Section 4.2: Cocycle isos for n > 0 and m = infty} gives us an isomorphism
	\[
	\varphi_{nt} = \varphi^{\infty,n}_{\infty,t}:(\opi_{n}^{G^nG^t})^{\ast}A_{\infty,n} \xrightarrow{\cong} (\opi_{t}^{G^nG^t})^{\ast}A_{\infty,t}.
	\]
	We now claim that $(\overline{\gamma_{m,n} \times \gamma_{s,t}})^{\ast}\varphi_{n,t}$ gives the desired isomorphism from 
	\[
	(\opi_{m,n}^{\Gamma_{m,n}\Gamma_{s,t}})^{\ast}A_{m,n} \xrightarrow{\cong} (\opi_{s,t}^{\Gamma_{m,n}\Gamma_{s,t}})^{\ast}A_{s,t}.
	\]
	To see this we first calculate that
	\begin{align*}
	(\overline{\gamma_{m,n} \times \gamma_{s,t}})^{\ast}(\opi_{n}^{G^nG^t})^{\ast}A_{\infty,n} &= \left(\opi_{n}^{G^nG^t} \circ \overline{\gamma_{m,n} \times \gamma_{s,t}}\right)^{\ast}A_{\infty,n} \\
	&= \left(\ogamma_{m,n} \circ \opi_{m,n}^{\Gamma_{m,n}\Gamma_{s,t}}\right)^{\ast}A_{\infty,n} \\
	&= (\opi_{m,n}^{\Gamma_{m,n}\Gamma_{s,t}})^{\ast}\ogamma_{m,n}^{\ast}A_{\infty,n} = (\opi_{m,n}^{\Gamma_{m,n}\Gamma_{s,t}})^{\ast}A_{m,n}.
	\end{align*}
	Similarly, we find
	\begin{align*}
	(\overline{\gamma_{m,n} \times \gamma_{s,t}})^{\ast}(\opi_{t}^{G^nG^t})^{\ast}A_{\infty,t} &= \left(\opi_{t}^{G^nG^t} \circ \overline{\gamma_{m,n} \times \gamma_{s,t}}\right)^{\ast}A_{\infty,t} \\
	&= \left(\ogamma_{s,t} \circ \opi_{s,t}^{\Gamma_{m,n}\Gamma_{s,t}}\right)^{\ast}A_{\infty,t} \\
	&= (\opi_{s,t}^{\Gamma_{m,n}\Gamma_{s,t}})^{\ast}\ogamma_{s,t}^{\ast}A_{\infty,t} = (\opi_{s,t}^{\Gamma_{m,n}\Gamma_{s,t}})^{\ast}A_{s,t}.
	\end{align*}
	Thus we define our isomorphisms by
	\[
	\varphi^{m,n}_{s,t} := (\overline{\gamma_{m,n} \times \gamma_{s,t}})^{\ast}\varphi_{nt}.
	\]
	
	We now prove the ABL-cocycle condition on these isomorphisms. For this we let $\Gamma_{m,n}, \Gamma_{s,t}, \Gamma_{i,j} \in \Lambda^X$ with $n \geq t \geq j \geq 1$. We now must show that
	\[
	\left(\opi_{nj}^{\Gamma_{m,n}\Gamma_{s,t}\Gamma_{i,j}}\right)^{\ast}\varphi^{mn}_{ij} = \left(\opi_{tj}^{\Gamma_{m,n}\Gamma_{s,t}\Gamma_{i,j}}\right)^{\ast}\varphi_{ij}^{st} \circ \left(\opi_{nt}^{\Gamma_{m,n}\Gamma_{s,t}\Gamma_{i,j}}\right)^{\ast}\varphi_{st}^{mn}.
	\]
	For this we consider the following three commutative diagrams:
	\[
	\xymatrix{
		\Gamma_{m,n}\times \Gamma_{s,t}\times \Gamma_{i,j} \times X \ar[rr]^-{\pi_{nt}^{\Gamma_{m,n}\Gamma_{s,t}\Gamma_{i,j}}} \ar[d]_{\gamma_{m,n} \times \gamma_{s,t} \times \gamma_{i,j} \times \id_X} & & \Gamma_{m,n}\times \Gamma_{s,t} \times X \ar[d]^{\gamma_{m,n}\times \gamma_{s,t} \times \id_X} \\
		G^n \times G^t \times G^j \times X \ar[rr]_-{\pi^{G^nG^tG^j}_{nt}} & & G^n \times G^t \times X
	}
	\]
	\[
	\xymatrix{
		\Gamma_{m,n}\times \Gamma_{s,t}\times \Gamma_{i,j} \times X \ar[rr]^-{\pi_{tj}^{\Gamma_{m,n}\Gamma_{s,t}\Gamma_{i,j}}} \ar[d]_{\gamma_{m,n} \times \gamma_{s,t} \times \gamma_{i,j} \times \id_X} & & \Gamma_{s,t}\times \Gamma_{i,j} \times X \ar[d]^{\gamma_{s,t}\times \gamma_{i,j} \times \id_X} \\
		G^n \times G^t \times G^j \times X \ar[rr]_-{\pi^{G^nG^tG^j}_{tj}} & & G^t \times G^j \times X
	}
	\]
	\[
	\xymatrix{
		\Gamma_{m,n}\times \Gamma_{s,t}\times \Gamma_{i,j} \times X \ar[rr]^-{\pi_{nj}^{\Gamma_{m,n}\Gamma_{s,t}\Gamma_{i,j}}} \ar[d]_{\gamma_{m,n} \times \gamma_{s,t} \times \gamma_{i,j} \times \id_X} & & \Gamma_{m,n}\times \Gamma_{i,j} \times X \ar[d]^{\gamma_{m,n}\times \gamma_{i,j} \times \id_X} \\
		G^n \times G^t \times G^j \times X \ar[rr]_-{\pi^{G^nG^tG^j}_{nj}} & & G^n \times G^j \times X
	}
	\]
	We use these induced relations to compute that
	\begin{align*}
	&\left(\opi_{tj}^{\Gamma_{m,n}\Gamma_{s,t}\Gamma_{i,j}}\right)^{\ast}\varphi_{ij}^{st} \circ \left(\opi_{st}^{\Gamma_{m,n}\Gamma_{s,t}\Gamma_{i,j}}\right)^{\ast}\varphi_{st}^{mn}  \\
	&= \left(\opi_{tj}^{\Gamma_{m,n}\Gamma_{s,t}\Gamma_{i,j}}\right)^{\ast}(\overline{\gamma_{s,t} \times \gamma_{i,j}})^{\ast}\varphi_{tj} \circ \left(\opi_{st}^{\Gamma_{m,n}\Gamma_{s,t}\Gamma_{i,j}}\right)^{\ast}(\overline{\gamma_{m,n} \times \gamma_{s,t}})^{\ast}\varphi_{nt} \\
	&= \left(\overline{\gamma_{s,t} \times \gamma_{i,j}} \circ \opi_{tj}^{\Gamma_{m,n}\Gamma_{s,t}\Gamma_{i,j}}\right)^{\ast}\varphi_{tj} \circ \left(\overline{\gamma_{m,n}\times \gamma_{s,t}} \circ \opi_{nt}^{\Gamma_{m,n}\Gamma_{s,t}\Gamma_{i,j}}\right)^{\ast}\varphi_{nt} \\
	&= \left(\opi_{tj}^{G^nG^tG^j} \circ \overline{\gamma_{m,n} \times \gamma_{s,t} \times \gamma_{i,j}}\right)^{\ast}\varphi_{tj} \circ \left(\opi_{nt}^{G^nG^tG^j} \circ \overline{\gamma_{m,n} \times \gamma_{s,t} \times \gamma_{i,j}}\right)^{\ast}\varphi_{nt} \\
	&= \left(\overline{\gamma_{m,n} \times \gamma_{s,t} \times \gamma_{i,j}}\right)^{\ast}\left(\left(\opi_{tj}^{G^nG^tG^j}\right)^{\ast}\varphi_{tj} \circ \left(\opi_{nt}^{G^nG^tG^j}\right)^{\ast}\varphi_{nt}\right) \\
	&= \left(\overline{\gamma_{m,n} \times \gamma_{s,t} \times \gamma_{i,j}}\right)^{\ast}\left(\opi_{nj}^{G^nG^tG^j}\right)^{\ast}\varphi_{nj}  =\left(\opi_{nj}^{G^nG^tG^j} \circ \overline{\gamma_{m,n} \times \gamma_{s,t} \times \gamma_{i,j}}\right)^{\ast}\varphi_{nj} \\
	&= \left( \overline{\gamma_{m,n}\times \gamma_{i,j}} \circ \opi_{nj}^{\Gamma_{m,n}\Gamma_{s,t}\Gamma_{i,j}}\right)^{\ast}\varphi_{nj} = \left(\opi_{nj}^{\Gamma_{m,n}\Gamma_{s,t}\Gamma_{i,j}}\right)^{\ast}\left(\overline{\gamma_{m,n} \times \gamma_{i,j}}\right)^{\ast}\varphi_{nj} \\
	&= \left(\opi_{nj}^{\Gamma_{m,n}\Gamma_{s,t}\Gamma_{i,j}}\right)^{\ast}\varphi_{ij}^{mn}.
	\end{align*}
	This shows the ABL-cocycle condition and concludes the proof.
\end{proof}
Lemma \ref{Lemma: Section 4.2: Structure Isos for n > 0} and Proposition \ref{Prop: Section 4.2: Cocycle maps for all n,t geq 1 and all m,s} together tell  us that if we try to produce an object $\Glue(\Lambda^X)$ from an object of $\Dbeqsimp{X}$ by sending $A$ to $(A_0, A_{m,n}, \theta_{m,n})$, we get that there are isomorphisms $\theta_{m,n}$ for all $m \in \N \cup \lbrace \infty \rbrace$ and for all $n \in \N$ whenever $n \geq 1$, and moreover that these objects satisfy the ABL-cocycle condition (at least when we restrict our attention to simplicial maps that all lie away from degree $0$). While this is an excellent start, and indeed where the meat of the work lies, we still must discuss what happens when $n = 0$. In particular, we must first prove that there are indeed isomorphisms $\theta_{m,0}:\quo_{\Gamma_{m,0}}^{\ast}A_{m,0} \xrightarrow{\cong} (\pi_2^{\Gamma_{m,0}})^{\ast}A_0$ for all $m \in \N$ and that these objects also satisfy the ABL-cocycle condition, both in the restricted to degree $0$ sense and the total sense. Our strategy will now be as follows: We first show that the isomorphisms $\theta_{m,0}$ all exist and that the $A_{m,0}$ all satisfy the ABL-cocycle condition with respect to the other objects $A_{n,0}$ concentrated in degree $0$. After accomplishing this, we will finally prove that the $A_{m,n}$ satisfy the ABL-cocylce ccondition in all cases, which allows us to finally deduce that we can produce objects in $\Glue(\Lambda^X)$ via the assignment $(A_n) \mapsto (A_0, A_{m,n}, \theta_{m,n})$.

\begin{lemma}\label{Lemma: Section 4.2: Isos quotAm0 with pi2A0}
	Let $m \in \N$ and $A \in \Dbeqsimp{X}_0$. Then there are isomorphisms
	\[
	\theta_{m,0}:\quo_{\Gamma_{m,0}}^{\ast}A_{m,0} \xrightarrow{\cong} (\pi_2^{\Gamma_{m,0}})^{\ast}A_0.
	\]\index[notation]{thetam0@$\theta_{m,0}$}
\end{lemma}
\begin{proof}
	Begin by considering the commutative diagrams:
	\[
	\xymatrix{
		\Gamma_{m,0} \times X \ar[r]^-{\pi_2^{\Gamma_{m,0}}} \ar[d]_{\quo_{\Gamma_{m,0}}} & X & \Gamma_{m,0} \times X \ar[r]^-{\pi_2^{\Gamma_{m,0}}} \ar[d]_{\quo_{\Gamma_{m,0}}} & X \ar[d]^{\psi} \\
		\quot{X}{\Gamma_{m,0}} \ar[r]_-{\ogamma_{m,0}} & \quot{X}{G} \ar[u]_{\psi^{-1}} & \quot{X}{\Gamma_{m,0}} \ar[r]_-{\ogamma_{m,0}} & \quot{X}{G} 
	}
	\]
	With this we calculate that
	\begin{align*}
	\quo_{\Gamma_{m,0}}^{\ast}A_{m,0}
	&= \quo_{\Gamma_{m,0}}^{\ast} \ogamma_{m,0}^{\ast}(\psi)^{-1}A_0 = (\ogamma_{m,0} \circ \quo_{\Gamma_{m,0}})^{\ast}(\psi^{-1})^{\ast}A_0 \\
	&= (\psi \circ \pi_{2}^{\Gamma_{m,0}})^{\ast}(\psi^{-1})^{\ast}A_0 = (\pi_2^{\Gamma_{m,0}})^{\ast}\psi^{\ast}(\psi^{-1})^{\ast}A_0 \\
	&= (\pi_{2}^{\Gamma_{m,0}})^{\ast}(\psi^{-1} \circ \psi)^{\ast}A_0 \\
	&= (\pi_2^{\Gamma_{m,0}})^{\ast}A_0.
	\end{align*}
	Thus we set our isomorphism as
	\[
	\theta_{m,0}:= \id_{\quo_{\Gamma_{m,0}}^{\ast}A_{m,0}},
	\]
	which proves the lemma.
\end{proof}
\begin{proposition}\label{Prop: Section 4.2: ABL cocylce for n = 0}
	Let $m,n \in \N$. Then there are isomorphisms
	\[
	\varphi_{mn}:\left(\opi_{m}^{\Gamma_{m,0}\Gamma_{n,0}}\right)^{\ast}A_{m,0} \xrightarrow{\cong} \left(\opi_{n}^{\Gamma_{m,0}\Gamma_{n,0}}\right)^{\ast}A_{n,0}
	\]
	which satisfy the ABL-cocycle condition, where the projection morphisms are given in the diagram below:
	\[
	\xymatrix{
		& \Gamma_{m,0} \times \Gamma_{n,0} \times X \ar[dr]^-{\pi_{n}^{\Gamma_{m,0}\Gamma_{n,0}}} \ar[dl]_{\pi_{m}^{\Gamma_{m,0}\Gamma_{n,0}}} \\
		\Gamma_{m,0} \times X \ar[dr]_{\pi_2^{\Gamma_{m,0}}} & & \Gamma_{n,0} \times X \ar[dl]^{\pi_2^{\Gamma_{n,0}}} \\
		& X &
	}
	\]
\end{proposition}
\begin{proof}
	We note that the diagram above factors as in the diagram below:
	\[
	\xymatrix{
		& \Gamma_{m,0} \times \Gamma_{n,0} \times X \ar[dr]^-{\pi_{n}^{\Gamma_{m,0}\Gamma_{n,0}}} \ar[dl]_{\pi_{m}^{\Gamma_{m,0}\Gamma_{n,0}}} \\
		\Gamma_{m,0} \times X \ar[dr]_{\pi_2^{\Gamma_{m,0}}} \ar[r]^{\ogamma_{n,0}} & \quot{X}{G} \ar[d]^{\psi^{-1}} & \Gamma_{n,0} \times X \ar[dl]^{\pi_2^{\Gamma_{n,0}}} \ar[l]_{\ogamma_{n,0}} \\
		& X &
	}
	\]
	With this we calculate that
	\begin{align*}
	\left(\opi_{m}^{\Gamma_{m,0}\Gamma_{n,0}}\right)^{\ast}A_{m,0} &=\left(\opi_{m}^{\Gamma_{m,0}\Gamma_{n,0}}\right)^{\ast}\ogamma_{m,0}^{\ast}(\psi^{-1})^{\ast}A_0 = \left(\ogamma_{m,0} \circ \opi_{m}^{\Gamma_{m,0}\Gamma_{n,0}}\right)^{\ast} (\psi^{-1})^{\ast}A_0 \\
	&= \left(\ogamma_{n,0} \circ \opi_{n}^{\Gamma_{m,0}\Gamma_{n,0}}\right)^{\ast}(\psi^{-1})^{\ast}A_0 \\
	&= \left(\opi_{n}^{\Gamma_{m,0}\Gamma_{n,0}}\right)^{\ast}\ogamma_{n,0}^{\ast}(\psi^{-1})^{\ast}A_0 = \left(\opi_{n}^{\Gamma_{m,0}\Gamma_{n,0}}\right)^{\ast}A_{n,0}.
	\end{align*}
	Thus we define $\varphi_{mn}$ to be the identity morphism; that this satisfies the cocycle condition is trivial. This proves the lemma.
\end{proof}

We now prove that these objects finally satisfy the remaining total ABL-cocycle condition. For this we need to construct isomorphisms between the pullbacks of the $A_{m,n}$ and the $A_{s,0}$ where $n \geq 1$. It will be helpful to keep the diagram of projections
\[
\xymatrix{
	\Gamma_{m,n}\times \Gamma_{s,0} \times X \ar[rr]^-{\pi_{m,n}^{\Gamma_{m,n}\Gamma_{s,0}}} \ar[d]_{\pi_{s,0}^{\Gamma_{m,n}\Gamma_{s,0}}} & & \Gamma_{m,n} \times X \ar[d]^{\pi_2^{\Gamma_{m,n}}} \\
	\Gamma_{s,0}\times X \ar[rr]_-{\pi_2^{\Gamma_{s,0}}} & & X
}
\]
in mind, as we will use it to build the isomorphisms
\[
\varphi^{mn}_{s0}:\left(\opi_{m,n}^{\Gamma_{m,n}\Gamma_{s,0}}\right)^{\ast}A_{m,n} \xrightarrow{\cong} \left(\opi_{s,0}^{\Gamma_{m,n}\Gamma_{s,0}}\right)^{\ast}A_{s,0}.
\]
We will also need some basic observations about how the functors $\overline{R}_{G^n}$, for $n \geq 1$, interact with the functors induced by $\psi:X \to \quot{X}{G}$ and $\psi^{-1}:\quot{X}{G} \to X$ that play the role of the quotient map for the scheme $X$. To be more explicit about this, fix an $n \geq 1$ and note that the projection map $\pi_2^{G^n}:G^n \times X \to X$ is simplicial; $\pi_2$ in this case factors as a composition
\[
\pi_2^{G^n} = d_0^{n-1} \circ \cdots \circ d_0^{0}
\]
which we have seen before. Furthermore, the action of $G$ through this map is trivial. As such, there exists a unique morphism $\tilde{d}:\quot{X}{G^n} \to \quot{X}{G}$ making the diagrams
\[
\xymatrix{
	G^n \times X \ar[r]^-{d} \ar[d]_{\quo_{G^n}} & X \ar[d]^{\psi} & G^n \times X \ar[r]^-{d} \ar[d]_{\quo_{G^n}} & X \\
	\quot{X}{G^n} \ar[r]_-{\tilde{d}} & \quot{X}{G} & \quot{X}{G^n} \ar[r]_-{\tilde{d}} & \quot{X}{G} \ar[u]_{\psi^{-1}}
}
\]
commute. Because of this a routine check using the various constructions above shows that the diagram
\[
\xymatrix{
	\DbQl{X} \ar[r]^-{d^{\ast}} \ar[d]_{(\psi^{-1})^{\ast}} & \DbQl{G^n \times X} \ar[d]^{\overline{R}_{G^n}} \\
	\DbQl{\quot{X}{G}} \ar[r]_-{\tilde{d}^{\ast}} & \DbQl{\quot{X}{G^n}}
}
\]
commutes as well. This will be our last essential technical ingredient in proving the proposition below.

\begin{proposition}\label{Prop: Section 4.2: Isos and ABL-cocycle for n > 0 and t = 0}
	Let $n,s \in \N$ with $n \geq 1$, let $m \in \N \cup \lbrace  \infty \rbrace$, and let $A \in \Dbeqsimp{X}_0$. Then there are isomorphisms
	\[
	\varphi^{mn}_{s0}:\left(\opi_{m,n}^{\Gamma_{m,n}\Gamma_{s,0}}\right)^{\ast}A_{m,n} \xrightarrow{\cong} \left(\opi_{s,0}^{\Gamma_{m,n}\Gamma_{s,0}}\right)^{\ast}A_{s,0}
	\]\index[notation]{phivarmns0@$\varphi^{m,n}_{s,0}$}
	which satisfy the ABL-cocycle condition.
\end{proposition}
\begin{proof}
	We note that since the morphism $\pi_2^{G^n} = d:G^n \times X \to X$ is simplicial and $A \in \Dbeqsimp{X}_0$, there is an isomorphism $\alpha_d:d^{\ast}A_0 \xrightarrow{\cong} A_d$ in $\DbQl{G^n \times X}$. As such, we take successive pullbacks to get an isomorphism
	\[
	\left(\opi_{m,n}^{\Gamma_{m,n}\Gamma_{s,0}}\right)^{\ast}\ogamma_{m,n}^{\ast} \overline{R}_{G^n}(\alpha_d)
	\]
	of the form
	\[
	\left(\opi_{m,n}^{\Gamma_{m,n}\Gamma_{s,0}}\right)^{\ast}\ogamma_{m,n}^{\ast} \overline{R}_{G^n}(d^{\ast}A_0) \xrightarrow{\cong} \left(\opi_{m,n}^{\Gamma_{m,n}\Gamma_{s,0}}\right)^{\ast}\ogamma_{m,n}^{\ast} \overline{R}_{G^n}(A_n)
	\]
	in $\DbQl{\quot{X}{\Gamma_{m,n} \times \Gamma_{s,0}}}$. Note that 
	\[
	\left(\opi_{m,n}^{\Gamma_{m,n}\Gamma_{s,0}}\right)^{\ast}\ogamma_{m,n}^{\ast} \overline{R}_{G^n}(A_n) = \left(\opi_{m,n}^{\Gamma_{m,n}\Gamma_{s,0}}\right)^{\ast}A_{m,n}.
	\]
	We will now prove that
	\[
	\left(\opi_{m,n}^{\Gamma_{m,n}\Gamma_{s,0}}\right)^{\ast}\ogamma_{m,n}^{\ast} \overline{R}_{G^n}(d^{\ast}A_0) = \left(\opi_{s,0}^{\Gamma_{m,n}\Gamma_{s,0}}\right)^{\ast}A_{s,0}
	\]
	in order to show that the above map induces our desired isomorphism. For this we begin by using the diagram
	\[
	\xymatrix{
		\DbQl{X} \ar[r]^-{d^{\ast}} \ar[d]_{(\psi^{-1})^{\ast}} & \DbQl{G^n \times X} \ar[d]^{\overline{R}_{G^n}} \\
		\DbQl{\quot{X}{G}} \ar[r]_-{\tilde{d}^{\ast}} & \DbQl{\quot{X}{G^n}}
	}
	\]
	and deduce that $\overline{R}_{G^n}(d^{\ast}A_0) = \tilde{d}^{\ast}(\psi^{-1})^{\ast}A_0$. Furthermore, from the diagram describing the projection morphisms and the diagram defining $\tilde{d}$, we observe that the diagram
	\[
	\begin{tikzcd}
	\quot{X}{\Gamma_{m,n} \times \Gamma_{s,0}} \ar[rr]{}{\opi_{s,0}^{\Gamma_{m,n}\Gamma_{s,0}}} \ar[d, swap]{}{\opi_{s,0}^{\Gamma_{m,n}\Gamma_{s,0}}} & & \quot{X}{\Gamma_{s,0}} \ar[d]{}{\ogamma_{s,0}} \\
	\quot{X}{\Gamma_{m,n}} \ar[rr, swap]{}{\tilde{d}\circ \ogamma_{m,n}} \ar[dr, swap]{}{\ogamma_{m,n}} & & \quot{X}{G} \\
	& \quot{X}{G^n} \ar[ur, swap]{}{\tilde{d}} &
	\end{tikzcd}
	\]
	commutes as well. We use this and the induced relations to calculate that
	\begin{align*}
&	\left(\opi_{m,n}^{\Gamma_{m,n}\Gamma_{s,0}}\right)^{\ast}\ogamma_{m,n}^{\ast}\overline{R}_{G^n}(d^{\ast}A_0) \\
&= \left(\opi_{m,n}^{\Gamma_{m,n}\Gamma_{s,0}}\right)^{\ast} \ogamma_{m,n}\tilde{d}^{\ast}(\psi^{-1})^{\ast}A_0 = \left(\opi_{m,n}^{\Gamma_{m,n}\Gamma_{s,0}}\right)^{\ast}\left(\tilde{d} \circ \ogamma_{m,n}\right)^{\ast}(\psi^{-1})^{\ast}A_0 \\
	&= \left(\tilde{d} \circ \ogamma_{m,n} \circ \opi_{m,n}^{\Gamma_{m,n}\Gamma_{s,0}}\right)^{\ast}(\psi^{-1})^{\ast}A_0 = \left(\ogamma_{s,0} \circ \opi_{s,0}^{\Gamma_{m,n}\Gamma_{s,0}}\right)^{\ast}(\psi^{-1})^{\ast}A_0 \\
	&= \left(\opi_{m,n}^{\Gamma_{m,n}\Gamma_{s,0}}\right)^{\ast}\ogamma_{s,0}^{\ast}(\psi^{-1})^{\ast}A_0 = \left(\opi_{m,n}^{\Gamma_{m,n}\Gamma_{s,0}}\right)^{\ast}A_{s,0}.
	\end{align*}
	Thus  we attain an isomorphism 
	\[
	\left(\opi_{m,n}^{\Gamma_{m,n}\Gamma_{s,0}}\right)^{\ast}\ogamma_{m,n}^{\ast}\overline{R}_{G^n}(\alpha_d):\left(\opi_{s,0}^{\Gamma_{m,n}\Gamma_{s,0}}\right)^{\ast}A_{s,0} \xrightarrow{\cong} \left(\opi_{m,n}^{\Gamma_{m,n}\Gamma_{s,0}}\right)^{\ast}A_{m,n},
	\]
	so setting
	\[
	\varphi^{mn}_{s0} := \left(\opi_{m,n}^{\Gamma_{m,n}\Gamma_{s,0}}\right)^{\ast}\ogamma_{m,n}^{\ast}\overline{R}_{G^n}(\alpha_d)^{-1}
	\]
	provides our desired isomorphism.
	
	We now prove the ABL-cocycle condition. Let $n, s, j \in \N$ with $n \geq 1$. Assume for the moment that $j \geq 1$ (without loss of generality with $n \geq j$) and let $m,i \in \N \cup \lbrace \infty \rbrace$. To check the ABL-cocycle condition we must show in this case that
	\[
	\left(\opi_{n0}^{\Gamma_{m,n}\Gamma_{i,j}\Gamma_{s,0}}\right)^{\ast}\varphi_{s0}^{mn} = \left(\opi_{j0}^{\Gamma_{m,n}\Gamma_{i,j}\Gamma_{s,0}}\right)^{\ast}\varphi_{s0}^{ij} \circ  \left(\opi_{nj}^{\Gamma_{m,n}\Gamma_{i,j}\Gamma_{s,0}}\right)^{\ast}\varphi_{ij}^{mn},
	\]
	which is equivalent to showing that
	\[
	\left(\opi_{n0}^{\Gamma_{m,n}\Gamma_{i,j}\Gamma_{s,0}}\right)^{\ast}(\varphi_{s0}^{mn})^{-1} = \left(\opi_{nj}^{\Gamma_{m,n}\Gamma_{i,j}\Gamma_{s,0}}\right)^{\ast}(\varphi_{ij}^{mn})^{-1} \circ \left(\opi_{j0}^{\Gamma_{m,n}\Gamma_{i,j}\Gamma_{s,0}}\right)^{\ast}(\varphi_{s0}^{ij})^{-1}.  
	\]
	To write out the definition of these morphisms, we introduce the following commutative diagram
	\[
	\begin{tikzcd}
	& & \quot{X}{\Gamma_{m,n} \times \Gamma_{i,j} \times \Gamma_{s,0}} \ar[dll,swap]{}{\opi_{nj}^{\Gamma_{m,n}\Gamma_{i,j}\Gamma_{s,0}}} \ar[d]{}[description]{\opi_{n0}^{\Gamma_{m,n}\Gamma_{i,j}\Gamma_{s,0}}} \ar[drr]{}{\opi_{j0}^{\Gamma_{m,n}\Gamma_{i,j}\Gamma_{s,0}}} & \\
	\quot{X}{\Gamma_{m,n} \times \Gamma_{i,j}} \ar[d, swap]{}{\opi_{mn}^{\Gamma_{m,n}\Gamma_{i,j}}} &  & \quot{X}{\Gamma_{m,n} \times \Gamma_{s,0}} \ar[dll, near end]{}{\opi_{m,n}^{\Gamma_{m,n}\Gamma_{i,j}}} \ar[drr, near start]{}{\opi_{s0}^{\Gamma_{m,n}\Gamma_{s,0}}} & & \quot{X}{\Gamma_{i,j}\times \Gamma_{s,0}} \ar[d]{}{\opi_{s,0}^{\Gamma_{i,j}\Gamma_{s,0}}} \\
	\quot{X}{\Gamma_{m,n}} \ar[d, swap]{}{\ogamma_{m,n}} & & \quot{X}{\Gamma_{i,j}} \ar[d]{}{\ogamma_{i,j}} & & \quot{X}{\Gamma_{s,0}} \ar[d]{}{\ogamma_{s,0}} \\
	\quot{X}{G^n} \ar[rr, swap]{}{\overline{d}} \ar[rrrr, bend right = 30, swap]{}{\tilde{d}^{\prime\prime}} & & \quot{X}{G^j} \ar[rr, swap]{}{\tilde{d}^{\prime}} & & \quot{X}{G} \\
	\ar[from = 2-1, to = 3-3, crossing over, near start]{}{\opi_{ij}^{\Gamma_{m,n}\Gamma_{i,j}}}
	\ar[from = 2-5, to = 3-3, crossing over, near end]{}{\opi_{ij}^{\Gamma_{i,j}\Gamma_{s,0}}}
	\end{tikzcd}
	\]
	where $d:G^n \times X \to G^j \times X$ is the simplicial projection induced by the fact that $n \geq j \geq 1$ and $d^{\prime}:G^j \times X \to X, d^{\prime\prime}:G^n \times X \to X$ are the total projections with maps $\tilde{d}^{\prime}, \tilde{d}^{\prime\prime}$ induced via the diagrams:
	\[
	\xymatrix{
		G^n \times X \ar[r]^-{d^{\prime\prime}} \ar[d]_{\quo_{G^n}} & X \ar[d]^{\psi} & G^j \times X \ar[d]_{\quo_{G^j}} \ar[r]^-{d^{\prime}} & X \ar[d]^{\psi} \\
		\quot{X}{G^n} \ar[r]_-{\tilde{d}^{\prime\prime}} & \quot{X}{G} & \quot{X}{G^j} \ar[r]_-{\tilde{d}^{\prime}} & \quot{X}{G}
	}
	\]
	With this we must prove that
	\begin{align*}
	\left(\opi_{n0}^{\Gamma_{m,n}\Gamma_{i,j}\Gamma_{s,0}}\right)^{\ast}(\varphi_{s0}^{mn})^{-1} = \left(\opi_{n0}^{\Gamma_{m,n}\Gamma_{i,j}\Gamma_{s,0}}\right)^{\ast}\left(\opi_{mn}^{\Gamma_{m,n}\Gamma_{s,0}}\right)^{\ast}\ogamma_{m,n}^{\ast}\overline{R}_{G^n}(\alpha_{d^{\prime\prime}})
	\end{align*}
	is equivalent to
	\begin{align*}
	&\left(\opi_{nj}^{\Gamma_{m,n}\Gamma_{i,j}\Gamma_{s,0}}\right)^{\ast}(\varphi_{ij}^{mn})^{-1} \circ \left(\opi_{j0}^{\Gamma_{m,n}\Gamma_{i,j}\Gamma_{s,0}}\right)^{\ast}(\varphi_{s0}^{ij})^{-1} \\
	&=\left(\opi_{j0}^{\Gamma_{m,n}\Gamma_{i,j}\Gamma_{s,0}}\right)^{\ast}\left(\opi_{m,n}^{\Gamma_{m,n}\Gamma_{i,j}}\right)^{\ast}\ogamma_{m,n}^{\ast}\overline{R}_{G^n}(\alpha_d) \\
	&\circ \left(\opi_{nj}^{\Gamma_{m,n}\Gamma_{i,j}\Gamma_{s,0}}\right)^{\ast}\left(\opi_{i,j}^{\Gamma_{i,j}\Gamma_{s,0}}\right)^{\ast}\ogamma_{i,j}^{\ast}\overline{R}_{G^j}(\alpha_{d^{\prime}}) 
	\end{align*}
	or, more explicitly, that
	\begin{align*}
	&\left(\opi_{n0}^{\Gamma_{m,n}\Gamma_{i,j}\Gamma_{s,0}}\right)^{\ast}\left(\opi_{mn}^{\Gamma_{m,n}\Gamma_{s,0}}\right)^{\ast}\ogamma_{m,n}^{\ast}\overline{R}_{G^n}(\alpha_{d^{\prime\prime}}) \\
	&=\left(\opi_{j0}^{\Gamma_{m,n}\Gamma_{i,j}\Gamma_{s,0}}\right)^{\ast}\left(\opi_{m,n}^{\Gamma_{m,n}\Gamma_{i,j}}\right)^{\ast}\ogamma_{m,n}^{\ast}\overline{R}_{G^n}(\alpha_d) \\
	&\circ \left(\opi_{nj}^{\Gamma_{m,n}\Gamma_{i,j}\Gamma_{s,0}}\right)^{\ast}\left(\opi_{i,j}^{\Gamma_{i,j}\Gamma_{s,0}}\right)^{\ast}\ogamma_{i,j}^{\ast}\overline{R}_{G^j}(\alpha_{d^{\prime}}).
	\end{align*}
	However, this follows mutatis mutandis to the proof of Proposition \ref{Lemma: Section 4.2: Structure Isos for n > 0} and so is omitted.
	
	We now set $j = 0$ as the final case of the proof and let $i \in \N$. To prove the ABL-cocycle condition in this case we must show that 
	\[
	\left(\opi_{ns}^{\Gamma_{m,n}\Gamma_{i,0}\Gamma_{s,0}}\right)^{\ast}(\varphi_{s0}^{mn})^{-1} = \left(\opi_{ni}^{\Gamma_{m,n}\Gamma_{i,0}\Gamma_{s,0}}\right)^{\ast}(\varphi_{i0}^{mn})^{-1} \circ \left(\opi_{is}^{\Gamma_{m,n}\Gamma_{i,0}\Gamma_{s,0}}\right)^{\ast}(\varphi_{s0}^{i0})^{-1}  
	\]
	as in the case above (note that we are now indexing $\Gamma_{i,0}$ by $i$ and $\Gamma_{s,0}$ by $s$ in the projections to remove ambiguity). However, by Proposition \ref{Prop: Section 4.2: ABL cocylce for n = 0} we have that $\varphi_{s0}^{i0}$ is an identity morphism, so the proof of this case reduces to showing that
	\begin{align*}
	\left(\opi_{ns}^{\Gamma_{m,n}\Gamma_{i,0}\Gamma_{s,0}}\right)^{\ast}(\varphi_{s0}^{mn})^{-1} = \left(\opi_{ni}^{\Gamma_{m,n}\Gamma_{i,0}\Gamma_{s,0}}\right)^{\ast}(\varphi_{i0}^{mn})^{-1}.
	\end{align*}
	We proceed as above and write the diagram below in order to produce the definition of these morphisms:
	\[
	\begin{tikzcd}
	& & \quot{X}{\Gamma_{m,n} \times \Gamma_{i,0} \times \Gamma_{s,0}} \ar[dll,swap]{}{\opi_{ni}^{\Gamma_{m,n}\Gamma_{i,0}\Gamma_{s,0}}} \ar[d]{}[description]{\opi_{n0}^{\Gamma_{m,n}\Gamma_{i,0}\Gamma_{s,0}}} \ar[drr]{}{\opi_{is}^{\Gamma_{m,n}\Gamma_{i,0}\Gamma_{s,0}}} & \\
	\quot{X}{\Gamma_{m,n} \times \Gamma_{i,0}} \ar[d, swap]{}{\opi_{mn}^{\Gamma_{m,n}\Gamma_{i,0}}} &  & \quot{X}{\Gamma_{m,n} \times \Gamma_{s,0}} \ar[dll, near end]{}{\opi_{m,n}^{\Gamma_{m,n}\Gamma_{i,0}}} \ar[drr, near start]{}{\opi_{s0}^{\Gamma_{m,n}\Gamma_{s,0}}} & & \quot{X}{\Gamma_{i,0}\times \Gamma_{s,0}} \ar[d]{}{\opi_{s,0}^{\Gamma_{i,0}\Gamma_{s,0}}} \\
	\quot{X}{\Gamma_{m,n}} \ar[d, swap]{}{\ogamma_{m,n}} & & \quot{X}{\Gamma_{i,0}} \ar[d]{}{\ogamma_{i,0}} & & \quot{X}{\Gamma_{s,0}} \ar[d]{}{\ogamma_{s,0}} \\
	\quot{X}{G^n} \ar[rr, swap]{}{\overline{d}} \ar[rrrr, bend right = 30, swap]{}{\tilde{d}} & & \quot{X}{G} \ar[rr, equals]{}{} & & \quot{X}{G} \\
	\ar[from = 2-1, to = 3-3, crossing over, near start]{}{\opi_{i0}^{\Gamma_{m,n}\Gamma_{i,0}}}
	\ar[from = 2-5, to = 3-3, crossing over, near end]{}{\opi_{i0}^{\Gamma_{i,0}\Gamma_{s,0}}}
	\end{tikzcd}
	\]
	As we can see from the diagram, in this case we actually have that $\tilde{d}^{\prime} = \id_{\quot{X}{G}}$ and the map $\tilde{d}$ in the $(m,n),(i,0)$ case and $\tilde{d}^{\prime\prime}$ in the $(m,n),(s,0)$ case actually coincide, as they are both induced from the total projection $d:G^n \times X \to X$. As such, a routine verification shows that
	\begin{align*}
	\left(\opi_{ns}^{\Gamma_{m,n}\Gamma_{i,0}\Gamma_{s,0}}\right)^{\ast}(\varphi_{s0}^{mn})^{-1} &= \left(\opi_{ns}^{\Gamma_{m,n}\Gamma_{i,0}\Gamma_{s,0}}\right)^{\ast}\left(\opi_{m,n}^{\Gamma_{m,n}\Gamma_{s,0}}\right)^{\ast}\ogamma_{m,n}\overline{R}_{G^n}(\alpha_d) \\
	&=\left(\opi_{ni}^{\Gamma_{m,n}\Gamma_{i,0}\Gamma_{s,0}}\right)^{\ast}\left(\opi_{m,n}^{\Gamma_{m,n}\Gamma_{i,0}}\right)^{\ast}\ogamma_{m,n}\overline{R}_{G^n}(\alpha_d) \\
	&= \left(\opi_{ni}^{\Gamma_{m,n}\Gamma_{i,0}\Gamma_{s,0}}\right)^{\ast}(\varphi_{i0}^{mn})^{-1}.
	\end{align*}
	This shows the ABL-cocycle condition in the $j=0$ case and hence proves the proposition.
\end{proof}
\begin{corollary}\label{Cor: Section 4.2: ABL Cocycle for all n}
	Let $A \in \Dbeqsimp{X}$. Then for any $j,n \in \N$ and any $i,m \in \N$ or $i,m \in \N \cup \lbrace \infty \rbrace$, there are isomorphisms
	\[
	\varphi_{ij}^{mn}:\left(\opi_{mn}^{\Gamma_{m,n}\Gamma_{i,j}}\right)^{\ast}A_{m,n} \xrightarrow{\cong} \left(\opi_{ij}^{\Gamma_{m,n}\Gamma_{i,j}}\right)^{\ast}A_{i,j}
	\]\index[notation]{phivarmnij@$\varphi_{ij}^{mn}$}
	which satisfy the ABL-cocycle condition.
\end{corollary}
\begin{proof}
	Apply Propositions \ref{Prop: Section 4.2: Cocycle maps for all n,t geq 1 and all m,s}, \ref{Prop: Section 4.2: ABL cocylce for n = 0}, and \ref{Prop: Section 4.2: Isos and ABL-cocycle for n > 0 and t = 0}.
\end{proof}

We now must prove how the functor $\overline{R}$ acts on morphisms. While the majority of our task is the construction of the objects, it is crucial that we do not forget how we manipulate morphisms.
\begin{proposition}\label{Prop: Section 4.2: Rhos commute with cocycle and structure isos}
	Let $\rho \in \Dbeqsimp{X}(A,B)$ and define
	\[
	\rho_{m,n} := \begin{cases}
	\ogamma_{m,0}^{\ast}(\psi^{-1})^{\ast}\rho_0 & \text{if}\,\, n = 0, m \in \N; \\
	\ogamma_{m,n}^{\ast}\overline{R}_{G^n}(\rho_{n}) & \text{if}\,\, n \geq 1, m \in \N \cup \lbrace \infty \rbrace.
	\end{cases}
	\]\index[notation]{rhomn@$\rho_{m,n}$}
	Then for all pairs $(m,n)$ and $(s,t)$ for which $\rho_{m,n}$ and $\rho_{s,t}$ are defined above, the diagrams
	\[
	\xymatrix{
		(\pi_{2}^{\Gamma_{m,n}})^{\ast}A_0 \ar[rr]^-{(\pi_2^{\Gamma_{m,n}})^{\ast}} \ar[d]_{\theta_{m,n}^{A}} & & (\pi_2^{\Gamma_{m,n}})^{\ast} \ar[d]^{\theta_{m,n}^{B}} \\
		\quo_{\Gamma_{m,n}}^{\ast}A_{m,n} \ar[rr]_-{\quo_{\Gamma_{m,n}}^{\ast}\rho_{m,n}} & & \quo_{\Gamma_{m,n}}^{\ast}B_{m,n}
	}
	\]
	\[
	\xymatrix{
		(\opi_{m,n}^{\Gamma_{m,n}\Gamma_{s,t}})^{\ast}A_{m,n} \ar[rr]^-{(\opi_{m,n}^{\Gamma_{m,n}\Gamma_{s,t}})^{\ast}\rho_{m,n}} \ar[d]_{(\quot{\varphi_{st}^{mn}}{A})} & & (\opi_{m,n}^{\Gamma_{m,n}\Gamma_{s,t}})^{\ast}B_{m,n} \ar[d]^{\quot{\varphi_{st}^{mn}}{B}} \\
		(\opi_{s,t}^{\Gamma_{m,n}\Gamma_{s,t}})^{\ast}A_{s,t} \ar[rr]_-{(\opi_{s,t}^{\Gamma_{m,n}\Gamma_{s,t}})^{\ast}\rho_{s,t}} & & (\opi_{s,t}^{\Gamma_{m,n}\Gamma_{s,t}})^{\ast}B_{s,t}
	}
	\]
	both commute.
\end{proposition}
\begin{proof}
	We first verify that all the $\theta_{m,n}$ diagrams commute. For this set $n = 0$ and let $m \in \N$. Then we find from Lemma \ref{Lemma: Section 4.2: Isos quotAm0 with pi2A0} that the isomorphism $\theta_{m, 0}$ is the identity morphis, so it suffices to show in this case that
	\[
	(\pi_2^{\Gamma_{m,0}})^{\ast}\rho_0 = \quo_{\Gamma_{m,0}}^{\ast}\rho_{m,0}.
	\]
	For this we recall the commuting diagram
	\[
	\xymatrix{
		\Gamma_{m,0} \times X \ar[rr]^-{\pi_2^{\Gamma_{m,0}}} \ar[d]_{\quo_{\Gamma_{m,0}}} & & X \\
		\quot{X}{\Gamma_{m,0}} \ar[rr]_-{\ogamma_{m,0}} & & \quot{X}{G} \ar[u]_{\psi^{-1}}
	}
	\]
	of schemes. With this we derive that
	\begin{align*}
	(\pi_2^{\Gamma_{m,0}})^{\ast}\rho_0 &= (\psi^{-1} \circ \ogamma_{m,0} \quo_{\Gamma_{m,0}})^{\ast}\rho_0 = \quo_{\Gamma_{m,0}}^{\ast}\ogamma_{m,0}^{\ast}(\psi^{-1})^{\ast}\rho_0 \\
	&= \quo_{\Gamma_{m,0}}^{\ast}\rho_{m,0},
	\end{align*}
	as was desired.
	
	We now establish that the $\rho_{m,n}$ commute suitably with the $\theta_{m,n}$ when $n \geq 1$. Thus let $n \geq 1$ and first set $m = \infty$. In this case we recall that the morphism $\pi_2:G^n \times X \to X$ is a simplicial morphism $d$, so there is an isomorphism
	\[
	\alpha_d:d^{\ast}A_0 \xrightarrow{\cong} A_n.
	\] 
	and from Lemma \ref{Lemma: Section 4: Quot pullback for An} that there is a natural isomorphism 
	\[
	\nu:\id_{\DbQl{G^n \times X}} \to \quo_{G^n}^{\ast} \circ \overline{R}_{G^n}.
	\] 
	By Lemma \ref{Lemma: Section 4.2: Structure Isos for n > 0} the isomorphism $\theta_{\infty,n}$ is given by
	\[
	\theta_{\infty,n} := \nu_{A_n} \circ \alpha_d.
	\]
	Observe that $\quo_{G^n}^{\ast}\overline{R}_{G^n}(A_n) = \quo_{G^n}^{\ast}A_{\infty,n}$. To establish the commuting diagram we first note that since $\rho \in \Dbeqsimp{X}$, the diagram
	\[
	\xymatrix{
		d^{\ast}A_0 \ar[r]^-{d^{\ast}\rho_0} \ar[d]_{\alpha_d^A} & d^{\ast}B_0 \ar[d]^{\alpha_d^B} \\
		A_{n} \ar[r]_-{\rho_n} & B_n
	}
	\]
	commutes. Similarly, from the naturality of $\nu$ we find that the diagram
	\[
	\xymatrix{4
		A_n \ar[rr]^-{\rho_n} \ar[d]_{\nu_{A_n}} & & B_n \ar[d]^{\nu_{B_n}} \\
		\quo_{G^n}^{\ast}\overline{R}_{G^n}(A_n) \ar[rr]_-{\quo_{G^n}^{\ast}\overline{R}_{G^n}(\rho_n)} & & \quo_{G^n}^{\ast}\overline{R}_{G^n}(B_n)
	}
	\]
	commutes as well. Using that $\overline{R}_{G^n}(\rho_n) = \rho_{\infty,n}$  and stacking the diagrams gives that
	\[
	\xymatrix{
		d^{\ast}A_0 \ar[rr]^-{d^{\ast}\rho_0} \ar[d]_{\alpha_d^A} & & d^{\ast}B_0 \ar[d]^{\alpha_{d}^{B}} \\
		A_n \ar[d]_{\nu_{A_n}} \ar[rr]_-{\rho_n} & & B_n \ar[d]^{\nu_{B_n}} \\
		\quo_{G^n}^{\ast}A_{\infty,n} \ar[rr]_-{\quo_{G^n}^{\ast}\rho_{\infty,n}} & & \quo_{G^n}^{\ast}B_{\infty,n}
	}
	\]
	commutes. However, the outer edges compose to
	\[
	\xymatrix{
		d^{\ast}A_0 \ar[rr]^-{d^{\ast}\rho_0} \ar[d]_{\theta_{\infty,n}^A} & & d^{\ast}B_0 \ar[d]^{\theta_{\infty,n}^{B}} \\
		\quo_{G^n}^{\ast}A_{\infty,n} \ar[rr]_-{\quo_{G^n}^{\ast}\rho_{\infty,n}} & & \quo_{G^n}^{\ast}B_{\infty,n}
	}
	\]
	which establishes the result for $m = \infty$ and $n \geq 1$. 
	
	For any $m \in \N$, we recall that $\theta_{m,n} =(\gamma_{m,n} \times \id_X)^{\ast}\theta_{\infty,n}$, where the morphism is given as in the diagram below:
	\[
	\begin{tikzcd}
	\Gamma_{n,m} \times X \ar[rrrr, bend left = 30]{}{\pi_{2}^{m,n}} \ar[rr]{}{\gamma_{m,n} \times \id_X} \ar[d, swap]{}{\quo_{\Gamma_{m,n}}} & & G^n \times X \ar[rr]{}{\pi_2^n} \ar[d]{}{\quo_{G^n}} & & X \\
	\quot{X}{\Gamma_{m,n}} \ar[rr, swap]{}{\overline{\gamma}_{m,n}} & & \quot{X}{G^n}
	\end{tikzcd}
	\]
	Applying $(\gamma_{m,n} \times \id_X)^{\ast}$ to the diagram
	\[
	\xymatrix{
		d^{\ast}A_0 \ar[rr]^-{d^{\ast}\rho_0} \ar[d]_{\theta_{\infty,n}^A} & & d^{\ast}B_0 \ar[d]^{\theta_{\infty,n}^{B}} \\
		\quo_{G^n}^{\ast}A_{\infty,n} \ar[rr]_-{\quo_{G^n}^{\ast}\rho_{\infty,n}} & & \quo_{G^n}^{\ast}B_{\infty,n}
	}
	\]
	gives the new commmuting diagram:
	\[
	\xymatrix{
		(\gamma_{m,n} \times \id_X)^{\ast}d^{\ast}A_0 \ar[rrr]^-{(\gamma_{m,n} \times \id_X)^{\ast}d^{\ast}\rho_0} \ar[d]_{\theta_{m,n}^A} &  & & (\gamma_{m,n} \times \id_X)^{\ast}d^{\ast}B_0 \ar[d]^{\theta_{m,n}^{B}} \\
		(\gamma_{m,n} \times \id_x)^{\ast}\quo_{G^n}^{\ast}A_{\infty,n} \ar[rrr]_-{(\gamma_{m,n} \times \id_X)^{\ast}\quo_{G^n}^{\ast}\rho_{\infty,n}} &  & & (\gamma_{m,n} \times \id_X)^{\ast}\quo_{G^n}^{\ast}B_{\infty,n}
	}
	\]
	This diagram is equal to our desired commuting diagram, so this establishes the result for the $\theta_{m,n}$.
	
	We now establish the result for the $\varphi_{st}^{mn}$ where without loss of generality $n \geq t$ for $n,t \in \N$. In all cases, it suffices to prove the desired commuting diagrams for the morphisms $(\varphi_{st}^{mn})^{-1}$, so this is what we will do. First consider the case where $t =0$ and $n \geq 1$. In this case we let $s \in \N$ and $m \in \N \cup \lbrace \infty \rbrace$. Now recall that by Proposition \ref{Prop: Section 4.2: Isos and ABL-cocycle for n > 0 and t = 0}
	\[
	(\quot{\varphi_{s0}^{mn}}{A})^{-1} = (\opi_{m,n}^{\Gamma_{m,n}\Gamma_{s,0}})^{\ast}\ogamma_{m,n}^{\ast}\overline{R}_{G^n}(\quot{\alpha_d}{A}),
	\] 
	where $d:G^n \times X \to X$ is the projection $\pi_2^{G^n}:G^n \times X \to X$ and the map $\quot{\alpha_d}{A}$ is the isomorphism asserted by the fact that $A \in \Dbeqsimp{X}_0$; the morphism $(\quot{\varphi_{s0}^{mn})^{-1}}{B}$ and $\quot{\alpha_d}{B}$ are defined similarly. Now using that $\rho$ is a morphism in $\Dbeqsimp{X}$ implies that the diagram
	\[
	\xymatrix{
		d^{\ast}A_0 \ar[r]^-{d^{\ast}A_0} \ar[d]_{\quot{\alpha_d}{A}} & d^{\ast}B_0 \ar[d]^{\quot{\alpha_d}{B}} \\
		A_n \ar[r]_-{\rho_n} & B_n
	}
	\]
	commutes, so applying the functor $(\opi_{m,n}^{\Gamma_{m,n}\Gamma_{s,0}})^{\ast} \circ \ogamma_{m,n}^{\ast} \circ \overline{R}_{G^n}$ and using the induced identities of Proposition \ref{Prop: Section 4.2: Isos and ABL-cocycle for n > 0 and t = 0} again we get that the diagram
	\[
	\xymatrix{
		(\opi_{s,0}^{\Gamma_{m,n}\Gamma_{s,0}})^{\ast}A_{s,0} \ar[rrr]^-{(\opi_{s,0}^{\Gamma_{m,n}\Gamma_{s,0}})^{\ast}\rho_{s,0}} \ar[d]_{\quot{(\varphi_{s0}^{mn})^{-1}}{A}} & & & (\opi_{s,0}^{\Gamma_{m,n}\Gamma_{s,0}})^{\ast}B_{s,0} \ar[d]^{\quot{(\varphi_{s0}^{mn})^{-1}}{B}} \\
		(\opi_{m,n}^{\Gamma_{m,n}\Gamma_{s,0}})^{\ast}A_{m,n} \ar[rrr]_-{(\opi_{m,n}^{\Gamma_{m,n}\Gamma_{s,0}})^{\ast}\rho_{m,n}} & & & (\opi_{m,n}^{\Gamma_{m,n}\Gamma_{s,0}})^{\ast}B_{m,n}
	}
	\]
	commutes as well. This is our desired commuting diagram and establishes the $n \geq 1, t = 0$ cases.
	
	Now let $n = 0 = t$ and let $m,s \in \N$. In this case we have from Proposition \ref{Prop: Section 4.2: ABL cocylce for n = 0} that the isomorphisms $\quot{\varphi_{s0}^{m0}}{A}$ are the corresponding identity morphisms, so in this case it suffices to check that the diagram
	\[
	\xymatrix{
		(\opi_{m,0}^{\Gamma_{m,0}\Gamma_{s,0}})^{\ast}A_{m,0} \ar@{=}[d] \ar[rrr]^-{(\opi_{m,0}^{\Gamma_{m,0}\Gamma_{s,0}})^{\ast}\rho_{m,0}} & & & (\opi_{m,0}^{\Gamma_{m,0}\Gamma_{s,0}})^{\ast}B_{m,0} \ar@{=}[d] \\
		(\opi_{s,0}^{\Gamma_{m,0}\Gamma_{s,0}})^{\ast}A_{s,0} \ar[rrr]_-{(\opi_{s,0}^{\Gamma_{m,0}\Gamma_{s,0}})^{\ast}\rho_{s,0}} & & & (\opi_{s,0}^{\Gamma_{m,0}\Gamma_{s,0}})^{\ast}B_{s,0}
	}
	\]
	commutes. However, this is a routine check by using the identities induced from the commuting diagram
	\[
	\xymatrix{
		& \quot{X}{\Gamma_{m,0}} \ar[dr]^{\ogamma_{m,0}} \\
		\quot{X}{\Gamma_{m,0} \times \Gamma_{s,0}} \ar[ur]^{\opi_{m,0}^{\Gamma_{m,0}\Gamma_{s,0}}} \ar[dr]_{\opi_{s,0}^{\Gamma_{m,0}\Gamma_{s,0}}} & & \quot{X}{G} \ar[r]^-{\psi^{-1}} & X \\
		& \quot{X}{\Gamma_{s,0}} \ar[ur]_{\ogamma_{s,0}}
	}
	\]
	and proceeding as in the proof of Proposition \ref{Prop: Section 4.2: ABL cocylce for n = 0}. This establishes the $n = 0 = t$ case.
	
	We finally proceed in the case $n, t \geq 1$ (without loss of generality with $n \geq t$). Let $s,m \in \N \cup \lbrace \infty \rbrace$ and recall that in this case we again have
	\[
	\quot{(\varphi_{st}^{mn})^{-1}}{A} = \left(\opi_{m,n}^{\Gamma_{m,n}\Gamma_{s,t}}\right)^{\ast}\ogamma_{m,n}^{\ast}\overline{R}_{G^n}(\quot{\alpha_d}{A})
	\]
	where $d:G^n \times X \to G^t \times X$ is the projection morphism (or identitiy if $n = t$) and the morphism $\quot{\alpha_d}{A}:d^{\ast}A_t \xrightarrow{\cong} A_n$ is given from the fact that $A \in \Dbeqsimp{X}_0$. From here proceeding as in the $n \geq 1, t = 0$ case establishes the commutativity of the diagram and completes the proof of the proposition.
\end{proof}
\begin{corollary}\label{Cor: Section 4.2: Functor from simplicial to glue}
	There is a functor $\overline{R}:\Dbeqsimp{X} \to \Glue(\Lambda^X)$\index[notation]{Roverline@$\overline{R}$} given by
	\[
	(A_n)_{n \in \N} \mapsto (A_0, A_{m,n}, \theta_{m,n})
	\]
	on objects and
	\[
	(\rho)_{n \in \N} \mapsto (\rho_0, \rho_{m,n})
	\]
	on morphisms.
\end{corollary}
\begin{proof}
	That this construction preserves composition and identities is immediate, so we only need to check that the object and morphism assignments are indeed objects of $\Glue(\Lambda^X)$. However, for objects this is verified by Lemma \ref{Lemma: Section 4.2: Structure Isos for n > 0}, Lemma \ref{Lemma: Section 4.2: Isos quotAm0 with pi2A0}, and Corollary \ref{Cor: Section 4.2: ABL Cocycle for all n}; for morphisms, this is Proposition \ref{Prop: Section 4.2: Rhos commute with cocycle and structure isos}.
\end{proof}
\begin{corollary}\label{Cor: Section 4.2: Functor from splicial to derived}
	There is a functor $R:\Dbeqsimp{X} \to \DbeqQl{X}$ given by the composition:
	\[
	\xymatrix{
		\Dbeqsimp{X} \ar[r]^-{\overline{R}} & \Glue(\Lambda^X) \ar[r]^-{\simeq} & \quot{\DbeqQl{X}}{ABL} \ar[r]^-{\simeq} & \DbeqQl{X}
	}
	\]\index[notation]{R@$R$}
\end{corollary}

We have (finally) developed up all the necessary machinery and techniques to prove the simplicial equivalence. This is presented briefly below, as it is a routine application of the material built up to this point.
\begin{Theorem}\label{Theorem: Section 4.2: Simplicial equivalence with derived category}
	There is an equivalence of categories
	\[
	\begin{tikzcd}
	\DbeqQl{X} \ar[r, bend left = 30, ""{name = U}]{}{L} & \Dbeqsimp{X} \ar[l, bend left = 30, ""{name = L}]{}{R} \ar[from = U, to = L, symbol = \simeq]
	\end{tikzcd}
	\]
	where $L$ is the functor of Proposition \ref{Prop: Section 4: The functor L from derived equiv cat to simplicial equiv cat} and $R$ is the functor of Corollary \ref{Cor: Section 4.2: Functor from splicial to derived}.
\end{Theorem}
\begin{proof}
	It suffices to prove that for any object $S \in \Dbeqsimp{X}_0$, there are pointwise isomorphisms
	\[
	L(R(S))_n \cong S_n
	\]
	for all $n \in \N$ in order to verify that $L(R(S)) \cong S$; note that this is because these isomorphisms are already sufficiently compatible with the simplicial pullbacks by construction, so providing their existence for each $n \in \N$ gives rise our natural isomorphism at $S$. Beginning this task, we note that when $n \geq 1$ that by construction upon chasing $R(S)$ through its definition, we have that the $\quot{X}{G^n}$-components of $R(S)$ are given by
	\[
	\quot{R(S)}{G^n} = \overline{R}_{G^n}(S_n).
	\]
	Consequently we calculate that
	\begin{align*}
	L(R(S))_n &= L_{G^n}(\quot{R(S)}{G^n}(S))_0 = L_{G^n}(\overline{R}_{G^n}(S_n))_0 = L_{G^n}(R_{G^n}(\underline{S}_n))_0 \\
	&= L_{G^n}(R_{G^n}(S_n)) \cong S_n.
	\end{align*}
	Alternatively, if $n = 0$ we find that
	\[
	L(R(S))_0 = \psi^{\ast}\quot{R(S)}{G} = (\psi)^{\ast}(\psi^{-1})^{\ast}S_0 = (\psi \circ \psi^{-1})^{\ast}S_0 = S_0,
	\]
	which in turn establishes the isomorphism for $n = 0$. Combining this with the construction above we get that
	\[
	(L \circ R)(S) \cong S
	\]
	and hence that $L \circ R \cong \id_{\Dbeqsimp{X}}$.
	
	Now fix an $A \in \DbeqQl{X}_0$. Because of how the objects are propagated through the functor $R$ by way of acyclic descent, in order to conclude that $R \circ L \cong \id_{\DbeqQl{X}}$ it suffices to prove that for all $\Gamma_{m,n} \times X \in \Lambda^X$, $\quot{A}{\Gamma_{m,n}} \cong \quot{(R \circ L)(A)}{\Gamma_{m,n}}$. For this we first let $n \geq 1$ and let $m \in \N \cup \lbrace \infty \rbrace$. In this case we recall that the morphism $\ogamma_{m,n}:\quot{X}{\Gamma_{m,n}} \to \quot{X}{G^n}$ is induced from the diagram:
	\[
	\xymatrix{
		\Gamma_{m,n} \times X \ar[rr]^-{\gamma_{m,n}\times \id_X} \ar[d]_{\quo_{\Gamma_{m,n}}} & & G^n \times X \ar[d]^{\quo_{G^n}} \\
		\quot{X}{\Gamma_{m,n}} \ar[rr]_-{\ogamma_{m,n}} & & \quot{X}{G^n}
	}
	\]
	Because $\gamma_{m,n}:\Gamma_{m,n}\to G^n$ is an $m$-acyclic resolution of $G^n$, it is routine to check that $\gamma_{m,n} \in \Sf(G)(\Gamma_{m,n},G^n)$ and $\gamma_{m,n} \times \id_X:\Gamma_{m,n} \times X \to G^n \times X$ is a morphism in $\SfResl_G(X)$. As such we have that there is an isomorphism
	\[
	\tau_{\gamma_{m,n}}^A:\ogamma_{m,n}^{\ast}(\quot{A}{G^n}) \xrightarrow{\cong} \quot{A}{\Gamma_{m,n}}
	\]
	in $\DbQl{\quot{X}{\Gamma_{m,n}}}$. Using this we calculate (also with the use of the natural isomorphism $R_{G^n} \circ L_{G^n} \cong \id_{\DbQl{\quot{X}{G^n}}}$) that
	\[
	\quot{R(L(A))}{\Gamma_{m,n}} = \ogamma_{m,n}^{\ast}R_{G^n}(L_{G^n}(\quot{A}{G^n})) \cong \ogamma_{m,n}^{\ast}(\quot{A}{G^n}) \cong \quot{A}{\Gamma_{m,n}},
	\]
	which establishes what was desired whenever $n \geq 1$. Now set $n = 0$. In this case we have that
	\[
	\quot{R(L(A))}{\Gamma_{m,n}} = \ogamma_{m,0}^{\ast}(\psi^{-1})^{\ast}L(A)_0 = \ogamma_{m,0}^{\ast}(\psi^{-1})^{\ast}\psi^{\ast}\quot{A}{G} = \ogamma_{m,0}^{\ast}\quot{A}{G}
	\]
	where $\ogamma_{m,0}$ is the morphism induced via the diagram:
	\[
	\xymatrix{
		\Gamma_{m,0} \times X \ar[r]^-{\pi_2} \ar[d]_{\quo_{\Gamma_{m,0}}} & X \ar[d]^{\psi} \\
		\quot{X}{\Gamma_{m,0}} \ar[r]_-{\ogamma_{m,0}} & \quot{X}{G}
	}
	\]
	We claim that $\quot{A}{\Gamma_{m,0}} \cong \ogamma_{m,0}^{\ast}(\quot{A}{G})$. However, from the construction of the isomorphism $\psi$ as a trivialization of the action of $G$ on $X$, together with the uniqueness of the morphism $\ogamma_{m,n}$, the equivariant descent induced by Theorem \ref{Thm: Section 4: EDC is ABL EDC}, and the commutativity of the diagram
	\[
	\xymatrix{
		& \quot{X}{G \times \Gamma_{m,0}} \ar[dr]^{\opi_{G}^{G\Gamma_{m,0}}} \ar[dl]_{\opi_{m,0}^{G\Gamma_{m,0}}} & \\
		\quot{X}{\Gamma_{m,0}} \ar[rr]_-{\ogamma_{m,0}} & & \quot{X}{G}
	}
	\]
	we derive that $\ogamma_{m,0}^{\ast}\quot{A}{G} \cong \quot{A}{\Gamma_{m,0}}$. This allows us to finally deduce that
	\[
	R \circ L \cong \id_{\DbeqQl{X}}
	\]
	and hence that $L$ and $R$ are inverse equivalences between the categories $\DbeqQl{X}$ and $\Dbeqsimp{X}$, as was to be shown.
\end{proof}

\chapter{The Stacky Equivariant Derived Category}\label{Section: Simplicial EDC is the stcky EDC}

We are now at the last climb before we can summit Mount Equivariant Derived Category; we only need to climb the equivalence of the stacky equivariant derived category at this point. To handle this final stretch, it will be particularly convenient to use a simplicial presentation of the category $\DbQl{[G \backslash X]}$ which we will show is compatible with the simplicial category $\Dbeqsimp{X}$. The way in which we prove this simplicial presentation of $\DbeqQl{[G \backslash X]}$ is equivalent to that of $\Dbeqsimp{X}$ will go through the language of coskeletal simplicial schemes and the presentation of the algebraic stack $[G \backslash X]$ as a quotient of a groupoid internal to the category of varieties. Many of the results here are known (cf.\@ \cite{OlssonSpacesStacks} and \cite{Champs} for textbook accounts), but we provide a quick exposition of some of these known results in a different (more categorical) style than how they are usually presented in Appendix \ref{Appendix: Stackification via Torsors} below. Primarily we will need that the stack $[G \backslash X]$ is the fppf-stackification of the internal groupoid in $\Var_{/\Spec K}$ induced by the $G$-action on $X$ (that is $[G \backslash X]$ is the stackification of the pseudofunctor $(\Gbb \times \Xbb)(-)$ with respect to the fppf topology) to deduce that the simplicial approximation of $[G \backslash X]$ is in fact a $2$-coskeletal simplicial scheme. This is a fundamental tool of ours, as it allows us to deduce an isomorphism of categories
\[
\Dbeqsimp{X} \cong \Dbeqstack{[G \backslash X]}
\]
(cf.\@ \ref{Cor: Section 4: Bounded derived simplicial cat is equiv}) which we will see is indispensable in our proof of equivalence of categories $\Dbeqsimp{X} \simeq \DbQl{[G \backslash X]}$. For readers unfamiliar with the algebraic stack $[G \backslash X]$, please see Appendix \ref{Appendix: Stackification via Torsors} for a basic definition and description.

We begin our study of the stacky equivariant derived category by defining the internal groupoid $\mathbb{G\times X}$ induced by the action of the smooth affine algebraic group $G$ on the left $G$-variety $X$. It is worth noting that by Corollary \ref{Cor: Section 4: Internal groupoid} and Proposition \ref{Prop: Section 4: Equivalence of quotient stack and symmetry stack} the fppf-stackification of the induced fibration of groupoids $(\Gbb \times \Xbb)^{+}$, is isomorphic to the quotient stack $[G \backslash X]$ so we can use the simplicial data that the object $\Gbb \times \Xbb$ has at hand to induce the simplicial information of $[G \backslash X]$ and make it more clear that such a construction is actually $2$-coskeletal (cf.\@ Lemma \ref{Lemma: Section 4: Simplicial presentation of a stack is 2coskeletal}).
\begin{definition}\label{Defn: Essentially alg theory internal gropuoid}\index{Internal Groupoid!Action Groupoid}\index[notation]{GtimesX@$\mathbb{G \times X}$}
	Let $G$ be a smooth affine algebraic group and let $X$ be a left $G$-variety. Then the collection
	\[
	\Gbb \times \Xbb := \begin{cases}
	(\Gbb \times \Xbb)_0 := X; \\
	(\Gbb \times \Xbb)_1 := G \times X; \\
	\Dom := \pi_2:G \times X \to X; \\
	\codom := \alpha_X:G \times X \to X; \\
	m = \mu \times \id_X: G \times G \times X \to G \times X; \\
	s_0 = \langle 1_G, \id_X \rangle:X \to G \times X; \\
	\inv = \langle -1 \circ \pi_1, \alpha_X \rangle:G \times X \to G \times X;
	\end{cases}
	\]
	where $-1$ is the inverse morphism of $G$, is an internal groupoid in $\Var_{/\Spec K}$ (cf.\@ Corollary \ref{Cor: Appendix A: Interal action gropoid}).
\end{definition}

The biggest benefit of using the internal groupoid presentation of $[G \backslash X]$ is that it makes it much clearer how to provide a simplicial approximation of $[G \backslash X]$. In particular, the natural simplicial presentation of $(\Gbb \times \Xbb)$ by the objects with $(\Gbb \times \Xbb)_0 = X$ and
\[
(\Gbb \times \Xbb)_{n} := \underbrace{(G \times X) \times_X (G \times X) \times_X \cdots \times_X (G \times x)}_{n\,\text{times}}
\]
for $n \geq 1$ (induced by the fact that $\Gbb \times \Xbb$ is a model of an essentially algebraic theory in $\Var_{/\Spec K}$) gives a presentation of the stack $[G \backslash X]$ by $K$-varieties.  Furthermore, that this simplicial presentation is $2$-coskeletal follows from the fact that the essentially algebraic theory of groupoids is itself $2$-coskeletal (in the sense that the only axiomatized coherence diagrams and morphisms all reside and involve diagrams only containing objects of types $G_0, G_1$, or $G_2$) gives us the structure we need on the simplicial presentation of the stack $[G \backslash X]$. We will also see below that this agrees with the simplicial presentation of the stack that we call $\underline{[G \backslash X]}_{\bullet}$ that is used in \cite{Behrend} and described below.

Let us now discuss the simplicial schematic presentations of the stack $[G \backslash X]$ as a simplicial algebraic space and how it compares with the simplical scheme $\underline{G \backslash X}_{\bullet}$. For this, we note from the general theory of fppf descent that for all $n \in \N$, there is a presentation $X \to [G \backslash X]$ as an {\'e}tale equivalence relation such that the algebraic space
\[
X_n = \underbrace{X \times_{[G \backslash X]} \cdots \times_{[G \backslash X]} X}_{n+1\, \text{terms}}
\]
is represented by a scheme (and hence we identify the algebraic space  $X_n$ with its representing scheme). In particular, it is well known (cf.\@ \cite{Illusie}, \cite{Champs}, and \cite{OlssonSpacesStacks}, for instance) that for all $n \in \N$ there is an isomorphism of schemes
\[
X_n \cong G^n \times X.
\]
\begin{definition}\index{Simplicial! Presentation of $[G \backslash X]$}\index[notation]{XmodGsimplicial@$\underline{[G \backslash X]}_{\bullet}$}
	Let $\underline{[G \backslash X]}_{\bullet}$ be the simplicial scheme given by the simplicial presentation of the stack $[G \backslash X]$. Explicitly:
	\[
	\underline{[G \backslash X]}_n = \underbrace{X \times_{[G \backslash X]} \cdots \times_{[G \backslash X]} X}_{n+1\, \text{terms}}
	\]
	for all $n \in \N$, the face maps $d_k^n:\underline{[G \backslash X]}_{n+1} \to \underline{[G \backslash X]}_{n}$ are the projection maps from the $k$-th coordinate, and the degeneracy maps $s_{t}^{n}:\underline{[G \backslash X]}_{n} \to \underline{[G \backslash X]}_{n+1}$ are
	\[
	s_{t}^{n} = \underbrace{\id_X \times_{[G \backslash X]} \cdots \times_{[G \backslash X]} \id_X}_{t\,\text{terms}} \times_{[G \backslash X]} \Delta \times_{[G \backslash X]} \underbrace{\id_{X} \times_{[G \backslash X]} \cdots \times_{[G \backslash X]} \id_X}_{n-t\,\text{terms}},
	\]
	when $0 \leq t \leq n$ and where $\Delta:X \to X \times_{{[}G \backslash X{]}} X$ is the diagonal map.
\end{definition}
We will now show that this construction provides a natural isomorphism of simplicial data between the simplicial presentation $\underline{[G \backslash X]}_{\bullet}$ and the simplicial scheme $\underline{G \backslash X}_{\bullet}$. To prove this, however, we will use some simplicial techniques. The idea of the proof is to use the fact that both $\underline{G\backslash X}_{\bullet}$ and $\underline{[G \backslash X]}_{\bullet}$ are $2$-coskeletal, so an isomorphism of simplicial schemes may be verified by proving the relevant simplicial diagrams are isomorphisms in degrees $0, 1,$ and $2$.

Let us now describe the tools we need for the proof of the isomorphism $\underline{[G \backslash X]}_{\bullet} \cong \underline{G \backslash X}_{\bullet}$. We will need the language of coskeletal simplicial objects. We follow the construction as described in \cite[Section V]{GoerssJardine} of coskeleta of simplicial objects in a category $\Cscr$. 

To describe what it means to be $n$-coskeletal we need to discuss the $n$-truncation functors and their right Kan extensions. Fix an $n \in \N$. There then is a full subcategory $\mathbf{\Delta}_{\leq n}$ of $\mathbf{\Delta}$\index[notation]{DeltaaTruncated@$\mathbf{\Delta}_{\leq n}$}\index{The Topologist's Simplex Category! Truncated} generated by saying that for any $k \in \N$, the object $[k] \in (\mathbf{\Delta}_{\leq n})_{0}$ if and only if $0 \leq k \leq n$. The inclusion $\mathbf{\Delta}_{\leq n} \to \mathbf{\Delta}$ induces the $n$-truncation functor functor defined via:
\begin{prooftree}
	\AxiomC{$\tr_n:[\mathbf{\Delta}^{\op},\Sch] \to [\mathbf{\Delta}_{\leq n}^{\op},\Sch]$}\doubleLine
	\UnaryInfC{$[\mathbf{\Delta}_{\leq n}^{\op} \hookrightarrow \mathbf{\Delta}^{\op}, \id_{\Sch}]:[\mathbf{\Delta}^{\op},\Sch] \to [\mathbf{\Delta}_{\leq n}^{\op},\Sch] $}
\end{prooftree}\index[notation]{Truncationn@$\tr_n$}
We then define the $n$-coskeleton functors as the right Kan extension\index{nCoSk@$n$-Coskeleton Functor}
\[
\cosk_n := \Ran_{\tr_n}\big(\id_{[\mathbf{\Delta}^{\op},\Sch]}\big):[\mathbf{\Delta}_{\leq n}^{\op},\Sch] \to [\mathbf{\Delta}^{\op},\Sch].
\]\index[notation]{cosk@$\cosk_n$}
Note that this exists by \cite[Theorem 3.9]{ConradCohomDesc}, as the category $\Sch_{/\Spec K}$ admits finite limits. Because $\cosk_n$ is a right Kan extension of the identity functor $\id_{[\mathbf{\Delta}^{\op},\Sch]}$ along $\tr_n$, it follows that we have an adjunction:
\[
\begin{tikzcd}
{[}\mathbf{\Delta}^{\op},\Sch{]} \ar[r, bend left = 30, ""{name = U}]{}{\tr_n} & {[}\mathbf{\Delta}_{\leq n}^{\op},\Sch{]} \ar[l, bend left = 30, ""{name = L}]{}{\cosk_n} \ar[from = U, to = L, symbol = \dashv]
\end{tikzcd}
\]
This adjunction allows us to determine when an object is $n$-coskeletal.
\begin{definition}[{\cite[Page 370, 371]{GoerssJardine}}]\index{Simplicial! Scheme! nCoskeletal@$n$-Coskeletal}
	A simplicial scheme $X_{\bullet}$ is said to be $n$-coskeletal for $n \in \N$ if the unit morphism
	\[
	X_{\bullet} \xrightarrow{\eta_{X_{\bullet}}} \cosk_n\big(\tr_n(X_{\bullet})\big)
	\]
	is an isomorphism.
\end{definition}
\begin{remark}
	A more concrete description of the $n$-truncation functors $\tr_n$ is as follows: Let $X_{\bullet}$ be a simplicial scheme. The $n$-trunction $\tr_nX_{\bullet}$ is given by restricting the simplicial scheme $X_{\bullet}$ to all information given in degrees $k \leq n$. Explicitly, $\tr_{\leq n}X_{\bullet}$ is the functor $\mathbf{\Delta}_{\leq n}^{\op} \to \Sch$ where $(\tr_nX_{\bullet})_k = X_k$ for all $0 \leq k \leq n$ and all face maps and degeneracy maps of $\tr_nX_{\bullet}$ are the same as those of $X_{\bullet}$, save for $s_{k}^{m}$ and $d_{\ell}^{s}$ are included in $\tr_nX_{\bullet}$ exactly when $0 \leq m, s < n.$ The trunction $\tr_nf_{\bullet}$ of a morphisms of simplicial schemes is easier to define; if we write $f_{\bullet} = (f_m)_{m \in \N}$ then $\tr_n f_{\bullet} = (f_{m} \; : \; 0 \leq m \leq n)$.
	
	For example, given the simplicial scheme $\underline{G \backslash X}_{\bullet}$ described above, $\tr_n(\underline{G \backslash X}_{\bullet})$ is given by taking only the terms $G^m \times X$ for $0 \leq m \leq n$ and taking face and degeneracy maps only in degrees between $0$ and $n$.
\end{remark}
\begin{remark}\label{Remark: Section 4: truncations and limits}
	Because $\tr_n$ is a left adjoint for all $n \in \N$, it preserves all colimits of simplicial schemes. It also follows that $\tr_n$ preserves all limits of simplicial schemes, and in fact does so pointwise. To see this, note that if $(X_{\bullet,i}; f_i)$ is a system of simplicial schemes with limit $X_{\bullet}$ in $[\mathbf{\Delta}^{\op},\Sch]$, then we can describe $X_{\bullet} = (X_{m})$
	\[
	X_m := \lim_{\substack{\longleftarrow \\ i \in I}} X_{m,i}
	\]
	for all $m \in \N$; the face and degeneracy maps induced via the universal properties of the relevant limits. Now upon $n$-truncating we find that $\tr_nX_{\bullet}$ may be described via
	\[
	(\tr_nX_{\bullet})_m = \lim_{\substack{\longleftarrow \\ i \in I}} X_{m,i}
	\]
	for all $0 \leq m \leq n$. This $n$-truncated simplicial scheme is a limit of the schemes $\tr_nX_{\bullet,i}$ in $[\mathbf{\Delta}_{\leq n}^{\op},\Sch]$. This implies in particular that if $X_{\bullet}$ is a limit of some simplicial schemes $X_i$ then $X_{\bullet}$ is $n$-coskeletal for all $n \in \N$, as right adjoints (and hence $\cosk_n$) preserve limits and give us the computation
	\begin{align*}
	\cosk_n(\tr_n(X_{\bullet}) &= \cosk_n\left(\tr_n\left(\lim_{\longleftarrow} X_{m,i}\right)\right) \cong \cosk_n\left(\lim_{\longleftarrow}\left(\tr_n(X_{m,i})\right)\right) \\
	&\cong \lim_{\longleftarrow} \cosk_n(\tr_n(X_{m,i})).
	\end{align*}
\end{remark}
\begin{lemma}\label{Lemma: Section 4: Simplicial presentation of quotient scheme is 2coskeletal}
	The simplicial scheme $\underline{G \backslash X}_{\bullet}$ is $2$-coskeletal.
\end{lemma}
\begin{proof}
We first show that for any $n \in \N$, $\underline{G \backslash X}_n$ is generated by limits of the objects and maps appearing in $\tr_2(\underline{G\backslash X}_{\bullet})$ (and similarly every face and degeneracy map is induced in some way from a limit map). We first note that this is trivially true for all $\underline{G \backslash X}_n$, $d_k^n$, and $s_t^n$ when $n \in \lbrace 0, 1 \rbrace$, $0 \leq k \leq n+1$, and $0 \leq t \leq n$. Furthermore, there is an isomorphism $G \times G \times X \cong (G \times X) \times_X (G \times X)$ fitting into the pullback diagram
\[
\xymatrix{
(G \times X) \times_X (G \times X) \pullbackcorner \ar[d]_{p_1} \ar[r]^-{p_2} & G \times X \ar[d]^{\alpha_X} \\
G \times X \ar[r]_-{\pi_2} & X
}
\]
and the corresponding maps give rise to $d_0^1 \cong p_1$ and $d_2^1 \cong p_2$. The map $d_1^1$, however, is not determined by pullback from degree $1$ truncation information.

Set $n \geq 3$. We then note that there is an isomorphism
\begin{align*}
\underline{G \backslash X}_n &= G^n \times X \cong \underbrace{(G \times X) \times_X \cdots \times_X (G \times X)}_{n\,\text{copies of}\,(G \times X)\,\text{in the pullback}} \\
&\cong \underbrace{\underline{G \backslash X}_1 \times_{\underline{G \backslash X}_0} \cdots \times_{\underline{G \backslash X}_0} \underline{G \backslash X}_1}_{n\,\text{copies of}\,\underline{G\backslash X}_1\,\text{in the pullback}}.
\end{align*}
while the outer face maps $d_0^n \cong p_1$ and $d_{n+1}^n \cong p_n$. In this case, however, the inner face maps all correspond to pulling back along a $d_1^1$ at the $(G \times X) \times_X (G \times X)$ level of the iterated pullback, and hence arise as limits of the degree $2$ truncation of $\underline{G \backslash X}_{\bullet}$. The construction of the degeneracy maps is trivial to verify, as they simply translate doing the diagonal insertion in a single component to insertion of the identity element of $G$ in the same component, which is also induced as a (higher) limit map of degree $2$ truncations of $\underline{G \backslash X}_{\bullet}$.
		
Now, because $\underline{G \backslash X}_{\bullet}$ is degreewise described by limits of the degree $2$ truncation of $\underline{G \backslash X}_{\bullet}$, the lemma follows from Remark \ref{Remark: Section 4: truncations and limits}.
\end{proof}
\begin{lemma}\label{Lemma: Section 4: Simplicial presentation of a stack is 2coskeletal}
	The simplicial scheme $\underline{[G \backslash X]}_{\bullet}$ is $2$-coskeletal.
\end{lemma}
\begin{proof}
	Because the algebraic space
	\[
	\Xfrak_n = \underbrace{X \times_{[G \backslash X]} \cdots \times_{[G \backslash X]} X}_{n+1\, \text{terms}}
	\]
	is represented by the scheme $G^n \times X$, we get that the degree $n$-component of simplicial scheme $\underline{[G \backslash X]}_{\bullet}$ is given by
	\begin{align*}
	\underline{[G \backslash X]}_n &= \underbrace{X \times_{[G \backslash X]} \cdots \times_{[G \backslash X]} X}_{n+1\, \text{terms}} \cong G^{n} \times X.
	\end{align*}
	Thus each component scheme of $\underline{[G \backslash X]}_{\bullet}$ is a limit for all $n \in \N$. From here the lemma follows from Remark \ref{Remark: Section 4: truncations and limits} as well.
\end{proof}
\begin{proposition}\label{Prop: Section 4: Two simplicial presentations are isomorphic}
	There is an isomorphism of simplicial schemes
	\[
	\underline{[G \backslash X]}_{\bullet} \cong \underline{G \backslash X}_{\bullet}.
	\]
\end{proposition}
\begin{proof}
	Because $\underline{G \backslash X}_{\bullet}$ and $\underline{G \backslash X}_{\bullet}$ are both $2$-coskeletal by Lemmas \ref{Lemma: Section 4: Simplicial presentation of quotient scheme is 2coskeletal} \ref{Lemma: Section 4: Simplicial presentation of a stack is 2coskeletal}, it suffices to show that $\tr_2\underline{[G \backslash X]}_{\bullet} \cong \tr_2\underline{G \backslash X}_{\bullet}$. We note that $\underline{[G \backslash X]}_0 = X_0 = X = \underline{G \backslash X}_0$ while $\underline{[G \backslash X]}_1 = X_1 = X \times_{[G \backslash X]} X \cong G \times X = \underline{G \backslash X}_1$. Notice now that by construction we have
	\begin{prooftree}
		\AxiomC{$s_0^0:\underline{[G \backslash X]}_0  \to \underline{[G \backslash X]}_{1}$}
		\UnaryInfC{$\Delta:X \to X \times_{[G \backslash X]} X$}
	\end{prooftree}
	Under the isomorphism $X \times_{[G \backslash X]} X \xrightarrow{\cong} G \times X$ we have that $\Delta$ corresponds to the morphism $\langle 1_G, \id_X\rangle:X \to G \times X$, which proves that the diagram
	\[
	\xymatrix{
		X \ar@{=}[r] \ar[d]_{s_0} & X \ar[d]^{\langle 1_G, \id_X\rangle} \\
		X \times_{[G \backslash X]} X \ar[r]_-{\cong} & G \times X
	}
	\]
	commutes. It follows similarly that the projection maps
	\[
	p_1, p_2:X \times_{[G \backslash X]} X \to X
	\]
	(which are in turn the morphisms $d_0^0, d_1^0:\underline{[G \backslash X]}_1 \to \underline{[G \backslash X]}_0$ in $\underline{[G \backslash X]}_{\bullet}$ respectively) correspond to the morphisms $\alpha_X,\pi_2:G \times X \to X$ and make the diagrams
	\[
	\xymatrix{
		X \times_{[G \backslash X]} X \ar[r]^-{\cong} \ar[d]_{p_1} & G \times X \ar[d]^{\alpha_X} & X \times_{[G \backslash X]} X \ar[d]_{p_2} \ar[r]^-{\cong} & G \times X \ar[d]^{\pi_2} \\
		X \ar@{=}[r] & X & X \ar@{=}[r]& X
	}
	\]
	commute as well.
	
	We now prove that the isomorphisms are compatible with the $d_k^1$ and $s_{\ell}^1$ morphisms. For this note that in the case of the $s_{\ell}^{1}$ morphisms, it suffices without loss of generality to prove this for $s_0^1$, as the $s_1^1$ case follows mutatis mutandis. In the $s_0^1$ case we have that $s_0^1$ is equal, by construction to the map:
	\begin{prooftree}
		\AxiomC{$s_0^1:[G \backslash X]_1 \to [G \backslash X]_2$} \\
		\UnaryInfC{$\Delta \times_{[G \backslash X]} \id_X:X \times_{[G \backslash X]} X \to X \times_{[G \backslash X]} X \times_{[G \backslash X]} X$}
	\end{prooftree}
	Under the isomorphisms $[G \backslash X]_1 \cong G \times X$ and  $[G \backslash X]_2 \cong G \times G \times X$, we find that this map corresponds to the insertion of the identity in the left-most coordinate (as we embed the $G \times X$ component along the neutral component of the left-most copy of $G$ when taking the diagonal along the left):
	\[
	\langle 1_G, \id_{G \times X}\rangle: G \times X \to G \times G \times X.
	\]
	In particular, this implies that the diagram
	\[
	\xymatrix{
		X \times_{[G \backslash X]} X \ar[r]^-{\cong} \ar[d]_{s_0^1} &  G \times X \ar[d]^{\langle 1_G, \id_{G \times X}\rangle} \\
		X \times_{[G \backslash X]} X \times_{[G \backslash X]} X \ar[r]_-{\cong} & G \times G \times X
	}
	\]
	commutes, as was desired. Following this mutatis mutandis with $s_1^1$ shows that the diagram
	\[
	\xymatrix{
		X \times_{[G \backslash X]} X \ar[r]^-{\cong} \ar[d]_{s_1^1} &  G \times X \ar[d]^{\langle \id_{G}, 1_G, \id_X\rangle} \\
		X \times_{[G \backslash X]} X \times_{[G \backslash X]} X \ar[r]_-{\cong} & G \times G \times X
	}
	\]
	commutes, which is our final degeneracy map to verify.
	
	We now verify the face maps. First, note that $d_0^1$ and $d_2^1$ follow, by the same arguments as in the $m = 0$ case, to be the maps which project away the left-most coordinate of $G$ and take the action of the right-most copy of $G$ with $X$, respectively. In particular, it is routine but tedious to verify that the diagrams
	\[
	\xymatrix{
		X \times_{[G \backslash X]} X \times_{[G \backslash X]} X \ar[r]_-{\cong} \ar[d]_{d_0^1} & G \times G \times X \ar[d]^{\pi_{23}}  \\
		X \times_{[G \backslash X]} X \ar[r]_-{\cong} & G \times X
	}	
\]
\[
\xymatrix{
 X \times_{[G \backslash X]} X \times_{[G \backslash X]} X \ar[r]_-{\cong} \ar[d]_{d_2^1} & G \times G \times X \ar[d]^-{\id_G \times \alpha_X} \\
 X \times_{[G \backslash X]} X \ar[r]_-{\cong} & G \times X
}
\]
	commute. We instead focus on proving that the morphism $d_1^1:[G \backslash X]_2 \to [G \backslash X]_1$ corresponds to the multiplication map $\mu \times \id_X:G \times G \times X \to G \times X$. For this note that the triple pullback $X \times_{[G \backslash X]} X \times_{[G \backslash X]} X$ is the limit in the diagram:
	\[
	\begin{tikzcd}
	& X \times_{{[}G \backslash X{]}} X \times_{{[}G \backslash X{]}} X \ar[dr]{}{{\pi}_{23}} \ar[d]{}[description]{{\pi}_{13}} \ar[dl, swap]{}{{\pi}_{12}} & \\
	X \times_{{[}G \backslash X{]}} X \ar[d, swap]{}{{\pi}_1}  & X \times_{{[}G \backslash X{]}} X \ar[dr, near end]{}[description]{{\pi}_2} \ar[dl, swap, near end]{}[description]{{\pi}_1} & X \times_{{[}G \backslash X{]}} X  \ar[d]{}{{\pi}_2} \\
	X \ar[dr] & X \ar[d] & X \ar[dl] \\
	& {[}G \backslash X{]} \ar[from = 2-3, to = 3-2, crossing over, near end]{}[description]{{\pi}_1} \ar[from = 2-1, to = 3-2, crossing over, near end]{}[description]{{\pi}_{2}}
	\end{tikzcd}
	\]
	As such, the map $d_1^1$, which corresponds to the project $\pi_{13}$, leaves the $X$-component alone but intermingles the $G$-components according to the natural $G$-action. However, this simply says that the diagram
	\[
	\xymatrix{
		X \times_{{[}G \backslash X{]}} X \times_{{[}G \backslash X{]}} \ar[r]^-{\cong} \ar[d]_{d_1^1} & G \times G \times X \ar[d]^{\mu \times \id_X} \\
		X \times_{{[}G \backslash X{]}} X \ar[r]_-{\cong} & G \times X
	}
	\]
	commutes, which is what we had to verify. This in turn proves that
	\[
	\tr_2\underline{[G \backslash X]}_{\bullet} \cong \tr_2\underline{G \backslash X}_{\bullet},
	\]
	which completes the proof of the proposition upon consideration of the isomorphisms
	\[
	\underline{[G \backslash X]}_{\bullet} \cong \cosk_2(\tr_2\underline{[G \backslash X]}_{\bullet}) \cong \cosk_2(\tr_2\underline{G \backslash X}_{\bullet}) \cong \underline{G \backslash X}_{\bullet}.
	\]
\end{proof}
\begin{corollary}\label{Cor: Section 4: Bounded derived simplicial cat is equiv}
	There is an isomorphism of bounded derived categories of simplicial sheaves
	\[
	D^{b}_{\text{eq}}(\underline{[G \backslash X]}_{\bullet};\overline{\Q}_{\ell}) \cong D^{b}_{\text{eq}}(\underline{G \backslash X}_{\bullet};\overline{\Q}_{\ell}).
	\]
\end{corollary}
As we move to show that $D^{b}_{\text{eq}}(\underline{G \backslash X}_{\bullet};\overline{\Q}_{\ell}) \simeq D^{b}_c([G \backslash X];\overline{\Q}_{\ell})$, we will need to briefly discuss the notions of fibred toposes and fibred sites as they appear in \cite[Expos{\'e} VI]{SGA4} so that we may explicitly describe the category $D_c^b([G \backslash X];\overline{\Q}_{\ell})$ of \cite{Behrend}. Because we need to talk about fibre functors within fibrations, we make the convention that all fibrations $p:\Ascr \to \Cscr$ are equipped with a chosen cleavage, as we will need to discuss pullback functors between categories of fibres.
\begin{remark}\index[notation]{Ascrrfpullback@$\Ascr(f)^{\ast}$}\index[notation]{Ascrobject@$\Ascr(U)$}
	As a point of notation: When given a (cloven) fibration $p:\Ascr \to \Cscr$, for any $U \in \Cscr_0$ and $f \in \Cscr_1$, we will write $\Ascr(U)$ to denote the category of fibres of $\Ascr$ above $U$ and $\Ascr(f)^{\ast}$ to denote the chosen pullback functor $\Ascr(f)^{\ast}:\Ascr(\Codom f) \to \Ascr(\Dom f)$.
\end{remark}
\begin{definition}[\cite{SGA4}]\index{Fibred (Grothendieck) Topos}
	A fibred (Grothendieck) topos over a category $\Cscr$ is a fibration $p:\Ecal \to \Cscr$ such that for each $U \in \Cscr_0$ and $f \in \Cscr_1$, $\Ecal(U)$ is a Grothendieck topos and $\Ecal(f)^{\ast}$ is the pullback of a geoemetric morphism $\Ecal(f):\Ecal(\Dom(f)) \to \Ecal(\Codom f)$.
\end{definition}
We will also need the notion of a total topos associated to a fibred topos over a category $\Cscr$; this comes from  \cite[Remarque VI.7.4.3]{SGA4}. Essentially the idea here is to take a fibration $p:\Ascr \to \Cscr$ whose corresponding pseudofunctor factors through the $2$-category of sites and then form a global version of the local sites (the sites on the fibre categories of $p$) called the total site of $\Ascr$. This total site associates a total topos of the fibration. Our goal is to take a fibred topos $q:\Ecal \to \Cscr$ and produce a fibred site over $\Cscr$, $p:\Ascr \to \Cscr$, in such a way that for all $U \in \Cscr$, $\Ecal(U) \simeq \Shv(\Ascr(U),J_U)$. If we can do this we'll also be able to realize the topos on the total site of $\Ascr$ as a global version of the sheaves which locally come from the sheaves on $\Ecal(U)$.
\begin{remark}\index{Morphism of Sites}\index{Site Morphism|see {Morphism of Sites}}
	Recall that if $(\Cscr,J)$ and $(\Dscr,K)$ are sites, then a functor $f:\Dscr \to \Cscr$ is a {morphism of sites $f:(\Cscr,J) \to (\Dscr,K)$} (in the style of \@ \cite[IV.4.9.1]{SGA4}) if $f$ is continuous, $f^{\ast}:\Shv(\Cscr,J) \to \Shv(\Dscr,K)$ is right exact, and there is a geometric morphism $g:\Shv(\Cscr,J) \to \Shv(\Dscr,K)$ for which the diagram
	\[
	\xymatrix{
		\Cscr \ar[d]_{\yon_J^{++}} \ar[r]^-{f} & \Dscr \ar[d]^{\yon_K^{++}} \\
		\Shv(\Cscr,J) \ar[r]_-{g^{\ast}} & \Shv(\Dscr,K)
	}
	\]
	commutes in $\Cat$, where $\yon^{++}$ is the sheafification of the Yoneda embedding in the relevant topology, i.e., $\yon_J^{++} := (-)_J^{++} \circ \yon_{\Cscr}$.
\end{remark}
\begin{definition}[\cite{SGA4}]\index{Fibred Site}
	A fibred site $\Ascr$ over a category $\Cscr$ is a cloven fibration $p:\Ascr \to \Cscr$ such that for all $U \in \Cscr_0$ the category $\Ascr(U)$ is a site $(\Ascr(U),J_U)$ and the morphism $\Ascr(f)^{\ast}:\Ascr(\Codom f) \to \Ascr(\Dom f)$ is a morphism of sites from $\Ascr(\Dom f)$ to $\Ascr(\Codom f)$.
\end{definition}
\begin{definition}[{\cite[Definition VI.7.4.1]{SGA4}}]\index{Fibred Site! Total Site}
	Let $p:\Ascr \to \Cscr$ be a fibred site and for each $U \in \Cscr_0$, let $a_!^U:\Ascr(U) \to \Ascr$ denote the inclusion of the category of $U$-fibres into $\Ascr$. We then define the total topology for $\Ascr$, denoted $\tot$, to be the finiest Grothendieck topology $J$ on $\Ascr$ for which all the functors $a_!^U:\Ascr(U) \to \Ascr$ are site-continuous. We define the total site of a fibred site $p:\Ascr \to \Cscr$ to be the pair $(\Ascr,\tot)$.\index[notation]{$\tot$}
\end{definition}
\begin{remark}
	An alternative, but equivalent, definition for fibred \\Grothendieck toposes and fibred sites over $\Cscr$ is to give pseudofunctors $S,E:\Cscr^{\op} \to \fCat$ for which:
	\begin{itemize}
		\item $S(U)$ is a site for all $U \in \Cscr_0$ and the pullback functor $S(f)$ is a site-continuous functor for all $f \in \Cscr_1$;
		\item $E(U)$ is a Grothendieck topos for all $U \in \Cscr_0$ and the pullback functor $E(f)$ is the inverse image of a geometric morphism $E(\Dom f) \to E(\Codom f)$ for all $f \in \Cscr_1$.
	\end{itemize}
\end{remark}
\begin{definition}\index{Fibred Site! Total Topos}\index{Total Topos|see {Fibred Site}}
	The total topos associated to a fibred site $p:\Ascr \to \Cscr$ is the topos of $\tot$-sheaves on $\Ascr$, i.e., $\Shv(\Ascr,\tot)$. We denote this topos by $\tottop(\Ascr)$.\index[notation]{$\tottop(\Ascr)$}
\end{definition}
\begin{remark}\label{Remark: Section 4: Total topos of a fibred topos}
	Let $p:\Ascr \to \Cscr$ be a fibred site. By construction, because all the morphisms $a_!^U:\Ascr(U) \to \Ascr$ are site-cocontinuous morphisms of sites $(\Ascr,\tot) \to (\Ascr(U),J_U)$ for all $U \in \Cscr_0$, there are induced geometric morphisms $\Shv(\Ascr(U),J_U) \to \tottop(\Ascr)$
	\[
	\begin{tikzcd}
	\tottop(\Ascr) \ar[rr, bend left =30, ""{name = U}]{}{(a_!^{U})^{\ast}} & & \Shv(\Ascr(U),J_U) \ar[ll, bend left = 30, ""{name = L}]{}{(a_!^U)_{\ast}} \ar[from = U, to = L, symbol = \dashv]
	\end{tikzcd}
	\]
	for all $U \in \Cscr_0$. Following \cite{Behrend}, we see from the fact that $\tot$ is the finest topology on $\Ascr$ for which all the maps $a_!^U$ induce the geometric morphisms above allows us to give the following equivalent description of the topos $\tottop(\Ascr)$:
	\begin{itemize}
		\item Objects: Pairs $(\Fscr,\Theta)$ where
		\[
		\Fscr := \lbrace \Fscr_U \; | \; U \in \Cscr_0 \rbrace
		\]
		and
		\[
		\Theta := \lbrace \theta_f \in \Shv(\Ascr(U),J_{U})(\Ascr(f)\Fscr_V, \Fscr_U) \; | \; f \in \Cscr_1, f:U \to V \rbrace
		\]
		subject to the cocycle condition 
		\[
		\theta_{g \circ f} = \theta_f \circ \Ascr(f)^{\ast}\theta_g.
		\]
		\item Morphism: A morphism $\alpha:(\Fscr,\Theta_f) \to (\Gscr,\Sigma_f)$ is a collection
		\[
		\alpha := \lbrace \alpha_U: \Fscr_U \to \Gscr_U \; | \; U \in \Cscr_0\rbrace
		\]
		such that for all $f \in \Cscr_1$ the diagram
		\[
		\xymatrix{
			\Ascr(f)^{\ast}\Fscr_{\Codom f} \ar[rr]^-{\Ascr(f)^{\ast} \alpha_{\Codom f}} \ar[d]_{\theta_{f}} & & \Ascr(f)^{\ast}\Gscr_{\Codom f} \ar[d]^{\sigma_f} \\
			\Fscr_{\Dom f} \ar[rr]_-{\alpha_{\Dom f}} & & \Gscr_{\Dom f}
		}
		\]
		commutes in $\Shv(\Ascr(\Dom f), J_{\Dom f})$.
		\item Composition: Pointwise.
		\item Identities: Pointwise identity morphisms.
	\end{itemize}
	This construnction works just as well if we simply have a fibred topos $p:\Ecal \to \Cscr$, even if we do not know that every geometric morphism in $\Ecal(U) \to \Ecal(V)$ is induced by site-continuous functors $\Ascr(U) \to \Ascr(V)$. In this way the total topos of \cite{Behrend} is a different notion than that of \cite[Remarque VI.7.4.3]{SGA4}; however, each of the topologies considered in \cite{Behrend} does have this property, so no harm is done. However, when given a general fibred topos $p:\Ecal \to  \Cscr$, we can define the total topos of $\Ecal$, $\tottop(\Ecal)$, to be the category generated by the $\Ecal$-analogue of the description of $\tottop(\Ascr)$ above.
\end{remark}
\begin{lemma}[\cite{Behrend}]\label{Lemma: Section 4: Cartesian objects of Total topos}
	Let $p:\Ecal \to \Cscr$ be a fibred topos, let $p_{\cart}:\Ecal_{\cart} \to \Cscr$\index[notation]{ECart@$\Ecal_{\cart}$} be the associated Cartesian fibration, let $\Scal_{\cart}$\index[notation]{SCart@$\Scal_{\cart}$} be the category of sections of $p_{\cart}$, and let $\tottop(\Ecal)$ be its total topos (cf.\@ the end of Remark \ref{Remark: Section 4: Total topos of a fibred topos} above). Then the category $\Scal_{\cart}$ is equivalent to the subcategory of $\tottop(\Ecal)$ generated by objects $(\Fscr,\Theta)$ whose transition morphisms $\theta_f:\Ecal(p)^{\ast}\Fscr_{\Codom f} \to \Fscr_{\Dom f}$ are all isomorphisms. Moreover $\Scal_{\cart}$ is a topos and there is a geometric morphism $\pi:\tottop(\Ecal) \to \Scal_{\cart}$ for which $\pi_{\ast}$ is a right-adjoint-right-inverse for $\pi^{\ast}$.
\end{lemma}
\begin{proof}
	The verification that $\Scal_{\cart}$ is equivalent to the subcategory described above of the above form comes from the fact that since every morphism in $\Ecal_{\cart}$ is Cartesian, each fibre category of $\Ecal_{\cart}$ is a groupoid. As such, upon taking sections of $p_{\cart}$ we find that the transition morphisms in $\Theta$ must all be isomorphisms.                                                                                                                                                                                                                                                
	
	We now describe the right-adjoint-right-inverse statement at the end of the lemma. For this define $\Fcal$ to be the subcategory of $\tottop(\Ecal)$ described above with $\Fcal \simeq \Scal_{\cart}$ and let $\pi^{\ast}:\Fcal \to \tottop(\Ecal)$ be the inclusion functor. Then we see that the limits and colimits in $\Fcal$ are given as limits and colimits in $\tottop(\Ecal)$, so $\pi^{\ast}$ creates and reflects all limits and colimits in $\Fcal$. This implies that $\pi^{\ast}$ admits a left and right adjoint by the Adjoint Functor Theorem, as $\Fcal$ and $\tottop(\Ecal)$ are both Grothendieck toposes. Finally, a routine calculation shows that $\pi_{\ast} \circ \pi^{\ast} = \id_{\Fcal}$.
\end{proof}

We now need one more concept involving Grothendieck toposes from \cite{Behrend}: That of a $c$-structure on a topos $\Ecal$.
\begin{definition}[{\cite[Definition 4.4.2]{Behrend}}]\index{c-Topos}\index{c-Structure|see {c-Topos}}
	Let $\Ecal$ be a Grothendieck topos. A c-structure on $\Ecal$ is given by a full subcategory $\overline{\Ecal}$ of $\Ecal$ such that:
	\begin{itemize}
		\item The category $\overline{\Ecal}$ is a Grothendieck topos;
		\item The inclusion functor $\pi^{\ast}:\overline{\Ecal} \to \Ecal$ is exact and has a right adjoint right inverse $\pi_{\ast}:\Ecal \to \overline{\Ecal}$;
		\item If $G$ is a group object in $\Ecal$ and $p:P \to B$ is a principal $G$-bundle in $\Ecal$, then if $B \in \overline{\Ecal}_0$ and $G \in \overline{\Ecal}_0$, $P \in \overline{\Ecal}_0$ as well.
	\end{itemize}
In such a case we call a topos $\Ecal$ with c-strucutre $\overline{\Ecal}$ a c-topos $(\Ecal, \overline{\Ecal})$.\index[notation]{EcalEcaloverline@$(\Ecal,\overline{\Ecal})$}
\end{definition}
\begin{remark}
	Recall that a principal $G$-bundle $p:P \to B$ in a topos $\Ecal$ is given by a left $G$-object $P$ and a morphism $p:P \to B$ such that the diagram
	\[
	\xymatrix{
		G \times P \ar[r]^-{\alpha_P} \ar[d]_{\pi_2} & P \ar[d]^{p} \\
		P \ar[r]_{p} & B
	}
	\]
	commutes and there exists a regular epimorphism $\gamma:U \to B$ (a cover in $\Ecal$ of $B$) such that in the pullback diagram
	\[
	\xymatrix{
		U \times_B P \pullbackcorner \ar[r] \ar[d] & P \ar[d]^{p} \\
		U \ar[r]_-{\gamma} & B
	}
	\]
	we have an isomorphism
	\[
	U \times_B P \cong G \times U.
	\]
\end{remark}
\begin{remark}[\cite{Behrend}]\label{Remark: Behrend c topos facts}
	Here are some general structural remarks about toposes with c-structures that are relevant for our application:
	\begin{itemize}
		\item A topos with c-strucutre $(\Ecal,\overline{\Ecal})$ is said to be Noetherian whenever $\overline{\Ecal}$ is Noetherian (cf.\@ Remark \ref{Remark: Noetherian Topos} for what we need from Noetherian toposes).
		\item A fibred c-topos $p:\Ecal \to \Cscr$ is a fibred topos $\Ecal$ such that there exists a fibred topos $\overline{p}:\overline{\Ecal} \to \Cscr$ with fully faithful functor $i:\overline{\Ecal} \to \Ecal$ making
		\[
		\xymatrix{
			\overline{\Ecal} \ar[rr]^-{i} \ar[dr]_{\overline{p}} & & \Ecal \ar[dl]^{p} \\
			& \Cscr &
		}
		\]
		commute such that for all $U \in \Cscr_0$, $(\Ecal(U),\overline{\Ecal}(U))$ is a c-topos.
		\item For any fibred topos $p:\Ecal \to \Cscr$, the pair $(\tottop(\Ecal),\Scal_{\cart})$ is a c-topos (cf.\@ \cite[Example 4.4.5]{Behrend}).
		\item If $(p:\Ecal \to \Cscr,\overline{p}:\overline{\Ecal} \to \Cscr)$ is a fibred c-topos, then $(\tottop(\Ecal),\tottop(\overline{\Ecal}))$ is a c-topos (cf.\@ \cite[Note 4.4.7]{Behrend}).
	\end{itemize}
\end{remark}

We now give the crucial construction which will allow us to describe stratifications of Grothendieck toposes. For this we follow \cite[Exercise IV.7.8.a]{SGA4}; note that background on open, closed, and locally closed subtoposes of a Grothendieck topos may be found in Appendix \ref{Appendix: Locally closed subtoposes}. Fix a Grothendieck topos $\Ecal$ and write
\[
\lvert \Ecal \rvert := {\Point(\Ecal)}_{/\text{iso}},
\]
i.e., the set of points of $\Ecal$ modulo isomorphism. Now fix an open subtopos of $\Ecal$ and sheaf representation $\Ecal \simeq \Shv(\Cscr,J)$ and $\Ucal \simeq \Shv(\Cscr,J)_{/\Fscr}$. Define now the set
\[
\lvert \Ucal \rvert := \lbrace [p] \in \lvert \Ecal \rvert \; : \; \Ucal_p \ne \emptyset \rbrace,
\]
i.e., the set of points $p$ of $\Ecal$ which do not send the open subtopos $\Ucal$ to the empty set. We then determine a topology on $\lvert \Ecal \rvert$ by setting
\[
\Open(\lvert \Ecal \rvert) := \lbrace \lvert \Ucal \rvert \; : \; \Ucal\, \text{is\,an\,open\,subtopos\,of\,} \Ecal\rbrace.
\]
This space has the property that it induces an isomorphism of posets
\[
\Open(\Ecal) \cong \Open(\lvert \Ecal \rvert),
\]
where the poset on the left is the poset of open subobject lattice of $\Ecal$, i.e.,
\[
\Open(\Ecal) = \mathbf{Sub}_{\text{Open}}(\Ecal).
\]\index[notation]{EcalAssociatedSpace@$\lvert\Ecal\rvert$}
\begin{definition}[{\cite[Exercise IV.7.8.a]{SGA4}}]\index{Topos! Associated Topological Space}
	Let $\Ecal$ be a Grothendieck topos. We then define the topological space associated to $\Ecal$ to be the space $(\lvert \Ecal \rvert, \Open(\lvert \Ecal \rvert))$
	constructed above.
\end{definition}
\begin{remark}\label{Remark: Stratifications of a topos}\index{Stratification! Of a Grothendieck Topos}
	A stratification of a Grothendieck topos $\Ecal$ is a finite collection $S$ of locally closed subtoposes $\Vcal$ of $\Ecal$ for which
	\[
	\Ecal \simeq \coprod_{\Vcal \in S} \Vcal
	\]
	and such that the closure of each $\Vcal$ is the union (pushout) of strata. Alternatively, we can describe a stratification of $\Ecal$ by finding a stratification of $\lvert \Ecal \rvert$ and then to each locally closed subset $\lvert \Vcal \rvert$ finding the corresponding locally closed subtopos. For details on locally closed subtoposes (such as definitions, constructions, and basic properties), we defer the reader to Appendix \ref{Appendix: Locally closed subtoposes} below.
\end{remark}

We now need a short digression into the theory of pre-L-stratifications of Grothendieck toposes discussed in \cite{Behrend}, as it will allow us to finally discuss stratifications (and hence constructible) of toposes with c-structures. The presentation of these notions here are quite brief, but for further details see Appendix \ref{Appendix: Stratifications and c-structures}. We now simply recall the definition of pre-L-stratifications (cf.\@ Definition \ref{Defn: Pre-L-stratification}) as given in \cite{Behrend}, save with different notational conventions. Fix a Noetherian Grothendieck topos $\Ecal$ and let $\mathsf{LC}(\Ecal)$ be the category of locally closed subtoposes of $\Ecal$, with a morphism $\Vcal \to \Wcal$ in $\mathsf{LC}(\Ecal)$ given by an inclusion of subtoposes.  A pre-L-stratification of $\Ecal$ is then a stratification $\Scal$ of $\Ecal$ toegether with a function $\Lcal$ which assigns to each strata $\Vcal \in \Scal$ a finite set of simple objects of (the hearts of the lcc-components of) the bounded derived fibration
\[
p:D^b \to \mathsf{LC}(\Ecal).
\]
For more details and justifications on this construction see \cite[Section 4.2, 4.3]{Behrend} for proofs and the full generality of the theory, or see Appendix \ref{Appendix: Stratifications and c-structures} for what we need and use in the scope of this paper. The fibration $p$ is defined by taking $D^b$ to be the category with
\[
D^b(\Vcal) := D^{b}(\Pro{\Vcal},\Z_{\ell})
\]
and cleaving the fibration via choosing pullback functors appropriately. Note that $\Pro{\Vcal}$ is the topos of projective sequences of $\Vcal$-sheaves and $D^{b}(\Pro{\Vcal},\Z_{\ell})$ is the bounded derived category of $\Z_{\ell}$-modules in $\Pro{\Vcal}$. The lcc-structure is given by saying that a sheaf $\Fscr = (\Fscr_n)$ is lcc each sheaf $\Fscr_n$ is locally isomorphic to a constant sheaf $\Delta_{X_n}$ for $X_n \in \FinSet_0$. The cd-structure (cf.\@ Definition \ref{Defn: Appendix Stratification: cd-structure}) is induced by defining the categories
\[
D^b_{\text{lcc}}(\Vcal) := D^b_{\text{lcc}}(\Pro{\Vcal},\Z_{\ell})
\]
for all strata $\Vcal$. This gives us the categories of constructible and $(\Scal,\Lcal)$-constructible $\ell$-adic sheaves on $\Ecal$ by setting
\[
D^{b}_{c}(\Pro{\Ecal};\overline{\Q}_{\ell}) =: \Dbb^{b}_{c}(\Ecal;\overline{\Q}_{\ell}),
\]\index[notation]{DboldEcal@$\Dbb^b_c(\Ecal;\overline{\Q}_{\ell})$}
\[ 
\qquad D^{b}_{(\Scal,\Lcal)}(\Pro{\Ecal};\overline{\Q}_{\ell}) =: \Dbb_{(\Scal,\Lcal)}^{b}(\Ecal;\overline{\Q}_{\ell}).
\]\index[notation]{DboldEcalSL@$\Dbb^b_{(\Scal,\Lcal)}(\Pro{\Ecal};\overline{\Q}_{\ell})$}
%

We now move to discuss stratifications of Grothendieck toposes with c-structures and how to produce the corresponding $\ell$-adic categories. For this we fix a c-structure $(\Ecal,\overline{\Ecal})$.
\begin{definition}[{\cite[Definition 4.4.9]{Behrend}}]\index{Pre-L-Stratification! Of a c-Topos}
	Fix a c-topos $(\Ecal,\overline{\Ecal})$. A (pre-L)-stratification of $(\Ecal,\overline{\Ecal})$ is a (pre-L)-stratification of $\overline{\Ecal}$.
\end{definition}
\begin{definition}[{\cite[Definition 4.4.9]{Behrend}}]\index{c-Topos! Constructible Sheaf}
	Fix a c-topos $(\Ecal,\overline{\Ecal})$. A sheaf $\Fscr \in \Ecal_0$ is constructible with respect to the c-structure if and only if $\Fscr \in \overline{\Ecal}_0$ and $\Fscr$ is a constructible sheaf in $\overline{\Ecal}$.
\end{definition}
\begin{remark}
	The terminology of calling an object of a Grothendieck topos $A \in \Ecal_0$ a sheaf $\Fscr$ of $\Ecal$ is simply an abuse of notation/terminology. We really mean that a constructible object is an object $A \in \Ecal_0$ such that in the sheaf representation $\Ecal \simeq \Shv(\Cscr,J); F_{\ast}:\Ecal \to \Shv(\Cscr,J)$, the corresponding sheaf $\Fscr = F_{\ast}A$ is constructible in $\Shv(\Cscr,J)$, but it is simply easier to say ``a sheaf in $\Ecal$'' instead.
\end{remark}

This allows us to build a derived category of $\ell$-adic sheaves over a c-topos $(\Ecal,\overline{\Ecal})$ when the topos $\overline{\Ecal}$ is Noetherian. In this case the technical constructions performed in \cite{Behrend} take the category of constructible sheaves in a c-topos $(\Ecal,\overline{\Ecal})$ allows us to perform constructions in such a way that we produce a category of constructible sheaves in the c-topos $(\Ecal,\overline{\Ecal})$ which is Artin-Rees $\ell$-adic with coefficients in our desired topos $\overline{\Ecal}$. 
\begin{definition}\index{c-Topos! Constructible Sheaves}
We define the category $\Dbb^{b}_{c}((\Ecal,\overline{\Ecal});\overline{\Q}_{\ell})$\index[notation]{DboldEcalEcal@$\Dbb_c^b\big((\Ecal,\overline{\Ecal});\overline{\Q}_{\ell}\big)$} to be the full subcategory of $\Dbb_c^b(\Ecal;\overline{\Q}_{\ell})$ where the complexes $A$ take values in $\overline{\Ecal}$ and are constructible as $\overline{\Ecal}$ objects.
\end{definition}

To construct the bounded derived category $\Dbb_c^b((\Ecal,\overline{\Ecal});\overline{\Q}_{\ell})$ we just need to take colimits over finer and finer stratifications. For this we will need to have the directed set $\Mcal$\index[notation]{Mcal@$\Mcal$} of stratifications of the c-topos $(\Ecal,\overline{\Ecal})$. Let $\Mcal$ be the collection of pre-L-stratifications of $(\Ecal,\overline{\Ecal})$ and note that these are themselves simply pre-L-stratifications of $\overline{\Ecal}$. By the fact that $\Mcal$ is isomorphic to the set $\Mcal^{\prime}$ of pre-L-stratifications of $\lvert \overline{\Ecal} \rvert$ and $\Mcal^{\prime}$ is a set by a routine size argument with the bound
\begin{align*}
&\sharp(\Mcal^{\prime}) \leq   \sharp(\Pcal(\lvert \Ecal \rvert))\cdot \sharp\left(\bigcup_{\substack{ \Scal \\ (\Scal, \Lcal) \in \Mcal^{\prime}}} \Set\left(\Scal,\lbrace Y \subseteq A_{\text{lcc}}^{b}(-)_0^{\text{simple}} \; | \; \sharp(Y) < \infty\rbrace\right)\right).
\end{align*}
The set $\Mcal$ is directed with respect to the partial order of Remark \ref{Remark: Appendix stratifications: Refinment order}. There is then a functor, for each pair of pre-L-stratifications $(\Scal,\Lcal), (\Tcal,\Jcal) \in \Mcal$ with $(\Scal,\Lcal) \leq (\Tcal,\Jcal)$,
\[
\Dbb^b_{(\Scal,\Lcal)}\big((\Ecal,\overline{\Ecal});\overline{\Q}_{\ell}\big) \to \Dbb^b_{(\Tcal,\Jcal)}\big((\Ecal,\overline{\Ecal});\overline{\Q}_{\ell}\big)
\]
which is given by inclusion; note that it is not necessarily full. 
\begin{definition}[\cite{Behrend}]\index{c-Topos! Derived Category of $\ell$-adic Sheaves}
The derived category of constructible $\ell$-adic sheaves on a c-topos is defined by
\[
\Dbb^b_c\big((\Ecal,\overline{\Ecal});\overline{\Q}_{\ell}\big) := \twoColim{(\Scal,\Lcal) \in \Mcal}\Dbb_{(\Scal,\Lcal)}^{b}\big((\Ecal,\overline{\Ecal});\overline{\Q}_{\ell}\big).
\]\index[notation]{DbbEcalEcalconstructible@$\Dbb^b_c\big((\Ecal,\overline{\Ecal});\overline{\Q}_{\ell}\big)$}
\end{definition}
In the case $\overline{\Ecal}$ is a topos of ({\'e}tale) sheaves on an algebraic stack $\Xfrak$, we write
\[
D_c^b(\Xfrak;\overline{\Q}_{\ell}) := \Dbb_c^b\big((\Ecal,\overline{\Ecal});\overline{\Q}_{\ell}\big).
\]\index[notation]{$D_c^b(\Xfrak;\overline{\Q}_{\ell})$}

We now need to discuss the explicit fibred topos that allows us to define and work with $\Dbb^b_c([G \backslash X];\overline{\Q}_{\ell})$. Before we jump in, however, let us motivate what we need. Because we are working with the simplicial scheme $\underline{[G \backslash X]}_{\bullet}$ as a simplicial presentation of $[G \backslash X]$, to describe an {\'e}tale c-topos on $[G \backslash X]$ we should be working with a fibration whose fibres look like the {\'e}tale toposes of $\underline{[G \backslash X]}_n$ and whose fibre functors are all pullbacks of the geometric morphisms
\[
d_k^n:\Shv(G^{n+1}\times X,\text{{\'e}t}) \to \Shv(G^n \times X,\text{{\'e}t}),
\]
\[
s_t^n:\Shv(G^{n}\times X,\text{{\'e}t}) \to \Shv(G^{n+1} \times X,\text{{\'e}t}),
\]
i.e., the fibre functors should be
\[
(d_k^n)^{\ast}:\Shv(G^n \times X,\text{{\'e}t}) \to \Shv(G^{n+1} \times X,\text{{\'e}t})
\]
and
\[
(s_t^n)^{\ast}:\Shv(G^{n+1}\times X,\text{{\'e}t}) \to \Shv(G^n \times X,\text{{\'e}t})
\]
for all $n \in \N$ with $0 \leq t \leq n, 0 \leq k \leq n+1$. Note that because the simplicial scheme $\underline{[G \backslash X]}_{\bullet}$ is already realized as a functor $\mathbf{\Delta}^{\op} \to \Sch$, and our desired construction of a fibred topos has morphisms going in the opposite direction of this, we need our fibration to be a pseudofunctor realized as a cosimplicial category $\Ecal:\mathbf{\Delta} \to \Cat$ for which each category $\Ecal([n])$ is a Grothendieck topos and each of the functors $\Ecal(\delta_k^n), \Ecal(\sigma_t^n)$ are pullbacks of geometric morphisms. this construction essentially forces the kind of fibrations we consider to define $\Dbb^b_c(\Xfrak;\overline{\Q}_{\ell})$: cosimplicial categories realized as functors $\mathbf{\Delta} \to \Cat$ which take values in Grothendieck toposes.
\begin{remark}
	In the formalism we have defined, it is possible to realize a cosimplicial category taking values as simultaneously a functor $\Ecal:\mathbf{\Delta} \to \Cat$ as we described above or, using the fact that each fibre functor $\Ecal(\delta_k^n)$ and $\Ecal(\sigma_t^n)$ admits a right adjoint, as a simplicial functor $\mathbf{\Delta}^{\op} \to \mathbf{GrTopos}$, where $\mathbf{GrTopos}$ is the $1$-category of Grothendieck toposes with morphisms geometric morphisms. Note, however, that doing this involves the Axiom of Choice: We have to choose right adjoints to each pullback functor, and while this gives an isomorphism between any two choices of translation, it is only unique up to a unique isomorphism. While the construction of our cosimplicial categories as simplicial toposes may be a more clean way of describing what how things are valued, because we have to take an extra step to throw away all of our pushforward functors when working with a simplicial topos, we have elected to use the cosimplicial language in this paper.
\end{remark}

Let $\Xfrak$\index[notation]{Xfrak@$\Xfrak$} be a finite type algebraic stack over $\Spec K$ and let $X \to \Xfrak$ be a presentation of the stack as a quotient of an {\'e}tale equivalence relation. Define the simplicial algebraic space $\Xfrak_{\bullet}$\index[notation]{XfrakBullet@$\Xfrak_{\bullet}$} by setting
\[
\Xfrak_{\bullet} := \underbrace{X \times_{\Xfrak} \cdots \times_{\Xfrak} X}_{n+1\, \text{terms}}
\]
and setting the face maps to the projections while setting the degeneracy maps to be insertion of the diagonals. 
\begin{definition}\label{Defn: Etale topos of simp alg space}\index{EtaleTopos@{\'E}tale Topos!Of a Simplicial Algebraic Space}
Let $X_{\bullet}$ be a simplicial algebraic space. Then the {\'e}tale topos on $X_{\bullet}$, $(X_{\bullet})_{\text{{\'e}t}}$\index[notation]{Xbulletet@$(X_{\bullet})_{\text{{\'e}t}}$} is the total topos of the fibration associated to the cosimplicial category $\Ecal:\mathbf{\Delta} \to \Cat$ defined by setting $\Ecal([n])$ to be
\[
\Ecal([n]) = \Ecal_n := \Shv(\Xfrak_n,\text{{\'e}t})
\]
and taking each pullback functor to be given by
\[
\Ecal(\delta_k^n) := (d_k^n)^{\ast}, \qquad \Ecal(\sigma_k^n) := (s_k^n)^{\ast}.
\]
\end{definition}
Taking the total topos of the fibration associated to the cosimplicial category $\Ecal:\mathbf{\Delta} \to \Cat$ induced by the simplicial algebraic space $\Xfrak_{\bullet}$ then gives rise to the {\'e}tale topos on $\Xfrak_{\bullet}$. While we gave this definition very explicitly above for any simplicial algebraic space, this is the particular construction we want to keep in mind for our applications (as the next remark in particular shows).

\begin{remark}
	Fix a finite type algebraic stack $\Xfrak$ over $\Spec K$ with presentation $X \to \Xfrak$ and let $\Xfrak_{\bullet}$ be the simplicial algebraic space described just before Definition \ref{Defn: Etale topos of simp alg space}. We will be particularly interested in the case that each algebraic space $\Xfrak_n$ is representable by a scheme $X_n$; in this case the simplicial algebraic space $\Xfrak_{\bullet}$ is actually isomorphic to a simplicial scheme $X_{\bullet}$. We have already seen that this this happens whenever the algebraic stack $\Xfrak = [G \backslash X]$ is a quotient stack of a variety $X$ by a smooth affine algebraic group $G$, although the unconvinced reader should see Appendix \ref{Appendix: Stackification via Torsors} for a more detailed justification of this fact.
\end{remark}
\begin{remark}
	There is an alternative {\'e}tale topos associated to the stack $\Xfrak$ described in \cite[Definition 5.2.2]{Behrend} which is called $\Xfrak_{\text{{\'e}t}}$\index[notation]{XfrakEt@$\Xfrak_{\text{{\'e}t}}$}; this is a notion that is well-behaved for all algebraic stacks $\Xfrak$, but we need it only for exactly two insights, so we do not define $\Xfrak_{\text{{\'e}t}}$ here. 
	
	The key insight that that this perspective gives are as follows: First, when $\Xfrak$ is an algebraic stack, there is an equivalence $\Xfrak_{\text{{\'e}t}} \simeq (\Xfrak_{\bullet})_{\cart} = \Scal_{\cart}$ so there is a c-structure on the pair of toposes $(\tottop(\Ecal),\Scal_{\cart}) \simeq (\tottop(\Ecal),\Xfrak_{\text{{\'e}t}})$; cf.\@ \cite[Section 5.4, Page 40]{Behrend} for the equivalence of toposes and Remark \ref{Remark: Behrend c topos facts} for the c-topos structure. Moreover, by \cite[Lemma 5.2.3]{Behrend}, the topos $\Xfrak_{\text{{\'e}t}}$ is Noetherian whenever $\Xfrak$ is of finite type, so whenever $\Xfrak$ is of finite type the c-topos $(\tottop(\Ecal),\Xfrak_{\text{{\'e}t}})$ is Noetherian (cf.\@ Remark \ref{Remark: Behrend c topos facts}).
\end{remark}
\begin{corollary}
	Let $\Xfrak = [G \backslash X]$ when $G$ is a smooth affine algebraic group and $X$ is a left $G$-variety. Then the c-topos $(\tottop(\Ecal),\Xfrak_{\text{{\'e}t}})$ is a Noetherian c-topos where $\Ecal$ is the {\'e}tale topos fibration on $\underline{G \backslash X}_{\bullet}$.
\end{corollary}
\begin{remark}
	Fix an algebraic stack $\Xfrak$. By \cite[Lemma 5.2.3]{Behrend}, if $\Xfrak$ is a finite type stack, then the topos $\Xfrak_{\text{{\'e}t}}$ is a Noetherian topos and so the c-structure $((\Xfrak_{\bullet})_{\text{\'e}t}, \Scal_{\cart}) \simeq ((\Xfrak_{\bullet})_{\text{{\'e}t}},\Xfrak_{\text{{\'e}t}})$ is Noetherian as well. In this case we get that there is a well-defined bounded derived category
	\[
	\Dbb_c^b\big(((\Xfrak_{\bullet})_{\text{{\'e}t}},\Xfrak_{\text{{\'e}t}});\overline{\Q}_{\ell}\big).
	\]
\end{remark}

\begin{definition}[{\cite[Section 5.4]{Behrend}}]\index{Derived $\ell$-adic Category! On a Finite Type Algebraic Stack}
	Let $\Xfrak$ be a finite type algebraic stack over a field $K$. The derived category of $\ell$-adic sheaves on $\Xfrak$ is defined as
	\[
	D_c^b(\Xfrak;\overline{\Q}_{\ell}) := \Dbb_c^b\big(((\Xfrak_{\bullet})_{\text{{\'e}t}},\Ecal_{\cart}); \overline{\Q}_{\ell}\big).
	\]\index[notation]{DbcXfrak@$D_c^b(\Xfrak;\overline{\Q}_{\ell})$}
\end{definition}
\begin{example}\index{Equivariant Derived Categrory! Stacky Version}
	Assume that $G$ is a smooth affine algebraic group over a field $K$ and that $X$ is a left $G$-variety. Then the algebraic stack $[G \backslash X]$ is finite type over $\Spec K$, so the bounded derived $\ell$-adic category on $[G \backslash X]$ is
	\[
	\DbQl{[G \backslash X]} = \Dbb_c^b\left(\big((\underline{[G \backslash X]}_{\bullet})_{\text{{\'e}t}},[G \backslash X]_{\text{{\'e}t}}\big);\overline{\Q}_{\ell}\right)
	\]\index[notation]{Dbcstack@$D^b_C([G \backslash X];\overline{\Q}_{\ell})$}
	and this is equivalent to the category
	\[
	\Dbb_c^b\left(\big((\underline{[G \backslash X]}_{\bullet})_{\text{{\'e}t}},\Scal_{\cart}\big)\right)
	\]
	for $\Scal_{\cart}$ the Cartesian sections of the the {\'e}tale topos fibration $\Ecal$ on $\underline{[G \backslash X]}_{\bullet}$.
\end{example}

With all the theory and constructions we have described above, we can finally prove that $\Dbeqstack{X} \simeq \DbQl{[G \backslash X]}$. The main technique that we will use here is to give a particularly useful $1$-category $\Cscr$ which is equivalent to the $2$-collimit that describes the $\ell$-adic category on $[G \backslash X]$ and then show that this $2$-colimit is equivalent to $\Dbeqstack{X}$, as the category $\Cscr$ is constructed as an incarnation of the fact that $\fCat$ is $2$-cocomplete and is particularly well-suited to comparisons with our category $\Dbeqstack{X}$.

\begin{Theorem}\label{Theorem: Section 4: Equivalence of Behrend's Derived Cat and Simplicial Derived Cat}
	There is an equivalence of categories
	\[
	D_{c}^{b}([G \backslash X];\overline{\Q}_{\ell}) \simeq D_{\text{eq}}^{b}(\underline{G \backslash X}_{\bullet};\overline{\Q}_{\ell}).
	\]
\end{Theorem}
\begin{proof}
	We first note that by Corollary \ref{Cor: Section 4: Bounded derived simplicial cat is equiv} it suffices to prove that 
	\[
	D_c^b([G \backslash X];\overline{\Q}_{\ell}) \simeq \eqstacksimplicial.
	\] 
	Observe that since $[G \backslash X]_{\text{{\'e}t}} \simeq \Ecal_{\cart}$, for $\Ecal$ the fibred topos $\Ecal \to \Delta^{\op}$ constructed above there is a nearly definitional equivalence of categories
	\[
	\Dbb^{b}_c\bigg((\tottop(\Ecal),[G \backslash X]_{\text{{\'e}t}});\overline{\Q}_{\ell}\bigg)  \simeq \Dbb^{b}_{c}\bigg((\tottop(\Ecal),\Ecal_{\cart});\overline{\Q}_{\ell}\bigg) = D^b_c([G \backslash X];\overline{\Q}_{\ell}),
	\]
	so it suffices to prove that the category  $\Dbb^{b}_c((\tottop(\Ecal),[G \backslash X]_{\text{{\'e}t}});\overline{\Q}_{\ell})$ is equivalent to $D^b_{\text{eq}}(\underline{[G \backslash X]}_{\bullet};\overline{\Q}_{\ell}).$ For this we recall from \cite{Behrend} that we have
	\[
	\Dbb_c^{b}\bigg((\tottop(\Ecal),[G \backslash X]_{\text{{\'e}t}});\overline{\Q}_{\ell}\bigg) := \twoColim{(\Scal,\Lcal) \in \Mcal}\Dbb_{(\Scal,\Lcal)}^{b}\bigg((\tottop(\Ecal),[G \backslash X]_{\text{{\'e}t}});\overline{\Q}_{\ell}\bigg)
	\]
	where $\Mcal$ is the filtered poset of pre-L-stratifications of the topos $[G \backslash X]_{\text{{\'e}t}}$. Because $\Mcal$, regarded as a category, has the property that there is at most one morphism between any two objects $(\Scal,\Lcal)$ and $(\Tcal,\Jcal)$, by \cite[Theorem A.3.4, Proposition A.3.6]{StackMicrolocalPerverse} we have that
	\[
	\twoColim{(\Scal,\Lcal) \in \Mcal}\Dbb_{(\Scal,\Lcal)}^{b}\bigg((\tottop(\Ecal),[G \backslash X]_{\text{{\'e}t}});\overline{\Q}_{\ell}\bigg) \simeq \Cscr,
	\]
	where $\Cscr$ is the category defined as follows:
	\begin{itemize}
		\item Objects: Pairs $(\Gscr,(\Scal,\Lcal))$ where $\Gscr \in \Dbb^b_{(\Scal,\Lcal)}((\tottop(\Ecal),[G \backslash X]_{\text{{\'e}t}});\overline{\Q}_{\ell})$ and $(\Scal,\Lcal) \in \Mcal$;
		\item Morphisms: For any two objects $(\Gscr,(\Scal,\Lcal)), (\Gscr^{\prime},(\Scal^{\prime},\Lcal^{\prime}))$, the morphisms are
		\begin{align*}
		&\Cscr\bigg((\Gscr,(\Scal,\Lcal)), (\Gscr^{\prime},(\Scal^{\prime},\Lcal^{\prime}))\bigg) \\
		&= \colim \Dbb_{(\Acal,\Bcal)}^{b}((\tottop(\Ecal),[G \backslash X]_{\text{{\'e}t}});\overline{\Q}_{\ell})\big((\Gscr,(\Scal,\Lcal)), (\Gscr^{\prime},(\Scal^{\prime},\Lcal^{\prime}))\big).
		\end{align*}
		where the colimit is taken over the indexing
		\[
		\colim_{\substack{(\Scal,\Lcal), (\Tcal,\Jcal) \leq (\Acal,\Bcal) \\ (\Acal, \Bcal) \in \Mcal}}.
		\]
		\item Identities: The colimit of the identity maps.
		\item Composition: The colimit of the compositions.
	\end{itemize}
	Our strategy is to now prove that $D^{b}_{\text{eq}}(\underline{[G \backslash X]}_{\bullet};\overline{\Q}_{\ell})$ is equivalent to $\Cscr$, as this will prove the claim.
	
	We now define our functors $F:D^{b}_{\text{eq}}(\underline{[G \backslash X]}_{\bullet};\overline{\Q}_{\ell}) \to \Cscr$ and $H:\Cscr \to D^{b}_{\text{eq}}(\underline{[G \backslash X]}_{\bullet};\overline{\Q}_{\ell})$. It is worth remarking here that the functor $F$ we define will rely in a crucial way on the Axiom of Choice but only in the sense that the well-definition of $F$ depends on Choice in a critical way, not in the sense that the particular decision made is any better than any alternate decision. 
	
	Let us begin by defining our functor $F$. Define $F$ by first noting that an object $\Fscr \in \eqstacksimplicial_0$ if $\Fscr = (\Fscr_n)$ is a sequence of (derived complexes of) sheaves $\Fscr_n \in D_c^{b}([G \backslash X]_n; \overline{\Q}_{\ell})$ such that for all simplicial maps $h \in \mathbf{\Delta}([m],[n])$, the induced morphism $\alpha_h:h^{\ast}\Fscr_n \to \Fscr_m$ is an isomorphism. Chasing through the construction of $[G \backslash X]_{\text{{\'e}t}}$, we see that an object $\Gscr \in D^{b}_{c}([G \backslash X]_{\text{{\'e}t}};\overline{\Q}_{\ell})_0$ is given by a sequence of objects $(\Gscr_n)$ for which given any simplicial map $h \in \mathbf{\Delta}([m],[n])$ the structure map $h^{\ast}\Gscr_n \to \Gscr_m$ is an isomorphism. Thus, in order to determine $F$ on objects we just need to show that $\Fscr$ is $(\Scal,\Lcal)$-constructible for some pre-L-stratification $(\Scal,\Lcal)$. However, since $\Fscr$ is constructible in the usual sense (so in particular each $\Fscr_n$ is constructible on $G^n \times X$), there is a stratification $\Scal$ of the topos $\Ecal_{\cart}$ which witnesses that $\Fscr$ is constructible; moreover, there is also a function $\Lcal$, by construction, which reads out how $\Scal$ locally trivializes $\Fscr$. Thus $\Fscr$ is $(\Scal,\Lcal)$-constructible for some pre-L-stratification. This is where we invoke the Axiom of Choice: {For every object $\Fscr \in D^b_{\text{eq}}(\underline{[G \backslash X]}_{\bullet};\overline{\Q}_{\ell})$, choose exactly one pre-L-stratification $(\Scal,\Lcal)$ for which $\Fscr$ is $(\Scal,\Lcal)$-constructible.} We thus define the functor $F:\eqstacksimplicial \to \Cscr$ by
	\[
	\Fscr \mapsto \big(\Fscr,(\Scal,\Lcal)\big)
	\]
	where $(\Scal,\Lcal)$ is the stratification chosen by our selection via the Axiom of Choice. To define $F$ on morphisms, we note that any morphism 
	\[
	\varphi \in \eqstacksimplicial(\Fscr,\Gscr)
	\] 
	is a morphism for the common refinement of the stratifications chosen for $\Fscr, \Gscr$. As such, we can then limit over all upper bounds to give the assignment on $\varphi$ by sending $\varphi$ to the limit $[\varphi_{(\Acal,\Bcal)}]$ in
	\[
	\colim_{\substack{(\Scal,\Lcal), (\Tcal,\Jcal) \leq (\Acal,\Bcal) \\ (\Acal, \Bcal) \in \Mcal}}\Dbb_{(\Acal,\Bcal)}^{b}((\tottop(\Ecal),[G \backslash X]_{\text{{\'e}t}});\overline{\Q}_{\ell})\big((\Fscr,(\Scal,\Lcal)),(\Gscr,(\Tcal,\Jcal))\big).
	\]
	That this assignment is functorial follows from the functorial nature of colimits.
	
	Let us now construct the functor $H:\Cscr \to \eqstacksimplicial$. On objects, this is straightforward: We send an object $(\Fscr,(\Scal,\Lcal)) \mapsto \Fscr$. By the same reasoning as above (the stratification of $\Ecal_{\cart}$ gives rise to a stratification of each $G^n \times X$ for which $\Fscr_n$ is constructible) what was given above, we see that $\Fscr \in \eqstacksimplicial_0$. To determine the functor on morphisms, let the morphism $[\varphi_{(\Acal,\Bcal)}]$ in
	\[
	 \colim_{\substack{(\Scal,\Lcal), (\Tcal,\Jcal) \leq (\Acal,\Bcal) \\ (\Acal, \Bcal) \in \Mcal}}\Dbb_{(\Acal,\Bcal)}^{b}((\tottop(\Ecal),[G \backslash X]_{\text{{\'e}t}});\overline{\Q}_{\ell})\big((\Fscr,(\Scal,\Lcal)),(\Gscr,(\Tcal,\Jcal))\big).
	\] 
	be given and define
	\[
	(\Acal_0,\Bcal_0) := (\Scal,\Lcal) \lor (\Tcal,\Jcal) = \inf\lbrace (\Acal,\Bcal) \in \Mcal \; | \; (\Scal,\Lcal), (\Tcal,\Jcal) \leq (\Acal,\Bcal)\rbrace.
	\]
	We then determine our functor on morphisms by the assignment sending the map $[\varphi_{(\Acal,\Bcal)}]$ to
	\[
	[\varphi_{(\Scal,\Lcal)}] \mapsto \varphi_{(\Acal_0,\Bcal_0)}.
	\]
	That this assignment is well-defined and preserves identities follows immediately from the fact that producing the join of $(\Scal,\Lcal)$ and $(\Tcal,\Jcal)$ is functorial. Thus we only need to prove that the assignment we gave on morphisms preserves composition to prove functoriality. However, this is routine and involves only the observation that if $\varphi$ is defined over the stratifications $(\Scal,\Lcal)$ and $(\Tcal,\Jcal)$ and $\psi$ is defined over the stratifications $(\Tcal,\Jcal)$ and $(\Ucal,\Ical)$, then the action of $(\psi \circ \varphi)_{(\Acal_2,\Bcal_2)}$ is the same as the action of $\psi_{(\Acal_1, \Bcal_1)} \circ \varphi_{(\Acal_0, \Bcal_0)}$ on the underlying complex $A$ on which the morphisms act.
	
	We now verify the equivalence of categories. On one hand, for each $\Fscr \in \Dbeqstack{X}_0$. Begin by noting that
	\[
	H(F(\Fscr)) = H(\Fscr,(\Scal,\Lcal)) = \Fscr.
	\]
	On the other hand let $(\Fscr,(\Scal,\Lcal)) \in \Cscr_0$. It then follows that
	\[
	F(H(\Fscr,(\Scal,\Lcal))) = F(\Fscr) = (\Fscr,(\Tcal,\Jcal));
	\]
	we claim that in $\Cscr$ that $(\Fscr,(\Scal,\Lcal)) \cong (\Fscr, (\Tcal,\Jcal))$. First fix an $(\Acal,\Bcal) \in \Mcal$ for which $(\Scal,\Lcal), (\Tcal,\Jcal) \leq \Mcal$. The image of $(\Fscr,(\Scal,\Lcal))$ in the category $\Dbb_{(\Acal,\Bcal)}((\tottop(\Ecal),[G \backslash X]_{\text{{\'e}t}});\overline{\Q}_{\ell})$ includes the complex $\Fscr$ into $\Dbb_{(\Acal,\Bcal)}$ and witness the constructibility via restrictions from $\Scal$ to $\Acal$. Similarly, the image of $(\Fscr,(\Tcal,\Jcal))$ includes the complex $\Fscr$ and witnesses the constructibility and trivialization via restrictions from $\Tcal$ to $\Acal$. Now consider that since $\Acal$ simultaneously refines $\Scal$ and $\Tcal$, the diagram
	\[
	\begin{tikzcd}
	\Fscr \ar[rr, ""{name = U}] \ar[d] & & \prod\limits_{V \in \Tcal} \Fscr|_{V} \ar[d] \\
	\prod\limits_{U \in \Scal} \Fscr|_{U} \ar[rr, ""{name = L}] & & \prod\limits_{W \in \Acal} \Fscr|_{W} \ar[from = U, to = L, Rightarrow, shorten <= 4pt, shorten >= 4pt]{}{\cong}
	\end{tikzcd}
	\]
	commutes up to invertible $2$-cell. As such, in the colimit we find that in the colimit that defines the hom-set we get that there is a morphism which acts as the identity on $\Fscr$ and records the $2$-morphisms that go between the different ways of trivializing $\Fscr$ either by starting at $(\Scal,\Lcal)$ or $(\Tcal,\Jcal)$, respectively; call these maps 
	\[
	\varphi:(\Fscr,(\Scal,\Lcal)) \to (\Fscr,(\Tcal,\Jcal))
	\] 
	and
	\[
	\psi:(\Fscr,(\Tcal,\Jcal)) \to  (\Fscr,(\Scal,\Lcal)).
	\]
	It then follows from construction that
	\[
	\psi \circ \varphi = \id_{(\Fscr,(\Scal,\Lcal))}
	\]
	and
	\[
	\varphi \circ \psi = \id_{(\Fscr,(\Tcal,\Jcal))}
	\]
	so we indeed get that
	\[
	(F \circ H)(\Fscr,(\Scal,\Lcal)) = F(\Fscr) = (\Fscr,(\Tcal,\Jcal)) \cong (\Fscr,(\Scal,\Lcal)).
	\]
	This proves that $F$ and $H$ are mutual inverse equivalences, which completes the proof of the theorem.
\end{proof}

Altogether, the work we have done allows us to produce the main theorem of Part \ref{Chapter: EDC Comp}: The four-way-equivalence of equivariant derived categories.
\begin{Theorem}\label{Theorem: Section 4.4: Main Theorem}
	For any smooth affine algebraic group $G$ and any left $G$-variety $X$, there is a four way equivalence of categories
	\[
	\DbeqQl{X} \simeq \quot{\DbeqQl{X}}{ABL} \simeq \Dbeqsimp{X} \simeq \Dbeqstack{X}
	\]
	which is compatable with the standard $t$-structures on each category.
\end{Theorem}
\begin{proof}
	The equivalence $\DbeqQl{X} \simeq \quot{\DbeqQl{X}}{ABL}$ is Corollary \ref{Cor: Section 4: EDCs are equiv}, the equivalence $\DbeqQl{X} \simeq \Dbeqsimp{X}$ is Theorem \ref{Theorem: Section 4.2: Simplicial equivalence with derived category}, and the equivalence $\Dbeqsimp{X} \simeq \Dbeqstack{X}$ is Theorem \ref{Theorem: Section 4: Equivalence of Behrend's Derived Cat and Simplicial Derived Cat}. From here composing against the various equivalences and inverse equivalences give the four-way-equivalence of the theorem. Finally, the statment involving the standard $t$-structures follows from Proposition \ref{Prop: Equiv of Standard t-structures on EDCs}, the fact that the standard $t$-structure is determined by standard chain complex cohomology on each category (and hence is preserved/determined by exact functors), and the fact that equivalences of categories are exact functors.
\end{proof}

\part{Conclusions, Appendices, and Back Matter}
\chapter{Conclusions and Future Work}

As we reach the finish line of this marathon of a thesis, we have taken an arduous journey through the Foothills of Categorical Algebraic Geometry. We've see a deep crawl up the hills of Grothendieck topologies, formal schemes, and adhesive functors (which hopefully have stuck with you); we've taken a long (far longer than advisable, really) run in the Valley of Equivariant Categories; and we finally ended our journey at the four-faced Mount Equivariant Derived Category. While it would be nice to say that the real treasure of this thesis was the friends we made along the way, the length of this thesis does not make friends so it is more appropriate to say that the treasure contained within this thesis is a certain amount of liquid with which to submerge walnuts of yet-to-be-known information. 

In discussing how this liquid may best be applied brings us to a strange aspect of mathematical writings: We should not be concluding a study, but merely stopping to rest. While it is true that eventually we need to wrap up projects and explain what we have found and learned, our goal should be continuing to map out the mathematical landscape when we return from our rest. I would like to treat this section not so much as a conclusion to studying and learning, but instead as a place to write some of the inclusions and suggestions for future works and future places to attempt mapping the mathematical landscape. With this in mind, let us proceed to recall what we did in each main chapter of this thesis and then present some suggestions for future work.

\section{Equivariant Categories}
In Part \ref{Chapter: EQCats} we introduced and studied in great depth equivariant categories over a (quasi-projective) variety $X$. The construction we provided is intimately steeped within $2$-categorical foundations, as this gave us tools and techniques to study what structure is really at hand in any given moment, as well as how it can be applied to deduce properties about equivariant categories. We were able to use these constructions to give careful and precise conditions when the equivariant category $F_G(X)$ has the following properties: $F_G(X)$ admits all (co)limits of specified shape $I$, $F_G(X)$ is additive, $F_G(X)$ is (symmetric) monoidal, and when $F_G(X)$ admits a subobject classifier. After studying the equivariant category itself we also saw how to define functors $F_G(X) \to E_H(X)$ between equivariant categories when there is a morphism $\varphi:H \to G$ of smooth algebraic groups. In studying these functors we were able to essentially prove that these functors more or less arise from pseudonatural transformations between the information expressed by each corresponding equivariant category, which in turn allowed us to use modifications to provide natural transformations between equivariant functors. These observations together allowed us to prove when equivariant functors preserve (co)limits, admit adjoint functors, and even when these lifted adjoints are equivalences. After studying this we were able to prove when $F_G(X)$ is a triangulated category and even examined some $t$-structures on $F_G(X)$ that arise from local $t$-structures varying pseudofunctorially through the resolution category $\SfResl_G(X)$.

An important but missing piece of Part \ref{Chapter Intro EqCats}, however, isdefinition of the equivariant $\infty$-category $\mathsf{F}_G(X)$ (associated to an $(\infty,1)$-pseudofunctor). While a construction of such an object will immediately give the equivariant derived $\infty$-category and $\infty$-category of equivariant perverse sheaves on a variety $X$,we still must rigorously construct such an object. Instead of doing this construction explicitly, I'll leave the following conjectures as a look towards some of the future work involving equivariant $\infty$-categories.

Before giving the conjecture explicitly, I would like to set the stage (albiet briefly and without much in the way of explanation). Let $G$ be a smooth algebraic group and let $X$ be a left $G$-variety over $\Spec K$. Write $\mathsf{Ner}(\SfResl_G(X)^{\op})$ for the simplicial nerve of the category $\SfResl_G(X)^{\op}$. Then $\mathsf{Ner}(\SfResl_G(X)^{\op})$ is the quasi-category where $0$-cells are objects of $\SfResl_G(X)^{\op}$ and for $n \geq 1$ the $n$-cells of $\mathsf{Ner}(\SfResl_G(X)^{\op})$ are given by (the opposite of) length $n$ paths of composable arrows
\[
\xymatrix{
\Gamma_0 \times X \ar[rr]^-{f_1 \times \id_X} & & \Gamma_1 \times X \ar[rr]^-{f_2 \times \id_X} & & \cdots \ar[rr]^-{f_n \times \id_X} & & \Gamma_n \times X
}
\]
in $\SfResl_G(X)$.

To continue the story, let $\mathfrak{QCat}$ denote the $(\infty,2)$-category of quasi-categories and let $\mathsf{Ner}_{\text{hc}}(\mathfrak{QCat})$ denote the homotopy coherent nerve of $\mathfrak{QCat}$. Then $\mathsf{Ner}_{hc}(\mathfrak{QCat})$ is a (large) quasi-category whose $0$-cells are $(\infty,1)$-categories. An $\infty$-pseudofunctor on $\SfResl_G(X)$ should then be an $\infty$-functor
\[
F:\mathsf{Ner}(\SfResl_G(X)^{\op}) \to \mathsf{Ner}_{\text{hc}}(\mathfrak{QCat}).
\]
By giving a straightforward generalization of Definition \ref{Defn: Prequivariant pseudofunctors}, we have pre-equivariant $\infty$-pseudofunctors as the $\infty$-functors $F:\mathsf{Ner}(\SfResl_G(X)^{\op}) \to \mathsf{Ner}_{\text{hc}}(\mathfrak{QCat})$ which factor as
\[
\xymatrix{
\mathsf{Ner}(\SfResl_G(X)^{\op}) \ar[drr]_{F} \ar[rr]^-{\mathsf{Ner}(\quo^{\op})} & & \mathsf{Ner}(\Var_{/\Spec K}^{\op}) \ar[d]^{\tilde{F}} \\
 & & \mathsf{Ner}_{\text{hc}}(\mathfrak{QCat})
}
\]
in the large $(\infty,2)$-category of quasi-categories $\mathfrak{QCAT}$. 
\begin{conjecture}\label{Conjecture: Con 1}
To a pre-equivariant $\infty$-pseudofunctor $F$ on $X$ we have an associated equivariant $\infty$-category $F_G(X)$ whose $0$-cells are given by pairs $(A,T_A)$ where $A$ is a collection
\[
A = \lbrace \AGamma \rbrace_{\Gamma \in \Sf(G)_0}
\]
where each $\AGamma$ is a $0$-cell in $F(\XGamma)$ and $T_A$ is a collection of transition equivalences
\[
T_A = \lbrace \tau_f^A:F(\of)(\AGammap) \xrightarrow{\simeq} \AGamma \rbrace_{f \in \Sf(G)_1}
\]
subject to a higher cocycle condition involving fillers of the desired diagrams.
\end{conjecture}
\begin{conjecture}
Assuming Conjecture \ref{Conjecture: Con 1}, the equivariant derived $\infty$-category on $X$ is given by the pre-equivariant $\infty$-pseudofunctor
\[
F:\mathsf{Ner}(\SfResl_G(X)^{\op}) \to \mathsf{Ner}_{\text{hc}}(\mathfrak{QCat})
\] 
induced by the assignment
\[
F(Y) := \mathsf{D}_c^b(Y),
\]
where $\mathsf{D}_c^b(Y)$ is the bounded derived $\infty$-category of constructible sheaves on a variety $Y$. Similarly, the equivariant $\infty$-category of perverse sheaves is induced by the pre-equivariant $\infty$-pseudofunctor
\[
P:\mathsf{Ner}(\SfResl_G(X)^{\op}) \to \mathsf{Ner}_{\text{hc}}(\mathfrak{QCat})
\]
induced by the assignment
\[
P(Y) := \mathsf{Ner}(\Per(Y)).
\]
\end{conjecture}
%
The following are suggestions of other ways in which we can take the equivariant category and use it to sketch a larger map of arithmetic algebraic geometry:
\begin{itemize}
	\item Give a foundation for equivariant $\infty$-category theory over a scheme. In particular, it would be of interest to work with equivariant $\infty$-categories as model independently as possible, but when restricting to the quasi-categorical model it is of particular interest to establish the following: If each fibre quasi-category $\mathsf{F}(\XGamma)$ is a stable $\infty$-category and each fibre $\infty$-functor $\mathsf{F}(\of)$ is additive or triangulated, is the category $\mathsf{F}_G(X)$ a stable $\infty$-category as well? I am reasonably certain that this will be true, so this leads us to follow up additional work:
	\item Establish that if the $\infty$-category $\mathsf{F}_G(X)$ is a stable $\infty$-category with a $t$-structure and if $\Acal = \mathsf{F}_G(X)^{\heartsuit}$, then is it true that $h\Acal \simeq \big(h(\mathsf{F}_G(X))\big)^{\heartsuit}$?
	\item Can we use the equivariant category $F_G(X)$ to give a foundation for equivariant crystalline cohomology? In particular, by setting $F$ to be the pseudofunctor induced by $F(\XGamma) = \mathbf{SatDieuComp}(\XGamma)$, i.e., $F(\XGamma)$ is the category of saturated Dieudonn{\'e} complexes of $\XGamma$ sheaves, do we get a well-behaved theory of equivariant crystalline cohomology? One particular step to sort out is to study if the derived d{\'e}calage functor $L\eta_p$ lifts to an equivariant version $L\eta_p:D_G^b(X) \to D_G^b(X)$.
	\item Can we extend the equivariant category to non-variety schemes? For instance, what happens when we try to define $F_G(X)$ for pro-varieties, quasi-separated quasi-compact schemes, general schemes, formal schemes, or adic spaces? Do we have to add extra schemes to $\Sf(G)$ and if so, what do we need to add?
	\item While it is extremely preliminary and poorly-formed at this time, it would be interesting (at least to me) to establish that if Conjecture \ref{Conjecture: Con 1} is true, is there an analogue to Remark \ref{Remark: Grothendieck Construction Trick} to equivariant $\infty$-categories by using the $\infty$-Grothendieck construction (cf.\@ \cite[Section 3.2]{LurieHTT}).
\end{itemize}

\section{The Four Flavours of Equivariant Derived Category}
In this section we saw that there is a four-way equivalence of categories
\[
\ABLDbeqQl{X} \simeq \DbeqQl{X} \simeq \Dbeqsimp{X} \simeq \DbQl{[G \backslash X]}.
\]
We did this very explicitly and in particular emphasized how to do computations and translations with the category $\Dbeqsimp{X}$, as it is the most straightforward category to define and work in directly. A large benefit to this result is that we now know that when we are using the equivariant derived category it suffices to work with the incarnation of the object that suits us best in the moment. For instance, if we just have a variety and a group acting on it and we need to get at the equivariant derived category for explicit computations, it is easiest to work with $\Dbeqsimp{X}$, as in this case we only need to define suitably compatible complexes of ($\ell$-adic) sheaves on the varieties $G^n \times X$; this is much more straightforward than either finding all $\Sf(G)$-varieties or a sufficiently nice family of acyclic resolutions of $X$. Alternatively, we have seen from our work on equivariant categories, that if we need to check or establish a categorical property of $\DbeqQl{X}$, it is easiest to use the formulation of Lusztig.

Much of the work we suggest doing with this four-way equivalence comes either from applications or seeing how to potentially extend such an equivalence, as we can see below.
\begin{itemize}
	\item See if we can establish the four-way equivalence at the full derived categorical level:
	\[
	D_G(X;\overline{\Q}_{\ell}) \simeq \quot{D_G(X;\overline{\Q}_{\ell})}{ABL} \simeq D_{\text{eq}}(\underline{G \backslash X}_{\bullet};\overline{\Q}_{\ell}) \simeq D_c([G \backslash X];\overline{\Q}_{\ell}).
	\]
	I expect this to work out because the first two equivalences were largely formal and did not depend on the boundedness of the derived categories in which we took local values, while the last equivalence should carry over because the main point of the derived category of \cite{Behrend} is that it does allow the unbounded derived category of a stack.
	\item Try to see if this equivalence can be maintained when we allow either $G$ or $X$ to be provarieties or if we allow $X$ to be a quasi-separated quasi-compact scheme instead of just varieties. Whether this can happen or not will depend on certain descent properties, such as whether or not our stack $[G \backslash X]$ remains of finite type or if the schemes $\XGamma$ remain quasi-separated and quasi-compact, so it would be interesting to see if any of the equivalences remain or if any break.
	\item Use the simplicial equivariant derived category to compute equivariant nearby and vanishing cycles functors explicitly and compute equivariant trace explicitly as well.
\end{itemize}

\section{Summary}

In a final summary, I hope that this thesis has contributed myriad categorical techniques with which to study arithmetic algebraic geometry. We introduced and performed a systematic study of equivariant categories over schemes. This systematic study should have many applications in representation theory and the Langlands Programme, as it provides a unified and careful description of why certain equivariant categories have the various properties they do. It also has the benefit of introducing equivariant $\infty$-categories and an equivariant $\infty$-categorical version of $\ell$-adic perverse sheaves in a way that will allow us to give a more in depth study of these objects. Finally, we showed that the four main flavours of bounded equivariant derived category, while different at a first taste, are actually secretly the same. This allows the benefit of being able to use the particular incarnation of the equivariant bounded derived category that best fits what we are studying and what we need to do. In particular, if we need to work with the information over smooth free $G$-varieties, we can use $\DbeqQl{X}$; if we have acyclic resolutions readily at hand, we can use $\ABLDbeqQl{X}$; if we need to do explicit calculations, we can use $\Dbeqsimp{X}$; and if we are in a stacky setting, we can use $\DbeqQl{[G \backslash X]}$.

With the future work and summary now wrapped up, we now come to the absolute strangest aspect of reading, writing, and generally experiencing math papers: They just end.

\appendix
\chapter{The {\'E}tale Cohomology Comparison Theorem and Stiefel Manifolds}\label{Appendix: Comparison Theorem and Stiefel Manifolds}

This appendix serves to compute the various results on the {\'E}tale Cohomology Comparison Theorem for fields $K$ and to compute the $\ell$-adic cohomology of the Stiefel manifolds over $\Spec K$ needed for Proposition \ref{Prop: Section 4: Existence of n acyclic maps}. This will also justify allowing the field $K$ to be of arbitrary characteristic in Part \ref{Chapter Intro Comparison}.
\section{The {\'E}tale Cohomology Comparison Theorem for Fields}
\begin{lemma}\label{Lemma: Base Change Galois}
Let $K/k$ be a Galois extension of fields and let $X$ be a scheme over $\Spec k$ with base change $X_K$. Then for any finite Abelian group $A$ of order coprime to the characteristic of $k$ and for any $n \in \Z$, there is an isomorphism
\[
H^n_{\text{{\'e}t}}(X;A) \cong H^n_{\text{{\'e}t}}(X_K;A).
\]
\end{lemma}
\begin{proof}
It suffices, after writing $K$ as a colimit of finite Galois extensions and then passing to the limit, to assume that $K/k$ is finite Galois. Apply the Smooth Base Change Theorem (cf.\@ \cite[Theorem I.7.3]{FreitagKiehl} or \cite[Theorem 20.1]{milneLEC}) to the diagram
\[
\xymatrix{
X_K \ar[r] \ar[d] \pullbackcorner & \Spec K \ar[d] \\
X \ar[r] & \Spec k
}
\]
where we observe that $\Spec K \to \Spec k$ is smooth because the extension $K/k$ is Galois and hence smooth of relative dimension $0$ in order to derive the isomorphism.
\end{proof}
\begin{corollary}\label{Cor: Cohomology is the same along Galois extns in ladic case}
Let $K/k$ be a Galois extension of fields, let $X$ be a scheme over $\Spec K$, and let $\ell$ be a prime distinct from the characteristic of $k$. Then for any $n \in \Z$ there are isomorphisms
\[
H^n(X;\Z_{\ell}) \cong H^n(X_K;\Z_{\ell})
\]
and
\[
H^n(X;\overline{\Q}_{\ell}) \cong H^n(X_K;\overline{\Q}_{\ell}).
\]
\end{corollary}
\begin{proof}
The trick is to show for the constant $\ell$-adic sheaf $\mathbbm{1}_X = (\mathbbm{1}_m)$ on $X$ that the result holds. However, this amounts to calculating, for $\nu_X:X \to \Spec k$ and $\nu_{X_K}:X_K \to \Spec K$
\[
H^n_{\text{{\'e}t}}(X;\mathbbm{1}(X)) = R^n(\nu_X)_{\ast}(\mathbbm{1}_m) \cong H^n_{\text{{\'e}t}}(\mathbbm{1}_{X_K};f^{\ast}(\mathbbm{1}_{X_K})_m(X_K)) = R^n(\nu_{X_K})(f^{\ast}\mathbbm{1}_m)
\]
where $f:X_K \to X$ is the base change morphism. Noting
\[
\mathbbm{1}_m(X) \cong \frac{\Z}{\ell^m\Z} \cong (\mathbbm{1}_{X_K})_m(X_K)
\]
applying Lemma \ref{Lemma: Base Change Galois} above gives the result for each $m \in \N$. Thus we have
\[
H^n(X;\Z_{\ell}) \cong H^n(X_K;\Z_{\ell});
\]
the result follows similarly in the $\overline{\Q}_{\ell}$ case.
\end{proof}
\begin{corollary}\label{Cor: Base Change to closure}
Let $k$ be a field and let $k^{\text{sep}}/k$ be a separable closure of $k$. Then for a scheme $X$ defined over $\Spec k$ and any $n \in \N$ there is an isomorphism
\[
H^n(X;\overline{\Q}_{\ell}) \cong H^n(X_{k^{\text{sep}}};\overline{\Q}_{\ell})
\]
\end{corollary}
\begin{proof}
This is Corollary \ref{Cor: Cohomology is the same along Galois extns in ladic case} specialized to the Galois extension $k^{\text{sep}}/k$.
\end{proof}
\begin{lemma}\label{Lemma: Lift to algebraically closed of char 0}
Let $k$ be a separably closed field and let $X$ be a smooth scheme over $\Spec k$. Then for any $n \in \Z$ there is a smooth scheme $\overline{X}$ defined over an algebraically closed field $K$ of characteristic zero such that there is an isomorphism
\[
H^n(X;\overline{\Q}_{\ell}) \cong H^n(\overline{X};\overline{\Q}_{\ell})
\]
where $\ell$ is coprime to the characteristic of $k$.
\end{lemma}
\begin{proof}
If $k$ is of characteristic $0$ then set $X = \overline{X}$ and note that the isomorphism is a tautalogy. Assume now that $k$ is has characteristic $p > 0$. We can then use a Witt vector construction to realize $k$ as a residue field for a complete discrete valuation ring $A$ of characteristic zero (cf.\@ \cite{Teichmüller1937}). Now use the specialization equivalence of sites
\[
\operatorname{sp}:(\Sch_{/\Spec k}, \text{f{\'e}t}) \xrightarrow{\simeq} (\Sch_{/\Spec R}, \text{f{\'e}t})
\]
to lift $X$ to a smooth scheme $\tilde{X} := \operatorname{sp}(X)$ over the trait $\Spec A$. Base changing along the generic point $j:\eta \to \Spec R$ gives a smooth scheme $\hat{X} = j^{\ast}(\tilde{X})$ over $K = \operatorname{Frac}(R)$. Finally, define $\overline{X}$ to be the base change of $\hat{X}$ to $\Spec \overline{K}$. The equivalence of categories $\operatorname{sp}$ implies that
\[
H^n(X;\overline{\Q}_{\ell}) \cong H^n(\tilde{X};\overline{\Q}_{\ell})
\]
and Smooth Base Change gives the isomorphisms
\[
H^n(\tilde{X};\overline{\Q}_{\ell}) \cong H^n(\hat{X};\overline{\Q}_{\ell}) \cong H^n(\overline{X};\overline{\Q}_{\ell}),
\]
and hence
\[
H^n(X;\overline{\Q}_{\ell}) \cong H^n(\overline{X};\overline{\Q}_{\ell}).
\]
\end{proof}
\begin{lemma}\label{Lemma: Complex isomorphism}
Let $k$ be a separably closed field and let $X$ be a smooth scheme over $k$. Then there is a finite type $\C$-scheme $\overline{X}_{\C}$ and an isomorphism
\[
H^n(X;\overline{\Q}_{\ell}) \cong H^n(\overline{X}_{\C};\overline{\Q}_{\ell}).
\]
\end{lemma}
\begin{proof}
Apply the proof of Lemma \ref{Lemma: Lift to algebraically closed of char 0} to give a scheme $\hat{X}$ over a field $K$ of characteristic $0$ with base change $\overline{X}$ to the algebraic closure $\overline{K}$ of $K$.  There is then an isomorphism
\[
H^n(X;\overline{\Q}_{\ell}) \cong H^n(\overline{X};\overline{\Q}_{\ell}),
\]
by Lemma \ref{Lemma: Lift to algebraically closed of char 0}.

Observe now that because both $\overline{K}$ and $\C$ are algebraically closed fields of characteristic $0$, there are equivalences of categories
\[
\begin{tikzcd}
\DbQl{\Spec \overline{K}} \ar[r, bend left = 30, ""{name = U}]{}{E_{\overline{K}}} & D^b_c(\QlVectfd) \ar[l, bend left = 30, ""{name = L}]{}{F_{\overline{K}}} \ar[from = U, to = L, symbol = \simeq]
\end{tikzcd}
\]
and
\[
\begin{tikzcd}
\DbQl{\Spec \C} \ar[r, bend left = 30, ""{name = U}]{}{E_{\C}} & D^b_c(\QlVectfd) \ar[l, bend left = 30, ""{name = L}]{}{F_{\C}} \ar[from = U, to = L, symbol = \simeq]
\end{tikzcd}
\]
where $\QlVectfd$ is the category of finite dimensional $\overline{\Q}_{\ell}$ vector spaces. Recall here that for a variety $Y$ over a field $F$, $\nu_Y:Y \to \Spec F$ is the structure map and $\mathbbm{1}_{Y}$ is the constant $\ell$-adic sheaf on $Y$. Find a complex $A \in D^b_c(\QlVectfd)_0$ such that 
\[
R(\nu_{\overline{X}})_{\ast}(\mathbbm{1}_{\overline{X}}) \cong F_{\overline{K}}(A)
\] 
and 
\[
E_{\overline{K}}\left(R(\nu_{\overline{X}})_{\ast}(\mathbbm{1}_{\overline{X}})\right) \cong A;
\]
such a complex exists by the equivalence of categories. Moreover, because these equivalences may be assumed to be adjoint equivalences (cf.\@ \cite[Proposition 4.4.5]{riehl2017category}) they preserve the standard $t$-structures on each category. As such there are $\overline{\Q}_{\ell}$-linear isomorphisms, for all $q \in \Z$,
\[
H^{q}_{\text{{\'e}t}}\left(R(\nu_{\overline{X}})_{\ast}(\mathbbm{1}_{\overline{X}});\overline{\Q}_{\ell}\right) \cong H^{q}_{\text{{\'e}t}}(F_{\overline{K}}(A);\overline{\Q}_{\ell}) \cong H^{q}(A;\overline{\Q}_{\ell}).
\]
Now let $\Fscr \in \DbQl{\Spec \C}_0$ such that $E_{\C}(\Fscr) \cong A$ and $F_{\C}(A) \cong \Fscr$. We then get $\overline{\Q}_{\ell}$-linear isomorphisms, for all $q \in \Z$,
\[
H^{q}(A;\overline{\Q}_{\ell}) \cong H^{q}_{\text{{\'e}t}}(F_{\C}(A);\overline{\Q}_{\ell}) \cong H^{q}_{\text{{\'e}t}}(\Fscr;\overline{\Q}_{\ell})
\]
which in turn gives an isomorphism
\[
H^{q}_{\text{{\'e}t}}\left(R(\nu_{\overline{X}})_{\ast}(\mathbbm{1}_{\overline{X}});\overline{\Q}_{\ell}\right) \cong H^{q}_{\text{{\'e}t}}(\Fscr;\overline{\Q}_{\ell}).
\]
Find a finite type scheme $\overline{X}_{\C}$ over $\Spec \C$ for which
\[
R(\nu_{\overline{X}_{\C}})_{\ast}(\mathbbm{1}_{\overline{X}_{\C}}) \cong \Fscr.
\]
We then get that
\[
H^q_{\text{{\'e}t}}\left(R(\nu_{\overline{X}_{\C}})_{\ast}(\mathbbm{1}_{\overline{X}_{\C}});\overline{\Q}_{\ell}\right) \cong H^q_{\text{{\'e}t}}(\Fscr;\overline{\Q}_{\ell})
\]
so
\[
H^{q}_{\text{{\'e}t}}\left(R(\nu_{\overline{X}})_{\ast}(\mathbbm{1}_{\overline{X}});\overline{\Q}_{\ell}\right) \cong H^q_{\text{{\'e}t}}\left(R(\nu_{\overline{X}_{\C}})_{\ast}(\mathbbm{1}_{\overline{X}_{\C}});\overline{\Q}_{\ell}\right);
\]
however, this implies that
\[
H^q(\overline{X};\overline{\Q}_{\ell}) \cong H^q(\overline{X}_{\C};\overline{\Q}_{\ell}),
\]
which finally shows that $H^{q}(X;\overline{\Q}_{\ell}) \cong H^{q}(\overline{X}_{\C};\overline{\Q}_{\ell})$ for all $q \in \Z$.
\end{proof}
\begin{Theorem}\label{Thm: Comparison}
Let $k$ be a field with separable closure $k^{\text{sep}}$, let $\ell$ be a prime which does not divide the characteristic of $k$, and let $X$ be a smooth scheme over $k$. Then there is a smooth scheme $\overline{X}$ over $\Spec \C$ for which there is an isomorphism
\[
H^n(X;\overline{\Q}_{\ell}) \cong H^n(\overline{X};\overline{\Q}_{\ell}) \cong H^n_{\text{sing}}(\overline{X}_{\text{an}};\overline{\Q}_{\ell})
\]
\end{Theorem}
\begin{proof}
Corollary \ref{Cor: Base Change to closure} and Lemma \ref{Lemma: Complex isomorphism} give the complex scheme $\overline{X}$ for which the isomorphism
\[
H^n(X;\overline{\Q}_{\ell}) \cong H^n(X_{k^{\text{sep}}};\overline{\Q}_{\ell}) \cong H^n(\overline{X};\overline{\Q}_{\ell})
\]
holds. Now appealing to the Comparison Theorem (cf.\@ \cite[Expos{\'e} XII]{sga1}) gives the isomorphism
\[
H^n(\overline{X};\overline{\Q}_{\ell}) \cong H^n_{\text{sing}}(\overline{X}_{\text{an}};\overline{\Q}_{\ell}),
\]
as was desired.
\end{proof}

\section{The $\ell$-adic Cohomology of Stiefel Manifolds}
We now show how how Lemma \ref{Lemma: Lift to algebraically closed of char 0} interacts with the Stiefel manifolds $W_{n,m}$ over a field $\Spec k$ by ultimately computing the $\ell$-adic cohomology of these schemes. The idea here is that we want to reduce the computation of $H(W_{n,m}^{\Spec k},\Z_{\ell})$ to that of $H^n_{\text{sing}}(W_{n,m}^{\Spec \C};\Z_{\ell})$, which is known by results of \cite{Borel}. Thus for this we will need to see how the variety $W_{n,m}^{k}$ changes, via specialization, to a reduced separated finite type scheme over $\Spec A$ and hence to a variety over $\Spec K$.

We proceed by setting the stage for what is to come by defining our Stiefel manifolds over arbitrary rings. Let $R$ be a commutative ring with identity and let $m,n \in \N$ with $m > n$. Write, for any $t \in \N$,
\[
\GL(t)_R := \Spec\left(\frac{R[x_{ij},y:1 \leq i,j \leq t]}{(\det(x_{ij})y-1)}\right).
\] 
Now embed $\GL(n)_R \to \GL(m)_R$ via the matrix embedding
\[
\begin{pmatrix}
a_{11} & \cdots & a_{1n} \\
\vdots & & \vdots \\
a_{n1} & \cdots & a_{nn}
\end{pmatrix} \mapsto \begin{pmatrix}
a_{11} & \cdots & a_{1n} & 0 \\
\vdots & & \vdots & 0 \\
a_{n1} & \cdots & a_{nn} & 0 \\
0 & \cdots & 0 & I_{m-n}
\end{pmatrix}
\]
where $I_{m-n}$ is the $(m-n) \times (m-n)$ identity matrix. Then there is a closed subgroup $P_R$ of $\GL(m)_R$ induced by the equations describing the span of matrices of the form
\[
\begin{pmatrix}
I_{n} & a_{1,n+1} & \cdots & a_{1m} \\
0 & \ast & \ast & \ast \\
0 & a_{m,n+1} & \cdots & a_{mm}
\end{pmatrix}.
\]
This subgroup $P_R$ of matrices is stabilized by conjugation by $\GL(n)_R$. We then define the Stiefel manifold over $\Spec R$, $W_{n,m}^{R}$, to be the quotient scheme
\[
W_{m,n}^{R} := \frac{\GL(m)_R}{P_R},
\]
which over $\R$ and $\C$ gives the space of $n$-frames in $\R^m$ or $\C^m$, respectively.
\begin{proposition}
	Let $k$ be a separably closed field with $m, n \in \N$ such that $m > n$. Then $\operatorname{sp}(W_{m,n}^{k}) \cong W_{n,m}^{A}$ (where $A$ is the complete DVR of characteristic zero with residue field $k$ given in Lemma \ref{Lemma: Lift to algebraically closed of char 0}). Furthermore, if $j:\eta \to \Spec A$ is the open embedding of the generic fibre, then $j^{\ast}(W_{m,n}^{A}) \cong W_{m,n}^{K}$ for $K = \operatorname{Frac}(A)$.
\end{proposition}
\begin{proof}
	Write We begin by making the specialization calculation. For this we recall from \cite[Expos{\'e} XIII.2]{SGA7} that the quasi-inverse to specialization is base change along the closed immersion $i:s \to \Spec A$ whihch embeds the special fibre $s$ of $\Spec A$. Thus, using that adjoints may always be promoted to adjoint equivalences (cf.\@ \cite[Proposition 4.4.5]{riehl2017category}), it suffices to show that $i^{\ast}(W_{m,n}^{A}) \cong W_{m,n}^{k}$. Note that this follows because from the adjunction we have for any scheme $X$ over $\Spec k$ and any scheme $Y$ over $\Spec A$
	\begin{prooftree}
		\AxiomC{$\operatorname{sp}(X) \cong Y$}
		\UnaryInfC{$X \cong i^{\ast}(Y)$}
	\end{prooftree}
	so since $i^{\ast} \circ \operatorname{sp} \cong \id_{\Sch_{/\Spec k}}$ and $\operatorname{sp} \circ i^{\ast} \cong \id_{\Sch_{/\Spec A}}$ we would get
	\[
	\operatorname{sp}(W_{m,n}^{k}) \cong \operatorname{sp}\left(i^{\ast}(W_{m,n}^{A})\right) \cong W_{m,n}^{A}.
	\]
	
	We now establish the desired computation. Note first that a routine calculation shows
	\begin{align*}
	i^{\ast}\GL(m)_{A} &\cong \GL(m)_A \times_{\Spec A} s \cong \GL(m)_A \times_{\Spec A} \Spec k \\
	&\cong \GL(m)_A \times_{\Spec A} \Spec\left(\frac{A}{\mfrak}\right)\\
	& \cong \Spec \left(\frac{A[x_{ij},y: 1\leq i,j \leq m]}{(\det(x_{ij})y - 1)}\right) \times_{\Spec A} \Spec \left(\frac{A}{\mfrak}\right) \\
	&\cong \Spec\left(\frac{k[x_{ij},y:1\leq i,j\leq m]}{(\det(x_{ij})y-1)}\right) \cong \GL(m)_k
	\end{align*}
	and similarly $i^{\ast}\GL(n)_A \cong \GL(n)_k$. From this we find that the embedding $\GL(n)_A \to \GL(m)_A$ which defines $P_A$ is preserved by $i^{\ast}$ and $i^{\ast}P_A$ is isomorphic to the algebraic subgroup of $\GL(m)_k$ generated by matrices of the form
	\[
	\begin{pmatrix}
	\pi_{\mfrak}(I_{m-n}) & \pi_{\mfrak}(a_{1,n+1}) & \cdots & \pi_{\mfrak}(a_{1,m}) \\
	0 & \ast & \ast & \ast & \\
	0 & \pi_{\mfrak}(a_{m,n+1}) & \cdots & \pi_{\mfrak}(a_{mm})
	\end{pmatrix}
	\]
	where $\pi_{\mfrak}:A \to A/\mfrak$ is the canonical surjection and each $a_{ij} \in A$. However, this cuts out exactly the subvariety $P_k$ of $\GL(m)_k$, so we conclude the existence of the isomorphism $i^{\ast}P_A \cong P_k$. Now we use the exactness of the functor $i^{\ast}$ to get
	\[
	i^{\ast}(W_{m,n}^{A}) = i^{\ast}\left(\frac{\GL(m)_{A}}{P_A}\right) \cong \frac{i^{\ast}\GL(m)_A}{i^{\ast}P_A} \cong \frac{\GL(m)_k}{P_k} = W_{m,n}^{k}.
	\]
	From our earlier observations we conclude that $\operatorname{sp}(W_{m,n}^{k}) \cong W_{m,n}^{A}.$
	
	We now establish the result about base changing along $j:\eta \to \Spec A$. Note that first
	\begin{align*}
	j^{\ast}\GL(m)_A &\cong \GL(m)_A \times_{\Spec A}  \eta \cong \GL(m)_A \times_{\Spec A} \Spec K \cong \GL(m)_K
	\end{align*}
	by an argument mutatis mutandis to the calculation of $i^{\ast}\GL(m)_A$. Similarly, we calculate that $j^{\ast}P_A \cong P_K$ by a standard localization of scalars argument. Using now that $j^{\ast}$ preserves quotient schemes (by virtue of $j:\eta \to \Spec A$ being both epic and monic in $\Sch$) we get
	\[
	j^{\ast}W_{m,n}^{A} = j^{\ast}\left(\frac{\GL(m)_A}{P_A}\right) \cong \frac{j^{\ast}\GL(m)_A}{j^{\ast}P_A} \cong \frac{\GL(m)_K}{P_K} = W_{m,n}^{K},
	\]
	as was desired.
\end{proof}
We now record an important proposition which gives the singular cohomology of $W_{m,n}^{\C}$ for $m > n$. This was computed in \cite{Borel} and presented in \cite{TopologyII}, but we give the results here for later use.
\begin{proposition}[\cite{Borel}, \cite{TopologyII}]\label{Prop: Cohomology Calc for Steifel}
Let $m,n \in \N$ with $m > n$. Then there is an isomorphism of cohomology rings
\[
H^{\bullet}_{\text{sing}}((W_{m,n}^{\C})_{\text{an}};\Z) \cong H^{\bullet}_{\text{sing}}\left(\prod_{i=1}^{n} \Sbb^{2m - 2(n-i) -1};\Z\right).
\]
\end{proposition}
As an immediate result, we use Proposition \ref{Prop: Cohomology Calc for Steifel} and the Universal Coefficient Theorem to derive the following.
\begin{proposition}\label{Prop: The bounds of calculation}
Let $m,n \in \N$ with $m > n$. Then:
\begin{itemize}
	\item For any $q \in \Z$, $H^{q}_{\text{sing}}((W_{m,n}^{\C})_{\text{an}};\Z_{\ell})$ is a free $\Z_{\ell}$-module;
	\item For any $q \in \Z$ with $1 \leq q \leq 2m - 2,$ $H^{q}_{\text{sing}}((W_{m,n}^{\C})_{\text{an}};\Z_{\ell}) \cong 0$.
\end{itemize}
\end{proposition}
\begin{proof}
For the first claim, we note that for any $q \in \Z$ we have
\begin{align*}
H^{q}((W_{m,n}^{\C})_{\text{an}};\Z) &\cong H^{q}_{\text{sing}}\left(\prod_{i=1}^{n} \Sbb^{2m - 2(n-i) -1};\Z\right) \\
&\cong \bigoplus_{q = a_1 + \cdots + a_n} \bigotimes_{i=1}^{n} H^{a_i}_{\text{sing}}(\Sbb^{2m - 2(n-i) - 1};\Z)
\end{align*}
by applying the K{\"u}nneth Formula and using that all cohomology is concentrated in nonnegative degrees, i.e., using that for any $s,t$ we have
\[
H^{s}(\Sbb^{t};\Z) \cong \begin{cases}
\Z & \text{if}\, s = 0, t; \\
0 & \text{else}.
\end{cases}
\]
With this observation in mind it follows that each group
\[
\bigoplus_{q = a_1 + \cdots + a_n} \bigotimes_{i=1}^{n} H^{a_i}_{\text{sing}}(\Sbb^{2m - 2(n-i) - 1};\Z)
\]
is free; as such, each $\Z_{\ell}$-module
\[
\bigoplus_{q = a_1 + \cdots + a_n } \bigotimes_{i=1}^{n} H^{a_i}_{\text{sing}}(\Sbb^{2m - 2(n-i) - 1};\Z_{\ell})
\]
is a free $\Z_{\ell}$-module, which gives the first claim.

For the second claim, fix $q \in \Z$ with $1 \leq a \leq 2m - 2$. Because each cohomology group $H^{s}(\Sbb^t;\Z)$ is nonzero only in nonnegative degrees, it suffices to show that each tensor product
\[
\bigotimes_{i=1}^{n} H^{a_i}_{\text{sing}}(\Sbb^{2m - 2(n-i) - 1};\Z_{\ell}) \cong 0
\]
with $a_1 + \cdots + a_n = q$ and where each $a_i \geq 0$. However, this follows from a straightforward but tedious combinatorial argument which culminates in showing that no $a_i$ may possibly contribute nonzero cohomology because of the dimension of the spheres $\Sbb^{2m - 2(n-i)-1}$ and the fact that $m > n$.
\end{proof}
\begin{Theorem}\label{Thm: Stiefel Computation}
Let $k$ be a field and let $m, n \in \N$ with $m > n$. Then
\[
H^{\bullet}(W_{m,n}^{k};\Z_{\ell}) \cong H^{\bullet}(W_{m,n}^{\C};\Z_{\ell}) \cong H^{\bullet}_{\text{sing}}((W_{m,n}^{\C})_{\text{an}};\Z_{\ell})
\]
and
\[
H^{\bullet}(W_{m,n}^{k};\overline{\Q}_{\ell}) \cong H^{\bullet}(W_{m,n}^{\C};\overline{\Q}_{\ell}) \cong H^{\bullet}_{\text{sing}}((W_{m,n}^{\C})_{\text{an}};\overline{\Q}_{\ell}).
\]
In particular, the cohomology $H^{q}(-;\Z_{\ell})$ above vanishes for all $q < 0, 1 \leq q \leq 2m - 2$ and is a free $\Z_{\ell}$-module otherwise.
\end{Theorem}
\begin{proof}
The isomorphisms between cohomology groups follow from Theorem \ref{Thm: Comparison} and the fact that in the proof of Lemma \ref{Lemma: Complex isomorphism} we can take $W_{m,n}^{\C}$ as the scheme over $\Spec \C$ which has the same $\ell$-adic cohomology as $W_{m,n}^{\overline{K}}$ (by construction, in fact). The statement about the bounds of the cohomology follow from Proposition \ref{Prop: The bounds of calculation} and further changing rings $\Z_{\ell} \to \overline{\Q}_{\ell}$ gives the final result.
\end{proof}
\chapter{Internal Categories, Internal Groupoids, and Stackification via Torsors}\label{Appendix: Stackification via Torsors}
This appendix is a quick introduction and tour of some pre-basic theory of internal categories and internal groupoids that can be found within categories. While we will not even need a smidgen of the full weight of the theory of essentially algebraic theories or their models, it is helpful to see exactly why it is that a groupoid scheme induces a pseudofunctor which takes values in the $2$-category $\mathfrak{Gpd}$\index[notation]{GpdFrak@$\mathfrak{Gpd}$} of groupoids. This will also elucidate why we can take such a stackification and why intuitively it matches up with a quotient stack $[X/G]$ when $X$ and $G$ are schemes (and in fact varieties).

For what follows we will need to recall the notion of a category internal to a $1$-category $\Cscr$. 
\begin{definition}\index{Internal Category}
	An internal category $\Cbb$ to a category $\Cscr$ is given by the following information:
	\begin{itemize}
		\item An object $C_0$ that we think of as the object-of-objects of $\Cbb$;
		\item An object $C_1$ that we think of as the object-of-morphisms of $\Cbb$;
		\item Morphisms $\dom, \codom:C_1 \to C_0$ and a morphism $s_0:C_0 \to C_1$ forw which the diagram
		\[
		\xymatrix{
			C_0 \ar[r]^-{s_0} \ar@/_/@{=}[dr] & C_1 \ar@<-0.5ex>[d]_{\dom} \ar@<.5ex>[d]^{\codom} \\
			& C_0	
		}
		\]
		commutes.
		\item The pullback
		\[
		\xymatrix{
			C_1 \times_{C_0}C_1 \ar[r]^-{p_2} \ar[d]_{p_1} \pullbackcorner & C_1 \ar[d]^{\codom} \\
			C_1 \ar[r]_-{\dom} & C_0
		}
		\]
		exists in $\Cscr$.
		\item There is an internal composition morphism $m:C_2 \to C_1$ such that the diagrams
		\[
		\xymatrix{
			(C_1 \times_{C_0} C_1) \times_{C_0} C_1 \ar[d]_{\cong} \ar[rr]^-{m \times \id_{C_1}} \ar[d]_{\cong} & & C_1 \times_{C_0}C_1 \ar[d]^{m} \\
			C_1 \times_{C_0} (C_1 \times_{C_0} C_1) \ar[dr]_{\id_{C_1} \times m} & & C_1 \\
			& C_1 \times_{C_0} C_1 \ar[ur]_{m}
		}
		\]
		and
		\[
		\xymatrix{
			C_0 \times_{C_0} C_1 \ar[rr]^-{s_0 \times \id_{C_1}} \ar[drr]_{p_2} & & C_1 \times_{C_0} C_1 \ar[d]_m & & C_1 \times_{C_0} C_0 \ar[ll]_-{\id_{C_0} \times s_0} \ar[dll]^{p_1} \\
			& & C_1
		}
		\]
		both commute.
	\end{itemize}
\end{definition}
When we need to refer to an internal category and all the data that determines the object, we will write $\Cbb$ as a sextuple $\Cbb = (C_0, C_1, \dom, \codom, s_0, m)$. If we are being really quite eexplicit, we will even name the domain and codomain of all the morphisms. 
\begin{definition}\index{Internal Groupoid}
An internal groupoid in a category $\Cscr$ is an internal category $\Gbb = (G_0, G_1, \dom, \codom, s_0, m)$ together with a structure morphism $\inv:G_1 \to G_1$ such that the diagrams
	\[
	\xymatrix{
		G_1 \ar[r] \ar[dr]_{\codom} \ar[r]^{\inv} & G_1 \ar[d]^{\dom} & G_1 \ar[dr]_{\dom} \ar[r]^{\inv} & G_1 \ar[d]^{\codom} \\
		& G_0 & & G_0
	}
	\]
	\[
	\xymatrix{
		G_1 \ar[d]_{\dom} \ar[r]^-{\Delta} & G_1 \times_{G_0} G_1 \ar[r]^-{\id_{G_1} \times \inv} & G_1 \times_{G_0} G_1 \ar[d]^{m} \\
		G_0 \ar[rr]_-{s_0} & & G_1
	}
	\]
	\[
	\xymatrix{
		G_1 \ar[d]_{\codom} \ar[r]^-{\Delta} & G_1 \times_{G_0} G_1 \ar[r]^-{\inv \times \id_{G_1}} & G_1 \times_{G_0} G_1 \ar[d]^{m} \\
		G_0 \ar[rr]_-{s_0} & & G_1
	}
	\]
	all commute in $\Cscr$.
\end{definition}
\begin{lemma}\label{Lemma: Section 4: internal category}
	Let $\Cscr$ be a category with $G$ a group object of $\Cscr$ and $X$ a left $G$-object in $\Cscr$. Then the sextuple $(X, G \times X, \pi_{2}:G \times X \to X, \alpha_X:G \times X \to X, \langle 1_G, \id_X\rangle:X \to G \times X, \mu \times \id_X: G \times G \times X \to G \times X)$ is an internal category $\mathbb{G\times X}$ in $\Cscr$.
\end{lemma}\index[notation]{GtimesX@$\mathbb{G \times X}$}
\begin{proof}
	We first must prove that there is a simultaneous splitting of $\alpha_X$ and $\pi_2$ by $\langle 1_G, \id_X \rangle$. However, on one hand this follows quickly from the fact that $\pi_2 \circ \langle 1_G, \id_X \rangle = \id_X$ while on the other hand it comes from the fact that $X$ is a $G$-object. Thus the identity point of $G$ acts as the identity on $X$, i.e., $\id_X = \alpha_X \circ \langle 1_G, \id_X \rangle$. Thus the diagram
	\[
	\begin{tikzcd}
	X \ar[ddr, equals] & G \times X \ar[r]{}{\alpha_X} \ar[l, swap]{}{\pi_2} & X \ar[ddl, equals] \\ \\
	& X \ar[uu]{}[description]{\langle 1_G, \id_X \rangle}
	\end{tikzcd}
	\]
	commutes, giving that $\langle 1_G, \id_X \rangle$ is simultaneously splits $\pi_2$ and $\alpha_X$. 
	
	We now consider the pullback
	\[
	\xymatrix{
		(G \times X) \times_{X} (G \times X) \pullbackcorner \ar[r]^-{p_1} \ar[d]_{p_2} & G \times X \ar[d]^{\alpha_X} \\
		G \times X \ar[r]_{\pi_2} & X
	}
	\]
	and note that this pullback is isomorphic to the diagram:
	\[
	\xymatrix{
		\pullbackcorner G \times G \times X \ar[rr]^-{\id_{G} \times \alpha_X} \ar[d]_{\pi_{23}} & & G \times X \ar[d]^{\pi_{2}}  \\
		G \times X \ar[rr]_{\alpha_X} & & X
	}
	\]
	in $\Cscr$. We will work with this representative of the pullback in order to define the internal category structure, as it streamlines the proof we present.
	
	We now need to define the composition on the object of morphisms of the internal category $\mathbb{G \times X}$. For this we define $m:G \times \times G \times X \to G \times X$ by
	\[
	\mu \times \id_X: G \times G \times X \to G \times X,
	\]
	where $\mu$ is the multiplication of the group. Note that because $G$ is an internal group to $\Cscr$, the multiplication on $\mathbb{G\times X}$ is associative and by construction the map $\langle 1_G, \id_X \rangle$ determines a unit for the composition of $\mathbb{G\times X}$. Thus we only need to verify that the source and target maps interact properly with the composition. In particular, we must verify that the diagrams
	\[
	\xymatrix{
		G \times G \times X \ar[rr]^-{\mu \times \id_X} \ar[d]_{\pi_{23}} & & G \times X \ar[d]^{\pi_{2}} & G \times G \times X \ar[rr]^-{\id_{X} \times \alpha_{X}} \ar[d]_{\mu \times \id_{X}} & & G \times X \ar[d]^{\alpha_X} \\
		G \times X \ar[rr]_-{\pi_2} & & X & G \times X \ar[rr]_-{\alpha_X} & & X
	}
	\]
	both commute. However, the left diagram commutes trivially and the right diagram commutes because $X$ is a left $G$-object in $\Cscr$, which proves the lemma.
\end{proof}
\begin{corollary}\label{Cor: Section 4: Internal groupoid}
	The internal category $\mathbb{G \times X}$ of Lemma \ref{Lemma: Section 4: internal category} is an internal groupoid with inverse morphism
	\[
	\inv := \langle -1 \circ \pi_1, \alpha_X \rangle: G \times X \to G \times X
	\]
	where $-1:G \to G$ is the inverse of the group.
\end{corollary}
\begin{proof}
	We first must establish that $\pi_2 \circ \inv = \alpha_X$ and that $\alpha_X \circ \inv = \pi_2$. However, for this we calculate that
	\[
	\pi_2 \circ \inv = \pi_2 \circ \langle -1 \circ \pi_1, \alpha_X \rangle = \alpha_X
	\]
	and that
	\[
	\alpha_X\circ \inv = \alpha_X \circ \langle -1 \circ \pi_1, \alpha_X \rangle = \pi_2,
	\]
	which was necessary to show.
	
	We now must verify that the diagrams 
	\[
	\xymatrix{
		G \times X \ar[d]_{\alpha_X} \ar[r]^-{\Delta} & (G \times X) \times_X (G \times X) \ar[rr]^-{\inv \times \id_{G\times X}} & & (G \times X) \times_X (G \times X) \ar[d]^{\mu \times \id_X} \\
		X \ar[rrr]_-{\langle 1_G, \id_X\rangle} & & & G \times X
	}
	\]
	and
	\[
	\xymatrix{
		G \times X \ar[d]_{\pi_2} \ar[r]^-{\Delta} & (G \times X) \times_X (G \times X) \ar[rr]^-{ \id_{G\times X}\times \inv} & & (G \times X) \times_X (G \times X) \ar[d]^{\mu \times \id_X} \\
		X \ar[rrr]_-{\langle 1_G, \id_X\rangle} & & & G \times X
	}
	\]
	both commute. However, upon converting this to the formulation of the pullback we use we find that the morphism
	\[
	(\inv \times \id_{G\times X}) \circ \Delta
	\]
	corresponds to the morphism
	\[
	\xymatrix{
		G \times X \ar[drrr] \ar[rrr]^-{\langle -1\circ \pi_1, \id_G\circ \pi_1, \id_X \rangle} & & & G \times G \times X \ar[d]^{\langle \id_{G\times G} \circ \pi_{12}, \alpha_X \circ \pi_{23}\rangle} \\
		& & & G \times G \times X
	}
	\]
	and the morphism
	\[
	(\id_{G \times X} \times \inv) \circ \Delta
	\]
	corresponds to the morphism
	\[
	\xymatrix{
		G \times X \ar[drrr] \ar[rrr]^-{\langle \id_G\circ \pi_1, -1\circ \pi_1, \id_X \rangle} & & & G \times G \times X \ar[d]^{\id_{G\times G \times X}} \\
		& & & G \times G \times X
	}
	\]
	However, with these constructions we find that the diagrams
	\[
	\begin{tikzcd}
	& & G \times G \times X \ar[dr]{}{\langle \id_{G\times G} \circ \pi_{12}, \alpha_X \circ \pi_{23}\rangle} \\
	G \times X \ar[d, swap]{}{\alpha_X} \ar[urr]{}{\langle \id_G\circ \pi_1, -1\circ \pi_1, \id_X \rangle} \ar[rrr]{}{(\inv \times \id_{G\times X}) \circ \Delta} & & & G \times G \times X \ar[d]{}{\mu \times \id_X} \\
	X \ar[rrr]{}{\langle 1_G, \id_X \rangle} & & & G \times X
	\end{tikzcd}
	\]
	and
	\[
	\begin{tikzcd}
	& & G \times G \times X \ar[dr]{}{\langle \id_{G\times G} \circ \pi_{12}, \alpha_X \circ \pi_{23}\rangle} \\
	G \times X \ar[d, swap]{}{\alpha_X} \ar[urr]{}{\langle \id_G\circ \pi_1, -1\circ \pi_1, \id_X \rangle} \ar[rrr]{}{(\inv \times \id_{G\times X}) \circ \Delta} & & & G \times G \times X \ar[d]{}{\mu \times \id_X} \\
	X \ar[rrr]{}{\langle 1_G, \id_X \rangle} & & & G \times X
	\end{tikzcd}
	\]
	both commute, which proves that $\mathbb{G \times X}$ is an internal groupoid to $\Cscr$.
\end{proof}
\begin{corollary}\label{Cor: Appendix A: Interal action gropoid}
	For any algebraic group $G$ and any left $G$-variety $X$, $\mathbb{G\times X}$ is a groupoid internal to $\Var_{/\Spec K}$.
\end{corollary}
\begin{example}\index{Internal Orbit Groupoid}
	When $\Cscr$ is $\Set$, so $G$ is a group in the usual sense and $X$ is a left $G$-set, the internal groupoid $\mathbb{G\times X}$ is a small groupoid where the elements of $X$ are the objects and there is a morphism $x \to y$ in $\mathbb{G\times X}$ if and only if there exists a $g \in G$ for which $y = gx$. In this way, we can think of the internal groupoid $\mathbb{G\times X}$ as an orbit groupoid for the action of $G$ on $X$.
\end{example}
\begin{remark}
For a $1$-category $\Cscr$, we will write $\Cat(\Cscr)$ for the category of internal categories of $\Cscr$. Objects in this category are internal categories and morphisms are internal functors (cf.\@ Remark \ref{Remark: Internal functor} for a construction of identities and composition).
\end{remark}

Our goal now is to build the quotient stack $[X/G]$ from the internal groupoid $\mathbb{G \times X}$. For this we first follow \cite{ElephantVol1}. To the internal category $\mathbb{G\times X}$ we associate the pseudofunctor $(\mathbb{G \times X})(-):\Sch_{/\Spec K}^{\op} \to \fCat$ defined below. Throughout this construction we let $U$ and $V$ be schemes over $\Spec K$.
\begin{itemize}
	\item For any scheme $U$ over $\Spec K$ we define the category $(\mathbb{G \times X})(U)$ by setting:
	\begin{itemize}
		\item The objects are given by $(\mathbb{G \times X})(U)_0 := \Sch_{/\Spec K}(U,X)$.
		\item The morphisms are given by $(\mathbb{G \times X})(U)_1 := \Sch_{/\Spec K}(U, G \times X)$.
		\item To define composition, we recall from Lemma \ref{Lemma: Section 4: internal category} that the morphism $\pi_2:G \times X \to X$ is the domain morphism and $\alpha_X:G \times X \to X$ is the codomain morphism of the internal category $\mathbb{G\times X}$. Thus two morphisms $g_1, g_2 \in (\mathbb{G \times X})(U)_1$ are composable in $(\mathbb{G \times X})(U)$ whenever there is a commuting square of $K$-schemes:
		\[
		\xymatrix{
			U \ar[d]_-{g_1} \ar[r]^-{g_2} & G \times X \ar[d]^{\pi_2} \\
			G \times X \ar[r]_-{\alpha_X} & X
		}
		\]
		Thus we define the composition $g_2 \bullet g_1$ in $(\mathbb{G \times X})(U)$ to be given by $(\mu \times \id_X) \circ c_{g_1,g_2}$, where $c_{g_1,g_2}$ is the unique morphism in the pullback diagram:
		\[
		\xymatrix{
			U \ar@/^/[drrr]^{g_2} \ar@/_/[ddr]_{g_1} \ar@{-->}[dr]^{\exists!c_{g_1,g_2}} & & \\
			& G \times G \times X \ar[d]_{\pi_{23}} \ar[rr]_-{\id_G \times \alpha_X} & & G \times X \ar[d]^{\pi_2} \\
			& G \times X \ar[rr]_{\alpha_X} & & X
		}
		\]
		\item For any object $x:U \to X$ in $(\mathbb{G \times X})(U)$, we define the identity $\id_{x}$ by the commuting diagram:
		\[
		\xymatrix{
			U \ar[dr]_{\id_{x}} \ar[r]^{x} & X \ar[d]^{\langle 1_G, \id_X \rangle} \\
			& G \times X
		}
		\]
	\end{itemize}
	\item We define the fibre functors $f^{\ast}:(\mathbb{G \times X})(V) \to (\mathbb{G \times X})(U)$, for any morphism $f \in \Sch_{/\Spec K}(V,U)$, by setting
	\[
	f^{\ast}(x:U \to X) := x \circ f:V \to X
	\]
	on objects $x:U \to X$ of $(\mathbb{G \times X})(U)$ and
	\[
	f^{\ast}(r:U \to G \times X) := r \circ f:V \to G \times X
	\]
	for any morphisms $r:U \to G \times X$ of $(\mathbb{G\times X})(U)$.
	\item The functoriality of the above fibre functors is strict, so we set all our coherence isomorphisms $\phi_{g_1,g_2}$ to be the identity natural transformations.
\end{itemize}
\begin{proposition}
	The pseudofunctor $(\mathbb{G \times X})(-)$ is fibred in groupoids, i.e., there is a factorization
	\[
	\xymatrix{
		\Sch_{/\Spec K}^{\op} \ar[rr]^-{(\mathbb{G \times X})(-)} \ar[dr]_{(\mathbb{G\times X})(-)} & & \fCat \\
		& \mathfrak{Gpd} \ar[ur]_{i}
	}
	\]
	of pseudofunctors where $i:\mathfrak{Gpd} \to \fCat$ is the inclusion of the $2$-category of groupoids into the $2$-category of categories.
\end{proposition}
\begin{proof}
	This is a routine verification that follows from the fact that $\mathbb{G\times X}$ is an internal groupoid in $\Sch_{/\Spec K}$ by Corollary \ref{Cor: Section 4: Internal groupoid}.
\end{proof}
Unfortunately, the pseudofunctor $(\mathbb{G \times X})(-)$ above is not a stack; it is only a prestack in general. However, what we will sketch here is an equivalence of fppf stacks
\[
(\mathbb{G\times X})(-)^{+} \simeq [X/G].
\]
\begin{definition}\label{Defn: Quot Stack}\index{Quotient Stack $[X/G]$}\index[notation]{XoverG$[X/G]$}
	For this we recall that the stack $[X/G]$ fibred in groupoids (in its incarnation as a pseudofunctor) is defined as follows (cf.\@ \cite{Neumann}):
	\begin{itemize}
		\item For any scheme $U$ over $\Spec K$, $[X/G](U)$ is the category of spans
		\[
		\xymatrix{
			T \ar[d]_{x} \ar[r]^{t} & U \\
			X
		}
		\]
		where $t$ admits fppf-local sections (cf.\@ Remark \ref{Remark: Local section in topology}) and $T$ is a $G$-torsor for which the morphism $x \in \GSch(T,X)$. Morphisms in $[X/G](U)$ are given by morphisms between $G$-torsors $f:T \to T^{\prime}$ that commute with the equivariant maps $x:T \to X$ and $x^{\prime}:T^{\prime} \to X$.
		\item For any scheme morphism $f^{\ast}:U \to U^{\prime}$ the fibre functor $f^{\ast}:[X/G](U^{\prime}) \to [X/G](U)$ is induced by the pullback of a $G$-torsor $T$ equipped with a morphism $t:T \to U^{\prime}$ against the morphism $f:U \to U^{\prime}$. 
		\item The coherence isomorphisms are induced by the universal property of the pullbacks above.
	\end{itemize}
\end{definition}
It is well-known (cf.\@ \cite{Neumann}, pages 51 and 52) that $[X/G]$ is an {\'e}tale stack on $X$ when $G$ is affine, and the case when $G$ is quasi-affine follows from this (as this implies that the maps $\pi_2,\alpha_X:G \times X \to X$ are quasi-affine and separated when $X$ is a variety). However, for general group schemes and non-quasi-affine algebraic groups, we only know in general that $[X/G]$ is an fppf stack which is a finite type algebraic stack when $G$ and $X$ are both varieties (cf.\@ \cite{Illusie}, \cite[Example 8.1.12]{OlssonSpacesStacks}, \cite[Pages 37, 38, Comment ii.a]{Champs}). However, we will present our stackification results in the fppf case; while the results are formal and follow mutatis mutandis from one site to the other, what is cruicial is being able to do effective descent along torsors, and this is always possible when $G$ is a smooth scheme and $X$ is a scheme in the fppf topology. We will sketch how to see the stackification of the prestack $(\mathbb{G\times X})(-)$ compares to the stack $[X/G]$ with a light touch on the (Grothendieck) topology used, save for occasionally mentioning that these are fppf stacks and fppf torsors. Note that in Remark \ref{Remark: Local section in topology} and Definition \ref{Defn: Groupoid Stackification} we recall what it means for a morphism to admit local sections in a pretopology and what the stackification of $(\mathbb{G \times X})(-)$ looks like. The reader already familiar with these notions can skip ahead to Proposition \ref{Prop: Section 4: Equivalence of quotient stack and symmetry stack}.

\begin{remark}\label{Remark: Local section in topology}\index{Local Sections in a Site}
	Recall that if $(\Cscr, J)$ is a site where the topology $J$ is induced by a pretopology $\tau$, we say that a morphism $f:X \to Y$ in $\Cscr$ admits $J$-local sections if there exists a $\tau$-covering family $\lbrace g_i:Y_i \to Y \; | \; i \in I \rbrace$ of $Y$ for which there exist morphism $f_i:Y_i \to X$ making the diagram
	\[
	\xymatrix{
		X \ar[r]^{f} & Y \\
		& Y_i \ar[u]_{g_i} \ar@{-->}[ul]^{\exists f_i}
	}
	\]
	commutative in $\Cscr$.
\end{remark}
\begin{definition}[\cite{AlgStacks}]\label{Defn: Section 4: Banal groupoid}\index{Banal Groupoid}
	Let $\Cscr$ be a category and let $f \in \Cscr(Y,X)$. We say that an internal groupoid 
	\[
	\mathbb{G} = (G_0, G_1, \id_{\mathbb{G}}, s, t, \inv, m)
	\] 
	in $\Cscr$ is a banal groupoid for $f$ if $G_0 = Y$, if the diagram
	\[
	\xymatrix{
		G_1 \ar[r]^-{s} \ar[d]_{t} & Y \ar[d]^{f} \\
		Y \ar[r]_-{f} & X
	}
	\]
	commutes, and if $\langle s, t \rangle: G_1 \to Y \times_X Y$ is an isomorphism.
\end{definition}

We now need to recall internal functors between internal categories, as this will let us talk about morphisms of groupoid schemes. In particular, this will allow us to also describe square morphisms of internal groupoids, which are particularly nice morphisms of internal groupoids.
\begin{remark}\label{Remark: Internal functor}\index{Internal Functor}
	Let 
	\[
	\Cbb = (C_0, C_1, \dom^{\Cbb}, \codom^{\Cbb}, s_0^{\Cbb},m^{\Cbb})
	\]
	and
	\[
	\Dbb = (D_0, D_1, \dom^{\Dbb}, \codom^{\Dbb},s_0^{\Dbb},m^{\Dbb})
	\]
	be internal categories in a $1$-category $\Cscr$. Recall that an internal functor $\underline{f}:\Cbb \to \Dbb$ is a pair of morphisms $\underline{f} = (f_0,f_1)$ where $f_0 \in \Cscr(C_0,D_0)$ and $f_1 \in \Cscr(C_1,D_1)$ such that the following diagrams commute:
	\begin{itemize}
		\item The diagram
		\[
		\xymatrix{
			C_0 \ar[r]^-{f_0} \ar[d]_{s_0^{\Cbb}} & D_0 \ar[d]^{s_0^{\Dbb}} \\
			C_1 \ar[r]_-{f_1} & D_1
		}
		\]
		commutes, which says that $\underline{f}$ preserves identity morphisms;
		\item The diagrams
		\[
		\xymatrix{
			C_1 \ar[d]_{\dom^{\Cbb}} \ar[r]^-{f_1} & D_1 \ar[d]^{\dom^{\Dbb}} & & C_1 \ar[r]^-{f_1} \ar[d]_{\codom^{\Cbb}} & D_1 \ar[d]^{\codom^{\Dbb}} \\
			C_0 \ar[r]_-{f_0} & D_0 & & C_0 \ar[r]_-{f_0} & D_0
		}
		\]
		commute, which together say that $\underline{f}$ preserves all relevant domains and codomains;
		\item The diagram
		\[
		\xymatrix{
			C_1 \times_{C_0} C_1 \ar[rr]^-{f_1 \times f_1} \ar[d]_{m^{\Cbb}} & & D_1 \times_{D_0} D_1 \ar[d]^{m^{\Dbb}} \\
			C_1 \ar[rr]_-{f_1} & & D_1
		}
		\]
		commutes, which says that $\underline{f}$ preserves composition.
	\end{itemize}
	Note that technically speaking we should show that the unique morphism
	\[
	\theta:C_1 \times_{C_0} C_1 \to D_1 \times_{D_0} D_1
	\]
	induced by the internal functor $\underline{f}$ is indeed the product map $f_1 \times f_1$, but this is reasonably immediate from the cube below:
	\[
	\begin{tikzcd}
	C_1 \times_{C_0} C_1 & & C_1 & \\
	& D_1 \times_{D_0} D_1 & & D_1 \\
	C_1 & & C_0 & \\
	& D_1 & & D_0
	\arrow[from = 1-1, to = 1-3]{}{p_1}
	\arrow[from = 1-1, to = 2-2, dashed]{}{\exists!\theta}
	\arrow[from = 1-1, to = 3-1,swap]{}{p_2}
	\arrow[from = 3-1, to = 3-3, near start]{}{\dom^{\Cbb}}
	\arrow[from = 3-1, to = 4-2, swap]{}{f_1}
	\arrow[from = 1-3, to = 2-4]{}{f_1}
	\arrow[from = 1-3, to = 3-3, near end,swap]{}{\codom^{\Cbb}}
	\arrow[from = 3-3, to = 4-4]{}{f_0}
	\arrow[from = 4-2, to = 4-4, swap]{}{\dom^{\Dbb}}
	\arrow[from = 2-4, to = 4-4]{}{\dom^{\Dbb}}
	\arrow[from = 2-2, to = 4-2, crossing over, near end]{}{q_2}
	\arrow[from = 2-2, to = 2-4, crossing over, near start]{}{q_1}
	\end{tikzcd}
	\]
	
	The composition of internal functors is done pointwise, i.e., 
	\[
	\underline{g} \circ \underline{f} = ( g_0 \circ f_0,g_1 \circ f_1)
	\]
	and the identity internal functor is $\underline{\id}_{\Cbb} = (\id_{C_0}, \id_{C_1})$.
\end{remark}
\begin{definition}[\cite{AlgStacks}]\label{Defn: Section 4: Square Morphism}]\index{Square Morphism of Groupoids}
	Let $\Gbb$ and $\Hbb$ be internal groupoids and let $\underline{f}:\Gbb \to \Hbb$ be an internal functor. We say that $\underline{f}$ is a square morphism of groupoids if the diagrams
	\[
	\xymatrix{
		G_1 \ar[r]^{f_1} \pullbackcorner \ar[d]_{\dom^{\Gbb}} & H_1 \ar[d]^{\dom^{\Hbb}} \\
		G_0 \ar[r]_{f_0} & H_0
	}
	\]
	and
	\[
	\xymatrix{
		G_1 \ar[r]^-{f_1} \ar[d]_{\codom^{\Gbb}} \pullbackcorner & H_1 \ar[d]^{\codom^{\Hbb}} \\
		G_0 \ar[r]_-{f_0} & H_0
	}
	\]
	are both pullbacks in $\Cscr$.
\end{definition}
\begin{definition}[\cite{AlgStacks}]\label{Defn: Groupoid Torsors}\index{Groupoid Torsor}
	Let $(\Cscr,J)$ be a site induced by a pretopology $\tau$ and let $\Gbb = (G_0, G_1, \cdots)$ be a groupoid internal to $\Cscr$. Then a $\Gbb$-torsor in $J$ is given by the following three pieces of information:
	\begin{itemize}
		\item A morphism $g:T^{\prime} \to T$ in $\Cscr$, a banal groupoid $\mathbb{B} = (T^{\prime},B_1,\cdots)$ for $g$, and a square morphism of groupoids $\underline{f}:\mathbb{B} \to \Gbb$ for which $g$ admits $J$-local sections.
	\end{itemize}
\end{definition}
\begin{definition}[\cite{AlgStacks}]\label{Defn: Groupoid Stackification}
	Let $(\Cscr,J)$ be a site induced by a pretopology $\tau$ and let $\Gbb = (G_0, G_1, \cdots)$ be a groupoid internal to $\Cscr$. Then the category of $\Gbb$-torsors in $J$, $[\Gbb]$,\index[notation]{Gtorsors@$[\Gbb]$} is defined as follows:
	\begin{itemize}
		\item Objects: $\Gbb$-torsors in $(\Cscr,J)$.
		\item Morphisms: Given $\Gbb$-torsors in $(\Cscr,J)$ of the form $T= (g,\mathbb{B}^g,\underline{f}^g)$ with $g:T^{\prime} \to T$ and $S = (h,\mathbb{B}^h,\underline{f}^h)$ with $h:S^{\prime} \to S$, a morphism $\alpha \in [\Gbb](S,T)$ is given as follows:
		\begin{itemize}
			\item A pair of morphisms $\rho:S \to T$ and $\psi:S^{\prime} \to T^{\prime}$ for which the diagram
			\[
			\xymatrix{
				S^{\prime} \ar[r]^-{\psi} \ar[d]_{h} & T^{\prime} \ar[d]^{g} \\
				S \ar[r]_-{\rho} & T
			}
			\]
			commutes. Note that this commutativity induces a unique morphism $\psi_1:B^h_1 \to B^g_1$ making the pair $\underline{\psi} = (\psi,\psi_1):\mathbb{B}^h \to \mathbb{B}^g$ into an internal functor in $\Cat(\Cscr)$ of internal categories in $\Cscr$.
			\item We require that the diagram
			\[
			\xymatrix{
				\mathbb{B}^h \ar[r]^-{\underline{\psi}} \ar[dr]_{\underline{f}^h} & \mathbb{B}^g \ar[d]^{\underline{f}^g} \\
				& \Gbb
			}
			\]
			commutes in the category $\Cat(\Cscr)$.
		\end{itemize}
		\item Composition: Pairwise composition, i.e., $(\rho, \psi) \circ (\rho^{\prime}, \psi^{\prime}) = (\rho\circ \rho^{\prime}, \psi \circ \psi^{\prime})$.
		\item Identities: $(\id_T, \id_{T^{\prime}})$ is the identitiy on $(g, \mathbb{B}, \underline{f})$ with $g:T^{\prime} \to T$.
	\end{itemize}
	When $\Gbb$ is the orbit groupoid of a smooth algebraic group $H$ acting on a $K$-variety $X$, we call\index{Stackification}
	\[
	[\Gbb] = [\Hbb \times \Xbb]
	\]
	the fppf-stackification of $\Gbb$ and write $[\Gbb] = (\Hbb \times \Xbb)^{+}$. Note also that the pseudofunctor description of $[\Gbb]$ gives the fibre categories $[\Gbb](U)$ as those $\Gbb$-torsors $(g,\mathbb{B},\underline{f})$ for which $\Codom(g) = U$ and induces the fibre functors by choosing a cleavage for the fibration$[\Gbb] \to \Cscr$ given by sending objects $(g,\mathbb{B},\underline{f})$ to $\Codom(g)$ and sending morphisms $(\psi,\psi_1)$ to $\psi$. Note that this is simply running the Grothendieck construction to describe $[\Gbb]$ as a pseudofunctor.
\end{definition}
\begin{remark}
	It follows from the theory of fppf descent (cf.\@ \cite[Pages 37, 38, Comment ii.a]{Champs} and \cite[Example 8.1.12]{OlssonSpacesStacks}) that when the internal groupoid $\Cbb = (\mathbb{G \times X})$ is the action groupoid of $G$ on $X$ for a smooth algebraic group $G$ and a variety $X$, the category $[\mathbb{G \times X}]$ is a finitely presented algebraic stack over $\Spec K$. This justifies our label of $[\Gbb \times \Xbb]$ as the fppf-stackification of $\Gbb \times \Xbb$.
\end{remark}

\begin{proposition}\label{Prop: Section 4: Equivalence of quotient stack and symmetry stack}
	There is an equivalence of stacks
	\[
	(\mathbb{G\times X})(-)^{+} \cong [X/G].
	\]\index[notation]{GtimesXplus@$(\Gbb\times\Xbb)^{+}$}
	where $(\mathbb{G \times X})(-)^{+}$ is the stackification of $(\mathbb{G \times X})(-)$.
\end{proposition}
\begin{proof}
	This is a routine unwinding of the definitions. In particular, because the morphism condition is the same between both categories, it suffices to check that groupoid torsors for $(\Gbb \times \Xbb)(-)^{+}$ form torsors (up to isomorphism) in $[X/G]$ and vice versa. For this we fix a scheme $U$ in $\Sch_{/\Spec K}$. Fix an arbitrary object $(T \xrightarrow{t} U, T \xrightarrow{x} X)$ in $[X/G](U)$. Then we note $T$ is a $G$-torsor with a span
	\[
	\xymatrix{
		T \ar[r]^-{t} \ar[d]_{x} & U \\
		X
	}
	\]
	where $t$ admits fppf-local sections and $x$ is $G$-equivariant. By taking the pullback of $T$ with itself we get a banal groupoid for $t$:
	\[
	\Tbb = \begin{cases}
	B_0 = T; \\ 
	B_1 =T \times_U T; \\
	\codom^{\Tbb} = \pi_1:T \times_U T \to T; \\
	\dom^{\Tbb} \pi_2:T \times_U T \to T;\\
	\mu = \pi_{13}:T \times_U T \times_U T \to T \times_U T;\\
	s_0 = \Delta:T \to T \times_U T; \\ 
	\inv = s:T \times_U T \to T \times_U T;
	\end{cases}
	\]
	where $s$ is the switching map $s:X \times Y \to Y \times X$. Let $\sigma:T \times_U T \to G$ be the morphism for which
	\[
	\xymatrix{
		T \times_U T \ar[rr]^-{\langle \sigma, \pi_2\rangle} & & G \times T
	}
	\] 
	is an isomorphism of schemes asserted by the fact that $T$ is a torsor. We now claim that $\underline{t} = (t, \langle \sigma, t \circ \pi_2\rangle)$ determines an internal functor $\Tbb \to \mathbb{G \times X}$. We show that $\underline{t}$ commutes with the domain and codomain morphisms, as this implies that it commutes with the multiplication and identities in turn. Note that the diagram
	\[
	\xymatrix{
		T \times_U T \ar[d]_{\langle \sigma, t\circ \pi_2\rangle} \ar[r]^-{\pi_1} & T \ar[d]^{t} \\
		G \times X \ar[r]_-{\alpha_X} & X
	}
	\]
	commutes by the equivariance of $t$ while
	\[
	\xymatrix{
		T \times_U T \ar[d]_{\langle \sigma, t \circ \pi_2\rangle} \ar[r]^-{\pi_2} & T \ar[d]^{t} \\
		G \times X \ar[r]_{\pi_2} & X
	}
	\]
	commutes by construction. Thus $\underline{t}$ is an internal functor; in fact, it is a routine observation to note that the two diagrams above are in fact pullbacks, so $\underline{t}$ is a square morphism of groupoids. This implies that $(t:T \to U, \Tbb,\underline{t})$ is in fact a ($\mathbb{G \times X}$)-torsor and provides a fully faithful functor $[X/G](U) \to (\Gbb \times \Xbb)(U)$.
	
	We now show that the functor $[X/G](U) \to (\Gbb \times \Xbb)(U)$ is fully faithful, as because this construction is pseudofunctorial in $U$ this will give rise to an isomorphism of fibrations and hence stacks. However, the fully faithfulness is almost immediate; if $(t,\mathbb{B},\underline{f})$ is an aribtrary $(\Gbb \times \Xbb)$ torsor in $(\Gbb \times \Xbb)(U)$, then $t$ is a morphism of schemes $t:T \to U$ for some $T$ for which $B_1 \cong T \times_U T$. However, this isomorphism induces an isomorphism of banal groupoids $\mathbb{B} \to \Tbb$, which makes the diagram
	\[
	\xymatrix{
		\mathbb{B} \ar[dr]_-{\underline{f}} \ar[r]^-{\cong} & \Tbb \ar[d]^{\underline{t}} \\
		& \Gbb \times \Xbb
	}
	\]
	commute in $\Cat(\Sch_{/\Spec K})$. Thus we get that $(g,\mathbb{B},\underline{f}) \cong (g,\Tbb,\underline{t})$ and so $[X/G](U) \simeq (\Gbb \times \Xbb)(U)$. This gives the isomorphism of stacks
	\[
	(\Gbb \times \Xbb)(-)^{+} \cong [X/G],
	\]
	as was desired.
\end{proof}
\chapter{Locally Closed Subtoposes of Grothendieck Toposes}\label{Appendix: Locally closed subtoposes}

We now recall from \cite[Expos{\'e} IV.6 and IV.7]{SGA4} some terminology and results involving Grothendieck toposes and the canonical topology on the isomorphism-classes of points of a topos. The reason why we will be discussing this is to present the core ideas that allow us to discuss stratifications of Grothendieck toposes; this will be necessary to define the derived category on $[X/G]$.

\begin{definition}[{\cite[Expos{\'e} IV.6]{SGA4}}]\index{Punctual Topos}\index{Terminal Topos|see {Punctual Topos}}
	The punctual topos is the topos
	\[
	P := \Shv(\ast) \cong \Set.
	\]
\end{definition}
\begin{remark}
	Note that $P$ is the terminal object in $\mathbf{GrTopos}$, the $1$-category of Grothendieck toposes. 
\end{remark}
\begin{definition}[{\cite[Definition IV.6.1]{SGA4}}]\index{Category of Points}\index[notation]{Toposfrak@$\fTopos$}
	Let $\Ecal$ be a Grothendieck topos. Then define the category
	\[
	\Point(\Ecal) := \fTopos(P,\Ecal)
	\]
	to be the category of geometric morphisms from the punctual topos $P$ to $\Ecal$. Note that $\fTopos$ is the $2$-category of toposes with $1$-morphisms geometric morphisms and $2$-morphisms geometric transformations.
\end{definition}
\begin{remark}[Open Subtoposes]\index{Subtopos! Open}
	Recall that if $\Ecal$ is a Grothendieck topos, then a subtopos $\Ucal$ of $\Ecal$ is described by a geoemtric morphism $j:\Ucal \to \Ecal$ for which the pushforward $j_{\ast}:\Ucal \to \Ecal$ is fully faithful. An argument given in \cite[\href{https://stacks.math.columbia.edu/tag/08LW}{Tag 08LW}]{stacks-project} explains that for a sheaf topos $\Shv(\Cscr,J)$, a slice topos $\Shv(\Cscr,J)_{/\Fscr}$ is a subtopos of $\Shv(\Cscr,J)$ if and only if the sheaf $\Fscr$ is a subterminal object, i.e., there is a monic $\Fscr \to \top_{\Shv(\Cscr,j)}.$ We define such toposes to be open subtoposes of $\Shv(\Cscr,J)$; more generally, we say that a subtopos $\hat{j}:\Ucal\to\Ecal$ of a Grothendieck topos $\Ecal$ is an open subtopos if there is a sheaf topos $\Shv(\Cscr,J)$ and a subterminal object $\Fscr$ for which the diagram
	\[
	\xymatrix{
		\Shv(\Cscr,J)_{/\Fscr} \ar[rr]^-{j} & & \Shv(\Cscr,J) \\
		\Ucal \ar[rr]_-{\hat{j}} \ar[u]^{\simeq} & & \Ecal \ar[u]_{\simeq}
	}
	\]
	commutes in $\Topos$.
\end{remark}
\begin{remark}\index[notation]{UcalP@$\Ucal_p$}
	Let $\Ecal$ be a Grothendieck topos and let $\Fscr \in \Ecal_0$. If $p \in \Point(\Ecal)$, i.e., for a geometric morphism $p:P \to \Ecal$, we write $\Fscr_p$ in place of the object $p^{\ast}\Fscr$. If $\Ucal \simeq \Shv(\Cscr,J)_{/\Fscr}$ is an open subtopos of $\Shv(\Cscr,J)$, we will sometimes write
	$\Ucal_p$ in place of $\Fscr_p$ as well, depending on context (cf.\@ Exercise IV.7.8.a).
\end{remark}

We now give the crucial construction which will allow us to describe stratifications of Grothendieck toposes. For this we follow \cite[Exercise IV.7.8.a]{SGA4}. Fix a Grothendieck topos $\Ecal$ and write
\[
\lvert \Ecal \rvert := {\Point(\Ecal)}_{/\text{iso}},
\]\index[notation]{EcalAssociatedSpace@$\lvert\Ecal\rvert$}
i.e., the set of points of $\Ecal$ modulo isomorphism. Now fix an open subtopos of $\Ecal$ and sheaf representation $\Ecal \simeq \Shv(\Cscr,J)$ and $\Ucal \simeq \Shv(\Cscr,J)_{/\Fscr}$. Define now the set
\[
\lvert \Ucal \rvert := \lbrace [p] \in \lvert \Ecal \rvert \; : \; \Ucal_p \ne \emptyset \rbrace,
\]\index[notation]{UcalAssociatedSpace@$\lvert\Ecal\rvert$}
i.e., the set of points $p$ of $\Ecal$ which do not send the open subtopos $\Ucal$ to the empty set. We then determine a topology on $\lvert \Ecal \rvert$ by setting
\[
\Open(\lvert \Ecal \rvert) := \lbrace \lvert \Ucal \rvert \; : \; \Ucal\, \text{is\,an\,open\,subtopos\,of\,} \Ecal\rbrace.
\]
This space has the property that it induces an isomorphism of posets
\[
\Open(\Ecal) \cong \Open(\lvert \Ecal \rvert),
\]
where the poset on the left is the poset of open subobject lattice of $\Ecal$, i.e.,
\[
\Open(\Ecal) = \mathbf{Sub}_{\text{Open}}(\Ecal).
\]
\begin{definition}[{\cite[Exercise IV.7.8.a]{SGA4}}]\index{Topos! Associated Topological Space}
	Let $\Ecal$ be a Grothendieck topos. We then define the topological space associated to $\Ecal$ to be the space $
	(\lvert \Ecal \rvert, \Open(\lvert \Ecal \rvert))$
	constructed above.
\end{definition}

We now describe closed subtoposes of a Grothendieck topos. While these are of parameterized by the closed subsets of $\lvert \Ecal \rvert,$ we would like to give a more concrete description of what it means to be the complement of an open subtopos $\Ucal$. For this find a sheaf representation $\Ecal \simeq \Shv(\Cscr,J)$ and find an open subtopos $\Ucal = \Shv(\Cscr,J)_{/\Fscr}$. We now want to define a subtopos $\Ccal_{\Ucal}$ of $\Shv(\Cscr,J)$ which acts as a complement to $\Ucal$ and hence determines a closed subtopos of $\Ecal$. For this we follow \cite[Expos{\'e} IV.9.3]{SGA4}. Fix $j:\Shv(\Cscr,J)_{/\Fscr} \to \Shv(\Cscr,J)$ the immersion of toposes. We define a functor $C_{\Ucal}:\Ecal \to \Ecal$ which we will use to build our closed complement of $\Ucal$. For this, for an arbitrary object $\Gscr \in \Ecal_0$, cosider the span
\[
\xymatrix{
\Fscr & \Fscr \times \Gscr \ar[r]^-{\pi_2} \ar[l]_-{\pi_1} & \Gscr
}
\]
in $\Ecal$ and note that such spans vary functorially in $\Gscr$. Taking the pushout of the diagram as in
\[
\xymatrix{
\Fscr \times \Gscr \ar[r]^-{\pi_2} \ar[d]_{\pi_1} & \Gscr \ar[d] \\
\Fscr \ar[r] & \Fscr \coprod\limits_{\Fscr \times \Gscr} \Gscr \pushoutcorner
}
\]
gives rise to a functor $C_{\Ucal}:\Ecal \to \Ecal$ given by sending an object $\Gscr$ to the pushout given above. Note that the functor $C_{\Ucal}$ admits the following alternate description: Because taking the spans $\Fscr \xleftarrow{\pi_1} \Fscr \times \Gscr \xrightarrow{\pi_2} \Gscr$ is functorial in $\Gscr$, it produces a functor $\Ecal \to [\Lambda_0^2,\Ecal]$ which sends an object $\Gscr$ to the span $\Fscr \xleftarrow{\pi_1} \Fscr \times \Gscr \xrightarrow{\pi_2} \Gscr$. Call this functor $F:\Ecal \to [\Lambda_0^2,\Ecal]$. Note here that $\Lambda_0^2$ is the category
\[
\xymatrix{
 & 1 \\
0 \ar[ur] \ar[rr] & & 2
}
\]
with the identity arrows left implicit; it is the category whose simplicial nerve coincides with the simplicial set $\Lambda_0^2,$ and hence justifies our naming convention. Because $\Ecal$ is cocomplete (as $\Ecal$ is a Grothendieck topos), there is a colimit functor
\[
\colim:[\Lambda_0^2,\Ecal] \to \Ecal
\] 
which sends a span in $[\Lambda_0^2,\Ecal]$ into the category $\Ecal$. The functor $C_{\Ucal}$ is then the composite
\[
\xymatrix{
\Ecal \ar[dr]_{F} \ar[rr]^-{C_{\Ucal}} & & \Ecal \\
 & [\Lambda_0^2,\Ecal] \ar[ur]_{\colim}
}
\]
in $\Cat$. Now, because $\Ecal$ is a Grothendieck topos and $C_{\Ucal}$ is cocontinuous, by the Adjoint Functor Theorem the functor $C_{\Ucal}$ is a left adjoint with right adjoint $R: \Ecal \to \Ecal$. Now define $\Ccal_{\Ucal}$ to be the category given as follows:
\begin{itemize}
	\item Objects: $\Gscr \in (\Ccal_{\Ucal})_0$ if and only if $\Gscr \in \Ecal_0$ and the projection $\Fscr \times \Gscr \xrightarrow{\pi_1} \Fscr$ is an isomorphism in $\Ecal$.
	\item Morphisms: As in $\Ecal$.
	\item Composition and identities: As in $\Ecal$.
\end{itemize}
\begin{proposition}[{\cite[Proposition IV.9.3.4]{SGA4}}]
The category $\Ccal_{\Ucal}$ is a subtopos of $\Ecal$ and there is an immersion $i:\Ccal_{\Ucal} \to \Ecal$ where $i^{\ast}:\Ecal \to \Ccal_{\Ucal}$ acts on objects by $i^{\ast}\Gscr := C_{\Ucal}(\Gscr)$ and $i_{\ast}:\Ccal_{\Ucal} \to \Ecal$ is the inclusion of categories.
\end{proposition}
\begin{definition}[{\cite[Section IV.9.3.5]{SGA4}}]\index{Subtopos! Closed}
The topos $\Ccal_{\Ucal}$ is the closed complenent of the open subtopos $\Ucal$. More generally, we say that a subtopos $\Vcal$ of $\Ecal$ is a closed subtopos if there exists an open subtopos $\Ucal$ for which $\Vcal \simeq \Ccal_{\Ucal}$.
\end{definition}
\begin{remark}
Closed subtoposes $\Vcal$ of $\Ecal$ correspond to closed subspaces of the space $\lvert \Ecal \rvert$ in exactly the same way that open subtoposes $\Ucal$ of $\Ecal$ correspond to open subspaces of $\lvert \Ecal \rvert$. This follows on one hand because in the lattice $\mathbf{Sub}(\Ecal)$, the closed complement $\Ccal_{\Ucal}$ at an open $\Ucal$ is complemented, i.e.,
\[
\Ccal_{\Ucal} \land \Ucal = \bot.
\]
\end{remark}

We now move to define locally closed subtoposes of a Grothendieck topos. These are subtoposes $\Vcal$ of $\Ecal$ which are dense in their closures, i.e., these are the subtoposes who are the meet (intersection) of an open subtopos and a closed subtopos in the subobject lattice of $\Ecal$.
\begin{definition}[{\cite[Exercise IV.9.4.9]{SGA4}}]\index{Subtopos! Locally Closed}
A subtopos $\Vcal$ of a Grothendieck topos $\Ecal$ is locally closed if it is the intersection of an open subtopos $\Ucal$ and a closed subtopos $\Fcal$, i.e.,
\[
\Vcal = \Ucal \cap \Fcal.
\]
\end{definition}

\begin{proposition}
Locally closed subtoposes of $\Ecal$ are equivalent to locally closed subspaces of $\lvert \Ecal \rvert$.
\end{proposition}
\begin{proof}
This is a routine consequence of the definition of $\lvert \Ecal \rvert$ and its topology.
\end{proof}
\chapter{Constructible Fibrations and Stratifications}\label{Appendix: Stratifications and c-structures}
This appendix seeks to provide a background for the notions of locally constructible fibrations, stratifications of Grothendieck toposes, and the theory of d0structures and cd-structures of a topological space defined in \cite{Behrend}.

We begin by briefly recalling the notion of a d-structure and a cd-structure over a topological space $X$, which is a concept essentially introduced (but not named) in \cite[Section 1.4]{BBD}; these concepts were named in \cite[Section 3.1]{Behrend}, and it is the latter exposition we follow. For the basics of manipulating these structures, see \cite[Section 3.1, 3.2]{Behrend} and \cite[Section 1.4]{BBD}. Ins this appendix we instead describe what such structures are and then move to show how they allow us to get at constructible (complexes of) sheaves.
\begin{definition}\index{Stratification! Of a Topological Space}
	A stratification of a topological space $X$ is a finite set $\Scal$ of nonempty locally closed subspaces of $X$ such that
	\[
	X = \bigsqcup_{V \in \Scal} V
	\]
	and such that for each $V \in \Scal$ there exists a subset $\Scal_V \subseteq \Scal$ such that
	\[
	\overline{V} = \bigcup_{W\in \Scal_V} W.
	\]
\end{definition}
\begin{definition}[{\cite[Definition 3.1.2]{Behrend}}]\index{d-Structure! On a Topological Space}
	Fix a topological space $X$ and let $L(X)$ be the category of locally closed subspaces of $X$. A d-structure on $X$ is then a fibration $D \to L(X)$ which satisfies the following:
	\begin{itemize}
		\item Each fibre category $D(V)$ has a $t$-structure with $D(\emptyset) = 0$ and for which each possible pullback functor $D(i_{WV}) = i_{WV}^{\ast}:D(W) \to D(V), i_{WV}:W \to V$, is exact.
		\item If $i: W\to V$ is a closed immersion in $L(X)$, then every possible pullback $i^{\ast}:D(W) \to D(V)$ fits into an adjunction diagram
		\[
		\begin{tikzcd}
		D(W) \ar[r, bend left = 40, ""{name = U}]{}{i^{\ast}} \ar[r, bend right = 40, swap, ""{name = L}]{}{i^{!}} & D(V) \ar[l, ""{name = M}]{}[description]{i_{\ast}} \ar[from = U, to = M, symbol = \dashv]{}{} \ar[from = M, to = L, symbol = \dashv]{}{}
		\end{tikzcd}
		\]
		where $i^{!}$ is exact and $i_{\ast}$ is fully faithful and $t$-exact.
		\item If $j:W \to U$ is an open immersion in $L(X)$ then every possible pullback $j^{\ast}:D(W) \to D(U)$ fits into an adjunction diagram
		\[
		\begin{tikzcd}
		D(W) \ar[r, bend left = 40, ""{name = U}]{}{j_{!}} \ar[r, bend right = 40, swap, ""{name = L}]{}{j_{\ast}} & D(U) \ar[l, ""{name = M}]{}[description]{j^{\ast}} \ar[from = U, to = M, symbol = \dashv]{}{} \ar[from = M, to = L, symbol = \dashv]{}{}
		\end{tikzcd}
		\]
		for which $j_!$ and $j_{\ast}$ are both exact and fully faithful while $j^{\ast}$ is $t$-exact.
		\item If $W = V \coprod U$ for an open immersion $j:U \to W$ and closed immersion $i:V \to W$, then $j^{\ast} \circ i_{\ast} = 0$. Moreover, for any $A \in D(V)_0$ there exist morphisms
		\[
		(i_{\ast} \circ i^{\ast})A \to (j_{!} \circ j^{\ast})A[1]
		\]
		and
		\[
		(j_{\ast} \circ j^{\ast})A \to (i_{\ast} \circ i^!)A[1]
		\]
		for which the triangles
		\[
		\xymatrix{
			(j_{!} \circ j^{\ast})A \ar[r] & A \ar[r] & (i_{\ast} \circ i^{\ast})A \ar[r] & (j_{!} \circ j^{\ast})A[1]
		}
		\]
		and
		\[
		\xymatrix{
			(i_{\ast} \circ i^{!})A \ar[r] & A \ar[r] & (j_{\ast} \circ j^{\ast})A \ar[r] & (i_{\ast} \circ i^{!})A[1]
		}	
		\]
		are distinguished in $D(V)$.
	\end{itemize}
\end{definition}
The next definition we need is that of a cd-structure on a topological space $X$. This is essentially given by d-structure $p:D \to L(X)$ together with a subfibration $q:D_{\text{lcc}} \to L(X)$ which witnesses the locally constant constructable objects in $D$.
\begin{definition}[{\cite[Definition 3.2.1]{Behrend}}]\label{Defn: Appendix Stratification: cd-structure}\index{cd-Structure! On a Topological Space}
	Fix a Noetherian topological space with category $L(X)$ of locally closed subspaces. A cd-structure on $X$ is a d-structure $p:D \to L(X)$ for which there exists a subfibration $q:D_{\text{lcc}} \to L(X)$ of $p$ which satisfies the following properties for all $V \in L(X)_0$:
	\begin{itemize}
		\item Write $A(V) := D(V)^{\heartsuit}$ and define $A_{\text{lcc}} := D(V)^{\heartsuit} \cap D_{\text{lcc}}(V)$. Then $A_{\text{lcc}}$ is a finite closed subcategory of $A(V)$.
		\item If $M \in D(V)_0$ and for all $n \in \Z$ we have that $H^n(M) \in A_{\text{lcc}}(V)_0$, then $M \in D_{\text{lcc}}(V)_0$.
		\item We have that $D_{\text{lcc}}(V) \cong D_{\text{lcc}}^{b}(V)$, i.e., for each object $M \in D_{\text{lcc}}(V)_0$, there exist $m,n \in \Z$ (wolog $m \leq n$) for which $H^{k}(M) \cong 0$ whenever $k <  m$ or $n < k$.
	\end{itemize}
\end{definition}
\begin{definition}[{\cite[Definition 3.2.2]{Behrend}}]\label{Defn: Pre-L-stratification}\index{Pre-L-Stratification! On a Topological Space}
	A pre-$L$-stratification of $X$ is a pair $(\Scal,\Lcal)$, where $\Scal$ is a stratification of $X$ and $\Lcal$ assigns to each stratum $V \in \Scal$ a finite set $\Lcal(V)$ of simple objects in the category $A_{\text{lcc}}(V)$.
\end{definition}
\begin{definition}[{\cite[Defniition 3.2.2]{Behrend}}]\index{Pre-L-Stratifications! Constructible Sheaves! with Respect to $(\Scal,\Lcal)$}\index{Pre-L-Stratifications! Constructible Sheaves}
	Fix a pre-L-stratification $(\Scal,\Lcal)$ of $X$. An object $M \in D(X)_0$ is $(\Scal,\Lcal)$-constructible if for every $V \in \Scal$ with immersion $i:V \to X$ and for every $n \in \Z$, the cohomology object $H^{n}(i^{\ast}M) \in A_{\text{lcc}}(V)_0$ and the Jordan-H{\"o}lder components of $H^{n}(i^{\ast}M)$ are all isomorphic to elements in $\Lcal(V)$.
\end{definition}
\begin{definition}[{\cite[Defniition 3.2.2]{Behrend}}]
	An object $M \in D(X)_0$ is constructible if there exists a pre-L-stratification $(\Scal,\Lcal)$ of $X$ such that $M$ is $(\Scal,\Lcal)$-constructible.
\end{definition}
\begin{definition}[{\cite[Definition 3.2.5]{Behrend}}]\index{Stratification! Refinement}
	Let $(\Scal,\Lcal), (\Tcal,\Jcal)$ be pre-L-stratifications of $X$. We say that $(\Scal,\Lcal)$ refines $(\Tcal,\Jcal)$ if:
	\begin{itemize}
		\item $\Scal$ refines $\Tcal$, i.e., if for all $V \in \Scal$ there exists a nonempty subset $\Tcal_V \subseteq \Tcal$ for which
		\[
		V = \bigcup_{W \in \Tcal_V} W;
		\]
		\item for every stratum $V \in \Scal$, every simple object $L \in \Lcal(V)$ is trivialized by $(\Tcal,\Jcal)$ when restricted to $V$.
	\end{itemize} 
\end{definition}
\begin{remark}\label{Remark: Appendix stratifications: Refinment order}
	Note that if $X$ is a topological space, there is a set $\Mcal$ of pre-L-stratifications of $X$. We can regard $\Mcal$ as a poset by saying that for two pre-L-stratifications $(\Scal, \Lcal) \leq (\Tcal,\Jcal)$ if and only if $(\Tcal,\Jcal)$ is a refinement of $(\Scal,\Lcal)$.
\end{remark}

\begin{definition}\index{Stratification! Of a Topos}
	A stratification of a Grothendieck topos $\Ecal$ is a finite collection $\Scal$ of locally closed subtoposes $\Vcal$ of $\Ecal$ such that
	\[
	\Ecal = \coprod_{\Vcal \in \Scal} \Vcal
	\]
	and such that the closure of each subtopos $\Vcal$, $\overline{\Vcal}$, is a union of strata, i.e., there is a subset $\Scal_{\Vcal} \subseteq \Scal$ for which
	\[
	\overline{\Vcal} = \bigcup_{\Wcal \in \Scal_{\Vcal}} \Wcal.
	\]
	Note that the union above is meant to mean the category generated by the union of the $\Wcal$.
\end{definition}
\begin{remark}
	Alternatively, we can describe a stratification of $\Ecal$ by finding a stratification of $\lvert \Ecal \rvert$ and then determining from each strata $V$ of $\lvert \Ecal \rvert$ the corresponding locally closed subtoposes $\Vcal$ (cf.\@ Remark \ref{Remark: Stratifications of a topos}).
\end{remark}
\begin{remark}\label{Remark: Noetherian Topos}\index{Noetherian Topos}\index{{Topos!Noetherian}|see {Noetherian Topos}}
Occasionally in this paper we will also need Noetherian toposes. These are essentially coherent toposes which satisfy the ascending chain condition on subtoposes. The definition can be found in \cite[Definition VI.2.11]{SGA4}, although the precise definition is not necessary for this paper. The big take-away for us are the following consequences of having a Noetherian topos $\Ecal$:
\begin{itemize}
	\item The subobject lattice of $\Ecal$ satisfies the ascending chain condition (cf.\@ \cite[Exercise 7.9]{JohnstoneToposTheory}).
	\item If the topos $\Ecal$ is Noetherian, so is $\lvert X\rvert$ (and vice versa)
	\item If the scheme $X$ is Noetherian, then $\Shv(X,\text{{\'e}t})$ is Noetherian.
\end{itemize}
In particular, this tells us that a Noetherian topos admits a notion of Noetherian induction on $\lvert X \rvert$.
\end{remark}

When given a Noetherian Grothendieck topos $\Ecal$, we can mimic the theory and definition of pre-L-stratifications verbatim by defining d-structures for $\Ecal$ as d-structures on $\lvert \Ecal \rvert.$ Then in this case that running the d-structure and lcc theory through $\lvert \Ecal \rvert$ gives rise to a notion of constructiblity that is particularly well behaved on $\Ecal$ and compatible with the $\ell$-adic formalism. In particular, running this theory on a Noetherian c-topos $(\Ecal,\overline{\Ecal})$ is necessary to define $\Dbb^b_c(\Xfrak;\overline{\Q}_{\ell})$; cf.\@ \cite[Definition 4.4.14]{Behrend}.

\backmatter

\nocite{*}
\printbibliography[heading=bibintoc]

\printindex[notation]
\printindex

\end{document}